\ifdefined\XeTeXversion\else\pdfoutput=1\fi % to make arXiv happy
\documentclass[10pt]{article}
\usepackage[utf8]{inputenc} 
\usepackage[T1]{fontenc}
\usepackage{amsmath, amsfonts, amsthm, amssymb, aliascnt}
\usepackage{a4wide}
\usepackage[normalem]{ulem}
\usepackage{etoolbox}

\providecommand{\coloneqq}{\mathrel{\mathop:}=}

\usepackage{float}  % in preamble
 
\usepackage{mathrsfs} % provides \mathscr

\usepackage{accents}
\DeclareMathAlphabet\mathbfcal{OMS}{cmsy}{b}{n}

\usepackage{enumerate, graphicx, xcolor, array, booktabs, microtype}
\usepackage{fancyhdr}
\usepackage[nottoc]{tocbibind}
  
\usepackage[pdfusetitle]{hyperref}

\hypersetup{linktoc = all}
\hypersetup{hidelinks}
\hypersetup{bookmarksnumbered}

\numberwithin{equation}{section}
\numberwithin{table}{section}
\numberwithin{figure}{section}
\makeatletter
\let\c@table\c@figure
\makeatother

\usepackage{footnotebackref}
 \definecolor{cadmiumgreen}{RGB}{0,107,60}

\theoremstyle{plain}
        \newtheorem{theorem}{Theorem}[section]
\newtheorem{definition}[theorem]{Definition}
        \newtheorem{lemma}[theorem]{Lemma}
        \newtheorem{corollary}[theorem]{Corollary}

\numberwithin{equation}{section}
\DeclareMathAlphabet\mathbfcal{OMS}{cmsy}{b}{n}
\newaliascnt{proposition}{theorem}
\newtheorem{proposition}[proposition]{Proposition}
\aliascntresetthe{proposition}

\newaliascnt{lemma}{theorem}
\aliascntresetthe{lemma}

\newaliascnt{corollary}{theorem}
\aliascntresetthe{corollary}

\newaliascnt{definition}{theorem}

\aliascntresetthe{definition}

\newaliascnt{remark}{theorem}
\aliascntresetthe{remark}

\usepackage{tikz}

\newsavebox{\savepar}
\newenvironment{boxit}{\begin{lrbox}{\savepar}
  \begin{minipage}[b]{5.5in}}
  {\end{minipage}\end{lrbox}\fbox{\usebox{\savepar}}}
\newcommand \bkey {\begin{boxit} \tt  {Notre point CLEF ici} \normalsize\\}
\newcommand \ekey {\end{boxit} \newline}

\newcommand \la 	\langle
\newcommand \ra 	\rangle

\newcommand \Dcal 	{\mathcal D}

\newcommand \bei 	{\begin{itemize}}
\newcommand \eei 	{\end{itemize}}
\newcommand \del	 \partial
\newcommand \auth 	{\textsc}   
\newcommand \Mcal 	{\mathcal M}   
   
\newcommand \Lcal 	{\mathcal L}   
\newcommand \Ecal 	{\mathcal E}

\newcommand \Fb 	{\overline F}

\newcommand \RR 	{\mathbb R}   
\newcommand \eps 	\epsilon  
\newcommand \be 	{\begin{equation}}
\newcommand \ee 	{\end{equation}}  

\newcommand \lam 	\lambda

\newcommand{\WARNONCE}{\gdef\WARNONCE{}\GenericWarning{}{There are still BLF/PLF comments in this file!}}

\renewcommand \th 	\theta

\newcommand \Bcal 	{\mathcal B}

\renewcommand \geq \geqslant
\renewcommand \leq \leqslant

\newcommand \bea  {\begin{eqnarray}}
\newcommand \eea  {\end{eqnarray}}

\newcommand \Jt {\widetilde {\particlecurrent}}

\newcommand \rhob {\raug}

\newcommand \Ocal {\mathcal O} 

\newcommand \Jtilde {\widetilde {\particlecurrent}}

\newcommand \Jcausal {{\particlecurrent}} 

\newcommand \Acal {\mathcal A}
\newcommand \Gb     {\overline G}
\newcommand \Jb     {\overline {\particlecurrent}} 
\newcommand \Tb     {\overline T}

\newcommand \bse {\begin{subequations}}
\newcommand \ese {\end{subequations}}

  \definecolor{cadmiumgreen}{rgb}{0.0, 0.42, 0.24}
  
  \newcommand \cof {\textit{cof} \hskip.03cm}
  
  \newcommand \R \RR

 \newcommand \Jh {\widehat {\particlecurrent}}

\newcommand \Cscr  {\mathscr C} 
  \newcommand \bCscr {\mathscr A}
  
  \newcommand \Xbb {\mathbb X}
  
\newcommand \bF {\overline F}
\newcommand \bG {\overline G}
\newcommand \bJ {\overline {\particlecurrent}}
\newcommand \bT {\overline T}
\newcommand \bB {\overline B}

\newcommand \hJ {\widehat {\particlecurrent}}

\newcommand \bXi {\overline \Xi}
\newcommand \brho {\raug}

\usepackage{newunicodechar}
\newunicodechar{−}{\ensuremath{-}}
\makeatletter
\patchcmd{\l@section}{1.0em}{0.45em}{}{\PackageError{entropy-UN}{Contents spacing patch failed}{}}
\makeatother

\newcommand{\Frest}{{F_{\mathrm{rest}}}}
\newcommand{\Grest}{{G_{\mathrm{rest}}}}
\newcommand{\Fsp}{{F_{\mathrm{sp}}}}
\newcommand{\Gsp}{{G_{\mathrm{sp}}}}
\newcommand{\particlecurrent}{N}
\newcommand{\rphys}{{\widehat\rho}}
\newcommand{\raug}{{\rho}}
\newcommand{\rrestnew}{{\bar\rho}}
\newcommand{\Cl}{\operatorname{\mathbf{Cl}}}
\newcommand{\slicenormal}{\mathbf n}
\newcommand{\smoothoperator}{J}

\begin{document}

% Keep displayed mathematics distinct from the surrounding prose.
\setlength{\abovedisplayskip}{4pt plus 2pt minus 1pt}
\setlength{\belowdisplayskip}{4pt plus 2pt minus 1pt}
\setlength{\belowdisplayshortskip}{3.5pt plus 2pt minus 1pt}

\title{\textbf{A theory of nonlinear relativistic elastodynamics\\  with shocks, cavitation, and rank loss}}
\author{Philippe G. LeFloch\thanks{Laboratoire Jacques-Louis Lions, Centre National de la Recherche Scientifique, Sorbonne Universit\'e, Paris, France. Email: {\tt contact@philippelefloch.org}.}}
\date{November 2026}
\hypersetup{pdfauthor={Philippe G. LeFloch},pdftitle={A theory of nonlinear relativistic elastodynamics with shocks, cavitation, and rank loss}}
\maketitle

\begin{center}\vskip-.5cm
\emph{In memory of Peter D. Lax}\footnote{Peter Lax was my mentor at the Courant Institute. He invited me to join the faculty as a Courant Instructor and welcomed me back, every time very warmly, over several decades. His pioneering work on hyperbolic conservation laws and shock waves has been a lasting inspiration for my own research. I remember his generosity, encouragement, and friendship with deep gratitude.
}
\end{center}

\begin{abstract}
We study the evolution of materials in the context of nonlinear (non-)relativistic elastodynamics, allowing for shock waves and cavitation (or vacuum) phenomena; our results cover compressible fluid, as well. First, we formulate the equations in terms of spacetime deformation minors, treated as independent augmented fields, define a notion of \emph{weak solutions with cavitation and rank loss}, and prove an \emph{Eulerian--Lagrangian equivalence theorem} in this class of solutions. 
Second, we introduce the notion of \emph{relativistic polyconvexity} and, under constitutive and acoustic hypotheses, we prove that the (negative of the) spacetime particle current can be interpreted as a \emph{mathematical entropy current}, whose density is convex in the conservative variables and extends to cavitation under constitutive continuation assumptions. 
Third, we state and prove a \emph{weak--strong stability theorem}, which compares the evolution given by a strong solution to that of an entropy-admissible weak solution. Under suitable conditions on the strong solution, the stability estimate allows for both cavitation and loss of rank in the elastic fields. 
Fourth, for constitutive laws with a positive quadratic elastic response, we establish a local-in-time \emph{existence theory with cavitation} in Minkowski spacetime. Following the fluid constructions of Makino, Ukai, and Kawashima and of LeFloch and Ukai, we show that the particle flow transports initially empty cavities with zero current and stress. The ranks of $F$ and $G$ remain constant along each particle trajectory. Finally, we analyze relativistic Ogden-type constitutive laws and justify the non-relativistic limit for suitably convergent and compatible initial data. Together, these results provide a comprehensive mathematical theory for the evolution of the elastic materials and compressible perfect fluids  in both the relativistic and non-relativistic regimes.
\end{abstract}

\makeatletter
{\small\setcounter{tocdepth}{1}\let\l@part\@gobbletwo\tableofcontents}
\makeatother

\clearpage

%============================================================================
 
\let\EntropySavedJt\Jt
\renewcommand{\Jt}{\widetilde{\particlecurrent}}
\let\EntropySavedJtilde\Jtilde
\renewcommand{\Jtilde}{\widetilde{\particlecurrent}}
\let\EntropySavedJcausal\Jcausal
\renewcommand{\Jcausal}{\particlecurrent}
\let\EntropySavedJb\Jb
\renewcommand{\Jb}{\overline{\particlecurrent}}
\let\EntropySavedJh\Jh
\renewcommand{\Jh}{\widehat{\particlecurrent}}
\let\EntropySavedbJ\bJ
\renewcommand{\bJ}{\overline{\particlecurrent}}
\let\EntropySavedhJ\hJ
\renewcommand{\hJ}{\widehat{\particlecurrent}}

%===========================================================================

\section{Introduction}
\label{section=1}

\subsection{\texorpdfstring{{Elastodynamics and fluid dynamics}}{Elastodynamics and fluid dynamics}}

\paragraph{{A new formulation.}}
We study the evolution of \emph{relativistic compressible fluids} and \emph{relativistic elastic materials}, allowing for shock waves, vacuum in fluids, and cavitation in materials. On a prescribed oriented and time-oriented Lorentzian spacetime $(\mathcal M,g)$ of dimension $3+1$, we formulate the equations in terms of the material current and the spacetime minors of the configuration gradient.  They serve both the weak formulation and the local existence theory.

The system of non-relativistic elastodynamics imposes differential constraints on the deformation gradient. Dafermos was a pioneer in studying the role of these constraints~\cite{Dafermos1981Elasticity,Dafermos1986Involutions}. Using the relative entropy method that he and DiPerna introduced~\cite{Dafermos1979Stability,DiPerna}, Dafermos proved a weak--strong uniqueness theorem, assuming uniform convexity of the entropy and suitable smallness bounds on the local and far-field behavior~\cite{Dafermos2013Involutions}. The generalization of these results to relativistic evolution proposed here requires us to work with spacetime differential forms. The negative of the particle current serves as the mathematical entropy current used to compare weak and strong solutions. In addition to considering relativistic dynamics, we also address elastic flows in the presence of cavitation. 

At the heart of this study is the notion of mathematical entropy, which goes back to the work of Peter Lax on weak solutions and shock waves for hyperbolic conservation laws~\cite{Lax1954Weak,Lax1957Hyperbolic,Lax1971ShockEntropy}.
We are also motivated by the seminal work of Godunov~\cite{Godunov1961Quasilinear}, who introduced a potential formulation of a class of nonlinear hyperbolic systems, and Friedrichs and Lax~\cite{FriedrichsLax1971}, who showed how an additional conservation law with a strictly convex density yields a symmetric hyperbolic formulation. For weak solutions, the corresponding entropy inequality provides an admissibility condition across shocks.    In this paper, we exhibit this entropy structure for relativistic elastic materials, allowing for cavitation and loss of rank of the deformation gradient. 

\paragraph{The augmented fields.}

More precisely, we propose to work with the \emph{configuration gradient} $F$, its \emph{quadratic spacetime minors} $G$, and the \emph{material current} $\particlecurrent$. The \emph{constitutive stress} $T$ also appears in the equations. The fields $F^A$ are one-forms and $G_A$ are two-forms, with $A=1,2,3$; spacetime indices range from $0$ to $3$, and repeated indices are summed. We denote by $d$ the exterior derivative and by $\nabla$ the Levi--Civita connection of $g$.  The spacetime volume form $\mathrm{vol}_g$ identifies the material three-form with the particle current $\particlecurrent$. For a homogeneous material, {with a material parametrization in which the reference volume form has constant components}, the covariant equations and the entropy inequality take the form
\be
\label{eq:intro-augmented-system}
 dF^A=0,\qquad dG_A=0,\qquad
 \nabla_\nu T^{\mu\nu}=0,\qquad -\nabla_\mu {\particlecurrent}^\mu\le0.
\ee
For fields obtained from a configuration map $f$, we have $F=df$ and the quadratic and cubic minor relations $G=G(F)$ and $\particlecurrent=\particlecurrent(F)$. In the augmented formulation, we treat $F,G,\particlecurrent$ as independent fields. We require $\particlecurrent$ to be future causal or zero and impose $i_{\particlecurrent}F^A=i_{\particlecurrent}G_A=0$. Its proper particle density is $n\coloneqq\sqrt{-g(\particlecurrent,\particlecurrent)}$.

Throughout the augmented theory, the polynomial identities defining $G(F)$ and $\particlecurrent(F)$ are imposed only when the set of exact spacetime minors is explicitly specified. The spatial one-form and dual two-form matrices, and the particle density, are otherwise independent constitutive variables. Their ranks and zero sets need not satisfy the relations of a single configuration. The exterior equations, orthogonality and constitutive admissibility remain part of the augmented system. The rank-three and timelike-kernel conditions used to introduce configurations in Section~\ref{section=2} below describe the initial classical formulation; they are relaxed in the augmented constitutive domain of Section~\ref{section=3} below.

We denote the \emph{particle-number current}\footnote{This is the current denoted by $\mathbb N$ in the fluid formulation of~\cite{TorusT2-1-2026,TorusT2-2-2026}. The formulation with particle production in those papers also imposes stress conservation and $- \nabla_\mu\mathbb N^\mu\le0$.} by $\particlecurrent$, and we choose the mathematical entropy current to be $-\particlecurrent$. Observe also that the stress conservation laws in \eqref{eq:intro-augmented-system} are the matter identities required by the contracted Bianchi identity when coupling to Einstein's equations. 

\paragraph{Constitutive law and cavitation.}

To close the system of equations~\eqref{eq:intro-augmented-system}, we prescribe an energy $\raug(F,G,\particlecurrent;g)$ that determines the stress $T$ in matter. When the current vanishes, the same limiting values of $(F,G,\particlecurrent)$  can give different stress tensors, so the constitutive description also specifies which stress limits are admissible. A nonzero stress persisting at zero current is called a \emph{residual stress}; it is distinct from an empty cavity, where both the current and stress vanish. To allow for shock waves and cavitation in this formulation, we impose the equations and the entropy inequality in the sense of distributions, with prescribed initial data sets. For a compressible fluid, the energy depends only on $\particlecurrent$ once the metric is fixed, and the equations for $F$ and $G$ can be omitted.

Loss of rank does not change the meaning of $F$ and $G$: they remain tensor fields, and either may lose rank while $\particlecurrent$ is timelike. The reconstruction of velocity is a separate question. For $n>0$, the current determines $u=\particlecurrent/n$; at zero current it determines no velocity. The local existence theory uses an \emph{extended material velocity}, whose values in an empty region are specified by the continued evolution equations. A nonzero lightlike current has a determined null direction but no unit timelike velocity and no material rest frame. Such constitutive limits must therefore be distinguished from the finite-velocity continuations used in the local existence theorem.

\subsection{\texorpdfstring{{Entropy structure and weak--strong stability}}{Entropy structure and weak--strong stability}}

\paragraph{Relativistic polyconvexity and mathematical entropy.}
The central condition proposed in the present work is \emph{relativistic polyconvexity}. We extend Ball's convexity in the minors~\cite{Ball1} by requiring convexity in the \emph{material rest frame}. For $n>0$, this is the orthonormal frame moving with the material velocity $u=\particlecurrent/n$. We denote by $\Frest $ the spatial components of $F$ in that frame and by $\Grest $ the spatial duals of the two-forms $G$. The \emph{rest energy} $\bar\rho(\Frest ,\Grest ,n)$ is $\raug(F,G,\particlecurrent;g)$ expressed in these variables. Relativistic polyconvexity requires $\bar\rho$ to be convex on a convex domain of independent variables $(\Frest ,\Grest ,n)$. In the orthonormal frame of an observer, ${\particlecurrent}^0$ is the particle density that observer measures. Its negative is the mathematical entropy density. The current can vanish while the material still carries stress, so we study the extension of this density to limits of the constitutive law.
For $n>0$, the \emph{constitutive enthalpy} of the elastic material is
\be
H(\Frest,\Grest,n)\coloneqq\partial_n\bar\rho(\Frest,\Grest,n),
\ee
where the derivative is taken with both rest-frame fields $\Frest$ and
$\Grest$ held fixed. We say that the law has \emph{positive enthalpy}
on a given positive-density domain if $H>0$ throughout that domain.

To express the density in the variables used for evolution, we take a smooth time function $t$ with compact spacelike slices $\Sigma_t$ without boundary. We denote by $\slicenormal$ their future unit normal and by $\iota_t$ their inclusion in $\mathcal M$. The observer with four-velocity $\slicenormal$ measures the particle density $r\coloneqq-g(\particlecurrent,\slicenormal)$. The \emph{total particle number} on a slice is
\be
\label{eq:intro-total-particle-number}
\mathcal N(t)\coloneqq\int_{\Sigma_t}-g(\particlecurrent,\slicenormal)\,d\mu_{\Sigma_t}
=\int_{\Sigma_t}r\,d\mu_{\Sigma_t}.
\ee
{Here $d\mu_{\Sigma_t}$ is the induced Riemannian volume measure.} For integrable fluxes, the entropy inequality makes $\mathcal N$ nondecreasing, while the total mathematical entropy is $-\mathcal N$. For exact minors of a locally Lipschitz configuration, the kinematic identity conserves $\mathcal N$, including through rank loss. The \emph{conservative variables} are the spatial form components and the energy--momentum density
\be
\label{eq:intro-conservative-variables}
 U\coloneqq\bigl(\iota_t^*F,\iota_t^*G,\mathcal P\bigr),\qquad
 \mathcal P^\mu\coloneqq-T^{\mu\nu}{\slicenormal}_\nu.
\ee
{Relative to the future unit normal $\slicenormal$, we call an entropy current \emph{future-convex} when its normal entropy density is convex as a function of these conservative variables on a specified convex domain.}
{In matter, on a portion of the constitutive law where $U$ determines the fields, the balance laws yield a first-order system with spatial compatibility equations. Schematically,}
\be
\label{eq:intro-first-order-system}
 \mathsf A^\mu(U;g,\slicenormal)\nabla_\mu U=\mathcal S_{\rm geom}(U;g,\slicenormal,\nabla {\slicenormal}),\qquad
 d_{\Sigma_t}(\iota_t^*F)=d_{\Sigma_t}(\iota_t^*G)=0.
\ee
{Here $\mathsf A^\mu$ are the matrices defining the principal differential operator, $d_{\Sigma_t}$ is the exterior derivative on the slice, and $\mathcal S_{\rm geom}$ contains lower-order terms from the prescribed geometry and the chosen normal. With positive enthalpy and a positive definite inertia matrix, a positive definite Hessian of the rest energy yields a symmetrization of the evolution in an admissible time direction and a locally strictly convex entropy density. At cavitation we impose these covariant equations in the sense of distributions and use the limiting values of the stress $T$ specified by the constitutive law.}

\paragraph{Weak solutions and Eulerian--Lagrangian equivalence.}

In addition to exhibiting an entropy structure, we also generalize Wagner's treatment of the vacuum in gas dynamics~\cite{Wagner0} as well as his approach based on differential forms in classical elasticity~\cite{Wagner1,Wagner}. 

\begin{theorem}[{Entropy structure of nonlinear elastodynamics with cavitation}]
\label{thm:intro-weak-equivalence}
Fix a covariant augmented law $\raug(F,G,\particlecurrent;g)$ in the independent variables $(F,G,\particlecurrent)$, under suitable orthogonality conditions. Assume that this law is relativistically polyconvex, in the sense that the rest energy, denoted by $\bar\rho(\Frest,\Grest,n)$, is convex on a convex domain of independent rest variables $(\Frest,\Grest,n)$.
\begin{itemize}
\item[(i)] \emph{{A future-convex mathematical entropy.}}
\begin{itemize}
\item[(a)] On the chosen open domain of positive-density rest variables, assume that 
% $\bar\rho$ is $C^2$, 
the enthalpy is positive, and the inertia matrix is positive definite\footnote{{The inertia matrix is the velocity block of the rest-frame time matrix. Together with the nondegeneracy conditions on the other blocks, its positivity gives a symmetric hyperbolic formulation in an admissible time direction.}}. In each unit spatial direction of the material rest frame, require the acoustic quadratic form on velocity amplitudes to be bounded above by the quadratic form defined by the inertia matrix throughout the chosen constitutive domain. On physical configurations this means that the rest-frame wave speeds do not exceed the speed of light.
\item[(b)] Fix spacelike slices $\Sigma_t$ with future unit normal $\slicenormal$. On an open convex domain of conservative variables $U=(\iota_t^*F,\iota_t^*G,\mathcal P)$ in matter, assume that $U$ determines the observed particle density 
% $r = 
$-g(\particlecurrent,\slicenormal)$. Then $\eta(U)=g(\particlecurrent,\slicenormal)$ 
% -r$ 
is a convex mathematical entropy density with entropy current $-\particlecurrent$. At finite constitutive limits approached within this domain where $r$ has a continuous trace, the extended entropy equals the negative of that trace. The entropy admissibility condition for weak solutions with cavitation is $-\nabla_\mu {\particlecurrent}^\mu\le0$ in the sense of distributions.
\end{itemize}

\item[(ii)] \emph{Weak solutions in the Eulerian and Lagrangian descriptions.}
Fix an orientation-preserving material parametrization $(t,x)\mapsto(t,\chi(t,x))$, where $\particlecurrent^\mu\del_\mu\chi^A=0$. Assume it is a locally bi-Lipschitz homeomorphism onto its image, with the map and inverse extending locally bi-Lipschitz to the initial slice, and that the fields, fluxes and source densities are locally integrable up to that slice. The Eulerian and Lagrangian descriptions then give \emph{equivalent notions of weak solution}, with the corresponding initial data. The mathematical entropy inequality, when imposed, is preserved in both directions, including at cavitation and loss of rank of $F$.

\item[(iii)] \emph{Perfect fluids.} These statements apply to relativistic perfect compressible fluids with vacuum under  corresponding hypotheses on the constitutive law and the material parametrization. The constitutive law depends only on the particle density, so the fields $F$ and $G$ and their equations may be omitted; the conservative variable is then $U=\mathcal P$.
\end{itemize}
\end{theorem}

Since $g(\particlecurrent,\slicenormal)$ is single-valued, these local densities agree on overlaps. Integrating the Hessian along segments gives convexity on the chosen convex domain, and Proposition~\ref{thm:lsc-entropy-extension}(iv) gives the boundary value. Theorem~\ref{un-prop:weak-coordinate} proves (ii) by transforming each balance law together with its geometric source. Proposition~\ref{prop:fluid-particle-number-reduction} gives the fluid reduction in (iii). Proposition~\ref{cor:weak-flow-coordinates-vacuum} constructs material parametrizations from a conserved current with a bounded velocity whose spatial Lipschitz constant is integrable in time, including where the current is null or zero. Section~\ref{sec:material-flux-measures} treats two separate test traces for collapses covered by its geometric and flux hypotheses.

For a locally Lipschitz configuration $f$, the exact minors of $df$ satisfy the closure equations and conserve the particle current. The first-order weak formulation therefore reduces to the stress balance laws, with $T$ chosen among the admissible constitutive  stress tensors. For configurations in $W^{2,1}_{\rm loc}\cap W^{1,\infty}_{\rm loc}$ with compact essential constitutive ranges and material momenta in $W^{1,1}_{\rm loc}$, {the Noether identity for the variational stress of the physical law gives equivalence} with the variational equations if the gradient has rank three almost everywhere (cf.~Proposition~\ref{prop:variational-EL-local04} below).

\paragraph{Weak--strong stability with cavitation.}

For weak--strong comparison, it is necessary to define the \emph{strong solution}, which must satisfy the balance laws and entropy equality with the regularity specified below. At this stage, no propagation law for an extended velocity in an empty cavity is needed. On the other hand, for the local existence theory below, we will introduce the notion of \emph{strong solutions with prescribed continuation}: they satisfy, in addition, the continued evolution system, including $\mathrm D_tV=0$ where the particle current and stress vanish. Such a solution serves as a strong reference in our approach below. 

Let $U_*$ denote the conservative variables of the strong solution. A finite supporting\footnote{A supporting covector gives an affine lower bound for the convex entropy, with equality at $U_*$. In the differentiable interior it equals $D\eta(U_*)$. At a corner there may be several choices, and a finite choice at the boundary is made.}  
covector $\lambda_*\in\del\eta(U_*)$ defines the \emph{relative entropy density} $D_\eta(U\mid U_*)$ as the difference between the entropy and its \emph{supporting affine function} at $U_*$. We write $\mathcal E(t)$ for its integral over the slice.
\bse
\be
\label{eq:intro-relative-entropy-density}
 {D_\eta(U\mid U_*)}\coloneqq\eta(U)-\eta(U_*)-\lambda_*\cdot(U-U_*)\ge0,
\ee
\be
\label{eq:intro-relative-entropy-integral}
 \mathcal E(t)\coloneqq\int_{\Sigma_t}D_\eta(U\mid U_*)\,d\mu_{\Sigma_t},
\ee
\ese
The value at $t=0$ uses the prescribed initial data sets.

{Write $Z=(F,G,\particlecurrent,T)$ for the fields and their stress, and $Z_*$ for the corresponding quantities of the strong solution. Pairing the differences of their balance and spatial compatibility laws with the reference coefficients gives a relative current with density $D_\eta(U\mid U_*)$ and divergence remainder $\mathcal R_*$ relative to $dt\,d\mu_{\Sigma_t}$. In coordinates adapted to the slices, Section~\ref{subsec:physical-weak-strong-vacuum} uses the densities $\mathscr Q^0=\sqrt{\det h_t}\,D_\eta$ and $\mathscr R_*=\sqrt{\det h_t}\,\mathcal R_*$, where $h_t$ is the induced spatial metric. Integration on the slices gives}
\be
 \frac{d}{dt}\mathcal E(t)
 \le\int_{\Sigma_t}\mathcal R_*(t,x;Z(t,x))\,d\mu_{\Sigma_t}
 \qquad\text{in the sense of distributions in time}.
\ee
We use the relative entropy method of Dafermos~\cite{Dafermos1979Stability} and DiPerna~\cite{DiPerna}, and specifically extend the argument for relativistic fluids with vacuum given in~\cite{LeFloch-JHDE}. Bounding the remainder by an integrable multiple of $D_\eta$ yields our main estimate. 

\paragraph{{The comparison theorem.}}
{We now specify what the comparison requires near cavitation.}
\begin{theorem}[Weak--strong stability for nonlinear elastodynamics with cavitation]
\label{thm:intro-main-results}
On a smooth periodic Cauchy slab with spacelike slices of future unit normal $\slicenormal$, fix a relativistically polyconvex law. Assume that its mathematical entropy extends as a lower semicontinuous convex function of $U$, finite and equal to $g(\particlecurrent,\slicenormal)$ for the unknowns $(F,G,\particlecurrent)$ and the stress $T$ given by the constitutive law and its specified continuation at cavitation, including the values prescribed by the initial data. Let $(F,G,\particlecurrent)$ be an entropy weak solution with integrable fluxes and prescribed initial terms, with (exact or) independent minors. One then write $Z=(F,G,\particlecurrent,T)$ for these fields together with their constitutive stress. Let $Z_*$ and its constitutive parameters be bounded and continuous up to the initial slice. One also requires $Z_*$ to satisfy the balance laws and entropy equality with the initial data set given by its values at $t=0$, and imposes the spatial compatibility equations for both solutions.
\begin{itemize}

\item[(i)] Choose bounded Lipschitz coefficients on the closed slab for the spatial compatibility equations and a bounded Lipschitz covector $\lambda_*$ for the entropy. At almost every point it gives a finite support at $U_*$ on the comparison set. On the initial slice the trace of $\lambda_*$ gives a support for the prescribed comparison values.

\item[(ii)] {Assume that the remainder $\mathcal R_*$ along the comparison solution satisfies}
\be
\label{eq:intro-weighted-remainder}
 \mathcal R_*(t,x;Z(t,x))\le a(t)\,D_\eta\bigl(U(Z(t,x))\mid U_*(t,x)\bigr)
\ee
{almost everywhere for some $a\in L^1(0,T)$. }
\item[(iii)] The initial relative entropy satisfies $\mathcal E(0)<+\infty$.
\end{itemize}
Then, for almost every $t$, one has
\be
\label{eq:intro-relative-entropy-bound}
 \mathcal E(t)\le \exp\left(\int_0^t a(s)\,ds\right)\mathcal E(0).
\ee
If $\mathcal E(0)=0$, then $D_\eta(U\mid U_*)=0$ almost everywhere. If zero relative entropy determines the three unknown fields $F,G,\particlecurrent$ and their associated stress $T$ uniquely, equal initial data sets imply equality of these fields almost everywhere. The estimate holds on the whole slice, including where the cavities of the two solutions differ, without a positive lower bound on density or uniform convexity of the entropy near the constitutive boundary.
\end{theorem}

\paragraph{A constitutive example for existence and stability.}

{Our main concrete example combines local existence and weak--strong stability for independent augmented fields. In units where the speed of light is one, fix $m,\kappa,\alpha,\beta>0$ and $\gamma>1$, and consider the quadratic elastic law
\be
\label{eq:intro-complete-quadratic-law}
 \bar\rho_{\rm q2}(\Frest,\Grest,n)
 ={m \, n}+\frac{\kappa}{\gamma}n^\gamma
 +\frac12 n^{2\gamma-1}
       \bigl(\alpha|\Frest|^2+\beta|\Grest|^2\bigr).
\ee
{This law extends the polytropic fluid energy by a quadratic elastic response in the independent rest-frame fields. When the elastic terms are omitted, its pressure is $p(n)=\kappa(\gamma-1)n^\gamma/\gamma$ and its squared sound speed is $\kappa(\gamma-1)n^{\gamma-1}/(m+\kappa n^{\gamma-1})$, which tends to zero at vacuum. The exponent $\gamma>1$ describes the density dependence; the quadratic degree concerns the elastic variables.}

{Theorem~\ref{thm:quadratic-second-weight-augmented-stability} below gives local existence with prescribed continuation on $\mathbb T^3$ or $\mathbb R^3$, and weak--strong comparison on a periodic slab, under its constitutive and acoustic bounds. The fields $F,G,\particlecurrent$ are independent and the two cavities may differ. The comparison includes a nonnegative initial entropy defect, which vanishes when the initial data lie in the same constitutive branch. A strong $L^1$ initial trace of the conservative variables also gives the estimate without the defect.}

The positive quadratic response is retained in Corollary~\ref{cor:nonlinear-second-weight} below, which adds nonlinear elastic terms to this complete augmented result. Two other extensions require different comparison hypotheses:
\begin{itemize}
\item For the broader nonlinear family introduced next, Theorem~\ref{uc:thm:integrated-3d} below estimates relative entropy together with the displacement of {corresponding material points} for every $\gamma>1$. It assumes exact spacetime minor relations and common affine periods. A pointwise estimate is available under the additional weighted divergence bound on the reference form multipliers in Theorem~\ref{uc:thm:common} below.
\item For Ogden-type laws without a vanishing density factor in the elastic terms, Theorem~\ref{og:theorem} below uses exact or mixed minor relations and {the required regularity of} the reference constitutive derivatives. Corollary~\ref{cor:Ogden-fully-augmented-stability} below permits independent augmented fields under its adjoint identity and dual Sobolev regularity assumptions.
\end{itemize}
}

\subsection{\texorpdfstring{{Local existence and applications}}{Local existence and applications}}

\paragraph{ Local existence with possible cavitation and rank loss.}
For local existence in Minkowski space, we consider \emph{elastic constitutive laws}. We write $\bar\rho_{\varepsilon_{\rm e}}$ for the rest energy, with elastic parameter $\varepsilon_{\rm e}\ge0$, rest mass per particle $m>0$ and light speed $c_{\rm light}$. {When the material volume form has normalized components,} an exact configuration satisfies $\Grest=\operatorname{cof}\Frest$ and $n=\det\Frest$ in the rest frame. In specifying the augmented law we use $(\Frest,\Grest,n)$ as independent variables.

We put $\theta\coloneqq(\gamma-1)/2>0$ and $q\coloneqq n^\theta$. The density-weighted constitutive family has energy
\be
\label{eq:intro-general-local-energy}
 \bar\rho_{\varepsilon_{\rm e}}(\Frest ,\Grest ,n)
 =n\bigl(mc_{\rm light}^2+\Phi(q)+\Psi(q,z)\bigr),
 \qquad z=\varepsilon_{\rm e}q(\Frest ,\Grest ).
\ee
The functions $\Phi(q)=q^2 A(q)$ and $\Psi(q,z)$ are $C^9$ on neighborhoods of the ranges under consideration. The function $\Psi$ is invariant under simultaneous spatial rotations of the rest vectors and satisfies $\Psi(q,0)=0$ and $\Psi_z(q,0)=0$. At $\varepsilon_{\rm e}=0$ the energy depends only on the particle density $n$ and defines a perfect fluid.

For the evolution on a Minkowski time slice, we use the \emph{regularized density} $r^\theta$, with $r$ the particle density measured on the slice. We call $\varepsilon_{\rm e}(\Fsp,\Gsp)$ the \emph{scaled elastic fields}, where $\Fsp$ is the spatial one-form matrix and $\Gsp$ is the matrix of spatial dual two-forms. {On the set of exact spacetime minors, $(\Fsp,\Gsp)=(D_xf,\operatorname{cof}D_xf)$.} The variables estimated at the highest Sobolev order are the regularized density, the velocity and the product of $r^\theta$ with the scaled minors. We write $C_tH^s$ for continuity in time with values in $H^s$ and $C_t^1H^{s-1}$ for continuous differentiability in time with values in $H^{s-1}$. For the fluid energy in units where the light speed is one, the enthalpy is $h=m+\Phi+\theta q\Phi'$ and the coefficient of compression in the regularized equations is $H_f=\Phi''+(1+1/\theta)\Phi'/q$, continued at $q=0$. The squared sound speed is $\theta^2q^2H_f/h$.

\begin{theorem}[Local existence for nonlinear elastodynamics with cavitation]
\label{thm:intro-local-evolution}
{Write $Q=r^\theta$ for the regularized observed density and $V_0$ for the prescribed initial velocity in units with light speed one.} Consider a constitutive law from the family above. On a bounded density interval starting at zero, assume positive lower bounds for $h$, $H_f$ and the elastic Hessian $\Psi_{zz}(q,0)$, and suppose that the fluid sound speed is at most the speed of light. On $\mathbb T^3$ or $\R^3$, prescribe a nonnegative regularized density and a velocity in $H^3$, with strict margins from the upper density bound and the speed of light. The density may vanish on open sets. For positive elastic parameter, prescribe independent scaled one-form and two-form fields satisfying their spatial closure constraints. No determinant or cofactor relation with the initial density is imposed. Require the scaled minors to be sufficiently small in supremum norm and to belong to $H^3$ on the torus or $W^{3,\infty}$ on $\R^3$.
\begin{itemize}
\item[(i)] The equations for the regularized variables extend continuously to zero density. For the prescribed initial data, one obtains a unique local solution of the normalized equations with prescribed continuation, with regularized density, velocity and density-weighted scaled minors in $C_tH^3\cap C_t^1H^2$. The increment of the scaled minors from their initial values belongs to $C_tH^2\cap C_t^1H^1$. The existence interval depends on the data bounds and strict margins, uniformly in the elastic parameter and down to zero density.
\item[(ii)] For positive elastic parameter, the particle flow transports the initial zero-density set and keeps the ranks of $F$ and $G$ constant along each particle trajectory. The current and stress vanish on the transported set. {For each $y$ in the initial zero-density set, the particle flow satisfies $X(t,y)=y+tV_0(y)$.}
\item[(iii)] At zero elastic parameter, one obtains the relativistic perfect-fluid equations in the regularized variables. The same regularity and local uniqueness hold for the prescribed density and velocity. The velocity is constant along the trajectories starting in vacuum.
\item[(iv)] {If $\varepsilon_{\rm e}>0$, $\Fsp_0=D_xf_0$ for an initial configuration $f_0$, and the initial data satisfy the exact minor relations}, the flow gives $f=f_0\circ X^{-1}$ and preserves $F=df$, $G=G(F)$ and $\particlecurrent=\particlecurrent(F)$. The two ranks and the density then satisfy the additional relations of these minors.
\end{itemize}
\end{theorem}
{Theorem~\ref{uf:thm:main} gives the precise constitutive and data bounds. Its proof is in Sections~\ref{uf:subsec:normalized-equations} to~\ref{uf:subsec:dependence-reconstruction}, from the positive energy to reconstruction at zero density.} For fixed positive elastic parameter, Proposition~\ref{uf:cor:fractional-existence} below gives local existence with the same transport of the cavity and minor ranks under weaker initial regularity. For $5/2<s\le3$, the initial regularized density, velocity and density-weighted scaled minors belong to $H^s$. The initial scaled minors without the density factor belong to $H^{s-1}$ on the torus or to uniformly local $H^{s-1}$ on $\R^3$. In particular the latter condition allows a nonintegrable background on $\R^3$.

\paragraph{{An Ogden-type example.}}
{
The following example belongs to~\eqref{eq:intro-general-local-energy}. Its integrated comparison uses exact spacetime minors, as in Theorem~\ref{uc:thm:integrated-3d}.
}

An explicit isotropic example of Ogden type is
\be
\label{eq:intro-local-energy}
\begin{aligned}
 \bar\rho_{\varepsilon_{\rm e}}(\Frest ,\Grest ,n)
 ={}&mc_{\rm light}^2n+
 n^\gamma\left(a_{\rm fl}+\tfrac{\varepsilon_{\rm e}^2}{2}
 (\alpha|\Frest |^2+\beta|\Grest |^2)\right)\\
 &+\tfrac{\varepsilon_{\rm e}^4}{4}n^{2\gamma-1}
 \bigl(\alpha_4 I_4(\Frest )+\beta_4 I_4(\Grest )\bigr),
 \qquad \gamma>1.
\end{aligned}
\ee
{The coefficients $a_{\rm fl},\alpha,\beta$ are positive and $\alpha_4,\beta_4\ge0$. Here $I_4(M)=\operatorname{tr}((M^TM)^2)=\sum_i\sigma_i(M)^4$, where the \emph{singular values} $\sigma_1(M)\ge\sigma_2(M)\ge\sigma_3(M)\ge0$ are the square roots of the eigenvalues of $M^TM$, counted with multiplicity. Taking $\Phi(q)=a_{\rm fl}q^2$ and $\Psi(q,z)=(\alpha|z_F|^2+\beta|z_G|^2)/2+(\alpha_4I_4(z_F)+\beta_4I_4(z_G))/4$ recovers~\eqref{eq:intro-local-energy} from~\eqref{eq:intro-general-local-energy}. In units where the light speed is one, $h=m+(1+2\theta)a_{\rm fl}q^2$, $H_f=2a_{\rm fl}(2+1/\theta)$ and $\Psi_{zz}(q,0)=\operatorname{diag}(\alpha I,\beta I)$. The required margins therefore hold on a sufficiently small density and scaled-form range for every $\gamma>1$, including at rank loss.

{
The current and stress have order $n$ on bounded subluminal ranges, with internal stress of order $n^\gamma$ (see~\eqref{uf:eq:stress-decomposition}). The data in~\eqref{uf:eq:moving-vacuum-data}--\eqref{uf:eq:moving-vacuum-flow} give a moving rank-two cavity.
}}

\paragraph{Relation to classical elastic laws.}
At unit volume the singular values of the inverse deformation gradient are the reciprocals of the principal stretches, and those of its cofactor are the stretches themselves. Quadratic powers give the Mooney--Rivlin invariants, while other powers give the corresponding terms in Ogden's energy~\cite{Ogden1972,DestradeDorfmannSaccomandi2022}. {The law~\eqref{eq:intro-local-energy} combines a positive quadratic elastic response with quartic Ogden-type terms and satisfies the local existence hypotheses. Corollary~\ref{cor:nonlinear-second-weight} below also treats other singular-value powers under {its regularity conditions}, retaining the positive quadratic response.}

{Without a vanishing density factor, the elastic energy and stress may remain nonzero at zero current. Section~\ref{subsec:spectral-elastic-energies} treats these Ogden-type laws for elastic and density exponents in $(1,2]$. In particular, when the power of $\Grest$ is below two, Theorem~\ref{og:theorem} requires the spatial $L^\infty$ norm of the divergence of the reference form multiplier to be integrable in time.}

\paragraph{{Classical local existence and the non-relativistic limit.}}
{We next give the limit of the relativistic solutions and draw a local existence conclusion for the classical equations.}
\begin{corollary}[\texorpdfstring{{Non-relativistic limit with cavitation}}{Non-relativistic limit with cavitation}]
\label{cor:intro-non-relativistic-limit}
On $\mathbb T^3$ or $\R^3$, fix the constitutive functions and elastic parameter of the density-weighted family as the speed of light tends to infinity. Assume the positive lower bounds for $H_f$ and $\Psi_{zz}(q,0)$ on the chosen density range. Use compatible initial data with uniformly bounded physical velocities, a common strict margin from the upper density bound, and uniformly small scaled elastic fields. Fix $5/2<s\le3$. Assume that the initial regularized density and physical velocity converge in $H^s$, and that the initial scaled elastic fields are bounded in $W^{3,\infty}$ and converge in $L^\infty$. At $s=3$, require $H^3$ convergence of the product of the initial regularized density and the initial scaled elastic minors.
\begin{itemize}
\item[(i)] {On a common local time interval, the regularized density, weighted elastic fields and velocity converge in $C_tH^s\cap C_t^1H^{s-1}$, including at zero density. The changes in the scaled elastic fields from their initial values converge in $C_tH^{s-1}\cap C_t^1H^{s-2}$. {The limit satisfies the classical augmented equations with cavitation. For positive elastic parameter and exact initial minors it satisfies the classical material equations. For $\varepsilon_{\rm e}=0$ it satisfies the compressible Euler equations with vacuum.} At one lower Sobolev order the rate is controlled by the initial differences, including the difference of the unweighted minors in $W^{2,\infty}$ on $\R^3$ or $W^{1,\infty}$ on $\mathbb T^3$, and by $c_{\rm light}^{-2}$.}
\item[(ii)] {The particle current in physical time $({\particlecurrent}^0,c_{\rm light}{\particlecurrent}^{\rm sp})$, the momentum density and the spatial stress converge uniformly. The energy density and flux converge uniformly after subtraction of the rest energy contribution.}
\end{itemize}
\end{corollary}

{Theorem~\ref{nr:thm:vacuum-limit} gives the rate~\eqref{nr:eq:lower-convergence}. Sections~\ref{nr:subsec:uniform-estimates} to~\ref{nr:subsec:physical-fluxes} prove the uniform existence, convergence at the initial Sobolev order, and convergence of the physical fluxes, respectively. Proposition~\ref{nr:prop:two-bounded-derivatives} gives the variant for exact minors on $\mathbb R^3$ with two bounded derivatives of the unweighted initial fields and specifies the endpoint conditions.}

\begin{corollary}[Local existence for classical elastodynamics with cavitation]
\label{cor:intro-classical-existence}
{Fix a positive elastic parameter and a law from the preceding corollary. Prescribe one compatible initial data set in the same Sobolev spaces with the same strict range bounds. The classical normalized equations admit a unique local strong solution with prescribed continuation at zero density. The particle flow transports the initially empty cavities, where the current and stress vanish, and preserves the ranks of the elastic fields. For exact initial minors, the solution also satisfies the classical material equations.}
\end{corollary}
{Corollary~\ref{cor:classical-existence-cavitation} proves this result by using the same spatial initial fields and observed density for all sufficiently large light speeds.}

\paragraph{ Relativistic and non-relativistic perfect fluids with vacuum.}
{When the elastic terms are omitted, the constitutive energy depends on the particle density alone. The four stress balance laws and the particle-number inequality then give the relativistic perfect-fluid formulation.}

\begin{corollary}[ Relativistic and non-relativistic perfect fluids with vacuum]
\label{cor:intro-perfect-fluids-vacuum}
{ Fix $m>0$, $\gamma>1$, $\theta=(\gamma-1)/2$, and a density energy from the family of Theorem~\ref{thm:intro-local-evolution}, with the elastic parameter equal to zero:
\be
\label{eq:intro-fluid-energy-pressure}
 \bar\rho_{\rm fl}(n)=mc_{\rm light}^2n+n\Phi(n^\theta),\qquad
 p(n)=\theta n^{1+\theta}\Phi'(n^\theta),\qquad
 \Phi(q)=q^2A(q).
\ee
Assume its constitutive regularity and positive lower bounds for $h$ and $H_f$ on the prescribed bounded density range. In the relativistic case impose the sound-speed bound and strict velocity margin of that theorem. Then the following conclusions hold.}
\begin{itemize}
\item[(i)] {Under the entropy, constitutive continuation and material-parametrization hypotheses in Theorem~\ref{thm:intro-weak-equivalence}, the relativistic weak formulation, mathematical entropy construction and Eulerian--Lagrangian equivalence apply with $U=\mathcal P$. The entropy construction uses a convex conservative domain, and the change of description a locally bi-Lipschitz material parametrization. Under the hypotheses of Theorem~\ref{thm:intro-main-results}, including a finite supporting covector for the strong reference and control of the remainder, its relative entropy estimate applies with vacuum.}
\item[(ii)] { Nonnegative regularized initial density and initial velocity in $H^3(\mathbb T^3)$ or $H^3(\mathbb R^3)$, with the strict range margins above, give a unique local strong solution with prescribed continuation. This conclusion holds for both the relativistic and classical fluid equations. The particle flow transports the initial zero-density set, and the velocity is constant along trajectories starting there.}
\item[(iii)] { For initial regularized observed densities $Q_0$ and physical velocities converging in $H^3$ with uniform bounds and strict margins, the relativistic solutions converge to the classical solution in $C_tH^3\cap C_t^1H^2$ on a local interval independent of $c_{\rm light}$ as $c_{\rm light}\to\infty$. The classical equations, with $n=r$ in this limit, are
\be
\label{eq:intro-classical-perfect-fluid}
 \del_t n+\operatorname{div}(nv)=0,\qquad
 \del_t({m \, n}v)+\operatorname{div}({m \, n}v\otimes v)+\nabla p(n)=0.
\ee
The regularized density and velocity require no additional derivative at vacuum.}
\end{itemize}
\end{corollary}
{ The pressure formula follows from $p=n\bar\rho_{\rm fl}'(n)-\bar\rho_{\rm fl}(n)$; the rest-mass term cancels. For $A(q)=a_{\rm fl}>0$, it gives $p(n)=a_{\rm fl}(\gamma-1)n^\gamma$, including the polytropic range $1<\gamma\le2$.}
\begin{proof}
{We set $\varepsilon_{\rm e}=0$ and omit the form rows in the constitutive and evolution systems. Proposition~\ref{prop:fluid-particle-number-reduction} below gives the reduced weak formulation; the entropy construction and comparison use {the hypotheses of (i)}. To apply the existence theorems through the elastic family, one may choose the auxiliary function $\Psi(q,z)=|z|^2/2$ in~\eqref{uf:eq:law}; it has the required positive $z$-Hessian and disappears from the fluid law at $\varepsilon_{\rm e}=0$. Theorem~\ref{uf:thm:main} below gives the relativistic local solution, and Theorem~\ref{nr:thm:vacuum-limit} below gives the classical solution and convergence in the regularized variables. {The flow formulas give the transport of vacuum and the prescribed velocity evolution there. In the fluid reduction $Y=Z=0$ and the specific stress in~\eqref{nr:eq:specific-stress} is $\theta Q\Phi'(Q)I$. Thus its contribution $r\theta Q\Phi'(Q)I$ in~\eqref{nr:eq:classical-fluxes} gives the pressure in~\eqref{eq:intro-classical-perfect-fluid}, with $r=Q^{1/\theta}=n$ at $\delta=0$.}}
\end{proof}

%-----------------------------------------------------------------------------------------------------------------------------------------

%-----------------------------------------------------------------------------------------------------------------------------------------

\subsection{Related results in the literature}

\paragraph{{Compressible fluids with vacuum.}}
{
The classical construction of Makino, Ukai, and Kawashima~\cite{MakinoUkaiKawashima1986} and the relativistic construction of LeFloch and Ukai~\cite{LeFlochU} provide the vacuum continuations used in the fluid reduction. The relativistic construction assumes bounds on the constitutive law and the velocity. The geometric comparison with vacuum in~\cite{LeFloch-JHDE} underlies our stability argument. Frauendiener~\cite{Frauendiener} and Walton~\cite{Walton} gave symmetric hyperbolic formulations of the relativistic fluid equations.
}

\paragraph{{Elastodynamics and related models.}}

For classical elastodynamics, Dafermos and Hrusa~\cite{DafermosHrusa1985} developed energy estimates and local existence, and Wagner~\cite{Wagner} gave a symmetric hyperbolic formulation and proved local existence for strictly polyconvex laws at positive determinant. {Augmented entropy formulations appear in  LeFloch~\cite{LeFloch-notes} and Qin~\cite{Qin}. The work of Demoulini, Stuart, and Tzavaras~\cite{DST0} introduces variational approximations to construct solutions described by Young measures~\cite{DST0} and proved weak--strong uniqueness for dissipative solutions described by Young measures and concentration measures~\cite{DST2012}. Brenier~\cite{Brenier2004} and Serre~\cite{Serre2004} developed augmented entropy structures for the Born--Infeld and nonlinear Maxwell systems. 

Our analysis starts from the geometric formulation of relativistic elasticity introduced by Carter and Quintana~\cite{CQ} and by Beig and Schmidt~\cite{BeigSchmidt} and Beig~\cite{Beig2023}. Beig and Schmidt obtained local existence (without allowing for cavitation) near unstrained configurations by symmetric hyperbolicity. We also refer to the work of Tahvildar-Zadeh \cite{Shadi} 
on the relativistic and non-relativistic elastodynamics with small shear strains. 
Christodoulou's book~\cite{Christodoulou000} provides a general variational and hyperbolic viewpoint relevant to continuum mechanics. Dependence on the particle current alone in affine spacetime gives the class introduced by Beig and LeFloch~\cite{BeigLeFloch}; dependence on the proper particle density gives the perfect-fluid equations. 

We also derive the balance laws and the mathematical entropy inequality across a given \emph{discontinuity} in the formalism of this paper. For a Lipschitz configuration the normal particle flux is continuous across the interface. The formation and subsequent evolution of these fronts remain open for general elastic constitutive laws in three dimensions. For fluids, shock formation is established for relativistic Euler flow by Christodoulou~\cite{Christodoulou2007Shocks} and, for classical flow with vorticity, by Buckmaster, Shkoller, and Vicol~\cite{BuckmasterShkollerVicol2023}. In two space dimensions Shkoller and Vicol~\cite{ShkollerVicol2024} describe the boundary of a classical development near the first singularity, and Buckmaster, Drivas, Shkoller, and Vicol~\cite{BuckmasterDrivasShkollerVicol2022} construct an entropy-admissible shock after the first singularity in azimuthal symmetry. It would be interesting to generalize the methods and results therein to the model of elastodynamics proposed in the present paper.

%-----------------------------------------------------------------------------------------------------------------------------------------

\subsection{Outline of this paper}

\paragraph{Part I: The entropy structure of relativistic elastodynamics.}
{Sections~2--6 proceed from the field equations to their mathematical entropy.
\begin{itemize}
\item {Section~\ref{section=2} derives the equations for the configuration gradient. Section~\ref{section=3} introduces the augmented balance laws and the propagation of the minor relations. Section~\ref{section=4-hyperbolicity} studies the characteristic directions and wave speeds of the two formulations.}
\item {Section~\ref{section=4-polyconvexity} introduces relativistic polyconvexity and compares it with Ball's condition. Section~\ref{sec:convex-entropy-relativistic-elastodynamics} proves convexity of the mathematical entropy density and extends it to the constitutive boundary for the laws of Theorem~\ref{thm:closed-vacuum-entropy}.}
\end{itemize}}

\paragraph{Part II: Weak--strong stability and local-in-time existence theory.}

\bei 

\item 
{Section~\ref{sec:weak-local04} defines weak solutions with cavitation, proves Eulerian--Lagrangian equivalence under the hypotheses of Theorem~\ref{un-prop:weak-coordinate}, and derives the jump conditions. Section~\ref{section=6} proves weak--strong stability under a bound on the remainder along the two solutions.}

\item {Section~\ref{sec:uniform-elastic-fluid-existence} states local existence and comparison results. Theorem~\ref{thm:quadratic-second-weight-augmented-stability} treats independent augmented fields. Theorem~\ref{uc:thm:integrated-3d} treats the broader family with exact minors and common affine periods. Section~\ref{sec:proofs-local-existence} gives the proofs.}
\eei

\paragraph{{Part III: Applications and technical details.}}

\bei

\item 
{Section~\ref{sec:non-relativistic-vacuum} proves the non-relativistic limit and gives classical local existence with cavitation (Corollary~\ref{cor:classical-existence-cavitation}). Section~\ref{sec:constitutive-comparison-applications} relates the Ogden-type laws to classical elastic energies and verifies the comparison hypotheses.}

\item {Section~\ref{Append-AA} derives the variational identities for the field equations. Section~\ref{app:integrated-3d} proves the estimates for the corrected relative entropy used in integrated stability. Section~\ref{subsec:constitutive-stability-verification} derives constitutive identities for the relative current and determines the stress and supporting covectors at zero current.}

\eei

%============================================================================

\part[\texorpdfstring{The entropy structure of relativistic elastodynamics}{The entropy structure of relativistic elastodynamics}]
  {\LARGE{The entropy structure of relativistic elastodynamics}}
\label{part:entropy-structure}

%============================================================================

\section{\texorpdfstring{A first-order formulation based on the configuration gradient}{A first-order formulation based on the configuration gradient}}
\label{section=2}

%-----------------------------------------------------------------------------------------------------------------------------------------

\subsection{Admissible configuration gradients}

\paragraph{Configuration map and its gradient.}
{
We now describe the admissible configuration gradients, including at zero particle density, and derive the balance laws for homogeneous materials. The calculations in this section assume sufficiently regular fields unless a different regularity is specified. The weak formulation will be presented later on in Section~\ref{sec:weak-local04}.\par}

Let $(\Mcal,g)$ be an oriented and time-oriented, $(3+1)$-dimensional Lorentzian manifold of signature $({-}{+}{+}{+})$. Throughout, we assume that the metric $g$ is sufficiently regular. We write $x=(x^\mu)$ for local spacetime coordinates. We use Einstein's summation convention over repeated indices: Greek spacetime indices $\mu,\nu,\ldots$ range from $0$ to $3$, spatial indices $i,j,k$ range from $1$ to $3$, and material indices $A,B,C$ range from $1$ to $3$. Lowercase indices $a,b,c$ have this same range when they denote material components.  Sums over antisymmetric spacetime indices use the full range unless independent components are explicitly chosen.  The scalar product of two tangent vectors $U,V$ is denoted by $g(U,V) = g_{\mu\nu} U^\mu V^\nu = U^\mu V_\mu$, where indices are lowered and raised using the metric $(g_{\mu\nu})$ and its inverse $(g^{\mu\nu})$, respectively.

On this background spacetime $(\Mcal,g)$, we describe the evolution of an elastic body by a \emph{configuration map}
\be
f: \Mcal \to \Bcal, \qquad  x \mapsto X= f(x) = (f^A(x)),  
\ee
where $X$ lies in the \emph{material manifold} $\Bcal$ ---a prescribed three-dimensional manifold representing the reference configuration of the body, endowed with a \emph{material volume form} $\Omega$ with components $\Omega_{ABC}$. The component balance laws for the independent augmented fields below use fixed volume coordinates on the material body, so $\Omega_{ABC}$ is constant in that chart. Such coordinates exist locally: if $\Omega=\omega(X)\,dX^1\wedge dX^2\wedge dX^3$ with $\omega>0$, we set $Y^1\coloneqq\int^{X^1}\omega(s,X^2,X^3)\,ds$, $Y^2\coloneqq X^2$, $Y^3\coloneqq X^3$; their Jacobian is $\omega$ and $\Omega=dY^1\wedge dY^2\wedge dY^3$. The homogeneous global coordinate results use fixed affine material coordinates, allowing periodic identifications with translational transition maps. A periodic gradient can come from a configuration on the universal cover that is the sum of an affine map and a periodic map. Under a nonlinear change of material coordinates, the equations for $G_A$ transform as bundle equations. The exact top-form current identity is valid for every smooth material volume form.

We call the differential of $f$ the \emph{configuration gradient} and write
\be
F \coloneqq  df, \qquad F^A_\mu \coloneqq  \del_\mu f^A.
\ee
On a spatial slice where the configuration map is invertible, the deformation map in material coordinates is its inverse and the two spatial gradients are inverse matrices (cf.~Section~\ref{sec:weak-coordinate-equivalence} below). At each point $x\in\Mcal$, it defines a linear map $F:T_x\Mcal\to T_{f(x)}\Bcal$. We first restrict the configuration gradient to the following rank-three, timelike-kernel regime, which permits the classical constitutive calculation. These restrictions define this part of the configuration domain, rather than the full class of augmented solutions. The extension below and the augmented formulation in Section~\ref{section=3} below allow lower ranks and the other constitutive limits specified there.
\begin{itemize}

\item[(A-1)]\underline{Maximal rank condition}: the linear map $df=F:T_x\Mcal\to T_{f(x)}\Bcal$ has maximal rank, that is, three.

\item[(A-2)] \underline{Timelike kernel condition}: the one-dimensional kernel $\ker F\subset T_x\Mcal$ is timelike with respect to the Lorentzian metric $g$.

\end{itemize}

\noindent We denote by $\Xbb$ the \emph{positive-density kinematic domain}: the pointwise open set of rank-three maps with timelike kernel, in the orientation component for which the current~\eqref{vmu} below is future directed. Its closure $\overline{\Xbb}$ contains rank-three gradients with a future null current. A gradient of rank at most two belongs to this closure precisely when its kernel contains a nonzero future causal vector (cf.~Section~\ref{section=3} below). The constitutive law separately specifies which gradients admit finite stress limits.

For $F\in\Xbb$, the fibres of $f$ are timelike curves representing material worldlines. Conditions (A-1)--(A-2) are used in the differential calculations in matter. The identity $F=df$ gives exterior compatibility at every rank. Global invertibility and non-interpenetration are separate from these local rank and kernel conditions. Theorem~\ref{un-prop:weak-coordinate} below uses an orientation-preserving locally bi-Lipschitz material parametrization.

%--------------------------------- 

\paragraph{Energy density over the material bundle.}
The equations of motion are derived from the least action principle applied to the \emph{action functional}
\be
\label{Acal}
\Acal(f) \coloneqq  \int_\Mcal \rphys[f, df; g] \, dV_g,
\ee
where $dV_g$ denotes the volume form associated with the spacetime metric $g$.  The Lagrangian $\rphys=\rphys[X,F;g]$ is the material energy density on the deformation bundle $\Dcal$. In addition, we assume \bei

\item[(A-3)] \underline{diffeomorphism invariance}, namely 
\be
\rphys \text{is invariant under spacetime diffeomorphisms}.
\ee
\eei {For a smooth spacetime diffeomorphism $\phi$ carrying admissible data to admissible data, we write $\phi^\star f=f\circ\phi$ and $\phi^\star g$ for the pull-back of the metric. The covariance identity is}
\be
\rphys[{\phi^\star} f, d({\phi^\star} f);{\phi^\star} g](x) = \rphys[f, df;g](\phi(x)), 
\qquad x \in \Mcal
\ee
for every admissible configuration map. The Noether identity gives the stress balance laws for stationary configurations (cf.~Proposition~\ref{prop:3} below).

By definition, the \emph{deformation bundle} $\Dcal$ has base $\Mcal\times\Bcal$, and the fiber over $(x,X)$ consists of all linear maps $F:T_x\Mcal\to T_X\Bcal$ satisfying (A-1)--(A-2) (both of which define an open subset).  The fiberwise closure $\overline{\Dcal}$ admits gradients of lower rank and rank-three gradients with null kernel. On $\Dcal$ we specify the differentiability of the energy in each result. At the boundary we use its finite constitutive limits. The stress may have more than one limiting value at the same degenerate gradient, so the weak solution includes the chosen stress field. The metric is prescribed and enters these constitutive data as a geometric parameter.

\begin{definition} 
\label{def-def0} 
A map $f:\Mcal\to\Bcal$ is an \underline{\emph{admissible configuration map}} if it is locally Lipschitz and $(f,df)$ takes values almost everywhere in $\overline{\Dcal}$, with the future orientation selected by the current~\eqref{vmu} below. The gradient has rank three and timelike kernel where $df\in\Dcal$. Its constitutive stress is specified by an admissible finite-flux extension and is required to be locally integrable when the balance laws are imposed.

An \underline{\emph{admissible constitutive law}} consists of an energy density $\rphys$ on $\Dcal$ satisfying covariance (A-3), with the differentiability required by the equations and a prescribed finite-flux completion at degenerate gradients.
\end{definition}
{The completion includes cavitation. We evaluate constitutive derivatives in matter and use the prescribed flux extension in the balance laws in the sense of distributions.}

\paragraph{{Material homogeneity.}}
Material homogeneity concerns dependence on $X$ throughout the constitutive domain.

\begin{definition}
\label{def:homogeneous-material}
A material is \underline{\emph{homogeneous in two material directions}}, or \underline{\emph{partially homogeneous}}, if, in the fixed material coordinates, its constitutive energy and finite-flux extension are independent of $X^1,X^2$. It is \underline{\emph{homogeneous}} if both are independent of all material coordinates. For the energy on its positive-density domain, the two cases are
\bse
\be
\label{eq:partial-homogeneity}
\rphys(X,F;g)=\rphys(X^3,F;g),\qquad \del_{X^1}\rphys=\del_{X^2}\rphys=0,
\ee
\be
\label{eq:homg}
\rphys(X,F;g)=\rphys(F;g),\qquad \del_{X^1}\rphys=\del_{X^2}\rphys=\del_{X^3}\rphys=0,
\ee
\ese
respectively.
\end{definition}

{In material volume coordinates, a perfect-fluid law of the form $\rphys=e(n)$ is homogeneous when $e$ is independent of the material coordinates.} The variational equation~\eqref{EL} below allows material dependence on $X$; homogeneity is imposed explicitly where the formulation uses $F$ without the configuration map.

%---------------------------------------------------------------------

\paragraph{The material current and loss of rank.}
To express the admissibility conditions, we use a \emph{material volume form} $\Omega_{ABC}$ prescribed on~$\Bcal$. For every configuration, including a rank-deficient one, the pull-back $3$-form on $\Mcal$ is
\be
\label{omeg}
\omega = f^* \Omega, 
\qquad 
\omega_{\mu\nu\rho} \coloneqq  F^A_\mu F^B_\nu F^C_\rho \, \Omega_{ABC} \circ f.
\ee
Taking its Hodge dual gives the \emph{material current}
\be
\label{vmu}
{\particlecurrent}^\mu \coloneqq  - \, \frac{1}{6} \omega_{\nu\rho\sigma} \, \eps^{\mu\nu\rho\sigma}. 
\ee
Here $\eps$ is the metric volume tensor, with indices raised by $g$. The current is future causal or zero on $\overline{\Xbb}$. On $\Xbb$ it is future timelike\footnote{In an oriented orthonormal Minkowski frame, with $\eps_{0123}=+1$, $g^{00}=-1$ and normalized material volume $\Omega_{123}=+1$, the sign in \eqref{vmu} gives ${\particlecurrent}^0>0$ for the rest configuration $f(t,x)=x$. Other rank-three gradients with timelike kernel can induce the opposite orientation and are excluded by this convention.}. The fields $\omega$ and $\particlecurrent$ contain equivalent information and may be used interchangeably.

To determine the direction of the current, we contract it with a component of the configuration gradient:
\be
{\particlecurrent}^\mu F^A_\mu
=
-\frac16\eps^{\mu\nu\rho\sigma}
F^A_\mu F^B_\nu F^C_\rho F^D_\sigma
\, \Omega_{BCD} \circ f,
\ee
which vanishes identically, since four vectors taken from the three-dimensional image of $F$ are necessarily linearly dependent and $\eps^{\mu\nu\rho\sigma}$ is totally antisymmetric. Hence we obtain the \emph{orthogonality condition}
\be
\label{ortho} 
{\particlecurrent}^\mu F^A_\mu = 0.
\ee
Thus $\particlecurrent$ lies in $\ker F$ at every rank. For $F\in\Xbb$, the kernel is one-dimensional and
\be
\ker F = \mathrm{span}\{\particlecurrent\}.
\ee
For a rank-three gradient, condition (A-2) is equivalent to a timelike current. In future-adapted coordinates its orientation gives
\be
g(\particlecurrent,\particlecurrent) = {\particlecurrent}^\mu {\particlecurrent}_\mu < 0, \qquad {\particlecurrent}^0 > 0.
\ee
On the entire causal state set, the \emph{particle number density} is
\be
n\coloneqq\sqrt{-g(\particlecurrent,\particlecurrent)}\ge0.
\ee
On $\{n>0\}$, the unit timelike \emph{material velocity} is
\be
\label{equa-uJ}
u^\mu \coloneqq  \frac{{\particlecurrent}^\mu}{n}. 
\ee
In matter, $g(u,u)=-1$ and $\particlecurrent=nu$. The relation between rank and current follows from the cubic minors in~\eqref{vmu}:
\be
\label{eq:rank-and-current}
 {\particlecurrent}(F)=0\quad\Longleftrightarrow\quad \operatorname{rank}F\le2 .
\ee
A nonzero null current consequently has $n=0$ and $\operatorname{rank}F=3$.  The admissible limits of the current and stress are specified by the constitutive graph of Section~\ref{section=3} below. A rank loss can have a nonzero limiting stress. We call an open region with $\particlecurrent=T=0$ an \emph{empty cavity}. A reference on a selected timelike branch may prescribe a Lipschitz continuation of $u$ to zero current.

In terms of $\omega$, the density identity holds on the whole causal set and the velocity identity on $\{n>0\}$.
\bse
\be
\label{eq:6699}
n^2  = - \frac{1}{36} g_{\mu\mu'} \omega_{\nu\rho\sigma} \omega_{\nu'\rho'\sigma'} \eps^{\mu\nu\rho\sigma} \eps^{\mu'\nu'\rho'\sigma'},
\ee
\be
\label{eq:uv}
u^\mu  = - \frac{1}{6 n} \omega_{\nu\rho\sigma} \, \eps^{\mu\nu\rho\sigma}.
\ee
\ese

The form $\omega(F)$ and the current $\particlecurrent(F)$ are cubic in $F$, and $n^2=-g(\particlecurrent,\particlecurrent)$ is of degree six. These polynomial expressions apply at every rank. The unit velocity $u=\particlecurrent/n$ is defined where $n>0$.

%-----------------------------------------------------------------------------------------------------------------------------------------

\subsection{Field equations for relativistic materials}

\paragraph{Least action principle.}
From now on, an admissible energy density $\rphys$ satisfying (A-3) is prescribed. We vary the configuration map in the action $\Acal$. For a material-indexed spacetime vector field $K_A$ in a fixed material chart, the spacetime divergence is
\be
\nabla_\mu K^\mu_A 
=
(- \det g)^{-1/2} 
\del_\mu \big( (-\det g)^{1/2} \, K_A^\mu \big).
\ee
{The first variation gives the following distributional equations.}

\begin{proposition}[Euler--Lagrange equations]
\label{prop-E-L}{Assume that $f\in W^{1,\infty}_{\rm loc}$ and that, on each compact set, its constitutive values belong almost everywhere to a compact subset of an open domain on which the energy is {$C^1$}.} {An admissible configuration map $f$ is stationary under all compactly supported variations of the extended action in that differentiability neighborhood if and only if, in the sense of distributions,}
\be
\label{EL}
\aligned
& \nabla_\mu \Big( \frac{\del \rphys}{\del F_\mu^A} \Big)  - \frac{\del \rphys}{\del f^A} = 0,
\endaligned
\ee
where $\del$ denotes differentiation with respect to the indicated argument.
\end{proposition}

{For a perfect fluid the energy depends on the proper particle density $n$ alone. The same variational statement applies to this choice under {its regularity assumptions}.} The Euler--Lagrange equations~\eqref{EL} are second order in the configuration map. For homogeneous materials, we use the stress balance laws and the compatibility equations for $F$ to obtain a first-order formulation.

%-------------------------------------------------------------

\paragraph{Euler equations for elastic materials.}
The \emph{energy--momentum tensor} $T=T(X,F;g)$ associated with the energy density is defined by
\be
\label{Tmunu}
T_{\mu \nu} 
\coloneqq 
2 \, \frac{\del \rphys}{\del g^{\mu\nu}} 
- \rphys \, g_{\mu\nu}.
\ee
This definition is consistent with the variational derivative of the action with respect to the spacetime metric and coincides with the standard Hilbert stress--energy tensor. For metric variations $g\mapsto g+\varepsilon h$, with $h$ a symmetric covariant tensor and the configuration map $f$ fixed, we obtain
\be
\left.\frac{d}{d\varepsilon}\right|_{\varepsilon=0}\Acal(f;g+\varepsilon h)
=-\frac12\int_\Mcal T^{\mu\nu}h_{\mu\nu}\,dV_g,
\ee
which justifies the definition~\eqref{Tmunu}.  Spacetime covariance relates its divergence to the first variation in $f$.

\begin{proposition}[Energy--momentum balance laws]
\label{prop:3}
Assume that the energy is $C^2$ on an open neighborhood of the configuration range, and assume a smooth configuration or the controlled regularity in Proposition~\ref{prop:variational-EL-local04} below. Consider any admissible energy density $\rphys$, which, in particular, is diffeomorphism-invariant. {If an admissible configuration $f: \Mcal \to \Bcal$ is stationary under all compactly supported variations of the extended action $\Acal$ in that differentiability neighborhood, then} the energy--momentum tensor \eqref{Tmunu} (evaluated at $f$ and its deformation tensor $F$) satisfies the balance laws
\be
\label{equa-main}
\nabla_\nu \bigl( T^{\mu\nu}(f,F;g) \bigr) = 0
\ee
or, equivalently, for every {smooth} vector field $\xi$ on $\Mcal$
\be
\nabla_\nu \big( \xi_\mu T^{\mu\nu}\big)
=
T^{\mu\nu} \nabla_{(\mu} \xi_{\nu)}
=
\pi^{(\xi)}_{\mu\nu} \, T^{\mu\nu},
\ee
where the \emph{deformation tensor} of a vector field $\xi$ is defined by
\be
\pi^{(\xi)}_{\mu\nu} 
\coloneqq  
\frac{1}{2}
\big( \nabla_\mu \xi_\nu + \nabla_\nu \xi_\mu \big).
\ee
\end{proposition}

{For a perfect fluid with rest energy $e(n)$ and timelike current $\particlecurrent=nu$, the tensor in these balance laws takes the familiar form $T^{\mu\nu}=(e+p)u^\mu u^\nu+pg^{\mu\nu}$, where $p=ne'(n)-e$.}

\begin{proof}
\bse
Independently of the Euler--Lagrange equations \eqref{EL} and for configuration maps $f$ in {this regularity class}, the identity
\be
\nabla^\nu T_{\mu\nu}
=
\Ecal_A F^A_\mu
\ee
holds, where $\Ecal_A$ denotes the Euler--Lagrange operator
\be
\Ecal_A 
\coloneqq 
\nabla_\mu \Big( \frac{\del \rphys}{\del F_\mu^A} \Big)
- \frac{\del \rphys}{\del f^A}.
\ee
We take $\xi$ to be a smooth compactly supported vector field on $\Mcal$ and write $\phi_\varepsilon$ for its flow, with $\phi_0=\mathrm{id}$ and $\partial_\varepsilon\phi_\varepsilon=\xi\circ\phi_\varepsilon$. The star denotes pull-back: $\phi_\varepsilon^\star f=f\circ\phi_\varepsilon$, and $\phi_\varepsilon^\star g$ is the pull-back of the metric. On a noncompact spacetime, the action in the next formula is restricted to a relatively compact domain containing the support of $\xi$. Diffeomorphism invariance gives
\be
0 = \left.\frac{d}{d\varepsilon}\right|_{\varepsilon=0}
\Acal(\phi_\varepsilon^\star f ; \phi_\varepsilon^\star g).
\ee
For smooth $f$, the variations of $f^A$ and $g_{\mu\nu}$ are $\xi^\mu F^A_\mu$ and $2\nabla_{(\mu}\xi_{\nu)}$. The configuration part of the first variation is $-\int \Ecal_A F^A_\mu\xi^\mu\,dV_g$, and the metric part is $-\int T^{\mu\nu}\nabla_\mu\xi_\nu\,dV_g$. Integration by parts gives $\nabla^\nu T_{\mu\nu}=\Ecal_A F^A_\mu$. {At this Sobolev regularity}, Proposition~\ref{prop:variational-EL-local04} below proves the same identity by the chain and product rules. Stationarity gives $\Ecal_A=0$, so~\eqref{equa-main} follows.
\ese
\end{proof}

%--------------------------------------------------------------------------------------------------------------------------------

%-----------------------------------------------------------------------------------------------------------------------------------------

\subsection{Evolution in the configuration gradient}

Interestingly, the dynamics of a homogeneous relativistic material can be described directly in terms of the configuration gradient $F$, from which the stress is computed using the metric. The four balance laws~\eqref{equa-main} are the evolution equations for this formulation. The twelve components of $F$ also satisfy the exterior compatibility equations that express its local recovery as $df$.

\begin{definition}
\label{def-firstorder}
The \underline{\emph{first-order balance laws}} of the dynamics of homogeneous materials (cf.~\eqref{eq:homg}) are the energy-momentum balance laws derived in \eqref{equa-main}, that is,
\be
\label{equa-main-deux}
\nabla_\nu T^{\mu\nu}=0.
\ee
Here $F=(F_\mu{}^A)$ ranges in the completed kinematic state set, including cavitation. Where $n>0$, the constitutive stress is given by \eqref{Tmunu}, equivalently \eqref{eq:useful} below. The constitutive continuation specifies the limiting stress there.
\end{definition}

In Section~\ref{sec:weak-local04} below we impose compatibility $dF=0$ and the initial terms as part of the weak formulation. For a perfect fluid, the same balance laws use the perfect-fluid stress, with the particle current  giving the additional continuity law or mathematical entropy inequality.
A change of spacetime basis acts on both $F$ and $g$. Differentiating covariance with respect to this change expresses the stress in terms of the derivative of the energy with respect to $F$.

\begin{proposition}[Structure of the energy-momentum tensor] 
\label{prop:stru}
{Let $\rphys=\rphys(X,F;g)$ be $C^2$ on an open positive-density constitutive domain preserved by simultaneous changes of spacetime basis.} The energy density is covariant under these changes if and only if
\be
\label{equa-Lag}
2\,\frac{\del \rphys}{\del g^{\mu\nu}}
=
g_{\alpha(\mu} F^A_{\nu)}
\frac{\del \rphys}{\del F^A_\alpha}
\ee
(where the parentheses denote symmetrization). Consequently, the energy--momentum tensor admits the decomposition
\be
\label{eq:useful}
T_\mu{}^\nu
=
- \rphys\, \delta_\mu{}^\nu
+
F^A_\mu
\frac{\del \rphys}{\del F^A_\nu}. 
\ee
{For a timelike current, $u^\mu F_\mu^A=0$ follows from~\eqref{ortho} and~\eqref{equa-uJ}. Contracting~\eqref{eq:useful} with $u^\mu$ therefore gives $T_\mu{}^\nu u^\mu=-\rphys u^\nu$.}  At $\particlecurrent=0$, the constitutive extension specifies the stress.
\end{proposition}

For a perfect fluid this eigenvalue relation is the usual identity $T_\mu{}^\nu u^\mu=-e(n)u^\nu$ on the timelike-current set.
Importantly, the assumed covariance under spacetime coordinate changes imposes \emph{no restriction} on the material symmetry group. It implies neither homogeneity (independence of $X$) nor material isotropy.

%============================================================================

\section{\texorpdfstring{An augmented formulation of relativistic elastodynamics}{An augmented formulation of relativistic elastodynamics}}
\label{section=3}

%-----------------------------------------------------------------------------------------------------------------------------------------

\subsection{A new set of unknowns}

\paragraph{The spacetime minor map.}
{
We first describe the new fields as functions of $F$, then use them as independent variables in the balance laws. The continuity equation and the transport of the forms will give conditions under which the initial minor relations persist at zero current.\par}

From this section onward, we assume a homogeneous material in the sense of Definition~\ref{def:homogeneous-material}, with the fixed material volume coordinates described above. The stress and the other fluxes in the augmented equations may still depend on the spacetime point through the prescribed metric. For an inhomogeneous energy $\rphys(X,F;g)$, the material equation~\eqref{EL} contains the source $\rphys_{X^A}(f,F;g)$. The chain rule for the energy adds $\rphys_{X^A}F^A_\mu$; differentiating other constitutive coefficients gives their corresponding material derivatives. These terms have lower differential order than the derivatives of $F$ in the principal part. Their treatment uses the configuration map, or transported material coefficients, together with regularity and relative-entropy estimates for the additional sources.

We augment the equations of Definition~\ref{def-firstorder} by the following quadratic and cubic minors, defined for every $F\in T^*\Mcal\otimes T\Bcal$ and at every rank:
\be
\label{eq-312}
G_{\mu\nu A} = F^B_{\ [\mu} F^C_{\ \nu]} \,\Omega_{ABC},
\qquad
{\particlecurrent}^\mu = -\frac{1}{6}\, F^A_{\ \nu} F^B_{\ \sigma} F^C_{\ \tau}
\Omega_{ABC}\,\eps^{\mu\nu\sigma\tau}, 
\ee
where the brackets denote antisymmetrization.

\begin{definition}
\label{def:spacetime-minor-map}
The \underline{\emph{spacetime minor map}} maps a configuration gradient to its first, second and third minors:
\be
\mathbf{\Xi}:\mathbb R^{12}\longrightarrow\mathbb R^{34},\qquad F\longmapsto(F,G(F),\particlecurrent(F)).
\ee
The \underline{\emph{set of exact spacetime minors}} is the polynomial graph
\be
\label{eq:spacetime-minor-manifold}
\mathcal C_{\rm min}\coloneqq\{(F,G,\particlecurrent):G=G(F),\ \particlecurrent=\particlecurrent(F)\}=\mathbf\Xi(\mathbb R^{12}).
\ee
Its kinematically admissible part and positive-density part are, respectively,
\be
\overline{\Cscr}\coloneqq\mathbf{\Xi}(\overline{\Xbb}),\qquad
\Cscr\coloneqq\mathbf{\Xi}(\Xbb).
\ee
\end{definition}

{The fields in $\mathcal C_{\rm min}$ satisfy $G=G(F)$ and $\particlecurrent=\particlecurrent(F)$ pointwise. This map is kinematic and is independent of the constitutive law. For a perfect fluid, the unknown is the particle current $\particlecurrent$ and the constitutive law gives its stress; the unused fields $F,G$ may be omitted.} Projection onto $F$ gives the inverse of $\mathbf\Xi$ on its image. Thus $\mathcal C_{\rm min}$ is smooth at every rank, and $\overline\Cscr$ is the finite-component closure of $\Cscr$. This graph gives the minors $(G,\particlecurrent)$ from $F$, while the constitutive law specifies $T$.

\begin{lemma}[Differential of the spacetime minor map]
\label{prop:constraint-rank} 
For every deformation tensor $F\in\mathbb R^{12}$, the differential
\be
D\mathbf{\Xi}(F): \R^{12}\to \R^{34}
\ee
is injective. Consequently, $\mathbf{\Xi}$ is an immersion on any open subset of its domain (in particular, on the rank-three and timelike-kernel locus), and for every \(\Xi=\mathbf{\Xi}(F)\in\Cscr\), the tangent space $T_\Xi\Cscr$ has dimension
\be
\dim T_{\Xi}\Cscr = 12
\quad
(\text{hence } \operatorname{codim}_{\R^{34}}\Cscr = 22).
\ee
\end{lemma}

\begin{proof}
\bse
Since the first component of the map \(\mathbf{\Xi}\) is \(F\) itself, we have
\be
D\mathbf{\Xi}(F)[\delta F]
=
\bigl(\delta F,\;DG(F)[\delta F],\;D\particlecurrent(F)[\delta F]\bigr).
\ee
In particular, if \(D\mathbf{\Xi}(F)[\delta F]=0\), then its first component gives \(\delta F=0\). Hence \(D\mathbf{\Xi}(F)\) is injective.  For completeness, the differentials of \(G\) and \(\particlecurrent\) are
\be
\aligned
DG(F)[\delta F]_{\mu\nu A}
& = \delta F^B_{\ [\mu} F^C_{\ \nu]} \,\Omega_{ABC}
+ F^B_{\ [\mu} \delta F^C_{\ \nu]} \, \Omega_{ABC},
\\
D\particlecurrent(F)[\delta F]^\mu
&= - \frac{1}{2} \, \delta F^A_{\ \nu} F^B_{\ \sigma} F^C_{\ \tau}
\Omega_{ABC}\eps^{\mu\nu\sigma\tau}. \qedhere
\endaligned
\ee
\ese
\end{proof}

%----------------------------------------------------------

\paragraph{Independent fields and constitutive energy.}
In the augmented formulation, the unknown fields are
\be
\bXi=(\bF,\bG,\bJ),
\ee
{Here $\bF$ consists of material one-forms, $\bG$ of antisymmetric material two-forms, and $\bJ$ is a vector. We impose the algebraic contractions with $\bJ$ and the exterior equations specified below.}

\begin{definition}
\label{def:constitutive-minor-extension}
{A \underline{\emph{constitutive energy in independent spacetime minors}} is a covariant $C^2$ {function $\raug(F,G,\particlecurrent;g)$ of the independent fields $(F,G,\particlecurrent)$ on an open neighborhood} $\mathcal N_g$ of the positive-density part $\Cscr$ of the set of exact spacetime minors. Its restriction to exact minors is the configuration energy}
\be
\label{eq:constitutive-minor-restriction}
\rphys(F;g)=\raug\bigl(F,G(F),\particlecurrent(F);g\bigr),\qquad F\in\Xbb,
\ee
\end{definition}

The current $\particlecurrent(F)$ uses the volume tensor of $g$, a dependence included when differentiating $\rphys$ in $F$ or $g$. For a perfect fluid one may choose $\raug=e(n)$ with $n=\sqrt{-g(\particlecurrent,\particlecurrent)}$; it is independent of the form variables $F$ and $G$ in the positive-density domain.} Restriction to $\Cscr$ gives the original energy; values away from $\mathcal C_{\rm min}$ are part of the chosen augmented law. This extends the classical minor formulation of LeFloch--Qin, motivated by Ball's polyconvexity, to a $3\times4$ configuration gradient: its quadratic and cubic spacetime minors replace spatial cofactors and determinant. The finite-flux graph specifies the stress in regions of cavitation without requiring $\raug$ to be differentiable at zero current.

The spacetime minors of gradients are \emph{weakly continuous} under the following integrability assumption. If $F^j=df^j\rightharpoonup F=df$ in $L^p_{\rm loc}$ with $3<p<+\infty$, then $G(F^j)\rightharpoonup G(F)$ in $L^{p/2}_{\rm loc}$ and $\particlecurrent(F^j)\rightharpoonup {\particlecurrent}(F)$ in $L^{p/3}_{\rm loc}$. Indeed, locally subtracting the means gives strong $L^p$ convergence of the potentials. Integration by parts in the exterior products gives distributional convergence of the quadratic minors and then of the cubic minors. Their bounds in $L^{p/2}$ and $L^{p/3}$ give {the asserted weak convergence} {(cf.~\cite{Ball1}). The stress limit requires additional constitutive estimates.
 
%-----------------------------------------------------------------------------------------------------------------------------------------

\subsection{The proposed augmented formulation}

\paragraph{Algebraic constraints.}
{We now replace the configuration gradient and its minors by independent fields and admit their specified constitutive limits. Conditions (A-1)--(A-2) from Section~\ref{section=2} are not imposed on this full augmented domain; its algebraic constraints and constitutive admissibility {are given below}.} All equations under consideration are \emph{standard balance laws}, that is, equations of the form
\be
\label{equa-bala}
\nabla_\mu X^{\mu\ldots}=Y^{\ldots},
\ee
for suitable tensor fields $X$ and $Y$. Eventually, such equations will be considered in the sense of distributions. For example, the compatibility condition
\be
\nabla_{[\mu} \bF^A_{\nu]}=0,
\ee
expressing that $\bF$ is (locally) a gradient, can be written in divergence form as
\be
\nabla_{[\mu} \bF^A_{\nu]}
=\frac12\nabla_\lambda\big(\delta^\lambda_\mu \bF^A_\nu-\delta^\lambda_\nu \bF^A_\mu\big)=0,
\ee
and hence it is itself a balance law of the form \eqref{equa-bala}. Throughout, we impose the antisymmetry condition
\be
\label{eq:449}
\bG_{\mu\nu A}=-\bG_{\nu\mu A}.
\ee
We observe that $\bF$ has $3\times 4=12$ independent components, each $\bG_{\mu\nu A}$ is antisymmetric in $\mu,\nu$ and thus carries $3\times \binom42=18$ components, and $\bJ$ is a spacetime $4$--vector with $4$ components. These add up to $12+18+4=34$ components, justifying the ambient space $\RR^{34}$. Motivated by the following lemma, we will impose the corresponding algebraic orthogonality and antisymmetry constraints on the augmented variables $\bF, \bG, \bJ$.

\begin{lemma}[Orthogonality and antisymmetry properties]
\label{lem:alg-constraints}
Let $F=(F^A{}_\mu)$ be a deformation tensor of arbitrary rank and consider
\be
G_{\mu\nu A}\coloneqq F_{[\mu}^B F_{\nu]}^C\,\Omega_{ABC},
\qquad
{\particlecurrent}^\mu\coloneqq -\frac16\,F_\alpha^A F_\beta^B F_\gamma^C\,\Omega_{ABC}\,\eps^{\mu\alpha\beta\gamma}.
\ee
The minors satisfy the two orthogonality identities and the antisymmetry relation
\be
\label{eq:450bis}
{\particlecurrent}^\mu F_\mu^A=0,
\qquad
{\particlecurrent}^\mu G_{\mu\nu A}=0,
\qquad
G_{\mu\nu A}=-G_{\nu\mu A}.
\ee
\end{lemma}

\begin{proof}
The antisymmetry of $G$ is immediate. The identity ${\particlecurrent}^\mu F_\mu^A=0$ follows from total antisymmetry of $\eps^{\mu\alpha\beta\gamma}$ and the fact that four vectors drawn from the three-dimensional image of $F$ are linearly dependent (equivalently, $\particlecurrent$ spans $\ker F$ when $\mathrm{rank}\,F=3$). The second identity follows from ${\particlecurrent}^\mu G_{\mu\nu A}=({\particlecurrent}^\mu F_\mu^B)F_\nu^C\Omega_{ABC}=0$.
\end{proof}

%---------------------------------------------------------------------

\paragraph{{The orthogonality condition.}}
{We impose on the independent forms the contraction identities satisfied by the exact minors.}
\begin{definition}
\label{def:transversality}
For each material index $A=1,2,3$, the one-form $\bF^A$ and the two-form $\bG_A$ satisfy the \underline{\emph{orthogonality condition}} with respect to a vector $\bJ$ if
\be
\label{eq:450bis-2}
\bJ^\mu\bF_\mu^A=0,\qquad \bJ^\mu\bG_{\mu\nu A}=0.
\ee
\end{definition}

The first equality means that $\bJ$ is orthogonal to the metric dual $(\bF^A)^\sharp$ for the Lorentzian metric $g$; the second means that $\iota_{\bJ}\bG_A=0$. Where $\bJ=nu$ with $n>0$, these conditions are equivalent to $\iota_u\bF^A=\iota_u\bG_A=0$. At $\bJ=0$ they hold automatically. If a velocity is continued into a zero-current region, orthogonality to that velocity is a separate condition wherever it is required.
Figure~\ref{fig:kinematic-constitutive-boundaries} illustrates these contractions and the causal character of the current.

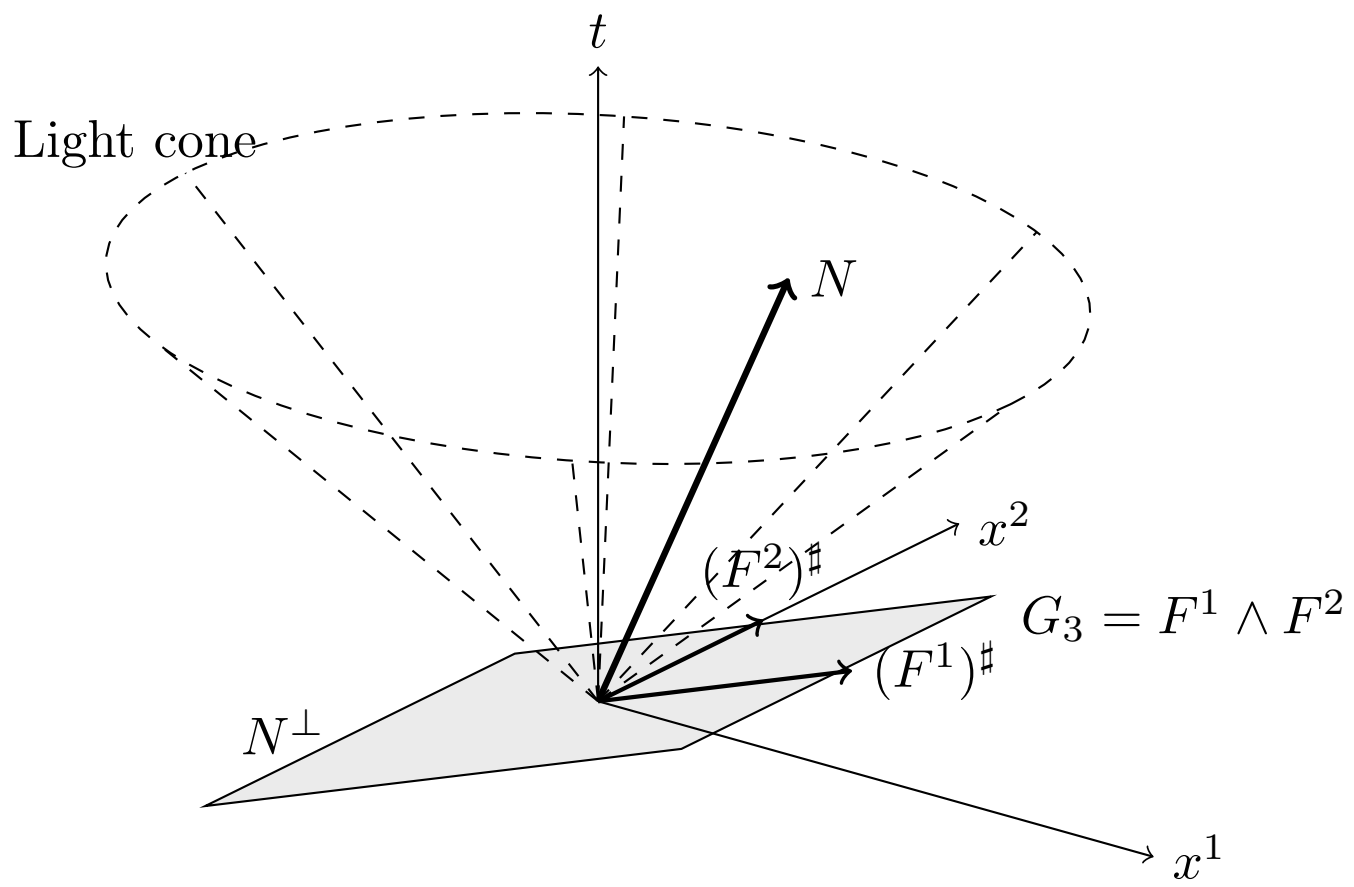
\begin{figure}[!ht]
\centering
\begin{tikzpicture}[x={(1cm,-.28cm)},y={(.65cm,.32cm)},z={(0cm,1cm)},scale=1.05]
 \draw[->] (0,0,0)--(3.5,0,0) node[right] {$x^1$};
 \draw[->] (0,0,0)--(0,3.5,0) node[right] {$x^2$};
 \draw[->] (0,0,0)--(0,0,4) node[above] {$t$};
 \draw[dashed] plot[domain=0:360,samples=90] ({2.6*cos(\x)},{2.6*sin(\x)},2.6);
 \foreach \a in {0,60,120,180,240,300} {
  \draw[dashed] (0,0,0)--({2.6*cos(\a)},{2.6*sin(\a)},2.6);
 }
 \fill[black,opacity=.08] (-1.5,-1.5,-.6)--(1.5,-1.5,.6)--(1.5,1.5,.6)--(-1.5,1.5,-.6)--cycle;
 \draw (-1.5,-1.5,-.6)--(1.5,-1.5,.6)--(1.5,1.5,.6)--(-1.5,1.5,-.6)--cycle;
 \draw[->,very thick] (0,0,0)--(1.2,0,3) node[right] {$\particlecurrent$};
 \draw[->,thick] (0,0,0)--(1.6,0,.64) node[right] {$(F^1)^\sharp$};
 \draw[->,thick] (0,0,0)--(0,1.6,0) node[above] {$(F^2)^\sharp$};
 \node[anchor=west] at (-2.7,-1.7,3.3) {Light cone};
 \node[anchor=east] at (-1.6,0,-.64) {${\particlecurrent}^\perp$};
 \node[anchor=west] at (1.7,1.3,.6) {$G_3=F^1\wedge F^2$};
\end{tikzpicture}
\caption{{A three-dimensional spacetime slice illustrating the orthogonality condition~\eqref{eq:450bis-2} in matter. The metric is $-dt^2+(dx^1)^2+(dx^2)^2$, and the omitted coordinate is $x^3$. For $\particlecurrent=r(\partial_t+v\partial_{x^1})$, $r>0$ and $|v|<1$, take $F^1=dx^1-vdt$ and $F^2=dx^2$. The displayed spacelike plane is ${\particlecurrent}^\perp$, spanned by the metric duals $(F^1)^\sharp$ and $(F^2)^\sharp$; the two-form $G_3$ is represented by their oriented plane. The drawing uses $v=0.4$. These are geometric representatives of forms, rather than three vector unknowns. As $|v|\to1$, $\particlecurrent$ approaches a null direction and this orthogonal plane becomes null. At $\particlecurrent=0$, the current alone no longer determines an orthogonal plane.}}
\label{fig:kinematic-constitutive-boundaries}
\end{figure}

{In the fluid reduction the current can be evolved without the auxiliary forms; if the latter are  included, the contractions are imposed under the same hypotheses as for an elastic material.}
The independent augmented fields satisfy the antisymmetry~\eqref{eq:449} and the conditions~\eqref{eq:450bis-2}.
These equations define the \emph{augmented algebraic constraint set}
\be
\bCscr\coloneqq
 \{(\bF,\bG,\bJ)\in\mathbb R^{34}: \text{\eqref{eq:449} and \eqref{eq:450bis-2} hold}\}.
\ee
{On the set of exact spacetime minors these relations follow from the exact minor identities. On $\bCscr$ they constrain the independent fields for $\particlecurrent\ne0$. At $\particlecurrent=0$, the constitutive closure below imposes further restrictions.}

\paragraph{Differential constraints.}
The exact quadratic minors are wedge products of closed one-forms. This gives the differential equations we impose on the independent two-forms.

\begin{lemma}[Closedness of the quadratic minor forms]
\label{lem:diff-constraints}
In material volume coordinates, assume that, for each $A$, the $1$--form $F^A\coloneqq F^A{}_\mu\,dx^\mu$ {belongs to $L^{p_A}_{\rm loc}$ and is closed in the sense of distributions, where $1\le p_A\le\infty$ and $p_B^{-1}+p_C^{-1}\le1$ for every distinct pair $B,C$, with $1/\infty=0$}, i.e.
\be
\nabla_{[\mu}F^A_{\nu]}=0.
\ee
By defining the quadratic minor forms $G(F)$ by
\be
G_{\mu\nu A}=F_{[\mu}^B F_{\nu]}^C\,\Omega_{ABC},
\ee
one also has {$G_A\in L^1_{\rm loc}$ and, in the sense of distributions,}
\be
\label{eq:451bis}
\nabla_{[\lambda}G_{\mu\nu]A}=0.
\ee
\end{lemma}

\begin{proof}
In fixed material volume coordinates, $\Omega_{ABC}$ is constant and $G_A=\tfrac12\Omega_{ABC}F^B\wedge F^C$. {For smooth closed forms,}
\be
dG_A=\tfrac12\Omega_{ABC}(dF^B\wedge F^C-F^B\wedge dF^C)=0.
\ee
{For the given closed forms, convolution on a smaller coordinate ball preserves closedness. The approximations converge strongly in each finite $L^{p_A}_{\rm loc}$ and, for $p_A=\infty$, almost everywhere with a common bound. H\"older's inequality and dominated convergence give convergence of each quadratic product in $L^1_{\rm loc}$, hence distributional closedness. The choice of material volume coordinates ensures $d(i_{\partial_A}\Omega)=0$. For a general volume density $\omega$, this derivative is $(\partial_A\omega)\,dX^1\wedge dX^2\wedge dX^3$.}
\end{proof}

%---------------------------------------------------------------

\paragraph{The augmented system.}
The constitutive chain rule gives the stress on the positive-density part of the set of exact spacetime minors $\Cscr$. We use the corresponding expression for the independent fields $(\Fb,\Gb,\Jb)$. The contraction identities used in this expression are established in Lemma~\ref{lem:algebraic-identities} below. {Under the transport and current-conservation hypotheses of Corollary~\ref{cor:physical-minors-local06} below, the polynomial relations propagate from the initial hypersurface, including at zero current. For an augmented weak solution, the exact minor relations are imposed only in the results that use them.}

Moreover, when $(\Fb,\Gb,\Jb)=(F,G(F),\particlecurrent(F))\in\Cscr$, the augmented tensor $\Tb$ is designed to be defined over $\bCscr$ and to coincide with the physical stress--energy tensor, namely
\be
\Tb^\mu_{\ \nu}(F,G(F),\particlecurrent(F)) = T^\mu_{\ \nu}(F),
\ee
so that $\Tb$ is indeed an \emph{extension of $T$ from $\Cscr$ to $\bCscr$}. Writing $\bCscr_+\coloneqq\{(F,G,\particlecurrent)\in\bCscr:\particlecurrent\text{is future timelike}\}$, on this positive-density part the augmented energy has the derivatives required by this formula. At the constitutive boundary, the chosen constitutive continuation specifies the stress through the limits below.

\begin{definition}
\label{def:constitutive-flux-graph}
For an augmented energy $\raug$, define the \underline{\emph{constitutive flux graph}} by
\be
\label{eq:finite-flux-state-space}
{\mathcal Z_{\raug}}\coloneqq
\{(F,G,\particlecurrent,T_{\raug}(F,G,\particlecurrent;g)):(F,G,\particlecurrent)\in\operatorname{Dom}(\raug)\cap\bCscr_+\}.
\ee
{We write $\Cl(\mathcal Z_{\raug})$ for its closure in the space of $(F,G,\particlecurrent,T)$, and call it the \underline{\emph{completed constitutive flux graph}}.}
\end{definition}

{Each graph consists of tensor values at a point $x$. We suppress $x$ and the prescribed metric $g$ in the notation. Their definition uses $g_x$ to determine the causal character of $\particlecurrent$, the proper density $n=\sqrt{-g_x(\particlecurrent,\particlecurrent)}$, and the constitutive stress. The fields and metric transform together.}

{
We write $Z=(F,G,\particlecurrent,T)$ for the unknowns and their constitutive stress. In matter the law determines $T$, and at the constitutive boundary its chosen continuation specifies the stress when several limits are possible. We use the same continuation for the initial stress. For a perfect fluid the constitutive flux graph and its closure use $(\particlecurrent,T)$ alone. Thus $\mathcal C_{\rm min}$ describes the exact polynomial minor relations, $\bCscr$ gives the constraints on the independent triples, and $\Cl(\mathcal Z_{\raug})$ specifies the admissible stresses.
}

\paragraph{{Fundamental and additional algebraic relations}}
Definition~\ref{propo=30DF} below uses antisymmetry, orthogonality and membership in the completed constitutive flux graph. These requirements have an algebraic consequence at every point of the completed constitutive flux graph: all the forms have zero contraction with the same nonzero future causal direction. For $\particlecurrent\ne0$ use $\particlecurrent$ itself. For $\particlecurrent=0$ normalize the timelike currents in an approximating sequence by a fixed time covector and pass to a subsequence. In an orthonormal observer frame the common direction can be written as $\ell=(1,v)$ with $|v|\le1$. Its contractions give $F_0^A=-v^iF_i^A$ and $G_{0j A}=-v^iG_{ij A}$, including at zero current. Hence the two-forms lie in the exterior powers of the same three-dimensional annihilator and satisfy
\be
\label{eq:minor-wedge-constraint}
G_A\wedge G_B=0,\qquad A,B=1,2,3.
\ee
{The common causal direction is an additional restriction on the limiting forms at zero current. The exact minors also satisfy the mixed identity}
\be
\label{eq:exact-mixed-minor}
 \Omega^{ABC}G_{\mu\nu A}G_{\rho\sigma B}
 ={\particlecurrent}^\lambda\bigl(\eps_{\lambda\rho\sigma\mu}F_\nu^C
                 -\eps_{\lambda\rho\sigma\nu}F_\mu^C\bigr).
\ee
{Here $\Omega^{123}=1$ in normalized material volume coordinates. To derive it, substitute $G_{\mu\nu A}=\Omega_{ADE}F_\mu^DF_\nu^E$ and use ${\particlecurrent}^\lambda\eps_{\lambda\rho\sigma\mu}=\Omega_{ABC}F_\rho^AF_\sigma^BF_\mu^C$. The identity holds whenever the exact minor relations hold. The mixed relation can fail for independent fields: with $\particlecurrent=\partial_0$, $F^1=dx^1$ and $F^2=F^3=G=0$, the orthogonality and wedge relations hold, but the mixed identity fails for $(C,\mu,\nu,\rho,\sigma)=(1,3,1,1,2)$.}

{The common causal direction is also sufficient for $(F,G,0)$ to belong to the algebraic closure $\overline{\bCscr_+}$. Normalize $\ell$ by a fixed time covector $\tau$, choose future timelike $u_j\to\ell$ with $\tau(u_j)=1$, and set $P_j=I-u_j\otimes\tau$. Then $(F\circ P_j,P_j^*G,\varepsilon_j u_j)\in\bCscr_+$ converges to $(F,G,0)$ as $\varepsilon_j\downarrow0$. For a nonzero future null current $\particlecurrent$, the same construction with $\ell=\particlecurrent/\tau(\particlecurrent)$ and ${\particlecurrent}_j=\tau(\particlecurrent)u_j$ shows that the orthogonality condition~\eqref{eq:450bis-2} suffices for membership of $\overline{\bCscr_+}$. A finite constitutive stress limit remains a separate requirement.

For exact gradients of rank at most two, the same causal-kernel condition characterizes $\overline{\Xbb}$. Necessity follows from normalized timelike kernels. For sufficiency, $F\circ P_j$ has rank at most two and converges to $F$. In bases of $\ker\tau$ and the material tangent space, write the restriction of $F$ as $L\operatorname{diag}(I_k,0)R$, where $k\le2$ and $L,R$ are invertible. Replace the zero diagonal entries by small nonzero entries, choosing the sign of the last one so that the resulting map $A_j$ gives a future current. Then $A_j\circ P_j$ has timelike kernel $\mathbb Ru_j$ and converges to $F$. This choice of sign is possible at every rank $k\le2$. At rank three, $F\in\overline{\Xbb}$ exactly when $\particlecurrent(F)$ is future causal. In the null case, $F\circ P_j$ has rank three and timelike kernel $\mathbb Ru_j$, because $F|_{\ker\tau}$ is invertible; its currents converge to $\particlecurrent(F)$ and have the future orientation for large $j$. These criteria concern gradients;  membership in ${\Cl(\mathcal Z_{\raug})}$ also requires a finite stress limit.
}

\paragraph{{The constitutive boundary.}}
{Taking the closure in both fields and stress includes the limiting stresses at zero current.}
\begin{definition}
\label{def:kinematic-constitutive-boundaries}
The \underline{\emph{kinematic boundary}} and the \underline{\emph{constitutive boundary}} are, respectively,
\be
\label{eq:kinematic-constitutive-boundaries}
 {\partial\Cscr}\coloneqq\overline{\Cscr}\setminus\Cscr,\qquad
 {\partial\mathcal Z_{\raug}}\coloneqq
 {\Cl(\mathcal Z_{\raug})}\setminus{\mathcal Z_{\raug}}.
\ee
\end{definition}
{Here $\partial\Cscr$ is the boundary relative to the set of exact minors, whereas $\partial\mathcal Z_{\raug}$ denotes the points added to the constitutive flux graph by closure. The latter need not be the topological boundary of $\mathcal Z_{\raug}$ in the ambient space.}

\begin{definition}
\label{def:residual-stress-zero-current}
{For a constitutive limit $(F,G,0,T)\in{\Cl(\mathcal Z_{\raug})}$, a nonzero $T$ is called a \underline{\emph{residual stress}}, that is, a stress persisting at zero particle current.}
\end{definition}
{Its value is specified by the constitutive continuation. An empty cavity has $\particlecurrent=T=0$, whereas $F$ and $G$  keep their tensor definitions and may be nonzero in the independent augmented formulation.} {On the other hand, the fields themselves and the quantities reconstructed from them have different behavior at the constitutive boundary. The following proposition makes this distinction precise.}
\begin{proposition}[ Velocity and stress at loss of rank]
\label{prop:boundary-fields-and-velocity}
{Let $F,G,\particlecurrent$ be constitutive fields with associated stress $T$ such that $(F,G,\particlecurrent,T)\in{\Cl(\mathcal Z_{\raug})}$, with $\particlecurrent$ future causal or zero, and let $n=\sqrt{-g(\particlecurrent,\particlecurrent)}$. One has the following alternatives.
\begin{itemize}
\item[(i)] The tensors $F,G$  keep their definitions at every rank. If $n>0$, the current determines the future unit velocity $u=\particlecurrent/n$, independently of their ranks.
\item[(ii)] If $\particlecurrent=0$, the relation $\particlecurrent=nu$ holds for every future unit timelike vector $u$ and determines none. A continued velocity must satisfy the additional conditions imposed by the chosen evolution or constitutive continuation.
\item[(iii)] If $\particlecurrent\ne0$ and $n=0$, the current is lightlike and cannot be written as $nu$ with a finite unit timelike vector $u$. There is no timelike rest frame for this current.
\end{itemize}}
\end{proposition}
{A loss of rank alone therefore requires no new unknown in place of $F$ or $G$. A separate choice can be necessary for the velocity or for the limiting stress. If several stress limits correspond to the same $(F,G,\particlecurrent)$, the constitutive continuation must select the stress, or specify parameters that do so. This is additional constitutive information, not a consequence of the ranks of $F,G$.}
\begin{proof}
{For $n>0$, normalization gives $g(u,u)=-1$ and fixes $u$ uniquely. At zero current and $n=0$, multiplication by $n$ gives zero for every such vector. A nonzero lightlike current has zero Lorentzian norm and remains nonzero, so it admits neither this normalization nor a frame with vanishing spatial current. None of these statements changes the tensor definitions of $F$ and $G$.}
\end{proof}

For a fluid described only by its current and stress, the constitutive closure is taken in those variables. Its boundary behavior depends on the pressure law and can include zero current or finite nonzero null-current limits.
The projection $\pi(F,G,\particlecurrent,T)=(F,G,\particlecurrent)$ satisfies $\pi({\Cl(\mathcal Z_{\raug})})\subseteq\overline{\pi({\mathcal Z_{\raug}})}$ by continuity; a limit of triples alone need not have a finite limiting stress.} For exact minors,~\eqref{eq:rank-and-current} and the timelike-kernel condition distinguish rank loss from a rank-three gradient with null current. The quadratic minor map $G(F)$ is the second exterior power of $F$ up to the material volume identification, so its rank on bivectors is $\binom{\operatorname{rank}F}{2}$. For independent augmented fields, the ranks of $F$ and $G$ can drop separately, even at positive density. {Such rank loss alone does not imply that the fields and their constitutive stress belong to the constitutive boundary ${\partial\mathcal Z_{\raug}}$.} The constitutive boundary specifies which limits of these fields and their stress are admitted by the law. {Table~\ref{tab:kinematic-constitutive-sets} summarizes these sets.}

\begin{table}[!ht]
\centering
\begingroup
\begin{tabular}{@{}p{.32\linewidth}p{.62\linewidth}@{}}
\toprule
Set & Meaning \\
\midrule
$\Xbb$ & \leavevmode{Gradients $F$ with $\operatorname{rank}F=3$ and future timelike $\particlecurrent(F)$} \\
$\mathcal C_{\rm min}$ & \leavevmode{Triples satisfying $G=G(F)$ and $\particlecurrent=\particlecurrent(F)$, at every rank} \\
$\Cscr$, $\overline\Cscr$ & \leavevmode{Images $\mathbf\Xi(\Xbb)$ and $\mathbf\Xi(\overline\Xbb)$ of the spacetime minor map} \\
$\bCscr$ & \leavevmode{Independent triples with $G$ antisymmetric and $\iota_{\particlecurrent}F=\iota_{\particlecurrent}G=0$} \\
${\mathcal Z_{\raug}}$ & \leavevmode{Triples in $\operatorname{Dom}(\raug)\cap\bCscr_+$ together with their constitutive stress} \\
${\Cl(\mathcal Z_{\raug})}$ & \leavevmode{Closure of the constitutive flux graph $\mathcal Z_{\raug}$} \\
${\partial\Cscr}$, ${\partial\mathcal Z_{\raug}}$ & \leavevmode{Kinematic and constitutive boundaries defined in~\eqref{eq:kinematic-constitutive-boundaries}} \\
\bottomrule
\end{tabular}
\caption{Spacetime minors and limiting values of the stress.}
\label{tab:kinematic-constitutive-sets}
\endgroup
\end{table}

\paragraph{Definition of the augmented system.}

We now specify how the constitutive law closes the balance laws.

\begin{definition}
\label{propo=30DF}
Fix a covariant \underline{\emph{constitutive energy}} $\raug(F,G,\particlecurrent;g)$ and its admissible positive-density domain $\operatorname{Dom}(\raug)\cap\bCscr_+$. On this domain, the constitutive relation $T=T_{\raug}(F,G,\particlecurrent;g)$ is {prescribed by the formula~\eqref{augment3} below. }At the constitutive boundary, prescribe the finite-flux graph ${\Cl(\mathcal Z_{\raug})}$ defined above. The unknowns of the \underline{\emph{augmented system}} are $(F,G,\particlecurrent)$. The stress $T$ is given by the constitutive relation in matter and its specified continuation at the constitutive boundary. Its equations are
\be
\label{augment1}
\nabla_\nu T^{\mu\nu}=0,\qquad (F,G,\particlecurrent,T)\in{\Cl(\mathcal Z_{\raug})},
\ee
\be
\label{augment2}
\nabla_{[\mu}\Fb^A_{\nu]} = 0,
\qquad
\nabla_{[\lambda}\Gb_{\mu\nu]A} = 0,
\ee
with the algebraic constraints defining $\bCscr$. 
\end{definition}

{For a perfect-fluid law, the stress and current equations form the reduced fluid formulation. The form equations can be  included as additional kinematic equations or omitted when these forms are not used.}
Choosing $\raug$ and its admissible flux extension closes the model, including in regions of cavitation, where the stress is specified by this constitutive continuation even when derivatives of $\raug$ with respect to $\particlecurrent$ do not exist.
The stress balance laws remain equalities in the sense of distributions. They are the conservation identity required by Einstein's equations when $T$ is their total matter source. For augmented weak solutions, particle production is instead prescribed by $\nabla_\mu {\particlecurrent}^\mu\ge0$ as the mathematical entropy condition of Definition~\ref{def:mathematical-entropy}. Below we prove conservation for exact minors and, under additional hypotheses, for sufficiently regular constitutive fields.

%------------------------------------------------------------

\paragraph{The augmented stress--energy tensor.}
{Next, we express the energy--momentum tensor in terms of the independent augmented fields.} Starting from the identity \eqref{eq:useful} expressed in terms of $\rphys$, namely
\be
\label{eq:tensor-elast}
T_\mu^{\ \nu}
= - \rphys\,\delta_\mu^{\ \nu}
+ F^A_\mu\,\frac{\del \rphys}{\del F^A_\nu},
\ee
{we use the constitutive energy $\raug(\bF,\bG,\bJ)$ in independent fields, with configuration energy $\rphys(F)=\raug(F,G(F),\particlecurrent(F))$. The chain rule gives}
\be
\frac{\del \rphys}{\del F^A_\nu}(F) = \frac{\del \rhob}{\del \bF^A_\nu}(F,G(F),\particlecurrent(F)) + \frac{\del \rhob}{\del \bG_{\alpha\beta B}}(F,G(F),\particlecurrent(F)) \frac{\del G_{\alpha\beta B}}{\del F^A_\nu}(F) + \frac{\del \rhob}{\del \bJ^\alpha}(F,G(F),\particlecurrent(F)) \frac{\del {\particlecurrent}^\alpha}{\del F^A_\nu} (F).
\ee
Here \(\frac{\del\rhob}{\del \bF}\), \(\frac{\del\rhob}{\del \bG}\), \(\frac{\del\rhob}{\del \bJ}\) denote partial derivatives of \(\rhob\) with respect to the corresponding blocks of variables. Substituting into \eqref{eq:tensor-elast} yields
\be
\aligned
T_\mu^{\ \nu}(F)  
&= - \rhob \,\delta_\mu^{\ \nu}
+ F^A_\mu\,\frac{\del \rhob}{\del \bF^A_\nu}
+ F^A_\mu\,\frac{\del \rhob}{\del \bG_{\alpha\beta B}} \frac{\del G_{\alpha\beta B}}{\del F^A_\nu}
+ F^A_\mu\,\frac{\del \rhob}{\del \bJ^\alpha} \frac{\del {\particlecurrent}^\alpha}{\del F^A_\nu}.
\endaligned
\ee

\begin{lemma}[Contraction identities]
\label{lem:algebraic-identities}
The following identities express the contraction of the derivatives of $G(F)$ and $\particlecurrent(F)$ with $F^A_\mu$ over the material index $A$. The remaining indices are spacetime indices:
\be
\label{eq:601}
\frac{\del G_{\mu'\nu' A'}}{\del F^A_\nu} F^A_\mu
=
-2\,\delta_{[\mu'}^\nu G_{\nu']\mu A'},
\qquad
\frac{\del {\particlecurrent}^{\mu'}}{\del F^A_\nu} F^A_\mu
=
-2\,\delta_\mu^{[\mu'} {\particlecurrent}^{\nu]}.
\ee
\end{lemma}

\begin{proof}
\bse
We use normalized antisymmetrization and the full-index-sum convention $d\rhob=P_A{}^\nu d\bF_\nu^A+Q^{\alpha\beta A}d\bG_{\alpha\beta A}+B_\alpha d\bJ^\alpha$, with $Q$ antisymmetric. Thus a derivative with respect to an independent component $G_{\alpha\beta A}$, $\alpha<\beta$, is $2Q^{\alpha\beta A}$. Differentiating the two factors of $G$ gives
\be
\frac{\del G_{\alpha\beta C}}{\del F_\nu^A}F_\mu^A
=\delta_\alpha^\nu G_{\mu\beta C}+\delta_\beta^\nu G_{\alpha\mu C}
=-2\delta_{[\alpha}^\nu G_{\beta]\mu C}.
\ee
For the cubic current, let $L$ be an invertible endomorphism of the spacetime tangent space at the point in question, with the metric volume held fixed. Pullback of the three-form gives $\particlecurrent(F\circ L)=(\det L)L^{-1}{\particlecurrent}(F)$. We differentiate at $L=I$ in the direction with entries $\delta L^\beta{}_{\alpha}=\delta_\mu^\beta\delta_\alpha^\nu$ to obtain
\be
\frac{\del {\particlecurrent}^\alpha}{\del F_\nu^A}F_\mu^A
=\delta_\mu^\nu {\particlecurrent}^\alpha-\delta_\mu^\alpha {\particlecurrent}^\nu
=-2\delta_\mu^{[\alpha}{\particlecurrent}^{\nu]}.
\ee
Its trace in $\mu=\nu$ is $3{\particlecurrent}^\alpha$, consistently with cubic homogeneity. This distinguishes the degree three of the polynomial from the coefficient two in the antisymmetrized contraction. Substitution in the constitutive chain rule determines both signs in the augmented stress formula.
\ese
\end{proof}

Substitution of the contraction identities into the constitutive chain rule gives the following stress tensor, agreeing with the original tensor on $\Cscr$. For an augmented fluid law depending only on $\particlecurrent$ and $g$, the derivatives below with respect to $F$ and $G$ vanish and this formula reduces to the perfect-fluid stress.

\begin{proposition}[Augmented stress--energy tensor]
The augmented stress--energy tensor $\bT=\bT(\bF, \bG, \bJ)$ is given by
\be
\label{augment3}
\Tb_\mu^{\ \nu}=-\rhob\,\delta_\mu^{\ \nu}
+\Fb^A_\mu\frac{\del\rhob}{\del\Fb^A_\nu}
-2\frac{\del\rhob}{\del\Gb_{\alpha\beta B}}\delta_{[\alpha}^{\nu}\Gb_{\beta]\mu B}
-2\frac{\del\rhob}{\del\Jb^\alpha}\delta_\mu^{[\alpha}\Jb^{\nu]}.
\ee
\end{proposition}

%-----------------------------------------------------------------------------------------------------------------------------------------

\subsection{Continuity equation}

\paragraph{Formulation in the deformation variable.}
For fields given by a configuration map, the current is the dual of the pullback of the material volume form. Closure of this three-form expresses particle conservation.

\begin{proposition}[Standard continuity equation]  
\label{prop-continuity}
For every admissible configuration map $f:\Mcal\to\Bcal$, its cubic minor current $\particlecurrent$ is conserved in the sense of distributions:
\be
\label{equa-cont}
\nabla_\mu {\particlecurrent}^\mu = 0.
\ee
{On the timelike-current set, the corresponding density--velocity form is}
\be
\label{equa-cont-2}
\nabla_\mu( n u^\mu) = 0.
\ee
\end{proposition}

For independent augmented fields, particle conservation requires a separate argument. The weak formulation then imposes the inequality for the current $-\particlecurrent$.

\begin{proof}
The current is determined by $i_{\particlecurrent}\mathrm{vol}_g=\omega=f^*\Omega$. Taking the exterior derivative gives $(\nabla_\mu {\particlecurrent}^\mu)\mathrm{vol}_g=d\omega=f^*(d\Omega)=0$, since $\Omega$ has top degree on $\Bcal$. The equality at Lipschitz regularity is proved in Proposition~\ref{prop:kinematic-weak-local} below.
\end{proof}

%--------------------------------------------------------

\paragraph{Conservation of the material current.}

For independent minors, we contract the stress balance laws with the particle current to obtain its conservation in matter. We then use the regularity of the current at its zero set. For a current $\bJ$ that is future timelike or zero, we write
\be
\Mcal_+\coloneqq\{\bJ\ne0\}.
\ee
On this open matter set, the energy $\brho$ has the derivatives used below. Its current-conjugate one-form is
\be
\bB_\mu\coloneqq\frac{\del\brho}{\del\bJ^\mu}(\bF,\bG,\bJ)
\qquad\hbox{on }\Mcal_+.
\ee
Writing $\bJ=nu$ with $n>0$, we vary its magnitude at fixed $u,\bF,\bG$. This gives $\bB_\mu\bJ^\mu=n\,\partial_n\bar\rho$ in the rest variables of Section~\ref{section=4-polyconvexity} below. Thus positive constitutive enthalpy implies the following nonvanishing condition in matter:
\be
\label{eq:5j0}
\bB_\mu\bJ^\mu\ne0\qquad\hbox{on }\Mcal_+.
\ee

\paragraph{{Continuation of the particle current.}}
{We now ask when the current is conserved across its zero set.}
\begin{proposition}[Augmented continuity at zero current]
\label{prop:aug-continuity}
Let $(\bF,\bG,\bJ,\bT)$ be a $C^1$ solution of the augmented balance laws, with $\bJ$ future timelike or zero. On $\Mcal_+$, assume the covariant constitutive law is {$C^1$}, its stress is~\eqref{augment3}, and~\eqref{eq:5j0} holds. The current satisfies
\be
\label{cont2}
\nabla_\mu\bJ^\mu=0
\ee
throughout the region, including the zero-current set. The same conclusion holds distributionally when the fields are continuous and belong to $W^{1,1}_{\rm loc}$, with the matter hypotheses satisfied almost everywhere on $\Mcal_+$. {It also holds for continuous fields in $BV_{\rm loc}$ when the constitutive law is {$C^1$} on $\Mcal_+$ and~\eqref{eq:5j0} holds at every point there; the balance laws and the conclusion are then understood as equalities of Radon measures.}
\end{proposition}

\begin{proof}
\bse
We first work on $\Mcal_+$ and suppress the bars on the fields. We set $P_A{}^\nu\coloneqq\del\brho/\del\bF^A_\nu$, $Q^{\alpha\beta A}\coloneqq\del\brho/\del\bG_{\alpha\beta A}$, and $B_\alpha\coloneqq\del\brho/\del\bJ^\alpha$, with $Q$ antisymmetric. We use full antisymmetric sums, so that
\be
d\raug=P_A{}^\nu\,dF_\nu^A+Q^{\alpha\beta A}\,dG_{\alpha\beta A}+B_\alpha\,d{\particlecurrent}^\alpha.
\ee
The constitutive algebraic identities give
\be
{\particlecurrent}^\mu T_\mu{}^\nu=-\raug {\particlecurrent}^\nu.
\ee
Closure and the orthogonality condition~\eqref{eq:450bis-2} give $\mathcal L_{\particlecurrent}F=\mathcal L_{\particlecurrent}G=0$ by Cartan's identity, or
\be
{\particlecurrent}^\lambda\nabla_\lambda F_\nu^A =-F_\lambda^A\nabla_\nu {\particlecurrent}^\lambda,\qquad {\particlecurrent}^\lambda\nabla_\lambda G_{\alpha\beta A} =-G_{\lambda\beta A}\nabla_\alpha {\particlecurrent}^\lambda -G_{\alpha\lambda A}\nabla_\beta {\particlecurrent}^\lambda.
\ee
Covariance, material homogeneity and $\nabla g=0$ permit the constitutive chain rule:
\be
{\particlecurrent}^\lambda\nabla_\lambda\raug
=-P_A{}^\nu F_\lambda^A\nabla_\nu {\particlecurrent}^\lambda
 -2Q^{\alpha\beta A}G_{\lambda\beta A}\nabla_\alpha {\particlecurrent}^\lambda
 +B_\alpha {\particlecurrent}^\lambda\nabla_\lambda {\particlecurrent}^\alpha.
\ee
\ese
Differentiating the identity for ${\particlecurrent}^\mu T_\mu{}^\nu$ and substituting~\eqref{augment3} cancels all derivative terms except
\be
\label{eq:JdivT-corrected-local}
{\particlecurrent}^\mu\nabla_\nu T_\mu{}^\nu
=-(B_\mu {\particlecurrent}^\mu)\nabla_\nu {\particlecurrent}^\nu
\qquad\hbox{on }\Mcal_+.
\ee
Stress balance laws and~\eqref{eq:5j0} therefore give $\nabla_\nu {\particlecurrent}^\nu=0$ on $\Mcal_+$.  If the energy depends explicitly on the material position $f$, the extra term in ${\particlecurrent}^\lambda\nabla_\lambda\raug$ is $\raug_{X^A}{\particlecurrent}^\lambda\partial_\lambda f^A$. It vanishes for a transported material position, in particular for $F=df$ by the orthogonality condition~\eqref{eq:450bis-2}. In the independent formulation used here the law is homogeneous and the labels are unnecessary for this identity. {This calculation uses the chain rule for $\raug$ and the product rule for ${\particlecurrent}^\mu T_\mu{}^\nu$; the conjugates $P,Q,B$ enter only as continuous coefficients. The regularity of $T$ is part of the field hypotheses, rather than a consequence of $C^1$ regularity of the law.}

We take an interior point $p$ with $\particlecurrent(p)=0$. In a future-adapted coordinate chart, causality gives ${\particlecurrent}^0\ge0$ and $|{\particlecurrent}^{\rm sp}|\le C {\particlecurrent}^0$ locally. The nonnegative $C^1$ function ${\particlecurrent}^0$ has a minimum at $p$, so all its first derivatives vanish there. The component bound then gives ${\particlecurrent}^i(p+h)=o(|h|)$, and hence every first derivative of $\particlecurrent$ vanishes at $p$. The connection terms vanish as well because $\particlecurrent(p)=0$. This proves~\eqref{cont2} on the zero set without evaluating any constitutive derivative there.
For continuous fields in $W^{1,1}_{\rm loc}$, the matter set is open and the product and constitutive chain rules hold almost everywhere on compact subsets of it. The stress balance laws therefore give $\nabla_\mu {\particlecurrent}^\mu=0$ almost everywhere on $\Mcal_+$. Every first weak derivative of a $W^{1,1}_{\rm loc}$ scalar function vanishes almost everywhere on each of its level sets. Applying this fact to the components of $\particlecurrent$ gives $\nabla {\particlecurrent}=0$ almost everywhere on $\{\particlecurrent=0\}$.  Thus~\eqref{cont2} holds almost everywhere and in the sense of distributions, including across the zero set of $\particlecurrent$.

{For continuous $BV_{\rm loc}$ fields, the product and chain rules give~\eqref{eq:JdivT-corrected-local} as an equality of measures on $\Mcal_+$. The continuous nonzero factor $B_\mu {\particlecurrent}^\mu$ and the balance laws imply $\operatorname{div}{\particlecurrent}=0$ there. For a continuous scalar $BV$ function $v$, we obtain vanishing of the derivative on $\{v=0\}$ by applying the chain rule to smooth $\phi_\epsilon\to0$ uniformly with bounded $\phi_\epsilon'\to\mathbf1_{\{0\}}$. Thus $D{\particlecurrent}^\mu=0$ on $\{\particlecurrent=0\}$, as are the connection terms, proving the result.}
This completes the proof of Proposition~\ref{prop:aug-continuity}.
\end{proof}

The continuity assumption in the $BV$ result matters at shocks. For an independent current, a finite-flux limit with nonzero null current also requires a separate conservation argument. In the positive-rest-mass class of Theorem~\ref{thm:closed-vacuum-entropy} below, zero proper density forces $\particlecurrent=0$ at finite-flux limits, where stress may persist. When the minors come from a locally Lipschitz configuration, current conservation follows directly from the minor relations (cf.~Proposition~\ref{prop:kinematic-weak-local} below).

\paragraph{Propagation of the minor relations.}
The transport equations also determine whether the algebraic identities imposed on the initial slice hold at later times.

\begin{corollary}[Propagation of the minor relations with cavitation]
\label{cor:physical-minors-local06}
Use fixed material volume coordinates and a smooth spacelike foliation, with time function $t$ and initial slice $\Sigma$. In coordinates adapted to the foliation, let $\beta$ be a bounded coordinate velocity with spatial Lipschitz constant locally integrable in time on flow tubes from $\Sigma$. Put $w=\del_t+\beta^i\del_i$ and assume that its flow from $\Sigma$ covers the region. Let $F^A,G_A,\particlecurrent$ be locally bounded fields with prescribed bounded initial data, satisfying
\be
dF^A=dG_A=0,\qquad
i_wF^A=i_wG_A=0,\qquad {\particlecurrent}={\particlecurrent}^0w,\qquad
\nabla_\mu {\particlecurrent}^\mu=0 .
\ee
The closure and conservation laws hold in the sense of distributions with initial tangential forms $F_{\rm in}^A,G_{{\rm in},A}$ and initial particle three-form $\omega_{\rm in}$ on $\Sigma$. Require orthogonality to $w$ also at zero current. If $G_{{\rm in},A}=\tfrac12\Omega_{ABC}F_{\rm in}^B\wedge F_{\rm in}^C$ and $\omega_{\rm in}=\tfrac16\Omega_{ABC}F_{\rm in}^A\wedge F_{\rm in}^B\wedge F_{\rm in}^C$, then:
\begin{itemize}
\item[(i)] \emph{Exact minors.} One has $G=G(F)$ and $\particlecurrent=\particlecurrent(F)$ almost everywhere in the flow region.
\item[(ii)] \emph{Ranks and zero current.} For almost every initial position, the ranks of $F$ and $G$ agree with their initial ranks at almost every time along the trajectory. The flow transports the zero set of $\omega_{\rm in}$ to the zero set of $\particlecurrent$, up to null sets on almost every slice.
\end{itemize}
In particular, the conclusion applies for a locally Lipschitz field $u$ with $dt(u)>0$, using $w=u/dt(u)$, and allows nonzero lightlike currents.
\end{corollary}

{The same kinematic conclusion applies to fluid configurations when the auxiliary forms $F,G$ are  included and the hypotheses of the corollary are satisfied; a fluid formulation may instead omit these form variables.}

\begin{proof}
Let $X_t$ be the spatial flow of $\beta$, and put $\Psi(t,y)=(t,X_t(y))$. The flow and its inverse are locally bi-Lipschitz up to the initial slice, by the spatial Gronwall bound and boundedness of $\beta$ (cf.~the proof of Proposition~\ref{cor:weak-flow-coordinates-vacuum} below). Its spacetime derivative sends $\del_t$ to $w$ almost everywhere. Theorem~\ref{un-prop:weak-coordinate} below transports the closure laws for $F,G$ to the coordinates $(t,y)$. The pullbacks have only spatial components, by their orthogonality to $w$. Their closure therefore makes their spatial components independent of $t$, with values given by the prescribed initial forms.

Similarly, the transformed vector density $|g|^{1/2}{\particlecurrent}$ has zero spatial components because $\particlecurrent={\particlecurrent}^0w$. Its conservation law in the flow chart makes its time component independent of $t$ with the prescribed initial value. Equivalently, $\omega_{\particlecurrent}=i_{\particlecurrent}\mathrm{vol}_g$ pulls back to the initial spatial current three-form. In Eulerian coordinates its time evolution reads
\be
\del_t(|g|^{1/2}{\particlecurrent}^0)+\del_i(|g|^{1/2}{\particlecurrent}^0\beta^i)=0,
\qquad |g|=|\det(g_{\mu\nu})|.
\ee

Pullback commutes with exterior products almost everywhere. Since the material volume coefficients are fixed, the initial quadratic and cubic relations give $\Psi^*G_A=G_A(\Psi^*F)$ and $\Psi^*\omega_{\particlecurrent}=\omega(\Psi^*F)$. Invertibility of the coordinate change gives the two claimed identities. On the slices the derivative of $X_t$ and its exterior powers are invertible almost everywhere. This preserves the ranks of the spatial restrictions of $F,G$, which are also their spacetime ranks by the orthogonality condition~\eqref{eq:450bis-2}. The spatial current three-form vanishes exactly where ${\particlecurrent}^0=0$, so its pullback preserves the zero set of $\particlecurrent$. All arguments concern the flow of $w$ and apply also when $\particlecurrent=0$. Closure gives local potentials for $F$; a global configuration requires compatible material coordinates.
\end{proof}

{The same transport argument also applies without initially imposing the exact minor relations. Under the kinematic and regularity hypotheses of Corollary~\ref{cor:physical-minors-local06}, set $\omega_{\particlecurrent}=i_{\particlecurrent}\varepsilon_g$ and consider the mixed tensor
\be
\label{eq:transport-mixed-minor-defect}
 D^C_{\mu\nu\rho\sigma}\coloneqq
 \Omega^{ABC}G_{\mu\nu A}G_{\rho\sigma B}
 -(\omega_{\particlecurrent})_{\rho\sigma\mu}F^C_\nu
 +(\omega_{\particlecurrent})_{\rho\sigma\nu}F^C_\mu.
\ee
In the flow coordinates $F^A$, $G_A$ and $\omega_{\particlecurrent}$ have only spatial components, all constant in time. Their tensor products have the same property. Thus if the spatial restriction of $D$ vanishes initially, $D$ vanishes throughout the flow region, including where $\particlecurrent=0$. This propagation does not require $G=G(F)$ or $\particlecurrent=\particlecurrent(F)$; the stability results below impose those relations {where explicitly required}.}

At zero current, the transport argument determines the minors from their initial values along the prescribed flow. The admissible stress at the transported fields is specified by the completed constitutive flux graph.

%============================================================================

% Restore the current-alias colors used in the remaining revision.
\let\Jt\EntropySavedJt
\let\Jtilde\EntropySavedJtilde
\let\Jcausal\EntropySavedJcausal
\let\Jb\EntropySavedJb
\let\Jh\EntropySavedJh
\let\bJ\EntropySavedbJ
\let\hJ\EntropySavedhJ

\section{\texorpdfstring{Hyperbolicity of the two proposed formulations}{Hyperbolicity of the two proposed formulations}}
\label{section=4-hyperbolicity}

%-----------------------------------------------------------------------------------------------------------------------------------------

\subsection{Hyperbolicity of the deformation gradient formulation}
\label{subsec:hyp-firstorder}

\paragraph{Observer decomposition.}
{
We now examine the characteristic directions of the two first-order formulations and compute their propagation speeds in matter. The observer decomposition makes it possible to compare the acoustic cones and to distinguish elastic waves from the modes associated with the constraints.\par}
To distinguish  the frequency measured by the observer from the spatial direction of propagation, we first choose an observer. Let $t_\mu$ be a timelike covector field, positive on future timelike vectors. The corresponding observer is the future unit vector $\slicenormal$, and $u$ denotes the material four-velocity. We set
\be
{\slicenormal}_\mu\coloneqq -\frac{t_\mu}{\sqrt{-g^{\alpha\beta}t_\alpha t_\beta}},
\qquad
g_{\mu\nu}{\slicenormal}^\mu {\slicenormal}^\nu=-1.
\ee
The $g$-orthogonal projection onto ${\slicenormal}^\perp\subset T\Mcal$ has components
\be
\pi^\mu{}_\nu\coloneqq \delta^\mu{}_\nu+{\slicenormal}^\mu {\slicenormal}_\nu.
\ee For any spacetime covector $V_\mu$ we use the splitting
\be
V_\mu = (-{\slicenormal}^\nu V_\nu)\,{\slicenormal}_\mu + \pi^\nu{}_\mu V_\nu,
\qquad
V_{\mathrm{time}}\coloneqq -{\slicenormal}^\nu V_\nu,
\qquad
(V_{\mathrm{space}})_\mu\coloneqq \pi^\nu{}_\mu V_\nu.
\ee
In particular, each material covector $F^A{}_\mu$ decomposes as
\be
F^A{}_\mu = (F_{\mathrm{time}})^A {\slicenormal}_\mu + (F_{\mathrm{space}})^A{}_\mu,
\qquad
(F_{\mathrm{time}})^A\coloneqq -{\slicenormal}^\nu F^A{}_\nu,
\qquad
(F_{\mathrm{space}})^A{}_\mu\coloneqq \pi^\nu{}_\mu F^A{}_\nu.
\ee
In an orthonormal frame with time vector $\slicenormal$, the components of the material four-velocity are $u=(1,v)/\sqrt{1-|v|^2}$, where $v$ is the spatial material velocity measured by $\slicenormal$.

\paragraph{Rank--one variations and acoustic forms.}
\label{subsubsec:rankone-F}

We fix an admissible state $(x,F)\in\Mcal\times\Xbb$ and a material direction $U\in T_X\Bcal$. Given a spacetime covector $\xi\in T_x^\ast\Mcal$, we consider the rank--one variation
\be
\delta F^A{}_\mu \coloneqq  U^A\,\xi_\mu.
\ee
We freeze the constitutive coefficients at $F$. The principal symbol is obtained by linearizing the field equations in the rank-one direction $U\otimes\xi$ and contracting with the phase covector $\xi$.

\begin{definition}
\label{def:acoustic-objects}
Assume $\rphys$ is $C^2$ near $F\in\Xbb$. The \underline{\emph{material acoustic matrix}} is
\bse
\be
\label{eq:material-acoustic-matrix}
\mathbb A_{AB}(\xi)\coloneqq\frac{\del^2\rphys}{\del F^A_\mu\,\del F^B_\nu}(F;g)\,\xi_\mu\xi_\nu.
\ee
For a nonzero material direction $U$, its \underline{\emph{acoustic quadratic form}} and \underline{\emph{acoustic tensor}} are defined by
\be
\label{eq:acoustic-form-and-tensor}
q_U(\xi)\coloneqq U^T\mathbb A(\xi)U=k(U,U)^{\mu\nu}\xi_\mu\xi_\nu.
\ee
\ese
The tensor $k(U,U)^{\mu\nu}$ is symmetric and contravariant. When it is nondegenerate, its inverse $k(U,U)_{\mu\nu}$ is the \underline{\emph{acoustic metric}} for that direction.
\end{definition}

{The matrix $\mathbb A(\xi)$ maps a material amplitude $U$ to the principal response for the phase covector $\xi$. In its definition the two spacetime indices of the constitutive Hessian are contracted with $\xi$, leaving the two material indices $A,B$.}
{For a perfect fluid, the energy depends only on the proper particle density $n$. The reduced fluid equations have the compression and zero-speed modes computed in Section~\ref{subsec:propagation-speeds} below.}
Substitution of $\delta F=U\otimes\xi$ into the principal field equations gives $\mathbb A(\xi)U=0$. Multiplication by $U^T$ gives $q_U(\xi)=0$. We use the signs of $q_U$ in time and spatial directions to express the acoustic conditions; characteristic covectors and propagation speeds are determined by $\det\mathbb A(\xi)=0$. The matrix equation must be satisfied by the wave amplitude, rather than only its scalar contraction.

%-----------------------------------------------------------------------------------------------------------------------------------------

\subsection{Hyperbolicity of the augmented formulation}
\label{subsec:hyp-augmented}

\paragraph{Observer decomposition for the augmented variables.}
\label{subsubsec:observer-split-Xi}

Consider the augmented system in the independent variables $(\bF,\bG,\bJ)\in\bCscr$. Hyperbolicity is governed by the principal symbol of the first-order balance laws \eqref{augment1}--\eqref{augment2}, hence by the Hessian of the augmented energy and by the linearized constraint map (propagation map).  Fix a timelike conormal with future unit normal $\slicenormal$ and spatial projection $\pi$. We write $(F,G,\particlecurrent)$ for the independent triple in the following observer decomposition.
\begin{itemize}
\item[(i)] \emph{Configuration gradient.}
\be
(F_{\mathrm{time}})^A\coloneqq -{\slicenormal}^\mu F^A{}_\mu,
\qquad
(F_{\mathrm{space}})^A{}_\mu\coloneqq \pi^\nu{}_\mu F^A{}_\nu.
\ee
\item[(ii)] {\emph{Current.}}
\be
{\particlecurrent}_{\mathrm{time}}\coloneqq -{\slicenormal}_\mu {\particlecurrent}^\mu,
\qquad
({\particlecurrent}_{\mathrm{space}})^\mu\coloneqq \pi^\mu{}_\nu {\particlecurrent}^\nu.
\ee
\item[(iii)] {\emph{Quadratic minors.} For each material index $A$, the two-form decomposes into}
\be
(G_{\mathrm{time}})_{\mu A}\coloneqq {\slicenormal}^\nu G_{\nu\mu A},
\qquad
(G_{\mathrm{space}})_{\mu\nu A}\coloneqq \pi^\alpha{}_\mu\pi^\beta{}_\nu G_{\alpha\beta A}.
\ee
\end{itemize}
We collect the components of the augmented fields into two groups:
\be
\Xi_{\mathrm{time}}\coloneqq (F_{\mathrm{time}},\,G_{\mathrm{time}},\,{\particlecurrent}_{\mathrm{space}}),
\qquad
\Xi_{\mathrm{space}}\coloneqq (F_{\mathrm{space}},\,G_{\mathrm{space}},\,{\particlecurrent}_{\mathrm{time}}).
\ee
The two groups have $15$ and $19$ components before the algebraic constraints are imposed. When $\particlecurrent=r(1,v)$ with $r>0$, orthogonality gives $F_0^A=-v^iF_i^A$ and $G_{0jA}=-v^iG_{ijA}$. The $19$ spatial fields and the three components of $v$ therefore determine the triple of augmented fields in matter.
\paragraph{Propagation map and acoustic tensor.}
\label{subsubsec:propmap-acoustic}

We fix a rank-three deformation gradient $F\in\Xbb$ and set $\Xi=\mathbf{\Xi}(F)=(F,G(F),\particlecurrent(F))\in\Cscr$. Given a material direction $U\in T_X\Bcal$ and a spacetime covector $\xi\in T_x^\ast\Mcal$, we consider the rank--one variation $\delta F^A{}_\mu\coloneqq U^A\xi_\mu$ and define the propagation map
\be
\Sigma(\xi,U)\coloneqq D\mathbf{\Xi}(F)[\delta F]
=\bigl(\delta F,\ DG(F)[\delta F],\ D\particlecurrent(F)[\delta F]\bigr)\in T_\Xi\Cscr\subset\RR^{34}.
\ee
Differentiating the polynomial minor formulas gives the components of this map:
\be
\aligned
&\Sigma_F(\xi,U)^A{}_\mu=U^A\,\xi_\mu,\qquad \Sigma_G(\xi,U)_{\mu\nu A}=2\,U^B\,\xi_{[\mu} F_{\nu]}^C\,\Omega_{ABC},\\
&\Sigma_{\particlecurrent}(\xi,U)^\mu=-\frac12\,U^A\xi_\nu F_\sigma^B F_\tau^C\,\Omega_{ABC}\,\eps^{\mu\nu\sigma\tau}.
\endaligned
\ee
We assume that the constitutive law factors through an augmented energy
\be
\rphys(F;g)=\raug(\Xi;g)\big|_{\Xi=\mathbf{\Xi}(F)},
\qquad \raug\in C^2(\mathcal E),
\ee
for some open neighborhood $\mathcal E\subset\RR^{34}$ of $\Cscr$. Along $F+sU\otimes\xi$, all minors are affine in $s$, since every higher variation alternates two identical rank-one factors. Thus $D^2\mathbf\Xi(F)[U\otimes\xi,U\otimes\xi]=0$, and the second chain rule gives
\be
\label{eq:acoustic-form-augmentation}
q_U(\xi)=D_F^2\rphys[U\otimes\xi,U\otimes\xi]
=D^2_{\Xi\Xi}\raug(\Xi)[\Sigma(\xi,U),\Sigma(\xi,U)].
\ee
Thus the acoustic forms agree on compatible rank-one variations. The augmented modes also involve the constraints. The signature and the value on a time covector give the following two conditions.

\begin{definition}
\label{def:LR1-firstorder}
\label{def:LR1-aug}
\label{def:uniform-time}
\begin{itemize}
\item[(i)] The \underline{\emph{Lorentzian acoustic condition}} holds if, for every $F\in\Xbb$ and every nonzero material direction $U$, the tensor $k(U,U)^{\mu\nu}$ has signature $(-,+,+,+)$\footnote{Spatial rank-one positivity is the classical Legendre--Hadamard condition; regular hyperbolicity additionally requires compatible time and spatial positivity~\cite{BeigSchmidt,Beig2023}. Here the terminology refers to the signature of the acoustic form for each fixed material direction.}. Its null cone is
\be
\mathcal C_U\coloneqq\{\xi\in T_x^*\Mcal\setminus\{0\}:k(U,U)^{\mu\nu}\xi_\mu\xi_\nu=0\}.
\ee

\item[(ii)] A $g$-timelike covector field $t$ is a \underline{\emph{uniform acoustic time direction}} on $\mathcal K\subset\Cscr$ if, for some $c>0$,
\be
k(U,U)^{\mu\nu}(\Xi)t_\mu t_\nu\le-c,\qquad
\Xi\in\mathcal K,\quad |U|_h=1,
\ee
where $h$ is a fixed auxiliary Riemannian metric on the material space.
\end{itemize}
\end{definition}

{The original formulation and its augmentation give the same condition by~\eqref{eq:acoustic-form-augmentation}. A plane-wave amplitude must also satisfy $\mathbb A(\xi)U=0$.}
{A perfect fluid has zero shear-wave speeds, so the strictly elastic signature in (i) is not imposed on its discarded shear directions. Its reduced density--velocity equations have speeds $0$ and $\pm c_s$ in the fluid rest frame.}
Condition (ii) selects a time direction uniformly timelike for the acoustic forms on $\mathcal K$. A symmetric-hyperbolic reduction also requires a positive time matrix for the constrained principal system.

%-----------------------------------------------------------------------------------------------------------------------------------------

\subsection{Propagation speeds and their multiplicities}
\label{subsec:propagation-speeds}

\paragraph{{Characteristic roots and acoustic speeds.}}
The determinant of the material acoustic matrix gives the characteristic covectors. We freeze the coefficients at a material state and use an observer orthonormal coframe. For a unit spatial covector $\widehat k$ and a plane wave with frequency $\omega$ and nonzero spatial frequency $k=|k|\widehat k$, the physical phase speed $s_{\rm ph}\coloneqq\omega/|k|$ satisfies
\be
\label{eq:physical-wave-speed-polynomial}
\det\mathbb A(-s_{\rm ph},\widehat k)=0.
\ee
Its order as a root gives the algebraic multiplicity. We now count the modes at a homogeneous rest configuration with positive density, exact minors and invertible spatial configuration gradient $\Frest =(F_i^A)$, nonzero constitutive enthalpy, a positive definite inertia matrix and positive spatial acoustic matrix. These are constitutive calculations on the positive-density material domain, not regularity assumptions on the weak solutions.

We assume that the energy is even in the time components of $F$ at fixed spatial gradient, as it is for the quadratic and {Ogden-type} laws below in the rest coframe. At rest, this symmetry makes the mixed time--space block vanish. For velocity perturbations we define
\be
\mathsf M\coloneqq-\Frest ^T\mathbb A(1,0)\Frest >0,\qquad
\mathsf K(\widehat k)\coloneqq \Frest ^T\mathbb A(0,\widehat k)\Frest >0.
\ee
The three positive squared acoustic speeds $c_\ell^2$, counted with multiplicity for $\ell=1,2,3$, are the eigenvalues of $\mathsf M^{-1/2}\mathsf K(\widehat k)\mathsf M^{-1/2}$. The physical branches are $s_{\rm ph}=\pm c_\ell$.

\begin{figure}[!htbp]
\centering
\includegraphics[width=.87\linewidth]{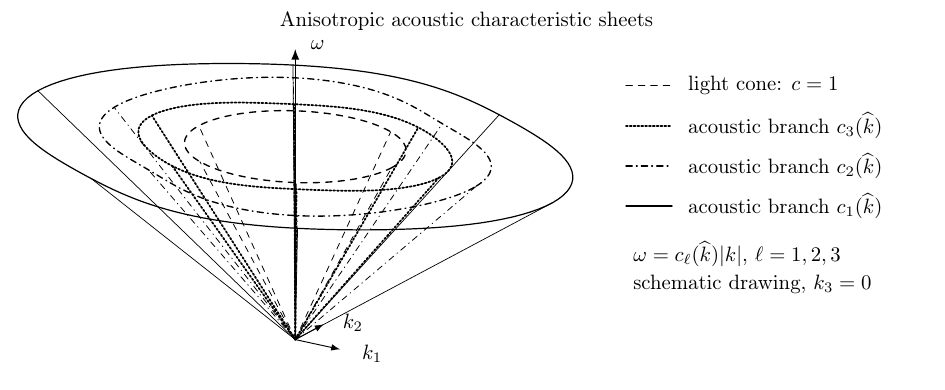}
\caption{Schematic positive-frequency characteristic sheets in the cotangent slice $k_3=0$. At a fixed direction $\widehat k$, each sheet has equation $\omega=c_\ell(\widehat k)|k|$, $\ell=1,2,3$, and the light cone has $\omega=|k|$. At fixed $\omega$, a slower branch lies farther from the frequency axis.  The cones illustrate direction-dependent wave speeds.}
\label{fig:anisotropic-acoustic-sheets}
\end{figure}

Figure~\ref{fig:anisotropic-acoustic-sheets} shows how direction-dependent phase speeds give characteristic sheets that need not be circular.

\paragraph{{Evolution symbols and constraint modes.}}
For the gradient evolution on $(F_i^A,v^i)$ there are $d=12$ variables. Before imposing the polynomial minor relations, the full augmented evolution uses the nine spatial one-form components, nine spatial two-form components, the density measured by the observer and three velocity components, hence $d=22$. For the augmented system at positive density, nonzero constitutive enthalpy allows the linearized energy row to be replaced by particle conservation. The time components of the exterior equations and the particle equation give the spatial increments associated with a velocity amplitude $b$:
\be
 \mathcal D_{12}(\widehat k)b=(\Frest b)\otimes\widehat k,\qquad
 \mathcal D_{22}(\widehat k)b=
 \bigl((\Frest b)\otimes\widehat k,\ 
      (\widehat k\cdot b)\Grest-(\Grest\widehat k)\otimes b,\ 
      n\,\widehat k\cdot b\bigr).
\ee
Here the second component uses the spatial dual of $G$, and $r=n$ at rest. Indeed, differentiating the spatial cofactor and determinant at $\Frest$ in the direction $(\Frest b)\otimes\widehat k$ gives the last two components. Let $\mathsf S_d(\widehat k)$ denote the linearized spatial momentum flux contracted with $\widehat k$. At the rest configuration the spatial Noether identity multiplies the linearized material equations by $\Frest^T$. A velocity variation gives $\delta F_0=-\Frest b$, so the time coefficient of the momentum equation is $\mathsf M$. For a spatial increment $(\Frest b)\otimes\widehat k$, the same identity gives $\mathsf S_d\mathcal D_d b=\Frest^T\mathbb A(0,\widehat k)\Frest b=\mathsf K b$. This calculation applies in both systems by~\eqref{eq:acoustic-form-augmentation}. The normalized evolution symbol and its characteristic polynomial are
\be
\label{eq:evolution-wave-multiplicities}
\mathsf P_d\coloneqq\begin{pmatrix}0&\mathcal D_d\\ \mathsf M^{-1}\mathsf S_d&0\end{pmatrix},\qquad
\det(s_{\rm ph}I_d-\mathsf P_d)
=s_{\rm ph}^{d-6}\prod_{\ell=1}^3(s_{\rm ph}^2-c_\ell^2),\qquad d\in\{12,22\}.
\ee
The Schur complement gives the determinant for $s_{\rm ph}\ne0$, then for every speed by polynomial identity. Before imposing the spatial equations, the gradient evolution has six zero roots and the augmented evolution has sixteen, in addition to the six acoustic roots.

{Since $\mathsf S_d\mathcal D_d=\mathsf K>0$, the zero eigenspace is $\ker\mathsf S_d\times\{0\}$, of dimension $d-6$, so these modes are semisimple. Put $\Pi_d\coloneqq\mathcal D_d\mathsf K^{-1}\mathsf S_d$ and $\mathsf H_d\coloneqq\mathsf S_d^T\mathsf K^{-1}\mathsf S_d+(I-\Pi_d)^T(I-\Pi_d)$. The projection $\Pi_d$ has range $\operatorname{ran}\mathcal D_d$ and kernel $\ker\mathsf S_d$; hence $\mathsf H_d>0$ and $\mathsf H_d\mathcal D_d=\mathsf S_d^T$. Thus $\operatorname{diag}(\mathsf H_d,\mathsf M)$ symmetrizes $\mathsf P_d$. To obtain a uniform lower bound at the fixed state, decompose an increment $h$ of the spatial fields as
\be
 h=\mathcal D_d b+h_0,\qquad b=\mathsf K^{-1}\mathsf S_d h,\qquad h_0=(I-\Pi_d)h.
\ee
The decomposition gives $h^T\mathsf H_d h=b^T\mathsf K b+|h_0|^2$. On the unit sphere, $\mathsf K$ has a positive minimum eigenvalue and $\mathcal D_d$ has bounded norm. This identity and the decomposition give a uniform lower bound by a positive multiple of $|h|^2$. The formula for $\mathsf H_d$ gives the upper bound on the same sphere. The frozen system is therefore strongly hyperbolic, including when acoustic speeds coincide. The symmetrizer here depends on the spatial direction.}

In the gradient formulation, spatial closure restricts the gradient increment to three amplitudes along the covector $\widehat k$. In the augmented formulation, spatial closure also requires $\delta\Gsp\,\widehat k=0$, leaving six two-form amplitudes and one density amplitude in addition to the three gradient amplitudes. This ten-dimensional field subspace contains $\operatorname{ran}\mathcal D_{22}$. The restricted momentum map onto $\mathbb R^3$ is surjective because $\mathsf S_{22}\mathcal D_{22}=\mathsf K>0$. Thus the augmented symbol with spatial closure has seven semisimple zero modes and the six acoustic modes. These seven field modes are transported with the material.

The linearized exact minor relations determine the two-form and density increments from the gradient increment. In both formulations the exact, spatially closed increments have the form $(\mathcal D_d a,b)$, with $a,b\in\mathbb R^3$. On these coordinates the symbol acts as $(a,b)\mapsto(b,\mathsf M^{-1}\mathsf K a)$. Thus each eigendirection of $\mathsf M^{-1}\mathsf K$ gives one pair of speeds. If the spatial acoustic matrix is only nonnegative and has a zero eigenvalue, its eigendirection in this sector gives a length-two Jordan block at speed zero. This conclusion concerns the deformation equations with exact minors. The reduced fluid equations omit the deformation variables.

\paragraph{Elastic materials and perfect fluids.}
{In an isotropic rest configuration, rotations about the propagation direction give the two transverse elastic modes the same speed. The longitudinal speed may differ. When the deformation breaks this symmetry, the three eigenvalues of $\mathsf M^{-1/2}\mathsf K(\widehat k)\mathsf M^{-1/2}$ need no longer coincide in the transverse directions. Proposition~\ref{prop:quadratic-vacuum-law} below verifies the acoustic bound for the quadratic elastic law in every direction. For a perfect fluid with $\bar\rho'(n)>0$ and $\bar\rho''(n)>0$, only the longitudinal pair $\pm\sqrt{n\bar\rho''(n)/\bar\rho'(n)}$ propagates relative to the fluid. The two transverse velocity modes have speed zero. For a polytropic pressure with positive rest mass the longitudinal speed also tends to zero at vacuum.}

%============================================================================

\section{\texorpdfstring{The notion of relativistic polyconvexity}{The notion of relativistic polyconvexity}}
\label{section=4-polyconvexity} 
\label{section=4-polycoercivity}

%-----------------------------------------------------------------------------------------------------------------------------------------

\subsection{Convexity in the rest frame and constitutive assumptions}
\label{subsec:transition-what}

\paragraph{The augmented energy and convexity at rest.}
{
{We now measure convexity in the material rest frame. The augmented energy $\raug(F,G,\particlecurrent;g)$ treats its arguments as independent on $\mathcal N_g$, which contains the positive-density exact minors under consideration. Restriction to $\Cscr$ gives the physical energy}}:
\be
\label{eq:augmented-energy-restriction}\rphys(F;g)=\raug\bigl(F,G(F),\particlecurrent(F);g\bigr).
\ee
In matter, let $n=\sqrt{-g(\particlecurrent,\particlecurrent)}>0$ and $u=\particlecurrent/n$. A \emph{rest frame} is an oriented orthonormal frame whose future time vector is the material velocity $u$. Write $a=1,\ldots,r_F$ and $b=1,\ldots,r_G$ for the selected one- and two-form indices, respectively. In the dual coframe, the \emph{rest variables} are given by the spatial components of the one-forms, the spatial duals of the two-forms, and the proper density:
\be\label{eq:rest-minor-notation}
 F_0^a=G_{0ib}=0,\qquad \Frest _{ai}\coloneqq F_i^a,\qquad
 G_{ijb}=\varepsilon_{ijk}\Grest _{bk},\qquad {\particlecurrent}=n u.
\ee
At zero current this decomposition is available on a selected constitutive branch with a future unit timelike limit of $u$. Depending on the law, the graph of finite fluxes can also include limits in which the spatial material speed tends to one. For fixed $g$, the \emph{rest energy} $\bar\rho(\Frest ,\Grest ,n)$ is the value of $\raug(F,G,\particlecurrent;g)$ on the rest data determined by $(\Frest ,\Grest ,n)$ through~\eqref{eq:rest-minor-notation}.
{The vanishing time components in~\eqref{eq:rest-minor-notation} follow from Definition~\ref{def:transversality}. The rest variables separate density and deformation from the velocity relative to an external observer.}

\begin{definition}
\label{active:def:dyn-poly-t}
\label{def:rest-polyconvexity}
\begin{itemize}
\item[(i)] {A covariant augmented law $\raug(F,G,\particlecurrent;g)$ is \underline{\emph{relativistically polyconvex}} if its rest energy} $\bar\rho(\Frest ,\Grest ,n)$ is convex on a specified convex domain $\mathcal D$ of the independent selected rest variables $(\Frest ,\Grest ,n)$. In matter, $\mathcal D$ is taken within the rest variables allowed by $\mathcal N_g$ and contains the rest values of the triples satisfying the exact minor relations under consideration.
\item[(ii)] A \underline{\emph{generalized fluid law}} depends on the augmented variables only through the material current $\particlecurrent$, so that $\raug(F,G,\particlecurrent;g)=\raug(\particlecurrent;g)$ with all background data fixed. A \underline{\emph{perfect fluid law}} depends on $\particlecurrent$ only through the proper particle density $n$:
\be
\label{eq:fluid-type-constitutive-law}
 \raug(\particlecurrent;g)=\bar\rho(n),\qquad n=\sqrt{-g(\particlecurrent,\particlecurrent)}.
\ee
\end{itemize}
\end{definition}
{The definition has the following consequences for admissible variations and time directions.}

\begin{itemize}
\item {\emph{Variables in which convexity holds.} Segments in $\mathcal D$ allow independent variations of the rest variables off the exact minors. Convexity in all spacetime components in a moving frame is a separate property.}
\item \emph{Convexity in matter and at the constitutive boundary.} The domain may include $n=0$, with the specified lower semicontinuous convex extension at the boundary. For the results that require strict positivity in matter, we assume that $\bar\rho$ is $C^2$ on the chosen positive-density rest domain and that $D^2\bar\rho$ is positive definite on the selected variables there. Continuity then gives $D^2\bar\rho\ge c_K I$, $c_K>0$, on each compact convex subset $K$ of the positive-density rest domain. We call this uniform convexity on $K$. No positive lower bound is imposed uniformly as vacuum or cavitation is approached.

\item \emph{Choice of observer.} For a prescribed foliation, the acoustic inequality is imposed on its time conormal. Once the enthalpy is positive and the inertia matrix is positive definite in matter, admissibility for every future timelike direction is equivalent to the bound of the acoustic speeds by the speed of light (cf.~Proposition~\ref{prop:observer-acoustic-coercivity} below).

\end{itemize}

\paragraph{Selected fields and polynomial minor relations.}
The dependence of the energy on the deformation determines which fields enter the constitutive equations.
\begin{itemize}
\item \emph{Selected components.} The numbers $r_F,r_G\in\{0,1,2,3\}$ count the selected one-forms $F^a$ and two-forms $G_b$. In fixed material coordinates, constant linear maps on the material indices of $F,G$, of ranks $r_F,r_G$, give the selected forms. We set
\be
\Xi_{\rm a}\coloneqq(F^a,G_b,\particlecurrent),\qquad
M\coloneqq4r_F+6r_G+4,\qquad {d_{\rm sel}}\coloneqq3r_F+3r_G+4.
\ee
{The augmented energy is a function $\raug(\Xi_{\rm a};g)$ of the selected fields and the metric.} The full elastic formulation has $(r_F,r_G)=(3,3)$ and $(M,d_{\rm sel})=(34,22)$; the fluid formulation has $(0,0)$ and $(4,4)$.  On the positive-density algebraic manifold $\mathcal B_{\rm a}$, the current is future timelike, the selected forms annihilate it, and $G_b$ is antisymmetric. Since ${\particlecurrent}^0>0$, the contractions determine the time components of the forms from the following $d_{\rm sel}$ coordinates:
\be
 Y=(F_i^a,G_{ijb},{\particlecurrent}^\mu),\qquad i<j\quad\text{in the two-form components}.
\ee
\item {\emph{Polynomial minor relations.} On the full set of exact minors at rest, $\Grest=\operatorname{cof}\Frest$ and $n=\det\Frest$. These polynomial relations concern pointwise gradients. The convexity condition uses independent rest variables and is invariant under spatial orthonormal changes of frame.}

\item {\emph{Acoustic forms.} On the exact minors at positive density, rest convexity makes the spatial acoustic form nonnegative. Positive inertia makes the time part negative. If the spatial form is positive definite for each nonzero material direction $U$, then $q_U$ has signature $(-,+,+,+)$ (Definition~\ref{def:LR1-firstorder}(i)). For the fluid speeds and possible loss of diagonalizability at zero acoustic speed, see Section~\ref{subsec:propagation-speeds}.}

\item \emph{Time directions at zero density.} For a fixed admissible time direction, the strict acoustic inequality is uniform on compact subsets of the positive-density range, as in Definition~\ref{def:uniform-time}(ii). In the fluid construction the regularized density and velocity equations extend to zero density, and the velocity is prescribed and evolved there as part of the complete reference data (cf.~Theorem~\ref{uf:thm:main} below).

\label{subsec:poly-motivation}

% Historical block A-015 transferred to the reserves.

\item \emph{Laws depending only on the current.} Beig and LeFloch~\cite{BeigLeFloch} allow an energy depending on the entire current in an affine spacetime. On the future timelike cone, a scalar formed from $\particlecurrent$ and $g$ alone depends only on $n$, as for a perfect fluid.
\end{itemize}

%-----------------------------------------------------------------------------------------------------------------------------------------

\subsection{The mathematical entropy current}
\label{subsec:poly-definition}

\paragraph{Supplementary conservation law, density, and flux.}
We impose the stress balance laws $\nabla_\beta T^{\alpha\beta}=0$. For sufficiently regular solutions under the constitutive assumptions of Proposition~\ref{prop:aug-continuity}, the particle current satisfies the supplementary conservation law $\nabla_\mu {\particlecurrent}^\mu=0$. {We use the opposite current to state the admissibility condition\footnote{{This sign convention is used in the formulation with particle production for the Einstein--Euler equations in}~\cite{TorusT2-1-2026,TorusT2-2-2026}.}.}
\begin{definition}
\label{def:mathematical-entropy}
The \underline{\emph{mathematical entropy current}} is $\mathcal J_\eta=-\particlecurrent$. In an orthonormal frame adapted to a spacelike slice its \underline{\emph{mathematical entropy density}} is $\eta=-{\particlecurrent}^0$, and its \underline{\emph{mathematical entropy flux}} is $(-{\particlecurrent}^i)$. An augmented weak solution is admissible for this current if $\nabla_\mu\mathcal J_\eta^\mu=-\nabla_\mu {\particlecurrent}^\mu\le0$ in the sense of distributions.
\end{definition}

{For a perfect fluid, $\particlecurrent=nu$ in matter, so $\mathcal J_\eta=-nu$ is the negative particle-number current used in the fluid formulation. For a sufficiently regular solution at positive density, with $p=n\rrestnew'(n)-\rrestnew(n)$, contraction of the stress equation gives
\be
\label{eq:regular-fluid-particle-conservation}
-u_\mu\nabla_\nu T^{\mu\nu}
=\rrestnew'(n)\nabla_\mu(nu^\mu).
\ee
If $\rrestnew'(n)>0$, stress conservation therefore gives particle conservation. This implication uses the constitutive chain rule. The weak formulation imposes stress conservation and the particle-number inequality by testing both against smooth functions.}

{
The scalar $\mathcal N(t)$ in~\eqref{eq:intro-total-particle-number} is the integral of the particle current on a slice. For compact slices without boundary and finite particle flux, testing with a function of time alone gives
\be
\label{eq:total-particle-number-inequality}
\int_0^T\mathcal N(t)\psi'(t)\,dt+
\mathcal N_{\rm in}\psi(0)\le0,
\qquad 0\le\psi\in C_c^\infty([0,T)),
\ee
where $\mathcal N_{\rm in}$ is computed from the prescribed initial particle data. Thus $\mathcal N$ has a nondecreasing representative on $(0,T)$ and $\mathcal N(t)\ge\mathcal N_{\rm in}$ for almost every $t>0$. Conservation replaces this inequality by equality for the exact-minor configuration class. The local entropy inequality remains the admissibility condition, including when there is a spatial boundary and particle flux through it.
\par}
In coordinates with $t=x^0$, the coordinate density and flux in the divergence inequality are $-|g|^{1/2}{\particlecurrent}^0$ and $(-|g|^{1/2}{\particlecurrent}^i)$. We prove below that the density is convex as a function of the conservative variables under {these constitutive and observer conditions}.
For the constitutive laws satisfying the conditions of Section~\ref{sec:closed-vacuum-entropy} below, the convex mathematical entropy extends to the boundary of its conservative domain. At the boundary the current is read from the fields and stress allowed by the constitutive law, without requiring a $C^1$ extension.

{In this section we fix a $C^2$ constitutive law $\raug$} {on the positive-density domain} and its selected variables, with the convexity at rest of Definition~\ref{def:rest-polyconvexity}. At a fixed spacetime point, choose a $t$-adapted coframe with $t=(1,0,0,0)$ and ${\particlecurrent}^0>0$. Write $Y$ for the selected coordinates of $(F,G,\particlecurrent)$ in this coframe. For the fixed law and observer $t$, the map $\Psi_t$ below sends $Y$ to the \emph{conservative variables} $V$ by taking the spatial form components and the time components of the constitutive stress. Proposition~\ref{prop:aug-continuity} gives conservation of the mathematical entropy current for sufficiently regular solutions under its constitutive nondegeneracy assumption. The density is convex under the constitutive and observer assumptions specified below:
\bse
\be
\label{eq:entropy-candidate-conservative}
V\coloneqq\Psi_t(Y)=(F_i^a,G_{ijb},T^{\mu0}),\qquad
\eta_t\bigl(\Psi_t(Y)\bigr)=e_t(Y)\coloneqq-{\particlecurrent}^0.
\ee
The goal is to determine the sign of the Hessian of $\eta_t$ in the conservative variables $V$. Where $\Psi_t$ is locally invertible, $e_t=\eta_t\circ\Psi_t$. We first compute a covector $\lambda$ such that $de_t=\lambda\cdot d\Psi_t$. In conservative coordinates this covector is $D_V\eta_t$. Differentiating the identity once more will give the quadratic form that represents the conservative Hessian in the original field coordinates. We set $P\coloneqq\del_F\raug$, $Q\coloneqq\del_G\raug$, $B\coloneqq\del_{\particlecurrent}\raug$, and $b\coloneqq B\cdot {\particlecurrent}\ne0$, with full antisymmetric sums as in Lemma~\ref{lem:algebraic-identities}. The mathematical entropy coefficients are
\be
\label{eq:lambda-explicit-local06}
\lambda_{F_i^a}\coloneqq({\particlecurrent}^0P_a{}^i-{\particlecurrent}^iP_a{}^0)/b,\qquad
\lambda_{G_{ijb}}\coloneqq2({\particlecurrent}^0Q^{ijb}-{\particlecurrent}^iQ^{0jb}+{\particlecurrent}^jQ^{0ib})/b,\qquad
\lambda_{T^{\mu0}}\coloneqq {\particlecurrent}_\mu/b.
\ee
\ese
{In the rest variables introduced above, $b=n\partial_n\bar\rho$. Indeed, scaling $\particlecurrent$ at fixed $F,G$ preserves its direction and the orthogonality condition~\eqref{eq:450bis-2}, while it scales $n$ and leaves the rest form components fixed. Thus the assumption $b\ne0$ is precisely $\partial_n\bar\rho\ne0$ in matter.}
Differentiating ${\particlecurrent}^\mu T_\mu{}^0=-\raug {\particlecurrent}^0$ at fixed metric, and differentiating the orthogonality condition~\eqref{eq:450bis-2}, gives
\bse
\be
-b\,d{\particlecurrent}^0={\particlecurrent}_\mu dT^{\mu0}
+\sum_{a,i}({\particlecurrent}^0P_a{}^i-{\particlecurrent}^iP_a{}^0)dF_i^a
+2\sum_{b,i<j}({\particlecurrent}^0Q^{ijb}-{\particlecurrent}^iQ^{0jb}+{\particlecurrent}^jQ^{0ib})dG_{ijb}.
\ee
\be
\label{eq:entropy-multiplier-differential}
de_t=\lambda\cdot d\Psi_t.
\ee
\ese
Indeed $F_\mu^a d{\particlecurrent}^\mu=-{\particlecurrent}^\mu dF_\mu^a$ and $G_{\mu\beta b}d{\particlecurrent}^\mu=-{\particlecurrent}^\mu dG_{\mu\beta b}$ cancel the time components of the exterior forms. Empty sums vanish, so the calculation includes the fluid formulation. The identity does not require invertibility of $\Psi_t$.

\paragraph{The Hessian of the mathematical entropy density.}
{The constitutive $C^2$ regularity makes $\lambda$ and $\Psi_t$ continuously differentiable. We define the continuous quadratic form using their first derivatives:}
\be
\label{eq:LP3-local-matrix}
{K_t\coloneqq(D_Y\lambda)^TD_Y\Psi_t.}
\ee
{For a smooth law, differentiating~\eqref{eq:entropy-multiplier-differential} gives
\be
 K_t=D_Y^2e_t-\sum_{\alpha=1}^{d_{\rm sel}}\lambda_\alpha D_Y^2\Psi_t^\alpha,
\ee
which proves symmetry. This calculation also applies to a $C^2$ law by approximation: near the compact rotation orbit of a fixed triple of rest variables, extend the objective rest energy and convolve with a rotation-invariant kernel. The smooth objective energies converge with their first two derivatives. The expressions for $K_t$, the rest quadratic form~\eqref{eq:rest-K-test} and the quadratic current~\eqref{eq:quadratic-observer-current} depend only on these derivatives and the algebraic change of frame, so the identities pass to the limit. The approximating laws need not satisfy convexity or the acoustic inequality; these assumptions are used for the limiting quadratic forms. Thus the formulas below require no third constitutive derivative.}
Where $\Psi_t$ has a local inverse, $K_t[\delta Y,\delta Y]=D_V^2\eta_t[D_Y\Psi_t\,\delta Y,D_Y\Psi_t\,\delta Y]$. Thus the sign of $K_t$ determines local convexity in the conservative variables there. The formula for $K_t$ itself also makes sense before inversion. In the rest frame we express it through the Hessian of the rest energy and the inertia matrix.

\begin{definition}
\label{def:constitutive-enthalpy}
The \underline{\emph{constitutive enthalpy}} is $H(\Frest ,\Grest ,n)\coloneqq\del_n\bar\rho(\Frest ,\Grest ,n)$, with $\Frest ,\Grest $ fixed in this derivative. A law has \underline{\emph{positive enthalpy}} on a prescribed positive-density rest domain if $H(\Frest ,\Grest ,n)>0$ at every point of that domain.
\end{definition}

Positive enthalpy enters the recovery of the particle density and the calculation of the entropy Hessian.
\begin{itemize}
\item[(i)] {At fixed rest minors, $H>0$ means that increasing the particle density increases the rest energy per unit volume. For a perfect fluid, the pressure satisfies $p=n\bar\rho'(n)-\bar\rho(n)$, so $H=\bar\rho'(n)=(\bar\rho+p)/n$ is the enthalpy per particle. A positive rest-mass term and a nondecreasing internal energy density ensure this sign.}
\item[(ii)] {For a fixed material velocity relative to an observer and fixed spatial form components, increasing the observed density increases the proper density while keeping the rest form components fixed. Positive enthalpy therefore makes the constitutive energy in the observer's variables strictly increasing with the observed density. This monotonicity is used to recover that density from the conservative variables in Section~\ref{sec:convex-entropy-relativistic-elastodynamics} below.}
\item[(iii)] {The identity $B\cdot {\particlecurrent}=nH$ gives the denominator of the entropy covector in the equation~\eqref{eq:lambda-explicit-local06}. In the rest frame, the factor $1/H$ also fixes the sign relating the rest-energy Hessian to the corresponding block of the entropy quadratic form.  We impose this sign in matter and examine the entropy limit at cavitation.}
\end{itemize}

To see what the inertia matrix measures, fix an orthonormal observer frame and write the current as $\particlecurrent=r(1,v)$, with observer density $r={\particlecurrent}^0$ and spatial velocity components $v^j$. At the rest state $v=0$ we have $r=n$. Let $P_i=T_i{}^0$ be the spatial momentum density in this frame. The \emph{inertia matrix} $\mathsf M$ is the derivative of $P$ with respect to $v$ at rest, with $r$ and the spatial form components held fixed. Thus a velocity perturbation gives $\delta P_i=\sum_{j=1}^3\mathsf M_{ij}\delta v^j$. In the linearized spatial momentum equations this contributes $\sum_{j=1}^3\mathsf M_{ij}\del_t\delta v^j$. The proof below derives this derivative from the constitutive energy.

\paragraph{{The velocity block at rest.}}
We write $\Frest _{Ai}$ for the $i$th spatial component of selected one-form $A$ and $\Grest _{Bi}$ for the $i$th component of the spatial dual of selected two-form $B$, as in~\eqref{eq:rest-minor-notation}. The material indices range over $1\le A\le r_F$ and $1\le B\le r_G$, while $i,j,k\in\{1,2,3\}$ are spatial indices. The derivatives of $\bar\rho(\Frest ,\Grest ,n)$ are taken with respect to the independent entries of $(\Frest ,\Grest ,n)$. With $\delta_{ij}$ the spatial Kronecker symbol, the entries of the inertia matrix are
\be
\label{eq:rest-inertia-test}
\begin{aligned}
\mathsf M_{ij}
={}& n\frac{\del\bar\rho}{\del n}\,\delta_{ij}
+\frac12\sum_{A=1}^{r_F}
\left(\Frest _{Ai}\frac{\del\bar\rho}{\del \Frest _{Aj}}+
      \Frest _{Aj}\frac{\del\bar\rho}{\del \Frest _{Ai}}\right)\\
&-\frac12\sum_{B=1}^{r_G}
\left(\Grest _{Bi}\frac{\del\bar\rho}{\del \Grest _{Bj}}+
      \Grest _{Bj}\frac{\del\bar\rho}{\del \Grest _{Bi}}\right)
+\delta_{ij}\sum_{B=1}^{r_G}\sum_{k=1}^3 \Grest _{Bk}\frac{\del\bar\rho}{\del \Grest _{Bk}} .
\end{aligned}
\ee
In matrix notation $(\bar\rho)_{\Frest}$ and $(\bar\rho)_{\Grest}$ denote the arrays of partial derivatives with respect to the entries of $\Frest $ and $\Grest $. The last sum in~\eqref{eq:rest-inertia-test} is the double contraction written $\Grest :(\bar\rho)_{\Grest}$ elsewhere. We write $\operatorname{sym}L=(L+L^T)/2$ for the symmetric part of a matrix.
The term $nH I_3$ is the density contribution to the inertia matrix. For a fluid it is the whole matrix, so $H>0$ in matter also gives a positive definite inertia matrix. For an elastic material the form derivatives contribute the other terms in the equation~\eqref{eq:rest-inertia-test}. We therefore require positive definiteness of the full matrix separately to control the velocity block of the entropy quadratic form.

\paragraph{{The velocity block and the entropy Hessian.}}
{In the rest frame the identities for the conservative variables separate the variations of the forms and density from those of the velocity.}
\begin{proposition}[Rest convexity and the velocity block of the mathematical entropy]
\label{prop:rest-entropy-test}
Consider a constitutive scalar built from the specified forms, current, metric and volume tensors, with no additional spacetime background. At {a triple of constitutive fields with positive density satisfying the orthogonality condition~\eqref{eq:450bis-2}}, choose a current rest coframe, with rest variables $(\Frest ,\Grest ,r)$ and $r=n>0$. The rows of $\Frest $ and $\Grest $ are the selected rest components, and $\mathsf M$ is the inertia matrix~\eqref{eq:rest-inertia-test}. {Assume $\bar\rho_r\ne0$, so that the entropy multiplier is defined.}
\begin{itemize}
\item[{(i)}] In the observer coordinates $(\Fsp,\Gsp,r,v)$ with $\particlecurrent=r(1,v)$, the quadratic form $K_{dt}$ for the rest time direction is given by
\be
\label{eq:rest-K-test}
b=r\bar\rho_r,\qquad
K_{dt}\big|_{v=0}
=\frac1{\bar\rho_r}\begin{pmatrix}D^2_{(\Frest,\Grest,r)}\bar\rho&0\\0&\mathsf M\end{pmatrix}.
\ee
\item[{(ii)}] If $D^2\bar\rho\ge0$ and $\mathsf M>0$, {the acoustic time form of the configuration formulation is negative in every nonzero material direction when the triple belongs to the set of exact spacetime minors}. With $\bar\rho_r>0$, the quadratic form $K_{dt}$ is nonnegative and is positive definite exactly when $D^2\bar\rho>0$. Conversely, under the same convexity and inertia matrix conditions, nonnegativity of $K_{dt}$ forces $\bar\rho_r>0$. For this coframe, strict positivity persists locally.
\end{itemize}
\end{proposition}

For a perfect fluid the form variables are absent and $\mathsf M=n\bar\rho'(n)I_3$. With positive enthalpy, the two blocks in~\eqref{eq:rest-K-test} are $\bar\rho''(n)/\bar\rho'(n)$ and $nI_3$.
We relate the sign of $K_t$ to the wave speeds for a general observer below.

\paragraph{{The Hessian calculation.}}
{We express the stress in the spatial fields and velocity, then differentiate the mathematical entropy density.}
\begin{proof}
We set $B(v)=(I-vv^T)^{1/2}$. The orthogonality condition~\eqref{eq:450bis-2} and a Lorentz boost to the rest frame give the reduced scalar
\bse
\be
\widetilde\rho(\Fsp ,\Gsp ,r,v)=\bar\rho\bigl(\Fsp B(v),\Gsp \operatorname{cof}B(v),r\sqrt{1-|v|^2}\bigr).
\ee
This formula is independent of the rest coframe's spatial rotation and even in $v$. The constitutive stress and the differential of the orthogonality condition~\eqref{eq:450bis-2} give
\be
P_i\coloneqq T_i{}^0=-\del_{v_i}\widetilde\rho,
\qquad E\coloneqq T^{00}=\widetilde\rho-v\cdot D_v\widetilde\rho.
\ee
{Indeed, $\del_{v_i}F_0^a=-F_i^a$, $\del_{v_i}G_{0j b}=-G_{ij b}$ and $\del_{v_i}{\particlecurrent}^j=r\delta_i^j$ give
\be
 \del_{v_i}\widetilde\rho
 =-F_i^aP_a{}^0-2Q^{0jb}G_{ij b}+rB_i.
\ee
The factor two uses the full antisymmetric sum. Substitution in~\eqref{augment3} gives $T_i{}^0=-\del_{v_i}\widetilde\rho$ and $T_0{}^0=-\widetilde\rho+v\cdot D_v\widetilde\rho$, proving both identities.} The expansions
\be
B(v)=I-\tfrac12vv^T+\Ocal(|v|^4),\qquad
\operatorname{cof}B(v)=I+\tfrac12(vv^T-|v|^2I)+\Ocal(|v|^4)
\ee
yield $-D_v^2\widetilde\rho|_0=\mathsf M$ with the entries in~\eqref{eq:rest-inertia-test}. Since $P=-D_v\widetilde\rho$, this equals $D_vP|_0$. The Hessian of the energy is
\be
D^2E\big|_0=\operatorname{diag}(D^2\bar\rho,\mathsf M).
\ee
\ese
The other conservative coordinates $\Fsp ,\Gsp $ and $e=-r$ are affine. At rest the stress multiplier is $\lambda_E=-r/b$ and $\lambda_{P_i}=0$. Formula~\eqref{eq:rest-K-test} follows from the definition of $K$ and $b=r\bar\rho_r$.

{For the acoustic conclusion, we take a triple on the positive-density part of the set of exact spacetime minors and set $\Frest_{\rm full}=(F_i^A)_{A,i=1}^3$. It is invertible because $\particlecurrent=\particlecurrent(F)$ is nonzero and timelike.} A time rank-one variation of the physical gradient leaves its spatial minors fixed and has $\delta v=-\Frest_{\rm full}^{-1}U$. A spatial rank-one variation induces a linear variation $\mathcal L_kU$ of the list $(\Fsp ,\Gsp ,r)$; its higher minor variations vanish. Evenness in $v$ eliminates the mixed term, so
\be
\label{eq:rest-rank-one-increment}
q_U(\tau,k)
=(\mathcal L_kU)^TD^2\bar\rho(\mathcal L_kU)
-\tau^2(\Frest_{\rm full}^{-1}U)^T\mathsf M(\Frest_{\rm full}^{-1}U).
\ee
Rest convexity makes the spatial part nonnegative. Invertibility of $\Frest_{\rm full}$ makes negativity in the time direction of the rest frame for every $U\ne0$ equivalent to $\mathsf M>0$. The remaining assertions follow from~\eqref{eq:rest-K-test}; its velocity block determines the positive sign of $\bar\rho_r$. Continuity preserves strict positivity in {the chosen chart}. {If $\bar\rho_r=0$, varying $r$ at fixed rest forms and velocity gives $de_{dt}=-dr\ne0$ but $d\Psi_{dt}=0$. Thus neither a finite differential multiplier nor a finite supporting covector exists for an entropy agreeing with $e_{dt}$ along these field values; a nonsmooth entropy without such a support remains possible.}
This completes the proof of Proposition~\ref{prop:rest-entropy-test}.
\end{proof}

The notation for the different expressions of the energy is summarized in Table~\ref{tab:energy-notation}. The mathematical entropy $\eta_t$ is a separate function of the conservative variables.

\begin{table}[htbp]
\centering
\renewcommand{\arraystretch}{1.1}
\begin{tabular}{@{}>{\raggedright\arraybackslash}p{.34\linewidth}>{\raggedright\arraybackslash}p{.59\linewidth}@{}}
\hline
Symbol & Meaning \\
\hline
$\rphys(F;g)$ & Energy of the configuration formulation \\
$\raug(F,G,\particlecurrent;g)$ & Constitutive energy in independent spacetime minors \\
$\bar\rho(\Frest ,\Grest ,n)$ & Constitutive energy in the material rest frame \\
$\widetilde\rho(\Fsp,\Gsp,r,v)$ & Expression of the energy in observer variables \\
$E=T^{00}$ & Energy measured in the chosen observer frame \\
$\eta_t(V)$ & Mathematical entropy in conservative variables \\
$\particlecurrent$ &  Particle-number current \\
$\mathcal N(t)$ & Total particle number on $\Sigma_t$ \\
$n,\ r=-g(\particlecurrent,\slicenormal)$ & Proper particle density and density measured by the slice observer \\
$\slicenormal$ & Future unit normal to the slice \\
\hline
\end{tabular}
\caption{Energy and mathematical entropy notation.}
\label{tab:energy-notation}
\end{table}

\paragraph{Local recovery of fields at positive density.}
The quadratic form computed in the rest frame can also be used to control the inverse change of variables.
\begin{proposition}[Local recovery of fields and bounds on the entropy Hessian]
\label{prop:LP3-inverse-local06}
Work on the positive-density constitutive domain where $\particlecurrent$ is timelike and $n>0$. All matrices in this proposition act on the $d_{\rm sel}$ selected variables. If $K_t$ is positive definite, $A_t=D_Y\Psi_t$ is invertible. If on a state set $K_t\ge kI$ and $\|D_Y\lambda_t\|\le L$, then $\|A_t^{-1}\|\le L/k$. Thus positive definiteness also gives a local inverse for the map to conservative variables. If also $\|A_t\|\le M$ and $\|K_t\|\le K_1$, the Hessian of the mathematical entropy expressed in conservative variables satisfies
\be
\frac{k}{M^2}I\le D_V^2(e_t\circ\Psi_t^{-1})
\le K_1\frac{L^2}{k^2}I.
\ee
\end{proposition}
\begin{proof}
If $A_tz=0$, then $z^TK_tz=(D_Y\lambda_t z)\cdot(A_tz)=0$, so positivity forces $z=0$. The $d_{\rm sel}\times {d_{\rm sel}}$ matrix is therefore invertible, and the {$C^1$} inverse function theorem applies. The same identity gives
\bse
\be
k|z|^2\le z^TK_tz\le L|z|\,|A_tz|,
\ee
proving the inverse bound. We write $\eta=e_t\circ\Psi_t^{-1}$. The differential identity identifies $D_V\eta=\lambda_t\circ\Psi_t^{-1}$. Differentiating once more gives
\be
D_V^2\eta=A_t^{-T}K_tA_t^{-1}>0.
\ee
\ese
The lower estimate uses $|A_t^{-1}y|\ge |y|/M$; the upper uses the inverse bound. Both hold {on this set} with its positive coercivity constants.
\end{proof}

{By~\eqref{eq:LP3-local-matrix}, changes of field coordinates transform $K_t$ by congruence.} Some elastic laws have a full augmented Hessian of signature $(15-,19+)$. Its rank drops when the energy is independent of $F,G$. Section~\ref{sec:convex-entropy-relativistic-elastodynamics} below treats convexity and extension to degenerate limits. 
For a perfect fluid the form variables are absent. In its rest frame the density coefficient in the entropy quadratic form is $\bar\rho^{\prime\prime}(n)/\bar\rho^\prime(n)$, so convexity in matter requires $\bar\rho^{\prime\prime}\ge0$ when $\bar\rho^\prime>0$.

\paragraph{Particle production and covariant relaxation.}
{A fluid relaxation model motivates the sign of the mathematical entropy inequality.\footnote{See~\cite{TorusT2-2-2026} for the particle-production convention.} A covariant potential $\chi(\beta,\zeta)$ determines the stress $T_\chi$, a relaxation current $A_\chi$ and a particle current $N_\chi$. Here $\beta$ is the rescaled fluid momentum and $\zeta$ is a symmetric tensor which vanishes at equilibrium. The stress is conserved, and a relaxation source $-\mathscr B$ in the equation for $A_\chi$ gives particle production:
\be
\label{eq:torus-particle-production-identity}
 \nabla_\alpha A_\chi^{\alpha\beta\gamma}=-\mathscr B^{\beta\gamma},
 \qquad \nabla_\alpha N_\chi^\alpha
      =\zeta_{\beta\gamma}\mathscr B^{\beta\gamma}\ge0 .
\ee
The principal matrices are symmetric third derivatives of $\chi$. Their contraction with a future time covector must be positive definite for symmetric hyperbolicity in that direction. Positivity for every future timelike covector bounds propagation by the speed of light. If the stress and particle current converge in the sense of distributions to their equilibrium values, stress conservation and the particle-number inequality pass to the limit. Identification of the limiting constitutive stress at cavitation requires a separate constitutive argument.}

\subsection{Admissible observers and convexity of the mathematical entropy density}
\label{subsec:observer-entropy}

\paragraph{Propagation directions and observers.}
{We work at a point with positive particle density and use the material rest frame. A unit spatial covector $e$ specifies the propagation direction of a wave, whereas a future timelike covector $t=(\tau,k)$, with $\tau>|k|$, specifies the time direction used by an observer. The quadratic form $K_t$ was defined in~\eqref{eq:LP3-local-matrix}; where the conservative variables can be inverted, its sign determines the convexity of the mathematical entropy density. We now relate this form to the wave speeds, first for one observer and then for all observers.}

\begin{proposition}[Wave speeds and admissible time directions]
\label{prop:observer-acoustic-coercivity}
{At a positive-density triple of rest variables, write $z=(\Frest ,\Grest ,r)$ with $r=n$. The symmetric matrix $A=D_z^2\bar\rho$ is the Hessian of the rest energy with respect to these variables, $H=\bar\rho_r$ is the enthalpy, and $\mathsf M$ is the inertia matrix defined in~\eqref{eq:rest-inertia-test}. Assume $A\ge0$, $H>0$ and $\mathsf M>0$. For a spatial covector $k$, let $\mathcal D_k$ map a velocity amplitude $b$ to the associated variations of the one-forms, dual two-forms and density:}
\be
\label{eq:retained-acoustic-map}
\mathcal D_k b=\bigl((\Frest b)\otimes k,\ 
 (k\cdot b)\Grest-(\Grest k)\otimes b,\ 
 r\,k\cdot b\bigr).
\ee
\begin{itemize}
\item[(i)] {\emph{Acoustic speeds.}} For a unit spatial covector $e$, {write $c_\ell(e)^2$ for the eigenvalues of the symmetric nonnegative matrix}
\be
\label{eq:physical-subluminal-speeds}
\mathsf M^{-1/2}\mathcal D_e^TA\mathcal D_e\mathsf M^{-1/2}.
\ee
{On the positive-density part of the set of exact spacetime minors, these eigenvalues are the squared acoustic speeds of the configuration formulation. The corresponding signed speeds are $\pm c_\ell(e)$ with $c_\ell(e)\ge0$,  including zero and repeated values.}
\item[(ii)] {\emph{A chosen observer.}} For a fixed observer $t=(\tau,k)$ with $\tau>|k|$, the quadratic form $K_t$ is nonnegative exactly when $\mathcal D_k^TA\mathcal D_k\le\tau^2\mathsf M$. {In this case, put $\mathsf C_t\coloneqq\tau^2\mathsf M-\mathcal D_k^TA\mathcal D_k$. The kernel of $K_t$ consists of the pairs
\be
 \delta v\in\ker \mathsf C_t,\qquad
 \delta z=-\tau^{-1}\mathcal D_k\delta v+w,\quad w\in\ker A,
\ee
and its nullity is $\dim\ker A+\dim\ker \mathsf C_t$. On the set of exact spacetime minors, the strict matrix inequality places $t$ inside the acoustic cone of time covectors.} Under that strict inequality, $K_t$ is positive definite exactly when $A$ is positive definite.
\item[(iii)] {\emph{All timelike observers.}} {The strict matrix inequality in (ii) holds for every future timelike $t$ if and only if}
\be
\label{eq:observer-acoustic-cone}
\mathcal D_e^TA\mathcal D_e\le\mathsf M
\qquad\text{for every unit spatial covector }e .
\ee
{Equivalently, every eigenvalue in (i) is at most one. On the set of exact spacetime minors, every signed acoustic speed then lies in $[-1,1]$, in units where the light speed is one.}
\end{itemize}
\end{proposition}

{For a perfect fluid the matrix in (i) has the compression eigenvalue $c_s^2=n\bar\rho''(n)/\bar\rho'(n)$ and two zero eigenvalues. Thus the condition in (iii) reduces to $0\le c_s^2\le1$ under {these convexity and positivity assumptions}.}

\phantomsection\label{app:observer-acoustic-proof}

\begin{proof}
We fix the spacetime and material points and a current rest coframe, using algebraic coordinates on the selected manifold. For a fixed base point $z_*=(\Frest_*,\Grest_*,r_*)$, $u_*=(1,0,0,0)$, we put $H_*=\bar\rho_r(z_*)$. {We write $A_*=A(z_*)$ and $\mathsf M_*=\mathsf M(z_*)$.}

{\bf (i).} \emph{Acoustic speeds.} On the set of exact spacetime minors at positive density, write $z_0=\Frest_{\rm full}^{-1}U$. The spatial gradient along the rank-one variation is $\Frest_{\rm full}(I+\lambda z_0\otimes k)$. Differentiating its affine minors gives $\mathcal L_k U=\mathcal D_k z_0$, including $\dot r=r_*\,k\cdot z_0$. Substitution in~\eqref{eq:rest-rank-one-increment} gives
\be
q_U(\tau,k)=z_0^T\bigl(\mathcal D_k^TA_*\mathcal D_k-\tau^2\mathsf M_*\bigr)z_0.
\ee

The characteristic equation in the rest frame is $\det(\mathcal D_k^TA\mathcal D_k-\omega^2\mathsf M)=0$. Congruence by $\mathsf M^{-1/2}$ proves~\eqref{eq:physical-subluminal-speeds}.

\paragraph{{Acoustic speeds for a prescribed observer.}}
\vskip.2cm \noindent {\bf (ii).} \emph{A prescribed observer.} With full antisymmetric derivatives, we use $P_*=\del_F\raug$ and $Q_*=\del_G\raug$ and define
\bse
\be
\label{eq:observer-current-coefficients}
w_{*\mu}=u_{*\mu}/H_*,\qquad C_{*a}^{\nu\alpha}=(u_*^\nu P_{*a}^\alpha-u_*^\alpha P_{*a}^\nu)/H_*,\qquad D_{*b}^{\nu\alpha\beta}=(u_*^\nu Q_{*b}^{\alpha\beta}-u_*^\alpha Q_{*b}^{\nu\beta}+u_*^\beta Q_{*b}^{\nu\alpha})/H_*.
\ee
We write $\Delta$ for the change of a flux from its value at this base point. With coefficients fixed, we differentiate the following affine combination of fluxes:
\be
\mathscr B^\nu=-\Delta {\particlecurrent}^\nu-w_{*\mu}\Delta T^{\mu\nu}
-C_{*a}^{\nu\alpha}\Delta F_\alpha^a
-D_{*b}^{\nu\alpha\beta}\Delta G_{\alpha\beta b}.
\ee
\ese
In any time-adapted coframe its time component is $e_t-e_{t*}-\lambda_{t*}\cdot(\Psi_t-\Psi_{t*})$ by the multiplier identity~\eqref{eq:entropy-multiplier-differential}; antisymmetry leaves only independent spatial form components. Its second variation at the base point gives a spacetime vector $j(\delta Y)$ quadratic in the perturbation. The first variation vanishes, so in any adapted coframe the time component of this quadratic current is $K_t$ in the corresponding field coordinates.

{We now calculate this quadratic current in the material rest frame.} We write $\particlecurrent=r(1,v)$ and $z=(\Fsp,\Gsp,r)$ in the fixed coframe at rest with respect to $u_*$. At the reference, $v_*=0$, so $\Fsp_*=\Frest_*$ and $\Gsp_*=\Grest_*$. The reduced energy of nearby field values is obtained by the boost to its own rest frame.
\bse
\be
B(v)=(I-vv^T)^{1/2},\qquad
\widetilde\rho(z,v)=\bar\rho(\Fsp B(v),\Gsp \operatorname{cof}B(v),r\det B(v)).
\ee
The stress identities used in Proposition~\ref{prop:rest-entropy-test} give $E=\widetilde\rho-v\cdot D_v\widetilde\rho$ and $\mathsf P=-D_v\widetilde\rho$. At the fixed base point the spatial momenta are $(\bar\rho_{\Frest})_*,(\bar\rho_{\Grest})_*$ and the two-form derivative has the factor $1/2$ corresponding to $G_{ij}=\varepsilon_{ijk}\Gsp _k$. In the following expansion, $X:Y$ is the sum of the products of corresponding entries of two spatial arrays, summed over the selected material indices and the spatial components. Thus $(\bar\rho_{\Frest})_*:(\Fsp -\Fsp _*)$ pairs the base point energy derivative with the change in the one-form components. Expanding the current gives
\be
\aligned
&H_*\mathscr B^0=E-\bar\rho_*-(\bar\rho_{\Frest})_*:(\Fsp -\Fsp _*)
-(\bar\rho_{\Grest})_*:(\Gsp -\Gsp _*)-H_*(r-r_*),\\
&H_*\mathscr B_{\rm sp}=\mathsf P-H_*rv-(\bar\rho_{\Frest})_*^T\Fsp\,v
+\Gsp ^T(\bar\rho_{\Grest})_*v-\bigl(\Gsp :(\bar\rho_{\Grest})_*\bigr)v.
\endaligned
\ee
\ese
For the last two terms we use $2\sum_j Q_*^{ijb}G_{\ell jb} =\delta_{i\ell}\,\Gsp :(\bar\rho_{\Grest})_* -(\Gsp ^T(\bar\rho_{\Grest})_*)_{i\ell}$.

Rotation invariance of the rest energy makes $\Frest^T\bar\rho_{\Frest}+\Grest^T\bar\rho_{\Grest}$ symmetric. At the reference, the first variations of the rest and spatial matrices agree because $DB(0)=D(\operatorname{cof}B)(0)=0$. The linear spatial and time coefficients vanish, so $D\mathscr B|_*=0$ on the whole selected algebraic tangent space. Evenness of $\widetilde\rho$ in $v$ and~\eqref{eq:rest-K-test} give its time Hessian; spatial differentiation gives the complete identities, with $\delta Y=(\delta z,\delta v)$,
\bse
\be
\label{eq:quadratic-observer-current}
H_*j^0(\delta Y)=\delta z^TA_*\delta z+\delta v^T\mathsf M_*\delta v,\quad
H_*j_{\rm sp}(\delta Y)\cdot k=2\,\delta z^TA_*\mathcal D_k\delta v,
\quad j^\nu(\delta Y)\coloneqq D^2\mathscr B^\nu|_*[\delta Y,\delta Y].
\ee
Differentiating the rotation identity gives
\be\begin{aligned}
&\delta\mathsf M-H_*\delta r\,I-(\bar\rho_{\Frest})_*^T\delta \Fsp +\delta\Gsp ^T(\bar\rho_{\Grest})_*-\bigl(\delta\Gsp :(\bar\rho_{\Grest})_*\bigr)I\\
&\qquad=(\delta \bar\rho_{\Frest})^T\Frest_* -\Gsp _*^T\delta \bar\rho_{\Grest}+\bigl(\Gsp _*:\delta \bar\rho_{\Grest}+r_*\delta \bar\rho_r\bigr)I.
\end{aligned}\ee
\ese
Multiplication by $k$ and $\delta v$ gives the second identity of~\eqref{eq:quadratic-observer-current}, including mixed derivatives between both form families and the density.

\paragraph{{The condition for a chosen observer.}}
{We next use the quadratic current to identify the condition for a chosen observer.} For a fixed time covector $t=(\tau,k)$ with $\tau>|k|$, the contraction $\sum_\nu t_\nu j^\nu$ gives the scalar quadratic form $K_t$ in the field variables
\bse
\be
H_*K_t[\delta Y,\delta Y]=\tau\bigl(\delta z+\tau^{-1}\mathcal D_k\delta v\bigr)^TA_*\bigl(\delta z+\tau^{-1}\mathcal D_k\delta v\bigr)+\tau^{-1}\delta v^T\bigl(\tau^2\mathsf M_*-\mathcal D_k^TA_*\mathcal D_k\bigr)\delta v.
\ee
Thus nonnegativity is equivalent to the chosen observer's matrix inequality: sufficiency uses $A_*\ge0$, and necessity sets $\delta z=-\tau^{-1}\mathcal D_k\delta v$. This algebraic coordinate change neither inverts the map to conservative variables nor introduces second-derivative terms, since $D\mathscr B|_*=0$.

For $|k|=1$ the decisive identity is
\be
\label{eq:observer-current-square}
H_*\bigl(j^0-j_{\rm sp}\cdot k\bigr)=(\delta z-\mathcal D_k\delta v)^TA_*(\delta z-\mathcal D_k\delta v)+\delta v^T(\mathsf M_*-\mathcal D_k^TA_*\mathcal D_k)\delta v.
\ee
\ese
{For $k\ne0$, we put $e=k/|k|$ and write $c_{\max}(e)$ for the largest of the speeds in~\eqref{eq:physical-subluminal-speeds}. Cauchy--Schwarz for $A_*\ge0$ and the generalized Rayleigh quotient give
\be
 |j_{\rm sp}\cdot e|\le c_{\max}(e)j^0,\qquad
 (\tau-|k|c_{\max}(e))j^0\le K_t\le
 (\tau+|k|c_{\max}(e))j^0.
\ee
For $c_{\max}(e)>0$, both bounds are attained by taking $\delta v$ to be a corresponding generalized eigenvector and $\delta z=\pm c_{\max}(e)^{-1}\mathcal D_e\delta v$. If $c_{\max}(e)=0$, the spatial component in direction $e$ vanishes. Thus these constants are sharp at the chosen field values. When $k=0$, $K_t=\tau j^0$.}

{When $\mathsf C_t\ge0$, the two terms in the completed square vanish separately; this gives the kernel and nullity in (ii).}

\vskip.2cm \noindent {\bf (iii).} \emph{All timelike observers.} Under~\eqref{eq:observer-acoustic-cone}, $c_{\max}(e)\le1$, so $j$ is future causal or zero. Also $\mathsf C_t\ge(\tau^2-|k|^2)\mathsf M>0$ for every timelike conormal, including when an acoustic speed equals one. Conversely, the matrix inequality for every timelike conormal gives~\eqref{eq:observer-acoustic-cone} by letting $\tau\downarrow|k|$. The physical characteristic covectors are then spacelike or null for $g$, an observer-independent property (cf.~\cite[Section~8]{Beig2023} for the acoustic and spacetime cones).

For a fluid, $A=\bar\rho''(n)$, $\mathsf M=n\bar\rho'(n)I$ and $\mathcal D_kb=n\,k\cdot b$. The cone bound is $n\bar\rho''(n)\le\bar\rho'(n)$, or $c_s^2\le1$ with positive enthalpy. The entropy quadratic form is positive definite under the strict conormal inequality and $\bar\rho''(n)>0$.
This completes the proof of Proposition~\ref{prop:observer-acoustic-coercivity}.
\end{proof}

\paragraph{Strict bounds on the acoustic speeds.}
{A strict inequality in~\eqref{eq:observer-acoustic-cone} for every unit direction gives a uniform gap on a compact positive-density rest range $K$, by compactness of $K\times\mathbb S^2$. If the elastic part already satisfies the acoustic bound, adding $m\,n$ leaves its Hessian unchanged and adds $m\,n I$ to the inertia matrix. The gap for the squared speeds is then at least $\inf_K m\,n/\lambda_{\max}(\mathsf M)>0$.
More generally let $\bar\rho=mc^2 n+\bar\rho_{\rm int}$, with smooth, objective and convex internal energy. On $K$, where $n\ge n_0>0$, the acoustic matrix for a unit direction is bounded independently of $c$, whereas the inertia matrix is $mc^2n I+\mathsf M_0$ with $\mathsf M_0$ bounded. For sufficiently large $c$, the enthalpy is positive and the inertia matrix strictly dominates the acoustic matrix on a neighborhood of $K$. This applies in particular to a smooth objective convex extension obtained from Ball's energy by Proposition~\ref{prop:Ball-perspective-new} below.}

\paragraph{One observer, a cone, or all observers.}
The rest direction satisfies the acoustic time inequality when the inertia matrix is positive definite. Positive enthalpy gives the sign of the entropy quadratic form in~\eqref{eq:rest-K-test}. For a general conormal, the principal matrices show how spatial compatibility yields a symmetric system.

Freeze the coefficients at fields at positive density and work in its rest coframe. For a spatial covector $k$, the spatial stress contracted with $k$ is $\mathcal D_k^T D_z\bar\rho-\bar\rho k$. Spatial closure makes each variation of a row of $\Frest$ parallel to $k$ and gives $\delta\Grest\,k=0$. Under these conditions, the term $(\delta\mathcal D_k^T)D_z\bar\rho$ is $(D_z\bar\rho\cdot\delta z)k$, canceling the derivative of $-\bar\rho k$. The momentum equation then has principal spatial term $\mathcal D_k^T A\delta z$. The closure evolution and particle conservation give $\mathcal D_k\delta v$ in the equations for $\delta z$. Multiplying these equations by $A$ gives the symmetric principal matrix at $t=(\tau,k)$:
\be
\label{eq:observer-symmetric-principal}
 \begin{pmatrix}
 \tau A & A\mathcal D_k\\
 \mathcal D_k^TA & \tau\mathsf M
 \end{pmatrix}.
\ee
Its quadratic form in the rest coordinates is $H K_t$. For $A>0$ and the strict conormal inequality in Proposition~\ref{prop:observer-acoustic-coercivity}(ii), the time matrix is positive definite. This calculation concerns the evolution with the spatial compatibility equations imposed, and applies to a perfect fluid by omitting the form blocks.

Let $\mathcal T_+^*$ denote the future unit timelike covectors, normalized by $g^{-1}(t,t)=-1$. A prescribed observer cone $\mathcal C_{\rm obs}$ has a nonempty compact unit section $\mathcal S_{\rm obs}\coloneqq\mathcal C_{\rm obs}\cap\mathcal T_+^*$. We distinguish three choices for the set on which to impose the observer conditions.
\begin{definition}
\label{def:polyconvex-observers}
\phantomsection\label{rem:polyconvexity-conventions}
\label{rem:polycoercivity-conventions}
{For a fixed relativistically polyconvex constitutive law $\raug$,} an observer $t\in\mathcal T_+^*$ is \underline{\emph{admissible}} on the prescribed positive-density domain if $K_t\ge0$ there and if, for every triple of exact spacetime minors and for every nonzero material amplitude $U$, the tensor of $q_U$ has one negative eigenvalue, its other eigenvalues are nonnegative, and $q_U(t)<0$. {The law is to be}
\begin{itemize}
\item[(i)] \emph{observer-invariant relativistically polyconvex} if every $t\in\mathcal T_+^*$ is admissible;
\item[(ii)] \emph{relativistically polyconvex on a prescribed observer cone} if every $t\in\mathcal S_{\rm obs}$ is admissible;
\item[(iii)] \emph{relativistically polyconvex relative to a prescribed observer} if one prescribed timelike direction is admissible.
\end{itemize}
\end{definition}

{For a perfect fluid, the admissible-observer conditions are evaluated on the current variables, omitting the form directions that do not occur in its law.}
Admissibility for all observers gives admissibility on every prescribed observer cone. Admissibility on a cone gives admissibility for each prescribed observer in its unit section. Positive rescaling preserves the signs of the quadratic currents. Under the positive enthalpy and positive definite inertia matrix assumptions of Proposition~\ref{prop:observer-acoustic-coercivity}, the light-cone bound~\eqref{eq:observer-acoustic-cone} is precisely the additional condition for admissibility of every observer. The observer condition therefore adds a restriction to rest-frame convexity.

Strict positivity of $K_t$ gives a uniform lower bound on compact positive-density state ranges and compact normalized observer families. A semidefinite $K_t$ can occur where recovery from the conservative variables fails.

%-----------------------------------------------------------------------------------------------------------------------------------------

\subsection{Comparison with Ball's polyconvexity}
\label{sec:Ball-comparison}
\paragraph{{Inverse deformation and polyconvexity.}}
Let $y$ be a locally invertible static deformation from material coordinates $X$ to spatial coordinates $x$, and set $f=y^{-1}$. In three spatial dimensions with reference particle density one, the two descriptions use the following minors:
\bse\label{eq:Ball-notation}
\be
 \nabla_X y,\qquad \operatorname{cof}(\nabla_X y),\qquad d\coloneqq\det(\nabla_X y)>0,
\ee
\be
 \Frest \coloneqq\nabla_x f=(\nabla_X y)^{-1},\qquad
 \Grest \coloneqq\operatorname{cof}\Frest =\frac{(\nabla_X y)^T}{d},\qquad
 n\coloneqq\det \Frest =\frac1d.
\ee
\ese
Write $A=\nabla_X y$ and $B=\operatorname{cof}A$ for the direct minors on a deformation. Ball's convex extension treats $(A,B,d)$ as independent. Since $dx=d\,dX$, multiplying the rest energy by $d$ gives the energy per unit reference volume. Ball's radial equilibria already include cavities under tension~\cite[Section~2.3]{BallOpen2002}. In the evolution studied here the particle flow can move and deform an initially empty cavity.

\paragraph{Conversion between current and reference volumes.}
Wagner~\cite{Wagner} gives the determinant-weighted change of variables in classical elasticity, {and Pedregal~\cite{Pedregal2008,Pedregal2024} studies the preservation of polyconvexity under the corresponding change of energy.} We express this correspondence in the relativistic rest variables. Ball's condition and positive definiteness of its Hessian are distinguished in~\cite{Ball1,BallOpen2002}. Comparisons with stress--strain monotonicity appear in~\cite{WollnerHolzapfelNeff2026}.

\begin{proposition}[Equivalence with Ball's polyconvexity at positive density]
\label{prop:Ball-perspective-new}
In three dimensions with reference particle density one, let ${\rho_B}(A,B,d)$ be a convex extension of a classical energy ${\widehat\rho_B}$ per reference volume:
\be
{\widehat\rho_B}(\nabla_X y)={\rho_B}\bigl(\nabla_X y,\operatorname{cof}(\nabla_X y),\det(\nabla_X y)\bigr).
\ee
On independent inverse-minor variables $(\Frest ,\Grest ,n)$, $n>0$, define the energy per current volume by
\bse
\be
\label{eq:Ball-perspective-new}
\bar\rho(\Frest ,\Grest ,n)=n{\rho_B}(\Grest ^T/n,\Frest ^T/n,1/n).
\ee
Then $\bar\rho$ is convex on its natural perspective domain. Conversely, a convex $\bar\rho$ gives the convex extension in the deformation variables
\be
{\rho_B}(A,B,d)=d\bar\rho(B^T/d,A^T/d,1/d).
\ee
\ese
These operations are inverse. They apply also to proper convex functions allowed to take the value $+\infty$ and preserve lower semicontinuity relative to $d>0$ and $n>0$. For $C^2$ extensions, the Hessians are congruent up to a positive factor and have the same rank at corresponding points. Thus the change of volume identifies the two convexity conditions at positive density. Approaching $n=0$ corresponds to $d\to+\infty$ in the direct variables.
\end{proposition}
{For reference density $n_0>0$, the relation $nd=n_0$ gives ${\rho_B}(A,B,d)=d\bar\rho(B^T/d,A^T/d,n_0/d)$. Conversely,~\eqref{eq:Ball-perspective-new} uses $n/n_0$ as its perspective parameter.}

\begin{proof}
By~\eqref{eq:Ball-notation}, the exact direct minors are $(\Grest^T/n,\Frest^T/n,1/n)$ and $dX=n\,dx$. On independent variables, the map $(z,t)\mapsto t{\rho_B}(z/t)$ is convex for $t>0$, by applying convexity with weights proportional to $\theta t_1$ and $(1-\theta)t_2$. Setting $(z,t)=((\Grest^T,\Frest^T,1),n)$ proves convexity, including convexity of the domain. The reverse formula uses the same argument. Lower semicontinuity follows from the continuous invertible change of variables and the positive continuous factors. For $s=(\Grest^T/n,\Frest^T/n,1/n)$ along an affine variation, $s''=-2(n'/n)s'$. The first-derivative terms cancel when we differentiate $n{\rho_B}(s)$ twice:
\be
D^2\bar\rho[\delta \Frest ,\delta \Grest ,\delta n]^2
=n\,D^2{\rho_B}[s',s'].
\ee
The map $(\Frest ,\Grest ,n)\mapsto s$ has invertible differential, including its last component $\delta(1/n)=-\delta n/n^2$. This proves equivalence of positive definiteness, as well as semidefinite convexity in the smooth case.
\end{proof}

At zero density the lower semicontinuous extension of the perspective is given by the recession function of the classical energy. Suppose ${\rho_B}$ is proper, convex and lower semicontinuous on $d>0$, and choose $W_0=(A_0,B_0,d_0)$ of finite energy. For $Z=(\Grest^T,\Frest^T,1)$, the value of the extended $\bar\rho$ at $(\Frest,\Grest,0)$ is
\be
\label{eq:Ball-recession-boundary}
 \lim_{\tau\to+\infty}
 \frac{{\rho_B}(W_0+\tau Z)-{\rho_B}(W_0)}{\tau}.
\ee
The quotients are nondecreasing by convexity. To prove the lower bound, let $n_j\downarrow0$ and $Z_j\to Z$, with last component one, and use convexity at $(1-\tau n_j)W_0+\tau Z_j$ for fixed $\tau>0$. Lower semicontinuity then bounds $\liminf_j n_j{\rho_B}(Z_j/n_j)$ below by the quotient at $\tau$. For the reverse bound, let $n=(d_0+\tau)^{-1}$ and $Z_n=n(W_0+\tau Z)$. Then $Z_n\to Z$, its last component is one and the transformed energies tend to the displayed limit. Thus the limit is independent of $W_0$, with $+\infty$ allowed.

We now examine recovery along exact minors. Let $F_0$ have rank two, and assume that the value of the extension at $(F_0,\operatorname{cof}F_0,0)$ is finite. Choose a rank-one matrix $R$ with $D\det(F_0)[R]>0$. The minors of $F_0+tR$ are affine in $t$, with positive determinant for $t>0$. If the energy is finite at one point of this ray, convexity bounds the upper limit along the ray by the value at $t=0$. Lower semicontinuity gives the reverse bound. At any rank, continuity at the boundary values of the fields also gives recovery when the law is defined on its positive-density neighborhood, by filling the null block of a singular value decomposition with small entries of the required orientation. Section~\ref{sec:convex-entropy-relativistic-elastodynamics} below treats the extension of the mathematical entropy in the conservative flux variables.

For a perfect fluid with rest energy per current volume $\bar\rho(n)$, pressure $p=n\bar\rho'(n)-\bar\rho(n)$, and positive reference density $n_0$, the energy per reference volume is $h(d)=d\bar\rho(n_0/d)$. Since $h'(d)=-p(n_0/d)$, this is convex in $d$ whenever $p$ is nondecreasing. The energy has no dependence on the shear variables.

%============================================================================

\section{\texorpdfstring{A convex mathematical entropy for relativistic elastodynamics}{A convex mathematical entropy for relativistic elastodynamics}}
\label{sec:convex-entropy-relativistic-elastodynamics}

%-----------------------------------------------------------------------------------------------------------------------------------------

\subsection{\texorpdfstring{{Future convexity on spacelike slices}}{Future convexity on spacelike slices}}

\paragraph{{The conservative variables.}}
{
We now recover the observed particle density from the conservative variables and construct affine supports for its negative. We then characterize its extension to the constitutive boundary.\par}

Let $t$ be a smooth time function with spacelike slices $\Sigma_t$, and let $\slicenormal$ be their future unit normal. In Section~\ref{section=4-polyconvexity} the material velocity defines the frame in which we measure the rest variables. Here the normal to the slice defines the frame for the densities in the evolution equations. At a fixed point choose an oriented orthonormal frame $(e_0,e_1,e_2,e_3)$ with $e_0=\slicenormal$, and denote its spatial dual covectors by $e^i$. The \emph{energy density} $E$ and \emph{momentum density} $P=(P^i)_{i=1}^3$ are
\be
\label{eq:slice-energy-momentum-density}
 E\coloneqq T({\slicenormal}^\flat,{\slicenormal}^\flat)=T^{00},\qquad
 P^i\coloneqq -T({\slicenormal}^\flat,e^i)=T^{0i}.
\ee
Write $z_{\rm s}=(\Fsp ,\Gsp )$ for the selected spatial components of the one-forms $F$ and the spatial duals of the two-forms $G$. The \emph{conservative variables} and the \emph{mathematical entropy density} on the slice are
\be
\label{eq:slice-variables-entropy-density}
 V=(z_{\rm s},E,P),\qquad r=-g(\particlecurrent,\slicenormal)={\particlecurrent}^0,\qquad \eta=-r=g(\particlecurrent,\slicenormal).
\ee
{We now introduce the corresponding convexity property of a spacetime current.}
\begin{definition}
\label{def:future-convex-current}
{Let $\mathcal J$ be a spacetime current whose normal density $-g(\mathcal J,\slicenormal)$ is a single-valued function of the conservative variables $V$ on a specified convex domain. The current is \underline{\emph{future-convex relative to $\slicenormal$}} on that domain if this density is convex as a function of $V$.}
\end{definition}
{For the mathematical entropy current $\mathcal J_\eta=-\particlecurrent$, this density is $\eta=-r$. The normal and conservative domain are fixed here. Constitutive limits and changes of observer require further hypotheses.}

The density $\eta$ is measured with respect to the induced volume on $\Sigma_t$. In coordinates with $x^0=t$, write ${\slicenormal}^\flat=-\alpha\,dt$ with lapse $\alpha>0$, and let $h$ be the induced metric. Then the corresponding coordinate density is $-|g|^{1/2}\,{\particlecurrent}^0=\sqrt{\det h}\,\eta$, with ${\particlecurrent}^0=dt(\particlecurrent)$ in this equality. At a fixed spacetime point the conservative components in the two conventions are related by an invertible linear map. Convexity is preserved by this map and multiplication by the positive factor $\sqrt{\det h}$.

\paragraph{{{Future convexity of the mathematical entropy current.}}}
The proof recovers the observed density from the stress and the spatial form components, then constructs a supporting plane for its negative.
\begin{theorem}[{Future convexity of the mathematical entropy at positive density}]
\label{prop:global-rest-current}
Fix a covariant augmented law with the selected forms orthogonal to the future timelike current $\particlecurrent$. The constitutive scalar is built from these fields, the metric and the volume tensors as in Proposition~\ref{prop:rest-entropy-test}.  Let the rest energy $\bar\rho(\Frest ,\Grest ,n)$ be defined on an open convex domain $\mathcal D$ of independent selected rest variables containing the triples satisfying the exact minor relations under consideration. The domain is a subset of the positive-density space
\be
\mathcal D\subseteq\mathbb R^{3r_F+3r_G}\times(0,\infty).
\ee
Impose the constitutive conditions below throughout the rest domain $\mathcal D$.
\begin{itemize}
\item[(i)] \emph{Rest convexity.} The energy $\bar\rho$ is convex, as required by relativistic polyconvexity in Definition~\ref{def:rest-polyconvexity}.
\item[(ii)] \emph{Regularity and positivity in matter.} The energy $\bar\rho$ is finite and {$C^2$}. Its constitutive enthalpy $H=\partial_n\bar\rho$ is positive, and its inertia matrix $\mathsf M$ in~\eqref{eq:rest-inertia-test} is positive definite.
\item[(iii)] \emph{Acoustic bound for all observers.} For every unit spatial covector $e$, the Hessian of the energy with respect to its independent rest-frame variables $A=D^2_{(\Frest,\Grest,n)}\bar\rho$ and the propagation map $\mathcal D_e$ of~\eqref{eq:retained-acoustic-map} satisfy $\mathcal D_e^TA\mathcal D_e\le\mathsf M$.
\end{itemize}
In an oriented orthonormal observer frame set $z=(\Fsp ,\Gsp ,r)$ and $\particlecurrent=r(1,v)$. For $|v|<1$, set $B(v)=(I-vv^T)^{1/2}$ and denote the corresponding rest variables by $\mathcal R_v z=(\Fsp B(v),\Gsp \operatorname{cof}B(v),r\det B(v))$. When $\mathcal R_v z\in\mathcal D$, the reduced energy and the stress components are
\be
\widetilde\rho(z,v)=\bar\rho(\mathcal R_v z),
\qquad E=\widetilde\rho-v\cdot D_v\widetilde\rho,
\qquad P=-D_v\widetilde\rho.
\ee
{Then the following conclusions hold.}
\begin{itemize}
\item[(a)] {\emph{Local recovery and convexity.}} For each such choice of fields, $D_v^2\widetilde\rho<0$ and the conservative map $\Psi:(z,v)\mapsto(\Fsp ,\Gsp ,E,P)$ has a {$C^1$} local inverse. On a sufficiently small convex conservative neighborhood of the image point, the mathematical entropy density $\eta_{\rm mat}(V)=-r(\Psi^{-1}(V))$ is $C^2$ and convex. At each point its Hessian has nullity $\dim\ker A$ at the corresponding triple of rest variables. In particular, it is positive definite wherever the Hessian of the energy with respect to its independent rest-frame variables is positive definite.

\item[(b)] {\emph{Convexity on a product domain.}} Suppose $Z$ is an open convex set of variables $z$ with $r>0$ and $W$ is an open convex set of velocities in $|v|<1$, with $\mathcal R_v z\in\mathcal D$ for every $(z,v)\in Z\times W$.
Then $\Psi$ restricts to a {$C^1$} diffeomorphism of $Z\times W$ onto its image $\Omega$. The density $\eta_{\rm mat}$ on $\Omega$ has a supporting plane at each point valid throughout $\Omega$, and is the restriction of a convex function of $V$. If $\mathcal D$ is the full rest domain, this conclusion applies with $Z=\mathbb R^{3r_F+3r_G}\times(0,\infty)$ and $W=\{|v|<1\}$.
\item[(c)] {\emph{A prescribed observer.}} For a fixed observer, (i) and (ii) suffice if one has
\be
 \tau^2\mathsf M-\mathcal D_k^T A\mathcal D_k>0,
\ee
at the chosen field values, where $(\tau,k)$ are the components of the observer's time conormal in the material rest frame. The product-domain conclusion holds if this strict inequality is satisfied at every point of that product.
\end{itemize}
\end{theorem}

{The inequality in (iii) is the light-cone condition of Proposition~\ref{prop:observer-acoustic-coercivity}. For a perfect fluid, the matter assumptions reduce to $\bar\rho'(n)>0$ and $0\le n\bar\rho''(n)/\bar\rho'(n)\le1$. When the law is defined on the full positive density range, the product domain is $(0,\infty)\times\{|v|<1\}$.}

\paragraph{{Local inversion and supporting planes.}}
{The sign of the velocity Hessian will give the supporting planes.}
\begin{proof}
{\bf 1.} {\emph{Concavity in the velocity.}} For fixed $v$, the map $\mathcal R_v$ is linear and invertible in $z$. Therefore $\widetilde\rho(\cdot,v)$ is convex on $\mathcal R_v^{-1}\mathcal D$, and its density derivative satisfies
\bse
\be
\widetilde\rho_r=\det B(v)\,\partial_n\bar\rho>0.
\ee
{The approximation argument following~\eqref{eq:LP3-local-matrix} applies to the present $C^2$ rest law. In particular, the quadratic-current identity~\eqref{eq:quadratic-observer-current} and the velocity-block identity below are understood through their continuous expressions in the first two constitutive derivatives.}

We first prove strict concavity in $v$, without assuming $D_z^2\bar\rho>0$. In the present observer chart the differential identity for $e=-r$ gives
\be
\lambda_E=-1/\widetilde\rho_r,\qquad
\lambda_P=v/\widetilde\rho_r,\qquad
\lambda_{(\Fsp ,\Gsp )}=D_{(\Fsp ,\Gsp )}\widetilde\rho/\widetilde\rho_r.
\ee
Indeed $dE-v\cdot dP=\widetilde\rho_r\,dr+D_{(\Fsp ,\Gsp )}\widetilde\rho\cdot d(\Fsp ,\Gsp )$ follows directly from $P=-D_v\widetilde\rho$. {For the smooth approximations, differentiating once more with $(\Fsp ,\Gsp ,r)$ fixed and passing to the limit gives}
\be
\label{eq:velocity-entropy-block}
K_{vv}=-\frac{D_v^2\widetilde\rho}{\widetilde\rho_r}.
\ee
\ese
To see that this block is positive definite, we use the quadratic current~\eqref{eq:quadratic-observer-current} of Proposition~\ref{prop:observer-acoustic-coercivity}. A nonzero velocity variation with $r,\Fsp ,\Gsp $ fixed changes the unit current direction: after a fixed Lorentz boost to the material rest frame its spatial velocity variation is nonzero. The rest time component of the quadratic current is therefore positive, since $\mathsf M>0$ and $A\ge0$. The light-cone bound in (iii) implies the strict matrix inequality for each timelike conormal. For a fixed observer the completed square in Proposition~\ref{prop:observer-acoustic-coercivity} uses only that observer's strict inequality: its second term is positive on the nonzero rest velocity increment and its first term is nonnegative. The contraction is $K_{vv}$ under the algebraic chart change, so the block in~\eqref{eq:velocity-entropy-block} is positive definite. The fixed-observer strict condition persists on a small neighborhood by continuity. The recovery and support arguments below use this positivity throughout the product under consideration.

\vskip.2cm \noindent {\bf 2.} {\emph{Recovery of the conservative variables.}} We set $L=-D_v^2\widetilde\rho>0$. The current block of the conservative Jacobian, in the order $(r,v)$ and $(E,P)$, is
\bse
\be
\begin{pmatrix}
\widetilde\rho_r+v\cdot P_r&v^TL\\
P_r&L
\end{pmatrix},
\qquad
\det D_{(r,v)}(E,P)=\widetilde\rho_r\det L>0.
\ee
Thus the full conservative map has a {$C^1$} local inverse, including when the rest-energy Hessian is semidefinite. For injectivity on $Z\times W$, we use the partial Legendre transform
\be
\Phi(z,P)\coloneqq\sup_{w\in W}\bigl(\widetilde\rho(z,w)+w\cdot P\bigr).
\ee
{The product-domain assumption ensures that} $\widetilde\rho(z,\cdot)$ is strictly concave on $W$ and $\widetilde\rho(\cdot,w)$ is convex on $Z$ for every $w\in W$. In particular, $\Phi$ is convex in $(z,P)$. At a point of $Z\times W$ with $P=-D_v\widetilde\rho(z,v)$, the stationary velocity $v$ is the unique maximizer, with $\Phi(z,P)=E$. We fix $(\Fsp ,\Gsp ,P)$ and take $r_1<r_2$ with $(\Fsp ,\Gsp ,r_i)\in Z$ for $i=1,2$. We suppose the supremum at $r_1$ is attained at $v_1\in W$. Evaluating the supremum at $r_2$ using $v_1$ and integrating $\widetilde\rho_r>0$ along the density segment gives
\be
\Phi(\Fsp ,\Gsp ,r_2,P)\ge\widetilde\rho(\Fsp ,\Gsp ,r_2,v_1)+v_1\cdot P
>\widetilde\rho(\Fsp ,\Gsp ,r_1,v_1)+v_1\cdot P=\Phi(\Fsp ,\Gsp ,r_1,P).
\ee
\ese
Thus for given $(\Fsp ,\Gsp ,E,P)$ there is at most one recovered $r$, and then $v$ is the unique maximizer. The injective local diffeomorphism is a diffeomorphism onto its open image $\Omega$.

\paragraph{{3. Supporting planes on a common domain.}}
We fix $V_0=(z_{{\rm s}0},E_0,P_0)\in\Omega$ with recovered density and velocity $r_0,v_0$. For $V\in\Omega$, we evaluate the supremum defining $\Phi(z_{\rm s},r,P)$ at $v_0$. The convexity of $\widetilde\rho(\cdot,v_0)$ gives
\bse
\be
\begin{aligned}
 E=\Phi(z_{\rm s},r,P)
 &\ge \widetilde\rho(z_{\rm s},r,v_0)+v_0\cdot P\\
 &\ge E_0+D_{z_{\rm s}}\widetilde\rho_0\cdot(z_{\rm s}-z_{{\rm s}0})
          +\widetilde\rho_{r0}(r-r_0)+v_0\cdot(P-P_0).
\end{aligned}
\ee
The subscript $0$ indicates evaluation at $(z_{{\rm s}0},r_0,v_0)$. Since $\widetilde\rho_{r0}>0$, isolating $-r$ gives
\be
\label{eq:entropy-affine-support-in-matter}
 \eta_{\rm mat}(V)\ge \eta_{\rm mat}(V_0)+\lambda(V_0)\cdot(V-V_0),
 \qquad
 \lambda(V_0)=\frac{(D_{z_{\rm s}}\widetilde\rho_0,-1,v_0)}{\widetilde\rho_{r0}}.
\ee
\ese
{On the local conservative image, the formula for $\lambda$ and the $C^1$ inverse show that $D_V\eta_{\rm mat}$ is $C^1$. The Hessian can also be read directly from the two Taylor remainders in the supporting inequality.  We denote by $D$ the difference between $\eta_{\rm mat}(V)$ and its supporting plane at $V_0$. Since $P=-D_v\widetilde\rho(z,v)$, one has
\be
\begin{aligned}
 \widetilde\rho_{r0}D
 &=\widetilde\rho(z,v_0)-\widetilde\rho(z_0,v_0)
      -D_z\widetilde\rho(z_0,v_0)\cdot(z-z_0)\\
 &\quad+\widetilde\rho(z,v)-\widetilde\rho(z,v_0)
      -D_v\widetilde\rho(z,v)\cdot(v-v_0).
\end{aligned}
\ee
The first remainder varies $z$ at fixed velocity, and the second varies the velocity at fixed $z$. Taylor expansion at $(z_0,v_0)$ gives the block diagonal quadratic form with blocks $D_z^2\widetilde\rho$ and $-D_v^2\widetilde\rho$, after multiplication by $\widetilde\rho_{r0}$. Equivalently, this form is the pullback of $\widetilde\rho_{r0}D_V^2\eta_{\rm mat}$ by the invertible Jacobian $D\Psi$. The velocity block is positive definite by Step~1, and $D_z^2\widetilde\rho=\mathcal R_{v_0}^T A\mathcal R_{v_0}$ at the base point. Thus the entropy Hessian has the same nullity as $A$.} The supremum of these affine functions is convex and agrees with $\eta_{\rm mat}$ on $\Omega$. Indeed, the inequality bounds every term from above, while the plane based at $V_0=V$ attains the bound. To obtain the local conclusion at any choice of fields with $\mathcal R_v z\in\mathcal D$, we choose a small product of open convex neighborhoods of $z$ and $v$ whose triples of rest variables lie in $\mathcal D$. Such a product exists by continuity of $(z,v)\mapsto\mathcal R_v z$ and openness of $\mathcal D$. Its conservative image is open and contains a convex neighborhood of the image point.
This completes the proof of Theorem~\ref{prop:global-rest-current}.
\end{proof}

{
\begin{corollary}[A criterion for a fixed observer]
\label{cor:fixed-observer-convex-entropy}
Let $Z$ be an open convex set of variables $z=(z_{\rm s},r)$ with $r>0$, and let $W$ be an open convex set of velocities. Suppose $\widetilde\rho\in C^1(Z\times W)$ is convex in $z$ at each fixed velocity, {concave} in $v$ at each fixed $z$, and satisfies $\widetilde\rho_r>0$. Define $P=-D_v\widetilde\rho$ and $E=\widetilde\rho+v\cdot P$.
\begin{itemize}
\item[(i)] \emph{Convex entropy.} The variables $(z_{\rm s},E,P)$ determine $r$. At each point of their image, $-r$ has an affine supporting plane valid on the whole image. It is therefore the restriction of a convex function.

\item[(ii)] \emph{Acoustic criterion.} For a $C^2$ rest law satisfying (i) and (ii) of Theorem~\ref{prop:global-rest-current} on a product $Z\times W$ as in that theorem, the concavity hypothesis follows from the non-strict inequality
\be
 \mathcal D_v^T A\mathcal D_v\le\mathsf M
 \qquad\text{at every triple of rest variables }\mathcal R_v z,\quad (z,v)\in Z\times W.
\ee
{The criterion also gives the entropy density when the conservative map has a singular derivative.}
\end{itemize}
\end{corollary}
{Section~\ref{og:proof-details} below uses the first-derivative criterion for unshifted spectral energies, whose second derivatives can diverge as a singular value tends to zero.}
\paragraph{{Proof of the fixed-observer criterion.}}
{We use the acoustic bound only in the direction of the chosen observer.}

\begin{proof}
Concavity suffices in Steps~2 and~3 of Theorem~\ref{prop:global-rest-current} to identify $E$ with the maximum of $\widetilde\rho(z,w)+w\cdot P$ on $W$, attained at the given velocity. The comparison of these maxima at two densities uses $\widetilde\rho_r>0$ and gives uniqueness of $r$. The support~\eqref{eq:entropy-affine-support-in-matter} then follows from convexity in $z$ at the reference maximizing velocity. Taking the supremum of these supports gives the convex extension.
For the last assertion, the observer's conormal in the rest frame has components $(\Gamma,\Gamma v)$, where $\Gamma=(1-|v|^2)^{-1/2}$. The completed square of Proposition~\ref{prop:observer-acoustic-coercivity} then gives $K_{vv}\ge0$. By~\eqref{eq:velocity-entropy-block} and $\widetilde\rho_r>0$, we have $D_v^2\widetilde\rho\le0$. This gives concavity on each convex velocity domain $W$. Convexity in $z$ and positivity of the density derivative follow from the rest energy by the linear transformation $\mathcal R_v$. For a perfect fluid the matrix inequality is $c_s^2|v|^2\le1$.
\end{proof}

}

The product hypothesis of Theorem~\ref{prop:global-rest-current} gives a supporting inequality between fields in the same conservative image, beyond the local convexity given by the acoustic bound.
For a fluid with $\bar\rho(n)={m \, n}+\kappa n^\gamma/\gamma$, $m\ge0$, $\kappa>0$ and $1\le\gamma\le2$, the Hessian of the energy with respect to its independent rest-frame variables is nonnegative, the inertia matrix is $n\bar\rho'(n)I_3>0$, and the squared sound speed is bounded by $\gamma-1$. The quadratic law in Section~\ref{subsec:three-invariant-class} below satisfies the hypotheses by the direct quadratic estimate. For the unshifted Ogden law we use Corollary~\ref{cor:fixed-observer-convex-entropy} in Section~\ref{og:proof-details} below.

%-----------------------------------------------------------------------------------------------------------------------------------------

\subsection{The mathematical entropy at degenerate constitutive limits}

\paragraph{Values of the convex extension.}
{We now use a lower semicontinuous convex extension to determine the mathematical entropy and its supporting covectors at the constitutive boundary of Definition~\ref{def:kinematic-constitutive-boundaries}.}

\begin{proposition}[The convex extension and its limiting values]
\label{thm:lsc-entropy-extension}
Let $\Omega$ be a nonempty subset of a finite-dimensional real vector space and let $f:\Omega\to\mathbb R$.
\begin{itemize}
\item[{(i)}] {\emph{Convex extension.}} The set
\be
\label{eq:entropy-closed-epigraph}
 \overline{\operatorname{co}}\{(V,s):V\in\Omega,\ s\ge f(V)\}
\ee
is the epigraph of a proper lower semicontinuous convex function $f_{\rm cl}$ if and only if $f$ has an affine lower bound on $\Omega$. The function $f_{\rm cl}$ is the greatest lower semicontinuous convex function bounded above by $f$ on $\Omega$, and it may take the value $+\infty$ outside $\Omega$. At any point of $\Omega$ at which an affine lower bound attains $f$, one has $f_{\rm cl}=f$.
\item[{(ii)}] {\emph{Boundary limits.}} If $\Omega$ is convex and $f$ is convex on $\Omega$, then $f$ has an affine lower bound. For every $V_0\in\overline\Omega$ and every $W$ in the relative interior of $\Omega$,
\be
 f_{\rm cl}(V_0)
 =\liminf_{\Omega\ni V\to V_0} f(V)
 =\lim_{\theta\downarrow0} f\big((1-\theta)V_0+\theta W\big).
\ee

\item[{(iii)}] {\emph{Entropy at a constitutive limit.}} Take $\Omega$ to be an image on which {Theorem}~\ref{prop:global-rest-current} gives supporting planes throughout $\Omega$. Then the mathematical entropy density $f=\eta_{\rm mat}=-r$ has an affine support at each point of $\Omega$. Thus $\eta_{\rm cl}\coloneqq f_{\rm cl}$ agrees with $\eta_{\rm mat}$ there. At a finite limit of constitutive fields and their stress at positive density with conservative variables $V_0$ and normal current component $r_0$, it is finite and $\eta_{\rm cl}(V_0)\le-r_0$. Equality holds if along one such sequence with variables $V_j$ the supporting covectors $\lambda_j$ satisfy $\limsup_j\lambda_j\cdot(V_0-V_j)\ge0$.
\item[{(iv)}] {\emph{Continuity at the boundary.}} {In the setting of (iii),} if $\Omega$ is convex, equality with $-r_0$ is equivalent to $\limsup_{\Omega\ni V\to V_0}r(V)=r_0$. A continuous trace of $r$ with value $r_0$ suffices. {More generally, let $\vartheta\colon[0,+\infty)\to[0,+\infty)$ be continuous and nondecreasing, with $\vartheta(0)=0$. If $\Omega$ is convex and $|r(V)-r_0|\le\vartheta(|V-V_0|)$ for $V\in\Omega$ near $V_0$, then $|\eta_{\rm cl}(V)+r_0|\le\vartheta(|V-V_0|)$ for $V\in\overline\Omega$ near $V_0$. In particular, a pointwise H\"older estimate for $r$ at $V_0$ holds with the same exponent for the extended entropy.}
\end{itemize}
\end{proposition}

{For a perfect fluid the proper density can tend to zero while the limiting current is lightlike and nonzero. The criterion in (iii) concerns the density measured by the observer. It allows unbounded supporting covectors, provided their pairing with $V_0-V_j$ has nonnegative upper limit. For the positive-rest-mass fluid the distinction between zero current and zero stress is described in Proposition~\ref{prop:vacuum-stress-subgradient} below.}

\paragraph{{Construction of the convex envelope.}}
{Vertical sections of the closed convex epigraph determine the extension.}
\begin{proof}
The set in~\eqref{eq:entropy-closed-epigraph} is nonempty, closed and convex and is invariant under upward translations of $s$. An affine lower bound for $f$ bounds this set from below. Each nonempty vertical section then has a finite minimum and these minima define $f_{\rm cl}$, with value $+\infty$ at empty sections. Conversely, if this set is the epigraph of a proper function, separate a point strictly below a finite vertical section from the epigraph. The coefficient of $s$ in the separating affine inequality is positive, since the epigraph contains a point with the same horizontal coordinate. The inequality therefore gives an affine lower bound for $f_{\rm cl}$ and hence for $f$ on $\Omega$.

The epigraph of every lower semicontinuous convex function bounded above by $f$ on $\Omega$ contains the set in~\eqref{eq:entropy-closed-epigraph}. This proves maximality. If an affine lower bound touches $f$ at a point of $\Omega$, it bounds $f_{\rm cl}$ from below there and gives equality. When $\Omega$ and $f$ are convex, a support at any point of the relative interior gives an affine lower bound on the affine hull of $\Omega$. Extending that affine function to the ambient space gives the bound in (i). The initial epigraph is convex, so its closure gives the lower limit formula. To see that one segment suffices, we choose $W$ in the relative interior. For $0<\theta\le1$, the point $(1-\theta)V_0+\theta W$ lies in the relative interior, where $f_{\rm cl}=f$ by the supporting-plane property. If $f_{\rm cl}(V_0)$ is finite, convexity bounds the upper limit along this segment by $f_{\rm cl}(V_0)$, and lower semicontinuity gives the opposite inequality for the lower limit. If $f_{\rm cl}(V_0)=+\infty$, lower semicontinuity gives divergence to $+\infty$ along the segment.

For a sequence with $V_j\to V_0$ and $r(V_j)\to r_0$, the points $(V_j,-r(V_j))$ belong to the epigraph of $\eta_{\rm cl}$. Closedness gives the upper bound at $V_0$. Each supporting plane at $V_j$ bounds the extended entropy from below, so
\be
 \eta_{\rm cl}(V_0)\ge -r(V_j)+\lambda_j\cdot(V_0-V_j).
\ee
Taking the upper limit proves the claimed equality whenever the upper limit of the last term is nonnegative. For convex $\Omega$, the existence of the supporting planes gives convexity of $\eta_{\rm mat}$ on $\Omega$. The lower limit formula is then the negative of the upper limit of $r$ in~$\Omega$. {The lower limit formula transfers {this modulus of continuity} to every point of $\overline\Omega$ near $V_0$.}
\end{proof}

\paragraph{Supporting covectors at the boundary.}
For weak--strong comparison, we will need a finite affine support of the entropy at the reference point.

\begin{definition}
\label{def:entropy-supporting-covector}
For a convex mathematical entropy $\eta$, a finite covector $\lambda_0$ is a \underline{\emph{supporting covector}} at $V_0$ if $\eta(V_0)$ is finite and $\eta(V)\ge\eta(V_0)+\lambda_0\cdot(V-V_0)$ on the domain of $\eta$.
\end{definition}

We write $\lambda_0\in\partial\eta(V_0)$ for this supporting inequality.
Let $V_0$ be a point of finite entropy in the closure of the convex conservative image $\Omega$. Put $V_\theta=(1-\theta)V_0+\theta W$ for $W$ in the relative interior of $\Omega$. There is a finite supporting covector at $V_0$ if and only if, for one such $W$,
\be
 \lim_{\theta\downarrow0}\frac{\eta(V_\theta)-\eta(V_0)}{\theta}>-\infty.
\ee
An affine support at $V_0$ bounds these slopes from below. Conversely, work in the affine hull and write $h=\eta-\eta(V_0)$. If the slopes are bounded below by $-L$, then $h(V_\theta)\ge-L\theta$. A closed relative ball of radius $\delta$ about $W$ lies in the relative interior and has $h\le M$ for some $M$. Convexity gives $h(V_\theta\pm\theta\delta e)\le\theta M$ for unit vectors $e$ in the affine hull. Testing any supporting covector $\lambda_\theta$ at these two points gives $|\lambda_\theta\cdot e|\le(M+L)/\delta$. A convergent subsequence of these covectors gives a finite support at $V_0$ by the convergence of the entropy along the segment (cf.~Proposition~\ref{thm:lsc-entropy-extension}(ii)). In particular, whenever a finite support exists at $V_0$,  the covectors in the affine hull stay bounded on these segments.

\begin{definition}
\label{def:continuous-constitutive-boundary}
Choose constitutive variables for a branch of the law and let $\mathcal D_b$ be a subset of the boundary of its positive-density domain in these variables. The law has a \underline{\emph{continuous constitutive extension to $\mathcal D_b$}} if its energy and stress possess limits at each point of $\mathcal D_b$ and these limits give continuous functions on the enlarged domain. The extension is \underline{\emph{locally H\"older continuous of exponent $\alpha$}}, $0<\alpha\le1$, if the extended energy and stress satisfy local H\"older estimates of exponent $\alpha$ in the chosen variables.
\end{definition}

{For a fluid, the required extension is a property of the energy and stress on the chosen enlarged domain. The limiting behavior of the current and the velocity must be specified when these variables are used at vacuum.}

The chosen constitutive variables may include a bounded velocity orthogonal to the forms at zero current.  Several limiting stress tensors can then occur at the same $(F,G,\particlecurrent)$. H\"older continuity of the energy and stress in those variables controls the map to $V$. To identify the entropy on a convex conservative image, we use continuity of the inverse density $r(V)$ at the boundary. In particular, all points of the completed constitutive flux graph over the same $V_0$ must agree on ${\particlecurrent}^0$ for this identification to hold on the whole graph. The existence of a finite supporting covector at $V_0$ is a separate condition used in Section~\ref{section=6} below for weak--strong stability.

%-----------------------------------------------------------------------------------------------------------------------------------------

\subsection{An explicit description of the constitutive boundary}

\paragraph{{Entropy at zero current.}}
\label{sec:closed-vacuum-entropy}

We now give a class of laws for which the conservative image, the entropy trace and the limiting stress can all be described. We use $z_{\rm s}=(\Fsp,\Gsp )$, $r={\particlecurrent}^0$ and $\particlecurrent=r(1,v)$ in the observer frame, with $q(v)\coloneqq\sqrt{1-|v|^2}$ and $n=rq(v)$. The other form components are determined by the orthogonality condition~\eqref{eq:450bis-2}. Consider the constitutive law
\be
\label{eq:massive-vacuum-class}
\raug(F,G,\particlecurrent;g)=\mathcal L(F,G,s;g)+\raug_{\rm d}(n),\qquad \raug_{\rm d}(n)={m \, n}+\mathcal U(n),\quad s=-g(\particlecurrent,\particlecurrent),\quad m\ge0.
\ee
The term $\mathcal L$ may be independent of $s$. Its domain contains the reconstructed fields orthogonal to the causal direction $(1,v)$ with $z_{\rm s}\in\mathbb R^{3r_F+3r_G}$, $r\ge0$ and $|v|\le1$. In these variables we put
\be
\ell(z_{\rm s},r,v)\coloneqq\mathcal L(F(z_{\rm s},v),G(z_{\rm s},v),r^2(1-|v|^2);g),
\qquad
\widetilde\rho=\ell+mrq+\mathcal U(rq).
\ee
{For the comparison at fixed spatial forms $z_{\rm s}$ and momentum $P$, we maximize over the closed velocity ball. We write $\Phi$ for this maximum at observed density $r$ and $\Phi_0$ for its value at $r=0$.}
\be
\label{eq:vacuum-hamiltonian}
\Phi(z_{\rm s},r,P)\coloneqq\max_{|v|\le1}
\{\widetilde\rho(z_{\rm s},r,v)+v\cdot P\},\qquad
\Phi_0(z_{\rm s},P)\coloneqq\Phi(z_{\rm s},0,P).
\ee
In $\Phi_0$, the velocity parametrizes limits from matter. We put $\Omega=\{(z_{\rm s},E,P):E>\Phi_0(z_{\rm s},P)\}$. To recover the particle density from the conservative variables, we need the maximizing velocity at positive observed density to lie inside the velocity ball. The following inward slope condition ensures this for every momentum $P$: for each $z_{\rm s}$, $r>0$ and $|e|=1$, one has
\be
\label{eq:boundary-inward-velocity}
 \limsup_{t\downarrow0}
 \frac{\widetilde\rho(z_{\rm s},r,(1-t)e)-\widetilde\rho(z_{\rm s},r,e)}{t}
 =+\infty.
\ee

\paragraph{{The conservative image and its entropy.}}
{The inward slope ensures that a positive observed density is represented by a timelike current.}
\begin{theorem}[Extension of the convex mathematical entropy to the constitutive boundary]
\label{thm:closed-vacuum-entropy}
Consider the law~\eqref{eq:massive-vacuum-class} on the full reconstructed domain. On the positive-density domain, require a convex $C^2$ rest energy with positive enthalpy, a positive definite inertia matrix and wave speeds bounded by the light speed in the sense of Theorem~\ref{prop:global-rest-current}(iii). Assume that the reduced elastic term is nonnegative and continuous up to the closed velocity ball and that $\mathcal U\in C^1([0,\infty))\cap C^2((0,\infty))$ is nondecreasing with $\mathcal U(0)=\mathcal U'(0)=0$. Require the infinite inward slope~\eqref{eq:boundary-inward-velocity} at each lightlike velocity with $r>0$.
\begin{itemize}
\item[(i)] {The positive-density image is the open convex set $\Omega$, with closure $E\ge\Phi_0(z_{\rm s},P)$.}
\item[(ii)] {In $\Omega$, the observed density $r(V)>0$ is determined by $\Phi(z_{\rm s},r(V),P)=E$. Its negative has the proper lower semicontinuous convex extension}
\be
\label{eq:closed-current-entropy}
\eta(V)=
\begin{cases}
-r(V),&E>\Phi_0(z_{\rm s},P),\\
0,&E=\Phi_0(z_{\rm s},P),\\
+\infty,&E<\Phi_0(z_{\rm s},P)
\end{cases}
\ee
{On $\overline\Omega$, this extension agrees with the supremum of the affine entropy supports in matter.}
\item[(iii)] {The conservative projection of the independent finite-flux graph is $\overline\Omega$. For every point of the completed constitutive flux graph with conservative projection $V$, one has $\eta(V)=-{\particlecurrent}^0$. Zero proper density implies zero current, and the fields with zero current project exactly onto $E=\Phi_0(z_{\rm s},P)$. Their stress may be nonzero.}
\end{itemize}
\end{theorem}

To verify~\eqref{eq:boundary-inward-velocity}, it suffices that $\ell(z_{\rm s},r,(1-t)e)-\ell(z_{\rm s},r,e)\ge-Ct$ and $\rho_{\rm d}(n)/n^2\to+\infty$ as $n\downarrow0$. Indeed, the proper density at $(1-t)e$ is $n=r\sqrt{2t-t^2}$, so the density contribution divided by $t$ tends to $+\infty$. This includes $m>0$, as well as $m=0$ with $\mathcal U(n)=a_{\rm fl}n^\gamma$, $a_{\rm fl}>0$ and $1<\gamma<2$, whenever the elastic term satisfies {this lower bound}. The convex extension agrees with the flux limits even if the supporting covectors in matter are unbounded near the boundary.

For exact minors of a gradient of rank $k$, the span of the quadratic forms has dimension $\binom{k}{2}$, and $\particlecurrent=0$ when $k\le2$. The extension includes the rank-two case with one independent quadratic minor as well as ranks zero and one, where they all vanish. For independent augmented variables, rank loss of either form family may also occur in matter, where the entropy is already defined.

At a boundary point $V_0$ with $E_0=P_0=0$, the hypotheses of the theorem already give a quantitative estimate. Indeed, evaluating the maximum at $v=0$ gives $E\ge\rho_{\rm d}(r)$ throughout $\overline\Omega$. Thus, if $\rho_{\rm d}(n)\ge c_0n^\gamma$ near zero for $c_0>0$ and $\gamma>1$, continuity of $r$ at $V_0$ gives, locally,
\be
 |\eta(V)-\eta(V_0)|=r(V)\le (E/c_0)^{1/\gamma}
 \le c_0^{-1/\gamma}|V-V_0|^{1/\gamma}.
\ee
When $m>0$, nonnegativity of $\ell$ and $\mathcal U$ also gives the global bound
\be
 E\ge\max_{|v|\le1}\{mr\sqrt{1-|v|^2}+v\cdot P\}
 =\sqrt{m^2r^2+|P|^2}.
\ee
In particular, $|\eta(V)-\eta(V_0)|\le |V-V_0|/m$ and $(0,-1/m,0)$ is a supporting covector at $V_0$ in the variables $(z_{\rm s},E,P)$.

For estimates at boundary points with nonzero energy or momentum, we also control the change of elastic energy as the density varies.
\paragraph{{Continuity at the boundary.}}
{We estimate the decay of the observed particle density at the boundary.}
\begin{corollary}[Continuity estimates at the constitutive boundary]
\label{cor:boundary-entropy-rates}
Assume the hypotheses of Theorem~\ref{thm:closed-vacuum-entropy} and a $C^{1,1}_{\rm loc}$ extension of $\ell$ to a neighborhood of the reconstructed domain, with $\ell_r(z_{\rm s},0,v)=0$.
\begin{itemize}
\item[(i)] If $m>0$, then for each compact set of $(z_{\rm s},P)$ there are $c,r_0>0$ such that
\be
 \Phi(z_{\rm s},r,P)-\Phi_0(z_{\rm s},P)\ge cr^2,\qquad 0\le r\le r_0.
\ee
Consequently the entropy restricted to $\overline\Omega$ has a pointwise H\"older estimate of exponent $1/2$ at every conservative boundary point.
\item[(ii)] If $\rho_{\rm d}(n)\ge c_0n^\gamma$ for small $n$, with $c_0>0$ and $1<\gamma<2$, the entropy on $\overline\Omega$ has a pointwise H\"older estimate of exponent $(2-\gamma)/(2\gamma)$ at every conservative boundary point. This includes $m=0$ and $\mathcal U(n)=c_0n^\gamma$. If the maximum defining $\Phi_0$ at a boundary point $V_0$ is attained at an interior velocity, the pointwise exponent improves to $1/\gamma$, with a constant depending on that velocity. {At a point admitting such an interior maximizer, suppose further that $m=0$ and $\rho_{\rm d}(n)=a n^\gamma+o(n^\gamma)$ with $a>0$. Let} $\mathcal M_0$ be the compact set of maximizers defining $\Phi_0$ at $V_0$, and put $q_{\max}=\max_{v\in\mathcal M_0}\sqrt{1-|v|^2}>0$. Along the inward energy direction one has
\be
 \eta(z_{{\rm s}0},E_0+\epsilon,P_0)
 =-\bigl(\epsilon/(a q_{\max}^\gamma)\bigr)^{1/\gamma}+o(\epsilon^{1/\gamma})
 \qquad(\epsilon\downarrow0).
\ee
Thus the exponent $1/\gamma$ is optimal at every such point.
\end{itemize}
\end{corollary}

\label{app:closed-vacuum-entropy}

\paragraph{{Construction and boundary values.}}
{We first recover positive observed densities by maximizing over the velocity, then pass to the constitutive boundary.}
\begin{proof}[Proof of Theorem~\ref{thm:closed-vacuum-entropy} and Corollary~\ref{cor:boundary-entropy-rates}]
{\bf 1.} \emph{The conservative image and the partial Legendre transform.} We work at a fixed spacetime point in the chosen observer frame. The reduced density is continuous on the closed velocity ball, locally uniformly in $(z_{\rm s},r)$. Hence the maximum defining $\Phi$ exists, is finite and depends continuously on $(z_{\rm s},r,P)$. For each fixed velocity, $\widetilde\rho$ is convex in $(z_{\rm s},r)$, including $r=0$ by continuity. Taking the maximum proves convexity of $\Phi$.

 We fix $r>0$ and $|e|=1$. Replacing $e$ by $(1-t)e$ changes the momentum term by $-te\cdot P$. Condition~\eqref{eq:boundary-inward-velocity} gives arbitrarily small positive $t$ for which the gain in reduced energy exceeds this loss. Thus every maximizing velocity is interior. Strict concavity with respect to the velocity, proved in {Theorem}~\ref{prop:global-rest-current}, makes the maximizer unique, with
\be
P=-D_v\widetilde\rho,\qquad E=\widetilde\rho+v\cdot P=\Phi(z_{\rm s},r,P).
\ee

The derivative $\widetilde\rho_r=q\,\partial_n\bar\rho>0$ gives strict increase of $\Phi$ between positive densities when we evaluate at the maximizing velocity for the smaller density. The same argument gives $\Phi(r)>\Phi(0)$ if the zero-density maximum has an interior maximizer. Otherwise, we choose a boundary maximizer $e$. Monotonicity at interior velocities and continuity give $\widetilde\rho(z_{\rm s},r,e)\ge\widetilde\rho(z_{\rm s},0,e)$; the strict gain just proved at density $r$ again gives $\Phi(r)>\Phi(0)$.  To bound the feasible densities, we evaluate at $v=0$ and use the supporting line of the convex rest energy at $r=1$. For $r\ge1$,
\be
\label{eq:uniform-density-from-energy}
 E\ge\Phi(z_{\rm s},r,P)\ge
 \bar\rho(z_{\rm s},1)+\partial_r\bar\rho(z_{\rm s},1)(r-1).
\ee
The slope is bounded below by a positive constant on compact ranges of $z_{\rm s}$. Continuity and strict monotonicity therefore give a unique $r>0$ for each $E>\Phi_0$, and no triple of constitutive fields with positive density below that threshold. This proves the description of the image of the map to conservative variables.

\vskip.2cm \noindent {\bf 2.} \emph{The entropy at the limiting field values.} For $V=(z_{\rm s},E,P)$ with $E\ge\Phi_0$, we define the largest feasible density
\bse
\be
r_{\max}(V)\coloneqq
\max\{r\ge0:\Phi(z_{\rm s},r,P)\le E\}.
\ee
It equals the positive inverse for $E>\Phi_0$ and zero for $E=\Phi_0$. The feasible set
\be
\{(V,r):r\ge0,\ E\ge\Phi(z_{\rm s},r,P)\}
\ee
\ese
is closed and convex, and its density coordinate is locally bounded for bounded $E,z_{\rm s}$ by the preceding estimate. The epigraph of $-r_{\max}$ is its projection with the additional linear inequality $s\ge-r$. It is convex, and is closed by that local bound on $r$. This proves properness, lower semicontinuity and convexity of~\eqref{eq:closed-current-entropy}.

At each boundary point $V_0=(z_{{\rm s}0},\Phi_0(z_{{\rm s}0},P_0),P_0)$, the recovered density $r(V)$ tends to zero along every sequence from $\Omega$. The local density bound gives convergent subsequences, and any subsequential limit $r_0$ satisfies $\Phi(z_{{\rm s}0},r_0,P_0)=\Phi_0(z_{{\rm s}0},P_0)$ by continuity of $\Phi$. Strict increase in $r$ gives $r_0=0$. The lower limit formula of Proposition~\ref{thm:lsc-entropy-extension} now gives~\eqref{eq:closed-current-entropy}, with an empty epigraph above conservative points outside $\overline\Omega$. We denote by $\eta_{\rm sup}$ the supremum of the affine functions supporting the entropy at triples of constitutive fields with positive density. Convexity gives $\eta_{\rm sup}\le\eta$, with equality on $\Omega$. At a boundary point $V_0$, we choose $W\in\Omega$, set $V_t\coloneqq(1-t)V_0+tW$, and write $g(t)\coloneqq\eta(V_t)$. Convexity and lower semicontinuity give $g(t)\to g(0)=0$. For small positive $t$ the tangent plane at $V_t$, evaluated at $V_0$, is $g(t)-tg'(t)$. The convex secant inequalities imply
\be
g(t)-g(0)\le tg'(t)\le g(2t)-g(t).
\ee
Thus $tg'(t)\to0$, and these supporting planes tend to zero at $V_0$. Consequently $\eta_{\rm sup}(V_0)=\eta(V_0)=0$. The two extensions agree on the physical conservative domain $\overline\Omega$.

\vskip.2cm \noindent {\bf 3.} \emph{Identification of limiting fields and stress.}
We consider triples of constitutive fields with positive density converging in $(F,G,\particlecurrent,T)$, with $r_j\to r_0$ and, along a subsequence, $v_j\to v_0$ in the closed velocity ball. Since $E_j=\Phi(z_{{\rm s}j},r_j,P_j)$, continuity gives $E_0=\Phi(z_{{\rm s}0},r_0,P_0)$. The maximizer identity also passes to the limit, so $v_0$ maximizes the objective at the limiting parameters. If $r_0>0$, all such maximizers are interior by Step~1. Hence $n_0=r_0\sqrt{1-|v_0|^2}>0$ and the entropy has its matter value $-r_0$. If $r_0=0$, then ${\particlecurrent}_0=0$ and $E_0=\Phi_0(z_{{\rm s}0},P_0)$. This is the zero entropy value and allows nonzero limiting stress.

\paragraph{{4. Approximation at the constitutive boundary.}}
We first prove a bound on the stress components that will allow us to approach every conservative boundary point by triples of constitutive fields with positive density. Nonnegativity of the reduced density and the partial Legendre formula give
\bse
\be
E=\sup_{|w|<1}\{\widetilde\rho(z,w)+w\cdot P\}\ge|P|.
\ee
This holds in every orthonormal observer frame. By covariance the stress is symmetric, so its energy current is future causal for every future timelike observer. The closed causal cone gives the same conclusion for the null directions $e_{\hat0}\pm e_j$. In the chosen frame,
\be
|T^{i0}\mp T^{ij}|\le T^{00}\mp T^{0j}\le2T^{00},
\qquad |T^{i0}|\le T^{00}.
\ee
The full stress tensor therefore satisfies
\be
\label{eq:general-dominant-bound}
|T|\le CE,
\ee
\ese
and these component inequalities pass to finite-flux limits.

It remains to recover every conservative boundary point. We fix $(z_{\rm s},P)$ and choose $r_j\downarrow0$ with maximizing velocities $v_j$. The energies satisfy $E_j=\Phi(z_{\rm s},r_j,P)\to\Phi_0$. The orthogonality condition~\eqref{eq:450bis-2} bounds all components of the original forms, and the currents tend to zero. The stress components are bounded by~\eqref{eq:general-dominant-bound}. The full fluxes therefore have a convergent subsequence. Its conservative projection is the prescribed boundary point. This proves recovery of $\overline\Omega$.

{\vskip.2cm \noindent {\bf 5.} \emph{A boundary estimate for positive rest mass.} \bse
We fix a compact set of $(z_{\rm s},P)$. {We apply Taylor's formula to the rest-frame restriction of $\mathcal L$, in the proper density $n$ with its rest form components fixed. Its first density derivative at $n=0$ vanishes by the assumption on $\ell_r$ at $v=0$, and {its first density derivative is Lipschitz on the compact range of rest forms}. Since those forms are independent of $r$ at fixed $(z_{\rm s},v)$ and $n=r\sqrt{1-|v|^2}$, this gives
\be
 |\ell(z_{\rm s},r,v)-\ell(z_{\rm s},0,v)|
 \le Cn^2=Cr^2(1-|v|^2).
\ee
Continuity includes $|v|=1$.} We choose $v_0$ to maximize the expression for $\Phi_0$ and set $v_\tau=(1-\tau)v_0$. The function $\ell(z_{\rm s},0,v)+v\cdot P$ is uniformly Lipschitz in $v$ on the closed ball, so its loss from $v_0$ to $v_\tau$ is at most $L\tau$. Since $\mathcal U\ge0$ and $1-|v_\tau|^2\ge\tau$ for $0<\tau\le1$, we obtain
\be
 \Phi-\Phi_0\ge mr\sqrt{1-|v_\tau|^2}-Cr^2(1-|v_\tau|^2)-L\tau.
\ee
\ese
For small $r$, the first two terms are at least $mr\sqrt{1-|v_\tau|^2}/2$. Taking $\tau$ to be a sufficiently small fixed multiple of $r^2$ {proves the required bound}, uniformly on the compact set. The maximum defining $\Phi_0$ is locally Lipschitz in $(z_{\rm s},P)$. Together with $r(V)\to0$ at the boundary and $\eta=-r$, the bound gives the pointwise $1/2$-H\"older estimate.

\paragraph{{A polytropic lower bound.}}
\vskip.2cm \noindent {\bf 6.} \emph{A polytropic lower bound at zero density.} We suppose $\rho_{\rm d}(n)\ge c_0n^\gamma$ for small $n$, with $1<\gamma<2$. For small $r$, {the quadratic density error from $\ell$ satisfies} $Cr^2(1-|v|^2)\le\frac{c_0}{2}(r\sqrt{1-|v|^2})^\gamma$. With the same $v_\tau$ and $L$ as above, we obtain
\be
 \Phi-\Phi_0\ge \tfrac{c_0}{2}r^\gamma\tau^{\gamma/2}-L\tau.
\ee
We choose $\tau=\kappa r^{2\gamma/(2-\gamma)}$, with $\kappa>0$ sufficiently small. Since $\gamma/2<1$, this yields $\Phi-\Phi_0\ge cr^{2\gamma/(2-\gamma)}$ for small $r$, uniformly on compact sets of $(z_{\rm s},P)$. The local Lipschitz bound for $\Phi_0$ therefore gives {the desired H\"older estimate}. }  To check sharpness for the class, we take the perfect-fluid energy $\rho_{\rm d}(n)=c_0 n^\gamma$ and fix $P\ne0$. We put $s=|v|$ and take $v=sP/|P|$, with $s\uparrow1$. The perfect-fluid stress gives $|P|=\gamma c_0 r^\gamma s(1-s^2)^{(\gamma-2)/2}$ and
\be
{E-|P|=\frac{|P|(1-s)}{s}\left(1-\frac{\gamma-1}{\gamma}(1+s)\right).}
\ee
{As $s\uparrow1$, the bracket tends to $(2-\gamma)/\gamma>0$. Hence $r$ is asymptotic to a positive constant times $(E-|P|)^{(2-\gamma)/(2\gamma)}$. This gives a nonzero boundary point where the exponent cannot be increased.}  For the improved estimate at a boundary point with an interior maximizer $v_0$, we set $q_0=\sqrt{1-|v_0|^2}>0$. Evaluating the maximum defining $\Phi$ at $v_0$, and using the quadratic estimate for $\ell(z_{\rm s},r,v_0)-\ell(z_{\rm s},0,v_0)$, gives for small $r$
\be
 E-E_0\ge c_0q_0^\gamma r^\gamma-Cr^2-C\bigl(|z_{\rm s}-z_{{\rm s}0}|+|P-P_0|\bigr).
\ee
Since $\gamma<2$, the $r^2$ term is absorbed for sufficiently small $r$. Continuity of $r$ at $V_0$ then yields $|\eta(V)|=r\le C|V-V_0|^{1/\gamma}$ on $\overline\Omega$ near $V_0$. If in addition $m=0$ and $\rho_{\rm d}(n)=a n^\gamma+o(n^\gamma)$, the quadratic estimate for $\ell$ gives, uniformly for $|v|\le1$,
\be
 \widetilde\rho(z_{{\rm s}0},r,v)+v\cdot P_0
 =\ell(z_{{\rm s}0},0,v)+v\cdot P_0
  +a r^\gamma(1-|v|^2)^{\gamma/2}+o(r^\gamma).
\ee
Evaluating at a maximizer for which $\sqrt{1-|v|^2}=q_{\max}$ gives the lower limit. Conversely, every limit of maximizers for positive $r$ belongs to $\mathcal M_0$, which gives the upper limit. These two limits give $\Phi(z_{{\rm s}0},r,P_0)-E_0=a q_{\max}^\gamma r^\gamma+o(r^\gamma)$. The entropy asymptotic follows by inversion. This completes the proof of Theorem~\ref{thm:closed-vacuum-entropy} and Corollary~\ref{cor:boundary-entropy-rates}.
\end{proof}

%-----------------------------------------------------------------------------------------------------------------------------------------

\paragraph{Stress and supporting planes at zero current.}
Under the additional assumptions on $\ell$ in Corollary~\ref{cor:boundary-entropy-rates}, fix $(z_{\rm s},P)$ and let $\mathcal M_0$ be the set of velocities attaining the maximum that defines $\Phi_0(z_{\rm s},P)$. The right derivative of $\Phi$ at zero density is
\be
 \lim_{r\downarrow0}\frac{\Phi(z_{\rm s},r,P)-\Phi_0(z_{\rm s},P)}{r}
 =m\max_{v\in\mathcal M_0}\sqrt{1-|v|^2}.
\ee
Indeed, $\ell_r(z_{\rm s},0,v)=0$ and $\mathcal U'(0)=0$, so the first-order term in the maximized expression is $mr\sqrt{1-|v|^2}$, with a uniform $o(r)$ remainder on the closed velocity ball. Evaluating at maximizers for $r=0$ gives the lower bound. Along any sequence $r\downarrow0$, a subsequence of maximizing velocities converges to an element of $\mathcal M_0$ and gives the upper bound.  Along the segment with fixed $(z_{\rm s},P)$ and increasing $E$, the entropy equals $-r$ at $E=\Phi(z_{\rm s},r,P)$. The inward entropy slope is finite from below precisely when the displayed derivative is positive, which requires $m>0$ and at least one timelike maximizer. {The same criterion characterizes finite supporting covectors in all conservative variables.} In particular, for $m=0$ the entropy is continuous at zero current and has an infinite inward slope. {We prove the following result in Section~\ref{app:vacuum-stress-subgradient} below.}

\begin{proposition}[Stress and entropy subgradients at zero current]
\label{prop:vacuum-stress-subgradient}
Under the hypotheses of Corollary~\ref{cor:boundary-entropy-rates} with $m>0$, also assume that $T_{\mathcal L}$ is continuous on the reconstructed domain in $(z_{\rm s},r,v)$ .
\begin{itemize}
\item[{(1)}] A zero-current state of the independent graph has the form
\be
\label{eq:vacuum-dust-graph}
T=T_{\mathcal L}(F,G,0)+D z\otimes z,\qquad
z=(1,v),\quad |v|\le1,\quad i_zF=i_zG=0,
\ee
where $D\ge0$ and $D=0$ if $|v|<1$. Conversely every such choice of fields and stress is a pointwise finite-flux limit of independent triples of constitutive fields with positive density. A positive $D$ is a residual rank-one null stress and has no finite entropy subgradient. {For such fields and stress with conservative variables $V_0=(z_{\rm s},E_0,P_0)$, at fixed $z_{\rm s},P_0$ the entropy satisfies
\be
 \eta(z_{\rm s},E_0+\delta,P_0)=-\frac{\sqrt{2D}}m\sqrt\delta+o(\sqrt\delta)
 \qquad(\delta\downarrow0).
\ee
Thus its pointwise H\"older exponent at $V_0$ is exactly $1/2$ on the independent graph.}
\item[{(2)}] At $V_0=(z_{{\rm s}0},\Phi_0(z_{{\rm s}0},P_0),P_0)$ a finite subgradient exists if and only if the maximization at zero current in~\eqref{eq:vacuum-hamiltonian} has an interior maximizer $v_0$. For each such maximizer,
\be
\label{eq:vacuum-entropy-subgradient}
\lambda_0=\frac{1}{m\sqrt{1-|v_0|^2}}
\bigl(\ell_{z_{\rm s}}(z_{{\rm s}0},0,v_0),-1,v_0\bigr)
\in\del\eta(V_0).
\ee
A timelike velocity orthogonal to the forms on the branch with $D=0$ therefore determines a supporting plane at zero current. The largest energy component among the supporting covectors is $-1/(m\max_{v\in\mathcal M_0}\sqrt{1-|v|^2})$. {At each such point $V_0$, the entropy is pointwise Lipschitz on $\overline\Omega$:
\be
 0\le-\eta(V)\le|\lambda_0|\,|V-V_0|.
\ee
The exponent one is optimal on the independent graph.}
\end{itemize}
\end{proposition}

%-----------------------------------------------------------------------------------------------------------------------------------------

\part[Weak-strong stability and local-in-time existence theory]{\LARGE Weak-strong stability and local-in-time existence theory}
\label{part:stability-and-existence}

%============================================================================

\section{The notion of weak solutions with cavitation}
\label{sec:weak-local04}

\subsection{Weak formulation and initial data sets}

\paragraph{Independent augmented fields and exact minor relations.}
{
{We now define entropy-admissible weak solutions and specify their initial terms. We then study material parametrizations through cavities and the shock jump conditions.\par}}
{The law may admit zero current with residual stress or nonzero null current. In an empty cavity $\particlecurrent=T=0$ and the weak formulation requires no velocity.}

\begin{definition}
\label{def:weak-local-complete}
Let $(\Mcal,g)$ be the prescribed oriented and time-oriented Lorentzian spacetime, with Levi--Civita connection $\nabla$ and volume measure $dV_g$. Fix the constitutive law $\raug$ and its finite-flux graph ${\Cl(\mathcal Z_{\raug})}$.
\begin{itemize}
{\item[(i)] An \underline{\emph{augmented weak solution}} {consists of measurable fields $(F,G,\particlecurrent)$ whose associated constitutive stress $T$ satisfies $(F,G,\particlecurrent,T)\in{\Cl(\mathcal Z_{\raug})}$ almost everywhere, with locally integrable flux components. The constitutive continuation at the constitutive boundary is specified as above.} Require the two-form $G$ to be antisymmetric, $\particlecurrent$ to be future causal or zero, and the orthogonality condition~\eqref{eq:450bis-2} expressed by ${\particlecurrent}^\mu F_\mu^A=0$ and ${\particlecurrent}^\mu G_{\mu\nu A}=0$ to hold almost everywhere. The closure equations and stress balance laws are
\be
dF^A=0,\qquad dG_A=0,\qquad \nabla_\nu T^{\mu\nu}=0
\quad\hbox{in }\mathcal D'(\Mcal).
\ee
}
{\item[(ii)] A \underline{\emph{configuration weak solution}} {consists of a locally Lipschitz configuration $f$ with locally integrable stress $T$ given by the constitutive law and its specified continuation at the constitutive boundary.} Set $F=df$, $G=G(F)$ and $\particlecurrent=\particlecurrent(F)$ almost everywhere, allowing rank loss. Require $(F,G,\particlecurrent,T)$ to belong almost everywhere to the admissible constitutive extension, $\particlecurrent$ to be future causal or zero, and $\nabla_\nu T^{\mu\nu}=0$ to hold in the sense of distributions. The current is conserved by the kinematic identity (cf.~Proposition~\ref{prop:kinematic-weak-local} below), so its mathematical entropy inequality follows.
}
\item[(iii)] Any given weak solution is \underline{\emph{entropy-admissible}} if it satisfies the mathematical entropy inequality in the sense of distributions:
\be
\label{entropy-ineq-null}
-\nabla_\mu {\particlecurrent}^\mu\le0.
\ee
\end{itemize}
\end{definition}

For example, the stress balance laws in the sense of distributions require $\int_\Mcal T^{\mu\nu}\nabla_\nu\varphi_\mu\,dV_g=0$ for every smooth compactly supported covector-valued test-function. The mathematical entropy inequality means that the particle-number current $\particlecurrent$ has nonnegative production.
Membership in ${\Cl(\mathcal Z_{\raug})}$ gives~\eqref{eq:minor-wedge-constraint} in both classes of solutions, including at zero current. The mixed identity~\eqref{eq:exact-mixed-minor} is used only for exact polynomial minors. Table~\ref{tab:weak-solution-terminology} summarizes the terminology.
\begin{table}[!htbp]
\centering
\begingroup
\begin{tabular}{@{}p{.42\linewidth}p{.54\linewidth}@{}}
\toprule
Term & Meaning \\
\midrule
{Augmented weak solution} & Independent fields satisfying {the balance laws of Definition~\ref{def:weak-local-complete}} \\
Configuration weak solution & Locally Lipschitz $f$, exact minors and stress balance laws \\
Solutions satisfying the exact minor relations & Augmented solutions with $G=G(F)$ and $\particlecurrent=\particlecurrent(F)$ a.e. \\
Entropy admissibility & $-\nabla_\mu {\particlecurrent}^\mu\le0$ in the sense of distributions \\
{Strong solution} & Fields satisfying Definition~\ref{def:entropy-strong-solution} below \\
Empty cavity & Open region on which $\particlecurrent=T=0$ \\
\bottomrule
\end{tabular}
\caption{Terminology for solutions and cavitation.}
\label{tab:weak-solution-terminology}
\endgroup
\end{table}

\paragraph{Initial data sets on a time slice.}
Consider a periodic coordinate slab $[0,T)\times\mathbb T^3$. Write $|g|=|\det(g_{\mu\nu})|$, so that $dV_g=|g|^{1/2}\,dt\,dx$, and denote the stress and form evolution equations by $\del_\mu\mathscr F_a^\mu=\mathscr S_a$, with $a$ indexing their rows. For smooth test-functions $\psi^a,\phi$ compactly supported in $[0,T)\times\mathbb T^3$, with $\phi\ge0$, the weak formulation with initial data is
\be
\begin{aligned}
\int_0^T\!\int_{\mathbb T^3}\bigl(\mathscr F_a^\mu\del_\mu\psi^a+\mathscr S_a\psi^a\bigr)\,dx\,dt
+\int_{\mathbb T^3}\mathscr F_{a,0}^{0}\psi^a(0,x)\,dx&=0,\\
\int_0^T\!\int_{\mathbb T^3}|g|^{1/2}{\particlecurrent}^\mu\del_\mu\phi\,dx\,dt
+\int_{\mathbb T^3}|g(0,x)|^{1/2}{\particlecurrent}_0^0\phi(0,x)\,dx&\le0.
\end{aligned}
\ee
The subscript $0$ denotes the prescribed initial data sets. The spatial exterior constraints hold in the sense of distributions on the slab.

\paragraph{Selected fields and the polynomial minor relations.}
The definition of weak solutions uses the selected fields $F^a,G_b,\particlecurrent$ and their constitutive stress. When $r_F=r_G=0$, it reduces to the four stress balance laws and the mathematical entropy inequality, as for generalized and perfect fluids. Any complementary kinematic fields obey their own equations and are outside the entropy estimate for the selected fields.

\begin{proposition}[Reduction to the particle-number fluid formulation]
\label{prop:fluid-particle-number-reduction}
Fix a perfect-fluid law $\raug(\particlecurrent;g)=\rrestnew(n)$, with $n=\sqrt{-g(\particlecurrent,\particlecurrent)}$,  and omit the form components. In matter assume $\rrestnew'(n)>0$ and put
\be
\label{eq:fluid-torus-density-dictionary}
\mu=\rrestnew(n),\qquad p=n\rrestnew'(n)-\rrestnew(n),\qquad
\particlecurrent=nu.
\ee
{Use the finite constitutive limits of this law at the constitutive boundary.}
\begin{itemize}
\item[(i)] {The entropy-admissible augmented weak solution of Definition}~\ref{def:weak-local-complete} { is exactly a solution of}
\be
\label{eq:fluid-particle-number-weak-system}
\nabla_\nu T^{\mu\nu}=0,\qquad
\nabla_\mu\particlecurrent^\mu\ge0,
\qquad T=(\mu+p)u\otimes u+p\,g^{-1}\quad(n>0)
\ee
{in the sense of distributions, with the same prescribed initial data sets and constitutive limits.}
\item[(ii)] {The stress formula holds in matter; its finite constitutive limits specify the stress at zero or null current.}
\end{itemize}
\end{proposition}
{The same normalization gives the particle current used for relativistic fluids.\footnote{In~\cite{TorusT2-1-2026,TorusT2-2-2026} the particle-density function satisfies $d\log N_{\rm fl}=d\mu/(\mu+p)$. Since $d\mu=\bar\rho'(n)\,dn$ and $\mu+p=n\bar\rho'(n)$, the choice $N_{\rm fl}=n$ identifies the two currents, with the same metric and constitutive limits.}}
\begin{proof}
{ We remove the absent form components from Definition}~\ref{def:weak-local-complete}{ . The remaining equations and inequality are those in}~\eqref{eq:fluid-particle-number-weak-system}{ , with their prescribed initial terms. The same constitutive law and its finite limits determine the stress in both descriptions, including at cavitation. In local coordinates, the particle inequality reads $\del_\mu(|g|^{1/2}\particlecurrent^\mu)\ge0$. For exact configuration minors, the kinematic identity makes this inequality an equality.}
\end{proof}

For a weak solution of Definition~\ref{def:weak-local-complete}, the restriction $(F,G,\particlecurrent)\in\mathcal C_{\rm min}$ to the manifold of spacetime minors means, almost everywhere,
\be
\label{eq:physical-weak-minors}
G=G(F),\qquad {\particlecurrent}=\particlecurrent(F).
\ee
The stress belongs to ${\Cl(\mathcal Z_{\raug})}$ and the prescribed initial data set satisfies the same minor relations. The closure equations, stress balance laws, entropy inequality and prescribed initial terms still apply. The minor relations make sense for unbounded gradients as well as for locally Lipschitz configurations. Their weak continuity requires separate hypotheses.

In fixed affine material coordinates, a closed bounded one-form $F^A$ has a locally Lipschitz potential. On a smaller ball, the potentials of smooth mollifications, fixed at a common base point, converge along a subsequence by their uniform Lipschitz bounds. The limit has gradient $F^A$. If the configuration itself gives the material parametrization, invertibility is a separate condition (cf.~Section~\ref{sec:weak-coordinate-equivalence} below).

\paragraph{{Current conservation for exact minors.}}
{The particle three-form is the pullback of the material volume form, including where the configuration gradient loses rank.}
\begin{proposition}[Distributional kinematic current identity]
\label{prop:kinematic-weak-local}
\begin{itemize}
\item[(i)] For every locally Lipschitz configuration and smooth material volume form, the current defined by~\eqref{vmu} obeys $\nabla_\mu {\particlecurrent}^\mu=0$ in the sense of distributions, at any rank. 
\item[(ii)] {If $f$ is locally Lipschitz up to an initial slice $t=0$, its initial particle density is the spatial minor $|g|^{1/2}{\particlecurrent}^0_{\rm in}=\omega(f_0)\det D_xf_0$, where $\Omega=\omega(X)\,dX^1\wedge dX^2\wedge dX^3$ in material coordinates and $f_0=f|_{t=0}$. The conservation law holds with this initial term, and the normal flux has a weakly-$*$ continuous representative in $L^\infty_{\rm loc}$ with this initial value.}
\item[(iii)] In fixed constant-volume material coordinates it also holds  when the forms $F^A\in L^{p_A}_{\rm loc}$ are closed in the sense of distributions, with $1\le p_A\le\infty$ and $\sum_{A=1}^3p_A^{-1}\le1$, using $1/\infty=0$, with $\particlecurrent$ its ordinary cubic minor.  The same conclusion holds for a smooth material volume form if the coordinate components of $f$ belong to the indicated $W^{1,p_A}_{\rm loc}$ spaces,  $f$ lies almost everywhere in one material chart and its current is locally integrable.
\end{itemize}
\end{proposition}
\paragraph{{Proof of the kinematic identity.}}
{We approximate the exterior products while retaining the spatial pullback on the initial slice.}
\begin{proof}
{\bf (i).} \emph{Conservation for Lipschitz configurations.} We work on a smaller compact spacetime chart whose Lipschitz image lies in a material chart. Smooth mollifications $f_\varepsilon$ converge uniformly to $f$, with uniformly bounded gradients converging strongly in every finite $L^p$. The smooth volume coefficients evaluated at $f_\varepsilon$ therefore converge uniformly, and the associated cubic minors and currents converge strongly in every finite $L^p$ there. Since $d(f_\varepsilon^*\Omega)=f_\varepsilon^*(d\Omega)=0$, each smooth current has zero divergence. Integrating the approximating currents against compactly supported test-functions and passing to the limit and localizing proves the geometric identity.

\vskip.2cm \noindent {\bf (ii).} \emph{Initial particle density.}  Reflect $f$ across $t=0$ and mollify in space at scale $\delta$ and then in time at scale $\delta^2$. The bulk gradients converge in every finite local $L^p$ with uniform bounds. At $t=0$, the spatial gradient differs from the spatially smoothed $D_xf_0$ by $\Ocal(\delta)$, by the Lipschitz bound in time. The smooth conservation laws pass to {the initial identity above}. For $j^\mu\coloneqq|g|^{1/2}{\particlecurrent}^\mu$, conservation makes the pairing of $j^0$ with each smooth compactly supported spatial test $\phi$ Lipschitz in time, with constant at most $\|j^{\rm sp}\|_{L^\infty}\|\nabla\phi\|_{L^1}$ on the test support. The initial identity, local boundedness and density of the tests give the asserted weak-$*$ continuous representative.

\vskip.2cm \noindent {\bf (iii).} \emph{Closed forms of lower regularity.}  For constant volume coefficients, we convolve the closed forms on a smaller coordinate ball. The smooth convolutions remain closed and converge strongly in each finite $L^{p_A}_{\rm loc}$. The infinite-exponent forms remain uniformly bounded and converge almost everywhere. We expand the difference between $F_\varepsilon^1\wedge F_\varepsilon^2\wedge F_\varepsilon^3$ and $F^1\wedge F^2\wedge F^3$ one factor at a time. H\"older's inequality and dominated convergence give strong convergence of the cubic minors in $L^1_{\rm loc}$, and the same distributional passage proves the claim. The case $p_1=p_2=p_3=p\ge3$ is included; for example, $(p_1,p_2,p_3)=(\infty,2,2)$ is also allowed.  For a variable material volume form, we write $\Omega=\omega(X)\,dX^1\wedge dX^2\wedge dX^3$ on the material chart. We choose smooth cutoffs $0\le\chi_k\le1$ with compact support in the chart and $\chi_k\to1$ pointwise. Each density $\chi_k\omega$ extends smoothly by zero to a bounded function on $\mathbb R^3$. Mollifying $f$ gives strong $L^1_{\rm loc}$ convergence of the cubic derivative products as above; the cutoff densities evaluated at the mollified maps converge almost everywhere and are uniformly bounded. The smooth pullback identity therefore gives conservation for the current with density $\chi_k\omega$. Finally these currents converge to $\particlecurrent$ in $L^1_{\rm loc}$ by dominated convergence, using the assumed integrability of $\particlecurrent$. This concerns strong approximation of fixed closed forms; it is not a weak compactness assertion for arbitrary $L^1$ sequences of minors.
\end{proof}

Part (ii) uses Lipschitz regularity up to the initial slice. Part (iii) concerns conservation in the interior.
H\"older's inequality gives local integrability of the unweighted cubic minors in Proposition~\ref{prop:kinematic-weak-local}(iii). For a variable material volume form it gives local integrability of the current whenever $\omega(f)$ is locally essentially bounded, for instance if the essential image of each compact spacetime subset lies in a compact subset of the material chart. A single material chart alone does not ensure this integrability. 

Multiplying the left-hand side of a distributional equation by a merely bounded gradient does not generally extend the smooth derivation of stress balance laws to the fields of a weak solution. Definition~\ref{def:weak-local-complete} imposes these laws in the sense of distributions. A nonzero current production measure is allowed in the independent augmented class but cannot arise from this proposition's exact Lipschitz configuration current.

%-----------------------------------------------------------------------------------------------------------------------------------------

%-----------------------------------------------------------------------------------------------------------------------------------------

\subsection{Weak formulation in Eulerian and material coordinates with cavitation}
\label{sec:weak-coordinate-equivalence}

\paragraph{Configuration maps and changes of variables for weak solutions.}
The configuration $f(t,x)$ assigns a point of the material manifold to each spacetime point. Where the assignment is invertible on each slice, its inverse is the direct deformation. There the spatial part of $F=df$ and the gradient of the direct deformation are inverse matrices. We use $\chi(t,x)$ for the map to material coordinates. It can be $f$ itself or a separate map, such as the inverse of the material flow in Section~\ref{sec:material-flow-at-cavitation} below.

For locally bi-Lipschitz material coordinates, we follow Wagner~\cite{Wagner1,Wagner} to transform the stress balance laws with their geometric sources, the closure equations, the mathematical entropy inequality, and their initial terms. The integrable fields may include jumps in the configuration gradient.

\begin{theorem}[Weak formulation in Eulerian and Lagrangian coordinates with cavitation]
\label{un-prop:weak-coordinate}
On a spacetime coordinate domain $\mathcal O$ with smooth prescribed metric $g$, let $\chi\in W^{1,\infty}_{\rm loc}(\mathcal O;\mathbb R^3)$ define an orientation-preserving  locally bi-Lipschitz homeomorphism $\Phi(t,x)=(t,\chi(t,x))$ onto its image. Let $\Psi(t,X)=(t,y(t,X))$ be its inverse and let $a\coloneqq\det D_x\chi>0$ almost everywhere. For a locally integrable vector density $P$ and source density $S$, set
\be
\label{eq:Piola-local04}
\widetilde P^b=\left(\frac{\del_\mu\Phi^b}{a}P^\mu\right)\circ\Psi,
\qquad \widetilde S=\left(\frac Sa\right)\circ\Psi,
\qquad b=0,1,2,3.
\ee
\begin{itemize}
\item[(i)] \emph{Balance laws and inequalities.} The equation $\del_\mu P^\mu=S$ is equivalent to $\del_b\widetilde P^b=\widetilde S$ in the sense of distributions, and either inequality sign is preserved in both directions.
\item[(ii)] \emph{Closure and minor relations.} For a locally integrable one-form or two-form $\alpha$, the equations $d\alpha=0$ and $d(\Psi^*\alpha)=0$ are equivalent in the sense of distributions, using the almost-everywhere derivative of $\Psi$ in the pullback. This pullback preserves the polynomial minor relations, with the particle current represented by its three-form, and the ranks of $F$ and $G$, including where the current vanishes or is nonzero and null.

\item[(iii)] \emph{Measure sources and initial terms.}  The balance-law and inequality conclusions also hold when $S$ is a signed Radon measure on the bulk: its transformed source is $\widetilde S=\Phi_\#S$. In particular, a nonnegative entropy-production measure is transported to its nonnegative pushforward. For the initial terms, assume the bulk measure is locally finite up to the initial slice and has no mass on that slice. If the maps extend as  mutually inverse locally bi-Lipschitz homeomorphisms to the initial slices and the bulk fields are locally integrable up to them, the equivalences include locally integrable prescribed initial densities and tangential forms. With $\Psi_0(X)\coloneqq y(0,X)$, the initial terms are given by
\be
\label{eq:material-initial-data}
\widetilde P_{\rm in}^0
 =\det D_X\Psi_0\, P_{\rm in}^0\circ\Psi_0,\qquad
 \widetilde\alpha_{\rm in}
 =\Psi_0^*\alpha_{\rm in}.
\ee
\end{itemize}
\end{theorem}

A source measure on the closed slab can also have an initial atom. In the convention
\be
 \int P^\mu\partial_\mu\phi+\int P^0_{\rm in}\phi(0,\cdot)+\int\phi\,dS=0,
\ee
{write $S=S_{>0}+\delta_{t=0}\otimes S_0$. Composing the initial test with $\Phi_0=\chi(0,\cdot)$ pushes the effective initial measure $P^0_{\rm in}\,dx+S_0$ forward by $\Phi_0$. This applies to arbitrary locally finite initial measures in both directions, including the initial atom of a nonnegative entropy-production measure.}

{Condition $\particlecurrent^\mu\partial_\mu\chi^A=0$ makes the spatial particle flux in material coordinates vanish by~\eqref{eq:Piola-local04}. With $\chi=f$ and exact minors, the spatial Jacobian is $|g|^{1/2}\particlecurrent^0$ in normalized material volume coordinates, so this map loses rank in empty cavities.}

\paragraph{Transport by material maps that lose rank.}
When the material map loses rank, the change of variables involves the pushforward of measures.
\begin{proposition}[Transport of the distributional formulation by material maps that may lose rank]
\label{prop:collapsed-material-coordinates}
Let $\Phi(t,x)=(t,\chi(t,x))$ be a proper  locally Lipschitz map between coordinate slabs, with source $\mathcal O$. For a locally integrable vector density $P$ and a signed Radon measure $S$, define Radon measures by
\be
\label{eq:singular-material-pushforward}
 \widetilde P^b\coloneqq\Phi_\#\bigl((\del_\mu\Phi^b)P^\mu\,dt\,dx\bigr),
 \qquad
 {\widetilde S\coloneqq\Phi_\#S}.
\ee
{For an integrable source density $s$, the corresponding measure is $S=s\,dt\,dx$.} For a locally integrable $k$-form $\alpha$, $k=1,2$, define the image form with measure coefficients by its action on smooth compactly supported $(4-k)$-forms $\zeta$,
\be
 \langle\widetilde\alpha,\zeta\rangle\coloneqq
 \int_{\mathcal O}\alpha\wedge\Phi^*\zeta .
\ee
\begin{itemize}
\item[(1)] An Eulerian balance law, inequality or closure equation gives the corresponding distributional relation for these measures. Each image measure arising from an integrable coefficient has an absolutely continuous part from the set where $\det D_x\chi\ne0$ and a singular part from the set where $\det D_x\chi=0$. The density of the absolutely continuous part is obtained by dividing each source coefficient in the pushforward by $|\det D_x\chi|$ and summing over the preimages where this determinant is nonzero.
\item[(2)] If $\Phi$ extends as a proper  locally Lipschitz map with the initial slices included, the bulk fields are locally integrable up to those slices, and $S$ is locally finite there without mass on the initial slice, the implication includes locally integrable prescribed initial densities and tangential forms. With $\Phi_0(x)\coloneqq\chi(0,x)$, the initial density is $\widetilde P_{\rm in}^0\coloneqq(\Phi_0)_\#(P_{\rm in}^0\,dx)$, and the initial form is defined by pairing $\alpha_{\rm in}$ with the pullback by $\Phi_0$ of the spatial test-form. The same decomposition holds for these initial measures, with the spatial determinant $\det D\Phi_0$.
\item[(3)] The image of a $k$-form has no contribution from points where $\operatorname{rank}D_x\chi<3-k$. The same threshold applies to its initial tangential form, with $D\Phi_0$ in place of $D_x\chi$. If $\chi=f$ is the configuration in normalized material-volume coordinates, the image of the densitized causal particle current with exact cubic minors is absolutely continuous, including across the image of a collapsed cavity. Its spatial components vanish, and its time component counts the preimages of the material label with positive spatial Jacobian, for almost every $(t,X)$. For these exact minors the image forms are $\operatorname{mult}\,dX^A$ and $\tfrac12 \operatorname{mult}\Omega_{ABC}\,dX^B\wedge dX^C$, where $\operatorname{mult}$ is that multiplicity. If the current is conserved, the multiplicity is constant on each connected time interval at almost every fixed material label. On a target product slab reaching the initial slice, it equals the transformed initial particle density under the hypotheses of (2). If $f$ is locally Lipschitz up to the initial slice, the initial terms in the weak formulation for the exact forms and current are given by the corresponding minors of $f(0,\cdot)$. Their images are the same material forms and particle density multiplied by the multiplicity of $f(0,\cdot)$.
\item[(4)] On a domain where $\Phi$ is an orientation-preserving  locally bi-Lipschitz homeomorphism with no other preimages, the transformations agree with those in Theorem~\ref{un-prop:weak-coordinate}. In general, the implications can fail in reverse.
\end{itemize}
\end{proposition}

%-----------------------------------------------------------------------------------------------------------------------------------------

\paragraph{\texorpdfstring{Proofs of the change-of-variables results.}{Proofs of the change-of-variables results.}}

Properness gives compact support for the pulled-back tests. A scalar test $\theta\circ\Phi$ is Lipschitz. We mollify it, with reflection at $t=0$ for tests meeting the initial slice. Uniform convergence gives the initial and source pairings, and bounded derivatives converging almost everywhere give the bulk pairing. A nonnegative smoothing kernel preserves the sign.

For form tests, reflect the map across $t=0$ and smooth at spatial scale $h$ and  time scale $\ell=o(h)$. The smoothed maps $\Phi_h$ converge uniformly and in every finite local $W^{1,p}$ with uniform Lipschitz bounds. Time smoothing changes the spatially smoothed derivative at zero by at most $C\ell/h$, by the time Lipschitz bound and differentiation of the kernel. Thus the tangential derivatives converge to $D_x\Phi_0$ in $L^p$ on the initial slice. With a cutoff equal to one near the pulled-back supports, $d(\Phi_h^*\theta)=\Phi_h^*(d\theta)$ passes by dominated convergence in the bulk and the initial pairing.

\begin{proof}[Proof of Theorem~\ref{un-prop:weak-coordinate}]
{The spacetime Jacobian is $a=\det D_x\chi$ because time is preserved. We apply the scalar pullback identity and change variables to obtain~\eqref{eq:Piola-local04}. The inverse transform follows by applying the same argument to $\Psi$, since $D\Psi\circ\Phi=(D\Phi)^{-1}$ almost everywhere. The source for a measure is $\Phi_\#S$, whose sign is preserved. The spatial change of variables on the initial slice gives the first identity in~\eqref{eq:material-initial-data}.

For the closure equations the form pullback identity gives the equivalence by applying it in each direction. The initial tangential form becomes $\Psi_0^*\alpha_{\rm in}$. At almost every point, pullback commutes with exterior products, including the quadratic and cubic minor operations, and preserves their ranks because $D\Psi$ and its exterior powers are invertible. When $F=df$ for a Lipschitz configuration, the chain rule gives $\Psi^*F=d(f\circ\Psi)$, and the current three-form becomes $(f\circ\Psi)^*\Omega$.}

{We next write the elastic equations in material form.} For the inverse spatial map $y$, we put $M\coloneqq D_Xy$, $j\coloneqq\det M$ and $v\coloneqq\del_ty$, all in material coordinates. The inverse identities give
\be
(D_x\chi)\circ\Psi=M^{-1},\qquad
(\del_t\chi)\circ\Psi=-M^{-1}v,\qquad
a\circ\Psi=j^{-1}.
\ee
Thus~\eqref{eq:Piola-local04} can equivalently be written
\be
\label{eq:material-relative-flux}
\widetilde P^0=j(P^0\circ\Psi),\qquad
\widetilde P^A=j(M^{-1})^A{}_i
\bigl(P^i\circ\Psi-v^i(P^0\circ\Psi)\bigr),\qquad
\widetilde S=j(S\circ\Psi).
\ee
The spatial material flux is relative to the motion of the chosen material coordinates. For coordinates transported by the material velocity, it is the flux relative to the moving particles. The change of variables preserves the entropy inequality and transports the initial entropy density by the same initial spatial Jacobian.

For the curved-spacetime stress equations, we apply this result row by row to
\be
P^\nu=|g|^{1/2}T^{\mu\nu},\qquad
S=-|g|^{1/2}\Gamma^\mu_{\alpha\beta}T^{\alpha\beta},
\ee
with $\mu$ fixed in the original smooth coordinates. This index specifies the row throughout the change of independent variables. The material metric $\widetilde g=(D\Psi)^T(g\circ\Psi)D\Psi$ is defined almost everywhere and is locally essentially bounded. It preserves the causal type of the transformed current by the pointwise change of basis. The stress balance laws differentiate the transformed fluxes and use the Christoffel symbols of the original smooth chart in their sources.
\end{proof}

\paragraph{{The pushforward by a map that loses rank.}}
{We pull back test functions and test forms from the image.}
\begin{proof}[Proof of Proposition~\ref{prop:collapsed-material-coordinates}]
{We apply the pullback identity above without using an inverse map. For each target test $\theta$, the definitions of the image measures give}
\bse
\be
\begin{aligned}
\int\del_b\theta\,d\widetilde P^b+\int\theta\,d\widetilde S
   +\int\theta(0,X)\,d\widetilde P_{\rm in}^0
&=\int P^\mu\del_\mu(\theta\circ\Phi)\,dt\,dx
  {+\int(\theta\circ\Phi)\,dS}\\
&\quad+\int \theta(0,\Phi_0(x))P_{\rm in}^0(x)\,dx .
\end{aligned}
\ee
This proves the material equation from the Eulerian one. The converse requires enough pulled-back tests to determine the Eulerian law. An arbitrary locally finite signed measure in the initial data set $M$ transforms as $(\Phi_0)_\#M$, since the pairing uses $\theta(0,\Phi_0(x))$. In particular $M=P^0_{\rm in}\,dx+S_0$ includes an initial source atom, with $S_{>0}$ used in the bulk. We write $\mathcal O_{\rm r}$ and $\mathcal O_{\rm s}$ for the differentiability points where $|\det D\Phi|$ is positive and zero, respectively. The area formula on a countable Lipschitz cover gives a null image for $\mathcal O_{\rm s}$. For a locally integrable coefficient $B$, the contribution from $\mathcal O_{\rm r}$ is
\be
 \frac{d\,\Phi_\#(B\mathbf 1_{\mathcal O_{\rm r}}\,dt\,dx)}{dt\,dX}(z)
 =\sum_{p\in\mathcal O_{\rm r}\cap\Phi^{-1}(z)}
       \frac{B(p)}{|\det D\Phi(p)|}.
\ee
Properness and the formula for $|B|$ make the sum locally integrable. The part from $\mathcal O_{\rm s}$ is singular. This applies to $B=(\del_\mu\Phi^b)P^\mu$, to an integrable source density, and to the contractions of $\alpha$ with minors of $D\Phi$. On an orientation-preserving bi-Lipschitz domain with no other preimages the sum becomes~\eqref{eq:Piola-local04}.

For a form, the bounded derivative on the compact inverse image of the test support bounds the pairing by the supremum norm of the test times the local $L^1$ norm of $\alpha$. Thus the image form has Radon measure coefficients. In coordinates with orientation $dt\wedge dx$, the distributional closure equation tested with the pullback of a smooth target $(3-k)$-form $\theta$ gives
\be
 \int_{\mathcal O}\alpha\wedge\Phi^*(d\theta)
 +(-1)^k\int_{\Sigma_0}\alpha_{\rm in}\wedge\Phi_0^*(\theta|_{t=0})=0 .
\ee
Here $\Sigma_0$ has the spatial coordinate orientation. The tangential approximation at $t=0$ justifies the initial pairing. The three-dimensional area formula for $\Phi_0$ gives the same decomposition with respect to $dX$, now for $P^0_{\rm in}$ and the contractions of $\alpha_{\rm in}$ with minors of $D\Phi_0$.

The rank identity $\operatorname{rank}D\Phi=1+\operatorname{rank}D_x\chi$ gives the threshold by the degree $4-k$ of the test-form. On the initial slice the test has degree $3-k$, yielding the same threshold. For $\chi=f$ in normalized material-volume coordinates, the cubic minor gives $|g|^{1/2}{\particlecurrent}^0=\det D_x f\ge0$. On the singular set causality then gives $\particlecurrent=0$. On the regular set the time component in~\eqref{eq:singular-material-pushforward} counts each preimage once, and the spatial components vanish by ${\particlecurrent}^\mu\del_\mu f^A=0$. The exact forms $F^A=\Phi^*(dX^A)$ and $G_A=\Phi^*(\tfrac12\Omega_{ABC}\,dX^B\wedge dX^C)$ pair with a pulled-back test through a top-degree pullback. Its integrand contains $\det D_xf$, so the same multiplicity $\operatorname{mult}$ multiplies the image forms. The transformed continuity equation gives $\del_t \operatorname{mult}=0$. By Fubini on a countable cover by product boxes, $\operatorname{mult}$ is constant along each connected time interval for almost every label. On a product slab the initial term specifies that constant.

\paragraph{{Initial terms and the reverse implication.}}
For the initial terms, extend $f$ to negative times by $f(0,\cdot)$. The Lipschitz extension has closed minors by Proposition~\ref{prop:kinematic-weak-local}. Its negative-time pairings give the initial tangential minors. Along slices tending to zero with nonnegative determinant, uniform convergence of the maps and boundedness of their gradients give distributional convergence of the minors by induction in $df^{A_1}\wedge\cdots\wedge df^{A_k}=d(f^{A_1}df^{A_2}\wedge\cdots\wedge df^{A_k})$. Thus the initial determinant is nonnegative, and the initial area formula gives the multiplicity for the exact forms and current.

\ese
{
The pushforward can lose the entire normal flux supported in a cavity collapsed to a hypersurface. A converse therefore needs additional information on the interior flux.
}
This completes the proof of Proposition~\ref{prop:collapsed-material-coordinates}.
\end{proof}

For the reverse implication, we use dual integrability of the transformed flux and the derivative of the inverse. Let $\Phi$ be a locally Lipschitz homeomorphism as above, with $\det D\Phi>0$ almost everywhere, and let its continuous inverse $\Psi$ belong locally to $W^{1,p}$ for $1\le p< +\infty$. Suppose the transformed vector-density coefficients are locally in $L^{p'}$, with $1/p+1/p'=1$. For the initial term, assume that $\Phi$ extends locally Lipschitz to the initial slice and that the extended maps are continuous inverses there. Require the same Sobolev and flux bounds up to that slice, and take the bulk source locally finite up to it, with zero mass on it. Use the pushforward initial measure of Proposition~\ref{prop:collapsed-material-coordinates}. The initial measure can be singular.

The area formula gives~\eqref{eq:Piola-local04}. The first-order expansion of $\Phi$ with invertible derivative gives $D\Psi(\Phi(x))=(D\Phi(x))^{-1}$ outside a set with null image. For smooth $\varphi$, the pullback $\varphi\circ\Psi$ is continuous and locally $W^{1,p}$. Reflect it for tests meeting the initial slice and mollify with a nonnegative kernel. Its derivatives converge in $L^p$ for the flux pairing, and uniform convergence gives the source and initial pairings. The inverse derivative formula then gives the Eulerian equation or inequality.

For scalar balance laws the Jacobian may vanish on a set of positive measure if the flux vanishes there. More precisely, allow $\det D\Phi\ge0$ and require $P=0$ almost everywhere on $\{\det D\Phi=0\}$, with the same assumptions on the continuous Sobolev inverse, transformed flux and initial map. The area formula then gives $\widetilde P$ by~\eqref{eq:Piola-local04} on the image of the full-rank set, with zero density elsewhere. In pulling a scalar test back by the inverse, the calculation recovers $\int P^\mu\partial_\mu\varphi$ over the full-rank set, which is the complete Eulerian flux pairing by the vanishing assumption. The source and initial measure terms are recovered on the whole domain by continuity of the inverse and the pushforward identities. An initial source atom can be included in the initial measure as above.

For $k=1,2$, suppose $\alpha=0$ almost everywhere where $\det D\Phi=0$. If $p\ge4-k$ and $\widetilde\alpha=\Psi^*\alpha$ on the image of the full-rank set, extended by zero elsewhere, belongs locally to $L^{p/(p-4+k)}$, the converse holds for bulk closure. Here and below the space is $L^\infty$ at a zero denominator. Mollify $\Psi$ near the compact pulled-back support, using a cutoff equal to one there. Uniform convergence of the maps and $L^p$ convergence of their derivatives give $L^{p/(4-k)}$ convergence of the pullbacks of $d\theta$ by multilinearity of the minors. H\"older's inequality gives the transformed pairing, and the area formula recovers the Eulerian pairing since $\alpha=0$ on the singular set.

For closure with an initial term, use the same map extensions and bulk bounds up to $t=0$. Require $\alpha_{\rm in}=0$ where $\det D\Phi_0=0$, $\Psi_0\in W^{1,p_0}_{\rm loc}$ for $3-k\le p_0< +\infty$, and the transformed initial form in $L^{p_0/(p_0-3+k)}_{\rm loc}$. Reflect $\Psi$ in time and smooth first in space at scale $h$. With $\omega$ the uniform time modulus of $\Psi$ on the compact neighborhood, use time scale $\ell(h)$ with $h^{-1}\omega(\ell(h))\to0$. This bounds the additional error in tangential derivatives at zero. The bulk convergence is still in $W^{1,p}$, and the initial tangential derivatives converge in $L^{p_0}$ to $D_X\Psi_0$. The initial pullback of degree $3-k$ therefore converges in $L^{p_0/(3-k)}$. H\"older's inequality and the spatial area formula give the prescribed Eulerian initial term.

%-----------------------------------------------------------------------------------------------------------------------------------------

\subsection{Material parametrizations through cavitation}
\label{sec:material-flow-at-cavitation}

\paragraph{The particle weight in material coordinates.}
In a bi-Lipschitz chart with $\chi=f$ in normalized material volume coordinates, we define the Eulerian coordinate velocity $\beta(t,x)\coloneqq v(t,f(t,x))$. The cubic-minor formula and the orthogonality condition~\eqref{eq:450bis-2} give
\be
|g|^{1/2}{\particlecurrent}^0=\det D_xf=a,\qquad
{\particlecurrent}^i={\particlecurrent}^0\beta^i,\qquad
\beta^i=u^i/u^0
\ee
on the timelike part. Formula~\eqref{eq:material-relative-flux} gives
\be
\widetilde{(|g|^{1/2}{\particlecurrent})}^0=1,\qquad
\widetilde{(|g|^{1/2}{\particlecurrent})}^A=0.
\ee
For independent fields the weight is $j(|g|^{1/2}{\particlecurrent}^0)\circ\Psi$. With exact minors and $\chi=f$, the Jacobian $a$ vanishes at zero current. Its reciprocal diverges when $a$ tends to zero through positive values. In the collapsed slab of Proposition~\ref{prop:collapsed-material-coordinates}, $a$ jumps directly from one to zero. A nonzero null current has $n=0<{\particlecurrent}^0$ and can admit nonsingular material coordinates on a slice.

The initial particle positions give an invertible material parametrization through a cavity. Their particle weight can vanish, and the construction uses the extended velocity specified below.

\begin{proposition}[Material flow coordinates through a cavity]
\label{cor:weak-flow-coordinates-vacuum}
On a coordinate slab $[0,T]\times\mathcal D$, with $\mathcal D=\mathbb T^3$ or $\mathbb R^3$, take smooth prescribed metric and spacelike slices up to $t=0$, with $t$ increasing toward the future. Let $\beta$ be a coordinate velocity with
\be
 \beta\in L^\infty((0,T)\times\mathcal D),\qquad
 \nabla_x\beta\in L^1(0,T;L^\infty(\mathcal D)),
\ee
and let $y$ be its spatial flow starting from the identity,
\be
\del_ty(t,X)=\beta(t,y(t,X)),\qquad y(0,X)=X.
\ee
The evolution equation holds for almost every time. Let $\particlecurrent$ be a future causal current, locally integrable up to the initial slice, allowing zero and nonzero null currents. Suppose it satisfies ${\particlecurrent}^i=\beta^i {\particlecurrent}^0$ almost everywhere and particle conservation with prescribed locally integrable initial data set. Then  the following conclusions hold.
\begin{itemize}
\item[(i)] The map $\Psi(t,X)\coloneqq(t,y(t,X))$ is bi-Lipschitz on compact flow slabs, and
\be
\label{eq:vacuum-flow-particle-weight}
\widetilde{(|g|^{1/2}{\particlecurrent})}^0(t,X)=\omega_0(X),\qquad
\widetilde{(|g|^{1/2}{\particlecurrent})}^A=0,\qquad
\omega_0(X)\coloneqq |g(0,X)|^{1/2}{\particlecurrent}_{\rm in}^0(X)\ge0.
\ee
\item[(ii)] The weight $\omega_0$ may vanish on an open set. The zero set of the current at almost every time is the image under $y(t,\cdot)$ of the zero set of $\omega_0$, up to sets of zero volume. A zero-current region may still carry constitutive stress.
\item[(iii)] If $\omega_0\in L^p(\mathcal D)$ for $1\le p\le\infty$, the densitized density $r\coloneqq|g|^{1/2}{\particlecurrent}^0$ satisfies
\be
 e^{-(1-1/p)A(t)}\|\omega_0\|_{L^p}
 \le\|r(t)\|_{L^p}\le
 e^{(1-1/p)A(t)}\|\omega_0\|_{L^p},\quad
 A(t)\coloneqq\int_0^t\|\operatorname{div}\beta(s)\|_\infty\,ds,
\ee
for almost every $t$, with $1/p=0$ when $p=\infty$. For $p< +\infty$, the density $r$ has a representative in $C([0,T];L^p)$ with $r(0)=\omega_0$, and the bounds hold at every time for this representative. For $p=\infty$, it has a weakly-$*$ continuous representative with the same initial value and bounds.  Theorem~\ref{un-prop:weak-coordinate} transfers the balance laws, closure equations and mathematical entropy inequality, with their initial terms under its integrability hypotheses.
\end{itemize}
\end{proposition}

For $\particlecurrent=nu$ with future unit $u$, take $\beta^i=u^i/u^0$. The hypotheses follow if $u$ is bounded, $u^0$ is bounded away from zero and $D_xu\in L_t^1L_x^\infty$. They concern the extension of $\beta$ through zero current, without prescribing a time derivative or an initial value of $u$. They also allow a null current, for instance $\particlecurrent=b(x^1-t)(\partial_t+\partial_1)$ in Minkowski coordinates with smooth nonnegative periodic $b$. Then $\beta=e_1$, $y(t,X)=X+te_1$, with vacuum on intervals where $b=0$.

\paragraph{{Proof of the material flow formula.}}
{We use the flow Jacobian to recover the initial density on each trajectory.}
\begin{proof}
For the spatially Lipschitz representative of $\beta$, Gronwall bounds the Lipschitz constants of the forward and inverse flows by $\exp(\int_0^t\|D_x\beta(s)\|_{\infty,\mathrm{op}}\,ds)$. Bounded velocity gives time Lipschitz continuity of both maps. We claim that the spatial Jacobian of the forward flow is
\be
 j(t,X)=\det D_X y(t,X)
 =\exp\left(\int_0^t\operatorname{div}_x\beta(\tau,y(\tau,X))\,d\tau\right)
\ee
for each $t$ and almost every $X$. To see this, convolve $\beta$ in space and denote the resulting flow by $y_\delta$. The common Lipschitz bounds and the integral flow equation give uniform convergence of the flows and their inverses. On compact flow slabs the smoothed divergences composed with $y_\delta$ converge in $L^1_{t,X}$ to $(\operatorname{div}\beta)\circ y$: approximate the divergence in $L^1$ by continuous functions and use the common change-of-variables bounds. The time integrals therefore converge in $C_tL^1_{\rm loc}$ and are bounded by $A(T)$. Their exponentials converge in the same space. Testing the smoothed Jacobians against $b\in C_c^\infty$ gives $\int b(y_\delta(t,\cdot)^{-1}(x))\,dx$. Uniform convergence of the inverse maps and bounded displacement give the limit $\int b(y(t,\cdot)^{-1}(x))\,dx$ on each slice. Change of variables by the orientation-preserving map $y(t,\cdot)$ identifies its Jacobian with the exponential above.

In Theorem~\ref{un-prop:weak-coordinate} take $\chi=y^{-1}$ on each slice. The identity ${\particlecurrent}^i={\particlecurrent}^0\beta^i$ makes the relative spatial flux zero. The transformed density is thus constant in time with prescribed value $\omega_0$, proving~\eqref{eq:vacuum-flow-particle-weight}. Future causality gives ${\particlecurrent}^0\ge0$ and $\particlecurrent=0$ on exactly its zero set. Since $j>0$, the density identity transports that set. For $1\le p< +\infty$, change of variables yields
\be
 \|r(t)\|_{L^p}^p=\int_{\mathcal D}|\omega_0(X)|^p j(t,X)^{1-p}\,dX.
\ee
The Jacobian formula gives $e^{-A(t)}\le j(t,\cdot)\le e^{A(t)}$, yielding the displayed bounds. The case $p=\infty$ follows from $r(t,y(t,X))=\omega_0(X)/j(t,X)$ and the bi-Lipschitz change of variables. The exponential formula makes $j$ continuous in time in $L^\infty$, with $j$ bounded away from zero. Composition by the inverse flow is strongly continuous on $L^p$ for $p< +\infty$: this follows first for continuous compactly supported functions from uniform convergence of the maps, and then by density and their uniform change-of-variables bounds. Thus $r(t)=(\omega_0/j(t))\circ y(t)^{-1}$ is continuous in $L^p$, including at $t=0$, which extends the norm bounds to every time. When $p=\infty$, the formula $\int r(t,x)\phi(x)\,dx=\int\omega_0(X)\phi(y(t,X))\,dX$ gives continuity against continuous compactly supported tests; their density in $L^1$ and the uniform $L^\infty$ bound give weak-$*$ continuity. The pointwise bounds on $j$ give both essential-supremum bounds at every time. The equivalence of the other balance laws, closure equations and the mathematical entropy inequality then follows, including the initial terms.
\end{proof}
If the prescribed initial particle data is initially given as a Radon measure $\mu_0$, the same argument forces it to have a locally integrable density. Indeed, the change of coordinates is the identity at $t=0$, so the transformed current has zero spatial components and prescribed initial term $\mu_0$. Its time component has slices $\mu_0$ for almost every $t$, by the distributional equation $\del_t\widetilde{(|g|^{1/2}{\particlecurrent})}^0=0$ with its initial term. The transformed bulk current is locally integrable, hence these slices are absolutely continuous.

In the inverse flow chart $\chi$, the first identity in~\eqref{eq:vacuum-flow-particle-weight} reads $|g|^{1/2}{\particlecurrent}^0=(\omega_0\circ\chi)\det D_x\chi$. The extended velocity determines the motion of the corresponding trajectories through the cavity. The strong solutions constructed below admit this flow. For weak solutions the spatial Lipschitz bound is an additional assumption.

\paragraph{Tests with independent traces across a collapsed cavity.}
\label{sec:material-flux-measures}
Wagner~\cite{Wagner0} uses discontinuous Lagrangian test-functions at vacuum in one-dimensional gas dynamics. Here we use independent traces on the two faces of a cavity, together with the equations in the sense of distributions inside it. For the planar collapse $X^1=\xi(x^1)\coloneqq\operatorname{sgn}(x^1)(|x^1|-1)_+$, with $X^2=x^2$ and $X^3=x^3$, the interval $-1<x^1<1$ is sent to the material plane $X^1=0$. The pairings in the weak formulation with these exterior tests, supplemented by the equations in the open cavity, determine the Eulerian formulation. We first treat the current when it is conserved in the cavity with its prescribed initial term, as happens for exact Lipschitz configurations.

Let the bulk fluxes, sources and forms be locally integrable up to $t=0$. Take locally integrable initial tangential forms. For scalar balance laws and the mathematical entropy inequality, allow the prescribed initial measure to be a locally finite measure. In the formulation with two traces, the interior laws use the restrictions of the same prescribed data to the initial cross-section of the cavity. A material test is smooth and compactly supported up to the plane on each half-space and allowed to meet $t=0$. Denote its two traces by $\theta_-$ and $\theta_+$. We use $\theta^{\rm ext}$ for the Eulerian test obtained by translating the two exterior pieces and extending between the faces by
\be
\theta^{\rm ext}(t,x)=\frac{1-x^1}{2}\,\theta_-(t,x^2,x^3)
                     +\frac{1+x^1}{2}\,\theta_+(t,x^2,x^3).
\ee
For test-forms interpolate the coefficients in the coordinate coframe. The lifted tests are Lipschitz, and the scalar lift preserves nonnegativity.  In the equations with two traces we use the pairings in the exterior weak formulation and add the interior pairings with this lift. Thus the formulation uses the fields in the cavity along with the exterior material fields. The bulk terms inside the slab are integrals in $x^1$ against the affine weights or their derivatives. A measure in the initial data set is paired with the continuous function $\theta^{\rm ext}(0,\cdot)$ on the whole initial slice, including both faces. The initial term for closure of a $k$-form is
\be
 (-1)^k\int_{\Sigma_0}\alpha_{\rm in}\wedge\iota_0^*\theta^{\rm ext},
\ee
where $\iota_0$ is the inclusion of the initial slice and $\Sigma_0$ has the spatial coordinate orientation. Thus the two traces at $t=0$ also pair with the initial data inside the cavity.

Reflection and mollification justify the lift as an Eulerian test by uniform convergence and bounded derivatives converging almost everywhere. Conversely, we subtract from a smooth Eulerian test the lift of its exterior pieces. The difference is Lipschitz, vanishes on both faces and is supported in the closed cavity. We compress it slightly in $x^1$, reflect in time and mollify within the spatial margin. The interior tests converge uniformly, with bounded first derivatives converging in each finite $L^p$ by continuity of dilations and convolution. Uniform convergence gives the initial pairing with measures too. Since the difference vanishes on both faces, the lift accounts for any initial face masses. An almost-everywhere subsequence gives the bulk pairing. Adding the interior identity to the equation with separate traces gives the Eulerian equation. For entropy the difference of the scalar tests is signed and its pairing is zero by interior current conservation. This argument also applies to the difference of any two Lipschitz lifts of the same exterior test, so the pairings are independent of the lift under the interior laws.

For an independent current, let $\mathcal I$ be the pairing defining the mathematical entropy inequality including the prescribed initial term. Suppose its restriction to the cavity, including the interior of its initial cross-section, is $-\mu_{\rm c}$, with $\mu_{\rm c}\ge0$ locally finite up to the faces. Extend this measure by zero on the faces. The two-trace entropy inequality is
\be
\label{eq:two-trace-interior-entropy-production}
 \mathcal I(\theta^{\rm ext})+\int\theta^{\rm ext}\,d\mu_{\rm c}\le0
\ee
for nonnegative tests. The Eulerian inequality implies it by removing the interior part of the production measure. Conversely, for a nonnegative Eulerian test $\varphi$ use the lift $\theta^{\rm ext}$ of its exterior pieces. The approximation just proved gives $\mathcal I(\varphi-\theta^{\rm ext})=-\int(\varphi-\theta^{\rm ext})\,d\mu_{\rm c}$, by uniform convergence for the measure term. Adding the pairings yields $\mathcal I(\varphi)\le-\int\varphi\,d\mu_{\rm c}\le0$. The difference of two lifts of the same exterior test satisfies this identity. Hence the corrected pairing in~\eqref{eq:two-trace-interior-entropy-production} is independent of the lift, including when $\mu_{\rm c}$ has an interior initial atom.

\paragraph{The interior laws and curved cavities.}
The interior stress and continuity laws hold if $\particlecurrent=T=0$ inside and their prescribed normal fluxes vanish on the initial cross-section. For a Lipschitz configuration up to $t=0$, the tangential minors of its trace give the initial terms for closure and current conservation at every rank (cf.~Proposition~\ref{prop:kinematic-weak-local} and the proof of Proposition~\ref{prop:collapsed-material-coordinates}).

The same equivalence holds for $K\circ\Phi_{\rm flat}\circ H^{-1}$, with $\Phi_{\rm flat}(t,s,z)=(t,\xi(s),z)$ and $z\in\R^2$. Here $H$ and $K$ preserve the time coordinate and spatial orientation, are bi-Lipschitz up to $t=0$, and are defined on the full reference strip over the test support. Apply Theorem~\ref{un-prop:weak-coordinate} to the fluxes, forms and initial data before taking the affine lift in $(t,s,z)$.

For moving Lipschitz faces $x^1=a_\pm(t,z)$ with separation bounded below, take $H$ affine in $s$ between them and a translation in $s$ outside. Its normal derivative inside is $(a_+-a_-)/2$, giving the inverse Lipschitz bound. Then $\Phi_{\rm flat}\circ H^{-1}$ collapses the segments at fixed $(t,z)$, and the two traces describe their endpoints.

{We now consider passage through a pinch.} Let $S$ be a closed set off which the local equivalence holds, and write $S_\epsilon$ for its spacetime $\epsilon$-neighborhood. Assume that the total variations of the bulk sources and of the prescribed initial measures give no mass to $S$. For integrable initial densities and tangential forms it suffices that $S$ meet the initial slice in a set of zero spatial volume. Suppose that each vector of flux or form coefficients $B$ satisfies
\be
 \label{eq:pinch-flux-decay}
 \lim_{\epsilon\downarrow0}\frac1\epsilon
       \int_{K\cap S_{3\epsilon}}|B|\,dt\,dx=0
 \qquad\text{for each compact subset $K$ of the coordinate slab}.
\ee
Then the Eulerian equations and entropy inequality, with their prescribed initial terms, hold across $S$.

We multiply each test by a smooth cutoff that vanishes on $S_\epsilon$, equals one off $S_{3\epsilon}$, takes values in $[0,1]$, and has first derivatives bounded by $C/\epsilon$. A partition of unity gives the distributional equations and entropy inequality for these tests. The terms with derivatives of the cutoff are bounded by $C/\epsilon$ times the integral in~\eqref{eq:pinch-flux-decay}. This condition also implies that the bulk $L^1$ coefficients give no mass to $S$. The source and initial pairings converge by dominated convergence and the cutoff preserves the sign of the entropy test. Once the entropy inequality has been recovered, use the complementary cutoff times a nonnegative test in its production measure. The cutoff equals one on $S$ after taking the complement, and the same bounds show that the measure gives no mass to $S$, including on the initial slice.

For example, suppose $|K\cap S_\epsilon|=\Ocal(\epsilon^d)$ with $d>0$ and the bulk sources are locally integrable. For $B\in L^p_{\rm loc}$ with $1<p< +\infty$, H\"older's inequality gives~\eqref{eq:pinch-flux-decay} when $d(1-1/p)\ge1$. At equality the $L^p$ norm over $K\cap S_{3\epsilon}$ tends to zero. For $p=\infty$ it suffices that $d>1$. Thus a smooth spacetime set of codimension two that meets the initial slice in a set of zero spatial volume is removable for locally square-integrable fluxes.

For exact Lipschitz configurations with the prescribed tangential minors at $t=0$, closure and particle conservation hold across $S$ by the polynomial minor identities. Only the stress fluxes then need~\eqref{eq:pinch-flux-decay}. Smooth metric coefficients are bounded on $K$, so it suffices to impose that condition on $T$. If $|K\cap S_\epsilon|=\Ocal(\epsilon^d)$ and $|T|\le C_K\operatorname{dist}(\cdot,S)^\alpha$ on $K$, it holds for $d+\alpha>1$. This includes rank two on a hypersurface where the stress vanishes at a positive H\"older rate.

%-----------------------------------------------------------------------------------------------------------------------------------------

\paragraph{\texorpdfstring{Rank loss and zero-density limits.}{Rank loss and zero-density limits.}}
\label{Append-BB}

{
Proposition~\ref{prop:boundary-fields-and-velocity} distinguishes the ranks and causal types for independent fields and exact minors. The same possibilities occur in the distributional formulation, with local integrability of the current, forms and constitutive stress as required in Definition~\ref{def:weak-local-complete}.
}

For a barotropic fluid in matter, write $\particlecurrent=nu$ and $\varepsilon=\rrestnew(n)/n$. The pressure is $p=n^2\varepsilon'=n\rrestnew'-\rrestnew$, and the stress has the equivalent forms
\be
\label{eq:main-current-fluid-stress}
T_{\mu\nu}=(\rrestnew+p)u_\mu u_\nu+p\,g_{\mu\nu}
=\frac{\rrestnew+p}{n^2}{\particlecurrent}_\mu {\particlecurrent}_\nu+p\,g_{\mu\nu}.
\ee
For $\rrestnew(n)=n^2/2$ the Hessian of $\raug(\particlecurrent;g)$ in $\particlecurrent$ is $-g$ even at $\particlecurrent=0$, where the derivative of the quadratic stress vanishes. With positive rest mass, the sound speed tends to zero for the polytropic laws with $\gamma>1$ of Theorem~\ref{uf:thm:main}.

%============================================================================

\subsection{Shock interfaces and entropy admissibility}
\label{sec:weak-shocks}

\paragraph{Jump conditions and the initial slice.}

The distributional equations determine the normal fluxes that can jump across a front. We first work in the interior of the spacetime slab.

\begin{proposition}[Jump conditions for a piecewise regular solution]
\label{def:shock-interface}
On a spacetime with a smooth metric, let $\Sigma=\{\varphi=0\}$ be a $C^1$ hypersurface, with a $C^1$ defining function $\varphi$ and nonzero conormal $\nu=d\varphi$. Orient the plus side by $\varphi>0$ and set $[A]\coloneqq A_+-A_-$.  {Assume that the augmented fields and their constitutive stress are piecewise $C^1$ up to $\Sigma$ and satisfy the balance laws off $\Sigma$. They form an entropy-admissible augmented weak solution on the whole neighborhood if and only if the ordinary divergence satisfies $-\nabla_\alpha\particlecurrent^\alpha\le0$ off $\Sigma$ and the traces satisfy}
\be
\label{eq:augmented-rankine-hugoniot}
\nu\wedge[F^A]=0,\qquad
\nu\wedge[G_A]=0,\qquad
\nu_\beta[T^{\alpha\beta}]=0,\qquad
-\nu_\alpha[{\particlecurrent}^\alpha]\le0.
\ee

\end{proposition}

For a reduced formulation, the jump conditions use the selected form families.
For a front with nonzero spatial conormal, rescale $\nu$ to write it as $-s\,dt+\xi$, with $\xi$ of unit spatial length. For a fixed set of left-hand traces and $\xi$, the first three relations in~\eqref{eq:augmented-rankine-hugoniot} are algebraic conditions on the right-hand traces and the normal speed $s$. A Hugoniot curve parametrizes a local branch of these relations issuing from the left-hand traces. The constitutive graph restricts the right-hand traces, and the current inequality tests entropy admissibility along the branch.

\begin{proof}
The implicit function theorem writes the $C^1$ interface locally as a graph. Integration by parts on its two sides gives $\del_\mu\mathbf1_{\{\varphi>0\}}=\delta(\varphi)\del_\mu\varphi$, where $\delta(\varphi)$ is the positive level-set measure normalized by the coarea formula in the coordinate chart. The smooth connection coefficients contribute no singular part. The singular parts are therefore given by the following formulas.
\be
\aligned
&(dF^A)_{\rm sing}=\delta(\varphi)\nu\wedge[F^A],\qquad
(dG_A)_{\rm sing}=\delta(\varphi)\nu\wedge[G_A],\\
&(\nabla_\beta T^{\alpha\beta})_{\rm sing}
=\delta(\varphi)\nu_\beta[T^{\alpha\beta}],\qquad
(-\nabla_\alpha {\particlecurrent}^\alpha)_{\rm sing}
=-\delta(\varphi)\nu_\alpha[{\particlecurrent}^\alpha].
\endaligned
\ee
{The ordinary divergence off $\Sigma$ gives the part of the entropy inequality absolutely continuous with respect to spacetime volume. The singular part is nonpositive precisely when the last inequality in~\eqref{eq:augmented-rankine-hugoniot} holds. The other three conditions cancel the singular parts of the balance laws.}\qedhere
\end{proof}

With respect to spacetime volume, the entropy production measure is $\nabla_\mu {\particlecurrent}^\mu\,dV_g$. Its interface part has density $\nu_\mu[{\particlecurrent}^\mu]$ relative to the positive measure $\delta(\varphi)\,dV_g$. Reversing the conormal and interchanging the sides leaves this density unchanged. Replacing $\varphi$ by another defining function with the same orientation rescales the density and the level-set measure by reciprocal factors. The formula also applies to null interfaces.
Suppose the two fields extend with bounded traces to the initial slice and the front meets that slice transversely. Integration by parts gives the interior interface terms above and the initial data on each initial half-domain. Their common edge has zero slice measure. Thus the prescribed tangential forms and normal stress flux must equal the corresponding traces almost everywhere. For the current, let ${\particlecurrent}_{\rm in}^0$ denote the prescribed time component in adapted coordinates. The initial term of the entropy pairing is the integral of $|g|^{1/2}({\particlecurrent}_{\rm in}^0-{\particlecurrent}^0(0+))$ against the initial test. The entropy inequality allows ${\particlecurrent}_{\rm in}^0\le {\particlecurrent}^0(0+)$, and current conservation requires equality. In the examples below the prescribed data are the actual traces.

\paragraph{Shocks of a continuous configuration.}
The closure condition fixes the common tangential derivative on the front. For a configuration given by two potentials it remains to match their values there. If $f$ is continuous across $\Sigma$ and has bounded $C^1$ derivatives up to each side, it is locally Lipschitz. Equality of tangential derivatives of the identical trace on both sides gives the Hadamard jump relation
\be
[F^A]=a^A\nu.
\ee
Its quadratic and cubic minor jumps are consequently determined by $F_-$, the material vector $a$ and $\nu$. In particular,
\be
\label{eq:exact-shock-current}
\nu_\alpha[\particlecurrent(df)^\alpha]=0.
\ee
In the expansion of the cubic current at $F_+=F_-+a\otimes\nu$, every increment contracted with $\nu$ contains two identical conormals in an alternating product and vanishes. Proposition~\ref{prop:kinematic-weak-local} also gives this conclusion without a rank assumption. Thus the current inequality holds as an equality at every such interface, and shock selection in this class requires a further admissibility criterion. A jump in the value of $f$ would instead give a singular measure in $df$. To describe thermodynamic entropy production one would also need a constitutive law with thermal variables.

At an interface with an empty cavity, the cavity traces are $\particlecurrent=T=0$. The exact-current jump identity and the stress jump condition therefore reduce on the matter side to $\nu_\alpha {\particlecurrent}^\alpha=0$ and $\nu_\beta T^{\alpha\beta}=0$. If the matter trace of $\particlecurrent$ is timelike, the first equality says that the interface moves with the normal component of the particle velocity and is timelike, since its tangent space contains that timelike current. If the current tends to zero on the matter side as well, that equality holds without determining the speed of the interface.

For a piecewise regular configuration {whose stress traces are the variational stresses of a $C^1$ physical energy}, write $p_A^\pm=P_A^{\nu,\pm}\nu_\nu$ for the material traction, with $P_A^\nu=\partial\rphys/\partial F^A_\nu$. Let $\langle B\rangle=(B_++B_-)/2$ denote the mean of two traces. The jump relation for $F$ and the stress identity used in Proposition~\ref{prop:variational-EL-local04} below give
\be
 \nu_\nu[T_\mu{}^\nu]
 =\langle F_\mu^A\rangle[p_A]
   +\nu_\mu\bigl(a^A\langle p_A\rangle-[\rphys]\bigr).
\ee
Here the only singular term in the material Euler--Lagrange expression is the traction jump. The pullback of the material volume form to $T\Sigma$ is nonzero exactly when $F|_{T\Sigma}$ has rank three. Since this pullback is also the restriction of the particle three-form to $T\Sigma$, that rank condition is equivalent to $\nu_\alpha {\particlecurrent}^\alpha\ne0$. In this case the tangential stress condition gives $[p_A]=0$. The remaining stress condition is the scalar equality $[\rphys]=a^A\langle p_A\rangle$. Conversely, these two relations imply the stress jump condition. Thus the variational traction condition alone does not give the complete stress jump condition. The decomposition applies also to a null hypersurface. At any rank, the tangential stress condition says that $[p]$ annihilates $F(T\Sigma)$, which is common to the two sides. To recover the remaining stress condition, choose any transverse vector $b$ with $\nu(b)=1$ and impose
\be
 \langle F^A\rangle(b)[p_A]+a^A\langle p_A\rangle-[\rphys]=0 ,
\ee
{which is independent of that choice once the tangential condition holds.}

There is a simple sufficient condition for traction continuity to give the complete stress condition. Assume that the segment $F_-+\theta a\otimes\nu$, $0\le\theta\le1$, lies in the domain of the energy at the common trace of $f$ and at the fixed metric. If the restriction $R(\theta)=\rphys(F_-+\theta a\otimes\nu;g)$ is convex, then $[p]=0$ implies $R'(0)=R'(1)=a^A p_A^-$. Monotonicity of $R'$ gives $[\rphys]=R(1)-R(0)=a^A\langle p_A\rangle$. If $R$ is strictly convex for $a\ne0$, a nontrivial gradient jump cannot satisfy traction continuity along this segment.  Under this strict convexity along the segment, stress balance therefore excludes a nontrivial exact gradient jump at an interface with nonzero normal particle flux.   If one side has $\particlecurrent=0$, the normal particle flux on the other side is zero by~\eqref{eq:exact-shock-current}. The stress then constrains the tangential projection of the traction jump and the scalar combination displayed above.

%-----------------------------------------------------------------------------------------------------------------------------------------

\section{A weak--strong stability theorem for relativistic elastodynamics}
\label{section=6}

%-----------------------------------------------------------------------------------------------------------------------------------------

\subsection{Strong solutions with cavitation}
\label{subsec:entropy-conditions-vacuum}
\label{subsec:entropy-strong-solutions}
\label{subsec:completed-polyconvexity}

\paragraph{Physical fluxes and reference parameters.}

We specify the constitutive parameters and support of the mathematical entropy used for comparison at cavitation.
Let $Z\in{\Cl(\mathcal Z_{\raug})}$ be a point of the completed constitutive graph and $Y=\pi(Z)$ the fields selected from it. When the model  includes both form families and the current, $Y=(F,G,\particlecurrent)$. The comparison set is the full completed constitutive domain unless additional relations are explicitly imposed. {For independent augmented fields, we take $\Theta\equiv0$ and $\mathcal Z_{\rm ph}(t,x)=\Cl(\mathcal Z_{\raug})$; {the fields satisfy the constitutive constraints without any additional requirement that they be exact minors}. When the comparison is restricted to the exact spacetime minors, we instead take}
\be
\label{eq:general-entropy-algebraic-constraints}
\mathcal Z_{\rm ph}(t,x)=\{Z\in{\Cl(\mathcal Z_{\raug})}:\Theta(t,x,\pi(Z))=0\},\qquad
\Theta(Y)=(G-G(F),\particlecurrent-\particlecurrent(F)).
\ee
The polynomial minors use the prescribed metric and fixed material volume coordinates. This restriction gives a particular comparison class within the augmented theory. The map $\Theta$ is $C^1$ near the prescribed reference range, with geometric dependence fixed by the model.

On a periodic coordinate slab with spacelike slices and time increasing toward the future, let $\mathscr F_a^\mu$ be the evolution fluxes of the stress and form equations, $\mathscr S_a$ their sources, and $\mathscr K_\ell^i$ the fluxes of the spatial involutions. Here $a=1,\ldots,d_{\rm sel}$, with $d_{\rm sel}=4+3r_F+3r_G$, indexes the evolution rows, and $\ell=1,\ldots,3r_F+r_G$ indexes the spatial constraints. A missing form family contributes no rows or constraints. We use densitized stress currents and include connection terms in $\mathscr S_a$. The mathematical entropy current is $\mathscr I^\mu=-|g|^{1/2}{\particlecurrent}^\mu$, with entropy source $\mathscr G=0$. We write $U=(\mathscr F_a^0)_a$. A reference parameter $\sigma_*$ consists of a physical constitutive state $Z_*$ and any boundary data needed to select a finite supporting covector for the entropy and the continued material momentum. For the laws with a positive rest-mass term considered above these are a finite-velocity constitutive branch and a future unit velocity orthogonal to the forms, also at ${\particlecurrent}_*=0$. The admissible parameter domain includes its finite-velocity zero-current limits. These are pointwise constitutive parameters, before any reference solution is chosen.

{Let $Z_*$ be the fields and stress selected by $\sigma_*$, and write $U_*=U(Z_*)$. For a finite covector $\lambda_*$ we define the relative density directly from the two sets of fluxes by
\be
\label{eq:general-entropy-first-remainder}
\mathscr Q^0(Z\mid Z_*)\coloneqq
\mathscr I^0(Z)-\mathscr I^0(Z_*)-\lambda_*\cdot\bigl(U(Z)-U_*\bigr).
\ee
We call $\lambda_*$ a supporting covector at $Z_*$ on a set of constitutive values if $\mathscr Q^0(Z\mid Z_*)\ge0$ for every $Z$ in that set. When the law provides a convex function $\eta$ of $U$ with $\mathscr I^0=\eta(U)$ on the set, this is the supporting inequality for $\eta$ restricted to the conservative image of that set, and $\mathscr Q^0(Z\mid Z_*)=D_\eta(U\mid U_*)$. For fully densitized variables $V_{\rm dens}=|g|^{1/2}V$, the entropy $|g|^{1/2}\eta(V_{\rm dens}/|g|^{1/2})$ has the same numerical subgradient as $\eta(V)$. When only the stress rows are densitized, the coefficients of the form rows acquire the factor $|g|^{1/2}$.}

{At a fixed spacetime point, let $Z$ denote the comparison fields and their stress, and let $Z_*$ denote the corresponding reference fields and stress. We write $Y_*=\pi(Z_*)$, $\Delta Y=Y-Y_*$, and, for each flux or source $\mathscr H$, $\Delta\mathscr H=\mathscr H(t,x;Z)-\mathscr H(t,x;Z_*)$.} We choose measurable constitutive linear maps $A_{a*}^\mu,B_{a*},B_{\eta*},M_{\ell*}^i$ on $\Delta Y$, with $B_{\eta*}=0$. They depend on $\sigma_*$ and the prescribed geometry, independently of reference derivatives and the comparison state, and are bounded on allowed bounded reference ranges, including the constitutive boundary. They agree with the constitutive differentials where $n>0$; the constitutive law specifies their finite affine continuations on the constitutive boundary. {We define the corresponding differences after subtraction of these linear terms by}
\be
\label{eq:general-entropy-constitutive-remainders}
\aligned
&r_{\mathscr F,a}^\mu=\Delta\mathscr F_a^\mu-A_{a*}^\mu\Delta Y,\qquad r_{\mathscr S,a}=\Delta\mathscr S_a-B_{a*}\Delta Y,\\
&r_{\mathscr G}=\Delta\mathscr G-B_{\eta*}\Delta Y,\qquad r_{\mathscr K,\ell}^i=\Delta\mathscr K_\ell^i-M_{\ell*}^i\Delta Y,\\
&r_\Theta=\Theta(t,x,Y)-\Theta(t,x,Y_*)-D_Y\Theta(t,x,Y_*)\Delta Y.
\endaligned
\ee {Thus $r_{\mathscr F}$ and $r_{\mathscr S}$ are the remainders of the evolution fluxes and sources, while $r_{\mathscr K}$ and $r_\Theta$ concern the differential and algebraic constraints. In the present formulation $\mathscr G=0$ and $B_{\eta*}=0$, so $r_{\mathscr G}=0$. These definitions do not by themselves assert a bound on the remainders at the constitutive boundary.}

\paragraph{{Constitutive compatibility with mathematical entropy.}}
 
Fix the metric and the observer used for the conservative variables at a point of spacetime. With the same volume convention for the fluxes, $U(Z)$ and $\mathscr I^0(Z)$ are functions of the constitutive fields and their stress. The following property concerns these functions on the chosen constitutive domain, including whichever finite limiting stresses are admitted at its boundary.

\begin{definition}
\label{def:entropy-extension-condition}
{A law is \underline{\emph{compatible with a convex mathematical entropy}} for the chosen observer on a specified set of constitutive values if there is a proper lower semicontinuous convex function $\eta$ of the conservative variables such that, for every $Z$ in this set, $\eta(U(Z))$ is finite and
\be
\label{eq:entropy-extension-identity}
 \eta(U(Z))=\mathscr I^0(Z).
\ee}
\end{definition}

{The constitutive property applies pointwise with the metric and observer at $(t,x)$, with a measurable choice of $\eta$. It gives the supporting inequality both in the slab and for the prescribed initial data.}

{For the local constitutive laws the form multipliers vanish at zero density (Proposition~\ref{uf:prop:entropy-multipliers}). Thus two sets of fields with $\particlecurrent=T=0$ have zero relative entropy even if their forms differ. Proposition~\ref{uf:prop:local-convex-entropy} shows that zero relative entropy determines the current and stress everywhere and the forms in matter.}

\paragraph{Strong solutions and supporting covectors.}

{The constitutive maps and their affine continuations are fixed before choosing the strong solution used for comparison.}

\begin{definition}
\label{def:entropy-strong-solution}
{In the selected formulation, a \underline{\emph{strong solution}} on the chosen comparison set consists of bounded continuous fields $(F_*,G_*,\particlecurrent_*)$, their associated continuous constitutive stress $T_*$, and a bounded continuous choice of constitutive parameters $\sigma_*$. Write $Z_*=(F_*,G_*,\particlecurrent_*,T_*)$ for the fields and stress, and $U_*=U(Z_*)$. They satisfy the following conditions:
\begin{itemize}
\item[(i)] The fields and stress take values in $\Cl(\mathcal Z_{\raug})$. They satisfy the evolution equations, spatial compatibility equations and mathematical entropy equality in the sense of distributions, with the initial data set given by their values at $t=0$.
\item[(ii)] At almost every point of the slab, $\lambda_*=\lambda(\sigma_*)$ is a finite supporting covector at $Z_*$ on the chosen comparison set, in the sense of~\eqref{eq:general-entropy-first-remainder}. The covector $\lambda_*$ and the involution coefficients $b_*^\ell=b^\ell(\sigma_*)$ are bounded and Lipschitz on the closed slab. Along the admissible reference branch these coefficients and the constitutive affine maps have finite limits and are bounded on bounded parameter ranges.
\item[(iii)] {The supporting inequality~\eqref{eq:general-entropy-first-remainder}, evaluated at the initial reference values and with the trace of $\lambda_*$, holds for the prescribed comparison data.}
\end{itemize}}
\end{definition}
{The parameters $\sigma_*$ may specify the continued velocity at zero current or a limiting constitutive stress. The supporting covector and regularity conditions apply also at nonzero lightlike current. At zero current the remainder bound may restrict the chosen reference parameters even though the definition imposes no velocity evolution there.}

For the laws of Proposition~\ref{prop:vacuum-stress-subgradient}, a timelike velocity orthogonal to the forms gives a finite support on the branch with no rank-one null stress, including at zero current with residual elastic stress. Positive rank-one null stress in that proposition has no finite supporting covector and can enter the comparison as a value of the weak solution. The regularity $C_tH^s\cap C_t^1H^{s-1}$, $s\ge3$, gives Lipschitz bounds for multipliers depending smoothly on the reference variables up to the boundary. Proposition~\ref{uf:prop:entropy-multipliers} below { verifies these bounds for the laws used in the local existence theorem}, while Theorem~\ref{og:theorem} below estimates  the derivatives of the Ogden-type energies directly. The criteria below use the adjoint identity or the continued material equations.

%-----------------------------------------------------------------------------------------------------------------------------------------

\subsection{Weak--strong stability with cavitation}
\label{subsec:physical-weak-strong-vacuum}

\paragraph{{The relative entropy current.}}
The relative current pairs the differences of the balance laws with the supporting covector and incorporates the spatial compatibility equations in the chosen foliation. Write $\Delta$ for the difference between the comparison fields and their stress $Z$ and the reference fields and their stress at the same spacetime point. The time component of the \emph{relative entropy current} is defined by~\eqref{eq:general-entropy-first-remainder}. Its spatial components in coordinate densities are
\be
\label{eq:general-entropy-relative-current}
\mathscr Q^i=\Delta\mathscr I^i-\lambda_*^a\Delta\mathscr F_a^i+b_*^\ell\Delta\mathscr K_\ell^i.
\ee
The spatial correction involving $b_*^\ell$ plays the role of the constraint term in Dafermos' relative entropy current~\cite{Dafermos2013Involutions}. The total remainder in the divergence is
\be
\label{eq:general-entropy-total-remainder}
\mathscr R_*(t,x;Z)\coloneqq\Delta\mathscr G-\lambda_*^a\Delta\mathscr S_a-(\del_\mu\lambda_*^a)\Delta\mathscr F_a^\mu+(\del_i b_*^\ell)\Delta\mathscr K_\ell^i.
\ee

 \paragraph{{Weak--strong stability.}}
{For uniqueness of the selected fields and stress, we will use the implication $\mathscr Q^0(Z\mid Z_*)=0\ \Longrightarrow\ Z=Z_*$ on the comparison set. In a reduced formulation, the implication concerns the fields retained in the balance laws and their stress. The omitted kinematic fields satisfy separate equations.}

\begin{theorem}[Weak--strong stability with cavitation]
\label{thm:physical-entropy-vacuum}
\label{thm:general-constitutive-entropy-comparison}
{On a smooth prescribed periodic Cauchy slab with spacelike slices, fix a constitutive law and the set of constitutive values used for comparison.} Fix a strong solution of Definition~\ref{def:entropy-strong-solution} and an entropy weak solution $(F,G,\particlecurrent)$ of Definition~\ref{def:weak-local-complete}, represented together with its constitutive stress by $Z$, with integrable fluxes in the chosen comparison set, with exact minors or independent augmented fields according to the chosen formulation. Use integrable prescribed initial data sets from the same set, with finite initial relative entropy, and impose the spatial involutions in the sense of distributions. {Assume that for some $a\in L^1(0,T)$ the remainder satisfies the following bound along the two solutions for almost every $(t,x)\in(0,T)\times\mathbb T^3$:}
\be
\label{eq:general-entropy-weighted-continuity}
{\mathscr R_*(t,x;Z(t,x))\le a(t)\mathscr Q^0(Z(t,x)\mid Z_*(t,x)).}
\ee
{This bound applies also on the zero-density sets of either solution.} Impose the evolution laws and entropy inequality with their initial terms.
\begin{itemize}
\item[(i)] The relative entropy satisfies, for almost every $t$,
\be
\label{eq:general-constitutive-entropy-stability}
\int_{\mathbb T^3}\mathscr Q^0(Z(t)\mid Z_*(t))\,dx\le\exp\left(\int_0^t a(s)\,ds\right)\int_{\mathbb T^3}\mathscr Q^0(Z_0\mid Z_{*,0})\,dx.
\ee
\item[(ii)] Zero initial relative entropy gives zero relative entropy almost everywhere in the slab.  It also gives the mathematical entropy equality with the prescribed initial term. Every physical quantity determined by this vanishing agrees with its reference value. {If zero relative entropy determines the selected fields and stress,} identical initial data for them imply equality of the selected physical fields and stress; in the full elastic formulation these are $F,G,\particlecurrent,T$. The physical velocity is then identified wherever $\particlecurrent$ is timelike. The estimate holds on the whole slice, including the zero-density sets of both solutions. The constitutive graph may include loss of rank of $F$ or $G$, zero current with zero stress or a residual stress as in Definition~\ref{def:residual-stress-zero-current}, and nonzero null current. At each such point of the strong reference, the hypotheses include the finite supporting covector and the regularity of Definition~\ref{def:entropy-strong-solution}. {Under this same implication,} equal initial data sets give agreement of the sets where $\particlecurrent=T=0$, up to null sets.
\end{itemize}
\end{theorem}

{We may verify~\eqref{eq:general-entropy-weighted-continuity} by proving the inequality for every admissible constitutive value at the given reference, as in Corollary~\ref{cor:PC3-sufficient-entropy-bound} below. The same conclusions hold under the weaker condition}
\be
\label{eq:general-entropy-integrated-continuity}
 \int_{\mathbb T^3}\mathscr R_*(t,x;Z(t,x))\,dx
 \le a(t)\int_{\mathbb T^3}\mathscr Q^0(Z(t,x)\mid Z_*(t,x))\,dx
\ee
{for almost every $t$. This allows spatial cancellations in $\mathscr R_*$. The comparison in Theorem~\ref{uc:thm:integrated-3d} below uses a corrected relative current and a displacement term. If a nonnegative source $\mathcal P_{\rm src}\in L^1(0,T)$ is added to the right-hand side, Lemma~\ref{thm:common-current-jhde} below gives the estimate with entropy production in~\eqref{eq:relative-entropy-production-bound}, including the possible initial atom. Its initial value is the relative density of the prescribed initial data.}

\paragraph{{Derivation of the relative entropy estimate.}}
{We apply the supporting inequality on the initial slice and use the spatial compatibility equations in the bulk.}
\begin{proof}[Proof of Theorem~\ref{thm:physical-entropy-vacuum} and the preceding observation]
{The supporting inequalities of Definition~\ref{def:entropy-strong-solution} give $\mathscr Q^0\ge0$ on the comparison solution, including at the constitutive boundary, and on its prescribed initial data.} We apply Lemma~\ref{thm:common-current-jhde} below with $w_*^a=\lambda_*^a$, $L^0=0$ and $L^i=b_*^\ell\Delta\mathscr K_\ell^i$. Each difference $\Delta\mathscr K_\ell$ is divergence-free in space. Testing with $b_*^\ell\varphi$ by smooth approximation gives $\del_iL^i=(\del_i b_*^\ell)\Delta\mathscr K_\ell^i$ after summation, since the coefficients are Lipschitz and the fluxes integrable. The lemma therefore gives the relative entropy inequality with remainder~\eqref{eq:general-entropy-total-remainder}.

The initial component is $\Delta\mathscr I^0_{\rm in}-\lambda_*^a(0,x)\Delta\mathscr F^0_{a,\rm in}=\mathscr Q^0(Z_0\mid Z_{*,0})$ by~\eqref{eq:general-entropy-first-remainder}. {Integrating the pointwise bound in space gives~\eqref{eq:general-entropy-integrated-continuity} with} {$\mathcal P_{\rm src}=0$}. Lemma~\ref{thm:common-current-jhde} below gives the estimate with the production measure and its initial atom. Dropping the nonnegative measure term gives~\eqref{eq:general-constitutive-entropy-stability}. If {$E_0=\mathcal P_{\rm src}=0$}, the lemma gives zero relative density and zero entropy production on the slab including the initial slice. {The implication from zero relative entropy then gives equality of the selected fields and stress.}\qedhere
\end{proof}

{For a perfect fluid the stress multiplier in our flux convention is $u_{*\mu}/\rrestnew'(n_*)$. The relation $\rrestnew+p=n\rrestnew'(n)$ gives the ratio used in~\cite{LeFloch-JHDE}. For $\rrestnew(n)=m\,n+\kappa n^\gamma/\gamma$ with $m>0$, the multiplier has a finite limit at vacuum when the reference velocity does. The supporting inequality and the remainder bound must still hold there.}

{Under the same remainder bound and the remaining reference conditions, the comparison applies row by row. For an evolution flux difference in $L^{p_a}$, $1<p_a\le\infty$, its bounded coefficient $\lambda_*^a$ belongs to $W^{1,p_a'}$, where $p_a'$ is conjugate to $p_a$. The trace of $\lambda_*$ gives nonnegative initial relative density. For a spatial compatibility row in $L^p$ with $1<p\le\infty$, we require $b_*^\ell$ bounded and $\nabla_x b_*^\ell\in L^{p'}$. Interior mollification gives the product rule by duality, and time cutoffs extend the spatial identity to the initial slice. The spatial correction has zero time component. Lemma~\ref{thm:common-current-jhde} gives the comparison with the Sobolev trace of $\lambda_*$.}

%-----------------------------------------------------------------------------------------------------------------------------------------

\subsection{Verification of the entropy-weighted bound}

\paragraph{The adjoint identity and constitutive remainders.}
To control the full constitutive remainder in Corollary~\ref{cor:PC3-sufficient-entropy-bound}, we also require a bounded measurable covector $\nu_*$ on the space of constraint values satisfying
\be
\label{eq:general-entropy-adjoint-identity}
B_{\eta*}-\lambda_*^aB_{a*}-(\del_\mu\lambda_*^a)A_{a*}^\mu+(\del_i b_*^\ell)M_{\ell*}^i=\nu_*\circ D_Y\Theta(t,x,Y_*).
\ee
{For the laws verified in Section~\ref{subsec:constitutive-stability-verification} below, the three material equations~\eqref{eq:physical-continued-material-EL} give this identity, with a uniform bound for $\nu_*$ as the reference density tends to zero. The equations are imposed on the continued material momentum $\mathcal P_{*A}^\nu$ throughout the slab.}

{
The identity cancels the first constitutive variation modulo the algebraic constraints. It is required at cavitation for the sufficient criterion, since at zero flux the balance laws alone leave the velocity parameter undetermined. The identity involves neither division by density nor inversion of $D_Y\Theta$. Theorem~\ref{thm:physical-entropy-vacuum} uses the entropy-weighted bound.
}

 \paragraph{{Control of constitutive remainders.}}
For every allowed bounded range of reference parameters there is a uniform constant $K$ such that, for almost every $(t,x)$ outside a single spacetime null set, every admissible reference parameter in the range, and every $Z\in\mathcal Z_{\rm ph}(t,x)$,
\be
\label{eq:general-entropy-Taylor-control}
\sum_{a,\mu}|r_{\mathscr F,a}^\mu|+\sum_a|r_{\mathscr S,a}|+|r_{\mathscr G}|+\sum_{\ell,i}|r_{\mathscr K,\ell}^i|+|r_\Theta|\le K \mathscr Q^0(Z\mid Z_*).
\ee
Here the same $K$ applies down to zero reference density and for all comparison fields in the prescribed range.

\begin{corollary}[A sufficient constitutive criterion]
\label{cor:PC3-sufficient-entropy-bound}
{Suppose that the law satisfies~\eqref{eq:general-entropy-Taylor-control}.} Let a strong solution in the sense of Definition~\ref{def:entropy-strong-solution} satisfy the adjoint identity~\eqref{eq:general-entropy-adjoint-identity} with bounded $\nu_*$. Then~\eqref{eq:general-entropy-weighted-continuity} holds with $a(t)=C$ for every admissible comparison. The constant depends only on $K$, the reference range, the Lipschitz norms of $\lambda_*,b_*$, $\|\nu_*\|_\infty$ and the geometry. Theorem~\ref{thm:physical-entropy-vacuum} gives the estimate with factor $e^{Ct}$ for every admissible weak solution with finite initial relative entropy, {and gives equality of fields and stress for identical initial data if zero relative entropy determines them.}

For a fixed reference, the constitutive remainder bound~\eqref{eq:general-entropy-Taylor-control} may be replaced by the following bound, with a nonnegative $k_*\in L^1(0,T)$, for  every collection of comparison fields and its constitutive stress on the admissible constitutive flux set outside a common spacetime null set:
\be
 \sum_{a,\mu}|r_{\mathscr F,a}^\mu|+\sum_a|r_{\mathscr S,a}|+|r_{\mathscr G}|
 +\sum_{\ell,i}|r_{\mathscr K,\ell}^i|+|\nu_*(r_\Theta)|
 \le k_*(t)\mathscr Q^0(Z\mid Z_*).
\ee
{In this reference-dependent variant, a measurable covector $\nu_*$ that is finite almost everywhere and satisfies the adjoint identity suffices.} The same conclusions then hold with exponential factor $\exp(C_*\int_0^t k_*(s)\,ds)$, where $C_*$ depends on the bounds for $\lambda_*,\partial\lambda_*,\partial b_*$. This condition controls only the component of $r_\Theta$ selected by the reference.
\end{corollary}

\paragraph{{Proof of the constitutive criterion.}}
{We use the adjoint identity to cancel the part of the flux difference that is linear in the fields.}
\begin{proof}
Both solutions satisfy the physical constraints. Their affine remainders and the adjoint identity~\eqref{eq:general-entropy-adjoint-identity} give
\bse
\be
D_Y\Theta(t,x,Y_*)\Delta Y=-r_\Theta,
\ee
\be
\label{eq:general-entropy-remainder-identity}
\mathscr R_*=-\nu_*(r_\Theta)+r_{\mathscr G}-\lambda_*^ar_{\mathscr S,a}-(\del_\mu\lambda_*^a)r_{\mathscr F,a}^\mu+(\del_i b_*^\ell)r_{\mathscr K,\ell}^i.
\ee
\ese
The constitutive remainder bound~\eqref{eq:general-entropy-Taylor-control} and the bounded reference coefficients yield $|\mathscr R_*|\le C \mathscr Q^0(Z\mid Z_*)$ for every collection of fields and constitutive stress in the chosen comparison set. The relative density is a combination of integrable fluxes with bounded coefficients, so the required product with the constant $C$ is integrable. Theorem~\ref{thm:physical-entropy-vacuum} applies.  Under the reference-dependent bound, the same identity gives $|\mathscr R_*|\le C_*k_*(t)D_\eta$. The factor $\nu_*$ is used only in the controlled contraction $\nu_*(r_\Theta)$ and is never differentiated. Since $k_*\in L^1$, the theorem gives {the asserted estimate}.
\end{proof}

\begin{corollary}[Independence of admissible auxiliary extensions]
\label{cor:physical-independence-vacuum-lift}
Fix the law, constitutive selection, set ${\Cl(\mathcal Z_{\raug})}$ of fields and  stress tensors, including their specified degenerate limits and geometric slab of Theorem~\ref{thm:physical-entropy-vacuum}, and {assume that zero relative entropy determines the selected fields and stress.} Let $(Z_1,\sigma_1)$ and $(Z_2,\sigma_2)$ be strong solutions with cavitation with identical initial data for $F,G,\particlecurrent,T$. Suppose the weighted bound~\eqref{eq:general-entropy-weighted-continuity} holds with the first solution as reference and the second as comparison. Then their selected physical fields and stress agree almost everywhere, even if their auxiliary extensions differ at zero current. In full elasticity this gives equality of $F,G,\particlecurrent,T$ and of the velocity in matter. Under the constitutive remainder bound~\eqref{eq:general-entropy-Taylor-control}, if the first reference also satisfies the adjoint identity~\eqref{eq:general-entropy-adjoint-identity} with bounded coefficient, the weighted bound follows from Corollary~\ref{cor:PC3-sufficient-entropy-bound}.
\end{corollary}
\begin{proof}
The second strong solution is an entropy weak solution with the initial data set given by its values at $t=0$. {The fields of both solutions and the reference supporting covector} are bounded on the compact slab, so their relative entropy is bounded and its product with $a\in L^1(0,T)$ is integrable. Identical initial data for $F,G,\particlecurrent,T$ give zero initial relative entropy, regardless of the auxiliary parameters. We apply Theorem~\ref{thm:physical-entropy-vacuum} to this pair. {Zero relative entropy then gives equality of the fields and constitutive stress.} At zero current the velocity parameters may differ.
\end{proof}

Theorem~\ref{prop:global-rest-current} gives supports in matter on a product of a convex range of spatial forms and densities and a convex velocity range. Choose the observer according to Definition~\ref{def:polyconvex-observers}, and impose the entropy extension and supporting inequality on the range of both solutions and their prescribed initial data. At the constitutive boundary, a finite mathematical entropy value need not admit a finite supporting covector (cf.~Proposition~\ref{prop:vacuum-stress-subgradient}). The weighted estimate for the Ogden-type laws, including $p=q=s=2$, is proved in Theorem~\ref{og:theorem} below.

\paragraph{Relative entropy in matter and at finite flux limits.}
For two solutions with positive density at the point of comparison, take the covariant coefficients $w_*,C_*,D_*$ of~\eqref{eq:observer-current-coefficients} at the reference. The covariant expression in an orthonormal observer frame, before passing to the coordinate densities and spatial correction of~\eqref{eq:general-entropy-relative-current}, is
\be
\label{eq:covariant-relative-current-matter}
 \mathscr Q^\nu=-\Delta {\particlecurrent}^\nu-w_{*\mu}\Delta T^{\mu\nu}
      -C_{*a}^{\nu\alpha}\Delta F_\alpha^a-D_{*b}^{\nu\alpha\beta}\Delta G_{\alpha\beta b}.
\ee
The multiplier identity~\eqref{eq:entropy-multiplier-differential} gives the relative entropy density as its time component. In the variables $(z,v)$ and under the full-rest-domain hypotheses of Theorem~\ref{prop:global-rest-current}, its observer density has the decomposition
\be
\label{eq:global-current-decomposition}
\begin{aligned}
\widetilde\rho_{r*}\mathscr Q^0
={}&\widetilde\rho(z,v_*)-\widetilde\rho(z_*,v_*)
-D_z\widetilde\rho(z_*,v_*)\cdot(z-z_*)\\
&+\widetilde\rho(z,v)-\widetilde\rho(z,v_*)
-D_v\widetilde\rho(z,v)\cdot(v-v_*).
\end{aligned}
\ee
Here $\widetilde\rho_{r*}>0$ and both lines on the right are nonnegative. On a smaller constitutive domain, the same density formula and its nonnegativity hold for each pair whose minor segment at fixed $v_*$ and velocity segment at fixed $z$ stay in the domains with {these convexity and concavity properties}. Consequently, under the hypotheses on the whole rest domain, $\mathscr Q$ is future causal or zero.
Write $z_{u_*}=(\Frest_*,\Grest_*,n_*)$ for the selected rest variables of the reference. If $D^2\bar\rho(z_{u_*})>0$, then zero relative density implies equality of the selected fields, current and stress. This also holds for finite-flux limits of the comparison fields and their stress.

Indeed, using the supporting coefficients of~\eqref{eq:entropy-affine-support-in-matter} in the observer frame gives
\be
 \widetilde\rho_{r*}\mathscr Q^0
 = E-E_*-v_*\cdot(P-P_*)-D_z\widetilde\rho(z_*,v_*)\cdot(z-z_*).
\ee
Substituting $E=\widetilde\rho+v\cdot P$ gives~\eqref{eq:global-current-decomposition}. Its two lines are the remainders from convexity in $z$ at fixed $v_*$ and concavity in $v$ at fixed $z$, respectively. Under the hypotheses on the whole rest domain these signs hold in every orthonormal observer frame, which proves that $\mathscr Q$ is future causal or zero.

To prove equality of the selected fields when the relative density vanishes, we use positive definiteness of the Hessian of the rest energy at the reference state, which makes the first remainder in~\eqref{eq:global-current-decomposition} strictly positive unless $z=z_*$. This remains a separating estimate for any convergent sequence of $z$: we integrate the Hessian on the initial part of its segment from $z_*$, where a uniform positive lower bound holds. Thus a finite-flux limit with relative density zero has $z_j\to z_*$. In particular $r_j\to r_*>0$. The second remainder can be written as
\be
\int_0^1 s\,(v-v_*)^TL(z,v_*+s(v-v_*))(v-v_*)\,ds.
\ee
For $z$ near $z_*$, $L=-D_v^2\widetilde\rho$ has a uniform positive lower bound on a fixed small ball about $v_*$. Since $|v-v_*|<2$, integrating on a fixed initial subinterval gives a positive constant times $|v-v_*|^2$. Hence $v_j\to v_*$ even if the approximating velocities had approached the light cone. The algebraic reconstruction and continuity at the resulting triple of constitutive fields with positive density identify all selected fields, the current and stress.

The covariant formula for $\mathscr Q$ uses finite products and differences, so it passes continuously to {these finite-flux limits}, including a reference boundary limit when its coefficients converge. The closed future causal cone preserves causality and positivity at those degenerate states.

%-----------------------------------------------------------------------------------------------------------------------------------------

\subsection{\texorpdfstring{Integration of the relative entropy current}{Integration of the relative entropy current}}

\paragraph{The comparison principle.}
The divergence inequality for the relative current yields a time estimate once its spatial flux has been integrated over the torus.

\begin{lemma}[Comparison by the relative entropy current]
\label{thm:common-current-jhde}
On $[0,T]\times\mathbb T^3$, let the currents $\mathscr F_a^\mu,\mathscr I^\mu$ and the sources $\mathscr S_a,\mathscr G$ be integrable and satisfy
\be
\del_\mu\mathscr F_a^\mu=\mathscr S_a,
\qquad \del_\mu\mathscr I^\mu\le\mathscr G,
\ee
with prescribed integrable initial data sets. Let a reference, denoted by a star, have integrable currents, sources and initial data sets and satisfy the same balance laws and the entropy equality. All currents and sources used below are evaluated on these two states at the same spacetime point. Let $w_*^a$ be Lipschitz on the closed slab.  Alternatively, suppose the evolution flux differences belong to $L^p$ for some $1<p\le\infty$. Use coefficients $w_*^a\in W^{1,p'}\cap L^\infty$ on the slab and their Sobolev traces at $t=0$, where $p'$ denotes the conjugate exponent. More generally, for each of the finitely many rows one may require $\mathscr F_a-\mathscr F_{a*}\in L^{p_a}$ and $w_*^a\in W^{1,p_a'}\cap L^\infty$, with $1<p_a\le\infty$ and the corresponding value at $t=0$ supplied by the Sobolev trace theorem. Let $L^0=0$ and let the integrable spatial correction $L^i$ have an integrable distributional divergence. Set
\bse
\be
\label{eq:common-current-jhde}
\mathscr Q^\mu=\mathscr I^\mu-\mathscr I_*^\mu
-w_*^a(\mathscr F_a^\mu-\mathscr F_{a*}^\mu)+L^\mu,
\ee
\be
\label{eq:common-remainder-jhde}
\mathscr R=\mathscr G-\mathscr G_*
-w_*^a(\mathscr S_a-\mathscr S_{a*})
-(\del_\mu w_*^a)(\mathscr F_a^\mu-\mathscr F_{a*}^\mu)
+\del_iL^i.
\ee
\ese
{Suppose $\mathscr Q^0\ge0$, its initial integral $E_0$ is finite, and, for $a,\mathcal P_{\rm src}\in L^1(0,T)$ with $\mathcal P_{\rm src}\ge0$,
\be
 \int_{\mathbb T^3}\mathscr R(t,x)\,dx
 \le a(t)\int_{\mathbb T^3}\mathscr Q^0(t,x)\,dx+\mathcal P_{\rm src}(t)
\quad\hbox{for almost every }t.
\ee
Putting $E(t)\coloneqq\int_{\mathbb T^3}\mathscr Q^0(t,x)\,dx$ and $A(t)\coloneqq\int_0^t a(s)\,ds$, one has
\be
 E(t)\le e^{A(t)}\left(E_0+\int_0^t e^{-A(s)}{\mathcal P_{\rm src}(s)}\,ds\right)
 \quad\hbox{for almost every }t.
\ee
The pointwise bound $\mathscr R\le a(t)\mathscr Q^0$ gives this condition with $\mathcal P_{\rm src}=0$. The function $E$ has a nonnegative right-continuous representative $E^+$ of bounded variation on $[0,T)$, with $E^+(0)\le E_0$, and the estimate holds for $E^+(t)$ at every $t\in[0,T)$. There is a nonnegative finite measure $\mu$ on $[0,T)$ such that, with $R(t)=\int_{\mathbb T^3}\mathscr R(t,x)\,dx$,
\be
 E^+(t)+\mu([0,t])=E_0+\int_0^t R(s)\,ds,\qquad 0\le t<T.
\ee
Its initial atom is $\mu(\{0\})=E_0-E^+(0)$, and
\be
\label{eq:relative-entropy-production-bound}
 E^+(t)+\int_{[0,t]}e^{A(t)-A(s)}\,d\mu(s)
 \le e^{A(t)}\left(E_0+\int_0^t e^{-A(s)}{\mathcal P_{\rm src}(s)}\,ds\right).
\ee
The measure $\mu$ is the time marginal of the nonnegative entropy production measure on $[0,T)\times\mathbb T^3$, defined using the prescribed initial mathematical entropy data. If $\mathcal P_{\rm src}=0$ and $E_0=0$, the entropy balance holds with equality for all tests, including those meeting the initial slice. If zero relative density implies equality of specified physical observables, equal initial observables giving $E_0=0$ imply their equality almost everywhere in the slab.} Neither invertibility of the map to conservative variables at the boundary of the state set nor a quadratic lower bound there is required.
\end{lemma}
\paragraph{{Entropy production and the initial value.}}
{We first derive the distributional inequality, then integrate it in time with the prescribed initial term.}
\begin{proof}
We integrate the difference of the balance laws against $w_*^a$ times a smooth scalar test-function and subtract the resulting identity from the difference of the entropy inequality and the reference equality. We extend the Lipschitz coefficients to a neighborhood of the closed slab and mollify. The smooth approximations converge uniformly on the slab, including at $t=0$, with uniformly bounded derivatives converging almost everywhere. Integrable fluxes, sources and prescribed initial data sets justify passage to the limit in the integral identities defining the weak formulation.  In the Sobolev case, reflect across the time endpoints before mollifying. The approximations are uniformly bounded, converge in $W^{1,p'}$, and converge in $L^1$ on the initial slice by the trace theorem. We extract a subsequence converging almost everywhere in the slab and on that slice. Strong convergence of their derivatives in $L^{p'}$ pairs with the $L^p$ evolution flux differences. For row-dependent exponents, this construction is made separately for each coefficient $w_*^a$; H\"older's inequality uses $(p_a,p_a')$, and the finite sum then converges. For the undifferentiated coefficients, dominated convergence applies also to the integrable sources and initial data sets. The spatial correction can be tested up to the initial slice by time cutoffs that vanish near $t=0$ and tend to one for $t>0$. Only spatial derivatives of the tests enter, and both $L^i$ and $\del_iL^i$ are integrable. The correction has zero time component, so the prescribed initial terms come from the entropy and evolution fluxes alone. The product rule gives
\be
\del_\mu\mathscr Q^\mu\le\mathscr R
\ee
in the sense of distributions, including the initial term $E_0$. Spatial integration against the constant function on the torus removes the flux.  The hypotheses give $\mathscr R\in L^1$. Writing $R(t)\coloneqq\int_{\mathbb T^3}\mathscr R(t,x)\,dx$, the distributional inequality and the initial term give
\be
 0\le E(t)\le E_0+\int_0^t R(s)\,ds\le E_0+\|R\|_{L^1(0,T)}
\ee
for almost every $t$: indeed $E(t)-\int_0^tR(s)\,ds$ has a nonincreasing representative whose initial value is at most $E_0$. Thus $|a|E$ and $|a|\mathscr Q^0$ are integrable. The nonincreasing representative of $E-\int_0^\cdot R$ has a right-continuous version. Adding back the continuous integral defines $E^+$, which is nonnegative and has bounded variation on $[0,T)$, with $E^+(0)\le E_0$. Its distributional derivative on $(0,T)$ satisfies $dE^+\le R\,dt$. The measure $\mu$ is $R\,dt-dE^+$ there, with the additional atom $E_0-E^+(0)$ at zero; it is finite and nonnegative, and its cumulative identity follows. The product rule for $E^+$ and the absolutely continuous factor $e^{-A}$, together with {$R\le aE+\mathcal P_{\rm src}$}, yields the displayed estimate including $\mu$. At $t=0$ the initial atom gives equality with $E_0$. The estimates for $E^+$ follow by dropping the nonnegative measure term.

We extend $\mathscr I$ and $\mathscr G$ by zero to $t<0$. Adding the prescribed initial mathematical entropy data times $\delta_{t=0}$ to the extended source and subtracting the divergence of the extended current defines a positive distribution $\mathscr D$. It is a Radon measure by the entropy inequality. The distributional product rule expresses $\mathscr D$ in the same way with $\mathscr R$, $\mathscr Q$ and the initial component of $\mathscr Q^0$. Testing with functions of time yields the marginal $\mu$, including its atom $E_0-E^+(0)$. If {$E_0=\mathcal P_{\rm src}=0$}, then $\mu=0$ and positivity gives $\mathscr D=0$. Vanishing of the relative density gives the asserted equality of the specified observables.
\end{proof}

\paragraph{Integration on compact Cauchy slices.}
A vector current $Q^\mu$ defines the flux three-form $\boldsymbol{\mathscr Q}\coloneqq i_Q\mathrm{vol}_g$.  In adapted coordinates put $\omega\coloneqq|g|^{1/2}$. The corresponding vector density is $\mathscr Q^\mu=\omega Q^\mu$. The slice flux is therefore $\mathscr Q^0\,dx$, and $d\boldsymbol{\mathscr Q}=\nabla_\mu Q^\mu\,\mathrm{vol}_g$. The local comparison with a weak solution extends to a compact slice whenever the relative fluxes agree on chart overlaps and the divergence bound holds in each chart.

\begin{proposition}[Relative entropy on compact slices]
\label{prop:compact-relative-flux}
On an oriented Cauchy slab with compact spacelike slices without boundary, let $\boldsymbol{\mathscr Q}$ be an integrable relative flux three-form with nonnegative slice density and finite prescribed initial integral $E_0$. Suppose its local distributional inequalities with initial terms give
\be
 d\boldsymbol{\mathscr Q}\le a(t)\,dt\wedge\boldsymbol{\mathscr Q},
 \qquad a\in L^1(0,T).
\ee
Assume that the right-hand side is integrable on every closed subslab $0\le t\le S<T$. {The slice integral has a nonnegative right-continuous representative $E^+$ of bounded variation on $[0,T)$, with $E^+(0)\le E_0$, and
\be
 E^+(t)\le E_0\exp\left(\int_0^t a(s)\,ds\right),\qquad 0\le t<T.
\ee
In particular, $E_0=0$ implies that the relative slice density vanishes almost everywhere in the slab.}
\end{proposition}
\begin{proof}
We sum the local inequalities using a partition of unity defined up to the initial slice. The terms containing its derivatives cancel because the three-form is global. For every nonnegative smooth time test $\varphi$ supported in $[0,T)$, we obtain
\be
-\int_0^T E(t)\varphi'(t)\,dt
\le E_0\varphi(0)+\int_0^T a(t)E(t)\varphi(t)\,dt,
\qquad E(t)\coloneqq\int_{\Sigma_t}\boldsymbol{\mathscr Q}.
\ee
Here the initial terms sum to $E_0\varphi(0)$ because the partition sums to one on $\Sigma_0$. On each $[0,S)$ the scalar argument of Lemma~\ref{thm:common-current-jhde}, with $R=aE$ and {$\mathcal P_{\rm src}=0$}, gives the asserted estimate. The right-continuous representatives agree on overlaps. Since $a\in L^1(0,T)$, the estimate bounds $E$ uniformly and gives $|a|E\in L^1(0,T)$. The nonincreasing function $E^+(t)-\int_0^t a(s)E(s)\,ds$ is therefore bounded on $[0,T)$, proving bounded variation on the whole interval. When $E_0=0$, nonnegativity gives vanishing of the slice density almost everywhere.
\end{proof}

In local coordinates, the covariant stress and current equations give
\be
\label{eq:intro-geometric-sources}
 \del_\nu(\omega T^{\mu\nu})=-\omega\Gamma^\mu_{\alpha\beta}T^{\alpha\beta},\qquad
 -\del_\mu(\omega {\particlecurrent}^\mu)\le0,
\ee
where $\Gamma^\mu_{\alpha\beta}$ are the Christoffel symbols. Contraction with a reference covector $w_{*\mu}$ combines the stress source and derivatives of $w_*$ into $-\omega(\nabla_\nu w_{*\mu})\Delta T^{\mu\nu}$. Thus the comparison includes the connection terms. The affine-periodic applications and local construction in Minkowski space specify their domains separately.

%============================================================================

\section{Local existence and weak--strong stability with cavitation}
\label{sec:uniform-elastic-fluid-existence}
\label{sec:local-elastic-stability}

%-----------------------------------------------------------------------------------------------------------------------------------------

\subsection{Constitutive laws and initial data}

\paragraph{Density and the rest-frame minors.}
{
{We now introduce elastic laws whose regularized equations extend to zero density. Section~\ref{sec:proofs-local-existence} will prove local existence and the associated comparison estimates.\par}}

{The construction uses timelike currents and their zero-current limits at bounded subluminal velocity. The ranks of $F$ and $G$ can drop independently at any density. The continued velocity transports the forms through an empty region, where the current and stress vanish.}

Set the speed of light to one and fix $\theta>0$, $p_\theta=1/\theta$. In Minkowski coordinates let $\Fsp,\Gsp$ be independent spatial one-form and dual two-form matrices, let $V$ be the velocity, $\Gamma=(1-|V|^2)^{-1/2}$ and $r={\particlecurrent}^0=n\Gamma$. For exact minors in normalized material-volume coordinates, $r=\det\Fsp$. We regard each row as a spatial vector and define
\bse
\be
\label{uf:eq:coordinates}
 q\coloneqq n^\theta,\qquad Q\coloneqq r^\theta=\Gamma^\theta q,\qquad P_V\coloneqq I+\frac{\Gamma^2}{\Gamma+1}V\otimes V.
\ee
\be
\label{uf:eq:rest-shapes}
 \Frest ^A\coloneqq P_V^{-1}\Fsp ^A,\qquad \Grest _A\coloneqq\Gamma^{-1}P_V\Gsp_A,\qquad
 z\coloneqq\varepsilon_{\rm e}q(\Frest ^A,\Grest _A)_{A=1}^3.
\ee
\ese
Here $\varepsilon_{\rm e}\ge0$ is the elastic parameter. The use of a density power follows Makino, Ukai, and Kawashima~\cite{MakinoUkaiKawashima1986} for classical fluids and LeFloch and Ukai~\cite{LeFlochU} for relativistic fluids.

\begin{definition}
\label{def:density-weighted-elastic-law}
Fix $m>0$ and $\theta>0$.  An \underline{\emph{elastic constitutive law}} for the local existence theory is an augmented law $\raug_{\varepsilon_{\rm e}}$ whose expression $\bar\rho_{\varepsilon_{\rm e}}$ in the rest frame is
\be
\label{uf:eq:law}
 \bar\rho_{\varepsilon_{\rm e}}(\Frest ,\Grest ,n)=n\bigl(m+\Phi(q)+\Psi(q,z)\bigr),\qquad
 q=n^\theta,\qquad z=\varepsilon_{\rm e}q(\Frest ,\Grest ),\qquad\varepsilon_{\rm e}\ge0.
\ee
{The functions $\Phi(q)$ and $\Psi(q,z)$ are $C^9$ on neighborhoods of the chosen coefficient ranges, with $z\in\mathbb R^{18}$.  Require $\Phi(0)=\Phi'(0)=0$. The function $\Psi$ is invariant under simultaneous spatial rotations of its six vector arguments and satisfies $\Psi(q,0)=\Psi_z(q,0)=0$.}
\end{definition}

{The conditions at $q=0$ give the Taylor factorization $\Phi(q)=q^2A(q)$, where $A(q)=\int_0^1(1-t)\Phi''(tq)\,dt$.} On bounded rest-minor ranges the internal energy $\bar\rho_{\varepsilon_{\rm e}}-{m \, n}$ is $\Ocal(n^\gamma)$, with $\gamma=1+2\theta>1$. The function $\Psi$ allows anisotropy and nonlinear coupling of the two form families. At $\varepsilon_{\rm e}=0$, the choice $\Phi(q)=a_{\rm fl}q^2$ gives the polytropic energy ${m \, n}+a_{\rm fl}n^\gamma$. 

For the Ogden-type polynomial law~\eqref{eq:intro-local-energy}, we check rest convexity directly. Fix $\varepsilon_{\rm e}>0$, put $Y=\varepsilon_{\rm e}(\Frest,\Grest)$, and write $B=\operatorname{diag}(\alpha I,\beta I)$. At $n>0$, set $X=\dot n/n$. After division by $n^\gamma$, the second variation of the rest energy at $Y=0$ is
\be
 a_{\rm fl}\gamma(\gamma-1)X^2+\dot Y^TB\dot Y .
\ee
For bounded $n$, this quadratic form remains positive definite for sufficiently small $|Y|$, uniformly down to $n=0$ in the scaled variations $(X,\dot Y)$. Indeed, the coefficients from the quadratic elastic energy are polynomial in $Y$, and those from the quartic energy are polynomials in $Y$ multiplied by $n^{\gamma-1}$. Their changes from the displayed form tend uniformly to zero as $|Y|\to0$. Convexity on the resulting product of a density interval and a ball in $Y$ extends to its face $n=0$ by continuity. At zero elastic parameter the same conclusion follows directly from convexity of ${m \, n}+a_{\rm fl}n^\gamma$.

{The quadratic example~\eqref{uf:eq:exact-polytropic-example} uses both Mooney--Rivlin invariants. The construction also allows nonlinear dependence on the forms. Section~\ref{subsec:spectral-elastic-energies} treats Ogden-type energies at rank loss.}

\paragraph{{The quadratic law and the continued velocity.}}
{We first choose a law for which the entropy controls the terms from independent forms in the comparison.}
\begin{definition}
\label{def:quadratic-second-density-weight-law}
Fix $m,\kappa,\alpha,\beta>0$ and $\gamma>1$.  The \underline{\emph{quadratic elastic law}} is
\be
\label{eq:quadratic-second-density-weight-law}
 \bar\rho_{\rm q2}(\Frest,\Grest,n)
 ={m \, n}+\frac{\kappa}{\gamma}n^\gamma
 +\frac{\alpha}{2}n^{2\gamma-1}|\Frest|^2
 +\frac{\beta}{2}n^{2\gamma-1}|\Grest|^2.
\ee
\end{definition}
{With $\theta=(\gamma-1)/2$ and $q=n^\theta$, set $z_F=q^2\Frest$ and $z_G=q^2\Grest$. Then the elastic internal energy equals $n(\alpha|z_F|^2+\beta|z_G|^2)/2$. Taking $\alpha=\beta=0$ gives the corresponding polytropic fluid law.}

\begin{definition}
\label{def:extended-material-velocity}
{On a region where $\particlecurrent$ is timelike or zero, an \underline{\emph{extended material velocity}} is a future unit timelike field $u$ satisfying $u=\particlecurrent/n$ wherever $n>0$, with values prescribed at zero current for the continued evolution.}
\end{definition}
{Here $u=\Gamma(1,V)$ and $\particlecurrent=n\Gamma(1,V)$. The prescribed velocity $V_0$ at $Q_0=0$ is transported according to $\mathrm D_tV=0$. Where $\particlecurrent=T=0$, the physical balance laws leave this continued velocity undetermined.}

\paragraph{Constitutive bounds and the initial elastic fields.}
On $0\le q\le q_*$, we assume
\bse
\be
\label{uf:eq:constitutive-bounds}
 h\coloneqq m+\Phi+\theta q\Phi'\ge h_0>0,\qquad
 H_f\coloneqq\Phi''+(1+1/\theta)\Phi'/q\ge\alpha_{\rm f}>0,\qquad
 \Psi_{zz}(q,0)\ge\alpha_{\rm e} I>0.
\ee
\be
\label{uf:eq:fluid-sound-bound}
 c_f(q)^2\coloneqq\theta^2q^2H_f(q)/h(q)\le1.
\ee
\ese
{The quotient $\Phi'(q)/q=\int_0^1\Phi''(tq)\,dt$ extends with at least seven continuous derivatives. For the perfect-fluid energy $\rho_{\rm f}(n)=n(m+\Phi(q))$, one has $\rho_{\rm f}'=h$ and $\rho_{\rm f}+p=nh$. Thus the variable $w$ of LeFloch and Ukai~\cite{LeFlochU}, calculated with the induced pressure law as a function of $\rho_{\rm f}$ and with light speed one, obeys}
\be
\label{uf:eq:fluid-regularized-coordinate}
 {
 w(\rho_{\rm f}(n))
 =\int_0^n\frac{c_f(\xi^\theta)}{\xi}\,d\xi
 =\int_0^q\sqrt{\frac{H_f(\zeta)}{h(\zeta)}}\,d\zeta .
 }
\ee
By the positive bounds in~\eqref{uf:eq:constitutive-bounds}, the derivative of $w$ with respect to $q$ is strictly positive and belongs to $C^7$ up to vacuum. Thus $q\mapsto w$ and its inverse are $C^8$, and $Q=\Gamma^\theta q$ is obtained by a further change of variable that is smooth in the velocity for $|V|<1$. Since $c_f(q)=\theta q\sqrt{H_f(q)/h(q)}$, the ratio $c_f/w$ tends to $\theta$ as $q\downarrow0$. In particular it satisfies the boundedness condition at vacuum in~\cite[(4.1)]{LeFlochU}. In vacuum our velocity satisfies $\mathrm D_tV=0$, the continuation in~\cite[Definition~5.1]{LeFlochU}. The change of density variable preserves the Sobolev order of the regularized pair. Our data are in global Sobolev spaces on $\R^3$ or $\mathbb T^3$, and the existence result of~\cite[Theorem~4.1]{LeFlochU} uses uniformly local spaces on $\R^3$. Their density variable before regularization is the rest energy $\rho_{\rm f}$, whereas our reconstruction uses the particle density $n$. The tame-solution definition also requires $\rho_{\rm f}$ to be $C^1$, which holds here for $1<\gamma\le3$ since $n=Q^{2/(\gamma-1)}/\Gamma$.

We fix $0\le v_0<v_*<1$ and $0<Q_{\rm in}<Q_*<q_*$. The coefficient construction below gives a number $s_*>0$ depending only on these ranges, the positive constants in the bounds~\eqref{uf:eq:constitutive-bounds}, and bounded constitutive derivatives. We choose $0<s_{\rm in}<s_{\rm out}<s_*$.

For the Sobolev estimates, we use the scaled elastic minors and their products with the regularized density:
\be
\label{uf:eq:scaled-fields}
 {Y\coloneqq\varepsilon_{\rm e}(\Fsp ^A,\Gsp _A)_{A=1}^3,\qquad Z\coloneqq QY,\qquad U\coloneqq(Q,Z,V),\qquad Z_0\coloneqq Q_0Y_0.}
\ee
{On $\mathcal D=\R^3$ or $\mathbb T^3$ prescribe $Q_0,V_0\in H^3(\mathcal D)$, with $0\le Q_0\le Q_{\rm in}$ and $|V_0|\le v_0$. The field $Y_0$, with vector components $Y_{F,0}^A=\varepsilon_{\rm e}\Fsp _0^A$ and $Y_{G,0,A}=\varepsilon_{\rm e}\Gsp _{0,A}$, belongs to $H^3(\mathbb T^3;\R^{18})$ on the torus and to $W^{3,\infty}(\R^3;\R^{18})$ on the whole space, with $\|Y_0\|_\infty\le s_{\rm in}$. The sum of the prescribed norms is bounded by $K_{\rm in}$. The algebra property on the torus and the multiplication estimate for $W^{3,\infty}$ fields on $\R^3$ give $Z_0\in H^3$.}
{For the augmented system, require the spatial one-forms to be closed and the spatial two-forms to be closed, equivalently $\partial_i(\Fsp_0)_{Aj}=\partial_j(\Fsp_0)_{Ai}$ and $\sum_i\partial_i(\Gsp_0)_{Bi}=0$, in the sense of distributions. These are differential constraints, not polynomial relations between the fields. The scalar $r_0=Q_0^{p_\theta}$ is independent of both matrices. For initial data on the set of exact spacetime minors with $\varepsilon_{\rm e}>0$, } the equations~\eqref{uf:eq:scaled-fields} hold initially with $\Fsp _0=Df_0$, $\Gsp _0=\cof \Fsp _0$ and $\det \Fsp _0=Q_0^{p_\theta}$. On the torus the configuration has affine periods and a periodic gradient. The periods may vary in a fluid-limit family and must match the lattice when the material torus is prescribed.  With fixed periods $A$ on the unit torus, the initial particle number is $\det A$ by the periodic determinant identity. At $\varepsilon_{\rm e}=0$, we set $Y_0=Z_0=0$ and prescribe only $Q_0,V_0$.

%-----------------------------------------------------------------------------------------------------------------------------------------

\subsection{Local existence and compatible initial data}

\paragraph{A local theorem uniform in the elastic parameter.}
{We now specify the additional evolution conditions used to construct a solution through cavitation.}
\begin{samepage}
\begin{definition}
\label{def:strong-prescribed-continuation}
{A \underline{\emph{strong solution with prescribed continuation}} on $[0,T]$,  for the constitutive laws and initial data specified above, consists of $U=(Q,Z,V)$ and the scaled elastic fields $Y$, with $Q\ge0$, $Z=QY$ and $|V|\le v_*<1$, such that:
\begin{itemize}
\item[(i)] For the chosen exponent $5/2<s\le3$, one has $U\in C([0,T];H^s)\cap C^1([0,T];H^{s-1})$. The fields attain the prescribed initial data, and $Y$ is transported from $Y_0$ by the particle flow according to~\eqref{uf:eq:flow-formulas} below.
\item[(ii)] The normalized equations~\eqref{uf:eq:normalized-evolution}--\eqref{uf:eq:shape-evolution} below hold with their prescribed continuous coefficient extension at $Q=0$.
\end{itemize}}
\end{definition}
\end{samepage}

{On $\{Q=0\}$ the continued equation for the extended material velocity is $\mathrm D_tV=0$.}
{For comparison, one must also verify the conditions of Definition~\ref{def:entropy-strong-solution} and Theorem~\ref{thm:physical-entropy-vacuum}.}
{This velocity continuation follows Makino, Ukai, and Kawashima~\cite{MakinoUkaiKawashima1986} and the relativistic-fluid construction of LeFloch and Ukai~\cite{LeFlochU}.}
\paragraph{{The local existence theorem.}}
{The independent spatial forms are transported by the same flow as the density.}
\begin{theorem}[ Local existence with possible cavitation and rank loss]
\label{uf:thm:main}
Fix a constitutive law~\eqref{uf:eq:law} satisfying~\eqref{uf:eq:constitutive-bounds}--\eqref{uf:eq:fluid-sound-bound}. Prescribe compatible initial fields~\eqref{uf:eq:scaled-fields} with $Q_0,V_0\in H^3$ and the bounds, elastic-field regularity and strict range margins of Section~\ref{sec:uniform-elastic-fluid-existence}. There exists $T>0$, depending only on those bounds and their strict margins, for which the following conclusions hold.
\begin{itemize}
{\item[(i)] For $\varepsilon_{\rm e}>0$, prescribe independent initial form matrices satisfying the spatial closure constraints and {the prescribed bounds}. The normalized system~\eqref{uf:eq:normalized-evolution}--\eqref{uf:eq:shape-evolution}, with its continuous coefficient extension at $Q=0$, admits a unique strong solution with prescribed continuation in the sense of Definition~\ref{def:strong-prescribed-continuation}, in the class
}\be
\label{uf:eq:solution-class}
 U\in C([0,T];H^3)\cap C^1([0,T];H^2),\qquad
 Y-Y_0\in C([0,T];H^2)\cap C^1([0,T];H^1).
\ee
{The particle flow $X$ is a $C^1$ diffeomorphism. It transports the two independent form families and the particle density. The ranks of $F$ and $G$ remain constant along each particle trajectory, without any relation between these ranks or the zero-density set. The reconstructed forms are closed, the current is conserved, and the constitutive stress satisfies its balance laws with the prescribed initial terms. The density and current are continuous; this assertion does not require differentiating $Q^{p_\theta}$ at zero.
}If $\mathcal V_0=\{Q_0=0\}$, then
\be
\label{uf:eq:inertial-vacuum-flow}
 \{Q(t)=0\}=X(t,\mathcal V_0),\qquad X(t,y)=y+tV_0(y)\quad(y\in\mathcal V_0),
\ee
with the periodic lift understood on the torus. There $\particlecurrent=T=0$ and the velocity is freely transported.  For every open set $\Omega_0\subset\mathcal V_0$ of finite volume, its image $\Omega_t=X(t,\Omega_0)$ has boundary $X(t,\partial\Omega_0)$ and Euclidean volume
\be
 |\Omega_t|=\int_{\Omega_0}\det(I+t\,DV_0(y))\,dy .
\ee
{Thus an initially empty region can move and change shape.
}\item[(ii)] {The existence time and regularized bounds are independent of $\varepsilon_{\rm e}$ and of a positive lower density.}
\item[(iii)] At $\varepsilon_{\rm e}=0$, set $Y_0=Z_0=0$. The same interval and regularity bounds apply to the reduced density--velocity equations. They give the relativistic barotropic fluid equations in the regularized variables, with the velocity freely transported on the zero-density set.
\item[(iv)] On the set of exact spacetime minors, assume additionally $\Fsp_0=Df_0$, $\Gsp_0=\operatorname{cof}\Fsp_0$ and $r_0=\det\Fsp_0$. Then $f=f_0\circ X^{-1}$ reconstructs the configuration, and all polynomial minor relations persist. For fixed positive elastic parameter, the density and current belong locally to $C_tH^2\cap C_t^1H^1$. The three material variational equations also hold in the sense of distributions. 
\end{itemize}
\end{theorem}

At $Y=Z=0$, the uniform time-matrix bounds still hold, and the form equations preserve $Y=Z=0$. The theorem also holds with $v_*^2\max_{0\le q\le q_*}c_f(q)^2<1$ in place of~\eqref{uf:eq:fluid-sound-bound}. This keeps the fluid Schur complement $h(I-c_f^2V\otimes V)$ uniformly positive at $Y_R=0$ for $|V|\le v_*$. Reducing $s_*$ preserves this positivity down to $q=0$. A light-speed bound for the elastic waves requires a further condition.

\paragraph{A complete comparison theorem for independent augmented fields.}
{For the quadratic law of Definition~\ref{def:quadratic-second-density-weight-law}, the factor $n^{2\gamma-1}$ makes the form multipliers vanish at the order needed for augmented comparison, including at zero density.}

\begin{theorem}[Local existence and augmented stability with cavitation]
\label{thm:quadratic-second-weight-augmented-stability}
Let $\mathcal D=\mathbb T^3$ or $\mathbb R^3$, and consider the law~\eqref{eq:quadratic-second-density-weight-law}.  Put $Q=r^\theta$, where $r=\particlecurrent^0$, and define the spatial fields
\be
 W_F=Q\Fsp,\qquad Z_F=QW_F,\qquad
 W_G=Q\Gsp,\qquad Z_G=QW_G,\qquad
 U_{\rm n}=(Q,Z_F,Z_G,V).
\ee
Choose bounded density and velocity ranges, with $|V|\le v_{\max}<1$, on which the reduced fluid wave speed satisfies
\be
 \frac{\kappa(\gamma-1)q^2}{m+\kappa q^2}\le1-\delta,
 \qquad\delta>0.
\ee
At zero elastic fields, the acoustic form is the reduced fluid form. For the fixed coefficients $\alpha,\beta>0$, restrict the weighted elastic fields to a sufficiently small ball so that the full acoustic bound holds with $1-\delta/2$ and
\be
\label{eq:quadratic-second-weight-rest-margin}
 (p_\theta+2)(p_\theta+4)
 \bigl(\alpha|q\Frest|^2+\beta|q\Grest|^2\bigr)
 <2\kappa p_\theta,
 \qquad p_\theta=1/\theta,
\ee
with a strict uniform margin. Choose a bounded comparison range for the unweighted spatial forms as well. Its transformations to a fixed observer's rest frame and the intervening density--form segments must stay in this convexity range; this holds after decreasing the upper density or form bounds. Prescribe initial fields with strict margins from the upper density, velocity and form bounds, with $Q_0\ge0$, satisfying the spatial closure constraints and
\be
 U_{{\rm n},0}\in H^3(\mathcal D),\qquad
 Z_{F0}=Q_0W_{F0}=Q_0^2\Fsp_0,\qquad
 Z_{G0}=Q_0W_{G0}=Q_0^2\Gsp_0.
\ee
{Prescribe the unweighted initial forms in $H^3$ on $\mathbb T^3$ or $W^{3,\infty}$ on $\mathbb R^3$; their products $W_{F0}=Q_0\Fsp_0$, $W_{G0}=Q_0\Gsp_0$ then belong to $H^3$ on either domain.} The local time and the constants below may depend on $\gamma$, the constitutive coefficients, the initial norm bounds and the strict range margins, but are independent of a positive lower bound for $Q_0$. Then the following conclusions hold on $[0,T]$.
\begin{itemize}
\item[(i)] The normalized evolution with the prescribed continuation at zero density has a unique local solution. Its fields satisfy the augmented balance laws and have the regularity
\be
 U_{\rm n}\in C_tH^3\cap C_t^1H^2,\qquad
 (\Fsp-\Fsp_0,W_F-W_{F0},\Gsp-\Gsp_0,W_G-W_{G0})
 \in C_tH^2\cap C_t^1H^1.
\ee
The particle flow transports the initial zero-density region, and the ranks of $F$ and $G$ remain constant along each particle trajectory.  On the transported cavity one has $\particlecurrent=T=0$.
\item[(ii)] On a smooth periodic slab, let the strong reference be the solution in (i), with the continued equation $D_tV_*=0$ on its cavity.  Let $D_\eta$ be the relative mathematical entropy density of the conservative variables, and let $D_{\eta,0}$ be its value on the prescribed initial data sets.  {For an augmented entropy weak solution of Definition~\ref{def:weak-local-complete} taking values in the same bounded constitutive branch, write $U_0$ and $\particlecurrent_0^0$ for its prescribed initial conservative variables and particle density. Its initial entropy defect is}
\be
 {\mathfrak d_0\coloneqq-\particlecurrent_0^0-\eta(U_0)\ge0.}
\ee
{The relative entropy satisfies, for almost every $t\in[0,T]$,}
\be
 {\int_{\mathbb T^3}D_\eta(t,x)\,dx
 \le e^{Ct}\int_{\mathbb T^3}\bigl(D_{\eta,0}(x)+\mathfrak d_0(x)\bigr)\,dx.}
\ee
{When $\mathfrak d_0=0$, in particular when the prescribed initial fields and stress belong to the same completed constitutive branch, one obtains}
\be
\label{eq:quadratic-second-weight-stability}
 \int_{\mathbb T^3}D_\eta(t,x)\,dx
 \le e^{Ct}\int_{\mathbb T^3}D_{\eta,0}(x)\,dx
\ee
{for almost every $t\in[0,T]$. The same estimate holds if the conservative variables have the strong initial trace $U(t)\to U_0$ in $L^1$ as $t\downarrow0$, even when $\mathfrak d_0>0$.}  No exact or mixed polynomial relation between $F,G,\particlecurrent$ is imposed.  The two form ranks may degenerate independently, the zero-density regions may differ, and the proof uses no derivatives of fields specific to the weak solution.  One imposes $\nabla_\mu\particlecurrent^\mu\ge0$.
\item[(iii)] {If $D_{\eta,0}=0$ and either $\mathfrak d_0=0$ or the above strong initial trace holds, the current and stress agree almost everywhere on $(0,T)\times\mathbb T^3$, and the fields $F$ and $G$ agree almost everywhere where the density is positive.}  The mathematical entropy does not identify independent fields $F$ and $G$ inside a common empty cavity.  {For the local strong solution, the exact polynomial minor relations propagate when they hold initially. Estimate~\eqref{eq:quadratic-second-weight-stability} also applies to weak solutions satisfying these relations when either of the initial conditions in (ii) holds.}
\end{itemize}
\end{theorem}

{For $\alpha=\beta=0$, omitting the form rows and variables gives the regularized polytropic fluid system with $Q,V$ in $H^3$.}

{\paragraph{{Nonlinear elastic response.}}
{We retain the positive quadratic part that gives positive time matrices and add nonlinear terms in the weighted rest fields.}
\begin{corollary}[Nonlinear elastic laws with cavitation]
\label{cor:nonlinear-second-weight}
Fix $m,\kappa>0$, $\gamma>1$ and $\theta=(\gamma-1)/2$. Let $\Psi$ be a $C^6$ function near the origin of the eighteen rest-form components, invariant under simultaneous spatial rotations, with $\Psi(0)=0$, $D\Psi(0)=0$ and $D^2\Psi(0)>0$. Consider
\be
\label{eq:nonlinear-second-weight-law}
 \bar\rho=n\left(m+\frac\kappa\gamma q^2+\Psi(z)\right),
 \qquad q=n^\theta,\qquad z=q^2(\Frest,\Grest).
\ee
Use the initial regularity and spatial compatibility conditions of Theorem~\ref{thm:quadratic-second-weight-augmented-stability}, and fix bounded density and velocity ranges with its strict fluid acoustic margin. For a sufficiently small bound on the unweighted forms, depending also on $\Psi$, one has:
\begin{itemize}
\item[(i)] the local existence and rank-transport conclusions of that theorem on $\mathbb T^3$ or $\mathbb R^3$, with the same prescribed continuation at zero density;
\item[(ii)] its augmented weak--strong comparison and uniqueness conclusions on the periodic slab, for comparison fields in the same bounded constitutive branch, including when the zero-density regions differ {and with the initial entropy defect or strong initial trace specified in that theorem}.
\end{itemize}
The same conclusions hold for the alternative
\be
\label{eq:second-weight-spectral-additions}
 \Psi(z_F,z_G)=\frac12\bigl(\alpha|z_F|^2+\beta|z_G|^2\bigr)
 +\sum_{\sigma\in\{F,G\}}\sum_{j=1}^{N_\sigma}
       \frac{a_{\sigma j}}{p_{\sigma j}}I_{p_{\sigma j}}(z_\sigma),
 \quad \alpha,\beta>0,\quad a_{\sigma j}\ge0,\quad p_{\sigma j}\ge4,
\ee
where $I_s(X)=\sum_{i=1}^3\sigma_i(X)^s$ and $\sigma_i(X)$ are the singular values. For this spectral alternative, one may instead fix $5/2<s<3$, replace $p_{\sigma j}\ge4$ by $p_{\sigma j}>s+1$, and prescribe $U_{{\rm n},0}\in H^s$, retaining {the hypotheses above} on the unweighted initial forms. In (i), the spaces $H^3,H^2,H^1$ are then replaced by $H^s,H^{s-1},H^{s-2}$, respectively; conclusion (ii) remains valid. No exact polynomial minor relations are required. The constants may depend on the positive quadratic margin and on the $C^6$ bounds in the first case, or the fixed coefficients and exponents in~\eqref{eq:second-weight-spectral-additions}.
\end{corollary}
Taking only the fourth powers in the sum in~\eqref{eq:second-weight-spectral-additions} gives quadratic and quartic terms with factors $n^{2\gamma-1}$ and $n^{4\gamma-3}$, respectively. The quadratic part gives the positive elastic Hessian at zero strain variables.
}

\paragraph{Weaker regularity of the initial minors.}
The weighted fields may have one more Sobolev derivative than the unweighted minors.

\begin{proposition}[Local existence with Sobolev-regular initial minors]
\label{uf:cor:fractional-existence}
Let $5/2<s\le3$ and $\mathcal D=\R^3$ or $\mathbb T^3$.  For the whole-space case, set
\be
 {\mathcal Y_{s,\rm uloc}\coloneqq H^{s-1}_{\rm uloc}(\R^3)},\qquad
 \|Y\|_{H^{s-1}_{\rm uloc}}\coloneqq\sup_{k\in\mathbb Z^3}\|\chi(\cdot-k)Y\|_{H^{s-1}},
\ee
where $\chi\in C_c^\infty(\R^3)$ equals one on a neighborhood of the unit cube. At $s=3$ this requires only a uniform local $H^2$ bound.
Under the constitutive assumptions, independent initial spatial closure constraints and strict range margins of Theorem~\ref{uf:thm:main}, replace the initial regularity by
\be
 Q_0,Z_0,V_0\in H^s(\mathcal D),\qquad Z_0=Q_0Y_0,\qquad
 Y_0\in\begin{cases}H^{s-1}(\mathbb T^3),&\mathcal D=\mathbb T^3,\\
 {\mathcal Y_{s,\rm uloc},}&\mathcal D=\R^3.\end{cases}
\ee
\begin{itemize}
\item[(i)] \emph{Local solution and particle flow.} For every fixed $\varepsilon_{\rm e}>0$, the normalized system has a unique local strong solution with prescribed continuation.  If $X$ is its particle flow and $\zeta=X^{-1}-\operatorname{Id}$, then
\be
 U,\zeta\in C_tH^s\cap C_t^1H^{s-1}.
\ee
\item[(ii)] \emph{Regularity of the elastic fields.} The scaled elastic fields are reconstructed by the flow formulas~\eqref{uf:eq:flow-formulas}.  On $0\le t\le\min\{T,1\}$ they satisfy
\be
 {\|Y(t)-Y_0\|_{L^\infty}\le Ct^{s-5/2},\qquad
 \mathcal D=\mathbb T^3\ \hbox{or}\ \mathbb R^3.}
\ee
where $C$ depends on {the initial data bounds} and margins. For each fixed  initial data set on the torus, the left-hand side is $o(t^{s-5/2})$ as $t\downarrow0$. On either domain, $Y_0\in W^{1,\infty}$ gives the stronger bound $Ct$. If $\mathcal D=\mathbb T^3$ or $Y_0\in W^{3,\infty}$, then $Y-Y_0\in C_tH^{s-1}\cap C_t^1H^{s-2}$.
\item[(iii)] \emph{Cavitation and conservation laws.} The particle flow transports the initial zero-density set and preserves the ranks of the reconstructed minors. The current and stress vanish on that set and satisfy their conservation laws with initial terms. If the initial data satisfy the exact minor relations, the three material equations also hold, and for fixed positive $\varepsilon_{\rm e}$ the density and current belong locally to $C_tH^{s-1}\cap C_t^1H^{s-2}$. The local interval depends on {the prescribed bounds} and strict margins.
\item[(iv)] \emph{Perfect fluids.} In the fluid case $\varepsilon_{\rm e}=0$, the reduced density--velocity system has the same $H^s$ regularity for $Q,V$, without elastic initial fields. Its flow transports the vacuum, and the current and stress satisfy their conservation laws with initial terms.
\end{itemize}
\end{proposition}

\paragraph{An elastic law with a polytropic density factor.}
Fix $\gamma=1+2\theta>1$ and $K,\alpha,\beta>0$. Write $z=(z_F,z_G)$ for the weighted one-form and two-form components in~\eqref{uf:eq:rest-shapes}. Taking $\Phi(q)=Kq^2/(\gamma-1)$ and $\Psi(q,z)=(\alpha|z_F|^2+\beta|z_G|^2)/2$ gives the rest energy
\be
\label{uf:eq:exact-polytropic-example}
\bar\rho(\Frest,\Grest,n)
={m \, n}+n^\gamma\left(\frac{K}{\gamma-1}+\frac{\varepsilon_{\rm e}^{\,2}}2
  \bigl(\alpha|\Frest |^2+\beta|\Grest |^2\bigr)\right).
\ee
At $\varepsilon_{\rm e}=0$, the pressure is $Kn^\gamma$ and $c_f^2=\gamma K n^{\gamma-1}/h$ with $h=m+\gamma K n^{\gamma-1}/(\gamma-1)$. The constants $H_f=2K(2+1/\theta)/(\gamma-1)$ and $\Psi_{zz}=\operatorname{diag}(\alpha I,\beta I)$ give the positivity in~\eqref{uf:eq:constitutive-bounds}. For every $\gamma>1$ we can take a sufficiently small density range for~\eqref{uf:eq:fluid-sound-bound}.  For comparison with the density coordinate $w$ of LeFloch--Ukai~\cite{LeFlochU}, write $\rho_{\rm d}(n)={m \, n}+Kn^\gamma/(\gamma-1)$ and $q=n^\theta$. Since $d\rho_{\rm d}/dn=h=(\rho_{\rm d}+p)/n$, their definition gives
\be
 w(\rho_{\rm d}(n))=\int_0^n\frac{\sqrt{dp/d\rho_{\rm d}}}{n'}\,dn'
 =\frac{2}{\sqrt{\gamma-1}}\operatorname{arsinh}\!\left(\sqrt{\frac{\gamma K}{m(\gamma-1)}}\,q\right),
\ee
{where the right-hand side is a smooth invertible function of $q$ near zero and gives the change of density variable.}

The Ogden-type law~\eqref{eq:intro-local-energy} has the same fluid coefficients with $a_{\rm fl}=K/(\gamma-1)$ and the same $\Psi_{zz}(q,0)$. Its fourth powers are polynomial in the weighted minors. The theorem applies with the elastic size bound chosen for all four coefficients. At unit volume the principal-stretch powers are $-2,2,-4,4$ (cf.~Section~\ref{subsec:spectral-elastic-energies} below).

\paragraph{A moving cavity with exact minors.}
{In a unit periodic cell take $f_0=(x^1\varphi(|x|),Ax^2,Bx^3)$ with $A,B>0$, and $V_0=\upsilon\sin(2\pi x^2)e_1$. Choose $0<R_1<R_2<1/2$ and a smooth nondecreasing $\varphi$ which vanishes on $[0,R_1]$, equals $\exp(-1/(r-R_1)^2)$ just outside $R_1$, and is a positive constant for $r\ge R_2$. Outside the larger ball $f_0$ is affine, so its gradient extends periodically. For $r=|x|>0$ one has
\be
\label{uf:eq:moving-vacuum-data}
 \det Df_0=AB\bigl(\varphi(r)+(x^1)^2\varphi'(r)/r\bigr).
\ee
The determinant vanishes on the inner ball, where the gradient has rank two. Its power $\theta$ is smooth across the sphere because the exponential dominates every inverse power of $r-R_1$, also after differentiation. Scaling $f_0$ and $\upsilon$ gives the strict initial margins. On the ball the continued equation $D_tV=0$ gives
\be
\label{uf:eq:moving-vacuum-flow}
 X(t,y)=y+t\upsilon\sin(2\pi y^2)e_1 .
\ee
Thus the cavity deforms. In the adjacent matter $\partial_t r(0)=-V_0\cdot\nabla r_0$ is nonzero wherever the shear has a nonzero component along the density gradient.}

The rest-frame minor identities extend using the transported velocity in~\eqref{uf:eq:rest-shapes}. The cofactor has rank one at rank two and vanishes at ranks one and zero. The two vanishing conditions on $\Psi$ make $n\Psi$ and its first two rest-minor derivatives $\Ocal(n^\gamma)$ on bounded rest-minor sets, regardless of rank.

%-----------------------------------------------------------------------------------------------------------------------------------------

\subsection{Entropy and stability for the local constitutive laws}
\label{uf:par:stability-scope}

\paragraph{A convex entropy on a bounded constitutive domain.}
{We return to the laws of Definition~\ref{def:density-weighted-elastic-law}. For comparison with exact minors, we either bound the divergences of the reference form multipliers or use the displacement of corresponding material points. Write $U=(Q,Z,V)$ for the normalized variables and $T^{0\rm sp}\coloneqq(T^{0i})_{i=1}^3$. The conservative variables are}
\be
\label{uf:eq:conservative-tuple}
 U_{\rm c}\coloneqq(\Fsp ,\Gsp ,T^{00},T^{0\rm sp}).
\ee
The density measured by this observer is $r={\particlecurrent}^0=Q^{p_\theta}$. For reference values we write $Q^\ast=r_*^\theta$ and $Z^\ast=Q^\ast Y_*$, distinguishing them from the upper bounds $q_*,Q_*$. The following construction allows independent minors.

\begin{proposition}[Convex entropy for the local constitutive laws]
\label{uf:prop:local-convex-entropy}
Fix a positive elastic parameter and the constitutive functions of Theorem~\ref{uf:thm:main},  with the constitutive bounds~\eqref{uf:eq:constitutive-bounds}. Choose $r_{\max}>0$ and $0<v_{\rm c}<1$, with $r_{\max}^{\theta}\le Q_*$. Instead of the sound-speed bound~\eqref{uf:eq:fluid-sound-bound}, assume that for some $\delta_{\rm c}>0$,
\be
 \sup_{\substack{0\le r\le r_{\max}\\|V|\le v_{\rm c}}}
 c_f\bigl((r/\Gamma(V))^\theta\bigr)^2|V|^2\le1-\delta_{\rm c},
 \qquad \Gamma(V)\coloneqq(1-|V|^2)^{-1/2}.
\ee
Choose $s_{\rm c}>0$ sufficiently small, depending also on $\delta_{\rm c}$. On the completed finite-velocity branch
\be
\label{uf:eq:bounded-comparison-branch}
 \mathcal K\coloneqq\{0\le r\le r_{\max},\ |Y|\le s_{\rm c},\ |V|\le v_{\rm c}\},\qquad Y=\varepsilon_{\rm e}(\Fsp ,\Gsp ),
\ee
{ the following conclusions hold.}
\begin{itemize}
\item[(i)] There is a proper convex entropy $\eta_{\mathcal K}$ with value $-r$ on the conservative image. The covectors selected by the constitutive multipliers extend continuously to zero density and support this entropy there. \item[(ii)] For every pair of states in this branch, its relative density satisfies
\be
\label{uf:eq:local-relative-entropy-comparison}
 D_{\eta_{\mathcal K}}(U_{\rm c}\mid U_{{\rm c}*})
 \asymp (r+r_*)\bigl(|Q-Q^\ast|^2+|Z-Z^\ast|^2\bigr)+r|V-V_*|^2.
\ee
The constants depend only on {the chosen coefficient ranges}, constitutive bounds and the observer and smallness margins; they do not depend on a positive lower density or on derivatives of either state. The corresponding rest spatial stress has its first-order remainder bounded by a constant times the rest-energy remainder. \item[(iii)] At zero elastic parameter the same conclusions hold in the reduced fluid variables, with $Y=Z=0$.
\end{itemize}
\end{proposition}
\paragraph{Pointwise comparison for the local constitutive laws.}
\label{uc:subsec:common-class}
Derivatives of the reference coefficients enter the bound on the remainder.

\begin{theorem}[Comparison under bounds on the form multipliers]
\label{uc:thm:common}
Fix a positive elastic parameter and a periodic reference in the class of solutions satisfying the exact minor relations of Theorem~\ref{uf:thm:main} or of Proposition~\ref{uf:cor:fractional-existence} with $5/2<s\le3$, and write $\gamma=1+2\theta$. Choose a bounded comparison domain of the form~\eqref{uf:eq:bounded-comparison-branch} containing {the range of the reference initial data} with strict upper bounds on density, velocity and elastic amplitude. Choose its elastic bound so that the orthogonal projections of affine segments of the original $(F,G,\particlecurrent)$ fields stay inside the constitutive domain where {the prescribed bounds} on the Hessian of the energy in its independent rest variables hold. On a possibly shorter interval the reference stays in this domain. Use the entropy $\eta_{\mathcal K}$ of Proposition~\ref{uf:prop:local-convex-entropy} and the reference form multipliers $C_*,D_*$ of~\eqref{eq:main-vacuum-finite-multipliers}. Assume, for some nonnegative $b_0\in L^1(0,T)$, that almost everywhere
\be
\label{uc:eq:geometric-coefficients}
 |\nabla_\nu C_*^{\nu\alpha}|+|\nabla_\nu D_*^{\nu\alpha\beta}|
 \le b_0(t)n_*^\gamma.
\ee
{The following estimates hold.}
\begin{itemize}
\item[(i)] Outside a single spacetime null set, every collection of comparison fields and its constitutive stress $\mathcal Z$ on this branch with the exact polynomial minor relations satisfies
\be
\label{uc:eq:weighted-bound}
 |\mathscr R_*(t,x;\mathcal Z)|\le C(1+b_0(t))D_{\eta_{\mathcal K}}(U_{\rm c}(\mathcal Z)\mid U_{{\rm c}*}(t,x)),
\ee
{where $C$ depends on the fixed elastic parameter, constitutive bounds, comparison domain and reference norms in the local existence result.}

\item[(ii)] For every entropy weak solution on this branch with exact minors and prescribed initial terms in the weak formulation,~\eqref{eq:general-constitutive-entropy-stability} holds with $a(t)=C(1+b_0(t))$. Zero initial relative entropy gives equality almost everywhere of current and stress, and of the configuration minors in matter.
\end{itemize}
\end{theorem}
\paragraph{Integrated stability for three-dimensional elastic data.}
\label{uc:subsec:integrated-3d}

The squared displacement of corresponding material points controls the configuration difference weighted by either density. Together with relative entropy, it covers the case where only one density vanishes.

\begin{theorem}[Integrated weak--strong stability without symmetry]
\label{uc:thm:integrated-3d}
Fix a density-weighted law of Definition~\ref{def:density-weighted-elastic-law} with positive elastic parameter and any $\gamma=1+2\theta>1$. Take as reference a local solution on $\mathbb T^3$ in the class of solutions satisfying the exact minor relations of Theorem~\ref{uf:thm:main} or Proposition~\ref{uf:cor:fractional-existence} at $s=3$. Let it and an entropy weak solution with exact spacetime minors lie in the bounded branch of Proposition~\ref{uf:prop:local-convex-entropy}, with the constitutive segment margins used in Theorem~\ref{uc:thm:common}, including for the values prescribed by the initial data. In normalized material-volume coordinates, the configuration potentials $f_{\mathrm w},f_*$ on the affine cover satisfy
\be
 f_j(t,x+k)=f_j(t,x)+Ak,\qquad k\in\mathbb Z^3,\quad j\in\{\mathrm w,*\}.
\ee
Both potentials satisfy $df_j=F_j$ and have the prescribed spatial gradients on the initial slice. Use the balance laws in the sense of distributions and entropy inequality with the initial terms of Definition~\ref{def:weak-local-complete}, and impose the spatial minor relations~\eqref{eq:physical-weak-minors} there. The zero-density regions may differ. Let $E(t)$ be the spatial integral of the relative entropy density $D_{\eta_{\mathcal K}}$ for almost every $t>0$. Define $E(0)$ by evaluating the same density on the prescribed initial data sets. For positive total particle number, at each time, including $t=0$, almost every material point $y$ has a unique preimage modulo the spatial lattice under each configuration. The lifted positions $X_{\mathrm w},X_*$, after subtraction of $A^{-1}y$, have Lipschitz time representatives with values in $L^\infty$ on the material torus. Their squared displacement is
\be
\label{eq:integrated-position-cost}
L(t)\coloneqq\frac12\int_{\mathbb R^3/A\mathbb Z^3}|X_{\mathrm w}-X_*|^2\,dy.
\ee
Take $L=0$ when the total particle number is zero. The constants below depend on $\gamma$, the bounded comparison branch, the constitutive coefficients, the fixed elastic parameter, the reference norms in the applicable local existence result, and the initial normalization of the potentials.
\begin{itemize}
\item[(i)] The relative entropy integral $E$ has a Lipschitz representative with initial value $E(0)$. There is a Lipschitz corrected integral $\mathcal E$ comparable to $E+L$ and satisfying
\be
\mathcal E(t)\le e^{Ct}\mathcal E(0){,\qquad 0\le t\le T}.
\ee
{The values at zero use the prescribed initial data sets and the normalized initial configurations.}  If $L(0)=0$, the continuous representative of $L$ satisfies $L(t)\le Ct^2e^{Ct}E(0)$ for every $t\in[0,T]$. In particular, if the normalized initial configurations agree, then $\|(r+r_*)^{1/2}(f_{\mathrm w}-f_*)\|_{L^2}\le Ct e^{Ct}\sqrt{E(0)}$ for almost every $t$, with $r={\particlecurrent}^0$ and $r_*={\particlecurrent}_*^0$.
\item[(ii)] In the fixed inertial coordinates, the currents and stress tensors satisfy, for every $1\le p< +\infty$ and almost every $t$,
\be
 \|\particlecurrent(t)-{\particlecurrent}_*(t)\|_{L^p(\mathbb T^3)}{+\|T(t)-T_*(t)\|_{L^p(\mathbb T^3)}}
 \le C_p\bigl(e^{Ct}(E(0)+L(0))\bigr)^{1/\max\{p,2,\gamma\}}.
\ee
The stress contribution $T_{\mathcal L}$ associated with the internal energy $\raug-{m \, n}$ satisfies, for every $1\le p< +\infty$,
\be
 \|T_{\mathcal L}(t)-T_{\mathcal L*}(t)\|_{L^p(\mathbb T^3)}
 \le C_p\bigl(e^{Ct}(E(0)+L(0))\bigr)^{1/\max\{p,2\}}.
\ee
For $1<\gamma\le2$, both $L^2$ bounds (current and total stress) have exponent $1/2$ in the initial error.
\item[(iii)] The configuration gradient and its quadratic minors satisfy
\be
 \int_{\mathbb T^3}(r+r_*)^\gamma
   \bigl(|F-F_*|^2+|G-G_*|^2\bigr)\,dx
 \le Ce^{Ct}\bigl(E(0)+L(0)\bigr),\qquad r={\particlecurrent}^0,\quad r_*={\particlecurrent}_*^0,
\ee
for almost every $t$.
\item[(iv)] If $E(0)+L(0)=0$, the particle current and stress agree almost everywhere, and the configuration minors agree almost everywhere on the common positive-density set. If the prescribed initial densities coincide, then $L(0)=0$ is equivalent to equality of the normalized initial configurations almost everywhere with respect to the common particle measure $r_0\,dx$. In particular, identical prescribed initial data sets give $E(0)=L(0)=0$ with a common normalization.
\end{itemize}
\end{theorem}

{In this exact-minor comparison, the particle-number inequality is an equality by the configuration identity. The general comparison theorem~\ref{thm:physical-entropy-vacuum} also allows nonnegative particle production for augmented weak solutions.}

\paragraph{{Consequences for cavities and initial velocities.}}
{We now identify which physical fields are independent of the initial velocity inside a cavity.}
\begin{corollary}[Independence of the physical evolution from initial data inside an empty cavity]
\label{cor:initial-vacuum-velocity}
Under the periodic constitutive and comparison hypotheses of Theorem~\ref{uc:thm:integrated-3d}, consider two local solutions in the class of solutions satisfying the exact minor relations of Theorem~\ref{uf:thm:main} or Proposition~\ref{uf:cor:fractional-existence} with $s=3$, with the same initial density coordinate $Q_0$. Suppose their normalized initial configurations and initial velocities agree on $\{Q_0>0\}$. On the interval where both solutions lie in the bounded comparison domain, their particle currents and stresses coincide on every time slice. Their velocities, normalized configurations and spacetime minors agree on the closure of the common positive-density region on each slice. The two particle flows agree on the closure of $\{Q_0>0\}$. In particular, if $\mathcal V_0=\{Q_0=0\}$, the common zero-density boundary is the image of $\partial\mathcal V_0$ under $y\mapsto y+tV_0(y)$, with either initial velocity on that boundary.
\end{corollary}
\begin{corollary}[Quantitative comparison of empty regions]
\label{cor:integrated-cavity-comparison}
Under the hypotheses of Theorem~\ref{uc:thm:integrated-3d}, let $r={\particlecurrent}^0$ and $r_*={\particlecurrent}_*^0$.
\begin{itemize}
\item[(i)] \emph{Disagreement of the zero-density sets.} For almost every $t$, one has
\be
\label{eq:integrated-cavity-comparison}
 \int_{\{r=0<r_*\}} r_*^\gamma\,dx + \int_{\{r_*=0<r\}} r^\gamma\,dx
 \le C e^{Ct}\bigl(E(0)+L(0)\bigr).
\ee
All density bounds in this corollary hold for the prescribed initial densities and, when both solutions are {local solutions constructed by the preceding existence theorems}, on every time slice. 
\item[(ii)] \emph{Transport and particle number.} For positive total particle number $M=\det A$, the quadratic transport distance between $r\,dx/M$ and $r_*\,dx/M$ on the spatial torus is at most $\sqrt{2L(t)/M}$. If $D(t)$ denotes the set where exactly one density is positive, its particle number satisfies
\be
 \int_{D(t)} (r+r_*)\,dx
 \le |D(t)|^{1-1/\gamma}\bigl(Ce^{Ct}(E(0)+L(0))\bigr)^{1/\gamma}.
\ee
For these times and any $\delta>0$, the portion of the disagreement set where the nonzero density is at least $\delta$ has volume at most $C\delta^{-\gamma}e^{Ct}(E(0)+L(0))$.

\item[(iii)] \emph{Regions of small relative density.} For every $0\le\lambda<1$ and for the same times,
\be
\label{eq:integrated-low-density-comparison}
 \int_{\{r\le\lambda r_*,\,r_*>0\}} r_*^\gamma\,dx
 +\int_{\{r_*\le\lambda r,\,r>0\}} r^\gamma\,dx
 \le \frac{C e^{Ct}}{(1-\lambda^\theta)^2}\bigl(E(0)+L(0)\bigr),
 \qquad \theta=\frac{\gamma-1}{2}.
\ee
The constant $C$ is independent of $\lambda$.
\item[(iv)] \emph{Density difference on the whole slice.} At the same times one has
\be
\label{eq:integrated-density-power-comparison}
 \|r^{\gamma/2}-r_*^{\gamma/2}\|_{L^2}^2
 \le C e^{Ct}\bigl(E(0)+L(0)\bigr).
\ee
 At the same times, for the bounded comparison domain of Theorem~\ref{uc:thm:integrated-3d},
\be
 \|r-r_*\|_{L^\gamma}^\gamma
 \le C\begin{cases}
 \bigl(e^{Ct}(E(0)+L(0))\bigr)^{\gamma/2}, &1<\gamma<2,\\
 e^{Ct}(E(0)+L(0)), &\gamma\ge2.
 \end{cases}
\ee
\end{itemize}
\end{corollary}

%============================================================================

\section{Proofs of existence and stability with cavitation}
\label{sec:proofs-local-existence}

%-----------------------------------------------------------------------------------------------------------------------------------------

\subsection{The equations and their positive energy}
\label{uf:subsec:normalized-equations}
\label{uf:subsec:positive-energy}

{
We prove Theorem~\ref{uf:thm:main} by constructing the positive energy in Steps~1--3, treating zero density in Steps~4--6 and recovering the physical fields in Steps~7--8.
}

{
{\bf 1.} \emph{Rest-frame matrices.} We set $y_R=(q,z)$ and $\kappa_\theta\coloneqq1+p_\theta$. For the independent fields, Taylor's integral formula gives the symmetric matrices
}
\bse
\be
\label{uf:eq:taylor-factors}
 A_1(q,z)\coloneqq\int_0^1\Psi_{zz}(q,tz)\,dt,\qquad
 B_1(q,z)\coloneqq\int_0^1(1-t)\Psi_{qzz}(q,tz)\,dt,
\ee
\be
 \Psi_z=A_1z,\qquad \Psi_q=z^TB_1z.
\ee
\ese
For $q,\varepsilon_{\rm e}>0$, let $\mathcal J_y$ be the Jacobian of $(n,\Frest,\Grest)\mapsto(q,z)$. The second variation of the rest energy~\eqref{uf:eq:law}, divided by $n$ and expressed in these variables, has matrix $\mathsf H\coloneqq n^{-1}\mathcal J_y^{-T}D^2\bar\rho\,\mathcal J_y^{-1}$. With $Y_R=z/q$ in matter and $Y_R$ prescribed at $q=0$, its entries are
\bse
\be
\label{uf:eq:rest-metric}
 \mathsf H_{qq}=H_f+\Psi_{qq}+\kappa_\theta z^TB_1Y_R-\kappa_\theta Y_R^TA_1Y_R,\qquad
 \mathsf H_{qz}=\Psi_{qz}+\kappa_\theta (A_1Y_R)^T,\qquad \mathsf H_{zz}=\Psi_{zz}.
\ee
{To calculate these blocks, we write $x=(\Frest ,\Grest )$ and differentiate $q=n^\theta$ and $z=\varepsilon_{\rm e}qx$ along an affine variation of $(n,x)$. The first and second derivatives satisfy
\be
 \frac{\dot n}{n}=\frac{\dot q}{\theta q},\qquad
 \ddot q=\frac{\theta-1}{\theta q}\dot q^{\,2},\qquad
 \ddot z=Y_R\ddot q+\frac{2\dot q}{q}(\dot z-Y_R\dot q).
\ee
Applying the chain rule to $n(m+\Phi+\Psi)$ gives the cross block $\Psi_{qz}+\kappa_\theta \Psi_z/q$ and the density block $H_f+\Psi_{qq}+\kappa_\theta \Psi_q/q-\kappa_\theta z\cdot\Psi_z/q^2$. The Taylor factors then give~\eqref{uf:eq:rest-metric}. This calculation is made for $q>0$; the displayed factors give its continuous extension to $q=0$ with bounded $Y_R$.} For the quadratic law~\eqref{uf:eq:exact-polytropic-example}, we write $S\coloneqq\operatorname{diag}(\alpha I,\beta I)$. The corresponding quadratic form is
\be
\label{uf:eq:quadratic-rest-schur}
 \binom{a}{b}^{\!T}\mathsf H\binom{a}{b}
 =(b+\kappa_\theta Y_Ra)^TS(b+\kappa_\theta Y_Ra)+\bigl(H_f-\kappa_\theta(\kappa_\theta+1)Y_R^TSY_R\bigr)a^2.
\ee
The square is positive definite precisely when $H_f>\kappa_\theta(\kappa_\theta+1)Y_R^TSY_R$. For $|Y_R|\le R$ and $H_f\ge\alpha_{\rm f}$, put $d_R\coloneqq\alpha_{\rm f}-\kappa_\theta(\kappa_\theta+1)\lambda_{\max}(S)R^2$. If $d_R>0$, then
\be
 \mathsf H\ge
 \frac{\min\{\lambda_{\min}(S),d_R\}}{(1+\kappa_\theta R)^2}\,I,
\ee
since the inverse of $(a,b)\mapsto(a,b+\kappa_\theta Y_Ra)$ has norm at most $1+\kappa_\theta R$. At $q=0$, $\Psi(q,0)=\Psi_z(q,0)=0$ gives $\Psi_{qq}(0,0)=\Psi_{qz}(0,0)=0$. Hence the same square and positivity criterion apply to every law in the family with $S=\Psi_{zz}(0,0)$ and $H_f=H_f(0)$.

\ese

For spatial vectors $k,b$, we put $t=k\cdot b$ and define $\ell_k(z)b$ as the list with entries $k(z_F^A\cdot b)$ and $z_{G,A}(k\cdot b)-b(z_{G,A}\cdot k)$. The kinematic map is
\bse
\be
\label{uf:eq:kinematic-column}
 D_kb\coloneqq\bigl(\theta q t,\ell_k(z)b+\theta zt\bigr).
\ee
The normalized inertia matrix $\mathsf M_R$ acts on the three velocity components. In matter it is the derivative of spatial momentum with respect to velocity at rest, divided by $n$. In the next formula, $\Psi_{z_F^A}$ and $\Psi_{z_{G,A}}$ are the three-component gradients with respect to the indicated one-form and dual two-form vectors. A dot pairs spatial components, and $x\otimes y$ has entries $x_i y_j$. The formula extends the normalized matrix to $q=0$:
\be
\label{uf:eq:rest-inertia}
 \mathsf M_R=\bigl(h+\Psi+\theta q\Psi_q+\theta z\cdot\Psi_z+\sum_A z_{G,A}\cdot\Psi_{z_{G,A}}\bigr)I+\sum_A\bigl(\operatorname{sym}(z_F^A\otimes\Psi_{z_F^A})-\operatorname{sym}(z_{G,A}\otimes\Psi_{z_{G,A}})\bigr),
\ee
\ese
where $\operatorname{sym}B\coloneqq(B+B^T)/2$.  {For an arbitrary matrix $B$, we consider the simultaneous variations $\delta \Frest=\Frest B$, $\delta n=n\operatorname{tr}B$ and $\delta \Grest=\Grest((\operatorname{tr}B)I-B^T)$. These variations define the material response on independent fields; they do not impose a determinant or cofactor relation between them.} The stress formula~\eqref{eq:useful} and the chain rule therefore give the spatial rest stress
\be
 {T_{\rm sp}=(n\bar\rho_n-\bar\rho+\textstyle\sum_A \Grest _A\cdot\bar\rho_{\Grest _A})I+\sum_A(\Frest ^A\otimes\bar\rho_{\Frest ^A}-\bar\rho_{\Grest _A}\otimes \Grest _A).}
\ee
The velocity derivative of the momentum density is $\bar\rho I+T_{\rm sp}=n\mathsf M_R$. Spatial rotation invariance makes this matrix symmetric. Differentiating the stress in a direction $\delta(n,\Frest ,\Grest )$ gives $L_k^TD^2\bar\rho\delta(n,\Frest ,\Grest )$, apart from
\bse
\be
\label{uf:eq:rest-constraint-terms}
 \sum_A\bigl(\delta \Frest ^A(\bar\rho_{\Frest ^A}\cdot k)-k(\bar\rho_{\Frest ^A}\cdot\delta \Frest ^A)-\bar\rho_{\Grest _A}(\delta \Grest _A\cdot k)\bigr),
\ee
where $L_kb=(n k\cdot b,k(\Frest ^A\cdot b),\Grest _A(k\cdot b)-b(\Grest _A\cdot k))$.  The variation of the scalar coefficient of the identity is $n\delta\bar\rho_n-\sum_A\bar\rho_{\Frest ^A}\cdot\delta \Frest ^A+\sum_A \Grest _A\cdot\delta\bar\rho_{\Grest _A}$, which gives~\eqref{uf:eq:rest-constraint-terms}. {These terms vanish for curl-free one-form variations and divergence-free dual two-form variations, without a cofactor relation.} Equivalently, their nonlinear expressions are $\bar\rho_{\Frest ^A,i}(\del_i \Frest ^A_j-\del_j \Frest ^A_i)-\bar\rho_{\Grest _A,j}\operatorname{div}\Grest _A$.  For $n>0$ and $\varepsilon_{\rm e}>0$, the Jacobian of the change to $(q,z)$ satisfies
\be
 \mathcal J_y L_k=D_k,\qquad
 L_k^T D^2\bar\rho\,\delta(n,\Frest ,\Grest )=nD_k^T\mathsf H\,\delta(q,z).
\ee
The first identity follows by differentiating $q=n^\theta$ and $z=\varepsilon_{\rm e}q(\Frest ,\Grest )$; the second is the definition of $\mathsf H$. The fluid reduction uses the density component alone. At zero density the matrices below are understood through their continuous extensions. Thus the rest equations, on the spatial constraints, have the symmetric matrices
\be
\label{uf:eq:rest-system}
 K_R^0=\operatorname{diag}(\mathsf H,\mathsf M_R),\qquad
 K_R(k)=\begin{pmatrix}0&\mathsf H D_k\\D_k^T\mathsf H&0\end{pmatrix}.
\ee
\ese
{For a perfect fluid the form blocks are omitted. At a point, in an inertial frame instantaneously comoving with the fluid, the normalized density--velocity system reduces to
\be
 H_f\bigl(\mathrm D_tq+\theta q\operatorname{div}V\bigr)=0,\qquad
 h\,\mathrm D_tV+\theta qH_f\nabla q=0.
\ee
The coefficients extend to $q=0$, where the chosen continued evolution gives $\mathrm D_tV=0$. This continuation is part of the strong-solution construction, not a condition imposed on all weak solutions.}

{For $|k|=1$ let $R=|Y_R|$. In~\eqref{uf:eq:rest-metric}, the changes from the fluid density and mixed blocks are $\Ocal(R^2)$ and $\Ocal(R)$, while the form component of $D_k$ is $\Ocal(qR)$. Since $\Psi(q,0)=\Psi_z(q,0)=0$, formula~\eqref{uf:eq:rest-inertia} gives the second matrix bound below. For the speeds in increasing order this yields
\bse
\be
\bigl\|D_k^T\mathsf H D_k-hc_f(q)^2 k\otimes k\bigr\|\le Cq^2R^2,
\qquad \|\mathsf M_R-hI\|\le Cq^2R^2 ,
\ee
\be
 |c_3(k)^2-c_f(q)^2|\le Cq^2R^2,
\qquad 0\le c_1(k)^2\le c_2(k)^2\le Cq^2R^2 .
\ee
\ese
Indeed, the uniform lower bounds on $h$ and $\mathsf M_R$ show that the Rayleigh quotient differs from $c_f(q)^2|k\cdot b|^2/|b|^2$ by at most $Cq^2R^2$. The min--max principle gives the speed estimates, hence squared speeds $\Ocal(q^2)$ at vacuum. A strict fluid acoustic margin persists for small $R$, in particular near vacuum for~\eqref{eq:intro-local-energy}. For a lower bound, the density component of $D_kb$ controls $q|k\cdot b|$. Subtracting $\theta z_F(k\cdot b)$ from the one-form component then controls $q\varepsilon_{\rm e}|\Frest b|$. The bounds on $\mathsf H$ and $\mathsf M_R$ give
\be
\label{uf:eq:acoustic-quotient-comparison}
\frac{b^TD_k^T\mathsf H D_k b}{b^T\mathsf M_R b}
\ge c q^2\frac{|k\cdot b|^2+\varepsilon_{\rm e}^2|\Frest  b|^2}{|b|^2},
\qquad b\ne0.
\ee
In particular $c_1(k)^2\ge c q^2\varepsilon_{\rm e}^2\sigma_{\min}(\Frest )^2$. For exact minors at positive density, $\det\Frest =n>0$, so all three physical speeds are positive at positive elastic parameter.}

For exact minors take $|b|=1$ and $k_\perp=k-(k\cdot b)b$. In each cofactor row the cross product of the row projections onto $b^\perp$ is parallel to $b$. Only terms with a component along $b$ contribute to its product with $k_\perp$, giving
\be
 |(\cof\Frest )k_\perp|\le C|\Frest |\,|\Frest  b|.
\ee
Substitution in $\ell_k(z)b$ gives the reverse of~\eqref{uf:eq:acoustic-quotient-comparison} up to a constant. Thus min--max compares the ordered squared speeds with $q^2$ times the ordered eigenvalues of $k\otimes k+\varepsilon_{\rm e}^2\Frest^T\Frest$. The constants are uniform for $|k|=1$ on the bounded ranges of the gradient and elastic parameter with the assumed matrix bounds, down to zero density. In the fluid case only the density component of $D_k$ remains.

\vskip.2cm \noindent {\bf 2.} \emph{Observer change and the coefficients at zero density.} We express the rest variables by the laboratory fields without division by $Q$. We define $\mathsf T(V)$ with blocks $\Gamma^{-\theta}P_V^{-1}$ on one-form vectors and $\Gamma^{-\theta-1}P_V$ on dual two-form vectors, and set
\bse
\be
\label{uf:eq:rest-lab-map}
 q=\Gamma^{-\theta}Q,\quad z=\mathsf T(V)Z,\quad Y_R=\Gamma^\theta\mathsf T(V)Y,\quad
 P_2\coloneqq Z\otimes Y,\quad P_{2R}\coloneqq\Gamma^\theta(\mathsf T\otimes\mathsf T)P_2.
\ee
\be
\label{uf:eq:boost-jacobian}
 C\coloneqq\Gamma P_V,\qquad \mathsf R\coloneqq\operatorname{diag}(\Gamma^{-\theta},\mathsf T),\qquad
 \mathsf X\coloneqq-D_VC,\qquad k_i\coloneqq P_V^{-1}e_i.
\ee
\ese
In $\mathsf X=-D_VC$, the subscript $V$ means that the covector $k$ in the map $D_k$ is $V$. It is not a velocity derivative. The boost and the rest frame at the base state are kept fixed while taking these variations. At fixed laboratory density coordinate $Q$, we have $\delta q=-\theta q V\cdot C\delta V$.  The corresponding linear map from laboratory variations to rest-frame variations is $\left(\begin{smallmatrix}\mathsf R&\mathsf X\\0&C\end{smallmatrix}\right)$. To see the form components, we vary a laboratory one-form $\Fsp $ and a dual two-form $\Gsp $ while keeping the base boost fixed. Their rest-frame variations are
\be
 \begin{aligned}
 \delta_{\rm fr}\Frest 
 &=P_V^{-1}\delta \Fsp -V(\Frest \cdot C\delta V),\\
 \delta_{\rm fr}\Grest 
 &=\Gamma^{-1}P_V\delta\Gsp-\Grest (V\cdot C\delta V)
 +(C\delta V)(V\cdot\Grest ).
 \end{aligned}
\ee
For the one-form, we apply the fixed boost to its time component $f_t=-\Fsp \cdot V$ and use $P_V-\Gamma V\otimes V=P_V^{-1}$. For the two-form, we represent the laboratory time components by $E=V\times \Gsp $ and apply the fixed boost to $\delta \Gsp $ and $\delta E=\delta V\times \Gsp +V\times\delta \Gsp $. Multiplication by $\varepsilon_{\rm e}q$ and variation of $q$ give the off-diagonal block $-D_VC$ in~\eqref{uf:eq:boost-jacobian}. These identities include $V=0$. The frame is fixed during the variation; differentiating a velocity-dependent choice of rest frame would also introduce a spatial rotation.

The product $\mathsf H D_k$ sends a velocity amplitude $b$ to the density and weighted form components coupled to it in direction $k$. It is smoother than $\mathsf H$. With $L_{k,b}z\coloneqq\ell_k(z)b$, both $A_1L_{k,b}$ and $P_{2R}$ are arrays on the $z$ components. The contraction $(A_1L_{k,b}):P_{2R}$ sums the products of corresponding entries over the two array indices. The two blocks of the extension $K_kb$ are
\be
\label{uf:eq:smooth-H-D}
\begin{aligned}
 (K_kb)_q&=\theta qH_ft+\theta q\Psi_{qq}t+\theta\kappa_\theta\Psi_qt+\Psi_{qz}\cdot\ell_k(z)b+\theta(\Psi_{qz}\cdot z)t+\kappa_\theta(A_1L_{k,b}):P_{2R},\\
 (K_kb)_z&=\theta q\Psi_{zq}t+\theta\kappa_\theta\Psi_z t+\Psi_{zz}\ell_k(z)b+\theta\Psi_{zz}zt.
\end{aligned}
\ee
The terms proportional to $z\cdot\Psi_z/q$ cancel. The remaining quotient $\Psi_z\cdot\ell_k(z)b/q$ is the contraction with $P_{2R}=z\otimes Y_R$, which is $\Ocal(q)$ for bounded $Y_R$. Thus $K_k$ is $C^6$ in $(q,z,P_{2R})$. In $F_{lk}\coloneqq D_l^T K_k$ the density row gains a factor $q$, with $qP_{2R}=z\otimes z$, and the form rows involve only $(q,z)$. Hence $F_{lk}$ is $C^6$ in these two variables.

We use the boost Jacobian to express the symmetric system in the observer variables. The blocks $\mathsf H_L$ and $\mathsf M_L$ multiply the time derivatives of $(Q,Z)$ and $V$, respectively. In direction $i$, $\mathsf B_i$ couples the velocity and density--form equations, while $\mathsf E_i$ acts within the velocity block.  In the base rest frame, we write $dt_R,dx_R$ for the coordinate covectors. The fixed Lorentz boost gives
\be
 dt=\Gamma(dt_R+V\cdot dx_R),\qquad
 dx^i-V_i\,dt=(P_V^{-1}e_i)\cdot dx_R.
\ee
Thus the rest-frame symbol corresponding to $dt$ is $\Gamma(K_R^0+K_R(V))$, and the symbol corresponding to $dx^i-V_i\,dt$ is $K_R(k_i)$. We introduce the matrix
\be
 \mathsf L\coloneqq\begin{pmatrix}\mathsf R&-D_VC\\0&C\end{pmatrix}.
\ee
Block multiplication gives the block-diagonal time matrix
\be
 \mathsf L^T\bigl(K_R^0+K_R(V)\bigr)\mathsf L
 =\operatorname{diag}\bigl(\mathsf R^T\mathsf H\mathsf R,
 C^T(\mathsf M_R-D_V^T\mathsf H D_V)C\bigr).
\ee
The same change of variables applied to $K_R(k_i)$ gives the spatial blocks below. These are algebraic identities also at $q=0$, using the extended rest matrices. The congruence of~\eqref{uf:eq:rest-system} gives
\bse
\be
\label{uf:eq:lab-matrices}
 \mathsf H_L=\Gamma\mathsf R^T\mathsf H\mathsf R,\quad
 \mathsf M_L=\Gamma C^T(\mathsf M_R-F_{VV})C,\quad
 \mathsf B_i=\mathsf R^TK_{k_i}C,\quad
 \mathsf E_i=-C^T(F_{V k_i}+F_{k_i V})C.
\ee
\be
\label{uf:eq:solved-matrices}
 \mathsf Q_i\coloneqq\mathsf M_L^{-1}\mathsf B_i^T,\qquad
 \mathsf R_i\coloneqq\mathsf M_L^{-1}\mathsf E_i.
\ee
\ese
We check uniform positivity first at zero density and then on the full density interval. At $q=0$ one has $D_k=0$ for every $k$ and $\mathsf M_R=mI$. Thus the spatial flux relative to the particle velocity vanishes, and the observer's time matrix has blocks $\Gamma\mathsf R^T\mathsf H\mathsf R$ and $m\Gamma C^TC$. On a compact range of rest shapes where $\mathsf H(0,0,Y_R)\ge h_{\rm v}I>0$, the constitutive formulas give, uniformly,
\be
 \mathsf H(q,qY_R,Y_R)=\mathsf H(0,0,Y_R)+\Ocal(q),\qquad
 \mathsf M_R-F_{VV}=mI+\Ocal(q^2).
\ee
For $q$ sufficiently small and $\Gamma\le\Gamma_{\max}$, it follows that
\be
 \mathsf H_L\ge\frac{h_{\rm v}}{2\Gamma_{\max}^{2\theta+1}}I,
 \qquad \mathsf M_L\ge\frac{m}{2}I .
\ee
Indeed, the least singular values of $\mathsf R$ and $C$ are bounded below by $\Gamma^{-\theta-1}$ and $\Gamma$, respectively.

On the full interval $0\le q\le q_*$ we start at $Y_R=0$, where the constitutive bounds give
\be
 \mathsf H=\operatorname{diag}(H_f,\Psi_{zz}(q,0))\ge\min(\alpha_{\rm f},\alpha_{\rm e})I,
 \qquad \mathsf M_R-F_{VV}=h(I-c_f^2 V\otimes V)\ge h_0(1-v_*^2)I.
\ee
Here $|V|\le v_*<1$ and $c_f^2\le1$. The eigenvalue of the second block along $V$ is $h(1-c_f^2|V|^2)$ and its transverse eigenvalues are $h$. Both bounds persist for sufficiently small $|Y_R|$ by continuity on the compact coefficient ranges. Since $P_V^{-1}$ and $\Gamma^{-1}P_V$ have norm at most one, $|Y_R|\le |Y|$. This gives the uniform upper and positive lower bounds for both time blocks whenever $|Y|\le s_*$, including at $Q=0$.

If $c_f^2\le1-\delta$ throughout the density interval, the same argument at $Y_R=0$ gives $\mathsf M_R-F_{kk}\ge h_0\delta I$ for every unit $k$. Compactness of the unit sphere gives a common positive lower bound for sufficiently small $s_*$. This additional margin gives strictly subluminal elastic wave speeds in matter.

\vskip.2cm \noindent {\bf 3.} \emph{The normalized evolution and spatial products.} We set $y\coloneqq(Q,Z)$ and $\mathrm D_t\coloneqq\del_t+V\cdot\nabla$. The upper spatial coefficients act on $b\in\R^3$ as
\be
\label{uf:eq:Pi}
 (\Pi_i b)_Q=\theta Qb_i,\qquad (\Pi_i b)_{Z_F^A}=e_i(Z_F^A\cdot b)+\theta Z_F^Ab_i,\qquad
 (\Pi_i b)_{Z_{G,A}}=(1+\theta)Z_{G,A}b_i-bZ_{G,A,i}.
\ee
The flux coefficient is $\Pi_i=\Gamma^{-1}\mathsf R^{-1}D_{k_i}C$.  For $b\in\mathbb R^3$, the one-form and two-form blocks of~\eqref{uf:eq:rest-lab-map} give
\be
 k_i\cdot Cb=\Gamma b_i,\qquad
 z_F^A\cdot Cb=\Gamma^{1-\theta}Z_F^A\cdot b,\qquad
 z_{G,A}\cdot k_i=\Gamma^{-1-\theta}Z_{G,A,i}.
\ee
Substitution gives~\eqref{uf:eq:Pi}, including at $Q=0$. The equations are
\bse
\be
\label{uf:eq:normalized-evolution}
 \mathrm D_t y+\sum_i\Pi_i\del_iV=0,\qquad
 \mathrm D_tV+\sum_i(\mathsf Q_i\del_i y+\mathsf R_i\del_iV)=0.
\ee
\be
\label{uf:eq:shape-evolution}
 \mathrm D_tY_F^A=-(DV)^TY_F^A,\qquad
 \mathrm D_tY_{G,A}=(DV)Y_{G,A}-(\operatorname{div}V)Y_{G,A}.
\ee
\ese
They preserve $Q\ge0$ and $Z=QY$ along the flow, including at zero. Their coefficients obey
\bse
\be
\label{uf:eq:symmetry}
 \mathsf H_L\Pi_i=\mathsf B_i,\qquad \mathsf M_L\mathsf Q_i=\mathsf B_i^T,\qquad
 \mathsf M_L\mathsf R_i=\mathsf E_i=\mathsf E_i^T.
\ee
The ${C^{9}}$ constitutive assumptions give $C^6$ dependence of every principal coefficient on $(Q,Z,P_2,V)$; $\Pi_i,\mathsf M_L,\mathsf E_i,\mathsf R_i$ do not need $P_2$. The time matrices are also $C^6$ in $(Q,Z,Y,V)$; only $\mathsf H_L$ depends on the unweighted $Y$. At $Q=0$, the product relation gives $Z=P_2=0$ and hence $\Pi_i=\mathsf Q_i=\mathsf R_i=0$, including when $Y\ne0$. The positive time blocks and the equation $\mathrm D_tV=0$ specify the evolution there. These coefficient formulas contain no $\varepsilon_{\rm e}$.

We write $W_s\coloneqq Y-Y_0$ and use a constant $K$ bounding $\|U\|_{H^3}$, $\|W_s\|_{H^2}$, and the prescribed norm of $Y_0$, for fixed coefficient ranges. The identity
\be
\label{uf:eq:weighted-product}
 d(P_2)_{ij}=Y_j\,dZ_i+Y_i\,dZ_j-Y_iY_j\,dQ,\qquad
 \|P_2\|_{H^3}\le C_K\|U\|_{H^3},
\ee
\ese
follows from $d(Z_iY_j)=Y_j\,dZ_i+Y_i(dZ_j-Y_j\,dQ)$, also at $Q=0$. Multiplication by $Y_0$ preserves $H^2$ on either domain, by the algebra estimate on the torus and the $W^{3,\infty}$ bound on $\mathbb R^3$. The same algebra estimate applies to $Y-Y_0$, while $\|P_2\|_2\le\|Y\|_\infty\|Z\|_2$. The principal coefficients thus belong to $H^3$ relative to their constant values, and are spatially Lipschitz. For the time weights we will use bounded material derivatives, allowing a nonintegrable background $\mathsf H_L(0,Y_0,0)$.

%-----------------------------------------------------------------------------------------------------------------------------------------

\subsection{Construction and estimates at zero density}
\label{uf:subsec:construction-vacuum}

\vskip.2cm \noindent {\bf 4.} \emph{Approximation and construction.} We choose a nonnegative even Schwartz function $\kappa$ with integral one, and, for integers $\ell\ge1$, denote by $\smoothoperator_{\ell}$ convolution with ${\ell}^3\kappa(\ell{}x)$, periodically extended on $\mathbb T^3$. It is self-adjoint on $L^2$. We introduce $w_{\ell}$ and put $V_{\ell}\coloneqq \smoothoperator_{\ell}w_{\ell}$ and $\mathrm D_{\ell}\coloneqq\del_t+V_{\ell}\cdot\nabla$. We use the equations for $y_{\ell},Y_{\ell}$ from~\eqref{uf:eq:normalized-evolution}--\eqref{uf:eq:shape-evolution}, with this velocity, and replace the velocity equation by
\be
\label{uf:eq:approximation}
 \mathrm D_{\ell}w_{\ell}=-\smoothoperator_{\ell}\mathcal F_{\ell},\qquad
 \mathcal F_{\ell}\coloneqq\sum_i(\mathsf Q_i\del_i y_{\ell}+\mathsf R_i\del_iV_{\ell}),\qquad
 w_{\ell}(0)=V_0.
\ee
All coefficients are evaluated at $(Q_{\ell},Z_{\ell},P_{2\ell},V_{\ell})$. We prescribe $y_{\ell}(0)=(Q_0,Q_0Y_0)$ and $Y_{\ell}(0)=Y_0$. In particular the density is not shifted away from zero.

For a fixed $\ell$, we denote by $X_{\ell}$ the flow of $V_{\ell}$, and set $B\coloneqq DX_{\ell}$ and $j\coloneqq\det B$. The upper equations have the explicit solutions
\bse
\be
\label{uf:eq:flow-formulas}
 Q_{\ell}\circ X_{\ell}=j^{-\theta}Q_0,\qquad
 (Z_F)_{\ell}\circ X_{\ell}=j^{-\theta}B^{-T}(Z_F)_0,\qquad
 (Z_G)_{\ell}\circ X_{\ell}=j^{-\theta-1}B(Z_G)_0.
\ee
\be
 (Y_F)_{\ell}\circ X_{\ell}=B^{-T}(Y_F)_0,\qquad
 (Y_G)_{\ell}\circ X_{\ell}=j^{-1}B(Y_G)_0.
\ee
\ese
These formulas transport $\{Q_0=0\}$ and preserve $Q_\ell\ge0$ and $Z_\ell=Q_\ell Y_\ell$. The formula for $P_{2\ell}$ applies the two form transformations to $P_{20}=Q_0Y_0\otimes Y_0$ and multiplies by $j^{-\theta}$; $P_{20},Q_0Y_0\in H^3$ by Step~3. For exact initial minors{ and $\varepsilon_{\rm e}>0$, the determinant and cofactor calculation of Step~8 applies to each smooth flow. }Convolution acts on the velocity across the cavity boundary; its continued equation is recovered on removal of the convolution.

At fixed $\ell$ the unknowns are $X_\ell-\mathrm{Id}\in H^6$ and $\mathcal W_\ell\coloneqq w_\ell\circ X_\ell\in H^3$. We keep $X_\ell$ in a small $C^1$ neighborhood of the identity and write $\mathcal A=X^{-1}$. Subtracting the two inverse identities gives
\bse
\be
\label{uf:eq:inverse-flow-difference}
 M(\mathcal A_1-\mathcal A_2)=-(X_1-X_2)\circ\mathcal A_2,\qquad
 M\coloneqq\int_0^1DX_1(\mathcal A_2+\tau(\mathcal A_1-\mathcal A_2))\,d\tau.
\ee
Differentiation of $X\circ\mathcal A=\mathrm{Id}$ bounds $\mathcal A-\mathrm{Id}$ in $H^6$. The intermediate maps remain bi-Lipschitz, so $M-I$ is bounded in $H^5$ and small in $L^\infty$. Multiplication by $M^{-1}$ bounds $\|\mathcal A_1-\mathcal A_2\|_{H^5}$ by $C\|X_1-X_2\|_{H^6}$. Integrating $Dg$ along the same segment gives
\be
\label{uf:eq:composition}
 \|g\circ\mathcal A\|_{H^3}\le C\|g\|_{H^3},\qquad
 \|g\circ\mathcal A_1-g\circ\mathcal A_2\|_{H^2}\le C\|g\|_{H^3}\|X_1-X_2\|_{H^6}.
\ee
Thus the reconstructed $y_\ell,P_{2\ell},w_\ell$ are bounded in $H^3$ with locally Lipschitz differences in $H^2$, since the factors formed from $DX_\ell$ and its determinant obey the same estimates. The force $\mathcal F_\ell$ is bounded in $H^2$ with differences in $H^1$. Fixed-scale convolution supplies $H^7$ for the velocity and $H^4$ for the force, allowing one derivative for composition by $X_\ell$. The equations
\be
\label{uf:eq:label-ODE}
 (X_{\ell}-\mathrm{Id},\mathcal W_{\ell})'=
 \bigl((\smoothoperator_{\ell}w_{\ell})\circ X_{\ell},-(\smoothoperator_{\ell}\mathcal F_{\ell})\circ X_{\ell}\bigr)
\ee
\ese
have a locally Lipschitz right-hand side in $H^6\times H^3$ and hence a unique local solution. At fixed $\ell$, a bound on $\|w_\ell\|_{H^3}$ bounds the higher derivatives of $V_\ell$ in the flow equation and hence $\|X_\ell-\mathrm{Id}\|_{H^6}$. While the $C^1$ and constitutive margins persist, the right-hand side and its local Lipschitz constants stay bounded. The solution therefore has a limit at a finite endpoint and can be restarted there.

\vskip.2cm \noindent {\bf 5.} \emph{The energy estimate without a positive lower density.} For $|\alpha|\le3$, we put $a_\alpha\coloneqq\del^\alpha y_{\ell}$, $b_\alpha\coloneqq\del^\alpha w_{\ell}$, and $c_\alpha\coloneqq \smoothoperator_{\ell}b_\alpha$. Their principal equations are
\bse
\be
\label{uf:eq:differentiated-approximation}
 \mathrm D_{\ell}a_\alpha+\sum_i\Pi_i\del_ic_\alpha=f_\alpha,\qquad
 \mathrm D_{\ell}b_\alpha+\smoothoperator_{\ell}\sum_i(\mathsf Q_i\del_ia_\alpha+\mathsf R_i\del_ic_\alpha)=g_\alpha.
\ee
The product estimate~\eqref{uf:eq:weighted-product} gives $H^3$ bounds for all differentiated principal coefficients. The commutator terms satisfy $\|f_\alpha\|_2+\|g_\alpha\|_2\le C_K\|(y_{\ell},w_{\ell})\|_{H^3}$. In the terms with three derivatives on a coefficient, its $L^2$ derivative multiplies a first field derivative in $L^\infty$; the intermediate products use $L^6/L^3$. These are the three-dimensional Sobolev products at order three.

We use the positive quadratic energy
\be
\label{uf:eq:energy}
 \mathcal E_3\coloneqq\sum_{|\alpha|\le3}\int_{\mathcal D}
 \bigl(a_\alpha^T\mathsf H_La_\alpha+b_\alpha^T\mathsf M_Lb_\alpha\bigr)\,dx,
 \qquad c_0\|(y_{\ell},w_{\ell})\|_{H^3}^2\le\mathcal E_3\le c_1\|(y_{\ell},w_{\ell})\|_{H^3}^2.
\ee
\ese
The constants $c_0,c_1>0$ depend on the fixed ranges and small shape bound. Self-adjointness of $\smoothoperator_{\ell}$ and the identities~\eqref{uf:eq:symmetry} combine the principal integrals into
\bse
\be
\label{uf:eq:paired-fluxes}
 \int(a_\alpha^T\mathsf B_i\del_ic_\alpha+c_\alpha^T\mathsf B_i^T\del_ia_\alpha)
 =-\int a_\alpha^T(\del_i\mathsf B_i)c_\alpha,\qquad
 \int c_\alpha^T\mathsf E_i\del_ic_\alpha=-\frac12\int c_\alpha^T(\del_i\mathsf E_i)c_\alpha.
\ee
\be
\label{uf:eq:metric-commutator}
 \|[\smoothoperator_{\ell},\mathsf M_L]b_\alpha\|_{H^1}\le C\|D\mathsf M_L\|_\infty\|b_\alpha\|_2,\qquad
 \Big\|\sum_i(\mathsf Q_i\del_ia_\alpha+\mathsf R_i\del_ic_\alpha)\Big\|_{H^{-1}}
 \le C_K(\|a_\alpha\|_2+\|b_\alpha\|_2).
\ee
\ese
The extra integral from commuting $\smoothoperator_{\ell}$ with $\mathsf M_L$ is bounded by the product of the two norms in the estimates~\eqref{uf:eq:metric-commutator}. The second estimate follows by integration by parts with Lipschitz coefficients. For the convolution kernel $K_{\ell}$ and a bounded Lipschitz matrix $M$,
\be
 [\smoothoperator_{\ell},M]b(x)=\int K_{\ell}(x-y)(M(y)-M(x))b(y)\,dy.
\ee
Differentiating in $x$ gives a term with $DK_{\ell}(x-y)(M(y)-M(x))$ and a term with $K_{\ell}(x-y)DM(x)$. Young's inequality bounds their $L^2$ norms by $C\|DM\|_\infty\|b\|_2$, since $\int |z||DK_{\ell}(z)|\,dz$ is independent of $\ell$. The undifferentiated term has the additional factor ${\ell}^{-1}$. This proves the $H^1$ bound in~\eqref{uf:eq:metric-commutator} for $\ell\ge1$, without an $L^2$ assumption on $M$ or its background.

The combined derivatives of the weights are also bounded uniformly. The physical velocity obeys
\bse
\be
\label{uf:eq:approximate-material-derivative}
 \mathrm D_{\ell}V_{\ell}=-\smoothoperator_{\ell}^2\mathcal F_{\ell}+[V_{\ell}\cdot\nabla,\smoothoperator_{\ell}]w_{\ell},
 \qquad \|\mathrm D_{\ell}V_{\ell}\|_{H^2}+\|\mathrm D_{\ell}y_{\ell}\|_{H^2}+\|\mathrm D_{\ell}Y_{\ell}\|_{H^2}\le C_K.
\ee
The shape derivative uses $DV_{\ell}\in H^2$ and multiplication by $Y_0+W_s$. Thus $\mathrm D_{\ell}\mathsf H_L$ and $\mathrm D_{\ell}\mathsf M_L$ are bounded in $L^\infty$. Indeed, for either weight $S=S(y_{\ell},Y_{\ell},V_{\ell})$ the chain rule along the approximate flow gives
\be
 \mathrm D_{\ell}S=S_y\mathrm D_{\ell}y_{\ell}+S_Y\mathrm D_{\ell}Y_{\ell}+S_V\mathrm D_{\ell}V_{\ell},
 \qquad \mathrm D_{\ell}Y_{\ell}=L(DV_{\ell})Y_{\ell}.
\ee
The derivatives of $S$ are bounded on the prescribed coefficient range, $H^2\subset L^\infty$, and $Y_{\ell}$ is bounded. These formulas control the time derivative of the energy weight without differentiating $Y_{\ell}$ in space or dividing by the density. Transport in the energy introduces only these derivatives and $\operatorname{div}V_{\ell}$. In the spatial integrals~\eqref{uf:eq:paired-fluxes}, we differentiate $\mathsf B_i$ and $\mathsf E_i$. By~\eqref{uf:eq:symmetry} and symmetry of $\mathsf M_L$, one has $\mathsf B_i=\mathsf Q_i^T\mathsf M_L$. Both factors have bounded spatial derivatives by~\eqref{uf:eq:solved-matrices} and~\eqref{uf:eq:weighted-product}, as does $\mathsf E_i$. These bounds hold on the entire density range, including $Q=0$. We conclude
\be
\label{uf:eq:energy-bound}
 \frac{d}{dt}\mathcal E_3\le C_K\mathcal E_3.
\ee
\ese

This calculation applies at the Sobolev regularity of the theorem. To justify it, we convolve $a_\alpha,b_\alpha$ with a second kernel $\smoothoperator_\delta$, independently of $\ell$. For a Lipschitz coefficient $A$, the integral formula for $[\smoothoperator_\delta,A]\del_i u$ bounds its $L^2$ norm by $C\|DA\|_\infty\|u\|_2$ and gives convergence to zero for fixed $u\in L^2$, first for smooth $u$ and then by density. These are the commutators in the equations solved for the time derivatives. The weights are multiplied afterward. Locally they have $H^2$ regularity, and their pullbacks by the Lipschitz flow have bounded time derivatives. For a smoothed field $u$ and either weight $S$, the flow formula and its Jacobian give
\be
\label{uf:eq:transport-weighted-energy}
 \frac{d}{dt}\int u^TSu\,dx
 =2\int u^TS\mathrm D_{\ell}u\,dx
 +\int u^T\bigl(\mathrm D_{\ell}S+(\operatorname{div}V_{\ell})S\bigr)u\,dx.
\ee
Thus the transport integration uses $\mathrm D_{\ell}S$ without requiring $DS\in L^\infty$. We combine the fluxes before taking $\delta\to0$. The additional integral controlled by the estimates~\eqref{uf:eq:metric-commutator} passes by $H^1$--$H^{-1}$ duality. On $\R^3$ a cutoff at infinity leaves an error bounded by its gradient times the finite energy, which tends to zero. This proves the estimate~\eqref{uf:eq:energy-bound} for the approximations just constructed.

For the increment $W_s=Y_{\ell}-Y_0$, the one-form and two-form actions $L$ in~\eqref{uf:eq:shape-evolution} give $\mathrm D_{\ell}W_s=L(DV_{\ell})W_s+L(DV_{\ell})Y_0-V_{\ell}\cdot\nabla Y_0$. The term linear in $W_s$ is bounded in $H^2$ by $C\|V_{\ell}\|_{H^3}\|W_s\|_{H^2}$, by the algebra property of $H^2$. This property also bounds the two terms involving $Y_0$ on $\mathbb T^3$, since $Y_0\in H^3$. On $\R^3$, multiplication by $Y_0$ and by $\nabla Y_0$ is bounded on $H^2$ by the prescribed $W^{3,\infty}$ norm. The commutator of two spatial derivatives with transport is bounded by $C\|V_{\ell}\|_{H^3}\|W_s\|_{H^2}$, using $L^\infty$ for $DV_{\ell}$ and $L^6L^3$ for the product of $D^2V_{\ell}$ and $DW_s$. This gives the $H^2$ transport estimate with $V_{\ell}\in H^3$ on both domains. We choose a bound $K$ for the field norms strictly larger than the initial bound after the constants $c_0,c_1$ in the energy bounds~\eqref{uf:eq:energy} have been accounted for, and choose a positive bound for $\|W_s\|_{H^2}$. On this set the estimates imply
\be
\label{uf:eq:bootstrap}
 \mathcal E_3(t)\le\mathcal E_3(0)e^{C_Kt},\qquad
 \|W_s(t)\|_{H^2}\le C_Kt e^{C_Kt},\qquad
 \|Y_{\ell}(t)\|_\infty\le s_{\rm in}e^{C_Kt}.
\ee
The positive kernel gives $|V_{\ell}(0)|\le v_0$. The equations bound $\del_t y_{\ell}$ and $\del_tV_{\ell}$ in $H^2$ uniformly, while $DV_{\ell}$ is bounded in $L^\infty$. We choose $T>0$ so that the three estimates in~\eqref{uf:eq:bootstrap} stay strictly inside the selected norm bounds, the shape bound stays below $s_{\rm out}$, the density stays below $Q_*$, the velocity stays below $v_*$, and the flow stays in the $C^1$ neighborhood of the identity. These are finitely many strict inequalities, satisfied for a sufficiently short time by the initial margins. The first-exit argument and the fixed-$\ell$ continuation criterion then give solutions on this same $[0,T]$. The constants involve neither $\ell$, $\varepsilon_{\rm e}$ nor a positive density minimum.

\vskip.2cm \noindent {\bf 6.} \emph{Passage to the limit and time continuity.} The equation for $V_{\ell}$ has the error
\bse
\be
\label{uf:eq:approximation-error}
 e_{\ell}\coloneqq(I-\smoothoperator_{\ell}^2)\mathcal F_{\ell}+[V_{\ell}\cdot\nabla,\smoothoperator_{\ell}]w_{\ell},
 \qquad \|e_{\ell}\|_{H^1}\le C_K{\ell}^{-1},\qquad \|e_{\ell}\|_{H^2}\le C_K.
\ee
For the first term use $\mathcal F_{\ell}\in H^2$ and the first kernel moment.  The even kernel gives $\|[V_{\ell}\cdot\nabla,\smoothoperator_{\ell}]w_{\ell}\|_{H^1}\le C{\ell}^{-2}\|V_{\ell}\|_{H^3}\|w_{\ell}\|_{H^3}$. Indeed, we write $\delta_z f(x)\coloneqq f(x)-f(x-z)$ and $\kappa_{\ell}(z)\coloneqq {\ell}^3\kappa(\ell{}z)$. Subtracting $Dw_{\ell}(x)$ in the commutator gives
\be
 [V_{\ell}\cdot\nabla,\smoothoperator_{\ell}]w_{\ell}
 =(V_{\ell}-\smoothoperator_{\ell}V_{\ell})\cdot Dw_{\ell}
 -\int\kappa_{\ell}(z)\,\delta_zV_{\ell}\cdot\delta_zDw_{\ell}\,dz.
\ee
The zero first moment gives $\|V_{\ell}-\smoothoperator_{\ell}V_{\ell}\|_{H^1}\le C{\ell}^{-2}\|V_{\ell}\|_{H^3}$. The two increments cost $|z|$ each in $H^2$ and $H^1$, respectively. The $H^2\times H^1$ product estimate and the second kernel moment bound the integral. The term $(I-\smoothoperator_{\ell}^2)\mathcal F_{\ell}$ contributes $\Ocal({\ell}^{-1})$ in $H^1$. Evenness also gives, for $0\le s\le3$,
\be
 \|w_{\ell}-V_{\ell}\|_{H^s}\le C {\ell}^{-\min\{2,\,3-s\}}\|w_{\ell}\|_{H^3}.
\ee
Indeed the Fourier multiplier $1-\widehat\kappa(\xi/\ell)$ is bounded by $C\min\{1,|\xi/\ell|^2\}$; the same estimate holds for periodic frequencies. In particular the two velocities differ by $\Ocal({\ell}^{-2})$ in $H^1$, without an additional derivative of the initial velocity. The $H^2$ bound follows by estimating the two transport expressions separately.

\ese

{
On the cone $Q\ge0$, $|Z|\le C_KQ$, the map $P_2=Z\otimes Z/Q$ is Lipschitz, with value zero at the vertex. Indeed, $|D_ZP_2|\le2C_K$ and $|\partial_QP_2|\le C_K^2$ away from the vertex. Subtraction in~\eqref{uf:eq:weighted-product} also gives
\bse
\be
\label{uf:eq:product-difference}
 \|\delta P_2\|_{H^1}\le C_K(\|\delta U\|_{H^1}+\|\delta Y\|_2).
\ee
Here the derivative difference is bounded by $C_K(\|D\delta U\|_2+\|D(Q_2,Z_2)\|_\infty\|\delta Y\|_2)$. The same $H^1$ estimate holds for differences of the principal coefficients and their products with the reference derivatives, since $H^2$ multiplies $H^1$. Pairing the fluxes by~\eqref{uf:eq:paired-fluxes} and using spatial convolution to justify the differentiated equations gives the lower-order difference estimate
\be
\label{uf:eq:cauchy}
 \frac{d}{dt}\mathcal E_{1,\ell M}
 \le C_K\mathcal E_{1,\ell M}+C_K\mathcal E_{1,\ell M}^{1/2}
 \bigl(\|Y_{\ell}-Y_M\|_2+\|e_{\ell}-e_M\|_{H^1}\bigr).
\ee
The transport equation for the unweighted forms gives
\be
 \frac{d}{dt}\|Y_{\ell}-Y_M\|_2
 \le C_K\bigl(\|Y_{\ell}-Y_M\|_2+\|V_{\ell}-V_M\|_{H^1}\bigr).
\ee
In the term $(V_\ell-V_M)\nabla Y_M$, we split $Y_M=Y_0+W_M$. We use $L^2\times L^\infty$ for $\nabla Y_0$ on $\mathbb R^3$ and $L^6\times L^3$ for $\nabla Y_0$ on the torus and for $\nabla W_M$ on both domains. The term with $D(V_\ell-V_M)$ uses $\|Y_M\|_\infty$. These inequalities also prove uniqueness with zero approximation error and equal data.
\ese

For convergence at a rate, the second kernel moment and the preceding commutator estimate give $\|e_\ell\|_2\le C_K\ell^{-2}$ and $\|V_{\ell}(0)-V_M(0)\|_2\le C_K(\ell^{-2}+M^{-2})$. At order zero the Lipschitz coefficient bound costs only $\|U_\ell-U_M\|_2$, since the reference derivatives are bounded. The paired energy and Gronwall then give
}
\be
 \|U_{\ell}-U_M\|_{C_tL^2}+\|P_{2\ell}-P_{2M}\|_{C_tL^2}
 \le C_{K,T}({\ell}^{-2}+M^{-2}).
\ee
Interpolating with the uniform $H^3$ bounds gives, for $0\le s<3$,
\be
 \|U_{\ell}-U_M\|_{C_tH^s}+\|P_{2\ell}-P_{2M}\|_{C_tH^s}
 \le C_{K,T}({\ell}^{-2}+M^{-2})^{1-s/3}.
\ee
The estimate for $U_{\ell}-U_M$ at $s=1$ and the shape transport inequality above give $\|Y_{\ell}-Y_M\|_{C_tL^2}\le C_{K,T}({\ell}^{-2}+M^{-2})^{2/3}$, since the initial shapes agree. Interpolation with the uniform $H^2$ bounds for $Y_{\ell}-Y_0$ yields, for $0\le r<2$,
\be
 \|Y_{\ell}-Y_M\|_{C_tH^r}\le C_{K,T}({\ell}^{-2}+M^{-2})^{(2-r)/3}.
\ee
Letting $M\to\infty$ gives these rates for the limiting fields; the estimate on $w_\ell-V_\ell$ gives the same $C_tH^s$ rate $\Ocal(\ell^{-2(1-s/3)})$ for $w_\ell$. With $s>5/2$ and $r>3/2$, the velocities converge in $C_tC^1$ and the unweighted forms uniformly. Thus $Q\ge0$ and $Z=QY$ on every slice, even for nonintegrable $Y_0$ on $\mathbb R^3$. All coefficients converge and~\eqref{uf:eq:approximation-error} tends to zero, giving~\eqref{uf:eq:normalized-evolution}--\eqref{uf:eq:shape-evolution}. The flow satisfies
\be
 \|X_{\ell}-X\|_{C([0,T];W^{1,\infty})}
 +\|Y_{\ell}-Y\|_{C_tL^\infty}\le C {\ell}^{-1/3}.
\ee
Indeed, splitting the Fourier integral or series at a frequency $R\ge1$ gives
\be
 \|f\|_\infty\le C\|f\|_2^{1/2}\|f\|_{H^3}^{1/2},\qquad
 \|Df\|_\infty\le C\|f\|_2^{1/6}\|f\|_{H^3}^{5/6}.
\ee
{
For the derivative bound, the low and high frequency contributions are $CR^{5/2}\|f\|_2$ and $CR^{-1/2}\|f\|_{H^3}$, respectively. Applying the two inequalities to $V_\ell-V$ gives errors $\Ocal(\ell^{-1})$ in $L^\infty$ and $\Ocal(\ell^{-1/3})$ for the derivative. Gronwall gives the same rates for $X_\ell-X$ and its derivative, using the uniform $1/2$-H\"older bound for $DV$ when comparing the flows. The inverse derivative formula and the reconstruction with Lipschitz $Y_0$ give the asserted bound on $Y_\ell-Y$. The same Fourier estimate for $Q_\ell-Q$ and the power map $r=Q^{p_\theta}$ give the current rate $\Ocal(\ell^{-\min\{1,p_\theta\}})$. With exact minors, $r_\ell=\det D_xf_\ell$ also gives a rate $\Ocal(\ell^{-1/3})$, so one may take the better of the two rates in this case.
}

{
The limit is weakly continuous in $H^3$ and the shape increment $W_s$ is weakly continuous in $H^2$. For a third derivative $u$ of $U$ the solved equations give $\partial_tu\in L^\infty_t H^{-1}$, so $\smoothoperator_\delta u$ is strongly continuous in $L^2$ at fixed $\delta$. Apply the grouped energy calculation to these smoothed fields. For each spatially Lipschitz solved matrix $A$, the commutator $[\smoothoperator_\delta,A]\partial_i u$ is bounded in $L^2$ by $C\|\nabla A\|_\infty\|u\|_2$ and tends to zero by density, as in Step~5. This gives convergence in $L_t^1L^2$. The energy uses the bounded material derivatives of the time weights. At each endpoint $\smoothoperator_\delta u(t_j)$ converges strongly to the weakly continuous representative $u(t_j)$ in $L^2$. The resulting energy identity between arbitrary times gives absolute continuity of the weighted $H^3$ energy.
}

The time weights are continuous in $L^\infty$ by the lower-order strong convergence and Sobolev embedding ($s>5/2$, $r>3/2$). Freezing the weight $S$ at $t_0$ changes the squared energy {by at most $\|S(t)-S(t_0)\|_\infty\|U(t)\|_{H^3}^2$}. Its continuity and weak $H^3$ continuity therefore give strong continuity in the fixed positive norm at $t_0$, including at zero. Transport gives $W_s\in C_tH^2$. In~\eqref{uf:eq:weighted-product}, $Y=Y_0+W_s$ is a continuous $H^2$ multiplier, hence $DP_2\in C_tH^2$ and $P_2=Z\otimes Y\in C_tL^2$. Thus $P_2\in C_tH^3$, and the equations give $U_t\in C_tH^2$, $(W_s)_t\in C_tH^1$. Estimate~\eqref{uf:eq:cauchy} with zero approximation error gives uniqueness.

%-----------------------------------------------------------------------------------------------------------------------------------------

\subsection{Dependence on the data and physical reconstruction}
\label{uf:subsec:dependence-reconstruction}

\vskip.2cm \noindent {\bf 7.} \emph{Continuous dependence on the initial data.} Suppose $U_0$ converges in $H^3$ and $Y_0$ in $L^\infty$ under the common bounds and margins of Theorem~\ref{uf:thm:main}. The difference estimates of Section~\ref{sec:non-relativistic-vacuum} apply with both parameters fixed at the relativistic value, using the positive time matrices of Step~2. On $\mathbb T^3$, interpolation gives convergence of $Y_0$ in $H^2$, {so the periodic difference estimate in the proof of Proposition~\ref{nr:prop:periodic-sobolev-data}, with both parameters fixed, gives} convergence of $U$ in $C_tH^3$ and $Y$ in $C_tH^2$. On $\mathbb R^3$, finite differences and the uniform $W^{3,\infty}$ bound give convergence of $Y_0$ in $W^{2,\infty}$, including on the initial vacuum. {At those same parameter values the estimates in the proof of Proposition~\ref{nr:prop:two-bounded-derivatives} give} convergence of $U$ and the inverse displacement in $C_tH^3$, and of the changes from the transported backgrounds in $C_tH^2$. Step~4 of the proof of Theorem~\ref{nr:thm:vacuum-limit} gives the latter convergence for the changes from the fixed initial fields. The equations give convergence of the respective time derivatives one order lower.

{\vskip.2cm \noindent {\bf 8.} \emph{Independent fields and exact minor relations.} Fix $\varepsilon_{\rm e}>0$. For the independent initial matrices the transport equations give
\bse
\be
\label{uf:eq:independent-flow-forms}
\Fsp(t,X(t,y))=\Fsp_0(y)(D_yX)^{-1},\qquad
\Gsp(t,X(t,y))=\Gsp_0(y)\operatorname{cof}((D_yX)^{-1}),
\ee
\be
\label{uf:eq:independent-flow-density}
r(t,X(t,y))=r_0(y)/\det D_yX(t,y),\qquad r=Q^{p_\theta}.
\ee
\ese
The invertible flow matrices preserve both form ranks, including at zero density. The formulas are pullbacks of spatial one- and two-forms, so they preserve spatial closure. Their spacetime extensions annihilating $(1,V)$ satisfy $dF=dG=0$ by transport. These identities follow for the $C^1$ flow by smooth approximation with uniform bounds on the first derivatives and the inverse Jacobian. At $t=0$ the flow is the identity.

The normalized momentum equation gives the three spatial rest-frame projections of stress conservation on $r>0$. The fourth projection follows from the constitutive identity
\be
\label{uf:eq:independent-stress-projection}
u_\mu\partial_\nu T^{\mu\nu}=-H\partial_\nu\particlecurrent^\nu,
\qquad H=\partial_n\bar\rho,
\ee
for closed forms annihilating $u$. Indeed, contracting the derivative of~\eqref{augment3} with $u$ and using $\mathcal L_uF=\mathcal L_uG=0$ cancels the form derivatives against the spatial stress by the chain rule. The remaining term is $-H(u\cdot\partial n+n\operatorname{div}u)$. In~\eqref{uf:eq:independent-flow-density}, change variables $x=X(t,y)$ and integrate $\partial_t[\varphi(t,X(t,y))]$ against $r_0(y)\,dy$. This gives $\operatorname{div}\particlecurrent=0$ with initial density $r_0$, also at vacuum. Hence all four stress equations hold in matter for independent forms. Their extension through the cavity follows from~\eqref{uf:eq:stress-decomposition} below.}

{
For exact initial minors set $f=f_0\circ X^{-1}$. Then $f_t=-D_xf\,V$ and~\eqref{uf:eq:independent-flow-forms} gives $\Fsp=D_xf$. Multiplicativity of the cofactor and determinant gives $\Gsp=\cof\Fsp$ and $\det\Fsp=r$, including at singular gradients. Orthogonality to $(1,V)$ determines the time components. For the singular-value bounds, put
}
\be
 K(t)\coloneqq\int_0^t\|DV(s)\|_{L^\infty,\mathrm{op}}\,ds,
\ee
where the norm is the essential supremum of the Euclidean operator norm. For $j=1,2,3$ and each initial position $y$, with singular values in decreasing order,
\be
 \begin{aligned}
 e^{-K(t)}\sigma_j(\Fsp _0(y))&\le\sigma_j(\Fsp (t,X(t,y)))\le e^{K(t)}\sigma_j(\Fsp _0(y)),\\
 e^{-2K(t)}\sigma_j(\cof \Fsp _0(y))&\le\sigma_j((\cof \Fsp )(t,X(t,y)))\le e^{2K(t)}\sigma_j(\cof \Fsp _0(y)).
 \end{aligned}
\ee
The flow bounds place the singular values of $DX$ and its inverse in $[e^{-K(t)},e^{K(t)}]$. Those of $\cof((DX)^{-1})$ are pairwise products. The min--max formula for products with these invertible matrices gives both estimates, even for a singular configuration gradient. The spacetime gradient acts by $F(\tau,w)=\Fsp(w-\tau V)$. The map $(\tau,w)\mapsto w-\tau V$ and its second exterior power are both surjective. Thus $F$ has rank $k=\operatorname{rank}\Fsp$ and its quadratic minors have rank $\binom{k}{2}$ (cf.~\eqref{eq:rank-and-current}). At fixed $\varepsilon_{\rm e}>0$ the gradient and cofactor belong locally to $C_tH^2\cap C_t^1H^1$. The identities $r=\det\Fsp$, $n=r/\Gamma$ and $\particlecurrent=r(1,V)$ give $n,\particlecurrent$ the same regularity by multiplication.

We identify the material fluxes at zero determinant as well as in matter. We define the specific Lagrangian
\bse
\be
\label{uf:eq:specific-Lagrangian}
 \ell(\Fsp ,q,V)\coloneqq\Gamma^{-1}\bigl(m+\Phi(q)+\Psi(q,\varepsilon_{\rm e}q(\Fsp P_V^{-1},\Gamma^{-1}(\cof \Fsp )P_V))\bigr),\qquad
 B_\ell\coloneqq\ell_V-\theta q\Gamma^2V\ell_q.
\ee
\be
\label{uf:eq:material-momenta}
 P^0\coloneqq-\cof(\Fsp )B_\ell,\qquad
 K_\ell\coloneqq(\ell+\theta q\ell_q)\cof(\Fsp )+r\ell_{\Fsp},\qquad
 P^{\rm sp}\coloneqq K_\ell+P^0\otimes V.
\ee
\ese
Each derivative of $\ell$ keeps $\Fsp ,q,V$ independent. At positive determinant, varying $f_t$ gives $\delta V=-\Fsp ^{-1}\delta f_t$ and $\delta q=-\theta q\Gamma^2V\cdot\delta V$; varying $\Fsp $ with $f_t$ fixed gives $\delta V=-\Fsp ^{-1}(\delta \Fsp )V$ and $\delta q=\theta q\delta r/r-\theta q\Gamma^2V\cdot\delta V$. Substitution in $\delta(r\ell)$ gives the material momenta in~\eqref{uf:eq:material-momenta}, using $r(\Fsp)^{-T}=\cof \Fsp $. These formulas extend without division to $r=0$.

We write $\ell=m/\Gamma+\ell_s$. Taylor's formula gives $\ell_s=q^2L(\Fsp ,q,V)$ with $L$ of class $C^6$ for fixed $\varepsilon_{\rm e}$ on compact ranges of $(\Fsp ,q,V)$. Thus the nonmass terms of $P^0,K_\ell$ are $q^2$ times $C^5$ functions of $(\Fsp ,q,V)$, except for a term containing
\bse
\be
\label{uf:eq:material-power}
 rq^2=Q^{p_\theta+2}\Gamma^{-2\theta},\qquad p_\theta+2>2.
\ee
They have spatial $H^2$ regularity and time derivatives locally in $H^1$: the power of $Q$ in~\eqref{uf:eq:material-power} has a $C^2$ extension through zero, and its time derivative contains $Q^{p_\theta+1}$, a spatially Lipschitz multiplier. Every component of these nonmass fluxes vanishes on $\{Q=0\}$. Locality of weak derivatives of an $H^1$ function then makes all their first derivatives zero almost everywhere on that set.

The mass terms are $P_m^0=m\Gamma\cof(\Fsp )V$ and $K_m=m\Gamma^{-1}\cof(\Fsp )$. The Piola identity and cofactor transport give
\be
\label{uf:eq:mass-material-equation}
 \del_tP_m^0+\operatorname{div}(P_m^0\otimes V+K_m)
 =m\cof(\Fsp )\mathrm D_t(\Gamma V).
\ee
\ese
The cofactor equation $\mathrm D_t\Gsp_A=(DV)\Gsp_A-\Gsp_A\operatorname{div}V$ cancels the transport-divergence terms. The Piola identity, proved by smoothing the $H^2$ spatial gradient and passing to the cofactor in $H^1$, leaves the term canceled by $\nabla\Gamma^{-1}=-\Gamma(DV)^TV$. At $Q=0$, $q=z=P_2=0$ and $\mathsf Q_i=\mathsf R_i=0$, so $\mathrm D_tV=0$ and~\eqref{uf:eq:mass-material-equation} vanishes. The scalar flow formula transports $\mathcal V_0$. Integration of the $C^1$ velocity along its trajectories yields~\eqref{uf:eq:inertial-vacuum-flow}. On $\Omega_0\subset\mathcal V_0$ it gives $DX=I+t\,DV_0$, and the orientation-preserving diffeomorphism gives the boundary and volume formulas.

On compact subsets of matter, $F$ is continuous and locally $W^{1,2}$ in spacetime. Smoothing $f$ in $W^{2,2}$ with uniform convergence of gradients gives the Noether identity $\del_\nu T_\mu{}^\nu=F_\mu^A\mathcal E_A$, with $\mathcal E_A=\del_tP_A^0+\del_iP_A^i$ in $L^2$. The normalized momentum equation gives the spatial projection of stress balance in the rest frame. The restriction of $F$ to the rest space is invertible, so $\mathcal E_A=0$ in matter. Its vanishing on the zero set and the locally integrable first weak derivatives of the material fluxes then give the material equations across the interface.

The flow also propagates any prescribed power of vanishing at the initial vacuum. We suppose $\mathcal V_0$ is nonempty and, on a neighborhood $\mathcal O_0$ of it, the initial density satisfies $c_{\rm v}\operatorname{dist}(y,\mathcal V_0)^\alpha\le r_0(y)\le C_{\rm v}\operatorname{dist}(y,\mathcal V_0)^\alpha$, with $c_{\rm v},C_{\rm v},\alpha>0$. We set $\mathcal V_t\coloneqq X(t,\mathcal V_0)$ and $K_{\rm div}(t)\coloneqq\int_0^t\|\operatorname{div}V(s)\|_\infty\,ds$. For $x=X(t,y)$ with $y\in\mathcal O_0$, the Jacobian formula and the Lipschitz bounds for $X$ and $X^{-1}$ give
\be
 c_{\rm v}e^{-K_{\rm div}(t)-\alpha K(t)}\operatorname{dist}(x,\mathcal V_t)^\alpha
 \le r(t,x)\le
 C_{\rm v}e^{K_{\rm div}(t)+\alpha K(t)}\operatorname{dist}(x,\mathcal V_t)^\alpha.
\ee
Here $K$ is defined above and distances are taken in $\R^3$ or on the flat torus. Since $\Gamma$ is bounded, the proper density has the same order of vanishing. We decompose the constitutive stress as $T=T_m+T_s$, with
\be
\label{uf:eq:stress-decomposition}
 T_m=rm\Gamma(1,V)\otimes(1,V),\qquad T_s=r\tau_s(q,z,V),\qquad
 |\tau_s|\le C(q^2+|z|^2),\quad |D_{q,z}\tau_s|\le C(|q|+|z|).
\ee
{
The bounded frame change gives $\tau_s\in C^6$, with quadratic vanishing also for its $V$ derivatives. The transport equations give $\mathrm D_t(q,z)=\Ocal(q)$ and $\mathrm D_tr=-r\operatorname{div}V$, so $T_s,\mathrm D_tT_s=\Ocal(n^\gamma)$. In matter the first derivatives of $T_s$ are bounded by $CQ^{p_\theta+1}(|dQ|+|dZ|+Q|dV|)$. At $Q=0$, the bound $|T_s|\le CQ^{p_\theta+2}$ gives differential zero, so $T_s$ is $C^1$ across the cavity. The continuity equation gives $\operatorname{div}T_m=r\mathrm D_t[m\Gamma(1,V)]$ in the sense of distributions, which vanishes at zero density. Thus stress balance holds throughout.
}

For fixed $\varepsilon_{\rm e}>0$ and exact minors, $r=\det\Fsp$. The smooth dependence of the stress and material momenta on $(r,q,z,V,\Fsp)$ in~\eqref{uf:eq:material-momenta} and~\eqref{uf:eq:stress-decomposition} transfers the uniform $\Ocal(\ell^{-1/3})$ rate of Step~6 to both. A time cutoff in their balance laws gives the initial terms $T^{0\nu}(0)$ and $P_A^0(0)$ computed from the prescribed configuration and velocity. On $\mathbb R^3$ the tests have compact spatial support. This completes the proof of Theorem~\ref{uf:thm:main}.

%-----------------------------------------------------------------------------------------------------------------------------------------

\subsection{\texorpdfstring{Extensions and reference equations at cavitation}{Extensions and reference equations at cavitation}}
\label{uf:subsec:extensions-multipliers}

\paragraph{Low density and Sobolev-regular initial minors.}
A product estimate for the transported minors gives the fractional Sobolev regularity.

\begin{proof}[Proof of Proposition~\ref{uf:cor:fractional-existence}]
The strict bound $s>5/2$ gives a spatially Lipschitz velocity and makes $H^{s-1}$ an algebra; both properties enter the flow and product estimates below. At the parameter value $\delta=1$, the constitutive law, normalized equations, and time matrices in the non-relativistic construction of Section~\ref{sec:non-relativistic-vacuum} below coincide with those of the present section. The positivity of the time matrices follows here from the acoustic and range assumptions above.  On $\R^3$, we set $\sigma=s-1>3/2$. A partition into unit cubes and the algebra estimate for $H^\sigma$ give
\be
 \|mf\|_{H^\sigma}\le C\|m\|_{H^\sigma_{\rm uloc}}\|f\|_{H^\sigma}.
\ee
If $A$ and $A^{-1}$ are bi-Lipschitz with $A-\operatorname{Id},A^{-1}-\operatorname{Id}\in H^s$, composition by either map is bounded on $H^2$. Indeed, the chain rule and change of variables bound the second derivative of $f\circ A$ by $C(\|D^2f\|_2+\|D^2A\|_3\|Df\|_6)$, and $H^{s-2}\subset L^3$. Interpolation gives the same conclusion on $H^\sigma$. Writing $M_mf=mf$ and $C_Af=f\circ A$, the identity $M_{m\circ A}=C_A M_m C_{A^{-1}}$ therefore makes $m\circ A$ a multiplier on $H^\sigma$. Since $DA-I\in H^\sigma$ and $H^\sigma$ is an algebra, the flow formula $Y=R(DA)(Y_0\circ A)$ gives a uniform multiplier bound for $Y$. The identity~\eqref{uf:eq:weighted-product} then gives $P_{20}\in H^s$ and $\|P_2\|_{H^s}\le C_K\|U\|_{H^s}$, including on $Q=0$.

We use~\eqref{uf:eq:approximation} with the flow formulas for $Q_{\ell},Z_{\ell},P_{2\ell}$. The multiplier bound closes the fractional energy and the transport estimate for $h_{\ell}=A_{\ell}-\operatorname{Id}$ with strict margins on a common interval. In the material derivative of the time weight we use $\mathrm D_{\ell}Y_{\ell}=L(DV_{\ell})Y_{\ell}$. In $H^1$, the error $(I-\smoothoperator_{\ell}^2)\mathcal F_{\ell}$ is $\Ocal({\ell}^{-(s-2)})$ by the $H^{s-1}$ bound for $\mathcal F_{\ell}$, and the transport commutator is $\Ocal({\ell}^{-1})$ by the first kernel moment.

\paragraph{{Differences of reconstructed forms.}}
For two approximations with the same $Y_0$, let $A_\tau=(1-\tau)A_2+\tau A_1$. These intermediate maps are uniformly bi-Lipschitz on the short interval and have bounded displacements. Since $s-2>1/2$, local Sobolev embedding gives $DY_0\in L^3_{\rm uloc}$. By change of variables, the $L^3$ norm of $DY_0\circ A_\tau$ on each unit cube is bounded by $C\|Y_0\|_{H^{s-1}_{\rm uloc}}$. Indeed, the image of a unit cube meets only a uniformly bounded number of such cubes. We apply H\"older's inequality cube by cube and bound the $L^6$ norms of $A_1-A_2$ by their $H^1$ norms on enlarged cubes with finite overlap. Integrating in $\tau$ gives
\be
 \|Y_0\circ A_1-Y_0\circ A_2\|_2
 \le C\|Y_0\|_{H^{s-1}_{\rm uloc}}\|A_1-A_2\|_{H^1}.
\ee
The calculation for smooth $Y_0$ passes to the given background by local approximation with uniform local norms and Fatou's lemma. The matrix factors in the reconstruction have bounded entries and $L^2$ differences controlled by $\|h_{\ell}-h_M\|_{H^1}$. Hence $\|Y_{\ell}-Y_M\|_2\le C\|h_{\ell}-h_M\|_{H^1}$ without a Lipschitz assumption on $Y_0$. Subtracting the inverse-flow equations and estimating $(V_{\ell}-V_M)D^2h_M$ in $L^2$ by $L^6\times L^3$ gives an $H^1$ transport estimate for $h_{\ell}-h_M$. Together with the coefficient-difference estimate~\eqref{uf:eq:product-difference} and the grouped $H^1$ energy, this proves the Cauchy property. Interpolation with the uniform $H^s$ bounds gives convergence in $C_tH^{s'}$ for $5/2<s'<s$. The flow formulas preserve the range margins and give $Y\in C_tL^\infty$, since the uniform local Sobolev bound makes $Y_0$ uniformly H\"older continuous with exponent $s-5/2$.

The limit is weakly continuous in $H^s$. Its time weights are continuous in $L^\infty$ and have bounded material derivatives by the equations for $U,Y$. The spatial flux matrices in the paired energy and the solved principal matrices are Lipschitz. The spatial mollification argument of Step~6, applied to $\Lambda^s U$, gives continuity of the weighted energy between arbitrary times (see also Proposition~\ref{nr:prop:two-bounded-derivatives}). Freezing the positive weight gives $U\in C_tH^s$ up to $t=0$. The transport energy gives the same continuity for $h$.

For a fixed background $Y_0$, the composition operators for $A(t)$ and $A(t)^{-1}$ are strongly continuous on $H^2$, first on smooth compactly supported functions and then by density and the uniform bounds above. Interpolation and the operator identity give strong continuity of multiplication by $Y_0\circ A(t)$ on $H^\sigma$. The flow formula gives the same property for multiplication by $Y$. The identity~\eqref{uf:eq:weighted-product} then yields $P_2\in C_tH^s$. The normalized and inverse-flow equations give $U,h\in C_t^1H^{s-1}$, while the flow and transport equations give $Y\in C_tH^{s-1}_{\rm loc}\cap C_t^1H^{s-2}_{\rm loc}$. The separate transport estimate gives the additional regularity of $Y-Y_0$ when $Y_0\in W^{3,\infty}$.

On $\mathbb T^3$ we apply transport directly to $Y$. With $\sigma=s-1>3/2$, the product and commutator estimates give
\be
 \frac{d}{dt}\|Y\|_{H^\sigma}\le C\|V\|_{H^s}\|Y\|_{H^\sigma}.
\ee
The algebra property and~\eqref{uf:eq:weighted-product} yield $\|P_2\|_{H^s}\le C_K\|U\|_{H^s}$, including initially. In the approximation, $Q_\ell\ge0$, $Z_\ell=Q_\ell Y_\ell$ and $P_{2\ell}=Q_\ell Y_\ell\otimes Y_\ell$ follow from~\eqref{uf:eq:flow-formulas}. Thus the same energy estimate closes with the bounds for $Y$ and the inverse displacement. In the difference equation for the forms, $\|(\Delta V)\nabla Y\|_2\le\|\Delta V\|_6\|\nabla Y\|_3$ since $H^{s-2}\subset L^3$. Together with~\eqref{uf:eq:product-difference} and the vanishing $H^1$ approximation error, this gives the $H^1\times L^2$ Cauchy estimate of Step~6 and uniqueness. The transport energy and equation yield $Y-Y_0\in C_tH^{s-1}\cap C_t^1H^{s-2}$.

The $H^s$ velocity bound gives $\|A(t)-\mathrm{Id}\|_\infty+\|DA(t)-I\|_\infty\le Ct$. Substituting in $Y=R(DA)(Y_0\circ A)$ gives
\be
 \|Y(t)-Y_0\|_\infty\le C\bigl(t+\omega_{Y_0}(Ct)\bigr),
\ee
where $\omega_{Y_0}$ is the modulus of continuity. Sobolev embedding, applied on overlapping unit cubes on $\mathbb R^3$, gives $\omega_{Y_0}(r)\le Cr^{s-5/2}$. On the torus, approximation of the fixed $Y_0$ by smooth functions in $H^{s-1}$ improves this to $o(r^{s-5/2})$. Lipschitz data give $\omega_{Y_0}(r)\le Cr$ on either domain.

{
Step~8 gives the rank and cavity transport and the balance laws with their initial terms. For exact minors, the determinant, cofactor and the factors $q^2$ and $Q^{p_\theta+2}\Gamma^{-2\theta}$ in~\eqref{uf:eq:material-momenta}--\eqref{uf:eq:stress-decomposition} preserve $C_tH^{s-1}_{\rm loc}\cap C_t^1H^{s-2}_{\rm loc}$, since $p_\theta+2>2$ and $H^{s-1}$ is an algebra. The material fluxes thus have locally integrable first derivatives, and their divergence vanishes on $\{Q=0\}$ by locality and~\eqref{uf:eq:mass-material-equation} with $\mathrm D_tV=0$. Omitting the form variables gives the fluid case.
}
This completes the proof of Proposition~\ref{uf:cor:fractional-existence}.
\end{proof}

\paragraph{Entropy coefficients of the local solutions.}
The density factor makes the entropy multipliers $C^1$ through cavitation.

\begin{proposition}[Regularity of the entropy multipliers]
\label{uf:prop:entropy-multipliers}
Fix a solution of Theorem~\ref{uf:thm:main} or a solution of Proposition~\ref{uf:cor:fractional-existence} on $\mathbb T^3$ or $\mathbb R^3$, with $5/2<s\le3$. Use the constitutive extension obtained by projecting the one- and two-forms onto the rest space of the current. The enthalpy at fixed rest minors is
\be
\label{uf:eq:entropy-enthalpy}
 H=h+\Psi+\theta q\Psi_q+\theta z\cdot\Psi_z\ge m.
\ee
\begin{itemize}
\item[(i)] \emph{Regularity on the closed slab.} The multipliers $w,C,D$ defined in the equation~\eqref{eq:main-vacuum-finite-multipliers} below, and the associated spatial involution coefficients, have bounded continuous first derivatives on the closed local slab, including its cavities. Their bounds depend only on the regularized norms in these results, the constitutive bounds and a bounded range of $\varepsilon_{\rm e}$. They require neither a positive density minimum nor a fourth derivative of the regularized density or velocity.
\item[(ii)] \emph{Bounds at zero density.}  Set $s=3$ for Theorem~\ref{uf:thm:main}, and choose $0<\alpha\le s-5/2$. With $\sigma=\alpha/(1+\alpha)$ and $\gamma=1+2\theta$, the form multipliers satisfy
\be
 |C|+|D|\le K n^\gamma,\qquad
 |\nabla C|+|\nabla D|\le K n^{1+\theta(1+\sigma)},
\ee
{Their material derivatives satisfy the sharper bound
\be
 |u^\mu\nabla_\mu C|+|u^\mu\nabla_\mu D|\le K n^\gamma.
\ee}
The constant $K$ also depends on $s$ and $\alpha$.
\item[(iii)] \emph{Use in comparison.} At $\varepsilon_{\rm e}=0$, they reduce to the fluid multipliers $w=u^\flat/h$, $C=D=0$. On $\mathbb T^3$,  if these covectors satisfy the supporting inequalities on the chosen comparison domain, the local solutions satisfy Definition~\ref{def:entropy-strong-solution}. Theorem~\ref{thm:physical-entropy-vacuum} also requires the entropy-weighted bound on the total remainder.
\end{itemize}
\end{proposition}
\begin{proof}
At fixed rest minors, $H=\bar\rho_n\to m$ as $n\downarrow0$. For positive elastic parameter,~\eqref{uf:eq:rest-metric} gives $\bar\rho_{nn}\ge0$ on the chosen range; in the fluid case $h_q=\theta qH_f\ge0$. Hence $H\ge m$. In the rest frame the forms have zero time components. A variation of the rest direction leaves their spatial components unchanged to first order, up to a rotation under which $\Psi$ is invariant. Thus $B=\del_{\particlecurrent}\raug=-Hu^\flat$ and $B\cdot\particlecurrent=nH$. For $P=\del_F\raug$, $Q=\del_G\raug$, the spatial derivatives in an oriented orthonormal rest frame are
\be
\label{uf:eq:rest-entropy-derivatives}
 P_A^i=\varepsilon_{\rm e}q^{p_\theta+1}\Psi_{z_F^{A,i}},\qquad
 2Q^{ijA}=\varepsilon_{\rm e}q^{p_\theta+1}\varepsilon^{kij}\Psi_{z_{G,A}^k},\qquad P_A^0=Q^{0iA}=0.
\ee
Formula~\eqref{eq:lambda-explicit-local06} gives $w=u^\flat/H$ and the form multipliers~\eqref{eq:main-vacuum-finite-multipliers}. Since $q,z,V\in C_tH^s\cap C_t^1H^{s-1}$ and $s>5/2$, these fields are $C^1$ with bounded derivatives. The boost is smooth on the velocity range, $H\ge m$, and $q^{p_\theta+1}$ is $C^1$ at zero. This proves (i) also for uniformly local initial minors.

For (ii), $\Psi_z(q,0)=0$ and the bounded rest minors give $|z|+|\Psi_z|\le Cq$. Sobolev embedding gives a spatial $C^{0,\alpha}$ bound for $D_x(q,z)$. For $q(x)>0$, Taylor expansion in the direction opposite to $\nabla q(x)$ at distance $q(x)^{1/(1+\alpha)}$, using $q\ge0$, gives $|\nabla q|\le Cq^\sigma$. At the same scale symmetric differences and $|z|\le Cq$ give $|\nabla z|\le Cq^\sigma$. The bounds extend by continuity at zero. Since $\mathrm D_t(q,z)=\Ocal(q)$ they hold for time derivatives as well. Differentiating~\eqref{uf:eq:rest-entropy-derivatives} gives (ii), with $q=n^\theta$. Along the flow, $\Psi_{zq}(q,z)=\Ocal(|z|)$, so $\mathrm D_t(q^{p_\theta+1}\Psi_z)=\Ocal(q^{p_\theta+2})=\Ocal(n^\gamma)$. All remaining factors have bounded material derivatives and $u^\mu\nabla_\mu=\Gamma\mathrm D_t$.

At zero elastic parameter $C=D=0$. On the periodic slab the supports assumed in (iii), the physical balance laws, current conservation and the prescribed data give Definition~\ref{def:entropy-strong-solution}. Theorem~\ref{uc:thm:common} supplies the separate remainder estimate under its comparison hypotheses.
\end{proof}

%-----------------------------------------------------------------------------------------------------------------------------------------

\paragraph{Entropy multipliers at zero current.}
\label{subsec:common-entropy-multipliers}
For comparisons with finite reference velocity, we use exact minors on the completed graph~\eqref{eq:finite-flux-state-space}. The extended unit velocity selects the reference stress, which may persist at zero current. For the positive-rest-mass law~\eqref{eq:massive-vacuum-class}, assume relativistic polyconvexity and the entropy conditions of Section~\ref{subsec:entropy-conditions-vacuum} on $\mathcal L(F,G,n^2;g)$ and $\mathcal U(n)$. At fixed rest forms the enthalpy is
\be
\label{eq:general-physical-enthalpy}
H=m+\mathcal U'(n)+2n\mathcal L_s(F,G,n^2;g).
\ee
The subscript $s$ denotes differentiation with respect to $s=-g(\particlecurrent,\particlecurrent)$. Let $F_*,G_*,u_*,H_*$ be bounded Lipschitz reference fields, with $H_*\ge m>0$, and let ${\particlecurrent}_*=n_*u_*$, $n_*\ge0$. We require
\be
\label{eq:main-vacuum-reference-kinematics}
i_{u_*}F_*=i_{u_*}G_*=0,\qquad
 dF_*=dG_*=0,\qquad \nabla_\mu(n_*u_*^\mu)=0.
\ee
For the smooth term $\mathcal L$, we use full antisymmetric sums and define
\be
\label{eq:main-vacuum-reference-forces}
P_A^\mu\coloneqq\del_{F_\mu^A}\mathcal L,\quad
Q^{\mu\nu b}\coloneqq\del_{G_{\mu\nu b}}\mathcal L,\quad
L_{*A}\coloneqq\nabla_\mu P_{*A}^\mu,\quad
K_*^{\lambda b}\coloneqq2\nabla_\mu Q_*^{\mu\lambda b},\quad
\Pi_{*\mu}{}^\nu\coloneqq\delta_\mu{}^\nu+u_{*\mu}u_*^\nu.
\ee
The covector combining acceleration and the enthalpy gradient is
\be
\label{eq:main-vacuum-reference-acceleration}
\mathcal A_{*\mu}\coloneqq H_*u_*^\nu\nabla_\nu u_{*\mu}
 +\Pi_{*\mu}{}^\nu\del_\nu H_*.
\ee
Since $H_*\ge m>0$, the entropy multipliers remain finite at zero density:
\be
\label{eq:main-vacuum-finite-multipliers}
w_{*\mu}\coloneqq u_{*\mu}/H_*,\qquad C_{*A}^{\nu\alpha}\coloneqq w_*^\nu P_{*A}^\alpha-w_*^\alpha P_{*A}^\nu,\qquad D_{*b}^{\nu\alpha\beta}\coloneqq w_*^\nu Q_{*b}^{\alpha\beta}-w_*^\alpha Q_{*b}^{\nu\beta}+w_*^\beta Q_{*b}^{\nu\alpha}.
\ee
With $\Delta$ denoting the difference from the reference, the relative current~\eqref{eq:covariant-relative-current-matter} now extends to zero density:
\be
\label{eq:main-vacuum-relative-current}
\mathscr Q^\nu\coloneqq-\Delta {\particlecurrent}^\nu-w_{*\mu}\Delta T^{\mu\nu}
-C_{*A}^{\nu\alpha}\Delta F_\alpha^A
-D_{*b}^{\nu\alpha\beta}\Delta G_{\alpha\beta b}.
\ee
The supporting inequality for the time component extends to zero reference current by Theorem~\ref{thm:closed-vacuum-entropy} and Proposition~\ref{prop:vacuum-stress-subgradient}. With only the stress rows densitized as in Section~\ref{subsec:entropy-conditions-vacuum}, the stress and form coefficients are $w_{*\mu}$, $|g|^{1/2}C_{*A}^{0i}$, and $2|g|^{1/2}D_{*b}^{0ij}$ for $i<j$, respectively. For the spatial equations $\del_iF_j^A-\del_jF_i^A=0$ and $\del_1G_{23b}+\del_2G_{31b}+\del_3G_{12b}=0$, we use
\be
\label{eq:physical-spatial-constraint-coefficients}
 b_{*,Aij}^F\coloneqq-|g|^{1/2}C_{*A}^{ij}\quad(i<j),\qquad
 b_{*,b}^G\coloneqq-2|g|^{1/2}D_{*b}^{123}.
\ee
By antisymmetry these coefficients give the spatial current~\eqref{eq:general-entropy-relative-current} for $|g|^{1/2}\mathscr Q$.

Testing the balance laws, including their initial terms, against these Lipschitz coefficients gives $\nabla_\nu\mathscr Q^\nu\le\mathscr R$, with $\mathscr R$ affine in the four flux differences. Multiplication by the volume density gives the inequality in Lemma~\ref{thm:common-current-jhde}.
%-----------------------------------------------------------------------------------------------------------------------------------------

\paragraph{Material equations with cavitation.}
\label{subsec:physical-entropy-regularity}

At positive density, the material momentum of the full elastic law is
\be
\label{eq:general-physical-material-momentum}
\mathcal P_A^\nu(F)=\frac{\del\rphys}{\del F_\nu^A}(F;g)=\frac{d}{d F_\nu^A}\raug(F,G(F),\particlecurrent(F);g).
\ee
The reference continuation specifies all three momenta at zero density. For the law~\eqref{eq:massive-vacuum-class}, use the derivatives in~\eqref{eq:main-vacuum-reference-forces} and set
\be
\label{eq:physical-current-piola-tensor}
\mathcal H_A^{\mu\nu}(F)\coloneqq
\frac{\del {\particlecurrent}^\mu(F)}{\del F_\nu^A}
=-\frac12\eps^{\mu\nu\rho\sigma}G_{\rho\sigma A}(F),
\qquad
\mathcal C_{*A}^{\nu}\coloneqq-u_{*\mu}\mathcal H_A^{\mu\nu}(F_*).
\ee
These polynomial expressions remain finite at rank loss and give the continued material momentum
\be
\label{eq:physical-continued-material-momentum}
\mathcal P_{*A}^{\nu}\coloneqq
P_{*A}^{\nu}
+2\Omega_{bAC}Q_*^{\nu\lambda b}F_{*\lambda}^C
-H_*u_{*\mu}\mathcal H_A^{\mu\nu}(F_*).
\ee
We impose the three material equations
\be
\label{eq:physical-continued-material-EL}
\nabla_\nu\mathcal P_{*A}^{\nu}=0
\qquad(A=1,2,3)
\ee
throughout the slab. These Euler--Lagrange equations in matter extend to zero current through the chosen velocity, without differentiating $n(F)$ at a degenerate gradient.

Closure and the Piola identities give the equivalent equation in the finite variables of~\eqref{eq:main-vacuum-reference-forces} and~\eqref{eq:main-vacuum-reference-acceleration}:
\be
\label{eq:physical-material-force-equation}
L_{*A}+\Omega_{bAC}F_{*\lambda}^C K_*^{\lambda b}
+\mathcal C_{*A}^{\mu}\mathcal A_{*\mu}=0.
\ee
For the constitutive continuations verified in Section~\ref{subsec:constitutive-stability-verification} below, the material equations give the adjoint identity~\eqref{eq:general-entropy-adjoint-identity}. The linear functional on increments $(\zeta_G,\zeta_{\particlecurrent})$ of the constraints, in the density convention of~\eqref{eq:physical-spatial-constraint-coefficients}, is
\be
\label{eq:physical-adjoint-normal-coefficient}
\nu_*(\zeta_G,\zeta_{\particlecurrent})=\frac{|g|^{1/2}}{H_*}\left(u_*^\mu\zeta_{G,\mu\lambda b}K_*^{\lambda b}-\zeta_{\particlecurrent}^\mu\mathcal A_{*\mu}\right).
\ee
The affine constitutive maps are fixed before the reference evolution. Lipschitz reference fields and conjugates give a bound on $\nu_*$ uniform down to $n_*=0$. Contracting the material equations with $F_{*\mu}^A$ gives the stress equations also at rank loss.

For conservative variables $V$ without volume factors, the time component of~\eqref{eq:main-vacuum-relative-current} is the relative density $\mathcal D_\eta(V\mid V_*)$ of~\eqref{eq:general-entropy-first-remainder}, with $\eta=-\particlecurrent^0$ on the completed graph. Section~\ref{subsec:constitutive-stability-verification} uses the continued material equations and exact minors to cancel the linear terms and isolate the density contribution. Once the remainder is controlled, Theorem~\ref{thm:general-constitutive-entropy-comparison} gives the integral estimate.

\subsection{Proofs of stability for the local constitutive laws}
\label{subsec:local-stability-proofs}

\paragraph{Convex entropy and pointwise comparisons.}
We construct supports at zero density and estimate the relative entropy.
\begin{proof}[Proof of Proposition~\ref{uf:prop:local-convex-entropy}]
{\bf 1.} \emph{The rest-energy remainder.} We put $Y_{\rm R}=\varepsilon_{\rm e}(\Frest,\Grest)$ and write
\bse
\be
 \bar\rho^\sharp(n,Y_{\rm R})\coloneqq {m \, n}+n\bigl(\Phi(n^\theta)+\Psi(n^\theta,n^\theta Y_{\rm R})\bigr),\qquad
 \gamma\coloneqq1+2\theta,
\ee
\be
 \mathcal B_{\bar\rho^\sharp}((n,Y_{\rm R})\mid(n_*,Y_{{\rm R}*}))\coloneqq
 \bar\rho^\sharp(n,Y_{\rm R})-\bar\rho^\sharp(n_*,Y_{{\rm R}*})-D\bar\rho^\sharp(n_*,Y_{{\rm R}*})\cdot(n-n_*,Y_{\rm R}-Y_{{\rm R}*}).
\ee
\ese
The derivative at zero is $D\bar\rho^\sharp(0,Y_{{\rm R}*})=(m,0)$. For $n>0$, the formula for the normalized matrix representing the second variation of the energy in the material rest frame in~\eqref{uf:eq:rest-metric} gives
\bse
\be
\label{uf:eq:weighted-rest-hessian}
 D^2\bar\rho^\sharp[(\dot n,\dot Y_{\rm R})]^2
 =nq^2\mathsf H\left[\left(\theta\frac{\dot n}{n},\dot Y_{\rm R}+\theta Y_{\rm R}\frac{\dot n}{n}\right)\right]^2
 \asymp n^{\gamma-2}\dot n^2+n^\gamma|\dot Y_{\rm R}|^2.
\ee
\be
\label{uf:eq:weighted-rest-remainder}
 \mathcal B_{\bar\rho^\sharp}\asymp B_\gamma(n\mid n_*)+(n+n_*)^\gamma|Y_{\rm R}-Y_{{\rm R}*}|^2
 \asymp(n+n_*)\bigl(|q-q^\ast|^2+|z-z^\ast|^2\bigr),
\ee
\ese
where $q^\ast\coloneqq n_*^\theta$, $z^\ast\coloneqq q^\ast Y_{{\rm R}*}$ and $B_\gamma(n\mid n_*)\coloneqq n^\gamma-n_*^\gamma-\gamma n_*^{\gamma-1}(n-n_*)$. Indeed, writing $X=\theta\dot n/n$ and $W=\dot Y_{\rm R}+Y_{\rm R}X$, the inverse change is $\dot n/n=X/\theta$, $\dot Y_{\rm R}=W-Y_{\rm R}X$. Both maps are uniformly bounded on the prescribed shape ball for each fixed $\theta>0$. The positive lower and upper bounds on $\mathsf H$ therefore give the comparison~\eqref{uf:eq:weighted-rest-hessian}; its constants may depend on $\theta$. Integrating on the affine $(n,Y_{\rm R})$ segment proves the first comparison. The shape coefficient obeys
\be
 \frac{2^{-\gamma-1}}{\gamma+1}(n+n_*)^\gamma
 \le\int_0^1(1-s)((1-s)n_*+sn)^\gamma\,ds
 \le\frac12(n+n_*)^\gamma.
\ee
For the lower bound, we use $(1-s)n_*+sn\ge(n+n_*)\min(s,1-s)$ and integrate. 

{For the two scalar comparisons put $S=n+n_*>0$ and $t=n/S$. By homogeneity the ratios
\be
 \frac{|n^\theta-n_*^\theta|}{S^{\theta-1}|n-n_*|},
 \qquad \frac{B_\gamma(n\mid n_*)}{S^{\gamma-2}|n-n_*|^2}
\ee
depend only on $t\in[0,1]$. Strict monotonicity of the power $n^\theta$ and strict convexity of $n^\gamma$ make them positive away from $t=1/2$, including at both endpoints. Taylor's formula at $1/2$ gives positive limits $\theta 2^{1-\theta}$ and $\gamma(\gamma-1)2^{1-\gamma}$, respectively. Compactness therefore gives the two-sided bounds for each fixed $\gamma>1$, independently of any positive lower bound for either density.}

For the shape comparison, we set $q=n^\theta$, $z=qY_{\rm R}$ and let $M$ be an upper bound on $|Y_{\rm R}|,|Y_{{\rm R}*}|$. The two identities
\be
 z-z^\ast=q(Y_{\rm R}-Y_{{\rm R}*})+Y_{{\rm R}*}(q-q^\ast)
 =q^\ast(Y_{\rm R}-Y_{{\rm R}*})+Y_{\rm R}(q-q^\ast)
\ee
show that $(q+q^\ast)|Y_{\rm R}-Y_{{\rm R}*}|\le2(|z-z^\ast|+M|q-q^\ast|)$, while the reverse bound follows from either identity. Since $q+q^\ast\asymp(n+n_*)^\theta$, this proves the shape part of the comparison without division by a density. The Hessian integral extends to either endpoint of zero density since $\gamma>1$; when both densities vanish, every term is zero.

{
The factors in~\eqref{uf:eq:taylor-factors} give $S=n^\gamma s(n^\theta,Y_{\rm R})$ for the spatial rest stress, with two bounded derivatives of $s$. In~\eqref{nr:eq:specific-stress} this follows from $\Phi'(q),\Psi_z(q,qY_{\rm R})=\Ocal(q)$ and $\Psi_q(q,qY_{\rm R})=\Ocal(q^2)$, with the corresponding smooth factors. Thus for $n>0$ the absolute second differential is bounded by $C(n^{\gamma-2}\dot n^2+n^\gamma|\dot Y_{\rm R}|^2)$. Integrate over the same affine segment as for~\eqref{uf:eq:weighted-rest-remainder}. The density weight is integrable at a zero endpoint since $\gamma>1$, and the first derivatives vanish at zero density. This gives
\be
\label{uf:eq:weighted-rest-stress-remainder}
 |S-S_*-DS_*\cdot(n-n_*,Y_{\rm R}-Y_{{\rm R}*})|\le C\mathcal B_{\bar\rho^\sharp}.
\ee
}
\paragraph{{Entropy in the observer frame.}}
\vskip.2cm \noindent {\bf 2.} \emph{The observer variables and supporting covectors.} We write $\mathsf A(V)$ for the linear transformation to rest-frame variables defined in~\eqref{uf:eq:rest-shapes}, and put $\Gamma=(1-|V|^2)^{-1/2}$. The observer Lagrangian is $L(r,Y,V)\coloneqq \bar\rho^\sharp(r/\Gamma,\mathsf A(V)Y)$. At fixed $V$, this transformation is linear and both form blocks contract the Euclidean norm. Thus $L$ is convex in $(r,Y)$. For the fluid part, $h(0)=m$ and $\del_q h=\theta qH_f\ge0$. Differentiation gives
\bse
\be
 -D_V^2L_{\rm fl}=rh\Gamma\bigl(I+\Gamma^2(1-c_f^2)V\otimes V\bigr).
\ee
{Its transverse eigenvalues are $rh\Gamma$, and its eigenvalue in the velocity direction is $rh\Gamma^3(1-c_f^2|V|^2)$, with $q=(r/\Gamma)^\theta$. Thus the observer condition in the proposition gives $-D_V^2L_{\rm fl}\ge m\delta_{\rm c}rI$. The light-speed bound~\eqref{uf:eq:fluid-sound-bound} is sufficient for this condition for every $v_{\rm c}<1$.}  The elastic perturbation vanishes faster than $r$ at zero density. Indeed, with $q=(r/\Gamma)^\theta$ and $Y_R=\mathsf A(V)Y$, Taylor's formula and $\Psi(q,0)=\Psi_z(q,0)=0$ give
\be
 L-L_{\rm fl}=\frac{rq^2}{\Gamma}\int_0^1(1-\tau)\,
 Y_R^T\Psi_{zz}(q,\tau qY_R)Y_R\,d\tau.
\ee
At fixed $r,Y$, the first two velocity derivatives of $q$ are $\Ocal(q)$, so the bounded constitutive derivatives and velocity range yield $|D_V^2(L-L_{\rm fl})|\le C_\Psi r^\gamma s_{\rm c}^2$. We fix $\mu\in(0,m\delta_{\rm c})$ and choose $C_\Psi r_{\max}^{\gamma-1}s_{\rm c}^2\le m\delta_{\rm c}-\mu$. The constant $C_\Psi$ depends on the fixed compact coefficient ranges and constitutive derivatives. We then have $\mu rI\le-D_V^2L\le CrI$ on the entire branch.

We write $E=T^{00}=L-V\cdot L_V$ and $P=T^{0\rm sp}=-L_V$. Differentiation gives $dE-V\cdot dP=L_r\,dr+L_Y\cdot dY$. Thus the entropy covector in the variables $(Y,E,P)$ is $(L_Y,-1,V)/L_r$, with the constant scaling of the form variables understood. Evaluating this covector at the reference and substituting the expressions for $E,P$ gives
\be
\label{uf:eq:observer-relative-decomposition}
\begin{aligned}
 L_{r*}D_{\eta_{\mathcal K}}
 &=L(r,Y,V_*)-L(r_*,Y_*,V_*)-L_{(r,Y)*}\cdot(r-r_*,Y-Y_*)\\
 &\quad+L(r,Y,V)-L(r,Y,V_*)-L_V(r,Y,V)\cdot(V-V_*).
\end{aligned}
\ee
The first remainder uses $V_*$ at both endpoints, so $n=r/\Gamma_*$ and $(q,z)=\Gamma_*^{-\theta}(Q,\mathsf A(V_*)Z)$ there. This change and its inverse are uniformly bounded, giving the density and form terms by~\eqref{uf:eq:weighted-rest-remainder}. The bounds on $-D_V^2L$ give the velocity term $r|V-V_*|^2$. Since $L_{r*}=H_*/\Gamma_*$ is bounded above and away from zero, these estimates give~\eqref{uf:eq:local-relative-entropy-comparison} through zero density.

Define $\eta_{\mathcal K}$ as the supremum of the affine supports at positive-density points of $\mathcal K$. By~\eqref{uf:eq:observer-relative-decomposition} every support lies below $-r$ on the conservative image and attains it at its base point. The conservative variables and supporting covectors are bounded at fixed $\varepsilon_{\rm e}>0$, so the supremum is finite and convex; passage to $r=0$ gives the entropy identity there.
\ese

\paragraph{{Supports at zero density.}}
At $U_{\rm c}^0=(\Fsp_0,\Gsp_0,0,0)$, taking $r\downarrow0$ with different $|V_*|\le v_{\rm c}$ gives the distinct supports $(0,-\Gamma_*/m,\Gamma_*V_*/m)$. The entropy is therefore nondifferentiable there. Its relative density $-r+(\Gamma_*/m)(E-V_*\cdot P)$ can depend on the continued reference velocity. On each bounded finite-velocity branch, the internal stress derivatives in the current and forms are respectively $\Ocal(n^{\gamma-1})$ and $\Ocal(n^\gamma)$ and extend by zero. The mass derivative is~\eqref{eq:nonlinear-fluid-affine-jet} with $h=m$. Together with the linear form fluxes these are the bounded affine continuations in Definition~\ref{def:entropy-strong-solution}; omitting the form blocks gives the fluid case.
This completes the proof of Proposition~\ref{uf:prop:local-convex-entropy}.
\end{proof}

At $n=0$, take $\Frest=\operatorname{diag}(0,\alpha,\beta)$ in the rest frame, with $\alpha,\beta>0$ in the small elastic range. Then $\Grest=\operatorname{cof}\Frest=\operatorname{diag}(\alpha\beta,0,0)$, $\particlecurrent=T=0$, $H=m$ and $C=D=0$. Comparing two such configurations gives
\be
\label{uf:eq:vacuum-constraint-remainder}
 D_{\eta_{\mathcal K}}(U_{\rm c}\mid U_{{\rm c}*})=0,\qquad
 r_{\Theta,G_{23\,1}}=-(\alpha-\alpha_*)(\beta-\beta_*).
\ee
Thus~\eqref{eq:general-entropy-Taylor-control} fails for the unweighted minors separately.

\paragraph{{Pointwise bounds on the constitutive remainder.}}
{We use the material equations to estimate the total remainder, distinguishing comparable from separated densities.}
\begin{proof}[Proof of Theorem~\ref{uc:thm:common}]
\noindent {\bf 1.} \emph{Reference coefficient bounds.} We put $\gamma=1+2\theta$.
{
We use the stress factorization just established. Writing $\mathcal L=\raug-m\,n$, both $\mathcal L$ and $T_{\mathcal L}$ are $n^\gamma$ times bounded $C^6$ functions of $(q,\varepsilon_{\rm e}F,\varepsilon_{\rm e}G,u)$ with the forms projected onto the rest space. Proposition~\ref{uf:prop:entropy-multipliers}, the normalized evolution and form transport give the following bounds for the coefficients defined in~\eqref{eq:main-vacuum-reference-forces}--\eqref{eq:main-vacuum-finite-multipliers}.
}
\bse
\be
\label{uc:eq:general-reference-bounds}
 |\mathcal A_*|\le Cq^\ast,\qquad
 |L_*|+|K_*|\le C n_*q^\ast,\qquad
 |u_*\cdot dH_*|\le C(q^\ast)^2,
\ee
\be
 |\nabla w_*|+|\nabla(m/H_*)|+|\nabla u_*|\le C.
\ee
\ese
Here $q^\ast=n_*^\theta$ is the reference density coordinate, as in the proof of Proposition~\ref{uf:prop:local-convex-entropy}. Indeed the relative spatial matrices vanish with $Q$ on $|Z|\le C Q$, so the normalized velocity equation gives $u_*\cdot\nabla u_*=\mathcal {\Ocal}(q^\ast)$. Also $H_*-m=\mathcal {\Ocal}((q^\ast)^2)$ and its first derivatives are $\mathcal {\Ocal}(q^\ast)$. Cartan transport gives $u_*\cdot dq^\ast,u_*\cdot dz_*=\mathcal {\Ocal}(q^\ast)$, proving the third bound. Differentiating $P=\varepsilon_{\rm e}q^{p_\theta+1}\Psi_{z_F}$ and its two-form counterpart, with $\Psi_z=\mathcal {\Ocal}(q)$, proves the force bound. These computations differentiate only $(q,z,V)$ once.  The unweighted reference minors belong to $C_tH^{s-1}\cap C_t^1H^{s-2}$ and need not be Lipschitz. Their derivatives are locally square integrable in spacetime. The product rule with the $C^1$ coefficients $P_*,Q_*,u_*,H_*$ and the distributional closure and Piola identities therefore yield the material force identity~\eqref{eq:physical-material-force-equation} almost everywhere, using the continued material equations.

\vskip.2cm \noindent {\bf 2.} \emph{The entropy bound and the remainder at comparable densities.} At a point of comparison, we use a constant Lorentz coframe at rest with respect to $u_*$. We write $n$ for the comparison proper density, $v$ for its relative velocity, $\Fsp $ for its spatial one-form matrix and $\Gsp =\cof \Fsp $. We write $\Delta$ for the difference from the reference value. We put $\Gamma_{\rm c}\coloneqq(1-v_{\rm c}^2)^{-1/2}$ and choose $\vartheta\in(0,1)$ with $\Gamma_{\rm c}\vartheta<1$. We shall use the following consequence of Proposition~\ref{uf:prop:local-convex-entropy}.  At $n=n_*=0$, the density term in this estimate is defined by continuity to be zero; indeed it is bounded above by $(n+n_*)^\gamma$ whenever $n+n_*>0$:
\bse
\be
\label{uc:eq:weighted-coercivity}
 \mathscr D\coloneqq D_{\eta_{\mathcal K}}(U_{\rm c}\mid U_{{\rm c}*})
 \ge c\bigl((n+n_*)^{\gamma-2}|n-n_*|^2
 +(n+n_*)^\gamma(|\Delta \Fsp |^2+|\Delta \Gsp |^2)+n|v|^2\bigr).
\ee
{
The relative density on the left uses the original observer. Write $r=n\Gamma(V)$ and $r_*=n_*\Gamma(V_*)$, with both Lorentz factors in $[1,\Gamma_{\rm c}]$. On this compact velocity range, $|v|\asymp|V-V_*|$. For comparable positive densities, the change of frame gives
\be
 |n-n_*|\le C(|r-r_*|+n|V-V_*|),\qquad
 |\Delta(\Fsp,\Gsp)|\le C\bigl(\varepsilon_{\rm e}^{-1}|Y-Y_*|+|V-V_*|\bigr).
\ee
The shape comparison in~\eqref{uf:eq:weighted-rest-remainder} applies to $Z=QY$ with $Q=r^\theta$. The additional velocity terms have weight at most $Cn^\gamma\le Cn$. For separated proper densities the ratio of the smaller observer density to the larger one is at most $\lambda=\Gamma_{\rm c}\vartheta<1$. Thus, with $R=\max(r,r_*)$, one has
}
\be
\label{uc:eq:separated-density-bound}
 (r+r_*)|r^\theta-r_*^\theta|^2
 \ge(1-\lambda^\theta)^2R^\gamma
 \ge2^{-\gamma}(1-\lambda^\theta)^2(r+r_*)^\gamma.
\ee
The scalar density remainder therefore controls $(n+n_*)^\gamma$, including either density equal to zero. Bounded shapes and the kinetic term give the bound~\eqref{uc:eq:weighted-coercivity} there also. If both densities are zero, both sides vanish by Proposition~\ref{uf:prop:local-convex-entropy}. This proves the bound~\eqref{uc:eq:weighted-coercivity} using the relative density of the original observer.

\paragraph{{Taylor expansion in spacetime variables.}}
For endpoints in $\mathcal K$, we write ${\particlecurrent}_\lambda=(1-\lambda){\particlecurrent}_*+\lambda {\particlecurrent}$. Wherever $r_\lambda={\particlecurrent}_\lambda^0>0$, its velocity is the average of the endpoint velocities with weights $(1-\lambda)r_*$ and $\lambda r$. The velocity bound gives the following comparison of densities
\be
 \Gamma(v_{\rm c})^{-1}r_\lambda\le n_\lambda\le r_\lambda\le r_{\max},
 \qquad \Gamma(v_{\rm c})\coloneqq(1-v_{\rm c}^2)^{-1/2}.
\ee
In particular, for comparable positive endpoint densities, $n_\lambda$ is comparable to $n_*$. The argument also gives a sufficient size for the elastic margin. In the Euclidean norm on independent form components, orthogonality at the endpoints bounds the full scaled spacetime forms by $\sqrt{1+v_{\rm c}^2}\,s_{\rm c}$. Indeed, each time component of a one-form is its spatial part contracted with the velocity, and the vector of time-space components of a two-form is obtained by contracting its spatial two-form with that velocity. The same norm bound holds along the affine form segment. The boost to the rest frame of ${\particlecurrent}_\lambda$ has operator norm at most $\Gamma(v_{\rm c})(1+v_{\rm c})$ on one-forms and two-forms, since its singular values on one-forms are $e^\xi,1,1,e^{-\xi}$ with $|\tanh\xi|\le v_{\rm c}$. In that frame, taking spatial components is an orthogonal projection in the Euclidean component norm. Hence it suffices that $\sqrt{1+v_{\rm c}^2}\,\Gamma(v_{\rm c})(1+v_{\rm c})s_{\rm c}$ lie strictly inside the radius chosen for the rest-minor bounds. For $n_*>0$ and $\vartheta n_*\le n\le n_*/\vartheta$, Taylor expansion along the segment bounds the second derivatives of the internal stress by
\be
 |D_{\particlecurrent}^2T_{\mathcal L}|\le Cn_*^{\gamma-2},\qquad
 |D_{\particlecurrent}D_{F,G}T_{\mathcal L}|\le Cn_*^{\gamma-1},\qquad
 |D_{F,G}^2T_{\mathcal L}|\le Cn_*^\gamma.
\ee
\ese
The ordinary internal stress remainder is therefore bounded by $C\mathscr D$, using the lower bound~\eqref{uc:eq:weighted-coercivity} and $|\Delta {\particlecurrent}|\le C(|n-n_*|+n|v|)$.

\paragraph{{Cancellation in the nonlinear source.}}
The joint Noether identity~\eqref{eq:physical-nonlinear-current-joint-Noether} applies at a positive-density reference to this smooth internal law. Its cancellation uses symmetry of the joint Hessian and does not require dependence on $n^2$. We put $\vartheta_*=m/H_*$ and $\Delta u=u-u_*$. Splitting off the mass stress gives the full identity
\bse
\be
\label{uc:eq:mass-internal-identity}
 \mathscr R_*=-\vartheta_*n\Delta u^\mu\Delta u^\nu\nabla_\nu u_{*\mu}
 +(\nabla\vartheta_*)\cdot\mathscr Q_m
 -(\nabla_\nu w_{*\mu})T_{\mathcal L,\mathrm{rel}}^{\mu\nu}
 +\mathscr R_{\rm force},\qquad
 \mathscr Q_m\coloneqq n(\chi-1)u,\quad\chi\coloneqq-g(u,u_*),
\ee
\be
 H_*\mathscr R_{\rm force}
 =b_A(\Delta\Gsp)_{Ai}\mathcal A_{*i}
  +b_A\Omega_{bAC}(\Delta\Fsp)_{Ci}K_*^{ib},\qquad b_A\coloneqq F_0^A=-(\Fsp)_{Aj}v_j.
\ee
\ese
{For the internal energy $\mathcal L=\raug-{m \, n}$, its current conjugate at the reference is $B_{\mathcal L*}=-(H_*-m)u_{*\flat}$. Its multiplier in~\eqref{eq:physical-nonlinear-current-joint-Noether} is therefore $\alpha_*=w_*\cdot B_{\mathcal L*}=1-\vartheta_*$. Thus $d\alpha_*=-d\vartheta_*$, and the terms differentiating these weights combine into $(\nabla\vartheta_*)\cdot\mathscr Q_m$. The acceleration terms from the mass and internal energies combine into $-H_*^{-1}\Delta {\particlecurrent}^\mu\mathcal A_{*\mu}$. Together with the form-force terms in that identity, the exact relations ${\particlecurrent}^i=-\Gsp _{Ai}b_A$ and $G_{0ib}=\Omega_{bAC}b_A(\Fsp)_{Ci}$ give the displayed formula for $\mathscr R_{\rm force}$ after subtracting the reference material equation.} The two mass terms in the identity~\eqref{uc:eq:mass-internal-identity} are bounded by $Cn|v|^2$. To bound the force terms, the exact minor relations and the reference material equations give the identities below.  In the reference rest frame, the cofactor relations are $\Fsp^T\Gsp=n\chi I_3$ and $\Fsp_*^T\Gsp_*=n_*I_3$; these identities give the density difference in the second equation.
\begin{samepage}
\bse
\be
 E_A\coloneqq L_{*A}+\Omega_{bAC}(\Fsp_*)_{Ci}K_*^{ib}
 =-(\Gsp_*)_{Ai}\mathcal A_{*i},
\ee
\be
\label{uc:eq:force-reorganization}
 H_*\mathscr R_{\rm force}
 =(n_*-n\chi)v\cdot\mathcal A_*
 -v_j(\Delta\Fsp)_{Aj}E_A
 +b_A\Omega_{bAC}(\Delta\Fsp)_{Ci}K_*^{ib}.
\ee
\ese
\end{samepage}
The cofactor identity gives $b_A(\Gsp)_{Ai}=-n\chi v_i$. Writing the other contraction as $-b_A(\Gsp_*)_{Ai}=v_j(\Fsp)_{Aj}(\Gsp_*)_{Ai}$ and expanding $\Fsp=\Fsp_*+\Delta\Fsp$ gives the term $n_*v_i$ and the remaining contraction with $\Delta\Fsp$. These are polynomial identities in the form components and require no division by the density.
The first force term before this reorganization must not be estimated separately. The bounds~\eqref{uc:eq:general-reference-bounds} give $|E|+|K_*|\le Cn_*q^\ast$ and $|\mathcal A_*|\le Cq^\ast$. Since $\gamma=1+2\theta$, weighted Young inequalities bound $|n-n_*|q^\ast|v|$ by $C(n_*^{\gamma-2}|n-n_*|^2+n|v|^2)$ and $n_*q^\ast|\Delta \Fsp ||v|$ by $C(n_*^\gamma|\Delta \Fsp |^2+n|v|^2)$.  The constants use the comparability bound $n_*\le n/\vartheta$. The correction $n(\chi-1)$ and the bounded comparison forms are controlled by the kinetic term. Thus $|\mathscr R_*|\le C\mathscr D$ for comparable densities.

\vskip.2cm \noindent {\bf 3.} \emph{The remainder when the densities are separated.} For noncomparable densities, we use the original affine remainder~\eqref{eq:physical-nonlinear-original-remainder}. Its internal stress terms are bounded by $C(n+n_*)^\gamma$. Bounded form increments and~\eqref{uc:eq:geometric-coefficients} bound the form terms by $Cb_0(t)n_*^\gamma$. For the mass stress $T_m={m \, n}\,u\otimes u$, direct contraction gives
\bse
\be
\label{uc:eq:mass-far-states}
 -(\nabla_\nu w_{*\mu})\Delta T_m^{\mu\nu}
 =-\frac{{m \, n}}{H_*}(\nabla_\nu u_{*\mu})u^\mu u^\nu
 -\frac{m}{H_*^2}(n\chi u^\nu-n_*u_*^\nu)\del_\nu H_*.
\ee
\be
 \bigl|(\nabla_\nu w_{*\mu})\Delta T_m^{\mu\nu}\bigr|
 \le C\bigl(n|v|^2+nq^\ast|v|+(n+n_*)(q^\ast)^2\bigr)
 \le C\bigl(n|v|^2+(n+n_*)^\gamma\bigr).
\ee
\ese
Indeed the reference-direction contraction uses $u_*\cdot dH_*=\mathcal {\Ocal}((q^\ast)^2)$, and the mixed contractions use $u_*\cdot\nabla u_*,dH_*=\mathcal {\Ocal}(q^\ast)$. Young's inequality gives the last bound since $(q^\ast)^2=n_*^{\gamma-1}$. For separated densities,~\eqref{uc:eq:weighted-coercivity} controls $(n+n_*)^\gamma$ and $n|v|^2$. Thus the full bound here is $C(1+b_0(t))\mathscr D$.

At zero reference density, the coefficient formulas $P=\varepsilon_{\rm e}q^{p_\theta+1}\Psi_{z_F}$ and their two-form counterparts have zero first derivatives at $(q,z)=(0,0)$. Thus $C_*,D_*$ and their first derivatives vanish, $dH_*=0$, and the normalized evolution gives $u_*\cdot\nabla u_*=0$. The mass estimate and the internal factor $n^\gamma$ give the weighted bound directly. If both states have zero density, the remainder and relative density are both zero, even when their kinematic labels differ. This proves~\eqref{uc:eq:weighted-bound}. Proposition~\ref{uf:prop:local-convex-entropy} gives the entropy and its supporting covectors, and Proposition~\ref{uf:prop:entropy-multipliers} gives the reference regularity. The local solution satisfies the physical balance equations, spatial involutions, continued material equations and prescribed initial conditions. Both solutions take their values in the bounded constitutive range on the periodic slab, so their fluxes are integrable.  Theorem~\ref{thm:physical-entropy-vacuum} gives the integrated estimate. If the initial relative entropy vanishes,~\eqref{uc:eq:weighted-coercivity} gives equality of the current and stress, and of the minors in matter.
This completes the proof of Theorem~\ref{uc:thm:common}.
\end{proof}

At $\varepsilon_{\rm e}=0$, the form multipliers vanish, so~\eqref{uc:eq:geometric-coefficients} is automatic in the reduced system. The normalized fluid equations give $\mathcal A_*=0$, including at zero density. On the bounded periodic comparison branch of this section, the density and mass estimates give comparison for these fluid laws  with the $H^s$ reference regularity of Proposition~\ref{uf:cor:fractional-existence}, for $5/2<s\le3$, including at vacuum.

%============================================================================

\paragraph{Construction of the corrected integral.}
{
The Taylor formula for $L-L_{\rm fl}$ above, at $V=0$, gives $n\Psi(n^\theta,n^\theta Y_{\rm R})/n^\gamma\to\tfrac12Y_R^T\Psi_{zz}(0,0)Y_R>0$ at every fixed $Y_R\ne0$. The internal energy has the same order because $2(2+\theta^{-1})A(0)=H_f(0)>0$.
}
We write \(f_*,f_{\mathrm w}\) for Lipschitz configuration potentials on the affine cover, \(\Fsp_*=D_x f_*\), \(r_*={\particlecurrent}_*^0\), \(r={\particlecurrent}_{\mathrm w}^0\), \(h=f_{\mathrm w}-f_*\), and \(\Delta\) for the difference between the weak solution and the reference. The common affine periods are \(A\), and the additive constants in the initial potentials are part of the normalization. We take the exterior closure and exact minor identities of the augmented formulation. On the finite-velocity branch, the spacetime forms are bounded through $r=0$. In matter their time components are given by $F_0^A=-V^iF_i^A$ and $G_{0iB}=-V^jG_{jiB}$, and these bounds pass to the completed branch. The time--space closure equations give the prescribed initial data set, in the sense of distributions, of each $F_i$. This trace is bounded and spatially closed by the boundedness and spatial closure of $F_i$ in the slab. Thus it has an affine-periodic Lipschitz potential $f_{\rm in}$ on the cover. For each solution set
\be
 f(t,x)=f_{\rm in}(x)+\int_0^t F_0(\tau,x)\,d\tau .
\ee
Time--space closure and its prescribed initial term give $\del_i f=F_i$ in the sense of distributions. The potential is Lipschitz, with uniform initial data set $f_{\rm in}$ and the same affine periods, since $F_0$ is periodic.

{
Put $R=r+r_*$ and let $\mathscr D$ be the relative entropy density and $E$ its spatial integral. The reference weights $q=n_*^\theta$ and $Q=r_*^\theta$ are comparable. Proposition~\ref{uf:prop:local-convex-entropy} and the shape comparison in~\eqref{uf:eq:weighted-rest-remainder} yield
\be
 R^{\gamma-2}|\Delta r|^2+R^\gamma\bigl(|\Delta F|^2+|\Delta G|^2\bigr)+r|V-V_*|^2\le C\mathscr D .
\ee
All terms vanish at $R=0$. For comparable positive densities, the time components follow from $F_0^A=-V^iF_i^A$ and $G_{0iA}=-V^jG_{jiA}$, and the factor $R^\gamma\le Cr$ controls the resulting velocity difference. For separated densities use~\eqref{uc:eq:separated-density-bound} and boundedness of the forms. Where $r=0$, products with $rV$ are read as $\particlecurrent^{\rm sp}$. Integration gives the asserted weighted bound.
}

\paragraph{Displacement at the same material point.}

{
Both particle numbers equal $M=\det A$. For $f=Ax+g$, with $g$ periodic, the integral of $\det(A+\tau Dg)$ is constant in $\tau$ by the Piola identity and periodic integration by parts. The formula passes from smooth $g$ to $W^{1,3}$ by convolution and convergence of determinants in $L^1$. If $M=0$, both densities vanish, so the current, stress, relative entropy and correction multipliers vanish on this branch, including for the initial data. Set $L=0$ in this case. If $M>0$, use the material torus $\mathcal B_A=\mathbb R^3/A\mathbb Z^3$ with its Euclidean volume.
}

For either configuration \(f\) with positive total particle number and \({\particlecurrent}^0=\det D_xf\ge0\), we have 
\be
 \int_{\mathbb T^3} g(f(x))\det D_x f(x)\,dx=\int_{\mathcal B_A}g(y)\,dy .
\ee
For smooth maps this follows by the homotopy from \(Ax\) to \(f\) and Stokes' theorem applied to the top-degree form \(g(y)\,dy\). Convolution of the periodic part gives uniform convergence of the maps and strong convergence of their derivatives in every finite \(L^p\), with a uniform bound on the derivatives. This proves the identity for Lipschitz maps and smooth \(g\).  The Lipschitz potential is defined on the closed slab. On any slice, take a sequence of times with $\det D_x f\ge0$ tending to that slice. The uniform convergence of the potentials and the bounded gradients give distributional convergence of the spatial minors by the exterior product identity used in Proposition~\ref{prop:collapsed-material-coordinates}. Hence $\det D_x f\ge0$ on every slice. The absolute Jacobian in the area formula is therefore this determinant, which equals $r$ almost everywhere in spacetime and gives the prescribed density at $t=0$. Applying the area formula in flat coordinate charts (cf.~Cibotaru and de Lira, \href{https://arxiv.org/pdf/1305.2968\#page=12}{Theorem 4.1(b)}), we obtain the following identity, where \({\operatorname{mult}}_f(y)\) counts preimages modulo the spatial lattice:
\be
 \int_{\mathcal B_A}g(y){\operatorname{mult}}_f(y)\,dy=\int_{\mathcal B_A}g(y)\,dy
\ee
for every smooth \(g\). Hence \({\operatorname{mult}}_f=1\) almost everywhere. The nondifferentiability and zero determinant sets have null images by the Lipschitz property and the area formula. We write \(X_f(t,y)\) for the lifted inverse with \(f(t,X_f(t,y))=y\). The displacement \(X_f-A^{-1}y\) is bounded, periodic, and measurable: it is the density of the pushforward of \((x-A^{-1}f)r\,dx\,dt\) relative to \(dy\,dt\). On the initial slice the same construction uses \(r_0\,dx\) and \(dy\). If the two initial densities coincide, the area formula sends a common full-particle-measure set to a full-measure set of material points. On these sets equality of the configurations is equivalent to equality of the lifted inverses. For zero total particle number both equalities are vacuous.

The exact-minor current is conserved with its initial data set by Proposition~\ref{prop:kinematic-weak-local}(ii), and ${\particlecurrent}^\mu\del_\mu f=0$ by the cubic minor identity. For a periodic smooth \(g\) on \(\mathcal B_A\), these two identities give, component by component, 
\be
 \del_\mu\big({\particlecurrent}^\mu g(f)(x-A^{-1}f)\big)=g(f){\particlecurrent}^{\rm sp}.
\ee
{The current has the initial value $r_{\rm in}=\det D_x f_{\rm in}$ by the spatial minor relations in the initial data set.}
\paragraph{{Evolution of the inverse material map.}}
In the continuity equation the periodic Lipschitz multiplier $g(f)(x-A^{-1}f)$ is admissible by smooth approximation. Extend the periodic part of $f$ constantly in time across $t=0$ and mollify it locally in time and periodically in space. The multipliers then converge uniformly up to $t=0$, their derivatives converge in $L^1$, and the derivatives remain uniformly bounded. Since $\particlecurrent$ is bounded, these properties justify the limit in both the bulk and initial terms. The area formula on the slices then gives
\be
 \del_t(X_f(t,y)-A^{-1}y)= V(t,X_f(t,y))
\ee
as a distribution in \((t,y)\) including the initial displacement prescribed by \(f_{\rm in}\). In particular, for a smooth periodic \(g\) and a smooth time test \(\psi\) vanishing near the final time, the distributional identity is
\be
 \int_0^T\!\int_{\mathcal B_A}g(y)\bigl[(X_f-A^{-1}y)\psi'
 +V(t,X_f)\psi\bigr]\,dy\,dt
 +\psi(0)\int_{\mathcal B_A}g(y)(X_f(0,y)-A^{-1}y)\,dy=0.
\ee
Integrating the bounded right-hand side with the prescribed initial position gives a representative that is Lipschitz in time in $L^\infty$, after subtracting $A^{-1}y$. For almost every material point it agrees with the inverse almost everywhere in time. Since $f$ and those trajectories are continuous, $f(t,X_f(t,y))=y$ for all times for those labels. On each slice the area formula gives uniqueness for almost every material point, so this representative also agrees there with the inverse constructed on that slice.

The area formula gives $\int_{\mathcal B_A}|(V-V_*)(t,X_{\mathrm w})|^2dy\le CE$. It also shows that the pushforward of $dy/M$ by $(X_{\mathrm w},X_*)$ modulo the spatial lattice has marginals $r\,dx/M$ and $r_*\,dx/M$. Its squared transport cost is at most $2L/M$, since the distance on the torus is bounded by the distance of the chosen lifts. This applies to the prescribed initial densities by the same area formula. The position equations make $\sqrt{2L}$ absolutely continuous. The Lipschitz bound on $V_*$ and Cauchy--Schwarz give, for almost every time,
\be
 \left|\frac{d}{dt}\sqrt L\right|\le C\sqrt E+C\sqrt L,
 \qquad |L'|\le C(E+L).
\ee
{
With the integrating factor $e^{-C_Vt}$ for the first inequality and $E(t)\le Ce^{C_Et}(E(0)+L(0))$, one obtains $L(t)\le Ct^2e^{2t\max\{C_V,C_E/2\}}E(0)$ when $L(0)=0$. The area formula also gives
}
\be
 \int (r+r_*)|h|^2\,dx\le C L .
\ee
Indeed \(h(X_{\mathrm w})=f_*(X_*)-f_*(X_{\mathrm w})\) and \(h(X_*)=f_{\mathrm w}(X_*)-f_{\mathrm w}(X_{\mathrm w})\). The area formula and the spatial Lipschitz bounds of the direct maps give separately
\be
 \int r|h|^2\,dx\le2\operatorname{Lip}(f_*)^2L,\qquad
 \int r_*|h|^2\,dx\le2\operatorname{Lip}(f_{\mathrm w})^2L.
\ee
In particular, when the initial configurations agree with the same normalization, $L(0)=0$ and
\be
 \left(\int_{\mathbb T^3}(r+r_*)|f_{\mathrm w}-f_*|^2\,dx\right)^{1/2}
 \le Ct e^{Ct}\sqrt{E(0)}
\ee
for almost every $t$. The two weights apply on their respective matter regions, even when the cavities differ. The initial normalizations and the spacetime gradient bounds control \(\|h\|_\infty\). With equal initial spatial gradients we normalize the two potentials so that \(h(0)=0\).

\paragraph{Coercivity of the corrected integral.}

We use the original relative current \(\mathscr Q\) and the reference coefficients \(w_*,C_*,D_*\) defined in~\eqref{eq:main-vacuum-finite-multipliers}. For the divergence coefficients and their linear combination we write 
\be
 c_A^\alpha=\del_\nu C_{*A}^{\nu\alpha},\qquad
 d_B^{\alpha\beta}=\del_\nu D_{*B}^{\nu\alpha\beta},\qquad
 s_A^\alpha=c_A^\alpha+2d_B^{\alpha\beta}\Omega_{BAC}F_{*\beta}^C .
\ee
Here \(\Omega\) is the material volume tensor, and antisymmetric indices use full sums. We use \(\beta_A=s_A^0/r_*\) in matter, continued by zero where $r_*=0$. Proposition~\ref{uf:prop:entropy-multipliers} gives \( |\nabla C_*|+|\nabla D_*|\le C r_*q\). Since \(F_*\) is bounded, this yields \(|d|+|c|\le C r_*q\) and \(|\beta|\le Cq\). In Section~\ref{app:integrated-3d} we construct a Lipschitz coefficient $\widetilde\beta$ with $|\widetilde\beta|\le Cq$ and $|\widetilde\beta-\beta|\le Cq^2$. We define the corrected time integral by
\be
 \mathcal H\coloneqq E+\int d_B^{0i}\Omega_{BAC}h^A\Delta F_i^C\,dx
                -\int\widetilde\beta_A\Delta r\,h^A\,dx .
\ee
For any \(\delta\in(0,1)\), the defining expression for $\mathcal H$, Young's inequality and the relative entropy estimate give the bound below. For the form correction, use \(r_*q\le C r_*^{\gamma/2}r_*^{1/2}\), which follows from \(q\le Cr_*^{(\gamma-1)/2}\). For the density correction, write \(q|\Delta r|=(r+r_*)^{(\gamma-2)/2}|\Delta r|\,q(r+r_*)^{(2-\gamma)/2}\) where \(r+r_*>0\). The square of the last factor is bounded by \(C r_*^{\gamma-1}(r+r_*)^{2-\gamma}\). Both corrections vanish where \(r_*=0\). Young's inequality therefore gives
\be
 |\mathcal H-E|\le \delta E+
 C_\delta\int\Big(r_*+r_*^{\gamma-1}(r+r_*)^{2-\gamma}\Big)|h|^2\,dx .
\ee
Where $R\coloneqq r+r_*>0$, the second weight equals $R(r_*/R)^{\gamma-1}\le R$; at $R=0$ it vanishes. The displacement estimate gives $|\mathcal H-E|\le\delta E+C_\delta L$. Fixing $\delta<1$ and taking $\Lambda>C_\delta$ gives both bounds on $\mathcal E\coloneqq\mathcal H+\Lambda L$ for every $\gamma>1$.

We also need a single estimate for the differences of the current and of the stress. Write $T=T_m+T_{\mathcal L}$ as in~\eqref{uf:eq:stress-decomposition}. On the bounded branch we claim
\be
\label{eq:local-weighted-current-and-stress}
 R^{\gamma-2}\bigl(|\Delta\particlecurrent|^2+|\Delta T|^2\bigr)+|\Delta T_{\mathcal L}|^2\le C\mathscr D ,
\ee
with the left-hand side defined to be zero at $R=0$. Indeed, the current and mass stress differences are bounded by $C(|\Delta r|+r|\Delta V|)$. The density weight in~\eqref{uf:eq:local-relative-entropy-comparison} controls $R^{\gamma-2}|\Delta r|^2$, and $rR^{\gamma-2}\le R^{\gamma-1}$ controls the velocity term. For comparable positive densities, integration of the constitutive derivatives along the allowed segments gives
\be
 |\Delta T_{\mathcal L}|\le C\bigl(R^{\gamma-1}|\Delta r|+R^\gamma(|\Delta Y|+|\Delta V|)\bigr).
\ee
After squaring, the factors relative to the entropy weights $R^{\gamma-2}, R^\gamma,r$ are bounded by $CR^\gamma, CR^\gamma, CR^{2\gamma-1}$. They stay bounded after multiplication by $R^{\gamma-2}$ since $\gamma>1$. For separated densities $|\Delta T_{\mathcal L}|\le CR^\gamma$ and $\mathscr D\ge cR^\gamma$, which bounds both squares, including when one density vanishes. This proves~\eqref{eq:local-weighted-current-and-stress}.

For $\gamma\le2$ the weighted square controls the ordinary square on this bounded branch. For $\gamma\ge2$ use $|\Delta\particlecurrent|+|\Delta T|\le CR$ to control the powers of order $\gamma$. Integration and interpolation with the uniform bounds on the torus give the asserted $L^p$ estimates; the unweighted internal-stress estimate starts at $p=2$.

Section~\ref{app:integrated-3d} proves the differential estimate for $\mathcal E$ with the prescribed initial data sets. The density term uses the same Taylor remainder of a power law as in the polytropic pressure estimate for classical fluids~\cite[Theorem~3.4 and (3.27)]{LeFloch-JHDE}.

\paragraph{{Consequences of the integrated estimate.}}
{The estimate controls the effect of the initial vacuum velocity and the size of the regions where the densities differ.}
\begin{proof}[Proof of Corollary~\ref{cor:initial-vacuum-velocity}]
{Where $Q_0>0$, equality of the initial configurations and velocities gives equality of the exact minors and the stress. Where $Q_0=0$, all terms in~\eqref{uf:eq:local-relative-entropy-comparison} vanish. Thus $E(0)=0$. The normalized initial configurations agree almost everywhere for the common particle measure $Q_0^{1/\theta}\,dx$. Part~(iv) of Theorem~\ref{uc:thm:integrated-3d} gives $L(0)=0$, so $E=L=0$ for almost every time.

The relative entropy estimate identifies $\particlecurrent,T$ and the minors in matter. To identify the configurations themselves, fix a time at which $L(t)=0$. For almost every material point $y$, the two lifted inverse positions coincide. At that common position both configurations take the value $y$. The area formula thus gives $f_{\mathrm w}=f_*$ almost everywhere for the particle measure $r(t,x)\,dx$. The configurations are continuous and $r>0$ on the open matter region, so they agree throughout that region. Spacetime continuity extends these equalities from almost every time to every slice, including the endpoints. Since $V={\particlecurrent}^{\rm sp}/{\particlecurrent}^0$ in matter, the velocities agree there as well. The velocities, minors and configurations are continuous in space, so they agree on the closure of the matter region. From any initial point in matter both trajectories stay in matter and solve the same Lipschitz velocity equation. The flows therefore coincide there and, by continuity, on its closure. The boundary formula follows from~\eqref{uf:eq:inertial-vacuum-flow} and the fact that each flow is a homeomorphism.}
\end{proof}

\begin{proof}[Proof of Corollary~\ref{cor:integrated-cavity-comparison}]
By~\eqref{uf:eq:local-relative-entropy-comparison}, the integral of $W\coloneqq(r+r_*)|r^\theta-r_*^\theta|^2$ is bounded by $CE(t)$. The corrected integral controls $E(t)$ in terms of $E(0)+L(0)$. On $\{r\le\lambda r_*,\,r_*>0\}$, we have $W\ge(1-\lambda^\theta)^2r_*^\gamma$, with the symmetric bound on the disjoint set where $r_*\le\lambda r$ and $r>0$. This gives~\eqref{eq:integrated-low-density-comparison}, including for the prescribed initial densities. Taking $\lambda=0$ gives~\eqref{eq:integrated-cavity-comparison}, and restricting to densities at least $\delta$ gives the volume bound. H\"older's inequality gives the particle number bound. For two local solutions, the integral of $W$ and the relative entropy integral are continuous in time. The bound on $\int W\,dx$ thus holds on every slice, and the pointwise bounds on the comparison sets apply on each slice.

{
By~\eqref{uf:eq:weighted-rest-remainder}, $W\asymp R^{\gamma-2}|\Delta r|^2\asymp |r^{\gamma/2}-r_*^{\gamma/2}|^2$, which proves~\eqref{eq:integrated-density-power-comparison}. For each $d\in\{|\Delta r|,|\Delta\particlecurrent|,|\Delta T|\}$, we have $d\le CR$ and $\int R^{\gamma-2}d^2\,dx\le CE$ by this comparison and~\eqref{eq:local-weighted-current-and-stress}. The integrands are zero at $R=0$. For $\gamma\ge2$, this gives $\|d\|_\gamma\le CE^{1/\gamma}$. For $1<\gamma<2$, H\"older's inequality with exponents $2/\gamma$ and $2/(2-\gamma)$ gives
\be
 \int d^\gamma\,dx\le
 \left(\int R^{\gamma-2}d^2\,dx\right)^{\gamma/2}
 \left(\int R^\gamma\,dx\right)^{1-\gamma/2}.
\ee
On $\{r\ge2r_*\}$, $W$ bounds $R^\gamma$, and elsewhere $R\le3r_*$. In particular,
\be
\label{eq:elastic-density-sensitive-rate}
 \|\Delta\particlecurrent(t)\|_{L^\gamma}+\|\Delta T(t)\|_{L^\gamma}
 \le C E(t)^{1/2}\bigl(E(t)+\|r_*(t)\|_{L^\gamma}^\gamma\bigr)^{1/\gamma-1/2}.
\ee
The same pointwise bounds apply to the initial fields. Substituting the estimate for $E(t)$ in terms of $E(0)+L(0)$ completes the proof.
}
\end{proof}

%============================================================================

{
\subsection{\texorpdfstring{The quadratic elastic law}{The quadratic elastic law}}
\label{subsec:proof-second-density-weight}

{We prove Theorem~\ref{thm:quadratic-second-weight-augmented-stability} and Corollary~\ref{cor:nonlinear-second-weight} together, from the normalized evolution to the entropy estimate and the treatment of rank loss.}

\vskip.2cm \noindent {\bf 1.} \emph{The Hessian in the rest frame and the inertia matrix.}
We put $p=p_\theta=1/\theta$ and use the rest fields $W_F=q\Frest$, $W_G=q\Grest$, so $z=qW$ also at $q=0$. For the next differentiation take $n>0$.  Let $a=\dot q$ and $b_F=\dot z_F$, $b_G=\dot z_G$ be induced by an affine variation of $(n,\Frest,\Grest)$.  Direct differentiation of $q=n^{1/p}$ and $z=q^2x$ gives
\be
 \ddot q=\frac{1-p}{q}a^2,
 \qquad
 \ddot z=\frac4q ab-\frac{2(p+2)}{q^2}za^2.
\ee
After division by $n$, the normalized second variation of the rest energy is
\be
\label{eq:quadratic-second-weight-rest-hessian}
\begin{aligned}
 &\alpha|b_F+(p+2)W_Fa|^2
 +\beta|b_G+(p+2)W_Ga|^2\\
 &\quad+\Bigl(2\kappa p-(p+2)(p+4)
 (\alpha|W_F|^2+\beta|W_G|^2)\Bigr)a^2.
\end{aligned}
\ee
Indeed, the fluid term $\kappa n^\gamma/\gamma=\kappa q^{p+2}/\gamma$ contributes $2\kappa p\,a^2$.  This follows either by differentiating $n^\gamma$ along the affine $n$--variation or by using $\ddot q=(1-p)a^2/q$ in the second derivative of $q^{p+2}$.  A form block contributes
\be
 \alpha|b+(p+2)(z/q)a|^2
 -\alpha(p+2)(p+4)|z/q|^2a^2.
\ee
{The strict margin~\eqref{eq:quadratic-second-weight-rest-margin} therefore gives a uniform positive lower bound for the normalized second variation of the energy with respect to the rest-frame variables.}  Substitution in the general inertia formula~\eqref{eq:rest-inertia-test}, followed by division by $n$, gives
\be
\label{eq:quadratic-second-weight-inertia}
\begin{aligned}
 \mathsf M_R={}&\left(m+\kappa q^2
 +\frac{\alpha(2\gamma-1)}2|z_F|^2
 +\frac{\beta(2\gamma+1)}2|z_G|^2\right)I_3\\
 &+\alpha z_F^Tz_F-\beta z_G^Tz_G.
\end{aligned}
\ee
Since $|z_Gb|^2\le|z_G|^2|b|^2$, we obtain $\mathsf M_R\ge(m+\kappa q^2)I_3$. The normalized Hessian and inertia matrix remain positive definite down to $q=0$. At $W_F=W_G=0$, both elastic transport columns vanish and the acoustic form is exactly the fluid form. At $q=0$ on the product cone, all columns vanish. The strict fluid acoustic margin and continuity on the compact density and velocity ranges give a common positive time matrix for sufficiently small elastic fields.

\vskip.2cm \noindent {\bf 2.} \emph{Normalized equations and local construction.}
{We set $Y'=(W_F,W_G)$ and $Z'=(Z_F,Z_G)=QY'$. For both results let $W=q(\Frest,\Grest)$ denote the corresponding fields in the material rest frame. The rest energy is
\be
 n\left(m+\frac\kappa\gamma q^2+\Psi(z)\right),\qquad
 z=(z_F,z_G)=qW.
\ee
The quadratic theorem uses $\Psi(z)=(\alpha|z_F|^2+\beta|z_G|^2)/2$. For the nonlinear corollary, define the symmetric matrices $A(z)=\int_0^1D^2\Psi(tz)\,dt$ and $S(z)=D^2\Psi(z)$, so $D\Psi(z)=A(z)z$. Put $c=p+2$, $h_0=2\kappa p$ and $t=k\cdot b$. The map $\ell_k(z)b$ has one-form entries $k(z_F^A\cdot b)$ and two-form entries $z_{G,A}t-b(z_{G,A}\cdot k)$. Direct differentiation gives
\bse
\be
\label{eq:q2-explicit-rest-blocks}
 \mathsf H=\begin{pmatrix}h_0-2cW^TA(z)W&cW^TA(z)\\cA(z)W&S(z)\end{pmatrix},
 \qquad D_kb=\binom{\theta qt}{\ell_k(z)b+2\theta zt}.
\ee
\be
\label{eq:q2-explicit-coupling}
 (\mathsf H D_kb)_q=\theta h_0qt+cW^TA(z)\ell_k(z)b,\qquad
 (\mathsf H D_kb)_z=S(z)\ell_k(z)b+\theta cD\Psi(z)t+2\theta S(z)zt.
\ee
\ese
For the quadratic law $A=S=\operatorname{diag}(\alpha I,\beta I)$, and the last two terms are $\theta(p+4)Szt$, recovering the quadratic contraction. The $C^6$ assumption on $\Psi$ gives $C^4$ dependence of $A$, $S$ and the resulting principal matrices. Three coefficient derivatives suffice for the $H^3$ commutators; a fourth controls their local Lipschitz dependence in $H^3$ and the auxiliary $H^4$ estimates used when removing the smoothing in the non-relativistic limit. No derivative of $\Psi$ beyond order six is used in these arguments. In the nonlinear case, $\mathsf H(0,0)=\operatorname{diag}(h_0,D^2\Psi(0))>0$, so a sufficiently small raw-form bound preserves a uniform positive margin. The inertia formula~\eqref{uf:eq:rest-inertia}, with $\theta z\cdot\Psi_z$ replaced by $2\theta z\cdot\Psi_z$, is $C^5$ in $(q,z)$ and equals $(m+\kappa q^2)I$ at $z=0$. The cancellation in the first component uses $z=qW$ exactly. Thus $K_k=\mathsf H D_k$ is $C^4$ in $(q,z,z\otimes W)$, while $D_l^TK_k$ is $C^4$ in $(q,z)$ because $q(z\otimes W)=z\otimes z$.}

The spatial constraints are imposed on the original forms. In particular, we do not impose closure on $W_F$ or $W_G$. At positive density the stress calculation in~\eqref{uf:eq:rest-constraint-terms} uses the Jacobian of $(n,F,G)\mapsto(q,q^2F,q^2G)$ and gives $\mathcal J_yL_k=D_k$. Its extra terms are the original curls of $F$ and divergences of the dual $G$, multiplied by their conjugates, so they vanish on the original closure constraints. The additional divergence terms in the weighted transport are already included in $2\theta zt$ in~\eqref{eq:q2-explicit-rest-blocks}.

For the observer change, we replace the form blocks of $\mathsf T$ in~\eqref{uf:eq:rest-lab-map} by $\Gamma^{-2\theta}P_V^{-1}$ and $\Gamma^{-2\theta-1}P_V$, and denote the resulting map by $\mathsf T_2$. With $\mathsf R_2=\operatorname{diag}(\Gamma^{-\theta},\mathsf T_2)$ and $C=\Gamma P_V$, the same block congruence gives
\bse
\be
\label{eq:q2-observer-blocks}
 \mathsf H_L=\Gamma\mathsf R_2^T\mathsf H\mathsf R_2,
 \quad\mathsf M_L=\Gamma C^T(\mathsf M_R-D_V^T\mathsf HD_V)C,
 \quad\mathsf B_i=\mathsf R_2^TK_{k_i}C,
 \quad\mathsf E_i=-C^T(D_V^TK_{k_i}+D_{k_i}^TK_V)C.
\ee
\be
\label{eq:q2-observer-transport}
 (\Pi_i b)_Q=\theta Qb_i,
 \quad(\Pi_i b)_{Z_F^A}=e_i(Z_F^A\cdot b)+2\theta Z_F^Ab_i,
 \quad(\Pi_i b)_{Z_{G,A}}=(1+2\theta)Z_{G,A}b_i-bZ_{G,A,i}.
\ee
\be
\label{eq:q2-observer-symmetry}
 \mathsf H_L\Pi_i=\mathsf B_i,
 \quad\mathsf M_L\mathsf Q_i=\mathsf B_i^T,
 \quad\mathsf M_L\mathsf R_i=\mathsf E_i=\mathsf E_i^T,
 \qquad\mathsf Q_i=\mathsf M_L^{-1}\mathsf B_i^T,
 \quad\mathsf R_i=\mathsf M_L^{-1}\mathsf E_i.
\ee
\ese
These formulas verify the symmetrization for the new weights. The principal matrices are $C^4$ in $(Q,Z_F,Z_G,Z\otimes W,V)$; only the density--form time block also uses $W$. All the relative spatial matrices vanish at $Q=Z=0$.

Let $D_t=\partial_t+V\cdot\nabla$.  The kinematic equations imply
\bse
\be
\label{eq:quadratic-second-weight-transport-a}
 D_tQ=-\theta Q\operatorname{div}V,\qquad
 D_t\Fsp=-(DV)^T\Fsp,\qquad
 D_t\Gsp=(DV)\Gsp-(\operatorname{div}V)\Gsp,
\ee
\be
\label{eq:quadratic-second-weight-transport-b}
 D_tW_F=-(DV)^TW_F-\theta(\operatorname{div}V)W_F,\qquad
 D_tW_G=(DV)W_G-(1+\theta)(\operatorname{div}V)W_G.
\ee
\ese
The equations for $Z_F=QW_F$ and $Z_G=QW_G$ follow by another multiplication by $Q$.  If $X$ is the particle flow, $B=DX$ and $j=\det B$, then
\bse
\be
\label{eq:quadratic-second-weight-flow-a}
 Q\circ X=j^{-\theta}Q_0,\quad
 \Fsp\circ X=B^{-T}\Fsp_0,\quad
 W_F\circ X=j^{-\theta}B^{-T}W_{F0},\quad
 Z_F\circ X=j^{-2\theta}B^{-T}Z_{F0},
\ee
\be
\label{eq:quadratic-second-weight-flow-b}
 \Gsp\circ X=j^{-1}B\Gsp_0,\quad
 W_G\circ X=j^{-1-\theta}BW_{G0},\quad
 Z_G\circ X=j^{-1-2\theta}BZ_{G0}.
\ee
\ese
These identities propagate both product relations and keep the ranks of $F$ and $G$ constant along particle trajectories.  {In the observer frame, write $W=(W_F,W_G)$ and $Z=(Z_F,Z_G)=QW$. The principal coefficients involve the entries of $Z\otimes W$, including the mixed products $Z_F\otimes W_G$ and $Z_G\otimes W_F$ when the energy couples the two form families. Differentiating $Z=QW$ gives, for all pairs of components and also at $Q=0$,}
\be
\label{eq:quadratic-second-weight-product}
 d(Z\otimes W)=dZ\otimes W+W\otimes dZ-W\otimes W\,dQ.
\ee
{ Since $H^2(\mathbb T^3)$ is an algebra, the right-hand side of~\eqref{eq:quadratic-second-weight-product} belongs to $H^2$ when $U_{\rm n}\in H^3$ and $W\in H^2$. The identity therefore gives the estimate}
\be
\label{eq:q2-product-tame-bound}
 \|Z\otimes W\|_{H^3}\le C_K\|Z\|_{H^3}+C_K\|Q\|_{H^3},
 \qquad \|W\|_{H^2}+\|U_{\rm n}\|_{H^3}\le K.
\ee
{On both domains, $W_0=Q_0(\Fsp_0,\Gsp_0)\in H^3$ by the initial product estimates. Thus $W\in H^2$ and the algebra estimate above gives $H^3$ principal coefficients. The transport sources $DV\,W_0$ and $V\cdot DW_0$ have $H^2$ norms bounded by $C\|V\|_{H^3}\|W_0\|_{H^3}$; the fixed-scale flow composition also uses only $W_0\in H^3$.} The material derivative of the time block requires only the transport equations for $W$ and $U_{\rm n}$, so it is bounded in $L^\infty$. {In the differentiated energy, the identities~\eqref{eq:q2-observer-symmetry} give exactly the paired fluxes~\eqref{uf:eq:paired-fluxes}, now with the second-weight matrices. Spatial integration by parts therefore differentiates $\mathsf B_i$ and $\mathsf E_i$; it does not require an $H^3$ bound for $\mathsf H_L$. For the convolved velocity equation, the extra term is paired in $H^1$--$H^{-1}$ using~\eqref{uf:eq:metric-commutator}. Its bound requires only the spatial Lipschitz norm of $\mathsf M_L$, which follows from its principal-coefficient bound.} The remaining commutators use the $H^3$ principal-coefficient bounds and $H^2$ lower-field bounds just established. Thus no fourth derivative of $Q$ or $V$ occurs.

We use the fixed-scale positive approximation of Section~\ref{uf:subsec:construction-vacuum}.  Only the velocity equation is convolved, while~\eqref{eq:quadratic-second-weight-flow-a}--\eqref{eq:quadratic-second-weight-flow-b} reconstruct the density and form variables.  The positive time matrix and~\eqref{eq:quadratic-second-weight-product} give
\be
 \frac d{dt}\mathcal E_3\le C\mathcal E_3,
 \qquad c\|U_{\rm n}\|_{H^3}^2\le\mathcal E_3\le C\|U_{\rm n}\|_{H^3}^2.
\ee
The transport equations give the accompanying $H^2$ estimates for the lower fields.  The maps $(Q,Z)\mapsto Z\otimes Z/Q$, extended by zero at $(0,0)$, are Lipschitz on every cone $|Z|\le KQ$.  Their difference estimates close the order-one Cauchy energy exactly as in~\eqref{uf:eq:product-difference}.  {The approximations therefore converge at the lower orders, including at $Q=0$; the limit belongs to the asserted class, with strong time continuity given by the endpoint argument.}  The flow formulas propagate spatial closure, the cavity and both ranks.  The normalized momentum equation, current conservation and the constitutive identity~\eqref{uf:eq:independent-stress-projection} give the stress balance laws.  On the cavity the internal stress and its first derivatives vanish, while the continued velocity equation reduces to $D_tV=0$.  This proves part~(i).

\vskip.2cm \noindent {\bf 3.} \emph{Coercivity and the constitutive remainder.}
{We construct the mathematical entropy directly on the bounded comparison branch. Write $Y=(\Fsp,\Gsp)$, $r={\particlecurrent}^0$, and let $\mathsf A(V)$ be the linear transformation from spatial to rest-frame form variables in~\eqref{uf:eq:rest-shapes}. For either law, set}
\be
 {L(r,Y,V)\coloneqq\bar\rho(r/\Gamma(V),\mathsf A(V)Y),\qquad
 E=L-V\cdot L_V,\quad P=-L_V.}
\ee
{The Hessian estimates in the material rest frame from Steps~1 and~2 and the strict margins on the comparison domain give on the density--form and velocity segments used in the comparison:}
\be
 {D^2_{(r,Y)}L[(\dot r,\dot Y)]^2\asymp
 r^{\gamma-2}|\dot r|^2+r^{2\gamma-1}|\dot Y|^2,\qquad
 -D_V^2L\asymp rI,\quad L_r\asymp1,\qquad r>0.}
\ee
{The velocity estimate follows from the strict fluid time-direction margin, uniformly on the velocity ball and its comparison segments after the elastic range is chosen sufficiently small. The bounds hold on a slightly enlarged range with the same properties. For a positive-density reference let $U_{\rm c}=(Y,E,P)$ and $\lambda_*=(L_{Y*},-1,V_*)/L_{r*}$. Direct substitution gives}
\be
 {\begin{aligned}
 L_{r*}\bigl(-r+r_*-\lambda_*\cdot(U_{\rm c}-U_{{\rm c}*})\bigr)
 ={}&L(r,Y,V_*)-L_*-D_{(r,Y)}L_*\cdot(r-r_*,Y-Y_*)\\
 &+L(r,Y,V)-L(r,Y,V_*)-L_V(r,Y,V)\cdot(V-V_*).
 \end{aligned}}
\ee
{The first remainder is nonnegative by convexity in $(r,Y)$ on {these segments}, and the second by concavity in $V$. Their integral formulas, with $R=r+r_*$, give the weights $R^{\gamma-2}|r-r_*|^2$, $R^{2\gamma-1}|Y-Y_*|^2$ and $r|V-V_*|^2$. These estimates extend to zero density: the density exponent $\gamma-2$ is integrable at an endpoint and the form energy vanishes there.}

{The supremum of the affine functions $-r_*+\lambda_*\cdot(U_{\rm c}-U_{{\rm c}*})$ over positive-density reference fields in the bounded branch is a convex function equal to $-r$ on its conservative image. At zero reference density, the limiting supporting covector is $(0,-\Gamma_*/m,\Gamma_*V_*/m)$; the bounded branch and the prescribed continuation select $V_*$. The same support inequality gives the value $0$ at every zero-current boundary field. Thus the entropy and the supporting covectors used in Theorem~\ref{thm:physical-entropy-vacuum} are defined on the completed branch, including the prescribed initial data.}

{Using $Z=Q^2Y$, $Q=r^\theta$, and the bounded form range converts these weights into the following comparison.} On the bounded branch, the relative mathematical entropy satisfies
\be
\label{eq:quadratic-second-weight-coercivity}
 D_\eta\asymp(r+r_*)
 \bigl(|Q-Q_*|^2+|Z_F-Z_{F*}|^2+|Z_G-Z_{G*}|^2\bigr)
 +r|V-V_*|^2.
\ee
For comparable positive densities, with $R=r+r_*$, this implies
\be
\label{eq:quadratic-second-weight-physical-coercivity}
 D_\eta\ge c\left(
 R^{\gamma-2}|r-r_*|^2+R^{2\gamma-1}
 (|F-F_*|^2+|G-G_*|^2)+R|V-V_*|^2\right).
\ee
For separated densities it controls $R^\gamma+r|V-V_*|^2$. {For both laws, with $X=(\Frest,\Grest)$, the reference conjugates and their differentials satisfy
\bse
\be
\label{eq:nonlinear-second-weight-conjugates}
 \bar\rho_X=n^\gamma D\Psi(z),\qquad
 d\bar\rho_X=n^\gamma\left(S(z)\,dz+\frac\gamma\theta A(z)W\,dq\right).
\ee
\be
\label{eq:quadratic-second-weight-multipliers}
 |C_*|+|D_*|+|\nabla C_*|+|\nabla D_*|\le Cn_*^\gamma.
\ee
\ese
The first identity follows from $z=q^2X$; its differential uses $D\Psi(z)=A(z)z$ and $z=qW$. Thus the multiplier bound uses only derivatives of normalized strong fields. {For a $C^6$ nonlinear energy,} Taylor's formula gives $n\Psi(q^2X)=n^{2\gamma-1}\chi(n^{\gamma-1},X)$, with bounded derivatives of $\chi$ on the chosen ranges. Taylor's integral formula for $H_a=n^a\phi(n^{\gamma-1},X)$, with $a\ge\gamma>1$ and bounded coefficient derivatives, gives}
\be
\label{eq:quadratic-second-weight-taylor}
 |H_a-H_{a*}-DH_{a*}\cdot(n-n_*,X-X_*)|
 \le C\bigl(R^{\gamma-2}|n-n_*|^2+R^a|X-X_*|^2\bigr).
\ee
When $1<\gamma<2$, the factor with exponent $\gamma-2$ is integrable along an affine segment ending at zero.  Formula~\eqref{eq:quadratic-second-weight-taylor} controls the fluid term with $a=\gamma$ and the elastic contribution with $a=2\gamma-1$.  In the comparable region the remaining mixed form fluxes satisfy
\be
 n_*^\gamma|V-V_*|\bigl(|F-F_*|+|G-G_*|\bigr)\le CD_\eta
\ee
by~\eqref{eq:quadratic-second-weight-physical-coercivity}. We now estimate the whole source, including the form terms. We use a fixed Lorentz coframe at rest with respect to the reference and the independent chart defined in~\eqref{eq:q2-independent-chart} below. We write $a=2\gamma-1$ and $\Delta X=(\Delta F,\Delta G,\Delta r,\Delta V)$. In the comparable-density region, the internal stress has the derivative bounds
\bse
\be
\label{eq:q2-internal-stress-derivatives}
 |\del_r^2T_{\rm int}|\le CR^{\gamma-2},\qquad
 |\del_rD_{F,G}T_{\rm int}|\le CR^{a-1},\qquad
 |D_{F,G}^2T_{\rm int}|\le CR^a.
\ee
\be
\label{eq:q2-internal-velocity-derivatives}
 |\del_rD_VT_{\rm int}|\le CR^{\gamma-1},\qquad
 |D_VD_{F,G}T_{\rm int}|\le CR^a,\qquad
 |D_V^2T_{\rm int}|\le CR^\gamma.
\ee
\ese
{These estimates follow by twice differentiating the fluid term and $n^a\chi(n^{\gamma-1},X)$, with $n=r/\Gamma(V)$ and the rest transformations linear in the spatial forms. The normalized positive Hessian gives~\eqref{eq:quadratic-second-weight-coercivity} for both laws; the observer velocity estimate follows by the same small-elastic-perturbation argument.} Their coefficients are uniform on a slightly enlarged bounded branch containing the chart segments. Weighted Young inequalities absorb the mixed derivatives: the relevant exponents satisfy $a\ge\gamma$ and $a-1\ge(a+\gamma-2)/2$. The form evolution fluxes are bilinear in $(F,V)$ or $(G,V)$, so their Taylor remainders are bounded by $C|\Delta V|(|\Delta F|+|\Delta G|)$. Their coefficients have size $\Ocal(n_*^\gamma)$ by~\eqref{eq:quadratic-second-weight-multipliers}, which is exactly the geometric mean of the form and velocity weights.

For the mass stress, we set $f_*(r,V)=-(\nabla_\nu w_{*\mu})T_m^{\mu\nu}(r,V)$. It is linear in $r$, its second velocity derivatives are bounded by $Cr$, and in the reference rest frame
\be
\label{eq:q2-mass-mixed-derivative}
 |\del_rD_Vf_*(r,V_*)|\le Cq_*.
\ee
Indeed this mixed derivative contains only the reference acceleration and the spatial derivative of $H_*$, both $\Ocal(q_*)$. At an intermediate velocity it is bounded by $C(q_*+|V-V_*|)$. The first term is controlled by the density and kinetic weights, since $q_*\simeq R^{(\gamma-1)/2}$; the second contributes at most $C|\Delta r||\Delta V|^2\le CR|\Delta V|^2$. Taylor's formula and the vanishing first differential proved in the next step therefore bound the complete remainder by~\eqref{eq:quadratic-second-weight-physical-coercivity}. No minor identity is used in this calculation.

In the separated region we estimate the original affine source. The internal stress and the reference form terms are bounded by $CR^\gamma$, using bounded unweighted comparison forms and~\eqref{eq:quadratic-second-weight-multipliers}. For the mass stress, the direct contraction~\eqref{uc:eq:mass-far-states} applies without a minor relation: $D_{u_*}u_*=\Ocal(q_*)$, $dH_*=\Ocal(q_*)$ and $D_{u_*}H_*=\Ocal(q_*^2)$ give $C(r|\Delta V|^2+R^\gamma)$. These are controlled by the separated-density entropy bound. At $r_*=0$, the free reference evolution gives the same estimate directly; at $r=r_*=0$ both the source and relative density vanish. Hence the complete relative-current remainder obeys
\be
\label{eq:quadratic-second-weight-remainder}
 |\mathscr R_*(t,x;\mathcal Z)|\le CD_\eta(\mathcal Z\mid\mathcal Z_*)
\ee
{ for every augmented comparison solution and its associated constitutive stress in the bounded constitutive range. The estimate uses the fields of the weak solution without differentiating them.}

\vskip.2cm \noindent {\bf 4.} \emph{Comparison and exact minor relations.}
The reference satisfies the augmented balance laws and spatial involutions of Definition~\ref{def:weak-local-complete}. We verify the first-order cancellation in the independent coordinates
\be
\label{eq:q2-independent-chart}
 X=(\Fsp,\Gsp,r,V),\qquad
 F_0^A=-\Fsp^A\cdot V,\qquad
 G_{0iA}=-V^jG_{jiA},\qquad
 \particlecurrent=r(1,V).
\ee
The map reconstructs the forms satisfying the orthogonality condition~\eqref{eq:450bis-2} directly and uses no inverse of a rank-dependent constraint. At $r>0$, the entropy multiplier identity~\eqref{eq:entropy-multiplier-differential} and its spatial counterparts, with the exterior equations included, read
\be
\label{eq:q2-chart-entropy-compatibility}
 d_X\mathscr I^\mu
 =\lambda^a d_X\mathscr F_a^\mu-b^\ell d_X\mathscr K_\ell^\mu,
 \qquad \mathscr K_\ell^0=0.
\ee
Here the involution coefficients are the spatial antisymmetric components of the form multipliers. The identity follows by applying the differential constitutive identity to each spacetime index; the remaining antisymmetric terms are precisely $dF$ and $dG$. Taking the exterior derivative in $X$, and contracting with $\del_\mu X_*$, gives
\be
\label{eq:q2-chart-adjoint-cancellation}
 (\del_\mu\lambda_*^a)d_X\mathscr F_{a*}^\mu
 -(\del_i b_*^\ell)d_X\mathscr K_{\ell*}^i
 =d_X\lambda_*^a\,\del_\mu\mathscr F_{a*}^\mu
 -d_Xb_*^\ell\,\del_i\mathscr K_{\ell*}^i=0.
\ee
The last equality uses the strong balance laws and spatial involutions. Hence the first differential of the total remainder vanishes in every independent chart direction. This calculation does not introduce the covector $\nu_*$ of~\eqref{eq:general-entropy-adjoint-identity}.

At reference vacuum the form multipliers and their first spacetime derivatives vanish by~\eqref{eq:quadratic-second-weight-multipliers}. Also $H_*=m$, $dH_*=0$ and $w_*=u_*^\flat/m$. The derivative of the internal stress in the chart vanishes because its density powers are $\gamma>1$ and $2\gamma-1>1$. The only possible first derivative is that of the mass stress in the $r$ direction, proportional to $u_*\otimes u_*$. Its contraction with $\nabla w_*$ vanishes by $u_*^\mu\nabla_\nu u_{*\mu}=0$. The first-order cancellation therefore holds also in the interior of the cavity, and extends to its boundary through the continuous coefficient formulas. The finite comparison remainder at such a reference is controlled directly by the continued equation $D_tV_*=0$ and the mass estimate, as in the separated-density calculation above.

Consequently,~\eqref{eq:quadratic-second-weight-remainder} is the total remainder required in Theorem~\ref{thm:physical-entropy-vacuum}. The supports and involution coefficients are continuous up to $t=0$: {the form coefficients are $n^\gamma A(z)z=q^{p+2}A(z)z$ times smooth functions of the velocity and $H^{-1}$, and $H=m+\kappa q^2+\Psi(z)+2\theta z\cdot D\Psi(z)\ge m$ on the chosen convex neighborhood of the origin.} Their first spacetime derivatives are bounded and vanish on the zero-density set for the form coefficients. {Thus they are valid test coefficients in Lemma~\ref{thm:common-current-jhde}, without differentiating the weak fields. We first identify its initial term without assuming the constitutive entropy identity for the weak initial data.}

{The bounded branch gives bounded fluxes. Testing their balance laws with a time cutoff supported in $[0,\varepsilon]$ shows that the time averages $\overline U_\varepsilon$ converge to $U_0$ weakly-* in $L^\infty$. For each nonnegative smooth spatial function, the particle inequality gives a lower bound for its pairing with $\overline{\particlecurrent^0}_\varepsilon$ in terms of its pairing with $\particlecurrent_0^0$, up to an error tending to zero. Every affine support $\ell$ in the construction of $\eta$ satisfies $\ell(U)\le-\particlecurrent^0$ in the slab. Averaging and passing to the initial terms gives $\ell(U_0)\le-\particlecurrent_0^0$. A countable family of supports has the same supremum, so $\mathfrak d_0\ge0$ almost everywhere. The initial relative current in the lemma is therefore}
\be
 {\int_{\mathbb T^3}\left(-\particlecurrent_0^0+\particlecurrent_{*,0}^0-\lambda_{*,0}\cdot(U_0-U_{*,0})\right)dx
 =\int_{\mathbb T^3}(D_{\eta,0}+\mathfrak d_0)\,dx.}
\ee
{The lemma and~\eqref{eq:quadratic-second-weight-remainder} give the estimate with the defect. The defect vanishes for initial fields and stress in the completed branch, proving~\eqref{eq:quadratic-second-weight-stability} in that case. If instead $U(t)\to U_0$ strongly in $L^1$, the bounded supporting covectors make $\eta$ Lipschitz, so the relative entropy integral converges to its initial value. Applying the comparison on $(\tau,t)$ and letting $\tau\downarrow0$ again gives~\eqref{eq:quadratic-second-weight-stability}. This alternative does not require the initial particle-production atom to vanish.}

{The coercivity above identifies the current and stress almost everywhere for positive times, and the fields $F$ and $G$ almost everywhere in matter, whenever the right-hand side of either applicable estimate vanishes. For the local strong solution, the flow formulas propagate the exact minor relations from initial data satisfying them, which proves part~(iii).} When $\alpha=\beta=0$, all form multipliers and form rows disappear, and the same argument uses only three Sobolev derivatives of the regularized fluid density and velocity. 
{\vskip.2cm \noindent {\bf 5.} \emph{Powers of the singular values at rank loss.}
We now check the alternative~\eqref{eq:second-weight-spectral-additions}, whose energy need not be $C^6$. At exponent four, $I_4(X)=\operatorname{tr}((X^TX)^2)$ is a polynomial, so each quartic summand may be kept unchanged in the approximation. It remains to treat a spectral exponent $s>4$. Put $\Phi_s(X)=I_s(X)/s$ and $g_s=D\Phi_s$. The contracted derivative
\be
\label{eq:spectral-contracted-flux-regularity}
 \mathcal K_s(X,L)=Dg_s(X)[XL]
 \quad\hbox{is }C^3\hbox{ on the full matrix space}.
\ee
Its bounds on bounded sets are independent of the smallest singular value. To prove this, use the symmetric dilation and the infinitesimal congruence
\bse
\be
 S_X=\begin{pmatrix}0&X\\X^T&0\end{pmatrix},\qquad
 \Omega_L=\operatorname{diag}(0,L^T),\qquad
 \Omega_LS_X+S_X\Omega_L^T
 =\begin{pmatrix}0&XL\\L^TX^T&0\end{pmatrix}.
\ee
\be
\label{eq:spectral-divided-difference-cancellation}
 t_j f^{[r]}(t_0,\ldots,t_r)
 =(tf)^{[r]}(t_0,\ldots,t_r)
 -f^{[r-1]}(t_0,\ldots,\widehat t_j,\ldots,t_r).
\ee
\ese
Here $f^{[r]}$ is the divided difference of order $r$, with its continuous value at repeated nodes; the hat in the last line means omission. For smooth $f$, the derivatives of $f(S)$ are sums of chains of matrix entries multiplied by these divided differences. This follows by differentiating polynomials and approximating the highest scalar derivative. In each derivative of $Df(S)[\Omega S+S\Omega^T]$, the distinguished matrix entry supplies one eigenvalue factor. Formula~\eqref{eq:spectral-divided-difference-cancellation} removes it. The repeated-integral bound for divided differences therefore gives
\be
 \bigl\|D_S^j\{Df(S)[\Omega S+S\Omega^T]\}\bigr\|
 \le C\|\Omega\|\bigl(\|f\|_{C^3}+\|tf\|_{C^4}\bigr),
 \qquad 0\le j\le3,
\ee
on a fixed bounded spectral interval. This bound also applies to differences of scalar functions. Take $f(t)=t|t|^{s-2}$ and approximate it by $f_\varepsilon(t)=t(t^2+\varepsilon^2)^{(s-2)/2}$. One has $f_\varepsilon\to f$ in $C^3$ and $tf_\varepsilon\to tf$ in $C^4$: near zero the highest derivatives are bounded by $C(|t|+\varepsilon)^{s-4}$, and away from zero convergence is smooth. The upper-right block gives~\eqref{eq:spectral-contracted-flux-regularity}, including at every rank loss.

The form row of $K_k$ is consequently $C^3$, since $D^2\Phi_s(X)X=(s-1)g_s(X)$. The density row requires a weighted product estimate. For the acoustic right action $XL$, put $J(X)=g_s(X):XL$. The identity $X^Tg_s(X)=(X^TX)^{s/2}$ shows that $J$ is $C^4$. With $z=qW$, the density row equals $cJ(z)/q$ in matter. The following formula defines and differentiates it also at $q=z=0$:
\bse
\be
 F_J=\int_0^1(1-t)D^2J(tz)[z,W]\,dt,\qquad
 B_J=\int_0^1D^2J(tz)\,dt,\quad
 C_J=\int_0^1(1-t)D^2J(tz)\,dt,
\ee
\be
\label{eq:spectral-weighted-density-flux}
 dF_J=B_J[W,dz]-C_J[W,W]\,dq.
\ee
\ese
Indeed, set $A_J(z)=\int_0^1DJ(tz)\,dt=C_Jz$. Then $F_J=A_J[W]$, $DA_J+C_J=B_J$, and the Sobolev product rule with $q\,dW=dz-W\,dq$ proves the last line without dividing by $q$. Since $B_J,C_J$ are $C^2$, its right-hand side is $H^2$ for $q,z\in H^3$ and $W\in H^2$. This gives the $H^3$ coefficient bound. For the fixed-scale construction, $F_J=A_J(z)[W]$ and $A_J\in C^3$ give local Lipschitz dependence in $H^2$; local Lipschitz dependence in $H^3$ is not needed here. The free Hessian is $C^2$, but this suffices for the spatial $H^2$ regularity and bounded material derivative of the time weight. The inertia and the contractions $D_l^TK_k$ are $C^3$; the observer formulas~\eqref{eq:q2-observer-blocks} therefore preserve the grouped energy estimates. These checks replace the generic $C^6$ coefficient argument; they do not assert $C^3$ dependence of every coefficient on the formal independent product variables.

For a smooth constitutive approximation use
\be
 \Phi_{s,\varepsilon}(X)
 =\frac1s\operatorname{tr}\bigl((X^TX+\varepsilon^2I)^{s/2}
                         -\varepsilon^sI\bigr).
\ee
It is convex and vanishes with its gradient at zero. The preceding divided-difference estimate gives uniform $C^3$ bounds and convergence for its contracted fluxes; its scalar contraction $J_\varepsilon$ converges in $C^4$. The positive quadratic response keeps the time margin uniform. Formula~\eqref{eq:spectral-weighted-density-flux} gives uniform $H^3$ bounds without a density floor. {We give the joint limit argument below for both Sobolev orders.}

For the fractional conclusion, fix $5/2<s<3$ and $p_{\sigma j}>s+1$. If no spectral term is present, the smooth quadratic argument applies directly. Otherwise, we choose
\be
 s-2<\mu<\min\bigl\{1,\min_{\sigma,j}(p_{\sigma j}-3)\bigr\}.
\ee
For a spectral exponent $3<a<4$, the preceding scalar approximation and divided-difference bound give
\be
 \|\mathcal K_{a,\varepsilon}-\mathcal K_a\|_{C^2}
 \le C\varepsilon^{a-3},\qquad
 \|\mathcal K_{a,\varepsilon}\|_{C^3}
 \le C(1+\varepsilon^{a-4}).
\ee
Put $\alpha=a-3$. For matrices at distance $h$, comparison with the approximation at $\varepsilon=h$ gives a $Ch^\alpha$ bound on the increment of $D^2\mathcal K_a$. For $D^2\mathcal K_{a,\varepsilon}$ the same bound follows by comparison with the limit when $h\ge\varepsilon$, and from $Ch\varepsilon^{\alpha-1}\le Ch^\alpha$ when $h<\varepsilon$. Thus every contracted flux is $C^{2,\mu}$ on the bounded matrix range, uniformly at all ranks; for exponents at least four this follows from polynomial smoothness or the $C^3$ bound. The scalar contraction $J$ is $C^{3,\mu}$ and the free Hessian is $C^{1,\mu}$. Approximation converges in any smaller H\"older exponent still greater than $s-2$.

These regularities suffice for the fractional Sobolev estimates. Indeed a bounded $C^\mu$ function is a multiplier on $H^{s-2}$ when $\mu>s-2$: in the double integral defining the fractional seminorm, its increment contributes an integrable factor $|x-y|^{2\mu-2(s-2)-3}$. Applying this bound to
\be
 D_x^2 F(u)=DF(u)D_x^2u+D^2F(u)[D_xu,D_xu]
\ee
gives the $H^s$ composition estimate for {$F(u)-F(0)$}, with $F\in C^{2,\mu}$ on the bounded range of $u\in H^s\subset W^{1,\infty}$. {On $\mathbb R^3$ the constant part is subtracted; on the torus it belongs to $H^s$.} The mean-value formula gives local differences in $H^{s-1}$. In~\eqref{eq:spectral-weighted-density-flux}, {the compositions $B_J(z)-B_J(0)$ and $C_J(z)-C_J(0)$ belong to $H^{s-1}$; their constant parts are bounded multipliers.} Hence that identity gives the density-row bound in $H^s$ with only $W\in H^{s-1}$. The fixed-scale source differences use $F_J=A_J(z)[W]$ as before.

{For the removal of the constitutive smoothing, this estimate also applies uniformly to fields depending on $\varepsilon$. For a spectral exponent $a=p_{\sigma j}$, fix the matrix $L$ and set $J(X)=g_a(X):XL$ and $J_\varepsilon(X)=D\Phi_{a,\varepsilon}(X):XL$. Suppose that $q_\varepsilon,z_\varepsilon$ are bounded in $H^s$ and $W_\varepsilon$ in $H^{s-1}$, with $z_\varepsilon=q_\varepsilon W_\varepsilon$ and with $z_\varepsilon$ in a fixed matrix ball $B_R$. Then, writing $F_J(z,W)$ for the integral defined immediately before~\eqref{eq:spectral-weighted-density-flux}, one has}
\be
 {\left\|F_{J_\varepsilon}(z_\varepsilon,W_\varepsilon)-F_J(z_\varepsilon,W_\varepsilon)\right\|_{H^s}
 \le C\|J_\varepsilon-J\|_{C^{3,\mu}(B_R)}\longrightarrow0.}
\ee
{Indeed, apply~\eqref{eq:spectral-weighted-density-flux} to $J_\varepsilon-J$ and use the $H^{s-1}$ algebra estimate, including the constant parts of $B_{J_\varepsilon-J}$ and $C_{J_\varepsilon-J}$ as bounded multipliers on $\mathbb R^3$. This argument uses $W_\varepsilon$ only in $H^{s-1}$ and does not divide by $q_\varepsilon$, so it remains valid at zero density. The same conclusion holds for a matrix $L_\varepsilon$ depending on the fields if $L_\varepsilon-L_{\rm c}$ is uniformly bounded in $H^s$ for a constant matrix $L_{\rm c}$: apply the $H^s$ product estimate componentwise.}

{We now remove the smoothing of the constitutive law and the velocity convolution together. Write $s=3$ for the integer case and retain $5/2<s<3$ for the fractional case, with the respective restrictions on the spectral exponents. For $s=3$ and $a>4$, the same integral identity gives $\|F_{J_\varepsilon}-F_J\|_{H^3}\le C_K\|J_\varepsilon-J\|_{C^4(B_R)}$ on these fields by the $H^2$ product estimate. In the smoothed equations the radial contraction is $Dg_{a,\varepsilon}(z)[z]=\mathcal K_{a,\varepsilon}(z,I)$. The homogeneity identity for $g_a$ applies only after taking the limit. Use the preceding contracted estimates for $\mathcal K_{a,\varepsilon}(z,L)$ in each form row and the weighted estimate for $F_{J_\varepsilon}$ in the density row. In $D_l^TK_k$ the density row is multiplied by a factor proportional to $q$, so this contribution involves $J_\varepsilon(z)$ itself. The rest inertia matrix involves $\Phi_{a,\varepsilon}$ and the entries of $z^Tg_{a,\varepsilon}(z)$, which are the same contractions with fixed elementary matrices.}

{In~\eqref{eq:q2-observer-blocks}, the changes of frame are smooth functions of $V$ on $|V|\le v_{\max}<1$. Their entries differ from their values at $V=0$ by uniformly bounded $H^s$ functions. Products give the same bounds, relative to constant matrices, for $\mathsf B_i,\mathsf E_i$ and $\mathsf M_L$. Inversion preserves these bounds by the uniform positive lower bound for $\mathsf M_L$. The upper time matrix $\mathsf H_L$ also contains the less regular field $W$. Along the particle flow its derivative is bounded by the ordinary chain rule in $(q,z,W,V)$, using the uniform first derivatives of the free Hessian, the transport equations and~\eqref{uf:eq:approximate-material-derivative} at order $s-1>3/2$. Spatial integration by parts differentiates only the paired flux matrices. The convolution commutator in~\eqref{uf:eq:metric-commutator} uses $\nabla\mathsf M_L\in L^\infty$. These bounds apply at $q=0$ and are uniform in both smoothing parameters.}

{At fixed $\varepsilon,\ell$, use the flow equations~\eqref{uf:eq:label-ODE} with $X-\mathrm{Id}\in H^6$ and $w\circ X\in H^s$. The composition difference uses one derivative of the transported initial factor, as in~\eqref{uf:eq:composition}, now from $H^s$ to $H^{s-1}$. This applies to $Q_0,Z_0,W_0$ on either domain. At fixed constitutive smoothing the resulting force differences are locally Lipschitz into $H^{s-2}$ by the product estimate with an $H^{s-1}$ multiplier. Spatial convolution supplies the derivatives required in the flow equations. The resulting local construction starts from the given initial fields.}

{Apply the paired energy estimate to $(y,w)$ at order $s$, where $y=(Q,Z)$, $w$ is the velocity before convolution and $V=\smoothoperator_\ell w$. Its spatial commutators use the $H^s$ bounds relative to the constant matrices and the Lipschitz bounds just obtained. Together with transport at order $s-1$ for $W$, this gives the common interval by the strict initial margins. In particular the force $\mathcal F$ of~\eqref{uf:eq:approximation} is uniformly bounded in $H^{s-1}$. In the equation for $V$, the first kernel moment gives}
\be
{
 \|(I-\smoothoperator_\ell^2)\mathcal F\|_2+
 \|[V\cdot\nabla,\smoothoperator_\ell]w\|_2
 \le C\ell^{-1}\bigl(\|\nabla\mathcal F\|_2+
              \|\nabla V\|_\infty\|\nabla w\|_2\bigr).}
\ee
{For the commutator, use $|V(x)-V(x-h)|\le\|\nabla V\|_\infty |h|$ under the convolution integral. Thus the error tends to zero independently of the constitutive smoothing.}

{Take any sequence $\varepsilon\downarrow0$, $\ell\to\infty$. The equations bound $\partial_t(y,V)$ in $H^{s-1}$ and $\partial_t W$ in $H^{s-2}$. Local compactness and interpolation give a subsequence converging in $C_tH^{s'}_{\rm loc}$ for $(y,V)$ and in $C_tH^{s'-1}_{\rm loc}$ for $W$, where $5/2<s'<s$. On $\mathbb R^3$ we exhaust space by balls. The normalized fields and their first spatial derivatives therefore converge locally uniformly, and $W$ converges locally uniformly as well. The integral defining $F_J$, the contracted fluxes and the change of frame identify every coefficient on the limit fields, also at $q=z=0$. The convolution error above tends to zero in $L^2$, so the normalized equations pass to the limit in the sense of distributions. The convergence includes $t=0$, where $V(0)=\smoothoperator_\ell V_0$ in the approximations. The flow formulas recover the unweighted forms from their prescribed initial values. In physical variables the flow transports $r_0\,dx$ and gives $r(1,V)$ its initial value in the sense of distributions. For the stress, the internal part is $r\tau(q,z,V)$ with $\tau=\mathcal O(q^2+|z|^2)$ and $D_{q,z}\tau=\mathcal O(q+|z|)$, as in~\eqref{uf:eq:stress-decomposition}. It is $C^1$ with zero differential where $Q=0$ since $|z|\le C Q$ and $r=Q^{p_\theta}$. The rest-mass part uses the transported density and the $C^1$ velocity. A time cutoff against a spatially compact test-function therefore recovers the prescribed initial stress also on $\mathbb R^3$.}

{For uniqueness one may use the paired difference energy at order zero. On the cone $|z|\le Cq$ the density coefficient $J(z)/q$ is Lipschitz in $(q,z)$: the bounds $|J(z)|\le C|z|^2$ and $|DJ(z)|\le C|z|$ control both derivatives in matter and give its continuous value zero at the vertex. The other coefficients of the equations solved for the time derivatives are Lipschitz on the same bounded range by the contracted estimates and the smooth change of frame. Their differences in $L^2$ are thus controlled by the difference of the normalized fields, and multiply bounded first derivatives of the reference. Gronwall's inequality gives uniqueness for these fields, and the flow determines the remaining fields. Finally the limit inherits the global Sobolev bounds by weak compactness. Applying spatial convolution to its differentiated equations gives continuity of the positive energy between arbitrary times as in Step~6, using the same paired fluxes and bounded material derivatives of the time weights. For the fractional order use $\Lambda^s$ in place of the third derivatives. The limiting equations give global time bounds one order lower, hence continuity by interpolation below the top orders, including $L^\infty$ continuity of the weights. The transport energy gives $W\in C_tH^{s-1}$, and applying it to the difference of each unweighted form from its prescribed initial value gives the same continuity for that difference. More explicitly, for a spatial form $Y$ let $Y^{\rm inc}=Y-Y_0$ and write $\mathcal L(DV)Y$ for the right-hand side of its transport equation in~\eqref{eq:quadratic-second-weight-transport-a}. On $\mathbb R^3$ one has
\be
 D_tY^{\rm inc}=\mathcal L(DV)Y^{\rm inc}
             +\mathcal L(DV)Y_0-V\cdot\nabla Y_0,\qquad Y^{\rm inc}(0)=0 .
\ee
The last two terms have $H^{s-1}$ norm bounded by $C\|Y_0\|_{W^{3,\infty}}\|V\|_{H^s}$ since $s-1\le2$. The global transport estimate therefore gives $\|Y^{\rm inc}(t)\|_{H^{s-1}}\le C t$ near zero; the initial form itself may have infinite $L^2$ norm. Weak continuity and continuity of the energy give $C_tH^s$ for the normalized fields up to the initial time, and their equations give $C_t^1H^{s-1}$.}

{Finally, for each spectral exponent $p_{\sigma j}>3$ in either alternative, we write $\nu_{\sigma j}\coloneqq1+p_{\sigma j}(\gamma-1)$ and use homogeneity to obtain}
\be
 {n\Phi_{p_{\sigma j}}(q^2X)=n^{\nu_{\sigma j}}\Phi_{p_{\sigma j}}(X),\qquad
 \nu_{\sigma j}>2\gamma-1>\gamma.}
\ee
{The form derivatives of $\Phi_{p_{\sigma j}}$ through order three are bounded on the chosen matrix range, including at rank loss. In observer variables the corresponding internal stress has the form $T_{\mathrm{int},\sigma j}=r^{\nu_{\sigma j}}\mathcal S_{\sigma j}(Y,V)$, where $\mathcal S_{\sigma j}$ and its derivatives through order two are bounded on the comparison range. Consequently, for $r>0$,}
\be
 {|\del_r^{i}D_{(Y,V)}^kT_{\mathrm{int},\sigma j}|\le C r^{\nu_{\sigma j}-i},\qquad i+k\le2.}
\ee
{These weighted estimates imply the stress derivative and Taylor bounds of Step~3; they allow density derivatives to be singular at vacuum. In the fractional range $s+1<p_{\sigma j}<4$ they also avoid requiring the auxiliary function $\chi$ of Step~3 to be $C^2$ at zero. The conjugate formula~\eqref{eq:nonlinear-second-weight-conjugates} still gives the $n_*^\gamma$ multiplier bound. The complete remainder, its cancellation at reference vacuum and its initial terms are thus controlled as before, including when the cavities differ. The positive quadratic part is used in coercivity throughout this argument.}

}
}

%============================================================================

\part[Applications and technical details]{\LARGE{Applications and technical details}}
\label{part:applications-and-technical-details}

\section{\texorpdfstring{The non-relativistic limit with cavitation}{The non-relativistic limit with cavitation}}
\label{sec:non-relativistic-vacuum}
 
\subsection{\texorpdfstring{Relativistic and classical constitutive energies}{Relativistic and classical constitutive energies}}

\paragraph{Physical time and the internal energy.}
{
{We take the non-relativistic limit of the normalized equations used in Theorem~\ref{uf:thm:main}. Uniform estimates down to zero density give existence on an interval independent of the light speed and convergence after subtraction of the rest energy.\par}}

{Beig and Schmidt~\cite{BeigSchmidt} obtained the classical equations in the limit of infinite light speed and described their symmetric hyperbolic formulation near unstrained configurations.} We denote the speed of light by $c_{\rm light}$ and use inertial coordinates $x^0=c_{\rm light}t$, with $\Fsp,\Gsp$ denoting the independent spatial form matrices. The relation $\Gsp=\operatorname{cof}\Fsp$ is imposed only for the restriction to exact spacetime minors. The physical velocity is denoted by $v$. We keep $m,\theta,\Phi,\Psi$ and the elastic parameter $\varepsilon_{\rm e}$ of the local construction, and define
\bse
\be
\label{nr:eq:parameter}
 \delta\coloneqq c_{\rm light}^{-2},\qquad
 \Gamma_\delta\coloneqq(1-\delta|v|^2)^{-1/2},\qquad
 P_\delta\coloneqq I+\frac{\delta\Gamma_\delta^2}{\Gamma_\delta+1}v\otimes v.
\ee
\be
\label{nr:eq:constitutive-law}
 \raug_{\delta,\varepsilon_{\rm e}}\coloneqq n\bigl(m/\delta+\Phi(q)+\Psi(q,z)\bigr),\qquad
 q=n^\theta,\quad z=\varepsilon_{\rm e}q(\Frest ^A,\Grest _A)_{A=1}^3,\quad \delta>0.
\ee
\ese
The $\delta$-independent functions $\Phi,\Psi$ satisfy Definition~\ref{def:density-weighted-elastic-law} and the positive lower bounds on $H_f,\Psi_{zz}(q,0)$ in~\eqref{uf:eq:constitutive-bounds} for $0\le q\le q_*$. We use the variables in~\eqref{uf:eq:scaled-fields} with physical velocity $v$:
\be
\label{nr:eq:unknowns}
 r\coloneqq Q^{p_\theta}=n\Gamma_\delta,\qquad
 Y=\varepsilon_{\rm e}(\Fsp ^A,\Gsp _A)_{A=1}^3,\qquad Z=QY,\qquad
 U_\delta\coloneqq(Q_\delta,Z_\delta,v_\delta),\qquad p_\theta=1/\theta.
\ee
The normalized equations at $\delta=0$ define $U_{\rm nr},Y_{\rm nr}$. We impose strict initial margins as in the local construction, using $|v_0|\le v_0^{\max}<v_*^{\max}<+\infty$. Taking $\delta_0(v_*^{\max})^2<1$ and reducing $\delta_0$ and the elastic bound gives positive time matrices on the common range.

\paragraph{Classical constitutive energies.}
The internal energy in the classical current-volume description is
\be
\label{nr:eq:rest-internal-energy}
 \bar\rho_{\rm nr}^{\rm int}(\Fsp ,\Gsp ,r)\coloneqq r\bigl(\Phi(r^\theta)+\Psi(r^\theta,\varepsilon_{\rm e}r^\theta(\Fsp ,\Gsp ))\bigr).
\ee
The Hessian calculation~\eqref{uf:eq:rest-metric} gives its convexity on convex portions of the local independent-minor domain; subtracting the affine rest-mass energy leaves that Hessian unchanged. For zero elastic parameter the absent form directions are omitted. On positive-density physical states, $d=\det(\nabla_X y)>0$ and $r=1/d$. The change of volume density in Proposition~\ref{prop:Ball-perspective-new} gives the energy per reference volume in the independent variables $(A,B,d)$ of Section~\ref{sec:Ball-comparison}:
\be
\label{nr:eq:ball-energy}
 {\rho_{B,\rm nr}^{\rm int}}(A,B,d)\coloneqq d\bar\rho_{\rm nr}^{\rm int}(B^T/d,A^T/d,1/d)
 =\Phi(d^{-\theta})+\Psi\bigl(d^{-\theta},\varepsilon_{\rm e}d^{-\theta}(B^T/d,A^T/d)\bigr).
\ee
This is the convex perspective in Ball's independent deformation minors, on the corresponding convex domains. The kinetic term is $m|v|^2/2$ per reference volume. The perspective comparison is made where the inverse configuration exists; the normalized evolution and the flux convergence in  Section~\ref{sec:non-relativistic-limit-evolution} below include $r=0$.

%-----------------------------------------------------------------------------------------------------------------------------------------

\subsection{\texorpdfstring{The convergence theorem}{The convergence theorem}}
\label{sec:non-relativistic-limit-evolution}

\begin{theorem}[Non-relativistic limit with cavitation]
\label{nr:thm:vacuum-limit}
Fix the constitutive functions in~\eqref{nr:eq:constitutive-law}, with the regularity, positivity and strict range margins specified in Section~\ref{sec:non-relativistic-vacuum}, and the elastic parameter $\varepsilon_{\rm e}\ge0$. Let $\mathcal D=\R^3$ or $\mathbb T^3$ and $5/2<s\le3$. Prescribe nonnegative $Q_{0,\delta}$ and physical velocities $v_{0,\delta}$ in a bounded subset of $H^s(\mathcal D)$, and scaled elastic fields $Y_{0,\delta}$ in a bounded subset of $W^{3,\infty}(\mathcal D)$, satisfying the smallness and strict range margins above. Put $Z_{0,\delta}=Q_{0,\delta}Y_{0,\delta}$. For $\varepsilon_{\rm e}>0$, impose the independent initial closure constraints of Theorem~\ref{uf:thm:main}; at $\varepsilon_{\rm e}=0$, set $Y_{0,\delta}=0$. Suppose, as $\delta\downarrow0$, that
\be
\label{nr:eq:initial-convergence}
 (Q_{0,\delta},v_{0,\delta})\longrightarrow(Q_{0,\rm nr},v_{0,\rm nr})\quad\hbox{in }H^s,
 \qquad Y_{0,\delta}\longrightarrow Y_{0,\rm nr}\quad{\hbox{in }L^\infty}.
\ee
At $s=3$, also suppose that $Z_{0,\delta}\to Z_{0,\rm nr}\coloneqq Q_{0,\rm nr}Y_{0,\rm nr}$ in $H^3$. There exist $\delta_0,T>0$ such that the following conclusions hold for $0<\delta\le\delta_0$.
\begin{itemize}
\item[(i)] {\emph{Existence on a common interval.}} The normalized relativistic solutions and the classical solution exist uniquely on $[0,T]$, with coefficients extended continuously to zero density and bounds uniform in $\delta$ in the class
\be
\label{nr:eq:regularity}
 U\in C_tH^s\cap C_t^1H^{s-1},\qquad Y-Y_0\in C_tH^{s-1}\cap C_t^1H^{s-2}.
\ee
The interval and bounds are independent of a positive lower density and of the elastic parameter under the fixed initial bounds.
\item[(ii)] {\emph{Convergence in the normalized variables.}} The density coordinate, weighted elastic fields and velocity, collected in $U=(Q,Z,v)$, and the increments of the scaled elastic fields converge as follows:
\be
\label{nr:eq:top-convergence}
 {U_\delta\longrightarrow U_{\rm nr}\quad\hbox{in }C_tH^s\cap C_t^1H^{s-1},}\qquad
 {(Y_\delta-Y_{0,\delta})\longrightarrow(Y_{\rm nr}-Y_{0,\rm nr})\quad\hbox{in }C_tH^{s-1}\cap C_t^1H^{s-2}.}
\ee
Writing $W_\delta\coloneqq Y_\delta-Y_{0,\delta}$, $W_{\rm nr}\coloneqq Y_{\rm nr}-Y_{0,\rm nr}$ and $b_j\coloneqq\|Y_{0,\delta}-Y_{0,\rm nr}\|_{W^{j,\infty}}$, one also has the bound
\be
\label{nr:eq:lower-convergence}
 \|U_\delta-U_{\rm nr}\|_{C_tH^{s-1}}+\|W_\delta-W_{\rm nr}\|_{C_tH^{s-2}}
 \le C\bigl(\|U_{0,\delta}-U_{0,\rm nr}\|_{H^{s-1}}+b_1^{3-s}b_2^{s-2}+\delta\bigr).
\ee
{At $s=3$ the interpolated term is $b_2$. On the torus it can be replaced by $b_1$. For the normalized variables alone, the bound improves to
\be
\label{nr:eq:normalized-lower-convergence}
 \|U_\delta-U_{\rm nr}\|_{C_tH^{s-1}}
 \le C\bigl(\|U_{0,\delta}-U_{0,\rm nr}\|_{H^{s-1}}+b_1+\delta\bigr).
\ee
}
\item[(iii)] {\emph{Convergence of the physical fluxes.}} The  particle current expressed in physical time $({\particlecurrent}^0,c_{\rm light}{\particlecurrent}^{\rm sp})$, momentum density, spatial stress, and energy density and flux after subtraction of the rest energy converge uniformly to their classical counterparts, defined in~\eqref{nr:eq:classical-fluxes} below. For $d_\delta\le1$, the uniform error is at most $Cd_\delta^{\min\{1,p_\theta\}}$, with $d_\delta=\|U_{0,\delta}-U_{0,\rm nr}\|_{H^{s-1}}+b_1+\delta$. For fixed positive $\varepsilon_{\rm e}$ and exact elastic data, the exponent can be replaced by $\max\{s-5/2,\min\{1,p_\theta\}\}$. Their conservation laws include the initial terms, even at zero density.
\item[(iv)] {\emph{Transport of cavitation and exact minors.}} The flow transports the zero-density set. For $\varepsilon_{\rm e}>0$, it keeps the ranks of $F$ and $G$ constant along each particle trajectory. On the zero-density set the current and stress vanish and $(\del_t+v\cdot\nabla)v=0$. If $\varepsilon_{\rm e}>0$ and the initial data satisfy the exact minor relations, the limiting fields satisfy these relations and the three classical material equations, with the continued fluxes~\eqref{nr:eq:classical-material} below, hold in the sense of distributions with their initial terms.
\end{itemize}
\end{theorem}
The rate uses the following bounds on the initial minors. The uniform bound and $L^\infty$ convergence give $Y_{0,\rm nr}\in W^{3,\infty}$. For $B_\delta=Y_{0,\delta}-Y_{0,\rm nr}$ with sufficiently small supremum norm, centered differences give
\be
 \|B_\delta\|_{W^{2,\infty}}\le C\|B_\delta\|_\infty^{1/3}.
\ee
With $b_0=\|B_\delta\|_\infty\le1$, differences at scale $h$ bound the first and second derivatives by $C(h^{-1}b_0+h^2)$ and $C(h^{-2}b_0+h)$, respectively. Taking $h=b_0^{1/3}$ gives $b_1\le Cb_0^{2/3}$ and $b_2\le Cb_0^{1/3}$. The minor term in~\eqref{nr:eq:lower-convergence} is then at most $Cb_0^{(4-s)/3}$ on $\mathbb R^3$ and $Cb_0^{2/3}$ on $\mathbb T^3$. The latter bound also suffices on both domains for~\eqref{nr:eq:normalized-lower-convergence} and $d_\delta$. For $s<3$, the small $H^2$ and bounded $H^3$ multiplication norms of $B_\delta$ give convergence to zero of the $H^s$ multiplication norm by interpolation, so~\eqref{nr:eq:initial-convergence} gives $Z_{0,\delta}\to Z_{0,\rm nr}$ in $H^s$. At $s=3$ this convergence is imposed.

%-----------------------------------------------------------------------------------------------------------------------------------------

{For fixed initial data, the same estimates give a local existence theorem for the classical equations.}

\begin{corollary}[Local existence for classical elastodynamics with cavitation]
\label{cor:classical-existence-cavitation}
{Fix $\varepsilon_{\rm e}>0$, $5/2<s\le3$ and $\mathcal D=\mathbb T^3$ or $\mathbb R^3$. Use the classical internal energy~\eqref{nr:eq:rest-internal-energy} with the constitutive hypotheses of Theorem~\ref{nr:thm:vacuum-limit}. Prescribe $Q_0\ge0$, $Q_0,v_0\in H^s(\mathcal D)$ and $Y_0\in W^{3,\infty}(\mathcal D)$, with the strict initial density and form bounds of that theorem. Require the two families of spatial forms in $Y_0$ to be closed in the sense of distributions. Put $Z_0=Q_0Y_0$. Then one has the following conclusions.
\begin{itemize}
\item[(i)] The classical normalized equations admit a unique local strong solution with prescribed continuation at zero density in the class~\eqref{nr:eq:regularity}. The existence time and bounds require no positive lower bound on $Q_0$.
\item[(ii)] The particle flow transports $\{Q_0=0\}$ and preserves the ranks of $F$ and $G$ along each trajectory. The current and stress vanish on the transported zero-density set, where $(\del_t+v\cdot\nabla)v=0$. The classical conservation laws hold in the sense of distributions with their initial terms.
\item[(iii)] If the initial data satisfy the exact minor relations, these relations persist and the three classical material equations with fluxes~\eqref{nr:eq:classical-material} hold in the sense of distributions with the prescribed initial data.
\end{itemize}}
\end{corollary}

\begin{proof}
{
We take $(Q_{0,\delta},v_{0,\delta},Y_{0,\delta})=(Q_0,v_0,Y_0)$ for all sufficiently small $\delta>0$. Thus both spatial form matrices and the observed density $r_0=Q_0^{p_\theta}$ are fixed. The two spatial closure constraints are independent of $\delta$. To reconstruct the initial spacetime forms, we extend the spatial forms by requiring them to annihilate $\partial_t+v_0^i\partial_i$ in physical time. This determines their time components also at $r_0=0$. In the inertial coordinates $x^0=c_{\rm light}t$, the initial current and proper density are
\be
\particlecurrent_{0,\delta}=r_0(1,\sqrt\delta\,v_0),\qquad
 n_{0,\delta}=r_0/\Gamma_\delta(v_0).
\ee
If the exact minor relations are prescribed, $D_x f_0=\Fsp_0$, $\Gsp_0=\operatorname{cof}\Fsp_0$ and $\det\Fsp_0=r_0$. The quadratic minors of the reconstructed $F$ then agree with $G$: they have the prescribed spatial components and annihilate the same velocity direction. The cubic current has time component $\det\Fsp_0$ and equals $\particlecurrent_{0,\delta}$. This calculation uses the polynomial minors and remains valid when the determinant vanishes.

Since $v_0\in L^\infty$, we choose $\delta_0$ such that $\delta_0(v_*^{\max})^2<1$ and~\eqref{nr:eq:inertia-margin} holds. The initial rest density variable is $q_{0,\delta}=\Gamma_\delta(v_0)^{-\theta}Q_0\le Q_0$. The two rest-frame transformations on the spatial form vectors, $P_\delta^{-1}$ and $\Gamma_\delta^{-1}P_\delta$, have operator norm at most one: $P_\delta$ has eigenvalue $\Gamma_\delta$ in the velocity direction and eigenvalue one in the orthogonal directions. The strict density and form bounds thus hold uniformly for these data. At $Q_0=0$ the stress and current vanish for each $\delta>0$, and the continued velocity equation uses the prescribed $v_0$.

Multiplication by $Y_0\in W^{3,\infty}$ is bounded on $H^s$ for $0\le s\le3$, by the product rule at integer orders and interpolation. Thus $Z_{0,\delta}=Q_0Y_0$ is fixed in $H^s$, including in $H^3$ when $s=3$. All initial convergences in Theorem~\ref{nr:thm:vacuum-limit} hold with zero difference. We apply its conclusions (i), (iii) and (iv) to obtain the classical normalized solution, its conservation laws with their initial terms, and, for exact initial minors, the classical material equations.
}
\end{proof}

%-----------------------------------------------------------------------------------------------------------------------------------------

\subsection{Uniform equations and estimates at zero density}
\label{nr:subsec:uniform-estimates}

{We now give a proof of Theorem~\ref{nr:thm:vacuum-limit}.} {\bf 1.} \emph{The powers of the light-speed parameter.} We first derive the matrices for positive $\delta$ and then extend them to zero. Dividing the physical constitutive law by $c_{\rm light}^2$ gives $n[m+\delta\Phi(q)+\delta\Psi(q,z)]$ in coordinates with normalized velocity $\sqrt\delta\,v$. The change from $x^0$ to $t$ has $\del_{x^0}=\sqrt\delta\,\del_t$. The normalized matrix representing the second variation of the energy in the material rest frame of this law is $\delta\mathsf H$, where $\mathsf H$ is the matrix in~\eqref{uf:eq:rest-metric}. Its rest inertia matrix is $mI+\delta\mathsf N$, where $\mathsf N\coloneqq\mathsf M_R-mI$ and $\mathsf M_R$ is defined in~\eqref{uf:eq:rest-inertia}. We keep the maps $D_k$, $K_k=\mathsf H D_k$ and $F_{lk}=D_l^T\mathsf H D_k$ of Section~\ref{sec:proofs-local-existence}; thus $F_{vv}$ denotes $F_{lk}$ evaluated at $l=k=v$.

The rest variables and frame matrices are smooth in $\delta$ at zero. With $y=(Q,Z)$, the matrix $\mathsf T_\delta$ converts the weighted form components $Z$ to their rest values $z$. The matrix $\mathsf R_\delta$ includes the density conversion, and $C_\delta$ acts on velocity variations. The resulting time blocks $\mathsf H_\delta,\mathsf M_\delta$ and the spatial blocks $\mathsf B_{i,\delta},\mathsf E_{i,\delta}$ have the same roles as in~\eqref{uf:eq:lab-matrices}:
\bse
\be
\label{nr:eq:rest-transform}
 \begin{gathered}
 \mathsf T_\delta\coloneqq\operatorname{diag}\bigl(\Gamma_\delta^{-\theta}P_\delta^{-1},\Gamma_\delta^{-\theta-1}P_\delta\bigr),\qquad
 q=\Gamma_\delta^{-\theta}Q,\quad z=\mathsf T_\delta Z,\\
 Y_R=\Gamma_\delta^\theta\mathsf T_\delta Y,\qquad
 \mathsf R_\delta\coloneqq\operatorname{diag}(\Gamma_\delta^{-\theta},\mathsf T_\delta),\qquad
 C_\delta\coloneqq\Gamma_\delta P_\delta,\qquad k_{i,\delta}\coloneqq P_\delta^{-1}e_i.
 \end{gathered}
\ee
\be
\label{nr:eq:time-matrices}
 \mathsf H_\delta\coloneqq\Gamma_\delta\mathsf R_\delta^T\mathsf H\mathsf R_\delta,\qquad
 \mathsf M_\delta\coloneqq\Gamma_\delta C_\delta^T\bigl(mI+\delta\mathsf N-\delta^2F_{vv}\bigr)C_\delta.
\ee
\be
\label{nr:eq:spatial-matrices}
 \mathsf B_{i,\delta}\coloneqq\mathsf R_\delta^TK_{k_{i,\delta}}C_\delta,\qquad
 \mathsf E_{i,\delta}\coloneqq-\delta C_\delta^T(F_{v k_{i,\delta}}+F_{k_{i,\delta}v})C_\delta.
\ee
\ese
The diagonal blocks of $\mathsf T_\delta$ act on the three one-form vectors and three two-form vectors, respectively. To verify the powers of $\delta$, we use $V=\sqrt\delta\,v$, $\delta\mathsf H$ and $mI+\delta\mathsf N$ in~\eqref{uf:eq:lab-matrices}. The change from variations of $V$ to variations of $v$ gives the congruence by $\operatorname{diag}(I,\sqrt\delta\,I)$, and $\del_{x^0}+V\cdot\nabla=\sqrt\delta\,\mathrm D_t$, where $\mathrm D_t\coloneqq\del_t+v\cdot\nabla$. The resulting symmetric system has a common factor $\delta^{3/2}$. Dividing by this factor gives~\eqref{nr:eq:time-matrices}--\eqref{nr:eq:spatial-matrices}. In particular the two boost directions in the time block give $\delta^2F_{vv}$, and the relative spatial velocity block carries the factor $\delta$. The displayed matrices extend to $\delta=0$.

\vskip.2cm \noindent {\bf 2.} \emph{The normalized system and its classical limit.} We use the maps $\Pi_i$ defined in~\eqref{uf:eq:Pi}, now acting on a physical velocity variation. The normalized equations and symmetry relations are
\bse
\be
\label{nr:eq:normalized-system}
 \mathrm D_ty+\sum_i\Pi_i\del_i v=0,\qquad
 \mathrm D_tv+\sum_i\mathsf M_\delta^{-1}\bigl(\mathsf B_{i,\delta}^T\del_i y+\mathsf E_{i,\delta}\del_i v\bigr)=0.
\ee
\be
\label{nr:eq:symmetry}
 \mathsf H_\delta\Pi_i=\mathsf B_{i,\delta},\qquad
 \mathsf E_{i,\delta}=\mathsf E_{i,\delta}^T.
\ee
\ese
The shape fields satisfy~\eqref{uf:eq:shape-evolution}, with $v$ in place of $V$, and their equations are independent of $\delta$. At zero parameter the frame matrices are identities, so the equations become
\be
\label{nr:eq:classical-normalized}
 \mathrm D_ty+\sum_i D_{e_i}\del_i v=0,\qquad
 m\mathrm D_tv+\sum_i K_{e_i}^T\del_i y=0.
\ee
All matrices in this formula are evaluated at $(q,z,Y_R)=(Q,Z,Y)$. These are the normalized classical elastic equations, extended to zero density. {In matter $K_{e_i}^T$ is the derivative of the specific momentum flux in direction $e_i$ after dropping the terms proportional to the spatial closure constraints and dividing the stress variation by $n$, as in~\eqref{uf:eq:rest-constraint-terms}--\eqref{uf:eq:rest-system}. The reconstruction satisfies those constraints, and the coefficients extend to zero density.} The limiting conservation laws will also follow directly from the physical fluxes below.

On the compact density range, the Taylor factors of~\eqref{uf:eq:rest-metric} give
$\|\mathsf H-\operatorname{diag}(H_f,\Psi_{zz}(q,0))\|\le C|Y_R|$ for $|Y_R|\le1$. The diagonal lower bounds give $\mathsf H\ge h_*I$ for small $|Y|$, since $|Y_R|\le|Y|$. At $q=0$ the shape-dependent Schur complement is the one in~\eqref{uf:eq:quadratic-rest-schur} with $S=\Psi_{zz}(0,0)$.

To estimate the second time block, we set $h_{{\rm f},\delta}=m+\delta(h-m)$ and define
$c_{f,\delta}^2=\delta\theta^2q^2H_f/h_{{\rm f},\delta}$, the squared fluid wave speed in coordinates with normalized light speed. Since $h_q=\theta qH_f$ and $h(0)=m$, we have $h_{{\rm f},\delta}\ge m$. At zero elastic minors the matrix inside the velocity block~\eqref{nr:eq:time-matrices} equals
$h_{{\rm f},\delta}(I-\delta c_{f,\delta}^2v\otimes v)$. At nonzero minors its difference from that matrix satisfies
\be
 \big\|mI+\delta\mathsf N-\delta^2 F_{vv}
       -h_{{\rm f},\delta}(I-\delta c_{f,\delta}^2v\otimes v)\big\|
 \le C\delta q^2|Y_R|^2(1+\delta|v|^2).
\ee
{This is the pair of rest acoustic estimates from Step~1 of Theorem~\ref{uf:thm:main}, using $F_{vv}$ quadratic in the direction $v$. The speed $c_{f,\delta}^2$ increases with $\delta$ at fixed $q$, since $\delta/[m+\delta(h-m)]$ does. Thus the uniform gap condition}
\be
\label{nr:eq:inertia-margin}
 \delta_0 (v_*^{\max})^2\sup_{0\le q\le q_*}c_{f,\delta_0}^2\le1-\zeta,
 \qquad \zeta>0,
\ee
{gives a lower bound $m\zeta/2$ for the matrix inside the velocity block once $|Y|$ is small enough. The condition holds for small $\delta_0$ and can also be used on a larger interval when the ranges satisfy it. The frame matrices and their inverses are uniformly bounded because $\delta_0(v_*^{\max})^2<1$. With $\Gamma_\delta\le\Gamma_*$, the singular values of $\mathsf R_\delta,C_\delta$ give}
\be
 \mathsf H_\delta\ge h_*\Gamma_*^{-2\theta-1}I,\qquad
 \mathsf M_\delta\ge (m\zeta/2)I .
\ee
{The compact coefficient ranges also give uniform upper bounds. All these estimates are independent of the elastic parameter and a positive density minimum and give the classical positive energy at $\delta=0$.}

The identities underlying the regularity of the principal matrices are unchanged by the parameter. We use $P_2=Z\otimes Y$ and $P_{2R}=\Gamma_\delta^\theta(\mathsf T_\delta\otimes\mathsf T_\delta)P_2$. The expression for $K_k$ in~\eqref{uf:eq:smooth-H-D} is $C^6$ in $(q,z,P_{2R})$, while $F_{lk}$ and $\mathsf N$ are $C^6$ in $(q,z)$. Thus the principal matrices of the equations solved for time derivatives are $C^6$ functions of $(\delta,Q,Z,P_2,v)$ on a neighborhood of the coefficient ranges. Their derivatives in $\delta$ have the same regularity bounds as the matrices themselves through the orders used below. At zero density, $Z=P_2=0$ and every relative spatial matrix vanishes; the equation there is $\mathrm D_tv=0$ for every $0\le\delta\le\delta_0$.

We put $r_s\coloneqq s-1$ and $W\coloneqq Y-Y_0$. We write $C_K$ for bounds depending on the fixed ranges and an upper bound $K$ on $\|U\|_{H^s}+\|W\|_{H^{r_s}}+\|Y_0\|_{W^{3,\infty}}$. In dimension three, $H^{r_s}$ is an algebra and $H^s$ embeds in $W^{1,\infty}$. The derivative identity~\eqref{uf:eq:weighted-product} yields $\|P_2\|_{H^s}\le C_K\|U\|_{H^s}$: each differentiated factor $\nabla Q,\nabla Z$ is in $H^{r_s}$, and $Y$ is the sum of an $H^{r_s}$ field and the bounded multiplier $Y_0$. In particular all principal matrices of the solved system are $H^s$ relative to their constant values at $U=P_2=0$, and are Lipschitz. These properties also hold for the symmetric spatial flux matrices $\mathsf B_{i,\delta},\mathsf E_{i,\delta}$ and the velocity time matrix $\mathsf M_\delta$. The upper time matrix $\mathsf H_\delta$ may depend on $Y$ without a density factor. Its material derivative is bounded using the equations for $(Q,Z,v,Y)$.

We give the commutator bounds used below. We write $\Lambda=(1-\Delta)^{1/2}$ and denote by $\|\widehat f\|_1$ the $L^1$ norm of the Fourier transform (the $\ell^1$ norm of Fourier coefficients on the torus). For a coefficient $A$ and a field $d$, at any order $a\ge0$, the multiplier of $[\Lambda^a,A]\del_i d$ is bounded, with $k$ the coefficient frequency and $\zeta$ the field frequency, by
\be
 \big|\langle k+\zeta\rangle^a-\langle\zeta\rangle^a\big|\,|\zeta_i|
 \le C_a
 \begin{cases}
 |k|\langle\zeta\rangle^a,& |k|\le\langle\zeta\rangle/2,\\
 \langle k\rangle^a|\zeta_i|,& |k|>\langle\zeta\rangle/2.
 \end{cases}
\ee
Young's inequality bounds the first contribution by $C\|\widehat{\nabla A}\|_1\|d\|_{H^a}$. For the second contribution, we can use either $C\|A\|_{H^a}\|\widehat{\nabla d}\|_1$ or $C\|\nabla A\|_{H^a}\|\widehat d\|_1$, since $|\zeta|\le2|k|$ there. Constant parts of $A$ are omitted. When $0\le a\le1$, the second multiplier bound is at most $C_a |k|\langle\zeta\rangle^a$, giving the first convolution bound for both contributions. The Fourier $L^1$ norm is controlled by $H^{r_s}$ by Cauchy--Schwarz.

{
We use~\eqref{uf:eq:approximation} for~\eqref{nr:eq:normalized-system}. Its flow preserves $Q\ge0$ and $Z=QY$. At fixed $\ell$, the argument of Step~4 for~\eqref{uf:eq:label-ODE} applies in $H^6\times H^s$: interpolation gives bounded composition on $H^s$ and $H^{s-1}$, and~\eqref{uf:eq:inverse-flow-difference} gives Lipschitz differences from $H^s$ to $H^{s-1}$ by integration along the segment of inverse maps. The reconstructed fields have the same bounds. Differentiation in $\mathcal F_\ell$ costs one further derivative, with multiplication by $H^{r_s}$ bounded on $H^{s-2}$. The fixed convolution restores the required regularity and bounds the $H^6$ flow as long as the $H^s$ velocity bound and the strict margins hold.
}

We apply $\Lambda^s$ to the equations solved for time derivatives, and put $a_{\ell}\coloneqq\Lambda^s y_{\ell}$, $b_{\ell}\coloneqq\Lambda^s w_{\ell}$. In the commutator estimate, $\nabla A$ and the first derivatives of the unconvolved fields are bounded in $H^{r_s}$. We use the energy
\be
\label{nr:eq:uniform-energy}
 \begin{gathered}
 \mathcal E_{s,\delta}\coloneqq\int_{\mathcal D}
 \big(a_{\ell}^T\mathsf H_\delta a_{\ell}+b_{\ell}^T\mathsf M_\delta b_{\ell}\big)\,dx,\\
 c_0\|(y_{\ell},w_{\ell})\|_{H^s}^2\le\mathcal E_{s,\delta}\le c_1\|(y_{\ell},w_{\ell})\|_{H^s}^2,\qquad
 \frac{d}{dt}\mathcal E_{s,\delta}\le C_K\mathcal E_{s,\delta}.
 \end{gathered}
\ee
We pair the spatial terms involving $a_{\ell},\smoothoperator_{\ell} b_{\ell}$ by~\eqref{nr:eq:symmetry}, as in~\eqref{uf:eq:paired-fluxes}. This differentiates $\mathsf B_{i,\delta},\mathsf E_{i,\delta}$. The convolution commutator with $\mathsf M_\delta$ uses the $H^1$--$H^{-1}$ pairing in~\eqref{uf:eq:metric-commutator}. Transport contributes the material derivatives of the weights. In the approximate velocity equation, both terms of~\eqref{uf:eq:approximate-material-derivative} are bounded in $H^{r_s}$ by the algebra estimate, so these material derivatives are uniformly bounded in $L^\infty$. The same argument applies to $\mathsf H_\delta$ since $\mathrm D_{\ell} Y_{\ell}$ is a product of $\nabla v_{\ell}$ and $Y_{\ell}$. For the equation of the increment $W=Y_{\ell}-Y_0$, the material derivative is $L(Dv_{\ell})(Y_0+W)-v_{\ell}\cdot\nabla Y_0$, where $L(Dv_{\ell})$ is the linear action in~\eqref{uf:eq:shape-evolution}. This source is bounded in $H^{r_s}$ by the algebra estimate and multiplication by $Y_0,\nabla Y_0$. The commutator at order $r_s$ uses $\|\nabla v_{\ell}\|_{H^{r_s}}$ and $\|\widehat W\|_1$. The energy estimates and the strict initial margins give a common interval by the first-exit argument following~\eqref{uf:eq:bootstrap}.

{
For the defect in~\eqref{uf:eq:approximation-error}, the same kernel estimates as in Step~6 give $\|(I-\smoothoperator_\ell^2)\mathcal F_\ell\|_2\le C_K\ell^{-(s-1)}$ from the $H^{s-1}$ bound. The even-kernel transport commutator costs $\Ocal(\ell^{-2})$ in $L^2$: the two velocity and gradient increments use $\|Dv_\ell\|_\infty,\ \|D^2w_\ell\|_2$ and the second moment. At order one the first term costs $\Ocal(\ell^{-(s-2)})$, and the transport commutator is $\Ocal(\ell^{-1})$ by the first moment and the algebra bound in $H^{s-1}$. The lower difference estimate~\eqref{uf:eq:product-difference}--\eqref{uf:eq:cauchy} applies with $H^{s-1}$ multiplying $H^1$ and $\nabla W\in L^3$ in the form transport equation.
}

We set $h_{\ell{}M}\coloneqq {\ell}^{-(s-1)}+M^{-(s-1)}$. On the cone $Q\ge0$, $|Z|\le C_KQ$, the solved coefficients are Lipschitz functions of $U$, as in~\eqref{uf:eq:product-difference}. The order-zero difference energy, with the common initial shape and the mollified initial velocities, gives $\|U_{\ell}-U_M\|_{C_tL^2}\le C_{K,T}h_{\ell{}M}$. Interpolation at order one and the shape transport estimate then give $\|Y_{\ell}-Y_M\|_{C_tL^2}\le C_{K,T}h_{\ell{}M}^{1-1/s}$. Interpolating with the uniform $H^s$ and $H^{s-1}$ bounds yields, for every $1<s'<s$,
\be
 \|U_{\ell}-U_M\|_{C_tH^{s'}}+\|Y_{\ell}-Y_M\|_{C_tH^{s'-1}}
 \le C_{K,T}h_{\ell{}M}^{(s-s')/s}.
\ee
The estimate includes the zero-density set and uses no additional derivatives of $Q$ or $v$. We can choose $s'>5/2$ to pass to the equations and the flow, so the coefficients and the time matrices converge uniformly.

At the highest order the limit has weakly continuous representatives. Spatial convolution of the differentiated equations justifies the energy identity between arbitrary times by the Lipschitz coefficient commutator as in the proof after~\eqref{uf:eq:cauchy}. The weights are continuous in $L^\infty$ and have bounded derivatives along the flow. The resulting continuity of the energy and weak continuity imply strong continuity of $U$ in $H^s$. Applying the same convolution argument to $\Lambda^{r_s}W$ gives strong continuity in $H^{r_s}$. The equations then give the time derivatives in~\eqref{nr:eq:regularity}, including at the initial time. The same lower difference estimate with zero error proves uniqueness. This proves part~(i).

%-----------------------------------------------------------------------------------------------------------------------------------------

\subsection{Convergence at the initial Sobolev regularity}
\label{nr:subsec:top-convergence}

\vskip.2cm \noindent {\bf 3.} \emph{An estimate at order $s-1$.} We denote by $\mathsf A_{i,\delta}$ the solved principal matrices, including advection. At fixed $(U,P_2)$ the parameter derivative is smooth, vanishes at $U=P_2=0$, and has bounded $H^j$ norms for $j=r_s,s$ by composition. In particular,
\be
\label{nr:eq:coefficient-parameter}
 \|\mathsf A_{i,\delta}(U,P_2)-\mathsf A_{i,0}(U,P_2)\|_{H^j}\le C_K\delta,\qquad j=r_s,s.
\ee
{On the torus, the lower difference estimate in the proof of Proposition~\ref{nr:prop:periodic-sobolev-data} below bounds the errors in $U$ in $H^{s-1}$ and in $Y$ in $H^{s-2}$ by the initial errors at these orders and by $C\delta$. The initial difference of $Y$ in $H^{s-2}$ is bounded by $Cb_1$, so the same bound holds for the increments $W=Y-Y_0$.

On $\mathbb R^3$, the lower estimate in the proof of Proposition~\ref{nr:prop:two-bounded-derivatives} below bounds the differences of $U$ and the inverse-flow displacements $h_j=A_j-\operatorname{Id}$ in $C_tH^{s-1}$ by $Cd_\delta$, where $d_\delta=\|U_{0,\delta}-U_{0,\rm nr}\|_{H^{s-1}}+b_1+\delta$. For a fixed background $Y_{0,\rm nr}$, the difference in the reconstructed fields $R(DA_j)(Y_{0,\rm nr}\circ A_j)$ is bounded in $H^{s-2}$ by $C\|h_\delta-h_{\rm nr}\|_{H^{s-1}}$. Indeed, the matrix difference has this bound by the polynomial formula for $R$, while the chain rule gives the same bound in $H^1$ for the difference of the composed backgrounds.

We write $B=Y_{0,\delta}-Y_{0,\rm nr}$. The remaining increment is $(R(DA_\delta)-I)(B\circ A_\delta)+B\circ A_\delta-B$. The first product has $H^{s-2}$ norm bounded by $Cb_1$ since $R(DA_\delta)-I$ is bounded in $H^{s-1}$ and $B\circ A_\delta$ in $W^{1,\infty}$. The last composition difference has $L^2$ norm at most $Cb_1$ and $H^1$ norm at most $Cb_2$ by the chain rule and the uniform $H^{s-1}$ bound for $h_\delta$. Interpolation at order $s-2$ bounds its norm by $Cb_1^{3-s}b_2^{s-2}$. These estimates prove~\eqref{nr:eq:lower-convergence} and~\eqref{nr:eq:normalized-lower-convergence}.}

\vskip.2cm \noindent {\bf 4.} \emph{Convergence at the highest initial regularity.}  We use the normalized convergence arguments in the proofs of Propositions~\ref{nr:prop:periodic-sobolev-data} and~\ref{nr:prop:two-bounded-derivatives} below. Those arguments use the uniform time-matrix estimates of Steps~1 and~2 and apply before the physical reconstruction of Steps~5 and~6. The initial convergence in~\eqref{nr:eq:initial-convergence} gives $U_{0,\delta}\to U_{0,\rm nr}$ in $H^s$ and $b_2\to0$, as proved after the theorem. On the torus it also gives $Y_{0,\delta}\to Y_{0,\rm nr}$ in $H^{s-1}$. The periodic argument gives~\eqref{nr:eq:top-convergence} there.

On $\mathbb R^3$, the argument with two bounded derivatives gives $U_\delta\to U_{\rm nr}$ and $h_\delta\to h_{\rm nr}$ in $C_tH^s\cap C_t^1H^{s-1}$. It gives convergence of $Y_j-Y_{0,j}\circ A_j$ in $C_tH^{s-1}$, also for $s=3$ because $b_2\to0$ implies the required local convergence of $D^2(Y_{0,\delta}-Y_{0,\rm nr})$. Here $A_j=\operatorname{Id}+h_j$. To recover the increments relative to the fixed initial background, we still need convergence of $Y_{0,j}\circ A_j-Y_{0,j}$ in $C_tH^{s-1}$. For the fixed background $Y_{0,\rm nr}\in W^{3,\infty}$, the chain rule up to order two bounds the difference of these compositions by $C\|h_\delta-h_{\rm nr}\|_{C_tH^2}$. For the varying difference $B_\delta=Y_{0,\delta}-Y_{0,\rm nr}$, the same chain rule bounds $B_\delta\circ A_\delta-B_\delta$ and its derivatives of order at most two in $L^2$ by $Cb_2$, except for $D^2B_\delta\circ A_\delta-D^2B_\delta$. This term is bounded pointwise both by $2b_2$ and by $C|h_\delta|$, using the uniform bound on $D^3 B_\delta$. Its $C_tL^2$ norm tends to zero by splitting into a fixed ball and its complement, since $h_\delta\to h_{\rm nr}$ in $C_tL^2$. This gives the required convergence in $C_tH^2$ for the background increments. Finally use
\be
 \del_t W=-v\cdot\nabla W+L(Dv)(Y_0+W)-v\cdot\nabla Y_0.
\ee
The convergence of $v,W$ and the $W^{2,\infty}$ convergence of the backgrounds imply convergence of $\del_tW$ in $C_tH^{s-2}$, by the same multiplication estimates as in Step~3.

%-----------------------------------------------------------------------------------------------------------------------------------------

\subsection{Classical fluxes and constitutive laws}
\label{nr:subsec:physical-fluxes}

\paragraph{The physical fluxes after rest-energy subtraction.}
The remaining step is to recover the physical equations, including the fluxes through the initial slice.
\vskip.2cm \noindent {\bf 5.} \emph{Convergence of current and stress.}  The flux formulas extend continuously to zero density. Write $z=(z_F,z_G)$, with $z_F^A$ and $z_{G,A}$ the spatial vectors of the weighted one- and two-forms, respectively. In the rest frame we write $I_{\rm int}\coloneqq\Phi(q)+\Psi(q,z)$ and define the symmetric specific spatial stress
\be
\label{nr:eq:specific-stress}
 \begin{aligned}
 \Sigma(q,z)&\coloneqq\bigl(\theta q(\Phi'+\Psi_q)+\theta z\cdot\Psi_z+
 \sum_Az_{G,A}\cdot\Psi_{z_{G,A}}\bigr)I\\
 &\quad+\sum_A\bigl(z_F^A\otimes\Psi_{z_F^A}-\Psi_{z_{G,A}}\otimes z_{G,A}\bigr).
 \end{aligned}
\ee
Rotation invariance of $\Psi$ gives symmetry. The rest stress is $n\Sigma$ and $\mathsf N=I_{\rm int}I+\Sigma$. With $q,z$ as in~\eqref{nr:eq:rest-transform}, the boosted energy is $T^{00}=r\Gamma_\delta(m/\delta+I_{\rm int}+\delta v\cdot\Sigma v)$. The physical-time fluxes are
\bse
\be
\label{nr:eq:finite-fluxes}
 \begin{gathered}
 {\particlecurrent}^0=r,\qquad c_{\rm light}{\particlecurrent}^{\rm sp}=rv,\\
 \mathcal P_\delta\coloneqq T^{0{\rm sp}}/c_{\rm light}
 =r\bigl(m\Gamma_\delta v+\delta(\Gamma_\delta I_{\rm int}v+P_\delta\Sigma v)\bigr),\\
 \mathcal S_\delta\coloneqq T^{{\rm sp},{\rm sp}}
 =r\bigl(\Gamma_\delta(m+\delta I_{\rm int})v\otimes v+
 \Gamma_\delta^{-1}P_\delta\Sigma P_\delta\bigr).
 \end{gathered}
\ee
\be
\label{nr:eq:renormalized-energy}
 \begin{gathered}
 e_\delta\coloneqq T^{00}-(m/\delta){\particlecurrent}^0
 =r\bigl(m\Gamma_\delta^2|v|^2/(\Gamma_\delta+1)
 +\Gamma_\delta(I_{\rm int}+\delta v\cdot\Sigma v)\bigr),\\
 f_\delta\coloneqq c_{\rm light}T^{0{\rm sp}}-(m/\delta)rv
 =r\bigl((m\Gamma_\delta^2|v|^2/(\Gamma_\delta+1)+\Gamma_\delta I_{\rm int})v+P_\delta\Sigma v\bigr).
 \end{gathered}
\ee
\ese
The identity $(\Gamma_\delta-1)/\delta=\Gamma_\delta^2|v|^2/(\Gamma_\delta+1)$ makes every flux a product of $Q^{p_\theta}$ and a smooth coefficient through $\delta=0$.

By~\eqref{nr:eq:normalized-lower-convergence} and $H^{s-1}\subset L^\infty$, the smooth factors differ from their limits by at most $Cd_\delta$. The common factor $r=Q^{p_\theta}$ obeys the same estimate with exponent $\min\{1,p_\theta\}$ by the modulus of continuity of the power function on the bounded nonnegative range.

The limiting fluxes, evaluated at the classical fields, are
\bse
\be
\label{nr:eq:classical-fluxes}
 \begin{gathered}
 r=Q^{p_\theta},\qquad \mathcal P_{\rm nr}=mrv,\qquad
 \mathcal S_{\rm nr}=mrv\otimes v+r\Sigma(Q,Z),\\
 e_{\rm nr}=r\bigl(m|v|^2/2+\Phi(Q)+\Psi(Q,Z)\bigr),\qquad
 f_{\rm nr}=e_{\rm nr}v+r\Sigma(Q,Z)v.
 \end{gathered}
\ee
\be
\label{nr:eq:classical-conservation}
 \del_t r+\operatorname{div}(rv)=0,\qquad
 \del_t\mathcal P_{\rm nr}+\operatorname{div}\mathcal S_{\rm nr}=0,\qquad
 \del_t e_{\rm nr}+\operatorname{div}f_{\rm nr}=0.
\ee
\ese
The flow gives $r\circ X=r_0/\det DX$ and hence continuity. Where $Q>0$, $Y=Z/Q$ is $C^1$ and the normalized equations give the rest-space projection of stress balance. The remaining projection follows from the independent constitutive identity~\eqref{uf:eq:independent-stress-projection} and current conservation. {Using the continuity equation in the product identity~\eqref{uf:eq:stress-decomposition} gives a mass-stress divergence proportional to $r$. This divergence vanishes at $Q=0$, where the internal stress and its first derivatives also vanish.} Time rescaling and subtraction of $(m/\delta)$ times the conserved current give the physical-time laws. Uniform convergence passes the bulk and initial terms to~\eqref{nr:eq:classical-conservation}.

For fixed positive $\varepsilon_{\rm e}$, the uniform initial convergence preserves the spatial closure constraints. If the initial data also lie on the set of exact spacetime minors, their polynomial minor relations persist. On the torus the affine periods converge by taking spatial means. To estimate the transported minors, put $\alpha=s-5/2$. At frequency cutoff $R\ge1$, the $L^\infty$ bound for $D(v_\delta-v_{\rm nr})$ is $Cd_\delta R^{1-\alpha}+CR^{-\alpha}$, by the difference bound in $H^{s-1}$ and the uniform $H^s$ bound. Taking $R=d_\delta^{-1}$ gives $Cd_\delta^\alpha$. The position difference is $\Ocal(d_\delta)$ by the velocity estimate. The differentiated flow equations use the $\alpha$-H\"older bound for $Dv_{\rm nr}$ and give an $\Ocal(d_\delta^\alpha)$ difference in $C_tW^{1,\infty}$, also for the inverse flows. Reconstruction with the uniformly Lipschitz initial minors gives the same rate for the minors and hence for their determinant, including at rank loss. Together with the power estimate this proves the improved flux exponent. At $Q=0$, the fluxes vanish and~\eqref{nr:eq:classical-normalized} gives $\mathrm D_tv=0$.

\paragraph{{6. The material equations at rank loss.}}
{We now restrict to initial data sets satisfying the exact minor relations and fix $\varepsilon_{\rm e}>0$. For these data, the same limit gives} the three equations obtained by varying the configuration. We subtract $(m/\delta)\det \Fsp $ from the Lagrangian density. This term contributes only $(m/\delta)\cof \Fsp $ to the spatial material flux, whose divergence is zero by the Piola identity; it leaves the equations unchanged, including at rank loss. The resulting specific Lagrangian, with $\Fsp ,q,v$ regarded as independent when differentiating, is
\bse
\be
\label{nr:eq:renormalized-Lagrangian}
 \begin{aligned}
 \ell_\delta^{\rm ren}(\Fsp ,q,v)&\coloneqq-\frac{m|v|^2}{1+\sqrt{1-\delta|v|^2}}\\
 &\quad+\Gamma_\delta^{-1}\bigl(\Phi(q)+\Psi(q,\varepsilon_{\rm e}q(\Fsp P_\delta^{-1},\Gamma_\delta^{-1}(\cof \Fsp )P_\delta))\bigr).
 \end{aligned}
\ee
\be
\label{nr:eq:renormalized-material-fluxes}
 \begin{gathered}
 B_\delta\coloneqq\ell_{\delta,v}^{\rm ren}-\theta q\delta\Gamma_\delta^2v\ell_{\delta,q}^{\rm ren},\qquad
 \mathsf p_\delta^0\coloneqq-(\cof \Fsp )B_\delta,\\
 \mathsf k_\delta\coloneqq(\ell_\delta^{\rm ren}+\theta q\ell_{\delta,q}^{\rm ren})\cof \Fsp +r\ell_{\delta,\Fsp }^{\rm ren},\qquad
 \mathsf p_\delta^{\rm sp}\coloneqq\mathsf k_\delta+\mathsf p_\delta^0\otimes v.
 \end{gathered}
\ee
\ese
The derivation of the material fluxes in~\eqref{uf:eq:material-momenta} applies with $\del_vq=-\theta q\delta\Gamma_\delta^2v$. In matter the stress balance laws and the Noether identity give $\del_t\mathsf p_\delta^0+\operatorname{div}\mathsf p_\delta^{\rm sp}=0$ in physical time, since the configuration gradient has rank three. The expressions above are finite at $\delta=0$, where
\bse
\be
\label{nr:eq:classical-material}
 \begin{gathered}
 \ell_{\rm nr}\coloneqq-m|v|^2/2+\Phi(Q)+\Psi(Q,\varepsilon_{\rm e}Q(\Fsp ,\cof \Fsp )),\\
 \mathsf p_{\rm nr}^0=m(\cof \Fsp )v,\qquad
 \mathsf p_{\rm nr}^{\rm sp}=(\ell_{\rm nr}+\theta Q\ell_{{\rm nr},Q})\cof \Fsp 
 +r\ell_{{\rm nr},\Fsp }+\mathsf p_{\rm nr}^0\otimes v.
 \end{gathered}
\ee
For the distributional equations at order $s$, the reconstructed spatial minors are locally $H^{r_s}$ in space and their time derivatives are locally $H^{s-2}$. As in the material calculation of Section~\ref{uf:subsec:dependence-reconstruction}, the internal terms of the fluxes are smooth multiples of $q^2$ or $rq^2$, with $rq^2=Q^{p_\theta+2}\Gamma_\delta^{-2\theta}$. Taylor's formula gives a factor $q^2$ in the internal Lagrangian, in its derivatives in $v$ and $\Fsp $, and in $q$ times its derivative in $q$. The remaining coefficients and their first weak derivatives are locally integrable on the bounded ranges. Since $p_\theta+2>1$, both powers have zero first derivative at $Q=0$. Thus all material fluxes have locally integrable weak first derivatives, and the derivatives of the internal terms vanish almost everywhere on the zero set. The kinetic terms in $\mathsf p_\delta^0$ and $\mathsf k_\delta$ are $m(\cof \Fsp )\Gamma_\delta v$ and $(m/\delta)(\Gamma_\delta^{-1}-1)\cof \Fsp $, respectively. For the Piola identity, smooth the configuration locally. Its gradients converge in $H^{s-1}_{\rm loc}$, and so do their cofactors because $H^{s-1}$ is an algebra in three dimensions. Passing to the limit in the smooth Piola identity gives $\operatorname{div}(\cof \Fsp )=0$ in the sense of distributions, including where $\det \Fsp =0$. Cofactor transport follows from the flow formulas. The divergence of the kinetic flux is $m(\cof \Fsp )\mathrm D_t(\Gamma_\delta v)$. Indeed, cofactor transport gives the term $m\Gamma_\delta(\cof \Fsp )(Dv)^Tv$, which cancels the derivative of $(m/\delta)(\Gamma_\delta^{-1}-1)$ in the spatial flux. The cancellation uses
\be
 \nabla\!\left(\frac{\Gamma_\delta^{-1}-1}{\delta}\right)
 =-\Gamma_\delta(Dv)^Tv,
\ee
\ese
whose two sides extend continuously to $\delta=0$. The remaining divergence vanishes on $Q=0$, where $\mathrm D_tv=0$. The material equations with these continued fluxes thus hold on the whole domain in the sense of distributions. The flow reconstruction gives $\Fsp ,\cof \Fsp \in C_tL^\infty_{\rm loc}$, while $Q,v\in C_tH^s\subset C_tL^\infty$. The flux formulas are smooth functions of $(r,\Fsp,q,v,\delta)$ on these bounded ranges, including $Q=0$ and $\delta=0$. The preceding flow and field estimates give uniform convergence to the classical material fluxes at rate $\Ocal(d_\delta^{s-5/2})$. The fluxes are continuous in time in $L^\infty_{\rm loc}$, so their values at $t=0$ give the initial terms. Uniform convergence gives convergence of the bulk and initial pairings.

{At zero density, the material momentum in the initial data is obtained by evaluating the same formula:}
\be
 \left.\mathsf p_\delta^0\right|_{t=0,\ Q_0=0}
 =m(\cof \Fsp _0)\Gamma_\delta(v_0)v_0,
 \qquad \Gamma_0(v_0)=1.
\ee
This normal material flux vanishes when $\Fsp _0=D_xf_0$ has rank at most one. At rank two it may be nonzero even though the initial particle current and stress--energy tensor both vanish.
This completes the proof of Theorem~\ref{nr:thm:vacuum-limit}.

\paragraph{{Periodic data.}}
{On the torus the unweighted initial elastic fields can converge in $H^{s-1}$, with convergence of the weighted fields still in $H^s$. We smooth $Q_0$ and $Y_0$ separately and define the weighted field by their smoothed product, which preserves the defining relation at zero density.}

\begin{proposition}[Non-relativistic limit for periodic Sobolev data]
\label{nr:prop:periodic-sobolev-data}
Let $5/2<s\le3$ and fix $\varepsilon_{\rm e}\ge0$. On $\mathbb T^3$, assume the constitutive bounds, physical initial minor relations and uniform strict range margins of Theorem~\ref{nr:thm:vacuum-limit}. Prescribe $Q_{0,\delta}\ge0$, $Z_{0,\delta}=Q_{0,\delta}Y_{0,\delta}$, and suppose that, as $\delta\downarrow0$,
\be
 (Q_{0,\delta},Z_{0,\delta},v_{0,\delta})\longrightarrow U_{0,\rm nr}\quad\hbox{in }H^s,
 \qquad Y_{0,\delta}\longrightarrow Y_{0,\rm nr}\quad\hbox{in }H^{s-1}.
\ee
At $\varepsilon_{\rm e}=0$, set $Y_{0,\delta}=Z_{0,\delta}=0$.
\begin{itemize}
\item[(i)] \emph{Convergence of the normalized fields.} The normalized relativistic and classical solutions exist on the same interval and satisfy
\be
 U_\delta\longrightarrow U_{\rm nr}\quad\hbox{in }C_tH^s\cap C_t^1H^{s-1},\qquad
 Y_\delta\longrightarrow Y_{\rm nr}\quad\hbox{in }C_tH^{s-1}\cap C_t^1H^{s-2}.
\ee
\item[(ii)] \emph{Convergence of the physical fluxes.} The current and renormalized physical fluxes of Theorem~\ref{nr:thm:vacuum-limit} converge uniformly, including at zero density, and their conservation laws hold with the initial terms. For fixed $\varepsilon_{\rm e}>0$, the physical density and current also converge in $C_tH^{s-1}\cap C_t^1H^{s-2}$.  If $1<\gamma\le3$, the physical density, current and renormalized physical fluxes converge in $C_tH^2\cap C_t^1H^1$ for every $\varepsilon_{\rm e}\ge0$, including the fluid case.
\item[(iii)] \emph{Cavitation and the minor relations.} The flow transports the initial zero-density set. For fixed positive elastic parameter, the limiting minors satisfy the polynomial relations, the flow preserves their ranks, and the three classical material equations hold with their initial terms.
\end{itemize}
\end{proposition}

\begin{proof}
We put $\sigma=s-1>3/2$. The periodic construction in Proposition~\ref{uf:cor:fractional-existence} uses only $U_0\in H^s$ and $Y_0\in H^\sigma$. The same estimates apply uniformly for small $\delta$ using the positive time matrices of~\eqref{nr:eq:inertia-margin}. They give a common interval and uniform bounds in $C_tH^s\times C_tH^\sigma$, including $\delta=0$. We prove convergence by smoothing a single classical initial data set and then removing the smoothing.

We write $\smoothoperator_h$ for convolution with a positive even smooth kernel of scale $h\in(0,1]$. For the classical initial data set $Q_h=\smoothoperator_hQ_0$, $v_h=\smoothoperator_hv_0$, $Y_h=\smoothoperator_hY_0$ and $Z_h=Q_hY_h$. The positivity and the relation $Z_h=Q_hY_h$ hold at every $h$. The positive kernel preserves the strict pointwise ranges of $Q_0,v_0,Y_0$. These auxiliary data need not satisfy the polynomial minor relations. The even kernel and the covariance identity for convolution give
\be
 \|(Q_h,Z_h,v_h)-U_{0,\rm nr}\|_{H^\sigma}
 +\|Y_h-Y_{0,\rm nr}\|_{H^{\sigma-1}}=o(h).
\ee
For the $Z$ term, the covariance is the average of $(Q_0(\cdot-y)-Q_h)(Y_0(\cdot-y)-Y_h)$; the first factor is $\Ocal(h)$ in $H^\sigma$ and the second is $o(1)$ there.  Since $Z_0=Q_0Y_0\in H^{\sigma+1}$, the remaining term $\smoothoperator_hZ_0-Z_0$ is $o(h)$ in $H^\sigma$. The algebra property therefore proves the displayed estimate without division by $Q_0$, including on its zero set. The data also converge at the higher orders $H^s\times H^\sigma$. To see this for $Z_h$, we use $Q_0\nabla Y_0=\nabla Z_0-Y_0\nabla Q_0\in H^\sigma$ and
\be
 Q_h\smoothoperator_h\nabla Y_0-\smoothoperator_h(Q_0\nabla Y_0)
 =(Q_h-Q_0)\smoothoperator_h\nabla Y_0+[Q_0,\smoothoperator_h]\nabla Y_0.
\ee
The first term tends to zero in $H^\sigma$ because $\|Q_h-Q_0\|_{H^\sigma}=o(h)$. The commutator is bounded on $H^\sigma$ by $C\|\nabla Q_0\|_{H^\sigma}\|Y_0\|_{H^\sigma}$ by integration by parts in the kernel; it tends to zero first for smooth fields and then by density.

The auxiliary higher norms obey
\be
 \|(Q_h,Z_h,v_h)\|_{H^{s+1}}+\|Y_h\|_{H^s}\le C h^{-1}.
\ee
For the term $Q_h\smoothoperator_h\nabla Y_0$, the preceding decomposition and the tame product estimate give this bound using
\be
 \|[Q_0,\smoothoperator_h]\partial_iY_0\|_{H^{\sigma+1}}
 \le Ch^{-1}\|Q_0\|_{H^{\sigma+1}}\|Y_0\|_{H^\sigma}.
\ee
To verify the commutator estimate, we write $\widehat\kappa$ for the Fourier multiplier of ${\smoothoperator_1}$ and $k,\ell$ for the frequencies of $Q_0,Y_0$. Its multiplier is $\ell_i(\widehat\kappa(h\ell)-\widehat\kappa(h(k+\ell)))$. If $|k|\le|\ell|/2$, the mean-value theorem and decay of $\widehat\kappa$ bound the $H^{\sigma+1}$-weighted multiplier by $Ch^{-1}|k|\langle\ell\rangle^\sigma$. Otherwise bound the two terms separately by $Ch^{-1}(\langle k\rangle^{\sigma+1}+|k|\langle\ell\rangle^\sigma)$.  For these two terms we use, respectively, $|\ell|\,|\widehat\kappa(h\ell)|\le Ch^{-1}$ and $\langle k+\ell\rangle|\widehat\kappa(h(k+\ell))|\le Ch^{-1}$. Fourier convolution and $H^\sigma\subset\mathcal F^{-1}\ell^1$ give the displayed estimate. The identity~\eqref{uf:eq:weighted-product} gives the same $\Ocal(h^{-1})$ bound for $P_{2h}=Z_h\otimes Y_h$ in $H^{s+1}$.  Indeed, the tame product estimate places the highest derivative on $Y_h$ or $(Q_h,Z_h)$, while their lower $H^\sigma$ and $H^s$ norms control $Y_h$ and $D(Q_h,Z_h)$ in $L^\infty$.  The grouped $H^{s+1}$ energy differentiates the spatial flux matrices and uses only $\mathrm D_t\mathsf H_\delta$, bounded through $\mathrm D_tY=L(Dv)Y$. The $H^s$ shape commutator is controlled by $\|\widehat{\nabla v}\|_1\|Y\|_{H^s}+\|\nabla v\|_{H^s}\|\widehat Y\|_1$, whose Fourier $L^1$ factors use the lower bounds. Gronwall's inequality propagates both $\Ocal(h^{-1})$ estimates, including at $Q=0$. In particular, at $s=3$ the algebra $H^2$ gives these estimates without a $W^{2,\infty}$ assumption on the periodic initial minors.

\paragraph{{Difference estimates and removal of smoothing.}}
For two normalized solutions, we write $d_\sigma=\|\Delta U\|_{H^\sigma}+\|\Delta Y\|_{H^{\sigma-1}}$ and $d_s=\|\Delta U\|_{H^s}+\|\Delta Y\|_{H^\sigma}$.  Since $H^\sigma$ multiplies $H^{\sigma-1}$ and $Z\otimes Z/Q$ is Lipschitz on the bounded cone, subtraction of~\eqref{uf:eq:weighted-product} gives
\be
 \|\Delta P_2\|_{H^\sigma}\le C d_\sigma,\qquad
 \|\Delta P_2\|_{H^s}\le C d_s.
\ee
{The lower estimate requires $\Delta Y$ only in $H^{\sigma-1}$ and holds at $Q=0$.} The grouped energy and shape transport equations therefore yield, for parameters $\delta_1,\delta_2$,
\be
 \|d_\sigma\|_{C_t}\le C\bigl(d_\sigma(0)+|\delta_1-\delta_2|\bigr),\qquad
 \|d_s\|_{C_t}\le C\bigl(d_s(0)+(1+M_2)(d_\sigma(0)+|\delta_1-\delta_2|)\bigr),
\ee
where $M_2=\|U_2\|_{C_tH^{s+1}}+\|Y_2\|_{C_tH^s}$. At the higher order the terms $\Delta A\,\nabla U_2$ and $\Delta v\cdot\nabla Y_2$ use $M_2d_\sigma$; all other factors use the uniform lower bounds. The dependence of the solved matrices on $\delta$ is smooth on the fixed ranges.

We compare each solution with the classical auxiliary solution. We put $U_{0,h}=(Q_h,Z_h,v_h)$ and define
\be
 \begin{aligned}
 \omega(h)\coloneqq{}&\|U_{0,h}-U_{0,\rm nr}\|_{H^s}
 +\|Y_h-Y_{0,\rm nr}\|_{H^\sigma}\\
 &+h^{-1}\bigl(\|U_{0,h}-U_{0,\rm nr}\|_{H^\sigma}
 +\|Y_h-Y_{0,\rm nr}\|_{H^{\sigma-1}}\bigr).
 \end{aligned}
\ee
The preceding smoothing estimates give $\omega(h)\to0$. {We use $d_s(0),d_\sigma(0)$ for the distances between the relativistic and classical initial data before smoothing. The difference estimates and the $\Ocal(h^{-1})$ bound for the smoothed solution give}
\be
 \|U_\delta-U_{\rm nr}\|_{C_tH^s}+\|Y_\delta-Y_{\rm nr}\|_{C_tH^\sigma}
 \le C\bigl(d_s(0)+\omega(h)+h^{-1}(d_\sigma(0)+\delta)\bigr).
\ee
First letting $\delta\downarrow0$ for fixed $h$, then $h\downarrow0$, proves both convergences. The solved matrices converge in $C_tH^\sigma$ by the lower difference estimate, and their products with $\nabla U_\delta$ give convergence of $\partial_tU_\delta$ in $C_tH^\sigma$. The shape transport equation similarly gives convergence of $\partial_tY_\delta$ in $C_tH^{\sigma-1}$.  The time-derivative estimates hold on the closed interval, including on the transported zero-density set.

For $\varepsilon_{\rm e}>0$, the initial physical minor relations pass to the limit and are transported by the flow. The spatial gradient and its cofactor are components of $Y/\varepsilon_{\rm e}$, so they converge in $C_tH^\sigma\cap C_t^1H^{\sigma-1}$. Their determinant gives $r=Q^{p_\theta}$ in the same class, even when the power is not differentiable at zero. The current, renormalized fluxes and material equations then follow as in the proof of Theorem~\ref{nr:thm:vacuum-limit}, including the initial terms.  For the improved density and flux convergence when $1<\gamma\le3$, we put $p=2/(\gamma-1)\ge1$. We first show that the map $w\mapsto w^p$ is continuous on the nonnegative cone of $H^2(\mathbb T^3)$. The cases $p=1$ and $p\ge2$ follow from the identity map and the $C^2$ chain rule, respectively.  For smooth nonnegative $w$, integration of the divergence of $|Dw|^2Dw/(w+\epsilon)$ gives
\be
 \int\frac{|Dw|^4}{(w+\epsilon)^2}\,dx
 =\int\frac{2D^2w(Dw,Dw)+\Delta w\,|Dw|^2}{w+\epsilon}\,dx.
\ee
Cauchy's inequality, positive mollification and Fatou's lemma yield
\be
 \left\|\mathbf1_{\{w>0\}}\frac{|Dw|^2}{w}\right\|_2\le C\|D^2w\|_2.
\ee
On $\mathbb R^3$, we integrate first against spatial cutoffs and let their radii tend to infinity at fixed $\epsilon$. The cutoff term tends to zero because $Dw\in L^3$. We then pass to the limit as $\epsilon\downarrow0$.
For $1<p<2$, the same inequality gives directly
\be
 \|\mathbf1_{\{0<w\le h\}}w^{p-2}|Dw|^2\|_2
 \le h^{p-1}\left\|\mathbf1_{\{w>0\}}\frac{|Dw|^2}{w}\right\|_2
 \le C h^{p-1}\|D^2w\|_2.
\ee
Together with $w^{p-1}|D^2w|\le h^{p-1}|D^2w|$ on this set, this controls both terms in the Hessian of $w^p$ near zero. Away from zero the usual chain rule applies. Approximation therefore gives $w^p\in H^2$ and continuity of the power map on the nonnegative cone of $H^2$, on both spatial domains. Since $w$ is bounded and $p>1$, its power is also in $L^2(\mathbb R^3)$. Applying this continuity to the established $C_tH^s$ convergence gives $Q_\delta^p\to Q_{\rm nr}^p$ in $C_tH^2$.

The continuity equation $\partial_t r+\operatorname{div}(rv)=0$, with $r=Q^p$, gives convergence of the time derivatives in $C_tH^1$. The Lorentz factor and the formulas~\eqref{nr:eq:finite-fluxes}--\eqref{nr:eq:renormalized-energy} then give the same convergence for the density, current and renormalized fluxes. The coefficient functions in these formulas depend smoothly on $(Q,Z,v,\delta)$ on {the chosen bounded ranges}, and their time derivatives converge in $C_tH^1$ by the normalized equations. This argument applies equally to the reduced fluid system and does not use the determinant of the elastic gradient.
This completes the proof of Proposition~\ref{nr:prop:periodic-sobolev-data}.
\end{proof}

{\paragraph{{The limit for the nonlinear elastic law.}}
{For the law of Corollary~\ref{cor:nonlinear-second-weight}, we use the square of the regularized density to weight the elastic fields in the non-relativistic limit.}
\begin{corollary}[Non-relativistic limit for the nonlinear elastic law]
\label{cor:nr-nonlinear-second-weight}
Let $\mathcal D=\mathbb T^3$ or $\mathbb R^3$, and fix the $C^6$ constitutive function $\Psi$ and parameters $m,\kappa,\gamma$ of Corollary~\ref{cor:nonlinear-second-weight}. Replace its rest-mass term by ${m \, n}/\delta$, with $\delta=c_{\rm light}^{-2}$. Use the physical velocity $v$. In this corollary $X$ collects the unweighted forms and $W=QX$ the singly weighted forms.
\be
\label{eq:nr-second-weight-variables}
 X=(\Fsp,\Gsp),\qquad Q=r^\theta,\qquad W=QX,\qquad Z=QW,\qquad U=(Q,Z,v).
\ee
Prescribe compatible initial data with the uniform norm bounds and strict range margins of that corollary. Suppose that $X_{0,\delta}$ is uniformly bounded in $H^3$ on the torus or $W^{3,\infty}$ on the whole space, and that
\be
\label{eq:nr-second-weight-initial-convergence}
 (Q_{0,\delta},Z_{0,\delta},v_{0,\delta})\longrightarrow U_{0,\rm nr}
 \quad\hbox{in }H^3,\qquad
 X_{0,\delta}\longrightarrow X_{0,\rm nr}\quad\hbox{in }L^\infty,
 \qquad Z_{0,\delta}=Q_{0,\delta}^2X_{0,\delta}.
\ee
For sufficiently small $\delta_0>0$, the relativistic and classical normalized solutions exist on a common interval and, as $\delta\downarrow0$, one has
\be
\label{eq:nr-second-weight-convergence}
 U_\delta\longrightarrow U_{\rm nr}\quad\hbox{in }C_tH^3\cap C_t^1H^2,
 \qquad W_\delta\longrightarrow W_{\rm nr}\quad\hbox{in }C_tH^2\cap C_t^1H^1.
\ee
The limit solves the classical normalized system for the law~\eqref{eq:nonlinear-second-weight-law}, with the same prescribed continuation at zero density. No positive lower density or rank bound is imposed.
\end{corollary}

\begin{proof}
We check the parameter bounds and the removal of smoothing.
\vskip.2cm \noindent {\bf 1.} \emph{The common interval.} In~\eqref{nr:eq:time-matrices}--\eqref{nr:eq:spatial-matrices}, use the Hessian and kinematic columns of~\eqref{eq:q2-explicit-rest-blocks}, and replace the form blocks of $\mathsf T_\delta$ by $\Gamma_\delta^{-2\theta}P_\delta^{-1}$ and $\Gamma_\delta^{-2\theta-1}P_\delta$. The physical-time scaling is unchanged. Formula~\eqref{eq:q2-explicit-coupling} shows that the solved matrices are $C^4$ in $(\delta,Q,Z,Z\otimes W,v)$, with parameter derivatives of class $C^3$. These bounds suffice for the coefficient compositions and differences below. Their time blocks are uniformly positive for small $\delta_0$ and {the chosen range of the forms}. The grouped estimates in the proof of Corollary~\ref{cor:nonlinear-second-weight} therefore give a common interval and uniform $H^3\times H^2$ bounds, including at $\delta=0$.

\vskip.2cm \noindent {\bf 2.} \emph{Convergence at the endpoint.} The initial products $W_{0,\delta}=Q_{0,\delta}X_{0,\delta}$ are uniformly $H^3$ and converge in $H^2$. On the torus this follows by interpolation of the bounded $H^3$ raw fields with their $L^\infty$ convergence. On $\mathbb R^3$, bounded third derivatives and $L^\infty$ convergence give $W^{2,\infty}$ convergence of $X_{0,\delta}$ by finite-difference interpolation, followed by the product rule.
For the fixed classical data, set $Q_h=\smoothoperator_hQ_0$, $v_h=\smoothoperator_hv_0$, $W_h=\smoothoperator_hW_0$ and $Z_h=Q_hW_h$. For $M_X\coloneqq\|X_0\|_\infty$, positivity of the kernel gives $|W_h|\le M_XQ_h$, which retains the range of the unweighted forms at zero density. The compatible smoothing estimates in the proof of Proposition~\ref{nr:prop:periodic-sobolev-data} apply with $Y=W$: they use only $Q_0,Z_0,v_0\in H^3$, $W_0\in H^2$ and $Z_0=Q_0W_0$. Since $W_0$ is globally $H^2$, their Fourier proof applies also on $\mathbb R^3$. The auxiliary data need not satisfy the original spatial closure constraints.

In those same high- and low-order comparison estimates, put $d_2=\|\Delta U\|_{H^2}+\|\Delta W\|_{H^1}$ and $d_3=\|\Delta U\|_{H^3}+\|\Delta W\|_{H^2}$. The product identity~\eqref{eq:quadratic-second-weight-product} gives the required coefficient-difference bounds; the additional divergence term in the transport of $W$ is linear in $Dv$ and has the same Sobolev estimates. The auxiliary $H^4\times H^3$ norm is $\Ocal(h^{-1})$. Thus, for a function $\omega(h)\to0$ defined by the compatible smoothing errors, one obtains
\be
 \|U_\delta-U_{\rm nr}\|_{C_tH^3}+\|W_\delta-W_{\rm nr}\|_{C_tH^2}
 \le C\bigl(d_3(0)+\omega(h)+h^{-1}(d_2(0)+\delta)\bigr).
\ee
First let $\delta\downarrow0$ at fixed $h$, then $h\downarrow0$. The initial convergences give $d_3(0)\to0$, proving the spatial convergence. The equations and the coefficient-difference estimates give the asserted time derivatives. The auxiliary higher norms disappear when smoothing is removed, so no fourth derivative of the initial fluid variables is required. The flow and product identities identify the classical continued system also at $Q=0$.
\end{proof}
}

\paragraph{Initial minors with two bounded derivatives.}
Measuring the unweighted minors relative to their transported initial values allows us to use initial minors with two bounded derivatives.

\begin{proposition}[Non-relativistic limit with two bounded derivatives of the initial minors]
\label{nr:prop:two-bounded-derivatives}
Let $\mathcal D=\R^3$ and $5/2<s\le3$. Assume the constitutive hypotheses, strict pointwise range margins and physical initial minor relations of Theorem~\ref{nr:thm:vacuum-limit}. Assume that $(Q_{0,\delta},Z_{0,\delta},v_{0,\delta})$ converges in $H^s$ as $\delta\downarrow0$, where $Z_{0,\delta}=Q_{0,\delta}Y_{0,\delta}$, and that $Y_{0,\delta}$ is bounded in $W^{2,\infty}$ and converges in $L^\infty$. The following conclusions hold.
\begin{itemize}
\item[(i)] There is a common time interval on which $U_\delta$ and the inverse particle displacement $h_\delta=X_\delta^{-1}-\operatorname{Id}$ belong to $C_tH^s\cap C_t^1H^{s-1}$ and converge in $C_tH^s\cap C_t^1H^{s-1}$ to their classical counterparts.   The normalized equations and flow reconstruction of Theorem~\ref{nr:thm:vacuum-limit} hold, including at cavitation. The physical fluxes specified there converge uniformly and satisfy the conservation laws with their initial terms. 
\item[(ii)] With $A_\delta\coloneqq X_\delta^{-1}$, the increments $Y_\delta-Y_{0,\delta}\circ A_\delta$ converge to $Y_{\rm nr}-Y_{0,\rm nr}\circ A_{\rm nr}$ in $C_tH^{s-1}$ if $s<3$, and in $C_tH^r$ for every $r<2$ if $s=3$.  At $s=3$ the convergence also holds in $C_tH^2$ when
\be
\label{nr:eq:initial-vacuum-hessians}
 \mathbf 1_{\{Q_{0,\rm nr}=0\}}D^2(Y_{0,\delta}-Y_{0,\rm nr})\longrightarrow0
 \quad\hbox{in }L^2_{\rm loc}(\mathbb R^3).
\ee
This condition is automatic if the initial zero-density set has zero volume. More precisely, put $\mathcal V_0=\{Q_{0,\rm nr}=0\}$ and write $R(DA)$ for the action of the inverse flow on the scaled spatial one-forms and two-forms, so that $Y=R(DA)(Y_0\circ A)$. At $s=3$, convergence of the increments in $C_tH^2$ is equivalent to
\be
\label{nr:eq:initial-vacuum-increment}
 \big\|\mathbf1_{\mathcal V_0}
 [(R(DA_{\rm nr})-I)\circ X_{\rm nr}]\,D^2(Y_{0,\delta}-Y_{0,\rm nr})\big\|_{C_tL^2}
 \longrightarrow0 .
\ee
\item[(iii)] {For each fixed positive elastic parameter, the physical density and current also belong locally to $C_tH^{s-1}\cap C_t^1H^{s-2}$, including at cavitation and in the classical limit. They converge to their classical counterparts in $C_tH^r_{\rm loc}$ for every $0\le r<s-1$. The endpoint $r=s-1$, together with convergence in $C_t^1H^{s-2}_{\rm loc}$, also holds when $s<3$, and when $s=3$ under~\eqref{nr:eq:initial-vacuum-hessians}. For $s=3$ and $1<\gamma\le3$, both endpoint convergences hold without that additional assumption. In each of these endpoint cases, the momentum, spatial stress, renormalized energy and energy flux of~\eqref{nr:eq:finite-fluxes}--\eqref{nr:eq:renormalized-energy} have the same local convergence.  For $1<\gamma\le3$ and every $\varepsilon_{\rm e}\ge0$, the physical density, current and these renormalized fluxes converge globally in $C_tH^2(\mathbb R^3)\cap C_t^1H^1(\mathbb R^3)$, including in the fluid case.}
\item[(iv)] {For every $\gamma>1$, the observed particle density, current, momentum density, spatial stress, and renormalized energy density and flux converge globally in $C_tH^{s-1}\cap C_t^1H^{s-2}$ if and only if the initial observed densities $r_{0,\delta}=Q_{0,\delta}^{p_\theta}$ converge in $H^{s-1}(\mathbb R^3)$. The same equivalence holds in the corresponding local Sobolev spaces.}
\end{itemize}
\end{proposition}

\paragraph{{The inverse particle map.}}
{We estimate derivatives of the transported forms through the displacement of the inverse flow.}
\begin{proof}
{\bf 1.} \emph{Control by the inverse particle map.} We write $r_s=s-1$. The displacement satisfies $(\del_t+v\cdot\nabla)h=-v$ with zero initial data, and
\be
 Y=R(DA)(Y_0\circ A),\qquad A=\operatorname{Id}+h ,
\ee
where the blocks of $R$ are the transpose and the adjugate of $DA$, as in the proof above. We use a common bound for $\|U\|_{H^s}+\|h\|_{H^s}+\|Y_0\|_{W^{2,\infty}}$ and the same strict coefficient margins. The field $R(DA)-I$ belongs to the algebra $H^{r_s}$ with a bound depending on $\|h\|_{H^s}$. The composed field $Y_0\circ A$ and its first derivatives are bounded. Its second derivatives are
\be
 D^2(Y_0\circ A)=(D^2Y_0\circ A)[DA,DA]+(DY_0\circ A)D^2h.
\ee
The first term is bounded and the second belongs to $L^2$. Thus multiplication by $Y_0\circ A$ is bounded on $L^2$ and $H^2$. In the second derivative of the product, the only term requiring more than bounded coefficients is estimated by
\be
 \|(DY_0\circ A)D^2h\,f\|_2
 \le C\|DY_0\|_\infty\|D^2h\|_2\|f\|_{H^2},
\ee
using $H^2(\R^3)\subset L^\infty$. By interpolation it is bounded on $H^{r_s}$. Using the algebra property and interpolating between $L^2$ and $H^{r_s}$, we conclude that $Y=R(DA)(Y_0\circ A)$ multiplies $H^a$ for all $0\le a\le r_s$. Consequently~\eqref{uf:eq:weighted-product} bounds $P_2$ in $H^s$ for $U\in H^s$. In particular, $P_{20}=Z_0\otimes Y_0$ belongs to $H^s$: the derivative identity expresses $DP_{20}$ using only $DQ_0,DZ_0$ multiplied by $Y_0$ and $Y_0\otimes Y_0$. Thus its initial $H^s$ bound follows from $Q_0,Z_0\in H^s$ and $Y_0\in W^{2,\infty}$, without a third derivative of $Y_0$.

\paragraph{{Approximation with the prescribed background.}}
{We use~\eqref{uf:eq:approximation} and the material-coordinate equation~\eqref{uf:eq:label-ODE} in $H^6\times H^s$. The flow reconstructs $(Q_\ell,Z_\ell,P_{2\ell})$ in $H^s$ from their initial values, with Lipschitz differences in $H^{s-1}$ locally in the flow variables. At fixed $\ell$, convolution compensates for the derivative in the source, whose differences lie in $H^{s-2}$. The multiplier bound uses $Z_0,P_{20}\in H^s$ and the $H^s$ transport estimate for the inverse-flow displacement $h_\ell$. The relation $\mathrm D_\ell Y_\ell=L(Dv_\ell)Y_\ell$ controls the material derivative of the upper time weight. Since $(U_\ell,P_{2\ell})$ are bounded in $H^s\subset W^{1,\infty}$, the spatial flux matrices have bounded first derivatives and the paired estimate~\eqref{nr:eq:uniform-energy} applies. The $H^s$ bounds of $(y_\ell,w_\ell,h_\ell)$, the flow formulas and the time derivative bounds close the strict range margins on an interval independent of $\ell$ and $\delta$.}

To pass to the limit in these approximations, we compare $h_{\ell},h_M$ in $H^1$ by their transport equations. Indeed, with $\Delta h=h_{\ell}-h_M$ and $\Delta v=v_{\ell}-v_M$, subtraction gives
\be
 (\del_t+v_{\ell}\cdot\nabla)\Delta h=-(I+Dh_M)\Delta v,\qquad \Delta h(0)=0.
\ee
The source is bounded by $C\|\Delta v\|_{H^1}$, since $Dh_M\in H^{r_s}$ multiplies $H^1$. The $H^1$ transport estimate then controls $\Delta h$ by the time integral of $\|\Delta v\|_{H^1}$. With the same $Y_0$ for both approximations, the reconstruction formula gives
\be
 \|Y_{\ell}-Y_M\|_2\le C\|h_{\ell}-h_M\|_{H^1}.
\ee
Indeed, the difference of the matrices $R$ is bounded in $L^2$ by $\|D(h_{\ell}-h_M)\|_2$, and the difference of the composed initial fields is bounded by the Lipschitz norm of $Y_0$ times $\|h_{\ell}-h_M\|_2$. The estimate~\eqref{uf:eq:product-difference} thus bounds the differences of the solved coefficients in $H^1$ by $C(\|U_{\ell}-U_M\|_{H^1}+\|h_{\ell}-h_M\|_{H^1})$. Multiplication by the derivative of either reference field is bounded on $H^1$, since that derivative belongs to $H^{s-1}$ and $s-1>3/2$. Pairing the difference equation with the positive time weights and combining it with the transport estimate for $\Delta h$ gives a differential inequality for an energy equivalent to $\|\Delta U\|_{H^1}+\|\Delta h\|_{H^1}$, with source $C(\|e_{\ell}\|_{H^1}+\|e_M\|_{H^1})$. In particular, this step uses $\Delta Y$ only in $L^2$ and does not require its derivative or an additional derivative of the fluid variables.
\be
 \|U_{\ell}-U_M\|_{C_tH^1}+\|h_{\ell}-h_M\|_{C_tH^1}
 \le C({\ell}^{-(s-2)}+M^{-(s-2)}).
\ee
The initial velocity difference is controlled by the same bound, while the initial displacement is zero. Interpolation gives convergence below order $s$, including in $C_t W^{1,\infty}$ since $s>5/2$, and the flow formula then gives uniform convergence of $Y_{\ell}$. The limit preserves $Q\ge0$, $Z=QY$ and the particle flow. Its strong convergence below order $s$ identifies the limit at $t=0$ with the prescribed initial data $U_0$; the uniform $H^s$ bounds and the equations give weak $H^s$ continuity up to $t=0$. In particular, $A_{\ell}=\operatorname{Id}+h_{\ell}$ converges in $C_tW^{1,\infty}$, so the density formula passes to the limit pointwise:
\be
 Q(t,x)=(\det DA(t,x))^\theta Q_0(A(t,x)).
\ee
The positive flow Jacobian gives $\{Q(t)=0\}=X(t,\{Q_0=0\})$. Lower-order convergence makes $A,DA$ continuous in time in the uniform norm. Since $Y_0$ is Lipschitz, the reconstruction gives $Y\in C_tL^\infty$ and, using $h(0)=0$,
\be
 \|Y(t)-Y_0\|_\infty\le C\|Y_0\|_{W^{1,\infty}}
 \|h(t)\|_{W^{1,\infty}}\longrightarrow0\qquad(t\downarrow0).
\ee
\paragraph{{Continuity at the initial Sobolev order.}}
{
The time weights are continuous in $L^\infty$ with bounded material derivatives, including where $Q=0$. We use the transport identity~\eqref{uf:eq:transport-weighted-energy} with the exact flow. For the differentiated equations, the Lipschitz commutator bound of Section~\ref{uf:subsec:construction-vacuum} removes spatial convolution, and the symmetric spatial fluxes are integrated by parts using their bounded first derivatives. The energy is continuous between arbitrary times by the argument of Step~6 in that section. Freezing the positive time weight and using weak continuity give $U,h\in C_tH^s$; the equations give $C_t^1H^{s-1}$.
}

\paragraph{{2. Differences of the principal matrices.}}
We write $\Delta$ for the difference between the relativistic and classical solutions, and set $b_1\coloneqq\|Y_{0,\delta}-Y_{0,\rm nr}\|_{W^{1,\infty}}$. The limit $Y_{0,\rm nr}$ belongs to $W^{2,\infty}$ by weak-* compactness. Applying a centered difference quotient to $Y_{0,\delta}-Y_{0,\rm nr}$ and using the uniform bound on its second derivatives gives $b_1\le C\|Y_{0,\delta}-Y_{0,\rm nr}\|_{L^\infty}^{1/2}$ for sufficiently small differences. Thus $b_1\to0$.  We write $g_\delta\coloneqq Y_{0,\delta}\circ A_\delta-Y_{0,\rm nr}\circ A_{\rm nr}$ for the difference of the composed initial fields. Subtracting~\eqref{uf:eq:weighted-product} gives
\begin{subequations}\label{nr:eq:Wtwo-difference}
\begin{align}
 \|\Delta P_2\|_{H^{r_s}}
 &\le C(\|\Delta U\|_{H^{r_s}}+\|\Delta h\|_{H^{r_s}}+b_1),\label{nr:eq:Wtwo-low}\\
 \|\Delta P_2\|_{H^s}
 &\le C(\|\Delta U\|_{H^s}+\|\Delta h\|_{H^s}+\omega_\delta),\label{nr:eq:Wtwo-top}\\
 \omega_\delta
 &\coloneqq \|g_\delta\otimes D(Q_{\rm nr},Z_{\rm nr})\|_{H^{r_s}}.
\end{align}
\end{subequations}
The zeroth-order estimate follows from the Lipschitz extension of $P_2$ on the cone. The differentiated difference is a sum of multipliers times $D\Delta U$ and products involving $\Delta Y$ times $D(Q_{\rm nr},Z_{\rm nr})$. For~\eqref{nr:eq:Wtwo-low}, the latter products are measured in $H^{r_s-1}$. The matrix difference $R(DA_\delta)-R(DA_{\rm nr})$ is bounded in $H^{r_s-1}$ by $C\|\Delta h\|_{H^{r_s}}$. The difference $Y_{0,\rm nr}\circ A_\delta-Y_{0,\rm nr}\circ A_{\rm nr}$ has $H^1$ norm at most $C\|\Delta h\|_{H^1}$ by the Lipschitz bounds for $Y_{0,\rm nr}$ and its first derivative. Since $r_s-1\le1$, multiplication by $DU_{\rm nr}\in H^{r_s}$ controls its product in $H^{r_s-1}$. The other part of $g_\delta$, namely $(Y_{0,\delta}-Y_{0,\rm nr})\circ A_\delta$, has $W^{1,\infty}$ norm at most $Cb_1$ and multiplies $H^{r_s-1}$ with this bound.

For~\eqref{nr:eq:Wtwo-top}, the matrix difference is bounded in $H^{r_s}$ by $C\|\Delta h\|_{H^s}$. The product $g_\delta\otimes D(Q_{\rm nr},Z_{\rm nr})$ is multiplied by $R(DA_{\rm nr})$ and, for the quadratic terms, by $Y_\delta$ or $Y_{\rm nr}$. All these multipliers are bounded on $H^{r_s}$. The same estimates for differences of the solved matrices include a term $C\delta$ from~\eqref{nr:eq:coefficient-parameter}. Their $L^\infty$ differences are bounded by $C(\|\Delta U\|_{H^{r_s}}+\delta)$ by the cone estimate.

We set $d_\delta=\|\Delta U(0)\|_{H^{r_s}}+b_1+\delta$. The order-$r_s$ energy and transport estimates bound $(\Delta U,\Delta h)$ in $C_tH^{r_s}$ by $C d_\delta$. In the equation for $\Delta h$, the source is $-(I+Dh_{\rm nr})\Delta v$ and is controlled by the $H^{r_s}$ algebra. The estimate~\eqref{nr:eq:Wtwo-low} also gives time continuity of the coefficients at order $r_s$ by comparing two times of a single solution.

\paragraph{{Convergence of the form increments.}}
We now show that $\omega_\delta\to0$ in $C_t$ and establish {the asserted convergence} of the increments. The lower estimate and the embedding $H^{r_s}\subset L^\infty$ give $\|g_\delta\|_{C_tL^\infty}\le C(b_1+\|\Delta h\|_{C_tL^\infty})\to0$. The uniform $W^{2,\infty}$ bound on the initial fields and the chain rule also give
\bse
\be
 \|g_\delta\|_{C_tW^{1,\infty}}
 \le C\bigl(b_1+\|h_\delta-h_{\rm nr}\|_{C_tW^{1,\infty}}\bigr)\longrightarrow0.
\ee
Here the convergence of the displacements in $W^{1,\infty}$ follows by interpolation between the lower $H^{r_s}$ convergence and their uniform $H^s$ bounds. By step~1, multiplication by $g_\delta$ is uniformly bounded on $H^2$, while its norm on $H^1$ is bounded by $C\|g_\delta\|_{W^{1,\infty}}$. Interpolating these two multiplication bounds, for $r_s<2$ we obtain
\be
 \|g_\delta f\|_{H^{r_s}}
 \le C\|g_\delta\|_{W^{1,\infty}}^{2-r_s}\|f\|_{H^{r_s}}.
\ee
Consequently $\omega_\delta\to0$ uniformly in time, with exponent $2-r_s=3-s>0$.

At $r_s=2$, we put $B_\delta\coloneqq Y_{0,\delta}-Y_{0,\rm nr}$ and $Q_0\coloneqq Q_{0,\rm nr}$. The initial weighted convergence gives
\be
 Q_0B_\delta=(Z_{0,\delta}-Z_{0,\rm nr})-(Q_{0,\delta}-Q_0)Y_{0,\delta}\longrightarrow0\quad\hbox{in }H^2.
\ee
The product rule makes the weighted control of the second derivatives explicit:
\be
 \|Q_0D^2B_\delta\|_2
 \le C\bigl(\|Q_0B_\delta\|_{H^2}+b_1\|Q_0\|_{H^2}\bigr)
 \longrightarrow0.
\ee
We set $f=D(Q_{\rm nr},Z_{\rm nr})$. For fixed $k>0$, the convergence $A_\delta\to A_{\rm nr}$ in $L^\infty$ gives $Q_0\circ A_\delta\ge k$ on $\{Q_0\circ A_{\rm nr}\ge2k\}$ for small $\delta$. With $E_k(t)\coloneqq\{Q_0\circ A_{\rm nr}(t)<2k\}$, the change of variables and the uniform bounds on $D^2B_\delta$ give
\be
 \|f(t)(D^2B_\delta)\circ A_\delta(t)\|_2
 \le Ck^{-1}\|f(t)\|_\infty\|Q_0D^2B_\delta\|_2
       +C\|f(t)\mathbf 1_{E_k(t)}\|_2.
\ee
\ese
The constants are uniform on the common time interval.  Here $Z_{\rm nr}=Q_{\rm nr}Y_{\rm nr}$ vanishes wherever $Q_{\rm nr}=0$, and the locality of weak derivatives gives $D(Q_{\rm nr},Z_{\rm nr})=0$ almost everywhere on that set. The flow formula identifies it with $X_{\rm nr}(t,\{Q_0=0\})$. To see that $t\mapsto f(t)\circ X_{\rm nr}(t)$ is continuous in $L^2$, we approximate $f(t_0)$ by a smooth compactly supported function. We use the uniform Jacobian bounds and the uniform convergence of the flow maps, followed by the $L^2$ continuity of $f$. Its image on the compact time interval is therefore compact in $L^2$. Since it vanishes on $\{Q_0=0\}$, its $L^2$ norms on $\{Q_0<2k\}$ tend to zero uniformly in $t$ as $k\downarrow0$. The other terms in the product rule tend to zero by $b_1\to0$ and the uniform $H^3$ bound for $h_\delta$.

For the fixed initial field $Y_{0,\rm nr}$, the difference $Y_{0,\rm nr}\circ A_\delta-Y_{0,\rm nr}\circ A_{\rm nr}$ tends to zero in $C_tH^2_{\rm loc}$. In its second derivative the term with $D^2Y_{0,\rm nr}$ follows by approximating this Hessian by smooth functions in $L^2$ on each ball containing the images of the spatial compact set. The approximation uses uniform Jacobian bounds and uniform convergence of the maps. The other terms use $\Delta h\to0$ in $C_tH^2\cap C_tW^{1,\infty}$. Multiplication by either composed field is uniformly bounded on $H^2$ by step~1. Local $H^2$ convergence therefore gives convergence after multiplication by every $f\in C_tH^2$ when we approximate its compact trajectory in $H^2$ by a finite set of smooth compactly supported functions and use the uniform multiplication bound for the errors. Together with the estimate for $(B_\delta\circ A_\delta)D(Q_{\rm nr},Z_{\rm nr})$, this proves $\omega_\delta\to0$ at $s=3$ using the weighted initial convergence alone.

For the unweighted increments at $s=3$, we write $M_j=R(DA_j)-I$ for $j=\delta,{\rm nr}$. After the top-order convergence of $h_j$ proved in Step~3, $M_\delta\to M_{\rm nr}$ in $C_tH^2$. The increment is $M_j(Y_{0,j}\circ A_j)$. The uniform multiplier bound and the fixed-field composition argument reduce convergence of the increment to that of $M_{\rm nr}(B_\delta\circ A_\delta)$ in $C_tH^2$. By $b_1\to0$, the lower-order terms in its first and second derivatives tend to zero in $C_tL^2$. The remaining second derivative is
\be
 M_{\rm nr}(D^2B_\delta\circ A_\delta)[DA_\delta,DA_\delta].
\ee
The two spatial slots are multiplied by uniformly invertible matrices. Changing variables by $X_\delta$ therefore reduces convergence to that of $(M_{\rm nr}\circ X_\delta)D^2B_\delta$ in $C_tL^2$. We may replace $X_\delta$ by $X_{\rm nr}$ in the multiplier, since $M_{\rm nr}\in C_tL^2$, the flow maps converge uniformly with uniformly bounded Jacobians and inverse Jacobians, and $D^2B_\delta$ is bounded. This follows by smooth approximation of the compact trajectory of $M_{\rm nr}$ in $L^2$.

We set $\mathcal V_0=\{Q_0=0\}$. On $\{Q_0\ge k\}$ the product is bounded in $L^2$ by $Ck^{-1}\|Q_0D^2B_\delta\|_2$, and on $\{0<Q_0<k\}$ its norm tends to zero uniformly as $k\downarrow0$ by compactness of the trajectory of $M_{\rm nr}\circ X_{\rm nr}$ in $L^2$. This proves the equivalence in~\eqref{nr:eq:initial-vacuum-increment}. The matrix acts on the minor components of the Hessian. Condition~\eqref{nr:eq:initial-vacuum-hessians} implies this criterion by local $L^2$ convergence on $\mathcal V_0$, boundedness of the Hessians and the uniform $L^2$ tail bound for the multiplier. The same condition gives the local convergence of the unweighted minors themselves at order two: $D^2B_\delta\to0$ in $L^2_{\rm loc}$ by splitting into $\{Q_0\ge k\}$, $\{0<Q_0<k\}$ and $\mathcal V_0$, using the local volume bound for the middle set as $k\downarrow0$. The chain rule with the common Jacobian bounds then gives $B_\delta\circ A_\delta\to0$ in $C_tH^2_{\rm loc}$.

\paragraph{{3. Comparison at order $s$.}}
We smooth the reference fields in space only. We write $\smoothoperator_\sigma$ for the convolution at length $\sigma>0$, and compare $U_\delta,h_\delta$ with $\smoothoperator_\sigma U_{\rm nr},\smoothoperator_\sigma h_{\rm nr}$. With $\mathsf A_{i,\delta}$ evaluated at the exact relativistic fields and $\mathsf A_{i,\rm nr}$ at the exact classical fields, the equation for $U_\delta-\smoothoperator_\sigma U_{\rm nr}$ has source
\be
 \sum_i(\mathsf A_{i,\rm nr}-\mathsf A_{i,\delta})\del_i \smoothoperator_\sigma U_{\rm nr}
       +\sum_i[\smoothoperator_\sigma,\mathsf A_{i,\rm nr}]\del_i U_{\rm nr}.
\ee
For the first term the tame product estimate uses~\eqref{nr:eq:Wtwo-top} times $\|DU_{\rm nr}\|_\infty$, and the $L^\infty$ coefficient difference times $\|D \smoothoperator_\sigma U_{\rm nr}\|_{H^s}\le C\sigma^{-1}$. The latter product tends to zero as $\delta\downarrow0$ at fixed $\sigma$ by the lower convergence.  For the commutator with $\smoothoperator_\sigma$, we set $C_i f\coloneqq[\Lambda^s,\mathsf A_{i,\rm nr}]\del_i f$. Since $\smoothoperator_\sigma$ commutes with $\Lambda^s$ and spatial derivatives, the exact identity is
\be
 \Lambda^s[\smoothoperator_\sigma,\mathsf A_{i,\rm nr}]\del_i U_{\rm nr}
 =[\smoothoperator_\sigma,\mathsf A_{i,\rm nr}]\del_i\Lambda^s U_{\rm nr}
 +(\smoothoperator_\sigma-I)C_iU_{\rm nr}+C_i(I-\smoothoperator_\sigma)U_{\rm nr}.
\ee
The first term tends to zero in $L_t^1L^2$ by the Lipschitz coefficient commutator estimate and approximation of $\Lambda^s U_{\rm nr}$ in $L^2$. The commutator estimate at order $s$ bounds the last term by $C\|(I-\smoothoperator_\sigma)U_{\rm nr}\|_{L_t^1H^s}$. The middle term tends to zero by strong convergence of convolution on $L^2$ and dominated convergence in time. Thus the source commutator tends to zero in $L_t^1H^s$, using only $s$ derivatives of $U_{\rm nr}$.

The displacement comparison has the same two terms, with transport coefficient $v$ and reference $h_{\rm nr}$, and the source $-(v_\delta-\smoothoperator_\sigma v_{\rm nr})$. The displacements vanish at $t=0$, and the initial-data convergence gives
\be
 \limsup_{\delta\downarrow0}\|U_\delta(0)-\smoothoperator_\sigma U_{\rm nr}(0)\|_{H^s}
 =\|(I-\smoothoperator_\sigma)U_{\rm nr}(0)\|_{H^s}\longrightarrow0\quad(\sigma\downarrow0).
\ee
{More precisely, we write $\mathcal D_{\delta,\sigma}(t)$ for the sum of the $H^s$ norms of the differences from the smoothed reference, and $\mathcal S_\sigma(t)$ for the sum of the $H^s$ smoothing errors of the reference. We write $\mathcal C_\sigma(t)$ for the sum of the $H^s$ norms of the smoothing commutators in the two equations. The integrated energy inequality gives
\be
 \mathcal D_{\delta,\sigma}(t)
 \le C\mathcal D_{\delta,\sigma}(0)+
 C\int_0^t\bigl(\mathcal D_{\delta,\sigma}+\mathcal S_\sigma+
     \mathcal C_\sigma+\omega_\delta+\sigma^{-1}d_\delta\bigr)(\tau)\,d\tau
 \qquad (0<\sigma\le1).
\ee
Here~\eqref{nr:eq:Wtwo-top} controls the $H^s$ coefficient difference multiplied by a uniformly bounded first derivative of the smoothed reference. The term $\sigma^{-1}d_\delta$ comes from the $L^\infty$ coefficient difference multiplied by the higher Sobolev norm of that derivative. The constant $C$ is independent of $\sigma$ and $\delta$. Gronwall's inequality, first taking $\delta\downarrow0$ at fixed $\sigma$, bounds the limit superior by $C(\|\mathcal S_\sigma\|_{C_t}+\|\mathcal C_\sigma\|_{L^1_t})$. Both terms tend to zero as $\sigma\downarrow0$, proving convergence in $C_tH^s$.}

Using~\eqref{nr:eq:Wtwo-low} in the evolution equations gives convergence of $\partial_t U_\delta$ in $C_tH^{s-1}$. The displacement equation
\be
 \partial_t h_\delta=-v_\delta-v_\delta\cdot\nabla h_\delta
\ee
and the algebra property of $H^{s-1}$ give the same convergence for $\partial_t h_\delta$. Both conclusions hold up to $t=0$ and throughout the transported zero-density set.

{For $1<\gamma\le3$, the power lemma in the proof of Proposition~\ref{nr:prop:periodic-sobolev-data} applies also on $\mathbb R^3$. It gives convergence of $r=Q^{p_\theta}$ in $C_tH^2$, and the continuity equation $\partial_t r+\operatorname{div}(rv)=0$ gives convergence of $\partial_t r$ in $C_tH^1$. Equations~\eqref{nr:eq:finite-fluxes}--\eqref{nr:eq:renormalized-energy} then give the global convergence of the physical fluxes, including for $\varepsilon_{\rm e}=0$. For every $\gamma>1$, uniform continuity of the power on the bounded range gives uniform convergence of these fluxes. At fixed $\varepsilon_{\rm e}>0$, the reconstructed minors also give uniform convergence in~\eqref{nr:eq:renormalized-material-fluxes}. Both limits preserve the distributional equations and their prescribed initial terms.}

For $s<3$, convergence of the increments follows from convergence of $R(DA)-I$ in $C_tH^{r_s}$ and the multiplier estimate for $g_\delta$ from Step~2. At $s=3$, the criterion in Step~2 gives the endpoint, while the uniform $H^2$ bound and interpolation give convergence below order two. The fixed-field composition argument gives continuity in time for each increment.

\paragraph{{Physical fields and their initial values.}}
On the set of exact spacetime minors, the flow formulas give $\Fsp\in C_tH^{r_s}_{\rm loc}$, and its transport equation gives $\partial_t\Fsp\in C_tH^{s-2}_{\rm loc}$. {In physical-time coordinates $(t,x)$, the reconstruction $F=(-\Fsp v,\Fsp)$ is continuous and locally $W^{1,2}$.} For positive elastic parameter, the determinant, $n=r/\Gamma_\delta$ and the current inherit these Sobolev bounds by multiplication. Interpolation with uniform convergence gives convergence below $r_s$. The composition estimates give the endpoint for $\Fsp$ when $s<3$ or when~\eqref{nr:eq:initial-vacuum-hessians} holds, and the transport equation gives convergence of its time derivative in $C_tH^{s-2}_{\rm loc}$. These conclusions pass to the determinant and current. For $1<\gamma\le3$, the power estimate above gives the density and current endpoint directly. Products with the smooth coefficients of the flux formulas yield the remaining assertions of (iii).

For (iv), we use $r_\delta=(r_{0,\delta}\circ A_\delta)\det DA_\delta$. Composition is bounded on $L^2$ by change of variables and on $H^2$ by the chain rule, using $D^2A_\delta\in L^3$ and $H^1\subset L^6$ for the differentiated test field. Interpolation gives uniform bounds on $H^{s-1}$, and smooth approximation gives strong convergence of the composed densities. The determinant minus one converges in $C_tH^{s-1}$. The continuity equation then gives convergence in $C_tH^{s-1}\cap C_t^1H^{s-2}$, and the smooth coefficients in the flux formulas give the same convergence for the physical fluxes. Locally, cut off the initial densities with a cutoff equal to one on a neighborhood of the images of the chosen spatial compact set under all the inverse flows. This neighborhood can be bounded because the velocities are uniformly bounded. In both cases evaluating the density convergence at $t=0$ gives the converse.
This completes the proof of Proposition~\ref{nr:prop:two-bounded-derivatives}.
\end{proof}

At zero elastic parameter, put $w(Q)=\int_0^Q\sqrt{H_f(\xi)/m}\,d\xi$ and $c(Q)=\theta Q\sqrt{H_f(Q)/m}$. The positive $C^7$ function $H_f$ gives a $C^8$ coordinate and inverse, including at $Q=0$. The classical system~\eqref{nr:eq:classical-normalized} reduces to
\be
 \mathrm D_tw+c(Q)\operatorname{div}v=0,\qquad
 \mathrm D_tv+c(Q)\nabla w=0.
\ee
For $p=Kr^\gamma$, $H_f$ is constant and $c(Q)=\theta w(Q)$. The classical mass density is $mr$, so $w$ is proportional to the density power used by Makino, Ukai, and Kawashima~\cite{MakinoUkaiKawashima1986} at constant thermodynamic entropy (see also~\cite{LeFlochU}). For $5/2<s\le3$ the resulting barotropic pair belongs to $C_tH^s\cap C_t^1H^{s-1}$. For $1<\gamma\le3$, the mass density $mQ^{2/(\gamma-1)}$ is also $C^1$. In vacuum $c(0)=0$ and the velocity obeys $\mathrm D_tv=0$.

The joint support of density and velocity is preserved on the interval of existence, both for this fluid system and for the relativistic fluid reduction. Indeed the transported density keeps its zero set along the particle paths, and on those paths the velocity is constant. The flow therefore carries the initial joint support to the current joint support. It fixes its complement pointwise since both fields initially vanish there. This gives the compact-support property in these fluid reductions (cf.~\cite{MakinoUkaiKawashima1986} and~\cite{LeFlochU}).

%============================================================================ 

\section{\texorpdfstring{Relativistic generalization of Ogden materials}{Relativistic generalization of Ogden materials}}
\label{sec:constitutive-comparison-applications}
\label{sec:vacuum-comparison} 

\subsection{\texorpdfstring{The quadratic case}{The quadratic case}}
\label{subsec:vacuum-law-examples}
\label{subsec:three-invariant-class}

{
We begin with a  quadratic Ogden-type law and then introduce relativistic Ogden-type laws using powers of the singular values and the proper density. We examine the convexity, propagation speeds, and behavior at rank loss of these laws, then verify the stability criterion under explicit assumptions on the strong reference.\par}

We construct isotropic laws from the squared norms of the spacetime minors and the proper density. For selected whole material components, the two form invariants are
\be
\label{eq:vacuum-law-quadratic-invariants}
I_F\coloneqq\sum_A g^{\mu\nu}F_\mu^AF_\nu^A,
\qquad
I_G\coloneqq\frac12\sum_b G_{\mu\nu b}G_b^{\mu\nu}.
\ee
In the rest frame these invariants are $|\Frest |^2$ and $|\Grest |^2$, where $\Grest $ is the spatial Hodge dual of $G$.

\begin{definition}
\label{def:quadratic-elastic-law}
For $m,\kappa>0$ and $\alpha,\beta>0$ on the selected form families, the \underline{\emph{quadratic Ogden-type law}} is the augmented constitutive law
\be
\label{eq:quadratic-vacuum-law}
\raug(F,G,\particlecurrent;g)\coloneqq {m \, n}+\frac\alpha2 I_F+\frac\beta2 I_G+\frac\kappa2 n^2,
\qquad n=\sqrt{-g(\particlecurrent,\particlecurrent)}.
\ee

\end{definition}
Its rest energy is quadratic in the independent form variables and the proper density, apart from the linear rest-mass term.
For the acoustic calculation, let $k\in\mathbb R^3$ be the spatial wave covector, identified with a vector by the Euclidean metric of the rest frame, and let $b\in\mathbb R^3$ be the velocity amplitude. The index $B$ specifies a material component of $G$. The two-form contribution satisfies the identity
\be
\label{eq:two-form-cone-identity}
|k|^2|\Grest _B\times b|^2-|\Grest _B(k\cdot b)-b(\Grest _B\cdot k)|^2
=[k\cdot(\Grest _B\times b)]^2
\ee

\begin{proposition}[Properties of the quadratic Ogden-type law]
\label{prop:quadratic-vacuum-law}
Consider the quadratic Ogden-type law of Definition~\ref{def:quadratic-elastic-law}.
\begin{itemize}
\item[(i)] \emph{Convexity and positivity.} The law is relativistically polyconvex, with positive definite Hessian of the energy with respect to its independent rest-frame variables, positive enthalpy and a positive definite inertia matrix $\mathsf M$ whenever $n>0$. The matrix $\mathsf M$ is defined in~\eqref{eq:rest-inertia-test}.\footnote{In the linearized momentum equations, $\mathsf M$ multiplies the time derivative of the velocity. Write $\mathsf A=D^2_{(\Frest,\Grest,n)}\bar\rho$. If $\mathsf M$ is positive definite and $\mathcal D_k^T\mathsf A\mathcal D_k$ is nonnegative, the squared acoustic speeds are the nonnegative eigenvalues of the symmetric matrix $\mathsf M^{-1/2}\mathcal D_k^T\mathsf A\mathcal D_k\mathsf M^{-1/2}$. Thus the wave speeds are real. For the laws considered here,  the positive definite Hessian of the energy with respect to its independent rest-frame variables and the acoustic bounds give the corresponding hyperbolicity properties, as specified in Lemma~\ref{lem:ogden-acoustic-speeds} below.}
\item[(ii)] \emph{Acoustic bound and entropy at the boundary.} The law satisfies the acoustic bound and the extension hypotheses of Theorem~\ref{thm:closed-vacuum-entropy} (nonnegativity and continuity of the reduced elastic energy up to the closed velocity ball, and infinite inward slope at its boundary).
\end{itemize}
\end{proposition}

{For the quadratic law, the change from the inverse to the direct deformation in Proposition}~\ref{prop:Ball-perspective-new} gives the classical energy per reference volume with reference particle density one:
\be
\label{eq:classical-quadratic-example}
{\rho_B}(A,B,d)=m+\frac{\alpha|B|^2+\beta|A|^2+\kappa}{2d},\qquad d>0.
\ee
{This energy is convex in the independent minors. At $d=1$ its two strain terms are the Mooney--Rivlin invariants $|A|^2$ and $|\cof A|^2$.}

\begin{proof}
The Hessian of the energy with respect to its independent rest-frame variables on the selected blocks is $\operatorname{diag}(\alpha I,\beta I,\kappa)$ and the enthalpy is $m+\kappa n$. For the acoustic inequality the one-form terms cancel between the inertia matrix and the spatial Hessian. The two-form terms leave the nonnegative squares in~\eqref{eq:two-form-cone-identity}, and the density term leaves $\kappa n^2|k\times b|^2$. The rest mass adds ${m \, n}|k|^2|b|^2$, giving strict subluminality for $n>0$. In the observer variables the quadratic part is nonnegative and continuous up to the closed velocity ball: the form invariants are the spatial squares minus the time-component squares, with the time components determined by orthogonality. The part depending on $r={\particlecurrent}^0$ other than the rest-mass term is even in $r$. At a unit boundary velocity $e$, moving to $(1-t)e$ adds the rest-mass contribution $mr\sqrt{2t-t^2}$, whereas the quadratic terms change by $\Ocal(t)$. For $r>0$, division by $t$ gives the infinite inward slope. Thus Theorem~\ref{thm:closed-vacuum-entropy} applies.
\end{proof}

%-----------------------------------------------------------------------------------------------------------------------------------------

\subsection{\texorpdfstring{Ogden-type elastic energies}{Ogden-type elastic energies}}
\label{subsec:spectral-elastic-energies}

\paragraph{\texorpdfstring{Isotropic energies and convexity at rest.}{Isotropic energies and convexity at rest.}}
{We now generalize the preceding quadratic law by considering isotropic functions of the singular values and the proper density. For a real $3\times3$ matrix $A$, write $\sigma(A)=(\sigma_1(A),\sigma_2(A),\sigma_3(A))$ for its singular values in nonincreasing order, as defined in the introduction after~\eqref{eq:intro-local-energy}, and put $\mathbb R_+=[0,\infty)$.}  Ball's spectral criterion for matrices and their minors~\cite{Ball1} suggests a convexity criterion in the two independent rest-frame matrices $\Frest ,\Grest $ and the proper density $n$. With unit reference density the exact minors satisfy $\Grest =\operatorname{cof}\Frest $ and $n=\det \Frest $.

\begin{definition}
\label{def:isotropic-rest-energy}
Let $\varphi:\mathbb R_+^3\times\mathbb R_+^3\times[0,\infty)\to\mathbb R$ be continuous and invariant under independent permutations of its first and second triples. The associated \underline{\emph{isotropic rest energy}} is
\be
\label{eq:spectral-convexity-rest}
\bar\rho(\Frest ,\Grest ,n)\coloneqq\varphi\bigl(\sigma(\Frest ),\sigma(\Grest ),n\bigr),
\qquad (\Frest ,\Grest ,n)\in\mathbb R^{3\times3}\times\mathbb R^{3\times3}\times[0,\infty).
\ee
\end{definition}
\begin{proposition}[Convexity of isotropic rest energies]
\label{prop:spectral-convexity-rest}
The energy of Definition~\ref{def:isotropic-rest-energy} has two properties.
\begin{itemize}
\item[(i)] The rest energy is convex on its closed domain if and only if $\varphi$ is convex and nondecreasing in each of its first six arguments. It is continuous there, including when either rest matrix loses rank and at $n=0$.
\item[(ii)] Under these convexity conditions, a covariant constitutive law with this rest energy and the positive-density constitutive regularity of Definition~\ref{def:constitutive-minor-extension} is relativistically polyconvex. If $\varphi$ is independent of the first two triples, the law is a perfect-fluid law, and the criterion reduces to convexity of $n\mapsto\varphi(n)$.
\end{itemize}
\end{proposition}

{The equivalence in (i) characterizes the particular extension~\eqref{eq:spectral-convexity-rest}, and a given constitutive law may admit other convex extensions.}

\paragraph{{Convexity as a function of the matrix entries.}}
{We first show how the convexity of a symmetric function of the singular values passes to the matrix entries.}
\begin{proof}
We take $\Frest _1,\Frest _2,\Grest _1,\Grest _2$ to be real $3\times3$ matrices, $n_1,n_2\ge0$ and $0\le\theta\le1$. We write $\Frest _\theta=\theta \Frest _1+(1-\theta)\Frest _2$, with analogous notation for $\Grest _\theta$ and $n_\theta$. Ky Fan's maximum principle~\cite{Fan1951}, applied in finite dimension, gives
\be
\label{eq:ky-fan-variational-principle}
\sum_{i=1}^k\sigma_i(A)
=\max_{U^TU=V^TV=I_k}\operatorname{tr}(U^TAV),
\qquad A\in\mathbb R^{3\times3},\quad U,V\in\mathbb R^{3\times k}.
\ee
Applying this identity to a convex combination and taking the maximum gives the following inequalities, for $k=1,2,3$:
\bse
\be
\label{eq:spectral-ky-fan}
\sum_{i=1}^k\sigma_i(\Frest _\theta)
\le\sum_{i=1}^k\bigl(\theta\sigma_i(\Frest _1)+(1-\theta)\sigma_i(\Frest _2)\bigr),
\ee
and the same inequalities hold for $\Grest $. Thus each singular-value triple on the left is weakly majorized by the corresponding convex combination on the right. A symmetric convex function that is nondecreasing in each coordinate is nondecreasing for weak majorization on the nonnegative orthant. Indeed, a weakly majorized vector is a convex combination of vectors dominated componentwise by permutations of the majorizing vector; monotonicity, symmetry and convexity give the assertion. We apply this property successively to the two triples and then use the convexity of $\varphi$:
\be
\label{eq:spectral-majorization-application}
\begin{aligned}
\bar\rho(\Frest _\theta,\Grest _\theta,n_\theta)
&\le\varphi\bigl(\theta\sigma(\Frest _1)+(1-\theta)\sigma(\Frest _2),\,
\theta\sigma(\Grest _1)+(1-\theta)\sigma(\Grest _2),\,n_\theta\bigr)\\
&\le\theta\bar\rho(\Frest _1,\Grest _1,n_1)+(1-\theta)\bar\rho(\Frest _2,\Grest _2,n_2).
\end{aligned}
\ee
\ese
{Conversely, if $\bar\rho$ is convex, we restrict it to pairs of diagonal matrices with nonnegative entries. The permutation invariance gives $\bar\rho(\operatorname{diag}x,\operatorname{diag}y,n)=\varphi(x,y,n)$, so $\varphi$ is convex. Keeping all other entries fixed, allow one diagonal entry of either matrix to vary over $\mathbb R$. The resulting function is convex and even, hence nondecreasing on $[0,\infty)$. This proves the monotonicity in each of the six arguments, also at $n=0$.} Singular values depend continuously on the matrix, so the rest energy is continuous also at repeated or zero singular values. The proof applies directly when either density is zero. The last statements follow from Definition~\ref{def:rest-polyconvexity}, with {the required constitutive regularity} when that definition is used for the evolution equations.
\end{proof}

\paragraph{The time direction and the acoustic cone.}
{The singular values use the positive spatial metric in the material rest space. For an observer measuring $\particlecurrent=r(1,v)$, with $r>0$ and $|v|<1$, let $\Fsp,\Gsp$ be the spatial one-form and dual two-form matrices of Proposition~\ref{prop:rest-entropy-test}. Recall $B(v)=(I-vv^T)^{1/2}$. The change of variables in Proposition~\ref{prop:rest-entropy-test} writes the scalar energy as
\be
\label{eq:spectral-observer-energy}
\widetilde\rho(\Fsp,\Gsp,r,v)=\varphi\bigl(\sigma(\Fsp B(v)),\sigma(\Gsp\operatorname{cof}B(v)),r\sqrt{1-|v|^2}\bigr).
\ee
Where the derivative exists, the energy measured by this observer is $E=\widetilde\rho-v\cdot D_v\widetilde\rho$. At fixed $v$, the reconstruction is linear in $(\Fsp,\Gsp,r)$, so Proposition~\ref{prop:spectral-convexity-rest} gives convexity of $\widetilde\rho$ in these variables. The dependence on $v$ encodes the orientation of the deformation relative to the observer. For the Ogden-type laws below we verify the additional acoustic inequality~\eqref{eq:observer-acoustic-cone} in matter.

}
\paragraph{\texorpdfstring{The relativistic Ogden-type family.}{The relativistic Ogden-type family.}}
In an orthonormal material basis, $\mathsf H_F$ and $\mathsf H_G$ are the Gram matrices of the one-forms and two-forms for the spacetime metric. The factor $1/2$ in the latter accounts for the two orders of each antisymmetric pair.
\be
\label{eq:spectral-material-grams}
(\mathsf H_F)_{AB}\coloneqq g^{\mu\nu}F_\mu^AF_\nu^B,
\qquad
(\mathsf H_G)_{AB}\coloneqq\tfrac12G_{\mu\nu A}G_B^{\mu\nu}.
\ee
In a current rest frame these are $\Frest\Frest^T$ and $\Grest\Grest^T$. In particular, $I_F=\operatorname{tr}\mathsf H_F$ and $I_G=\operatorname{tr}\mathsf H_G$, linking these matrices to the scalar invariants of the quadratic law. For a symmetric positive definite matrix, its real power is defined by the spectral decomposition.

\begin{definition}
\label{def:ogden-power-law}
Let $m>0$ be the rest mass, let $\alpha,\beta,\kappa>0$ be constitutive coefficients, and let $p,q,s>1$ be real exponents. Then the \underline{\emph{relativistic Ogden-type elastic law}} has rest energy
\be
\label{eq:spectral-unshifted-vacuum}
\bar\rho(\Frest ,\Grest ,n)={m \, n}+\frac\alpha p\sum_{i=1}^3\sigma_i(\Frest )^p
+\frac\beta q\sum_{i=1}^3\sigma_i(\Grest )^q+\frac\kappa{s} n^{s}.
\ee
\end{definition}
{The exponents $p$ and $q$ act on the singular values of the two rest matrices, while $s$ governs the density contribution.  The choice $p=q=s=2$ recovers Definition~\ref{def:quadratic-elastic-law}.}

\paragraph{\texorpdfstring{Weak--strong stability for the Ogden-type laws.}{Weak--strong stability for the Ogden-type laws.}}
{We now verify weak--strong stability for the Ogden-type law~\eqref{eq:spectral-unshifted-vacuum}, whose elastic terms have no vanishing density factor. This is a separate application of the general comparison theorem. The following statements specify the reference regularity and the minor or adjoint identities used at rank loss; they do not assert local existence for these laws.} Let $\mathcal L$ be the sum of the two elastic terms in~\eqref{eq:spectral-unshifted-vacuum}, expressed in the covariant Gram matrices~\eqref{eq:spectral-material-grams}. The form conjugates are $P_*=\del_F\mathcal L$ and $Q_*=\del_G\mathcal L$. In an orthonormal rest frame each conjugate uses the derivative of the corresponding spectral energy, obtained by replacing each singular value $\sigma$ by $\sigma^{s_0-1}$ for $s_0=p,q$. This defines the conjugates continuously at rank loss. Contractions involving the two-form derivative sum over both orders of each antisymmetric spacetime index pair. We define the divergence coefficient by
\be
\label{og:spectral-div}
K_*^{\lambda b}\coloneqq2\nabla_\nu Q_*^{\nu\lambda b}.
\ee
In an orthonormal rest frame of the reference, the two relations used in the comparison are
\be
\label{og:eq:mixed-minor-relations}
 {\particlecurrent}^i+(\Gsp)_{Ai}F_0^A=0,\qquad
 G_{0ib}-\Omega_{bAC}F_0^A(\Fsp)_{Ci}=0 .
\ee
Exact minors satisfy both identities. One may also impose just these relations on independent fields in that frame. For instance, if both currents are at rest in this frame and the time components of both form families vanish, the relations place no restrictions on their spatial minors.

\paragraph{{Statement of weak--strong stability.}}
{The material equations for the strong solution give the cancellation needed to estimate its relative entropy remainder.}
\begin{theorem}[Weak--strong stability for relativistic Ogden-type laws with cavitation]
\label{og:theorem}
Let $1<p,q,s\le2$ and consider the law of Definition~\ref{def:ogden-power-law} on a smooth periodic Cauchy slab. Let $F_*,G_*,u_*,H_*=m+\kappa n_*^{s-1}$ be bounded Lipschitz reference fields on the closed slab, with exact polynomial minors, ${\particlecurrent}_*=n_*u_*$, future unit $u_*$ and the constitutive stress of the finite-velocity branch. Require the orthogonality condition $\iota_{u_*}F_*^A=\iota_{u_*}G_{*A}=0$ of Definition~\ref{def:transversality}, the exterior equations, current and stress conservation, and the continued material equations~\eqref{eq:physical-continued-material-EL} throughout the slab, including at $n_*=0$, with their prescribed initial data sets. Define $P_*,Q_*$ by the continuous derivatives of the powers of singular values and $K_*$ by~\eqref{og:spectral-div}.  If $q<2$, assume $K_*\in L^1(0,T;L^\infty_x)$ in the sense that $K_*$ is represented by a measurable field with integrable spatial essential supremum.
\bei
\item {Bound by the relative entropy.} Outside one spacetime null set, the absolute value of the  remainder in the divergence of the relative entropy current is bounded by $a(t)D_\eta(U\mid U_*)$ for every point of ${\Cl(\mathcal Z_{\raug})}$ satisfying~\eqref{og:eq:mixed-minor-relations}, where $a(t)=C(1+\|K_*(t)\|_\infty)$. 
\item {Comparison with the strong reference.} Every entropy weak solution satisfying~\eqref{og:eq:mixed-minor-relations} almost everywhere, with integrable fluxes, initial data in the completed constitutive graph and finite initial relative entropy, satisfies~\eqref{eq:general-constitutive-entropy-stability} with this function $a$. The two zero-current regions may differ.
\item {Weak--strong uniqueness.} Equal initial data sets imply equality of $F,G,\particlecurrent,T$.
\eei

\end{theorem}

{This includes the quadratic law of Definition~\ref{def:quadratic-elastic-law} at $p=q=s=2$.} {The prescribed velocity also specifies the reference branch at zero current. For $p=q=s=2$, the bound on $K_*$ follows from the Lipschitz regularity of $G_*$.}

{
\paragraph{{Comparison with independent augmented fields.}}
{For independent fields we impose the adjoint identity on the strong solution and specify the integrability of the form multipliers.}
\begin{corollary}[Weak--strong stability for augmented Ogden-type laws]
\label{cor:Ogden-fully-augmented-stability}
Let $1<p,q,s\le2$ and consider the law of
Definition~\ref{def:ogden-power-law}. Put
\be
 \mathsf Z_a(X)=D\bigl(\operatorname{tr}((X^TX)^{a/2})/a\bigr),
 \qquad a'=a/(a-1).
\ee
Here $X$ is a real $3\times3$ rest matrix.  At zero current the prescribed
finite-velocity continuation selects the rest frame.
Let $Z_*=(F_*,G_*,\particlecurrent_*,T_*)$ collect bounded reference fields with cavitation and their constitutive stress. Suppose that these fields and stress satisfy the balance laws,
mathematical entropy equality and conditions for attainment of the initial data of
Definition~\ref{def:entropy-strong-solution}. Require its supporting covector
to be finite, with the Sobolev regularity of the {variant given after}
Theorem~\ref{thm:physical-entropy-vacuum} for the elastic multipliers.
Assume the following properties:
\begin{itemize}
\item its finite-velocity continuation satisfies all the augmented balance
laws and the adjoint identity, in the sense of distributions,
\eqref{eq:general-entropy-adjoint-identity} throughout the slab, including on
the constitutive boundary, with zero right-hand side;
\item { the derivatives of the Ogden-type energies} satisfy
\be
 \mathsf Z_p(\Frest _*)\in W^{1,p'},\qquad
 \mathsf Z_q(\Grest _*)\in W^{1,q'}
\ee
on the spacetime slab and have bounded representatives and the prescribed initial data sets; the changes between the continued rest frame and the chosen
observer, as well as all remaining entropy and involution coefficients, are
bounded Lipschitz functions on the closed slab.
\end{itemize}
Fix an entropy-admissible augmented weak solution with finite initial relative entropy. Assume that its fields and associated stress, including their prescribed initial values, belong to the same completed constitutive graph.
Then the estimate~\eqref{eq:general-constitutive-entropy-stability} holds.

No exact-minor relation $G=G(F)$,
$\particlecurrent=\particlecurrent(F)$, or mixed exact-minor identity is
required. The orthogonality and wedge identities defining the augmented
admissible set are  imposed. The ranks of $F$ and $G$, and the zero-current
regions, may differ between the two solutions. No derivative of a field of
the weak solution is assumed, and the particle-number inequality need not be
an equality. Zero initial relative entropy implies equality of
$F,G,\particlecurrent,T$ almost everywhere. In the quadratic case
$p=q=s=2$,  the derivatives of the Ogden-type energies are linear, so their Sobolev
regularity follows from the Lipschitz regularity of the strong reference.
\end{corollary}

\begin{proof}
The assumed Sobolev regularity of the rest-frame conjugates, together with
the Lipschitz changes of frame, is dual to the $L^p$ and $L^q$ control of
the corresponding form fluxes.  The row-dependent Sobolev variant following
Theorem~\ref{thm:physical-entropy-vacuum}, proved by mollification and
H\"older's inequality in Lemma~\ref{thm:common-current-jhde}, therefore
applies with exponents $(p,p')$ to the $F$-rows and $(q,q')$ to the $G$-rows,
including at rank loss. For the prescribed initial data in this graph, $\mathscr I^0_0=-|g(0,x)|^{1/2}N_0^0=\eta_{0,x}(U_0)$, so the initial term in the comparison is precisely the integral of $D_{\eta,0}$. The form evolution fluxes,
the spatial involution fluxes and the particle-number mathematical entropy
current are affine in the augmented fields, so their constitutive remainders
vanish.

The estimates~\eqref{og:eq:stress-power} and~\eqref{og:eq:boost-control},
together with the large-field argument in the proof of
Theorem~\ref{og:theorem}, control the only nonaffine flux remainder, namely
the stress remainder, by the relative entropy.  Here $\Theta=0$ because no
exact-minor constraint is imposed; orthogonality and the wedge identities
remain part of the completed constitutive graph.  The adjoint identity
cancels the complete first variation, so
\eqref{eq:general-entropy-weighted-continuity} holds independently of a
positive lower density or rank bound.  Theorem~\ref{thm:physical-entropy-vacuum}
gives the estimate, and {the same coercive bounds determine the fields and stress when the relative entropy is zero.}
\end{proof}

For exponents below two the Sobolev assumption on the constitutive derivatives applies at deficient rank as well. For the quadratic law the derivatives are linear in the rest matrices.
}

\paragraph{\texorpdfstring{Convexity and propagation speeds.}{Convexity and propagation speeds.}}
{We first verify the constitutive convexity and then quantify the acoustic bound for this family.}
{
\begin{proposition}[Relativistic polyconvexity of Ogden-type laws]
\label{prop:ogden-power-properties}
{The relativistic Ogden-type law of Definition~\ref{def:ogden-power-law}, including its quadratic case $p=q=s=2$, has the following properties.}
\begin{itemize}
\item[(i)] For $p,q,s>1$, its rest energy is convex and continuous on the closed domain of independent rest variables, including at zero density and when either rest matrix loses rank. Wherever its covariant constitutive law has the regularity required in Definition~\ref{def:rest-polyconvexity}, the law is relativistically polyconvex.
\item[(ii)] If $1<p,q,s\le2$, the inertia matrix is positive definite at every triple of independent fields with $n>0$. The Hessian of the energy in its independent rest variables is positive definite wherever the selected matrices with exponent below two are invertible. A matrix with exponent two may have any rank. On the set of exact spacetime minors, $n>0$ implies that both full rest matrices are invertible, so both conclusions apply there.
\end{itemize}
\end{proposition}
}
\begin{proof}

For a singular-value decomposition $A=U\operatorname{diag}(\sigma_i)V^T$ of an invertible matrix, put $h_{ij}=(U^THV)_{ij}$ for a variation $H$. For a smooth function $g$ on $(0,\infty)$ write $f(t)=g(\sqrt t)$. Differentiating the spectral function first at distinct singular values gives
\be
\label{eq:spectral-rest-hessian}
\begin{aligned}
D^2\operatorname{tr}f(AA^T)[H,H]
={}&\sum_i g''(\sigma_i)h_{ii}^2\\
&+\tfrac12\sum_{i<j}\frac{g'(\sigma_i)-g'(\sigma_j)}{\sigma_i-\sigma_j}(h_{ij}+h_{ji})^2
+\tfrac12\sum_{i<j}\frac{g'(\sigma_i)+g'(\sigma_j)}{\sigma_i+\sigma_j}(h_{ij}-h_{ji})^2.
\end{aligned}
\ee
The divided difference in the first sum over pairs takes its continuous value at repeated positive singular values. The energy~\eqref{eq:spectral-unshifted-vacuum} is convex on the closed domain of independent rest variables by Proposition~\ref{prop:spectral-convexity-rest}. Assume first that $1<p,q,s\le2$. First take independent invertible rest matrices and $n>0$. Their Gram matrices remain positive definite on a neighborhood, where their real powers define the smooth spacetime energy. The Hessian of the energy with respect to its independent rest-frame variables $\mathsf A=D^2_{(\Frest,\Grest,n)}\bar\rho$ is positive definite there. Indeed, the coefficients in the spectral Hessian formula~\eqref{eq:spectral-rest-hessian} are positive for $g(t)=t^{s_0}/s_0$, $s_0=p,q$, since $g'(t),g''(t)>0$ at positive singular values, including the limiting divided differences at repeated values. At exponent two the unit-coefficient Hessian block is the identity even at deficient rank. The corresponding block of $\mathsf A$ is $\alpha I$ for the one-forms or $\beta I$ for the two-forms. For an exponent strictly between one and two, the matrix energy has a continuous first derivative at deficient rank; a perturbation in a missing singular direction contributes $|t|^{s_0}/s_0$ and precludes a finite Hessian there. The density block is $\kappa(s-1)n^{s-2}>0$. The form contributions to the inertia matrix are nonnegative, and its density contribution is $n(m+\kappa n^{s-1})I_3>0$. This proves the proposition. The acoustic estimates and their limits at deficient rank are established in Lemmas~\ref{lem:ogden-acoustic-speeds} and~\ref{lem:ogden-rank-loss}.
\end{proof}

{
\paragraph{{Acoustic speeds and their bounds.}}
{We can bound the eigenvalues using the quadratic form of the spatial symbol and the inertia matrix.}
\begin{lemma}[Bounds on the acoustic speeds]
\label{lem:ogden-acoustic-speeds}
Consider the law of Definition~\ref{def:ogden-power-law} with $1<p,q,s\le2$, including the quadratic case $p=q=s=2$. At a triple of independent fields with $n>0$, assuming invertibility only for  selected rest matrices whose exponent is below two, let $\mathsf A=D^2\bar\rho$ be the Hessian of the energy with respect to its independent rest-frame variables and $\mathsf M$ the inertia matrix. For a unit spatial wave covector $k$, let $\mathcal D_k$ be the acoustic variation map of~\eqref{eq:retained-acoustic-map}. One has:
\begin{itemize}
\item[(i)] The three squared acoustic speeds are the generalized eigenvalues of $(\mathcal D_k^T\mathsf A\mathcal D_k,\mathsf M)$. For each such eigenvalue,
\be
\label{eq:ogden-acoustic-speed-bound}
0\le c_{\rm ac}^2\le
1-\frac{{m \, n}+\kappa(2-s)n^{s}}{\lambda_{\max}(\mathsf M)}<1.
\ee
\item[(ii)] {The mathematical entropy form $K_t$, computed at rest in~\eqref{eq:rest-K-test}, is positive definite for every future timelike observer when $n>0$.}
\end{itemize}
\end{lemma}
}
\begin{proof}
{To bound the wave speeds, put $\bar\rho_{s_0}(X)=\operatorname{tr}((XX^T)^{s_0/2})/s_0$ for either form matrix, with $s_0=p$ or $q$. We use the \emph{concavity of trace powers}: for positive definite symmetric matrices $K_0,K_1$,
\be
\label{eq:ogden-trace-concavity}
\operatorname{tr}\bigl((\lambda K_0+(1-\lambda)K_1)^\theta\bigr)
\ge\lambda\operatorname{tr}(K_0^\theta)+(1-\lambda)\operatorname{tr}(K_1^\theta),
\qquad 0<\theta\le1,\quad0\le\lambda\le1.
\ee
For completeness, write $\Psi_\theta(K)=\operatorname{tr}(K^\theta)$, and let $\lambda_i>0$ be the eigenvalues of $K$. In an orthonormal eigenbasis, a symmetric increment $L$ satisfies
\be
\label{eq:ogden-trace-second-variation}
D^2\Psi_\theta(K)[L,L]
=\theta(\theta-1)\sum_i\lambda_i^{\theta-2}L_{ii}^2
+2\theta\sum_{i<j}\frac{\lambda_i^{\theta-1}-\lambda_j^{\theta-1}}{\lambda_i-\lambda_j}L_{ij}^2\le0,
\ee
where divided differences take their continuous values at repeated eigenvalues. Both coefficients are nonpositive, proving the concavity. With $K=XX^T$, $L=XY^T+YX^T$ and $\theta=s_0/2$, the chain rule gives
}
\paragraph{{Acoustic consequences of trace concavity.}}
\bse
\be
\label{eq:ogden-trace-concavity-step}
{D^2\bar\rho_{s_0}(X)[Y,Y]
=\frac1{s_0}D^2\Psi_\theta(K)[L,L]+\operatorname{tr}(K^{\theta-1}YY^T).}
\ee
{The first term is nonpositive and vanishes at $s_0=2$. This proves the first estimate below. The second uses the acoustic variation $\mathcal D_k b$ and the inertia matrix $\mathsf M$ already defined.}
\be
\label{eq:unshifted-shape-hessian}
D^2\bar\rho_{s_0}(X)[Y,Y]\le\operatorname{tr}\bigl((XX^T)^{s_0/2-1}YY^T\bigr),
\ee
\be
\label{eq:unshifted-observer-gap}
|k|^2 b^T\mathsf M b-(\mathcal D_k b)^T\mathsf A\,\mathcal D_k b
\ge {\bigl({m \, n}+\kappa(2-s)n^{s}\bigr)}|k|^2|b|^2,
\qquad k,b\in\mathbb R^3.
\ee
\ese
For the one-form variation $Y=(\Frest b)\otimes k$, the right-hand side of~\eqref{eq:unshifted-shape-hessian}, multiplied by $\alpha$, equals its contribution to $|k|^2 b^T\mathsf M b$. For the two-form variation $Y=\Grest(k\cdot b)-(\Grest k)\otimes b$, apply~\eqref{eq:two-form-cone-identity} to each row of $(\Grest\Grest^T)^{(q/2-1)/2}\Grest$. After multiplication by $\beta$, the same expression is bounded by the two-form contribution to the inertia matrix. The density terms leave
${m \, n}|k|^2|b|^2+\kappa n^{s}\bigl(|k|^2|b|^2-(s-1)(k\cdot b)^2\bigr)$.
Since $s\le2$, Cauchy's inequality bounds the last bracket below by $(2-s)|k|^2|b|^2$, proving~\eqref{eq:unshifted-observer-gap} for the unshifted powers.  In particular, for $|k|=1$ every squared acoustic speed $c_{\rm ac}^2$, given by a generalized eigenvalue of $(\mathcal D_k^T\mathsf A\mathcal D_k,\mathsf M)$, satisfies~\eqref{eq:ogden-acoustic-speed-bound}.
Since $\bar\rho_n=m+\kappa n^{s-1}>0$, Proposition~\ref{prop:observer-acoustic-coercivity} gives a positive definite form $K_t$ for every future timelike observer in matter. The gap from the light speed may tend to zero as $n\downarrow0$.

For $p=q=s=2$, the computation gives an explicit symmetric matrix whose eigenvalues determine the speeds.  We write $\mathsf C_F=\Frest^T\Frest$ and $\mathsf C_G=\Grest^T\Grest$, and denote by $g_B$ row $B$ of $\Grest$, viewed as a spatial vector. For $|k|=1$, one has
\bse
\be
\label{eq:ogden-quadratic-inertia}
\mathsf M=({m \, n}+\kappa n^2)I_3+\alpha \mathsf C_F+
\beta\bigl((\operatorname{tr}\mathsf C_G)I_3-\mathsf C_G\bigr),
\ee
\be
\label{eq:ogden-quadratic-acoustic-gap}
\mathsf M-\mathcal D_k^T\mathsf A\mathcal D_k
={m \, n}I_3+\kappa n^2(I_3-k\otimes k)
 +\beta\sum_B(k\times g_B)\otimes(k\times g_B).
\ee
\ese
Thus the squared speeds are the eigenvalues of
\be
\label{eq:ogden-quadratic-speed-matrix}
I_3-\mathsf M^{-1/2}
\left[{m \, n}I_3+\kappa n^2(I_3-k\otimes k)
+\beta\sum_B(k\times g_B)\otimes(k\times g_B)\right]\mathsf M^{-1/2}.
\ee
The matrix in brackets is bounded below by ${m \, n}I_3$ in matter. Thus the propagation speeds remain strictly below the speed of light also at the endpoint $p=q=s=2$.

This completes the proof of Lemma~\ref{lem:ogden-acoustic-speeds}.
\end{proof}

\paragraph{\texorpdfstring{Rank loss and zero current.}{Rank loss and zero current.}}
{The acoustic estimate has a limit even when the inertia matrix loses invertibility; the limiting stress also depends on how the constitutive boundary is approached.}

{
\begin{lemma}[Acoustic forms and the inertia matrix at rank loss]
\label{lem:ogden-rank-loss}
{For the law of Definition~\ref{def:ogden-power-law} with $1<p,q,s\le2$, consider independent rest matrices $\Frest,\Grest$ and $n\ge0$ in a fixed timelike rest frame. One has:
\begin{itemize}
\item[(i)] The acoustic forms on the variations $XL$ extend continuously at every loss of rank of either matrix. Here $L=b\otimes k$ for $X=\Frest$ and $L=(k\cdot b)I_3-k\otimes b$ for $X=\Grest$. The inertia matrix extends continuously, and~\eqref{eq:unshifted-observer-gap} holds at all these limits.
\item[(ii)] { Recall the Gram matrices $\mathsf C_F=\Frest^T\Frest$ and $\mathsf C_G=\Grest^T\Grest$. The inertia matrix is given by}
\be
\label{eq:ogden-independent-inertia}
\mathsf M=n(m+\kappa n^{s-1})I_3+\alpha \mathsf C_F^{p/2}+\beta\bigl(\operatorname{tr}(\mathsf C_G^{q/2})I_3-\mathsf C_G^{q/2}\bigr).
\ee
For $n>0$ this matrix is positive definite, irrespective of either form rank. At $n=0$, its kernel is
\be
\label{eq:ogden-independent-kernel}
\ker\mathsf M=\ker\Frest\cap\mathcal K_G,\qquad
\mathcal K_G=\begin{cases}
\mathbb R^3,&\operatorname{rank}\Grest=0,\\
\operatorname{range}\Grest^T,&\operatorname{rank}\Grest=1,\\
\{0\},&\operatorname{rank}\Grest\ge2.
\end{cases}
\ee
On $(\ker\mathsf M)^\perp$, when nonzero, the acoustic quadratic form and the quadratic form defined by the inertia matrix determine squared propagation speeds in $[0,1]$. The acoustic form vanishes on $\ker\mathsf M$.
\item[(iii)] On the set of exact spacetime minors, $\Grest=\operatorname{cof}\Frest$ and $n=\det\Frest\ge0$. The alternatives are
\be
\label{eq:ogden-exact-rank-cases}
\begin{array}{c|cccc}
\operatorname{rank}\Frest&3&2&1&0\\
\operatorname{rank}\Grest&3&1&0&0\\
n&>0&0&0&0.
\end{array}
\ee
At $n=0$ the kernel in (ii) reduces to $\ker\Frest$. The upper acoustic endpoint can be attained. On a finite-velocity branch, $F=G=\particlecurrent=0$ gives $T=0$.
\end{itemize}
All conclusions include $p=q=s=2$. A nonzero null current has no timelike rest frame; the present classification concerns finite-velocity branches. Limits with null current or residual stress are distinguished below.}
\end{lemma}
}
\paragraph{{Verification at rank loss.}}
{The acoustic variations are products of each rest matrix with a matrix depending on the propagation direction and amplitude.}
\begin{proof}
{{\bf 1.} \emph{Continuous limits of the acoustic and inertia forms.}}

The acoustic variations of the two rest matrices have the form $XL$ with $L=b\otimes k$ and $L=(k\cdot b)I-k\otimes b$, respectively. In these directions the acoustic quadratic form extends continuously through loss of rank. Indeed, in~\eqref{eq:spectral-rest-hessian} with $Y=XL$ one has $h_{ij}=\sigma_i(V^TLV)_{ij}$. Each pair term is bounded by $C\max(\sigma_i,\sigma_j)^{s_0}\|L\|^2$ for the power $1<s_0\le2$. A pair with one zero singular value contributes $\sigma_i^{s_0}|(V^TLV)_{ij}|^2$, with $\sigma_i>0$ and $\sigma_j=0$. These are the terms obtained by differentiating $\bar\rho_{s_0}(X(I+tL))$ on the fixed range of $X$. Compactness of the orthogonal factors in the singular-value decomposition therefore gives the asserted continuity, independently of the bases in the null spaces. The form contributions to the inertia matrix are continuous as well, since $X^T(XX^T)^{s_0/2-1}X=(X^TX)^{s_0/2}$ in matter. Thus~\eqref{eq:unshifted-observer-gap} passes to rank loss as a quadratic-form inequality. {\vskip.2cm \noindent {\bf 2.} \emph{Independent ranks and the kernel of the inertia matrix.} Formula~\eqref{eq:ogden-independent-inertia} follows by differentiating each matrix power in~\eqref{eq:rest-inertia-test}. All three summands are nonnegative. The kernel of their sum is the intersection of their kernels. If $\mu_1,\mu_2,\mu_3\ge0$ are the eigenvalues of $\mathsf C_G^{q/2}$, the eigenvalues of $\operatorname{tr}(\mathsf C_G^{q/2})I_3-\mathsf C_G^{q/2}$ are $\sum_{j\ne i}\mu_j$. This proves~\eqref{eq:ogden-independent-kernel}, including all independent ranks. By (i), the nonnegative acoustic form is bounded above by $|k|^2\mathsf M$, so it vanishes on $\ker\mathsf M$ and induces a form on its orthogonal complement. The generalized speeds there lie in $[0,1]$.

\vskip.2cm \noindent {\bf 3.} \emph{Restriction to the set of exact spacetime minors.} }{For exact minors at $n=0$, the kernel is $\ker\Frest$. Indeed, the one-form contribution to~\eqref{eq:rest-inertia-test} is $\alpha(\Frest^T\Frest)^{p/2}$. At rank two, the cofactor has the form $\Grest=a\otimes e$, with $e$ a unit vector spanning $\ker\Frest$; its inertia contribution is a nonnegative multiple of the projection onto $e^\perp$. At ranks one and zero the cofactor vanishes. These cases give the asserted kernel by~\eqref{eq:ogden-independent-kernel}.

For $b\in\ker\Frest$, the one-form acoustic variation vanishes. The two-form variation is zero at ranks one and zero because $\Grest=0$; at rank two, writing $b=te$ gives
\be
\label{eq:ogden-rank-two-variation-zero}
 \Grest\bigl((k\cdot b)I_3-k\otimes b\bigr)
 =a\otimes\bigl((k\cdot b)e-(e\cdot k)b\bigr)=0.
\ee
Together with the zero density variation, this also verifies the acoustic kernel directly. Part~(ii) gives the speeds on its orthogonal complement. To attain the upper endpoint, take $\Frest=\operatorname{diag}(\sigma,0,0)$, $\sigma>0$, $\Grest=0$, $n=0$, $k=e_2$ and $b=e_1$. The energy along the acoustic variation is $\sigma^p(1+t^2)^{p/2}/p$, whose second derivative is $\sigma^p$ at zero; the acoustic and inertia forms both equal $\alpha \sigma^p$, so $c^2=1$. At rank zero, all rest fields vanish, and their finite-velocity reconstructions give $F=G=\particlecurrent=T=0$.}
\end{proof}

\paragraph{{Entropy and stress on the boundary.}}
Continuity of $B(v)$ and $\operatorname{cof}B(v)$ for $|v|\le1$ gives a continuous scalar energy in~\eqref{eq:spectral-observer-energy} at $r=0$. At zero current a finite-velocity reference specifies a rest frame, and the admissible flux includes the limiting stress.

For $m>0$ and $1<p,q,s\le2$, the positive-density image in $V=(z_{\rm s},E,P)$ is $E>\Phi_0(z_{\rm s},P)$, including singular forms by continuity. With $\delta E\coloneqq E-\Phi_0(z_{\rm s},P)$, the inward velocity variation of Section~\ref{og:proof-details} below gives
\be
\label{og:eq:current-boundary-trace}
 m^2({\particlecurrent}^0)^2\le \delta E\,(\delta E+2|P|).
\ee
Local Lipschitz continuity of $\Phi_0$ gives a pointwise $1/2$-H\"older bound for ${\particlecurrent}^0$ at each boundary point $V_0$, and a Lipschitz bound when $P_0=0$ since both factors on the right are $\Ocal(|V-V_0|)$. Thus ${\particlecurrent}^0\to0$ at the boundary and Proposition~\ref{thm:lsc-entropy-extension} gives $\eta_{\rm cl}(V)=-{\particlecurrent}^0$ on the graph of finite fluxes, including exact minors and exponents below two.

Rank two can leave positive elastic energy at zero current. At rank zero the finite-velocity branch gives $F=G=\particlecurrent=T=0$. A different limit at the same $(F,G,\particlecurrent)$ can carry null stress even with exact minors along the approximating sequence. For $D>0$, choose $\max\{1/p,2/q-1\}<a<1$ and put $b=(1-a)/2$. Take $\Frest=\operatorname{diag}(n^a,n^b,n^b)$, with determinant $n$ and $\Grest=\operatorname{cof}\Frest$, and boost in direction $e_1$ with Lorentz factor $\sqrt{D/({m \, n})}$. The boosted forms tend to zero: the amplified entries are of orders $n^{a-1/2}$ and $n^{a+b-1/2}$, respectively. The boosted current is $\Ocal(n^{1/2})$. For the stress with the rest-mass term removed, the sum of rest energy and normal stress is $\alpha n^{ap}+2\beta n^{q(a+b)}+\kappa n^{s}$ by~\eqref{augment3}. Both rest quantities tend to zero and the sum is $o(n)$. Multiplication of this sum by the squared Lorentz factor therefore still gives zero in the limit. The transverse  stress components are unchanged by the boost and tend to zero. The remaining mass stress converges to $D(1,e_1)\otimes(1,e_1)$. At its conservative projection $(0,D,De_1)$ the entropy has no finite supporting covector on the full conservative domain. Indeed, hold $z_{\rm s}=0$ and $P=De_1$ fixed. The maximizing velocity in~\eqref{eq:vacuum-hamiltonian} is $\lambda e_1$ for $0<\lambda<1$. Write $\ell=\sqrt{1-{\lambda}^2}$. Stationarity gives $D=r\lambda(m+\kappa(r\ell)^{s-1})/\ell$, so $\ell/r\to m/D$ as $r\downarrow0$.  The energy then satisfies
\be
 E-D=mr\ell+\frac{\kappa}{s}(r\ell)^{s}-D(1-\lambda)
     =\frac{m^2}{2D}r^2+o(r^2).
\ee
The inward entropy slope $-r/(E-D)$ tends to $-\infty$, and the boundary H\"older exponent $1/2$ is sharp there. At $F=G=\particlecurrent=T=0$, the bound $E\ge mr$ gives the supporting covector $(0,-1/m,0)$ for the empty limit on the finite-velocity branch.

\paragraph{Relation to Ogden's laws.}
{Let $y$ be the direct deformation and use the variables $(A,B,d)$ of Section~\ref{sec:Ball-comparison}, with $A=\nabla_X y$, $B=\operatorname{cof}A$ and $d=\det A>0$. Proposition~\ref{prop:Ball-perspective-new} gives} the following convex energy per reference volume for~\eqref{eq:spectral-unshifted-vacuum}:
\be
\label{eq:classical-spectral-example}
 {\rho_B}(A,B,d)=m+
 \frac{\alpha}{p}\,d^{1-p}\sum_i\sigma_i(B)^p+
 \frac{\beta}{q}\,d^{1-q}\sum_i\sigma_i(A)^q+
 \frac{\kappa}{s}\,d^{1-s}.
\ee
{At $p=q=s=2$, this formula is exactly~\eqref{eq:classical-quadratic-example}.} {If $\lambda_i=\sigma_i(A)$ are the principal stretches, the singular values of $B$ are $d/\lambda_i$, up to ordering. At $d=1$ the elastic terms thus give the Ogden powers $-p$ and $q$, up to the additive normalization of the energy~\cite{Ogden1972,DestradeDorfmannSaccomandi2022}. The factors in $d$ specify its volume dependence. At zero density, formula~\eqref{eq:Ball-recession-boundary} gives $\alpha\sum_i\sigma_i(\Frest)^p/p+\beta\sum_i\sigma_i(\Grest)^q/q$ on the independent boundary, including matrices with deficient rank.  For the elastic laws of Section~\ref{sec:uniform-elastic-fluid-existence}, Section~\ref{sec:non-relativistic-vacuum} proves convergence of the evolution. Theorem~\ref{og:theorem} treats these Ogden-type laws for references satisfying its continued equations.}

\paragraph{Consequences of the stability estimate.}

For independent comparison fields the two identities~\eqref{og:eq:mixed-minor-relations} can be replaced by a bound on their contribution to the divergence. Let $\delta_{\particlecurrent}^i,\delta_{G,ib}$ be their left-hand sides and let $b_+\coloneqq\max\{b,0\}$ for a scalar $b$. The integrated positive contribution of these algebraic defects is denoted by $\mathcal P_{\rm alg}$.  We assume that this contribution satisfies
\be
\label{og:eq:minor-defect-source}
 {\mathcal P_{\rm alg}(t)}\coloneqq\int_{\mathbb T^3}\frac{|g|^{1/2}}{H_*}
       \bigl(K_*^{ib}\delta_{G,ib}-\mathcal A_{*i}\delta_{\particlecurrent}^i\bigr)_+\,dx
 \quad\hbox{belongs to }L^1(0,T),
\ee
then the comparison holds with the source $\mathcal P_{\rm alg}$ in~\eqref{eq:relative-entropy-production-bound}. The contraction is unchanged under spatial rotations of the reference frame. The homogeneous estimate thus holds if the contraction is nonpositive almost everywhere.

Write $w_*=u_*/H_*$, $C_*=w_*\wedge P_*$ and $D_*=w_*\wedge Q_*$, with the component conventions of the equation~\eqref{eq:main-vacuum-finite-multipliers}. With $\omega=|g|^{1/2}$ the remainder in the theorem is
\be
\label{og:eq:grouped-remainder}
\mathscr R_*=-\omega\big((\nabla_\nu w_{*\mu})\Delta T^{\mu\nu}+(\nabla_\nu C_{*A}^{\nu\alpha})\Delta F^A_\alpha+(\nabla_\nu D_{*b}^{\nu\alpha\beta})\Delta G_{\alpha\beta b}\big).
\ee
The material equations and~\eqref{og:spectral-div} bound the divergences of $C_*,D_*$ by $a(t)$. For $q=2$, $K_*$ is bounded because $Q_*$ is linear in $G_*$. For $1<q<2$ the bound on $K_*$ is an additional hypothesis. The constant in $a(t)$ depends on the geometry, constitutive coefficients {and the Lipschitz norms of the reference fields}. The dependence on $K_*$ is displayed in $a(t)$. The constant is independent of a positive lower density or lower singular-value bound. The comparison of the future causal relative current is uniform on compact sets of normalized timelike observers.

Consider an open region where the reference gradient has rank at most two. Its rest cofactor matrix $X=\Grest_*$ has rank at most one. For $1<q<2$, a sufficient local condition for the required bound on $K_*$ is
\be
 |DX|\le b(t)|X|^{2-q}\quad\text{almost everywhere on }\{X\ne0\},
 \qquad b\in L^1(0,T).
\ee
Here $D$ denotes the spacetime derivatives in a Lipschitz orthonormal reference rest frame, and $|X|$ is the Frobenius norm. In rank at most one the spectral derivative is $|X|^{q-2}X$, continued by zero. For $\epsilon>0$ the derivative of $(\epsilon^2+|X|^2)^{(q-2)/2}X$ is bounded by $C_q b(t)$ where $X\ne0$, and vanishes almost everywhere on $\{X=0\}$ by locality of weak derivatives. Uniform convergence gives the same weak derivative bound for $|X|^{q-2}X$. The change from the rest frame to a smooth spacetime chart adds only bounded terms to the divergence. This verifies the local bound through the loss of the remaining cofactor rank, with no positive lower bound on $|X|$.

The estimate applies when either solution has zero current in the other's matter. In a region with $\particlecurrent=T=0$, evaluation at $v=0$ in the maximum defining $\Phi_0$ bounds both spatial spectral energies by $E=T^{00}=0$. The bound of the time components by the spatial ones passes from the timelike graph to its finite-flux limits. Hence $F=G=0$ there. In the rest frame of the reference the time component of the relative current, before multiplication by the volume density, then satisfies
\be
H_*\mathscr Q^0
=\alpha(p-1)\bar\rho_p(\Frest_*)+\beta(q-1)\bar\rho_q(\Grest_*)
 +\frac{\kappa(s-1)}{s}n_*^{s} .
\ee
At a reference with ${\particlecurrent}_*=T_*=0$, the form conjugates vanish and $H_*=m$, giving $m\mathscr Q^0=T^{00}-m{\particlecurrent}^0$ in its rest frame. The continued material equations and the bound on $K_*$ for $q<2$ remain the reference conditions.

The stationary rank-two configuration $f=(0,bx^2,cx^3)$, $b,c>0$, has zero current and positive rest energy $\alpha(b^p+c^p)/p+\beta(bc)^q/q$. Its stress and continued material fluxes are constant. The rest gradients $\operatorname{diag}(\epsilon,b,c)$ recover this value as $\epsilon\downarrow0$. In the density-weighted example~\eqref{uf:eq:moving-vacuum-flow}, the rank-two cavity is transported and has $\particlecurrent=T=0$.

%-----------------------------------------------------------------------------------------------------------------------------------------

\subsection{\texorpdfstring{Proof of stability for Ogden-type laws}{Proof of stability for Ogden-type laws}}
\label{sec:elastic-constitutive-estimates}
\label{app:constitutive-proofs}

\label{og:proof-details}

{We now give a proof of Theorem~\ref{og:theorem}.} {We first estimate the matrix power remainders and the stress in a reference rest frame. We then construct the convex entropy at finite-flux limits, use the material equations to control the source, and justify the integral identities defining the weak formulation with their initial terms.}

{\bf 1.} \emph{The matrix power remainder.} For $1<s_0\le2$, we set $W_{s_0}(X)=\operatorname{tr}((X^TX)^{s_0/2})/s_0$, $Z_{s_0}=DW_{s_0}$, $S_{s_0}(X)=(X^TX)^{s_0/2}$, and denote the first-order remainder of $W_{s_0}$ at $X_*$ by $B_{s_0}(X\mid X_*)$. The gradient $Z_{s_0}$ replaces each singular value $\sigma$ by $\sigma^{s_0-1}$ in the singular value decomposition. With Frobenius norms and $\Delta X=X-X_*$, we claim that
\bse
\be
\label{og:eq:stress-power}
|S_{s_0}(X)-S_{s_0}(X_*)-DS_{s_0}(X_*)[\Delta X]|\le C_{s_0} B_{s_0}(X\mid X_*),
\ee
\be
\label{og:eq:stress-dual}
|(XX^T)^{(2-s_0)/4}(Z_{s_0}(X)-Z_{s_0}(X_*))|^2\le C_{s_0} B_{s_0}(X\mid X_*).
\ee
\ese
At $s_0=2$ the zeroth matrix power is the identity.  In fact, $Z_2(X)=X$ and $S_2(X)=X^TX$, so
\be
 S_2(X)-S_2(X_*)-DS_2(X_*)[\Delta X]=\Delta X^T\Delta X,
 \qquad |Z_2(X)-Z_2(X_*)|^2=2B_2(X\mid X_*).
\ee
Both identities hold at every rank. To prove the first inequality, we take real symmetric matrices $A,B$ with orthonormal eigenpairs $(a_i,e_i),(b_k,h_k)$ and set $u_{ik}=e_i\cdot h_k$ and $f(t)=|t|^{s_0}$. The derivative of $f(A)$ multiplies matrix entries in the eigenbasis of $A$ by the first divided differences $f[a_i,a_j]$ (first for polynomials, then by polynomial approximation in $C^1$). For $R_f=f(B)-f(A)-Df(A)[B-A]$ this gives
\bse
\be
\label{og:eq:spectral-trace}
\operatorname{tr}R_f=\sum_{ik}u_{ik}^2 B_f(b_k\mid a_i),
\ee
\be
\label{og:eq:spectral-entries}
(R_f)_{ij}=\sum_k u_{ik}u_{jk}(b_k-a_i)(b_k-a_j)f[a_i,b_k,a_j].
\ee
\ese
Here $B_f$ is the scalar first-order remainder. The second identity subtracts the affine interpolant at $a_i,a_j$, using the tangent when they coincide. A summand with a zero product vanishes. For $M=\max(|a|,|b|,|c|)>0$, we have $|f[a,b,c]|\le C_{s_0} M^{s_0-2}$ by homogeneity and continuity on $M=1$ (triple coincidences occur only at $\pm1$). Similarly $B_f(b\mid a)\ge c_{s_0}\max(|a|,|b|)^{s_0-2}|b-a|^2$ for $(a,b)\ne(0,0)$ by strict convexity and Taylor expansion on the nonzero diagonal.  For the three arguments $(a_i,b_k,a_j)$ put $M_i=\max(|a_i|,|b_k|)$, $M_j=\max(|a_j|,|b_k|)$ and $M=\max(M_i,M_j)$. Since $s_0\le2$,
\be
M^{s_0-2}\le(M_iM_j)^{(s_0-2)/2}\qquad(M_iM_j>0).
\ee
The two scalar lower bounds for $B_f$ therefore control the product multiplying $u_{ik}u_{jk}$ in~\eqref{og:eq:spectral-entries}\newline by $C_{s_0}\sqrt{B_f(b_k\mid a_i)B_f(b_k\mid a_j)}$. If either maximum is zero, its corresponding difference factor is zero and the summand vanishes. Cauchy's inequality gives $|(R_f)_{ij}|\le C_{s_0}\sqrt{d_i d_j}$ with $d_i=\sum_k u_{ik}^2 B_f(b_k\mid a_i)$, hence $|R_f|\le C_{s_0}\operatorname{tr}R_f$. We apply this to the symmetric matrices with off-diagonal blocks $X_*,X_*^T$ and $X,X^T$ and zero diagonal blocks. A diagonal block of their absolute $s_0$th power is $S_{s_0}$ and its full trace is $2s_0W_{s_0}$. This proves~\eqref{og:eq:stress-power}. In the trace formula each scalar factor $\max(|a_i|,|b_k|)^{s_0-2}$ is bounded below by $(|X|+|X_*|)^{s_0-2}$ when its difference factor is nonzero. Summing and using the Frobenius norm of the difference of the symmetric block matrices gives $B_{s_0}(X\mid X_*)\ge c_{s_0}(|X|+|X_*|)^{s_0-2}|\Delta X|^2$, with value zero when both matrices vanish. In particular the comparison is quadratic on bounded balls.

\vskip.2cm \noindent {{\bf 2.} \emph{The dual weighted estimate.}}
We put $t=s_0/(s_0-1)$, $Z=Z_{s_0}(X)$ and $Z_*=Z_{s_0}(X_*)$. Since $DW_t(Z)=X$ and $W_t(Z)+W_{s_0}(X)=Z:X$, we have $B_{s_0}(X\mid X_*)=B_t(Z_*\mid Z)$. For $g(x)=|x|^t/t$, homogeneity gives $B_g(b\mid a)\ge c_t |a|^{t-2}|b-a|^2$. We use~\eqref{og:eq:spectral-trace} for $g$ with $A,B$ formed from $Z,Z_*$, respectively. The equality $(b_k-a_i)u_{ik}=e_i\cdot(B-A)h_k$ shows that the trace remainder controls $\operatorname{tr}(|A|^{t-2}(B-A)^2)$. Its block on material indices gives~\eqref{og:eq:stress-dual}, since $(ZZ^T)^{(t-2)/2}=(XX^T)^{(2-s_0)/2}$.

\vskip.2cm \noindent {\bf 3.} \emph{Stress in the rest frame of the reference.} We write $X$ for the spatial matrix of one form family of the weak solution, taking the spatial dual for two-forms, and use $Y$ for the components with a time index (a column for one-forms). We set $K=XX^T-YY^T$, $L=XX^T$, $\theta=s_0/2$, $P=K^{\theta-1}$, $P_0=L^{\theta-1}$ and $e=\operatorname{tr}(Y^TPY)$. We first take $K>0$. The energy measured in this frame is $E_{s_0}=\operatorname{tr}(K^\theta)/s_0+e$ for unit elastic coefficient. Concavity of the trace power gives $E_{s_0}-W_{s_0}(X)\ge e/2$. Order reversal for powers between $-1$ and zero gives $D=P-P_0\ge0$ by the resolvent formula. Positivity, trace monotonicity and Cauchy--Schwarz give
\bse
\be
\label{og:eq:gram-control}
\operatorname{tr}(X^TDX)=\operatorname{tr}(K^\theta-L^\theta)+e\le e,\qquad \operatorname{tr}(Y^TDY)\le e,\qquad |Y^TDX|\le e .
\ee
At $Y=0$,  the spatial stress components are linear combinations of $S_{s_0}(X)$, $\operatorname{tr}(S_{s_0}(X))I$ and $W_{s_0}(X)I$. At general $Y$, apart from the common isotropic term $-\operatorname{tr}(K^\theta)I/s_0$, the spatial stress is $X^TPX$ for one-forms and $\operatorname{tr}(X^TPX)I-X^TPX-Y^TPY$ for two-forms. Its change from $(X,0)$ is at most $Ce$ by~\eqref{og:eq:gram-control} and trace concavity. The momentum is a fixed contraction of $Y^TPX$, whose reference linearization uses $Y^T Z_{s_0}(X_*)$. In addition to $Y^TDX$, the remainder involves $Y^T(Z_{s_0}(X)-Z_{s_0}(X_*))$, bounded by $C\sqrt{e B_{s_0}(X\mid X_*)}$ by~\eqref{og:eq:stress-dual} and $P_0\le P$. Thus the linear stress map $A_{s_0*}$ just computed satisfies
\be
\label{og:eq:boost-control}
|\Delta T_{s_0}-A_{s_0*}(\Delta X,Y)|\le C d_{s_0},\quad d_{s_0}\coloneqq E_{s_0}-W_{s_0}(X_*)-Z_{s_0}(X_*):\Delta X\ge B_{s_0}(X\mid X_*)+e/2.
\ee
\ese
The map uses $DS_{s_0}(X_*)$, $Z_{s_0}(X_*)$ and their contractions and is bounded at rank loss.  In~\eqref{og:eq:gram-control}, the mixed bound follows from $|Y^TDX|\le\|D^{1/2}Y\|_{\rm F}\|D^{1/2}X\|_{\rm F}$. On bounded spatial ranges, $0<K\le XX^T\le R^2 I$ implies $P\ge R^{s_0-2} I$, so $e\ge R^{s_0-2}|Y|^2$. At a fixed timelike velocity, we approximate the rest matrices by invertible matrices and boost them back with the same velocity. The energy and all stress components then converge by continuity of the nonnegative matrix powers in the rest frame. Thus the bound by $d_{s_0}$ in~\eqref{og:eq:boost-control}, with the term $e/2$ replaced by $R^{s_0-2}|Y|^2/2$ on bounded spatial ranges, extends by this approximation.

\vskip.2cm \noindent {\bf 4.} \emph{Convex entropy and finite fluxes.} We use~\eqref{eq:spectral-observer-energy}, with $\ell=(1-|v|^2)^{1/2}$ and $P$ the momentum density. The two Gram matrices are $A(I-vv^T)A^T$ and $C((1-|v|^2)I+vv^T)C^T$. They are concave in $v$ in the matrix order: the subtracted quadratic terms have quadratic forms $(x\cdot v)^2$ and $|x\times v|^2$, respectively, before multiplication by the form matrices. The trace powers are increasing and concave on positive semidefinite matrices. The scalar terms $mr\ell$ and $(\kappa/s)(r\ell)^{s}$ are concave in $v$, with strict concavity from the mass term for $r>0$. At fixed $v$, convexity in $(z_{\rm s},r)$ follows from the linear reconstruction and rest convexity. Consequently the maximum $\Phi$ defined by~\eqref{eq:vacuum-hamiltonian} is continuous and convex in $(z_{\rm s},r,P)$. Evaluation at $v=0$ gives $\Phi\ge mr$.

Along $v=\lambda e$ with $|e|\le1$ and $0\le {\lambda}\le1$, both Gram matrices increase as $\lambda$ decreases. The elastic energies therefore increase in this inward direction, including for singular matrices. If $r>0$ and a maximizing velocity $e$ lay on the unit sphere, replacing it {by $\lambda e$} with $\lambda<1$ would gain at least $mr\sqrt{1-{\lambda}^2}$ in scalar energy and lose at most $(1-\lambda)|P|$ in the momentum term. This proves that the unique maximizing velocity at positive density is interior. We now choose a maximizing velocity $e$ for $\Phi_0(z_{\rm s},P)$, with $|e|\le1$. Evaluating $\Phi(z_{\rm s},r,P)$ {at $\lambda e$} and using $\sqrt{1-{\lambda}^2|e|^2}\ge\sqrt{1-{\lambda}^2}$ gives
\be
 \Phi(z_{\rm s},r,P)-\Phi_0(z_{\rm s},P)
 \ge \max_{0\le {\lambda}\le1}\{mr\sqrt{1-{\lambda}^2}-(1-\lambda)|P|\}
 =\sqrt{m^2r^2+|P|^2}-|P|.
\ee
In particular, $\Phi(z_{\rm s},r,P)>\Phi_0(z_{\rm s},P)$ for $r>0$. Evaluating at the interior maximizer for a smaller positive density also gives strict increase between any two positive densities, since $\widetilde\rho_r=\ell(m+\kappa(r\ell)^{s-1})>0$ there.

At fixed timelike velocity,  the stress contributions from the forms extend continuously to singular matrices. Indeed, in the material frame their entries are combinations of the energy and the matrix powers in Step~3, and the change to the observer frame is a fixed Lorentz transformation. The stress identities $P=-D_v\widetilde\rho$ and $E=\widetilde\rho+v\cdot P$ extend there as well, since the spectral energies are $C^1$ in the spatial matrices for $p,q>1$ and the reconstruction is smooth for $|v|<1$. Thus $E=\Phi(z_{\rm s},r,P)$ on all these triples with timelike current. Conversely, for every $E>\Phi_0(z_{\rm s},P)$ there is a unique $r>0$ giving this $E$. Its maximizing velocity is interior, and the first-order condition gives the prescribed $P$ in the reconstructed fields. This proves that the positive-density image with singular forms included is the convex set $E>\Phi_0(z_{\rm s},P)$. At fixed $r,v$, triples with invertible spatial matrices approximate the added timelike triples. Their entropy values therefore lie in the closure of the epigraph for invertible matrices.

With $E=\Phi(z_{\rm s},r,P)$ this bound is~\eqref{og:eq:current-boundary-trace}. The contractions in the spatial reconstruction and the bounded spectral derivatives on bounded sets give local Lipschitz continuity of $\Phi_0$ in $(z_{\rm s},P)$, uniformly under the velocity maximum. Hence $r(V)\le C|V-V_0|^{1/2}$ near each finite boundary point $V_0$. The identity $E=\Phi(z_{\rm s},r,P)$ passes to each finite-flux limit by continuity. If $r>0$ in the limit, any limiting velocity still maximizes the continuous objective on the closed ball. The maximum is interior, so the limiting current is timelike. In particular, zero proper density at a finite-flux limit implies $r=0$. At the triples with timelike current the fixed-velocity supporting inequality of Theorem~\ref{prop:global-rest-current} applies directly, using convexity in $(z_{\rm s},r)$ and their interior maximizing velocity. Proposition~\ref{thm:lsc-entropy-extension} now gives the closed convex entropy~\eqref{eq:closed-current-entropy} and its equality with $-{\particlecurrent}^0$ on the finite-flux graph, including the values prescribed by the initial data.

At a boundary point $V_0=(z_{{\rm s}0},\Phi_0(z_{{\rm s}0},P_0),P_0)$, a finite supporting covector exists precisely when the maximum defining $\Phi_0$ is attained at some $|v_0|<1$. Indeed, with $(z_{{\rm s}0},P_0)$ fixed, let $v_r$ maximize $\Phi$ at density $r>0$ and put $\ell_r=\sqrt{1-|v_r|^2}$. If all maximizers at zero density lie on the unit sphere, compactness and continuity give $\ell_r\to0$. Evaluation of the zero-density objective at $v_r$ gives
\be
 0<\Phi(z_{{\rm s}0},r,P_0)-\Phi_0(z_{{\rm s}0},P_0)
 \le mr\ell_r+(\kappa/s)(r\ell_r)^{s}=o(r).
\ee
Along this ray the entropy decreases by $r$, which rules out a finite supporting slope in the energy direction. Conversely, at an interior maximizer the fixed-velocity energy is convex and differentiable in $(z_{\rm s},r)$ up to $r=0$, also for singular form matrices. Its tangent inequality gives the supporting covector $(D_{z_{\rm s}}\widetilde\rho,-1,v_0)/(m\sqrt{1-|v_0|^2})$ at $V_0$. For the prescribed finite-velocity reference, $v_*$ is a maximizer by continuity from matter, and the same tangent inequality gives $(D_{z_{\rm s}}\widetilde\rho,-1,v_*)/(\ell_*H_*)$ throughout the branch. These are the time coefficients of~\eqref{eq:main-vacuum-relative-current} by the stress identities and continuity on the reference branch. Its nonnegative density in every observer frame gives a future causal current and the comparison $(\tau-|k|)\mathscr Q^0\le t_\nu\mathscr Q^\nu\le(\tau+|k|)\mathscr Q^0$ in a rest coframe, with $t=(\tau,k)$ and $\tau>|k|$.

In the reference rest coframe $H_*\mathscr Q^0=\alpha d_p+\beta d_q+H_*\mathscr Q_{\rm fl}^0$, with nonnegative summands. For the density contribution, the estimates~\eqref{eq:nonlinear-fluid-affine-jet}--\eqref{eq:nonlinear-fluid-jet-remainder} apply with $\gamma=s$. Along with~\eqref{og:eq:boost-control} they bound the stress remainder by $C\mathscr Q^0$ and the squared differences of $F,G,\particlecurrent$ by $C_R\mathscr Q^0$ on bounded field ranges. At large spatial fields or large ${\particlecurrent}^0$, the spatial power remainders and $B_{s}({\particlecurrent}^0\mid n_*)$ grow superlinearly, uniformly on the bounded reference range. The orthogonality condition~\eqref{eq:450bis-2} bounds the time components by the spatial ones. The total measured energy equals $H_*\mathscr Q^0$ plus an affine function of the spatial forms and ${\particlecurrent}^0$ with bounded coefficients. Each elastic stress is dominant: in its material rest frame  the spatial stress eigenvalues lie between $-W_{s_0}$ and $(s_0-1)W_{s_0}$, by the spatial stress formulas above.  Indeed, the positive matrix $X^TP_0X$ has trace $s_0W_{s_0}$; the one-form and two-form stress formulas give these bounds for each of its eigenvalues. Since $1<s_0\le2$, the absolute value of each principal stress is at most the energy density $W_{s_0}$. The fluid stress is also dominant:  its pressure satisfies ${0\le p\le m\,n+(\kappa/s)n^s}$ for $1<s\le2$ and $m\ge0$. Thus the energy controls all stress components and $1+|F|+|G|+|\particlecurrent|+|T|\le C\mathscr Q^0$ on this complementary range.

At an empty reference on this branch, $F_*=G_*={\particlecurrent}_*=T_*=0$. In its rest frame the form terms in $\mathscr Q^0$ control $|\Fsp|^p$ and $|\Gsp|^q$ by~\eqref{og:eq:boost-control}. The density-dependent part of the stress gives the remaining bound: for $r={\particlecurrent}^0$ and $\chi=u^0\ge1$, the scalar energy minus $mr$ equals
\be
 mr(\chi-1)+\frac{\kappa}{s}r^{s}\chi^{-s}
       \bigl(1+s(\chi^2-1)\bigr)
 \ \ge\ \frac{\kappa}{s}r^{s} .
\ee
Indeed, concavity of $x^{s/2}$ for $1<s\le2$ gives $\chi^{s}\le1+\frac{s}{2}(\chi^2-1)$. The orthogonality condition controls the time components of the forms and causality gives $|\particlecurrent|\le\sqrt{2}r$ in this frame. Thus $\mathscr Q^0$ controls $|F|^p+|G|^q+|\particlecurrent|^{s}$ at the empty reference, including for finite-flux limits by continuity.

We have obtained the estimates on the independent timelike graph, using fixed-velocity approximation for singular form matrices. Along any sequence of quadruples of fields and constitutive stress $(F_j,G_j,{\particlecurrent}_j,T_j)$ converging to a finite-flux limit, the relative current and stress remainder converge because they are affine in these fluxes with the reference fixed. The field bounds pass to the limit as well. For exact limiting field values this also gives the constitutive remainder bound~\eqref{eq:general-entropy-Taylor-control}: the polynomial constraint remainder is quadratic on bounded ranges and elsewhere equals $-D\Theta_*\Delta(F,G,\particlecurrent)$. The comparison below uses the total divergence remainder directly.  If $\mathscr Q^0=0$, the estimates identify the forms, current and stress. They include constitutive fields of the weak solution with residual stress at $\particlecurrent=0$.

\vskip.2cm \noindent {\bf 5.} \emph{Material equations and the constitutive source.} We use~\eqref{eq:general-vacuum-joint-Noether} along the given reference. To justify this linear identity at rank loss, we first replace $|x|^{s_0}/s_0$ in its elastic terms by $((\epsilon^2+x^2)^{s_0/2}-\epsilon^{s_0})/s_0$, keeping $F_*,G_*,w_*=u_*/H_*$ fixed. The identity uses closure and the orthogonality condition~\eqref{eq:450bis-2} before the material equations are imposed. The shifted conjugates converge uniformly on bounded reference ranges by convergence of $x(\epsilon^2+x^2)^{s_0/2-1}$ to $\operatorname{sgn}(x)|x|^{s_0-1}$. The spatial stress differentials converge there as well: $x^2(\epsilon^2+x^2)^{s_0/2-1}\to |x|^{s_0}$ in $C^1$ on compact intervals, by the uniform derivative bound $C|x|^{s_0-1}$ near zero. First divided differences give convergence of the matrix differentials.  Indeed, the scalar functions and their derivatives converge uniformly on bounded intervals, so their first divided differences converge uniformly there, including on the diagonal. The spectral formula for the matrix differential then gives convergence also at repeated or zero eigenvalues. The momentum differentials in the rest frame use the conjugates alone, as in Step~3. Covariance gives convergence in the chosen chart. We test the linear identity with smooth compactly supported increments and integrate its divergence terms by parts. A factor such as $w_*^\alpha\nabla_\nu P_*^\nu$ is understood as $\nabla_\nu(w_*^\alpha P_*^\nu)-(\nabla_\nu w_*^\alpha)P_*^\nu$. Uniform convergence and the Lipschitz bound on $w_*$ give the unshifted identity distributionally.

We now use the material equations~\eqref{eq:physical-continued-material-EL} for the actual reference. In differentiating $Q_*^{\nu\lambda b}F^C_{*\lambda}$, we may approximate $F_*$ uniformly with bounded derivatives converging almost everywhere, since $Q_*$ is continuous and $K_*$ is integrable on the slab. The term $Q_*^{\nu\lambda b}\nabla_\nu F^C_{*\lambda}$ vanishes by closure.  The condition $\iota_{u_*}G_*=0$ of Definition~\ref{def:transversality} and the Piola identity give $\mathcal H_A^{\mu\nu}=u_*^\mu\mathcal C_{*A}^\nu-u_*^\nu\mathcal C_{*A}^\mu$ and $\nabla_\nu\mathcal H_A^{\mu\nu}=0$. Thus the divergence of the last term in~\eqref{eq:physical-continued-material-momentum} is $\mathcal C_{*A}^\mu\mathcal A_{*\mu}$. The material equation now gives~\eqref{eq:physical-material-force-equation},  hence $|L_*|\le C(1+|K_*|)$ for $L_*=\nabla_\nu P_*^\nu$. The limiting linear identity on a basis of increments gives the same bound for the divergences of $C_*,D_*$ in~\eqref{eq:main-vacuum-finite-multipliers}. The identity holds outside one null set for every increment. Reference current conservation gives $\nabla_* H_*=-n_*h'(n_*)\nabla_\mu u_*^\mu$ for $h(n)=m+\kappa n^{s-1}$, including almost every point where $n_*=0$, since $H_*=m$ there and $nh'(n)$ extends by zero.

In~\eqref{og:eq:grouped-remainder}, the fluid cancellation~\eqref{eq:forced-fluid-remainder} and the material relation now give~\eqref{eq:physical-nonlinear-additive-force-cancellation} on exact triples with timelike current. The products involving $F_0$ contain two field differences. The remaining term is bounded on bounded ranges using $n|h(n)-H_*|^2\le C_R B_{s}(n\mid n_*)$, $|{\particlecurrent}^{\rm sp}|^2=n^2(\chi^2-1)$ and $n(\chi^2-1)\le C\mathscr Q_{\rm fl}^0$.  For the first inequality, homogeneity gives the stronger scalar bound
\be
 n\bigl|n^{s-1}-n_*^{s-1}\bigr|^2
 \le C_{s}(n+n_*)^{s-1}B_{s}(n\mid n_*),\qquad n,n_*\ge0.
\ee
Indeed, after setting $n+n_*=1$, the ratio of the two sides has finite limits at both endpoints and at $n=n_*$; its diagonal limit is $2^{2-s}(s-1)/s$. This proves the estimate for every fixed $s>1$, including when either density is zero. On a bounded density range the factor $(n+n_*)^{s-1}$ is bounded. For independent timelike approximants, let $\delta_{\particlecurrent}^i$ and $\delta_{G,ib}$ denote the two left-hand sides in~\eqref{og:eq:mixed-minor-relations}. Subtracting the right-hand side of~\eqref{eq:physical-nonlinear-additive-force-cancellation} from the actual force term gives $H_*^{-1}(K_*^{ib}\delta_{G,ib}-\mathcal A_{*i}\delta_{\particlecurrent}^i)$. The coefficient of $\delta_{\particlecurrent}$ is independent of the enthalpy of the weak solution: the terms containing $h$ cancel in the subtraction. Both defects are polynomials in the approximating fields, and their coefficients are bounded by $a(t)$. Along the sequence the polynomials and the original affine remainder converge. On bounded field ranges the resulting upper bound for the remainder is $a(t)D_\eta$ plus the nonnegative density integrated in~\eqref{og:eq:minor-defect-source}. When the polynomial defects vanish, the bound holds also in absolute value. At large fields the bound by $a(t)D_\eta$ alone follows from Step~4.  Both on bounded ranges and at large fields the estimate is linear in the bounds on the divergences and the material forces, {giving the asserted formula for $a(t)$}. For the prescribed foliation the coordinate density is $D_\eta(U\mid U_*)=|g|^{1/2}\,\mathscr Q^0$, where $\mathscr Q$ denotes the vector in~\eqref{eq:main-vacuum-relative-current}. Causality compares this density with the density in the reference rest frame, uniformly on the prescribed compact range of time directions.

\vskip.2cm \noindent {{\bf 6.} \emph{Distributional identities and initial data.}}
We justify testing the exterior equations when the divergences are only integrable in time. We write $\slicenormal$ for the future unit normal to the slices. The dominant stress bounds of Step~4 imply that $-\omega {\slicenormal}_\mu T^{\mu0}$ is nonnegative and controls $|T|$ on the compact slab. Testing the stress balance laws with $-\varphi(t){\slicenormal}_\mu$ gives an absolutely continuous representative for the spatial integral of this density. Its derivative is bounded by a constant times that integral, and its value at zero is the prescribed initial integral. Hence $\sup_t\|T(t)\|_{L^1_x}< +\infty$. In the slice normal frame the elastic energies measured by the observer bound their spatial spectral powers by Step~3, and the mass energy is at least $m{\particlecurrent}^0$. The orthogonality condition and causality then bound $\|F(t)\|_{L^1_x}+\|G(t)\|_{L^1_x}+\|\particlecurrent(t)\|_{L^1_x}$ uniformly by the slice energy and volume.

In each coordinate chart the covariant divergences of the antisymmetric tensors are ordinary divergences after multiplication by $\omega$. Localize the tensor densities in space and smooth first in space. At fixed spatial scale the smoothed divergences are in $L^1_t C_x$ on compact sets, since spatial convolution maps $L^1_x$ continuously to $C_x$. Smooth next in time using values at positive times tending to the given time. As the time scale tends to zero, the resulting divergences converge in $L^1_tC_x$, while the continuous coefficients converge uniformly up to $t=0$. The slice bounds just proved and integrability of the initial data sets justify this limit in the integral identities defining the weak formulation. Remove spatial smoothing next: for almost every time the divergences converge almost everywhere in space, and their spatial essential suprema are bounded by $C a(t)$. Dominated convergence first in space and then in time gives the required pairing with $F,G$. The partition of unity terms cancel upon summation. Testing the stress balance with the Lipschitz $w_*$ then gives the relative entropy inequality with the prescribed initial term. On each slice the contribution of the defects is at most $\mathcal P_{\rm alg}(t)$. The stress estimate and the slice bounds make the original divergence remainder integrable. The scalar argument of Lemma~\ref{thm:common-current-jhde} now gives~\eqref{eq:relative-entropy-production-bound}, with the relative entropy of the prescribed initial data sets. Taking $\mathcal P_{\rm alg}=0$ gives~\eqref{eq:general-constitutive-entropy-stability} and the uniqueness conclusion.
This completes the proof of Theorem~\ref{og:theorem}.

%============================================================================

\section{Technical identities and estimates}
\label{sec:technical-identities-estimates}
\setlength{\abovedisplayskip}{4pt plus 2pt minus 1pt}
\setlength{\belowdisplayskip}{4pt plus 2pt minus 1pt}
\setlength{\belowdisplayshortskip}{3.5pt plus 2pt minus 1pt}

%============================================================================

\subsection{Variational identities for the field equations}
\label{Append-AA}

\paragraph{Variational equations.}
\label{append-A}

{We derive the field equations from the variational and covariance identities and specify when stress balance is equivalent to the material equations. The Noether identity below holds for Lipschitz configurations with locally integrable second derivatives and compact essential constitutive ranges.\par}

\begin{proposition}[Variational equivalence at controlled regularity]
\label{prop:variational-EL-local04}
Let $f\in W^{2,1}_{\rm loc}\cap W^{1,\infty}_{\rm loc}$, and suppose that its constitutive state belongs almost everywhere on each compact set to a compact subset of the admissible constitutive domain. Assume that the scalar physical law is $C^1$ on an open neighborhood of those values and that the evaluated material momentum $P_A{}^\nu=\partial\rphys/\partial F_\nu^A$ belongs to $W^{1,1}_{\rm loc}$. This last condition follows, in particular, if the law is $C^{1,1}_{\rm loc}$. {Then the variational stress $T_\mu{}^\nu=-\rphys\delta_\mu^\nu+F_\mu^AP_A{}^\nu$ satisfies, in the sense of distributions,}
\be
\nabla_\nu T_\mu{}^\nu
=F_\mu^A\left(\nabla_\nu\frac{\del\rphys}{\del F_\nu^A}-\frac{\del\rphys}{\del f^A}\right).
\ee
The stress balance laws hold exactly when the material Euler--Lagrange expression annihilates the image of $F$ almost everywhere. If $F$ has rank three almost everywhere, these laws are equivalent to~\eqref{EL} and to stationarity under compactly supported variations of the extended action.
\end{proposition}
\begin{proof}
The compact constitutive range bounds $P$ and the first derivatives of $\rphys$. The Sobolev chain rule gives $\rphys\in W^{1,1}_{\rm loc}$, and the product rule for the bounded Sobolev fields $F,P$ gives $T_\mu{}^\nu=-\rphys\delta_\mu^\nu+F_\mu^AP_A{}^\nu\in W^{1,1}_{\rm loc}$. If the law is $C^{1,1}_{\rm loc}$, its localized first derivatives are Lipschitz on the compact constitutive ranges, so the chain rule gives the required regularity of $P$. The chain rule gives $\nabla_\mu\rphys=(\del\rphys/\del f^A)F_\mu^A+P_A{}^\nu\nabla_\mu F_\nu^A$, with metric terms accounted for by scalar covariance and $\nabla g=0$. Subtracting the derivative of $\rphys$ from the divergence of $F_\mu^AP_A{}^\nu$ leaves the material expression contracted with $F$ and the term $P_A{}^\nu(\nabla_\nu F_\mu^A-\nabla_\mu F_\nu^A)$. The latter vanishes in $L^1_{\rm loc}$ because $F=df$. The contraction with $F$ vanishes precisely when the material covector $\mathcal E$ vanishes on the image of $F$. At rank three this image is the full material tangent space. {For this variational stress the same identity holds at rank three when the particle current is nonzero and lightlike.}

The first-variation calculation below identifies the distributional material equations with stationarity of the extended action.
\end{proof}

\begin{corollary}[Continuous gradients of bounded variation]
\label{cor:variational-continuous-bv}
Let $f$ be a $C^1$ configuration with $F=df\in BV_{\rm loc}$, and suppose its constitutive state lies locally in the admissible domain of a $C^1$ scalar law. Assume that the evaluated momentum $P$ belongs to $BV_{\rm loc}$, as is the case for a $C^{1,1}_{\rm loc}$ law. Define the material Euler--Lagrange expression as the vector-valued Radon measure
\be
 \mathcal E_A\coloneqq\nabla_\nu\frac{\partial\rphys}{\partial F^A_\nu}-\frac{\partial\rphys}{\partial f^A}.
\ee
The Noether identity of Proposition~\ref{prop:variational-EL-local04} holds as an equality of measures using the continuous gradient. If the stress balance laws hold, the polar covector of $\mathcal E$ annihilates the image of $F$ at $|\mathcal E|$-almost every point. In particular, $\mathcal E$ is supported on the closed set where $F$ has rank below three. If $F$ has rank three $|\mathcal E|$-almost everywhere, the stress balance laws and the variational equations are equivalent.
\end{corollary}
\begin{proof}
Work locally in spacetime and material charts. The fields $F$ and $P$ are continuous. Mollification of $f$ gives uniform convergence of its gradient on compact sets and weak-* convergence of the second derivatives as measures with bounded total variation. We apply the chain rule for the $C^1$ energy to these smooth configurations. Its coefficients converge uniformly on compact sets, so the limit gives the measure chain rule for $\rphys$, with the metric terms as in Proposition~\ref{prop:variational-EL-local04}. Independently mollifying the continuous $BV$ fields $F$ and $P$ gives the product rule for $FP$: in each derivative term one factor converges uniformly and the derivative of the other converges weak-* as a measure. The covariant curl of $F=df$ vanishes as a measure, so the two terms containing derivatives of $F$ cancel. This proves $\nabla_\nu T_\mu{}^\nu=F_\mu^A\mathcal E_A$ with the continuous representative of $F$. If the stress divergence vanishes, the polar decomposition $\mathcal E=e|\mathcal E|$ gives $F^Te=0$ $|\mathcal E|$-almost everywhere in these coordinates. At points of rank three, $F^T$ is injective. Since $|e|=1$ $|\mathcal E|$-almost everywhere, the open set of rank-three points has zero $|\mathcal E|$-measure. Thus the material equations hold throughout the domain whenever the rank-deficient set also has zero $|\mathcal E|$-measure. The converse and the equivalence with stationarity follow as in Proposition~\ref{prop:variational-EL-local04}.
\end{proof}

\begin{proof}[Proof of Proposition~\ref{prop-E-L}]
{For a locally Lipschitz configuration, localize a compactly supported variation to a material chart. The compact essential range lies inside the differentiability domain, so the extended action is defined along $f+\varepsilon\eta$ for small $|\varepsilon|$. The local bounds on the first constitutive derivatives justify differentiation under the integral. At a degenerate gradient these variations need not remain in the closed physical domain.} Differentiating the extended action, with $P_A{}^\mu=\partial\rphys/\partial F_\mu^A$, gives
\be
\delta\mathcal A[\eta]=\int\bigl(\rphys_{f^A}\eta^A+P_A{}^\mu\del_\mu\eta^A\bigr)\,dV_g.
\ee
Distributional integration by parts and a partition of unity give~\eqref{EL} and its converse. Proposition~\ref{prop:variational-EL-local04} gives the Noether identity {under the additional regularity assumption}.

For an auxiliary torsion-free material connection $\Gamma^B{}_{CA}$, the variation of $df$ is $D_\mu\eta^B=\del_\mu\eta^B+\Gamma^B{}_{CA}F_\mu^C\eta^A$. The horizontal material derivative and covariant divergence are
\bse
\be
\rphys_{\mid A}=\rphys_{f^A}-P_B{}^\mu\Gamma^B{}_{CA}F_\mu^C,
\ee
\be
D_\mu P_A{}^\mu=|g|^{-1/2}\del_\mu(|g|^{1/2}P_A{}^\mu)
-\Gamma^B{}_{CA}F_\mu^CP_B{}^\mu.
\ee
\ese
The connection terms cancel in $D_\mu P_A{}^\mu-\rphys_{\mid A}$, giving the same equation in bundle notation.
\end{proof}

%-----------------------------------------------------------------------------------------------------------------------------------------

\paragraph{Covariance and the stress--energy tensor.}
\label{append-B}

We now prove Proposition~\ref{prop:stru}. Spacetime covariance means that the Lagrangian density, viewed as a functional $\rphys[f;g]$, satisfies
\be
\rphys[{\phi^\star} f;{\phi^\star} g](x)
=
\rphys[f;g](\phi(x))
\ee
for every diffeomorphism $\phi:\Mcal\to\Mcal$. Let $\phi_t$ be a one–parameter group generated by an arbitrary vector field $\xi$ on $\Mcal$. Under the Lie derivative along $\xi$,
\be
\Lcal_\xi F_\lambda^A=\xi^\mu\del_\mu F_\lambda^A+F_\mu^A\del_\lambda\xi^\mu,
\qquad
(\Lcal_\xi g^{-1})^{\mu\nu}=\xi^\lambda\del_\lambda g^{\mu\nu}
-2g^{\lambda(\mu}\del_\lambda\xi^{\nu)}.
\ee
In the derivative of the covariance identity, the terms proportional to $\xi$ give $\xi^\mu\del_\mu\rphys$ by the chain rule and cancel. The remaining terms give
\be
\left(F_\mu^A\frac{\del\rphys}{\del F_\lambda^A}
 -2g^{\lambda\nu}\frac{\del\rphys}{\del g^{\mu\nu}}\right)\del_\lambda\xi^\mu=0.
\ee
Since $\del_\lambda\xi^\mu$ is arbitrary at a point, this proves~\eqref{equa-Lag}.

Conversely, we assume that $\rphys$ is $C^2$ in $(F,g)$ on its open positive-density constitutive domain, preserved by simultaneous changes of spacetime basis, and satisfies~\eqref{equa-Lag}. We write $\gamma\coloneqq(g_{\mu\nu})$ and regard $F$ as a $3\times4$ matrix. For constant symmetric matrices $H,L$, we define
\bse
\be
\mathcal X_H(\gamma,F)\coloneqq(2H,F\gamma^{-1}H).
\ee
The equation~\eqref{equa-Lag} gives $\mathcal X_H\rphys=0$, and direct differentiation gives
\be
[\mathcal X_H,\mathcal X_L](\gamma,F)
=(0,F[\gamma^{-1}L,\gamma^{-1}H]).
\ee
\ese
The matrices $\gamma^{-1}H$ are the $\gamma$-self-adjoint matrices. Their commutators span the $\gamma$-skew-adjoint matrices: in a diagonal frame for $\gamma$, we commute a diagonal matrix with distinct entries with self-adjoint matrices supported on each pair of indices. Since $\rphys$ is $C^2$, its differential therefore vanishes on every infinitesimal simultaneous change of basis. Integration along paths in $\mathrm{GL}^+(4)$ proves covariance under positive-determinant changes. Reflection $S$ in the timelike line $\ker F$ satisfies $FS=F$ and $S^T\gamma S=\gamma$. For a negative-determinant matrix $A$, the positive-determinant matrix $SA$ thus produces the same transformed pair $(FA,A^T\gamma A)$. This proves the converse where $n>0$; covariance of the finite-flux closure follows by continuity of the change of basis on its field and stress components.

\bse
We substitute the covariance identity into
\be
T_{\mu\nu} = 2 \, \frac{\del \rphys}{\del g^{\mu\nu}} - \rphys g_{\mu\nu} 
\ee
and obtain the mixed stress components
\be
T_\mu^{\ \nu}
= - \rphys \delta_\mu^{\ \nu} + F^A_\mu \frac{\del \rphys}{\del F^A_\nu},
\ee
that is, \eqref{eq:useful}. Contracting the mixed stress with $u^\mu$ and using $u^\mu F_\mu^A=0$ gives $u^\mu T_\mu{}^\nu=-\rphys u^\nu$. Symmetry of $T_{\mu\nu}$ gives the same eigenvector identity with either index raised. This completes the proof of Proposition~\ref{prop:stru}.
\ese

%-----------------------------------------------------------------------------------------------------------------------------

% Split source interval: app:observer-acoustic-proof

%============================================================================

\subsection{The corrected estimate for integrated stability}
\label{app:integrated-3d}

{
{We now give a proof of Theorem~\ref{uc:thm:integrated-3d}. We estimate} the constitutive source and the coefficients near zero reference density, then identify the limit at $t=0$ of the corrected integral evaluated on the initial data set.\par}

We differentiate the corrected integral $\mathcal E$, using $q=n_*^\theta$ and the reference fields $(Q,Z,Y)$ of Section~\ref{sec:uniform-elastic-fluid-existence}. The coefficient averaged at spatial radius proportional to $q^2$ will allow us to bound the divergence in terms of $E$ and $L$.

{\bf 1.} \emph{Estimate of the grouped constitutive source.}
We use the corrected current
\be
 \widetilde{\mathscr Q}^{\,\mu}
 =\mathscr Q^\mu+d_B^{\mu\nu}\Omega_{BAC}h^A\Delta F_\nu^C
       -\widetilde\beta_A\Delta {\particlecurrent}^\mu h^A.
\ee
Its time component integrates to $\mathcal H$. We introduce the source with the original coefficient $\beta$,
\be
 \mathcal R_0=-(\del_\nu w_{*\mu})\Delta T^{\mu\nu}
               -s_A^\mu\Delta F_\mu^A-\beta_A\Delta {\particlecurrent}^\mu\Delta F_\mu^A .
\ee
In the first correction, $\del_\mu d_B^{\mu\nu}=0$ and antisymmetry cancels the second derivative of $h$. The term $d_B^{\mu\nu}\Omega_{BAC}\Delta F_\mu^A\Delta F_\nu^C$ cancels the quadratic part of $-d_B^{\mu\nu}\Delta G_{\mu\nu B}$, leaving the contribution with $2d_B^{\mu\nu}\Omega_{BAC}F_{*\nu}^C$.  Both currents are conserved with their prescribed initial terms by the identity for Lipschitz configurations in Proposition~\ref{prop:kinematic-weak-local}. Using the stress and compatibility equations as in Lemma~\ref{thm:common-current-jhde}, we obtain
\be
\label{eq:integrated-relative-rate}
 \del_\mu\widetilde{\mathscr Q}^{\,\mu}
 {=} \mathcal R_0+(\beta_A-\widetilde\beta_A)\Delta {\particlecurrent}^\mu\Delta F_\mu^A
       -(\del_\mu\widetilde\beta_A)\Delta {\particlecurrent}^\mu h^A .
\ee
{The antisymmetric coefficient $d_B^{\mu\nu}$ is continuous up to the initial slice. Indeed, $C_*,D_*$ have the form $q^{p_\theta+2}$ times smooth functions of $(U,Y)$, with $p_\theta=1/\theta>0$. For a spatial or time derivative, the identity $Z=QY$ gives
\be
 q^{p_\theta+2}\del_\alpha Y
 =q^{p_\theta+1}\Gamma_*^{-\theta}
       (\del_\alpha Z-Y\del_\alpha Q),\qquad \alpha=0,1,2,3.
\ee
The right-hand side is continuous because $U\in C_tH^3\cap C_t^1H^2$ and $Y\in C_tH^2$ on the torus. The other terms in the first derivatives of $C_*,D_*$ have the same continuity and vanish at $q=0$. Thus these derivatives are continuous even when the derivatives of $Y$ are not.}
We mollify the antisymmetric coefficient and the potentials in the interior of the slab. The mollified coefficient has zero divergence there and converges locally uniformly. The gradients of the mollified potentials are uniformly bounded and converge strongly in every finite $L^p$ on compact subsets. This gives the first product rule in the interior. The first correction has time component $d_B^{0i}\Omega_{BAC}h^A\del_i h^C$, with $d_B^{0i}=\del_jD_{*B}^{j0i}$. As $t\downarrow0$, $d_B^{0i}h$ converges uniformly and $\del_i h(t)$ converges weak-$*$ to $\del_i h(0)$. This gives its product trace when the interior identity is tested with time cutoffs.

For the second product rule, we use $\widetilde\beta_Ah^A$ in the continuity equation by smooth approximation, uniform at $t=0$ and with derivatives converging in $L^1$. Its initial term uses the prescribed density difference.  For each normal flux component $P$ in $(r,F_i,G_{ij},T^{\mu0})$, its balance law expresses $\partial_t P$ as a spatial divergence of bounded fluxes. Testing in space and integrating the resulting scalar identity from the prescribed initial term gives
\be
 \|P(t)-P_{\rm in}\|_{(W^{1,1}(\mathbb T^3))^*}\le Ct
\ee
on a common full-measure set of times. The slice fluxes and the prescribed fluxes are bounded by the branch hypothesis, so density extends the convergence to every $L^1$ test. The time components of the original and corrected currents are affine in these normal fluxes: the entropy term itself is $-r$. All their coefficients, including the products with $h$, converge uniformly at $t=0$. Thus $E$ and $\mathcal H$ have the traces obtained by evaluating the same affine expressions with the prescribed initial data sets. For positive total particle number the Lipschitz time representatives of the positions give $|L(t)-L(0)|\le Ct$, with the prescribed initial positions. For zero total particle number $L=0$. Thus the comparison of $\mathcal E$ with $E+L$ extends to $t=0$.

We prove $|\mathcal R_0|\le C\mathscr D$ for $\gamma>1$. We use the constitutive derivatives $P_*,Q_*$ and the forces $L_*,K_*$ of~\eqref{eq:main-vacuum-reference-forces}. The constitutive formulas give $P_*,Q_*=\Ocal(n_*^\gamma)$. Their derivatives along $u_*$ have the same bound, since $D_t(q,z)$ and $D_tV_*$ are $\Ocal(q)$ by the normalized equations. Differentiation of the antisymmetric products \(C_*=w_*\wedge P_*\) and \(D_*=w_*\wedge Q_*\) gives 
\bse
\be
 c_A^\mu=-w_*^\mu L_{*A}+\Ocal(r_*^\gamma),\qquad
 d_B^{\mu\nu}=\tfrac12(-w_*^\mu K_*^{\nu B}+w_*^\nu K_*^{\mu B})+\Ocal(r_*^\gamma).
\ee
All terms in which the derivative acts on \(w_*\) have size \(\Ocal(r_*^\gamma)\) by its bounded first derivatives. Thus \(s_A^\mu\) is a multiple of \(u_*^\mu\) up to \(\Ocal(r_*^\gamma)\), since \(u_*^\nu F_{*\nu}=0\). Its component along the spatial slices after subtraction of \(\beta_A {\particlecurrent}_*^\mu\) satisfies 
\be
 t_A^\mu\coloneqq s_A^\mu-\beta_A{\particlecurrent}_*^\mu=\Ocal(r_*^\gamma),\qquad t_A^0=0 .
\ee
\ese

For comparable positive proper densities, the original remainder obeys \(|\mathscr R_*|\le C\mathscr D\) by the comparable-density argument in the proof of Theorem~\ref{uc:thm:common}, which does not require~\eqref{uc:eq:geometric-coefficients}. That argument uses the physical material equations of the reference and the Taylor remainder of its internal stress. By the quadratic minor identity, \(\mathcal R_0-\mathscr R_*\) equals \(d_B^{\mu\nu}\Omega_{BAC}\Delta F_\mu^A\Delta F_\nu^C-\beta_A\Delta {\particlecurrent}^\mu\Delta F_\mu^A\). The orthogonality condition~\eqref{eq:450bis-2} gives $|u_*\cdot\Delta F|\le C|V-V_*|$. The expansion of $d$ bounds its contribution by
\bse
\be
 C r_*^\gamma|\Delta F|^2+C r_* q |V-V_*||\Delta F|\le C\mathscr D.
\ee
In the observer's coordinates we have \(\Delta {\particlecurrent}^\mu\Delta F_\mu=(r-r_*)(\Delta F_0+V_*\cdot\Delta F_{\rm sp})+r(V-V_*)\cdot\Delta F_{\rm sp}\) for each material index. Its product with \(\beta\) is bounded by \(Cq|r-r_*||V-V_*|+Cqr|V-V_*||\Delta F|\). We put $R=r+r_*$ and use $q=n_*^\theta\simeq R^\theta$ on this region. We factor the first product using the density weights
\be
 q|r-r_*||V-V_*|
 \le C\bigl(R^{(\gamma-2)/2}|r-r_*|\bigr)
        \bigl(\sqrt r|V-V_*|\bigr),\qquad
 qr|V-V_*||\Delta F|
 \le C\bigl(R^{\gamma/2}|\Delta F|\bigr)
        \bigl(\sqrt r|V-V_*|\bigr).
\ee
Each squared factor is controlled by $\mathscr D$, by the density, kinetic and form bounds of Section~\ref{subsec:local-stability-proofs}. This proves the required estimates with constants independent of a positive lower density, including for $1<\gamma<2$.

In the separated region, $\min(n,n_*)\le\vartheta\max(n,n_*)$ with the ratio chosen before~\eqref{uc:eq:weighted-coercivity}. The observer densities then satisfy $\min(r,r_*)\le\kappa\max(r,r_*)$, where $\kappa=\Gamma_{\rm c}\vartheta<1$.  The separated-density estimate~\eqref{uc:eq:separated-density-bound} applies with $\lambda=\kappa$. The entropy bound~\eqref{uf:eq:local-relative-entropy-comparison} therefore controls $(r+r_*)^\gamma+r|V-V_*|^2$, including when one density vanishes. The two form terms in \(\mathcal R_0\) must be combined. The exact orthogonality condition~\eqref{eq:450bis-2} gives
\be
 -s_A^\mu\Delta F_\mu^A-\beta_A\Delta {\particlecurrent}^\mu\Delta F_\mu^A
   =-t_A^\mu\Delta F_\mu^A+\beta_A {\particlecurrent}^\mu F_{*\mu}^A .
\ee
Indeed, substituting $s_A=t_A+\beta_A{\particlecurrent}_*$ combines the two current terms into $-\beta_A{\particlecurrent}^\mu\Delta F_\mu^A$. Since ${\particlecurrent}^\mu F_\mu^A=0$, this leaves $\beta_A {\particlecurrent}^\mu F_{*\mu}^A$. In the observer frame $F_{*0}^A=-V_*^iF_{*i}^A$, so this contraction is $\beta_A r(V^i-V_*^i)F_{*i}^A$ wherever $r>0$ and vanishes at $r=0$. Together with the bound on $t_A$ this gives $C(r_*^\gamma+rq|V-V_*|)$. Young's inequality gives
\be
 rq|V-V_*|\le\tfrac12r|V-V_*|^2+Crq^2
 \le\tfrac12r|V-V_*|^2+C(r+r_*)^\gamma,
\ee
\ese
since \(q^2\le Cr_*^{\gamma-1}\). The entropy controls both terms on this region. This includes points where the weak solution has zero density and the strong solution has positive density. For the stress contraction, the internal stress is \(\Ocal(n^\gamma)\), and the mass stress calculation~\eqref{uc:eq:mass-far-states} bounds its contribution by \(C(r|V-V_*|^2+(r+r_*)^\gamma)\). This calculation uses only \(D_{u_*}u_*=\Ocal(q)\), \(dH_*=\Ocal(q)\) and \(D_{u_*}H_*=\Ocal(q^2)\), valid for the local references specified in Theorem~\ref{uc:thm:integrated-3d}. At reference vacuum the form coefficients and their first derivatives vanish. The free transport equation gives the mass stress bound at $q=0$. The grouped source vanishes where both solutions are empty.

\vskip.2cm \noindent {\bf 2.} \emph{Regularization of the coefficient at the density scale.}
In addition to $U=(Q,Z,V_*)\in C_tH^3\cap C_t^1H^2$, we have $Y\in C_tH^2\cap C_t^1H^1$ at $s=3$ in the local existence results. The smooth factors in $C_*=q^{p_\theta+2}\mathcal C(U,Y)$ and $D_*=q^{p_\theta+2}\mathcal D(U,Y)$ come from $z=qY_R$, with $Y_R$ the scaled rest minors, and the identity
\be
 \Psi_z(q,qY_R)=q\int_0^1\Psi_{zz}(q,\sigma qY_R)Y_R\,d\sigma .
\ee
The derivatives in \(c^0,d^{0i}\) are spatial by antisymmetry. Differentiating \(q^{p_\theta+2}\mathcal C\) and dividing by \(r_*=\Gamma_*q^{p_\theta}\) give, up to the bounded factor \(\Gamma_*^{-1}\), the terms \((p_\theta+2)q(\del_iq)\mathcal C\) and \(q^2\del_i\mathcal C\). The same calculation applies to \(\mathcal D\). With \(q\del_iY=\Gamma_*^{-\theta}(\del_i Z-Y\del_iQ)\), this yields \(\beta=q B(U,Y,\nabla U)\), where \(B\) is smooth in the fields and linear in the last argument. This formula defines \(B\) also at \(q=0\). Sobolev products give \(B\in C_tH^2\cap C_t^1H^1\).  The smoothing below uses its uniform spatial $W^{1,6}$ bound, which also gives the $C^{0,1/2}$ bound.

We set \(D_t=\del_t+V_*\cdot\nabla\). The normalized equations give \(D_tU=q A^i(U,Y)\del_i U\) with smooth coefficients, and \(D_tY\) is algebraic in \(Y,\nabla V_*\). The density row of~\eqref{uf:eq:normalized-evolution} gives \(D_tQ=-\theta Q\,\operatorname{div}V_*\) on the full slab, and the velocity equation bounds \(D_t\Gamma_*\). Thus \(|D_tq|\le Cq\) by \(q=Q\Gamma_*^{-\theta}\). The terms with \(\nabla Y\) are bounded because
\be
 q\nabla Y=\Gamma_*^{-\theta}(\nabla Z-Y\nabla Q).
\ee
{Differentiating the normalized equation gives
\be
 D_t(\del_iU)=qA^j\del_i\del_jU+(\del_iq)A^j\del_jU
 +q(\del_iA^j)\del_jU-(\del_iV_*)\cdot\nabla U.
\ee
The product $q\nabla A^j$ is bounded by the identity for $q\nabla Y$ above. Since $D_tY$ is bounded and $D^2U\in L^\infty_tL^6_x$, the chain rule gives $D_tB=B_0+qB_1$, with $B_0\in L^\infty$ and $B_1\in L^\infty_tL^6_x$. The calculation holds almost everywhere, including where $q=0$, without division by $q$.} 

We choose a nonnegative smooth unit-integral function \(\varphi\) supported in the unit ball and put
\bse
\be
 S_\ell B(t,x)\coloneqq\int\varphi(y)B(t,x+\ell y)\,dy,
 \qquad \ell_\epsilon\coloneqq c(q+\epsilon)^2,\qquad
 \beta_\epsilon\coloneqq q S_{\ell_\epsilon}B .
\ee
{We choose $c>0$ so that every averaging ball lies in periodic coordinates. At radius $\ell$, the $C^{0,1/2}$ and $W^{1,6}$ bounds give an error $Cq\ell^{1/2}$ and a derivative term $Cq\ell^{-1/2}$. Their balance at $\ell=cq^2$ gives an error of order $q^2$ and a bounded derivative. With $5/2<s<3$, the corresponding estimates have $\alpha=s-5/2<1/2$. They would require $1/\alpha\le\kappa\le1/(1-\alpha)$ for $\ell=cq^\kappa$, an empty interval. Here the $L^6$ averaging factor is $C\ell_\epsilon^{-1/2}$. Since $|\nabla\ell_\epsilon|\le C(q+\epsilon)$, the spatial chain rule gives}
\be
 |\nabla\beta_\epsilon|
 \le C|\nabla q|\,\|B\|_{L^\infty}
 +Cq\bigl(1+|\nabla\ell_\epsilon|\bigr)\ell_\epsilon^{-1/2}
       \|\nabla B\|_{L^6}
 \le C.
\ee
{Here $q\ell_\epsilon^{-1/2}\le c^{-1/2}$ also at $q=0$.} In the material derivative of the average, one term is \(S_{\ell_\epsilon}D_tB\). On the averaging ball, \(q(t,x+\ell_\epsilon y)\le C(q(t,x)+\epsilon)\) by the Lipschitz bound of \(q\). This gives the uniform bound
\be
 |S_{\ell_\epsilon}(qB_1)|
 \le C(q+\epsilon)\ell_\epsilon^{-1/2}\|B_1(t)\|_{L^6}
 \le C\|B_1(t)\|_{L^6}.
\ee The other term is
\be
 \int\varphi(y)\bigl(V_*(t,x)-V_*(t,x+\ell_\epsilon y)
          +(D_t\ell_\epsilon)y\bigr)\cdot\nabla B(t,x+\ell_\epsilon y)\,dy .
\ee
The factor in parentheses is \(\Ocal(\ell_\epsilon)\) because \(V_*\) is Lipschitz and \(|D_t\ell_\epsilon|\le C\ell_\epsilon\). The same averaging estimate bounds this term. For fixed positive \(\epsilon\), we first differentiate \(S_\ell B(t,x)\) with \(t,x,\ell\) as independent variables and \(\ell>0\), then substitute \(\ell=\ell_\epsilon(t,x)\). Convolution at positive radius permits these derivatives with the regularity of \(B\) above. We obtain \(|D_t\beta_\epsilon|\le Cq\), uniformly as \(\epsilon\downarrow0\).

Set $\widetilde\beta=qS_{cq^2}B$ for $q>0$ and zero for $q=0$. The spatial H\"older bound gives $|\beta_\epsilon-\widetilde\beta|\le Cq|\ell_\epsilon-cq^2|^{1/2}\le C\epsilon^{1/2}$ on the closed slab. The uniform derivative bounds give Lipschitz continuity up to $t=0$, where the same averaging formula applies. Testing the inequalities for $D_t\beta_\epsilon$ against nonnegative $L^1$ functions passes the bounds to the weak-star limit. The H\"older estimate at radius $cq^2$ gives the second of the resulting bounds.
\be
\label{eq:regularized-material-coefficient}
 |\widetilde\beta|\le Cq,\qquad |\beta-\widetilde\beta|\le Cq^2,\qquad
 |\nabla\widetilde\beta|\le C,\qquad |D_t\widetilde\beta|\le Cq .
\ee
\ese

\vskip.2cm \noindent {\bf 3.} \emph{Relative entropy and displacement estimates.}
{
The algebraic error in~\eqref{eq:integrated-relative-rate} is bounded by $Cq^2|\Delta\particlecurrent||\Delta F|$. For $R=r+r_*>0$, the weighted estimates in Section~\ref{subsec:local-stability-proofs} bound $R^{(\gamma-2)/2}|\Delta\particlecurrent|$ and $R^{\gamma/2}|\Delta F|$ by $C\sqrt{\mathscr D}$. The product of their weights is $R^{\gamma-1}\ge c q^2$, so the error is at most $C\mathscr D$. It vanishes at $R=0$.
}

To estimate the term with the differentiated coefficient, we decompose the current difference as $\Delta {\particlecurrent}=(\Delta r)(1,V_*)+(0,r(V-V_*))$. The spatial factor $r(V-V_*)={\particlecurrent}^{\rm sp}-rV_*$ is defined on the whole slab, including where the weak solution has zero density. Hence the term is
\bse
\be
 -h^A\bigl((\Delta r)D_t\widetilde\beta_A+r(V-V_*)\cdot\nabla\widetilde\beta_A\bigr).
\ee
Since $q^2\le C(r+r_*)^{\gamma-1}$, the density bound gives $q|\Delta r|\le C\sqrt{\mathscr D}\sqrt{r+r_*}$. The kinetic bound gives $\sqrt r\,|V-V_*|\le C\sqrt{\mathscr D}$. Consequently the absolute value of this term is at most $C\sqrt{\mathscr D}\sqrt{r+r_*}\,|h|$, including when either density vanishes. Cauchy--Schwarz and the displacement estimate of Section~\ref{subsec:local-stability-proofs} bound its integral by $CE^{1/2}L^{1/2}$. Together with the bound for $\mathcal R_0$ and $|L'|\le C(E+L)$ this gives
\be
 \mathcal E'\le C(E+L)\le C'\mathcal E
\ee
\ese
The uncorrected current also has a divergence equality with bounded right-hand side, by conservation of the exact-minor currents and the product rule of Lemma~\ref{thm:common-current-jhde}. After spatial integration this gives a Lipschitz representative of $E$ with trace $E(0)$ by the prescribed initial term. For the corrected current,~\eqref{eq:integrated-relative-rate} and the bound for $L'$ give $|\mathcal E'|\le C(E+L)$. Thus $\mathcal E$ is Lipschitz with the trace computed in Step~1, and Gronwall's inequality gives the comparison. Parts~(ii) and~(iii) imply (iv).
This completes the proof of Theorem~\ref{uc:thm:integrated-3d}.

\subsection{Relative entropy identities and stress at zero current}
\label{subsec:constitutive-stability-verification}

\paragraph{Identities for the relative entropy current.}

{This subsection gives the algebraic identities used in Sections~\ref{subsec:local-stability-proofs} and~\ref{sec:elastic-constitutive-estimates}.} In computing the divergence, we use the exterior and material equations to cancel the terms linear in the minor increments, then estimate the density terms.
For the multipliers in~\eqref{eq:main-vacuum-finite-multipliers}, the source in the  inequality for the relative entropy current is
\be
\label{eq:physical-nonlinear-original-remainder}
\mathscr R=-(\nabla_\nu w_{*\mu})\Delta T^{\mu\nu}
-(\nabla_\nu C_{*A}^{\nu\alpha})\Delta F_\alpha^A
-(\nabla_\nu D_{*A}^{\nu\alpha\beta})\Delta G_{\alpha\beta A}.
\ee
In an orthonormal reference rest frame, write $q_A\coloneqq F_0^A$, $\Fsp _{Ai}\coloneqq F_i^A$, $\Gsp \coloneqq\operatorname{cof}\Fsp $ and $\Gsp _*=\operatorname{cof}\Fsp _*$. The exact minors and the orthogonality condition~\eqref{eq:450bis-2} for the reference give
\be
\label{eq:physical-minor-time-identities}
G_{0ib}=\Omega_{bAC}q_A(\Fsp)_{Ci},\qquad
{\particlecurrent}^i=-(\Gsp)_{Ai}q_A,\qquad
F_{*0}=G_{*0i}=0.
\ee
With the forces and acceleration from~\eqref{eq:main-vacuum-reference-forces} and~\eqref{eq:main-vacuum-reference-acceleration}, the material equation~\eqref{eq:physical-material-force-equation} becomes $L_{*A}+\Omega_{bAC}(\Fsp_*)_{Ci}K_*^{ib}+(\Gsp_*)_{Ai}\mathcal A_{*i}=0$.

For a smooth constitutive energy $\mathcal L(F,G)$ independent of the current, its stress is
\be
T_{\mathcal L}^{\mu\nu}
=F^{A\mu}P_A^\nu+2G^{\mu\lambda A}Q^\nu{}_{\lambda A}
-\mathcal Lg^{\mu\nu},
\qquad P\coloneqq\del_F\mathcal L,\quad Q\coloneqq\del_G\mathcal L,
\ee
with full antisymmetric sums. Closure gives $\nabla_\nu T_{\mathcal L}^{\mu\nu}=F^{A\mu}L_A+G^{\mu\lambda A}K_{\lambda A}$. To display the mixed-form cancellation, we evaluate every undotted field in the next identity at the reference and let a dot denote an arbitrary variation of the form components. Differentiation of $C=w\wedge P$ and $D=w\wedge Q$, followed by Cartan's identities and the exterior closure equations for the reference, gives
\be
\label{eq:general-vacuum-joint-Noether}
\aligned
&-\nabla_\nu w_\mu\,\dot T_{\mathcal L}^{\mu\nu}-(\nabla_\nu C_A^{\nu\alpha})\dot F_\alpha^A-(\nabla_\nu D_A^{\nu\alpha\beta})\dot G_{\alpha\beta A}-w^\alpha\dot F_\alpha^AL_A-w^\alpha\dot G_{\alpha\lambda A}K^{\lambda A}\\
&=w^\mu\bigl((\nabla_\mu F_\alpha^A)\dot P_A^\alpha-(\nabla_\mu P_A^\alpha)\dot F_\alpha^A+(\nabla_\mu G_{\alpha\beta A})\dot Q^{\alpha\beta A}-(\nabla_\mu Q^{\alpha\beta A})\dot G_{\alpha\beta A}\bigr)=0.
\endaligned
\ee
The material tensors are fixed, so the last equality follows from the symmetry of the joint Hessian of $\mathcal L$ in $(F,G)$ and $\nabla g=0$. To check the two-form terms, antisymmetry and full sums give
\be
 (\nabla_\nu D^{\nu\alpha\beta})\dot G_{\alpha\beta}
 =(\nabla_\nu w^\nu)Q^{\alpha\beta}\dot G_{\alpha\beta}
 +w^\nu(\nabla_\nu Q^{\alpha\beta})\dot G_{\alpha\beta}
 -2(\nabla_\nu w^\alpha)Q^{\nu\beta}\dot G_{\alpha\beta}
 -w^\alpha K^\beta\dot G_{\alpha\beta}.
\ee
The third term cancels the corresponding variation of the stress. Differentiating $i_wG=0$ and using $dG=0$ gives $2w^\mu(\nabla_\nu G_{\mu\lambda})\dot Q^{\nu\lambda}=w^\mu(\nabla_\mu G_{\nu\lambda})\dot Q^{\nu\lambda}$, which yields the remaining terms in the equation above.

For a smooth constitutive term $L$ depending on the original fields $(F,G,\particlecurrent)$, put $B_L=\del_{\particlecurrent} L$ and $\alpha_*=w_*^\mu B_{L*\mu}$. For its conjugates and their divergences use~\eqref{eq:main-vacuum-reference-forces}; the multipliers $C_*,D_*$ are formed using the same $w_*$. Including the current in the linear variation gives
\be
\label{eq:physical-nonlinear-current-joint-Noether}
\aligned
&-\nabla_\nu w_{*\mu}\,\dot T_L^{\mu\nu}
-(\nabla_\nu C_*^{\nu\alpha})\dot F_\alpha
-(\nabla_\nu D_*^{\nu\alpha\beta})\dot G_{\alpha\beta}
-(\del_\nu\alpha_*)\dot {\particlecurrent}^\nu\\
&=w_*^\mu\dot F_\mu L_*
+w_*^\mu\dot G_{\mu\lambda}K_*^\lambda
+w_*^\mu\dot {\particlecurrent}^\nu(dB_{L*})_{\mu\nu}.
\endaligned
\ee
{We sum over material indices. The exterior closure equations and $i_{u_*}F_*=i_{u_*}G_*=0$ give the form terms as above. For the current terms, ${\particlecurrent}_*=n_*H_*w_*$ and $\nabla_\nu {\particlecurrent}_*^\nu=0$ imply
\be
 w_*^\nu\nabla_\nu {\particlecurrent}_*^\mu
 ={\particlecurrent}_*^\nu\nabla_\nu w_*^\mu-{\particlecurrent}_*^\mu\nabla_\nu w_*^\nu.
\ee
This identity requires no division by $n_*$. The current conservation is the exterior equation for its dual three-form. More explicitly, after subtracting the current-force term on the right-hand side of~\eqref{eq:physical-nonlinear-current-joint-Noether}, the current contributions on the left reduce to
\be
\begin{aligned}
 &({\particlecurrent}_*^\nu\nabla_\nu w_*^\mu-{\particlecurrent}_*^\mu\nabla_\nu w_*^\nu)\dot B_{L*\mu}
   -w_*^\nu(\nabla_\nu B_{L*\mu})\dot {\particlecurrent}^\mu\\
 &\qquad=w_*^\nu\bigl((\nabla_\nu {\particlecurrent}_*^\mu)\dot B_{L*\mu}
                  -(\nabla_\nu B_{L*\mu})\dot {\particlecurrent}^\mu\bigr).
\end{aligned}
\ee
The derivatives of $w_*\cdot B_{L*}$ cancel those multiplying the undifferentiated $B_{L*}$ in the stress variation. Adding the corresponding form contributions displayed above leaves the pairing of the joint Hessian of $L$ with the two variations $\nabla_{w_*}(F_*,G_*,{\particlecurrent}_*)$ and $(\dot F,\dot G,\dot {\particlecurrent})$, in opposite orders. Its symmetry makes their difference zero and proves~\eqref{eq:physical-nonlinear-current-joint-Noether}. If its coefficients and the stress differential extend continuously to a selected branch on the constitutive boundary, the identity holds there by continuity.} With $u_*$ in place of $w_*$ and the corresponding $C_*,D_*,\alpha_*$, writing $B_{L*}=-h_{L*}u_{*\flat}$ gives $u_*^\mu(dB_{L*})_{\mu\nu}=-h_{L*}u_*^\mu\nabla_\mu u_{*\nu}-\Pi_{*\nu}{}^\lambda\del_\lambda h_{L*}$.

For the polytropic density term with $1<\gamma\le2$ and positive rest mass $m>0$, put $r\coloneqq n_*$, $h(n)\coloneqq m+\kappa n^{\gamma-1}$, $H=H_*=h(r)$ and $\mathscr Q_{\rm fl}\coloneqq-\Delta {\particlecurrent}-u_*\cdot\Delta T_{\rm fl}/H$.  Put $\chi=-g(u,u_*)$, $a=\kappa/\gamma$, $\rho_{\rm d}(n)={m \, n}+a n^\gamma$, $\theta_*=\nabla_\nu u_*^\nu$ and $c_r^2=r h'(r)/H$, continued by zero at $r=0$. In the reference rest frame, current conservation gives $\nabla_0 H=-Hc_r^2\theta_*$. This holds almost everywhere at $r=0$ as well, since the nonnegative Lipschitz function $H-m$ has zero derivative there. The mixed derivatives of $w_*=u_*^\flat/H$ satisfy $\nabla_i w_{*0}+\nabla_0 w_{*i}=\mathcal A_{*i}/H^2$. Substitution into the stress contraction gives:
\be
\label{eq:forced-fluid-remainder}
-\nabla_\nu w_{*\mu}\,\Delta T_{\rm fl}^{\mu\nu}
=\mathscr R_{\rm fl}^{0}
-\frac{nh}{H^2}(u-u_*)^\mu\mathcal A_{*\mu},
\qquad |\mathscr R_{\rm fl}^{0}|\le C\mathscr Q_{\rm fl}^0.
\ee
{Write $\Delta u=u-u_*$, $D_*H=u_*^\nu\nabla_\nu H$, and $\Delta p=p-p_*$. To see this separation, use $u_*^\mu\nabla_\nu u_{*\mu}=0$  in the contraction of the stress difference. The terms linear in $\Delta u$ from differentiating $u_*/H$ combine into $-nh\Delta u\cdot( H\nabla_*u_*+dH)/H^2$. Since $H\nabla_*u_*+dH=\mathcal A_*-u_*^\flat D_*H$, the part containing $\mathcal A_*$ is the term displayed above.  The remaining terms involving time derivatives combine with it to give $-(nh(2\chi-1)-rH-\Delta p)D_*H/H^2$. Using $nh=\rho_{\rm d}+p$ and the reference continuity identity, we obtain
\be
\mathscr R_{\rm fl}^{0}
=-\frac{nh}{H}\Delta u^\mu\Delta u^\nu\nabla_\nu u_{*\mu}
-\frac{nh}{H^2}(\chi-1)\Delta u\cdot dH
+\frac{\theta_*}{H}\left[-\frac{(\gamma-1)m}{H}B_{\rho_{\rm d}}(n\mid r)
+2c_r^2nh(\chi-1)\right].
\ee
The pressure term follows from $p-p_*-c_r^2(\rho_{\rm d}-\rho_{{\rm d}*})=(\gamma-1)m B_{\rho_{\rm d}}/H$. To bound the remainder, put $X=n\chi$ and $B_\gamma(x\mid y)=x^\gamma-y^\gamma-\gamma y^{\gamma-1}(x-y)$. Direct expansion gives
\be
H\mathscr Q_{\rm fl}^0
={m \, n}\chi(\chi-1)+aB_\gamma(X\mid r)+a n^\gamma\bigl(\gamma\chi^2-(\gamma-1)-\chi^\gamma\bigr).
\ee
The last bracket bounds $(\gamma-1)(\chi^2-1)$ since $\chi^\gamma\le\chi^2$. By homogeneity, $B_\gamma(x\mid y)$ is comparable to $|x^{\gamma/2}-y^{\gamma/2}|^2$: the ratio is positive and continuous after normalizing by $x+y$, with a positive limit at equal arguments by Taylor expansion. The triangle inequality and $(\chi^{\gamma/2}-1)^2\le\chi^2-1$ thus bound $B_{\rho_{\rm d}}(n\mid r)$ by $C\mathscr Q_{\rm fl}^0$. Finally $|\Delta u|_*^2=2\chi(\chi-1)$ and $(\chi-1)|\Delta u|_*\le\sqrt2\chi(\chi-1)$ give the remainder bound.} The force contribution is
\be
\mathscr R_{\rm force}
=\frac{q_AL_{*A}+G_{0ib}K_*^{ib}}H
-\frac h{H^2}{\particlecurrent}^i\mathcal A_{*i}.
\ee
The polynomial minor relations and material equation give
\be
\label{eq:physical-nonlinear-additive-force-cancellation}
\mathscr R_{\rm force}
=\frac{H-h}{H^2}{\particlecurrent}^i\mathcal A_{*i}
+\frac1Hq_A(\Delta\Gsp)_{Ai}\mathcal A_{*i}
+\frac1Hq_A\Omega_{bAC}(\Delta\Fsp)_{Ci}K_*^{ib}.
\ee
The same cancellation applies to the spectral law by passing to the continuous conjugates and their time-integrable divergence bounds in the proof of Theorem~\ref{og:theorem}.

The scalar stress differential is also used at zero current. Write $\rrestnew(n)={m \, n}+\kappa n^\gamma/\gamma$ with $1<\gamma\le2$, $m,\kappa\ge0$, $m+\kappa>0$, $h=\rrestnew'$, $p=nh-\rrestnew$, $T_{\rm f}=nh\,u\otimes u+pg^{-1}$, and $h_r=h(r)$. Differentiation in a current increment $j$ gives
\be
\label{eq:nonlinear-fluid-affine-jet}
A_{{\rm f}*}j=h_r\bigl(u_*\otimes j+j\otimes u_*+g(u_*,j)u_*\otimes u_*\bigr)
-rh'_r g(u_*,j)\bigl(u_*\otimes u_*+g^{-1}\bigr).
\ee
The coefficients have finite limits $h_r\to m$ and $rh_r'\to0$ as $r\downarrow0$. For $r_{\rm f}\coloneqq\Delta T_{\rm f}-A_{{\rm f}*}\Delta {\particlecurrent}$, write $u=(\chi,b)$ in the reference rest frame and define the scaled scalar density $Q_{\rm f}\coloneqq-h_r\Delta {\particlecurrent}^0-u_{*\mu}\Delta T_{\rm f}^{\mu0}$. When $h_r>0$, this is $h_r\mathscr Q_{\rm fl}^0$ for $\mathscr Q_{\rm fl}=-\Delta {\particlecurrent}-u_*\cdot\Delta T_{\rm f}/h_r$. Put $c_r^2\coloneqq rh_r'/h_r$ where $h_r>0$, and take its limit at $r=0$ when $h_r=0$. Direct substitution gives
\bse
\be
\label{eq:nonlinear-fluid-jet-remainder}
r_{\rm f}^{00}=Q_{\rm f},
\qquad r_{\rm f}^{0i}=\bigl(h\chi-h_r\bigr){\particlecurrent}^i,
\ee
\be
r_{\rm f}^{ij}=nh\,b_i b_j+
\bigl[P_{\rm rel}+c_r^2\bigl(Q_{\rm f}-nh(\chi^2-1)\bigr)\bigr]\delta_{ij}.
\ee
\ese
{For any $\gamma>1$, put $B_{\rrestnew}(n\mid r)=\rrestnew(n)-\rrestnew(r)-h_r(n-r)$, and define $B_p$ by the same first-order subtraction for $p$. Since $p=(\gamma-1)(\rrestnew-{m \, n})$ and $p'(r)=c_r^2 h_r$, the pressure remainder is
\be
\label{eq:polytropic-pressure-remainder}
 P_{\rm rel}=p-p(r)-c_r^2\bigl(\rrestnew-\rrestnew(r)\bigr)
 =B_p-c_r^2 B_{\rrestnew}
 =\frac{(\gamma-1)m}{h_r} B_{\rrestnew}(n\mid r).
\ee
Here the last expression applies for $h_r>0$, in particular at $r=0$ if $m>0$.} The preceding calculation of the polytropic contribution gives $B_{\rrestnew}(n\mid r)+nh\chi(\chi-1)\le C\mathscr Q_{\rm fl}^0$. Homogeneity gives $n|h-h_r|^2\le C B_{\rrestnew}(n\mid r)$ {on bounded density ranges}, by normalizing with $\max(n,r)$ and taking the finite limit of the ratio at $n=r>0$. To estimate the mixed component, split it as $(h-h_r){\particlecurrent}^i+h(\chi-1){\particlecurrent}^i$ and use $|{\particlecurrent}^{\rm sp}|=n\sqrt{\chi^2-1}$. These estimates bound $r_{\rm f}$ by $C\mathscr Q_{\rm fl}^0$ when $m>0$.

%-----------------------------------------------------------------------------------------------------------------------------------------

\paragraph{Stress and entropy at zero current.}
\label{app:vacuum-stress-subgradient}

{We now give a proof of Proposition~\ref{prop:vacuum-stress-subgradient}.}

{\bf 1.} \emph{Finite fluxes at zero current.} We first consider a positive-density sequence tending to a zero-current state. We write $z_j\coloneqq(1,v_j)$, $q_j\coloneqq\sqrt{1-|v_j|^2}$ and $n_j\coloneqq r_jq_j$. Its mass and scalar fluid stress tensors are
\bse
\be
T_{m,j}=D_jz_j\otimes z_j,\qquad D_j=\frac{mr_j}{q_j},\qquad
T_{\mathcal U,j}
=\frac{\mathcal U'(n_j)}mT_{m,j}
+\bigl(n_j\mathcal U'(n_j)-\mathcal U(n_j)\bigr)g^{-1}.
\ee
{The stress $T_{\mathcal L,j}$ converges by its assumed continuity}. The scalar fluid contribution to the observer energy is
\be
T_{\mathcal U,j}^{00}
=n_j\mathcal U'(n_j)\bigl((u_j^0)^2-1\bigr)+\mathcal U(n_j)\ge0.
\ee
\ese
Finite total energy thus bounds $D_j$. Since $n_j\to0$ and $\mathcal U'(0)=0$, the scalar stress tends to zero. Passing to subsequences gives $v_j\to v$, $D_j\to D\ge0$, and~\eqref{eq:vacuum-dust-graph}. The orthogonality condition~\eqref{eq:450bis-2} passes to $i_zF=i_zG=0$. If $|v|<1$, then $D_j\to0$ because $r_j\to0$.

Conversely, we fix the spatial forms and reconstruct their time components using timelike $v_j\to v$. They converge to the prescribed forms because those annihilate $(1,v)$. For $|v|<1$, we take $v_j=v$ and any $r_j\downarrow0$. For $|v|=1$ and $D>0$, we take $r_j=(D/m)q_j$; for $D=0$, we take $r_j=q_j^2$. These choices give ${\particlecurrent}_j\to0$, $mr_j/q_j\to D$ and {convergence of $T_{\mathcal L}$}. The estimate for the scalar stress term then proves convergence of every component of the prescribed fields and stress. The construction concerns independent augmented fields; it imposes no polynomial-minor relations between them.

\vskip.2cm \noindent {\bf 2.} \emph{Finite supporting covectors at zero current.} We now examine supporting planes at $V_0$. If $v_0$ is an interior maximizer at zero current, then $P_0=-\ell_v(z_{{\rm s}0},0,v_0)$. States with fixed $(z_{{\rm s}0},v_0)$ and $r\downarrow0$ tend to $V_0$, while
\be
\widetilde\rho_r\longrightarrow m q(v_0)>0,\qquad
D\eta=\frac{(\widetilde\rho_{z_{\rm s}},-1,v)}{\widetilde\rho_r}
\longrightarrow\lambda_0.
\ee
Passing the supporting inequality at positive density to the limit proves~\eqref{eq:vacuum-entropy-subgradient}. {Since $\eta(V_0)=0$ and $\eta\le0$ on $\overline\Omega$, the supporting inequality gives the pointwise Lipschitz bound. Along the same family with $(z_{{\rm s}0},v_0)$ fixed, the flux formulas give $|V-V_0|=\Ocal(r)$, while $-\eta(V)=r$. Thus no exponent larger than one can hold at $V_0$.}

{By the right derivative of $\Phi$ computed before the proposition, $\Phi(z_{{\rm s}0},r,P_0)-\Phi_0(z_{{\rm s}0},P_0)=m q_{\max}r+o(r)$, where $q_{\max}$ is the largest value of $\sqrt{1-|v|^2}$ among the zero-current maximizers. If $q_{\max}=0$, the inward slope of $\eta$ is $-\infty$, which precludes a finite supporting covector.} {For $q_{\max}>0$, evaluating} any support on the vertical ray gives $\lambda_E\le-1/(m q_{\max})$. Formula~\eqref{eq:vacuum-entropy-subgradient} attains this bound at a maximizer with $q(v_0)=q_{\max}$.  More generally, for $a>0$,
\be
 (\lambda_z,-a,\lambda_P)\in\partial\eta(V_0)
 \quad\Longleftrightarrow\quad
 (\lambda_z/a,1/a,\lambda_P/a)\in\partial\Phi(z_{{\rm s}0},0,P_0),
\ee
where the support inequality for $\Phi$ is taken on $r\ge0$. Indeed, we substitute $E=\Phi(z_{\rm s},r,P)$ and $\eta=-r$ into the support inequality for $\eta$ and divide by $a$. These points cover the domain of finite entropy. Every supporting covector of $\eta$ has a negative energy component by the vertical slope bound.

Finally, a state with $D>0$ in~\eqref{eq:vacuum-dust-graph} has $|v|=1$ and $P_0=-\ell_v(z_{\rm s},0,v)+Dv$. The concave objective at zero current $f(w)=\ell(z_{\rm s},0,w)+w\cdot P_0$ consequently satisfies
\be
f(w)\le f(v)+Dv\cdot(w-v)
=f(v)-D(1-v\cdot w)<f(v)\qquad(|w|<1).
\ee
There is no interior maximizer, and hence no finite entropy subgradient. {We next determine the entropy asymptotics. We write $q(w)=\sqrt{1-|w|^2}$ and put $\delta_r=\Phi(z_{\rm s},r,P_0)-\Phi_0(z_{\rm s},P_0)$. {The quadratic density estimate in the proof of Theorem~\ref{thm:closed-vacuum-entropy} gives $|\ell(z_{\rm s},r,w)-\ell(z_{\rm s},0,w)|\le Cr^2q(w)^2$ uniformly in $w$}, and $\mathcal U(rq)=o(rq)$ uniformly as $r\downarrow0$. The concavity bound above and $1-v\cdot w\ge q(w)^2/2$ therefore give
\be
 \delta_r\le\sup_{0\le q\le1}\left\{-\left(\frac D2-Cr^2\right)q^2+(m+o(1))rq\right\}
 \le\frac{m^2}{2D}r^2+o(r^2).
\ee
For the reverse bound choose $q=mr/D$ and $w=\sqrt{1-q^2}\,v$. Taylor's formula at $v$, where $D_w f(v)=Dv$, gives $f(w)-f(v)=-Dq^2/2+\Ocal(q^4)$. Since $\mathcal U\ge0$, this choice gives $\delta_r\ge m^2r^2/(2D)+o(r^2)$. Inverting the resulting relation proves the displayed entropy asymptotics. The general $1/2$-H\"older bound of Corollary~\ref{cor:boundary-entropy-rates} and this asymptotic show that the exponent is optimal.} These zero-current states with nonzero rank-one null stress remain admissible fields and constitutive stress with finite fluxes.

This completes the proof of Proposition~\ref{prop:vacuum-stress-subgradient}.

%============================================================================

\paragraph{Acknowledgments}

The questions addressed in this paper have been a longstanding interest of mine; over the years I have returned to them repeatedly and written many (unpublished) notes. I am grateful to Bruno Le Floch (Sorbonne) for his many insightful observations that helped finalize this work. For numerous discussions and encouragement over many years, I am    very grateful to Costas M. Dafermos (Brown). I also thank Bobby Beig (Vienna) for many discussions which led to our joint paper~\cite{BeigLeFloch} on a class of generalized relativistic fluids.  Part of this research was done during yearly visits as a research fellow at the Courant Institute of Mathematical Sciences, New York University.  

\paragraph{Funding.} During the completion of this project, I was partially supported by the research project ANR-23-CE40-0010-02: Einstein-PPF, funded by the Agence Nationale de la Recherche, and by the MSCA Staff Exchange Project 101131233: Einstein-Waves, funded by the European Union.

\paragraph{Note on AI use}
After the first version was written, I used ChatGPT~6 to improve its presentation and produce all of the figures. This was followed by a human-only review to ensure mathematical validity and clarity.

\subsection*{Statements and declarations}

The author has no conflict of interest. There is no data associated with this research.

\end{document}